%% file: main.tex
\documentclass[10pt]{amsproc}
\usepackage[utf8]{inputenc}
\input{header}
\title[$\ku$-theoretic spectral decompositions for $S^n$ and $\PP^n$]{$\ku$-theoretic spectral decompositions for spheres and projective spaces}
\author{S. K. Devalapurkar}
\address{1 Oxford St, Cambridge, MA 02139}
\email{sdevalapurkar@math.harvard.edu, \today}
\thanks{Part of this work was done when the author was supported by NSF DGE-2140743. \newline 
The present article is far from being free of typos/errors. Any feedback would be greatly appreciated! I'll post major updates to the arXiv, but I'll upload minor edits to my website at \url{https://sanathdevalapurkar.github.io/files/hyperboloid_spectral_decomp.pdf}; so please see there for the most up-to-date version.}

\begin{document}

\maketitle

\begin{abstract}
    Ben-Zvi--Sakellaridis--Venkatesh described a conjectural extension of the geometric Satake equivalence to spherical varieties, whose spectral decomposition is described by Hamiltonian varieties. The goal of this article is to study their conjecture, especially in the case of spherical varieties of relative rank $1$, using tools from homotopy theory. Our discussion relates their conjecture to classical topics in homotopy theory such as the EHP sequence and Hopf fibrations, as well as more modern topics such as Hochschild (co)homology. We will also study an analogue of the derived geometric Satake equivalence and of the Ben-Zvi--Sakellaridis--Venkatesh conjecture with coefficients in connective complex K-theory. In this generalized setting, the dual \textit{group} (\`a la Langlands, Gaitsgory--Nadler, Sakellaridis--Venkatesh, Knop--Schalke) remains unchanged, but the specific dual ``representation'' of the dual group changes. On the spectral/Langlands dual side, we expect that the appropriate replacement of Hamiltonian varieties are given by what we term ``$\ku$-Hamiltonian varieties''; this is a notion interpolating between Hamiltonian and quasi-Hamiltonian varieties (\`a la Alekseev--Malkin--Meinrenken).
    Finally, we suggest possible generalizations to more exotic cohomology theories such as complex cobordism.
\end{abstract}

\tableofcontents

\newpage

\section{Introduction}
\input{intro/intro}
\input{intro/notation}
\input{intro/thanks}

\newpage

\section{Equivariant connective K-theory}\label{sec: ku big section}
\input{connective-equiv-k-theory/shearing}
\input{connective-equiv-k-theory/equivariant-k-theory}
\input{connective-equiv-k-theory/equivariant-little-ku}
\input{connective-equiv-k-theory/k-homology}

\newpage

\section{(Derived) geometric Satake and variants}
\input{equivalences/full-faithfulness}
\input{equivalences/geometric-satake-review}
\input{equivalences/spherical-varieties}
\input{equivalences/whittaker-induction}
\input{equivalences/G-to-H-equiv-homology}
\input{equivalences/sheaves-equivalences}
\input{equivalences/ku-satake}

\newpage

\input{each-type/big-table}
\section{Case-by-case analysis}\label{sec: case by case}
\input{each-type/case-by-case-summary}
\input{each-type/homology-of-loops-sphere}
\input{each-type/type-An}
\input{each-type/type-Bn}
\input{each-type/type-Cn}
\input{each-type/type-Dn}
\input{each-type/type-F4}
\input{each-type/type-G2}
\input{each-type/type-B3prime}

\newpage

\section{Structures on the spectral side}\label{sec: everything deforms}
\input{ku-hamiltonian/beta-hamiltonian}
\input{ku-hamiltonian/duals-are-lagrangian}
\input{ku-hamiltonian/S1-equiv-MU}
\input{ku-hamiltonian/langlands-over-sphere}

\newpage
\appendix
\counterwithin{lemma}{section}
\input{appendices/shv-loops}
\newpage
\input{appendices/comparison-cohen-jones-yan}
\newpage
\input{appendices/open-qns}

\newpage
\bibliographystyle{alphanum}
\bibliography{main}
\end{document}

%% file: header.tex
\usepackage[dvipsnames]{xcolor}
\usepackage{todonotes}

\usepackage{etoolbox}
\usepackage{rotating}

\makeatletter
\patchcmd{\@thm}{\let\thm@indent\indent}{\let\thm@indent\noindent}{}{}
\patchcmd{\@thm}{\thm@headfont{\scshape}}{\thm@headfont{\bfseries}}{}{}

\usepackage[T1]{fontenc}

\usepackage{hyperref}

\usepackage{mathtools}

\let\U\undefined
\let\G\undefined

\usepackage{multirow}
\usepackage{longtable}

\usepackage{amsmath}
\usepackage{amsfonts}
\usepackage{amssymb}
\usepackage{float}
\usepackage{amsthm}
\usepackage{comment}
\usepackage{amscd}
\usepackage{amsxtra}
\usepackage{epsfig}
\usepackage{epigraph}
\usepackage{pgfplots}
\usepackage{pdflscape} 
\usepackage{afterpage}

\usepackage{stmaryrd}

\usepackage{listings}
\usepackage{color}

\definecolor{dkgreen}{rgb}{0,0.6,0}
\definecolor{gray}{rgb}{0.5,0.5,0.5}
\definecolor{mauve}{rgb}{0.58,0,0.82}

\usepackage{color}
\usepackage[all]{xy}
\usepackage{verbatim}

\usepackage{eucal}
\usepackage{mathrsfs}

\usepackage{graphicx}
\usepackage{tikz-cd}

\usepackage{url}
\usepackage{verbatim}

\usepackage{xy}
\xyoption{all}
\newcommand{\xar}[1]{\xrightarrow{{#1}}}

\usepackage{amsthm}
\usepackage{hhline}
\usepackage{enumitem}

\usepackage[normalem]{ulem}

\usepackage{caption} 
\makeatletter
\def\@tocline#1#2#3#4#5#6#7{\relax
  \ifnum #1>\c@tocdepth 
  \else
    \par \addpenalty\@secpenalty\addvspace{#2}%
    \begingroup \hyphenpenalty\@M
    \@ifempty{#4}{%
      \@tempdima\csname r@tocindent\number#1\endcsname\relax
    }{%
      \@tempdima#4\relax
    }%
    \parindent\z@ \leftskip#3\relax \advance\leftskip\@tempdima\relax
    \rightskip\@pnumwidth plus4em \parfillskip-\@pnumwidth
    #5\leavevmode\hskip-\@tempdima
      \ifcase #1
       \or\or \hskip 1em \or \hskip 2em \else \hskip 3em \fi%
      #6\nobreak\relax
    \hfill\hbox to\@pnumwidth{\@tocpagenum{#7}}\par
    \nobreak
    \endgroup
  \fi}
\makeatother

\usepackage[noabbrev,capitalize]{cleveref} 

\crefformat{nul}{(#2#1#3)}
\Crefformat{nul}{(#2#1#3)}

\crefname{section}{\S}{\S\S}
\crefname{subsection}{\S}{\S\S}
\crefname{axioms}{Axiom}{Axioms}
\crefname{exercise}{Exercise}{Exercises}
\crefname{exercisenum}{Exercise}{Exercises}
\crefname{construction}{Construction}{Constructions}
\crefname{problem}{Problem}{Problems}
\crefname{theorem}{Theorem}{Theorems}
\crefname{definition}{Definition}{Definitions}
\crefname{prop}{Proposition}{Propositions}
\crefname{lemma}{Lemma}{Lemmas}
\crefname{example}{Example}{Examples}
\crefname{examplealph}{Example}{Examples}
\crefname{corollary}{Corollary}{Corollaries}
\crefname{nonexample}{Nonexample}{Nonexamples}
\crefname{equation}{}{}
\crefname{summary}{Summary}{Summaries}
\crefname{recollection}{Recollection}{Recollections}
\Crefname{recollection}{Recollection}{Recollections}
\Crefname{nonexample}{Nonexample}{Nonexamples}
\Crefname{corollary}{Corollary}{Corollaries}
\Crefname{corollary}{Corollary}{Corollaries}
\Crefname{axioms}{Axiom}{Axioms}
\Crefname{exercise}{Exercise}{Exercises}
\Crefname{exercisenum}{Exercise}{Exercises}
\Crefname{construction}{Construction}{Constructions}
\Crefname{problem}{Problem}{Problems}
\Crefname{theorem}{Theorem}{Theorems}
\Crefname{definition}{Definition}{Definitions}
\Crefname{prop}{Proposition}{Propositions}
\Crefname{lemma}{Lemma}{Lemmas}
\Crefname{example}{Example}{Examples}
\Crefname{examplealph}{Example}{Examples}
\Crefname{section}{\S}{\S\S}
\Crefname{subsection}{\S}{\S\S}

\DeclareMathOperator{\Fun}{Fun}
\DeclareMathOperator{\Hom}{Hom}

\DeclareMathOperator{\Ext}{Ext}

\DeclareMathOperator{\spec}{Spec}
\DeclareMathOperator{\spf}{Spf}

\DeclareMathOperator{\supp}{Supp}

\DeclareMathOperator{\picsp}{\mathfrak{pic}}

\DeclareMathOperator{\fib}{fib}

\DeclareMathOperator{\QCoh}{\mathrm{QCoh}}

\DeclareMathOperator{\Tr}{Tr}

\DeclareMathOperator{\End}{End}
\DeclareMathOperator{\Pic}{Pic}

\DeclareMathOperator{\colim}{colim}

\DeclareMathOperator{\Tot}{Tot}
\DeclareMathOperator{\im}{im}

\DeclareMathOperator{\Sym}{Sym}
\DeclareMathOperator{\LSym}{LSym}
\DeclareMathOperator{\proj}{Proj}

\let\Bar\relax
\DeclareMathOperator{\Bar}{\mathrm{Bar}}

\numberwithin{equation}{section}
\newtheorem{lemma}{Lemma}[subsection]

\newtheorem{corollary}[lemma]{Corollary}

\newtheorem{theorem}[lemma]{Theorem}

\newtheorem*{corno}{Corollary}
\newtheorem*{thmno}{Theorem}
\newtheorem*{propo}{Proposition}
\newtheorem{prop}[lemma]{Proposition}
\newtheorem{conjecture}[lemma]{Conjecture}
\newtheorem*{conjectureno}{Conjecture}

\newtheorem*{conj-moore}{Conjecture~\ref{moore-splitting}}

\theoremstyle{definition}

\newtheorem{question}[lemma]{Question}
\newtheorem{assume}[lemma]{Assumption}
\newtheorem{definition}[lemma]{Definition}
\newtheorem{warning}[lemma]{Warning}
\newtheorem{example}[lemma]{Example}

\newtheorem{construction}[lemma]{Construction}
\newtheorem{remark}[lemma]{Remark}
\newtheorem{notation}[lemma]{Notation}
\newtheorem{variant}[lemma]{Variant}
\newtheorem{setup}[lemma]{Setup}

\newtheorem{recall}[lemma]{Recollection}

\newtheorem{slogan}[lemma]{Slogan}
\newtheorem{observe}[lemma]{Observation}

\newtheorem{strategy}[lemma]{Strategy}
\newtheorem{expect}[lemma]{Expectation}
\newtheorem{hypothesis}[lemma]{Hypothesis}
\newtheorem{goal}[lemma]{Goal}

\newtheorem*{notno}{Notation}

\renewcommand{\AA}{\mathbf{A}}
\newcommand{\FF}{\mathbf{F}}
\newcommand{\Z}{\mathbf{Z}}
\newcommand{\QQ}{\mathbf{Q}}
\newcommand{\cc}{\mathbf{C}}
\newcommand{\cJ}{\mathcal{J}}
\newcommand{\cC}{\mathcal{C}}
\newcommand{\cd}{\mathcal{D}}
\newcommand{\RR}{\mathbf{R}}
\newcommand{\Deltab}{\mathbf{\Delta}}

\newcommand{\Sp}{\mathrm{Sp}}

\newcommand{\Aff}{\mathrm{Aff}}

\newcommand{\Msfg}{\M^\mathrm{s}_\mathrm{FG}}

\newcommand{\co}{\mathcal{O}}

\newcommand{\PP}{\mathbf{P}}

\newcommand{\GG}{\mathbf{G}}

\newcommand{\LMod}{\mathrm{LMod}}

\newcommand{\Mod}{\mathrm{Mod}}
\newcommand{\coMod}{\mathrm{coMod}}

\newcommand{\cL}{\mathcal{L}}

\newcommand{\CP}{{\mathbf{C}P}}
\newcommand{\HHP}{\mathbf{H}P}
\newcommand{\RP}{\mathbf{R}P}
\newcommand{\OP}{\mathbf{O}P}

\newcommand{\Map}{\mathrm{Map}}

\newcommand{\dX}{\mathcal{X}}

\newcommand{\dZ}{\mathfrak{Z}}

\newcommand{\cg}{\mathcal{G}}

\newcommand{\cA}{\mathcal{A}}

\newcommand{\Ab}{\mathrm{Ab}}

\newcommand{\Alg}{\mathrm{Alg}}
\newcommand{\Sq}{\mathrm{Sq}}

\newcommand{\N}{\mathrm{N}}

\newcommand{\Eoo}{{\mathbf{E}_\infty}}

\newcommand{\CAlg}{\mathrm{CAlg}}

\newcommand{\cI}{\mathcal{I}}
\newcommand{\cf}{\mathcal{F}}

\newcommand{\univ}{\mathrm{univ}}

\newcommand{\Def}{\mathrm{Def}}

\newcommand{\Lie}{\mathrm{Lie}}

\newcommand{\fin}{\mathrm{fin}}
\newcommand{\ev}{\mathrm{ev}}

\newcommand{\cH}{\mathcal{H}}

\newcommand{\cN}{\mathcal{N}}

\newcommand{\ul}[1]{\underline{#1}}
\newcommand{\ol}[1]{\overline{#1}}

\newcommand{\Perf}{\mathrm{Perf}}

\newcommand{\E}[1]{{\mathbf{E}_{{#1}}}}
\newcommand{\mmod}{/\!\!/}

\newcommand{\id}{\mathrm{id}}

\newcommand{\gl}{\mathfrak{gl}}
\newcommand{\SL}{\mathrm{SL}}
\newcommand{\GL}{\mathrm{GL}}
\newcommand{\BGL}{\mathrm{BGL}}

\newcommand{\MO}{\mathrm{MO}}
\newcommand{\MU}{\mathrm{MU}}
\newcommand{\BU}{\mathrm{BU}}
\renewcommand{\O}{\mathrm{O}}
\newcommand{\SO}{\mathrm{SO}}
\newcommand{\BO}{\mathrm{BO}}
\newcommand{\KU}{\mathrm{KU}}

\newcommand{\ku}{\mathrm{ku}}
\newcommand{\ko}{\mathrm{ko}}

\newcommand{\U}{\mathrm{U}}
\newcommand{\SU}{\mathrm{SU}}

\renewcommand{\H}{\mathrm{H}}

\newcommand{\BSO}{\mathrm{BSO}}
\newcommand{\BSU}{\mathrm{BSU}}

\newcommand{\MSO}{\mathrm{MSO}}

\newcommand{\Spin}{\mathrm{Spin}}

\newcommand{\Conf}{\mathrm{Conf}}

\newcommand{\Ind}{\mathrm{Ind}}

\newcommand{\F}{\mathrm{F}}

\newcommand{\HH}{\mathrm{HH}}
\newcommand{\HC}{\mathrm{HC}}

\newcommand{\Coh}{\mathrm{Coh}}

\newcommand{\Cat}{\mathrm{Cat}}

\newcommand{\Sch}{\mathrm{Sch}}
\newcommand{\op}{\mathrm{op}}

\newcommand{\cB}{\mathcal{B}}
\newcommand{\dR}{\mathrm{dR}}

\newcommand{\fr}[1]{\mathfrak{#1}}
\newcommand{\g}{\mathfrak{g}}
\newcommand{\G}{\mathrm{G}}

\newcommand{\cP}{\mathcal{P}}

\newcommand{\pr}{\mathrm{pr}}

\newcommand{\cU}{\mathcal{U}}

\newcommand{\Rep}{\mathrm{Rep}}

\newcommand{\gr}{\mathrm{gr}}
\newcommand{\ds}{\mathrm{ds}}
\newcommand{\fil}{\mathrm{fil}}

\newcommand{\Free}{\mathrm{Free}}

\newcommand{\can}{\mathrm{can}}

\newcommand{\Vect}{\mathrm{Vect}}

\newcommand{\cR}{\mathcal{R}}

\newcommand{\coCAlg}{\mathrm{coCAlg}}

\newcommand{\cV}{\mathcal{V}}

\newcommand{\cl}{\mathrm{cl}}

\newcommand{\unit}{\mathbf{1}}

\newcommand{\und}{\mathrm{und}}
\newcommand{\cM}{\mathcal{M}}

\newcommand{\et}{\mathrm{et}}
\newcommand{\Shv}{\mathrm{Shv}}

\newcommand{\SH}{\mathrm{SH}}

\newcommand{\Bun}{\mathrm{Bun}}
\renewcommand{\ss}{\mathrm{ss}}

\newcommand{\IndCoh}{\mathrm{IndCoh}}

\newcommand{\Loc}{\mathrm{Loc}}

\newcommand{\Ad}{\mathrm{Ad}}

\renewcommand{\log}{\mathrm{log}}

\newcommand{\M}{\mathcal{M}}

\newcommand{\PGL}{\mathrm{PGL}}

\renewcommand{\det}{\mathrm{det}}

\newcommand{\threerightarrows}{{\substack{\rightarrow\\[-1em] \rightarrow \\[-1em] \rightarrow}}}

\newcommand{\Rees}{\mathrm{Rees}}

\newcommand{\IC}{\mathrm{IC}}

\newcommand{\per}{\mathrm{per}}

\newcommand{\bX}{\mathbb{X}}
\newcommand{\reg}{\mathrm{reg}}

\newcommand{\cS}{\mathcal{S}}

\renewcommand{\S}{Section }

\newcommand{\Prst}{\mathrm{Pr}^\mathrm{L,st}}

\newcommand{\WW}{\mathbf{W}}

\newcommand{\Top}{{\mathcal{S}}}

\renewcommand{\sqcup}{\amalg}

\makeatletter
\providecommand{\leftsquigarrow}{%
  \mathrel{\mathpalette\reflect@squig\relax}%
}
\newcommand{\reflect@squig}[2]{%
  \reflectbox{$\m@th#1\rightsquigarrow$}%
}
\makeatother

\newcommand{\bull}{\bullet}

\newcommand{\pw}[1]{[\![#1]\!]}
\newcommand{\ls}[1]{(\!(#1)\!)}

\newcommand{\act}{\circlearrowleft}
\newcommand{\LinCat}{\mathrm{LinCat}}

\newcommand{\Cn}{{\Z/n}}

\newcommand{\rank}{\mathrm{rank}}

\newcommand{\sh}{\mathrm{sh}}

\renewcommand{\sl}{\mathfrak{sl}}

\renewcommand{\tilde}{\widetilde}

\renewcommand{\min}{\mathrm{min}}

\newcommand{\Gr}{\mathrm{Gr}}

\newcommand{\diag}{\mathrm{diag}}

\newcommand{\oblv}{\mathrm{oblv}}

\newcommand{\Syn}{\mathrm{Syn}}

\newcommand{\Frac}{\mathrm{Frac}}

\newcommand{\rot}{\mathrm{rot}}

\newcommand{\ld}[1]{{\check{{#1}}}}

\usepackage{slashed}

\usepackage{relsize}
\usepackage[bbgreekl]{mathbbol}
\usepackage{amsfonts}
\DeclareSymbolFontAlphabet{\mathbb}{AMSb} 
\DeclareSymbolFontAlphabet{\mathbbl}{bbold}

\newcommand{\Ghat}{\hat{\mathbf{G}}}

\newcommand{\pdb}[1]{\langle{#1}\rangle}

\newcommand{\pgl}{\mathfrak{pgl}}

\newcommand{\std}{\mathrm{std}}

\newcommand{\Bl}{\mathcal{B}}

\newcommand{\Perv}{\mathcal{P}\mathrm{erv}}

\newcommand{\arth}{\mathrm{Arth}}

\newcommand{\Orb}{\mathrm{Orb}}

\newcommand{\Sat}{\mathrm{Sat}}

\newcommand{\punc}[1]{{#1}^\circ}

\newcommand{\WInd}{\mathrm{WhitInd}}
\newcommand{\olQell}{{\ol{\QQ_\ell}}}

\newcommand{\glob}{\mathrm{glob}}
\newcommand{\Lat}{\mathrm{Lat}}

\newcommand{\mutil}{\tilde{\mu}}

\newcommand{\bV}{\mathbf{V}}

\newcommand{\Mir}{\mathrm{Mir}}

\newcommand{\faux}{\mathrm{faux}}

%% file: intro/intro.tex
The goal of this article is multifold: we will study the local geometric Langlands conjecture of Ben-Zvi--Sakellaridis--Venkatesh (\cite[Conjecture 7.5.1]{bzsv}) in the case of affine homogeneous spherical varieties of relative rank $1$ (under \cref{hypothesis: rank 1 weakly placid}); we will also show through examples that this conjecture admits a generalization with coefficients in connective complex K-theory\footnote{I would like to emphasize at the outset that the $\infty$-category defined as ``constructible sheaves with $\ku$-coefficients'' in this article is rather contrived; it is constructed essentially by naively ``closing up'' the subcategory of local systems (on this subcategory, the definition is OK, at least upon rationalization). I hope to return to the foundational question of providing a better definition in later work.}; and finally, we will suggest some broader ideas in stable homotopy theory which attempts to relate the aforementioned Langlands duality conjecture with the main results of \cite{grg-reg}. Each of these steps will be explained in detail later in the introduction, but let me say here that (I believe) the main contribution of this article is the proposal of an approach to the aforementioned Ben-Zvi--Sakellaridis--Venkatesh conjecture for affine homogeneous varieties via an analogue of the regular centralizer group scheme.

\begin{remark}
    Despite the title of this article, knowledge of complex K-theory/stable homotopy theory is \textit{not} necessary to understand the proof of our main result \cref{intro: thm bzsv rk 1}. The reader uninterested in complex K-theory/stable homotopy theory should obviously ignore the relevant sections, and should note that any result stated in the body of this article for connective complex K-theory $\ku$ (often localized at some even integer) will imply an analogous result over $\Z$ (localized at the aforementioned integer) simply by setting the symbol $\beta$ to zero.
\end{remark}

Since the goals of this article involve a few different areas of mathematics, we will give a somewhat long-winding introduction to motivate the results presented here; discussion of work which is new to this article starts around \cref{intro: thm bzsv rk 1}. Apologies in advance to experts for the prolonged introduction!

\subsection{Spherical harmonics and its geometrization}

Broadly speaking, the Langlands program aims to study spectral decompositions of automorphic representations. The historically first example of such a spectral decomposition is the Fourier transform: 
\begin{example}
    The Fourier transform says that any $f\in L^2(S^1; \cc)$ can be expressed in terms of the spherical functions $\exp(inx)$ (which are eigenvectors for the derivative operator). One can view Fourier analysis as providing a decomposition $L^2(S^1; \cc) \cong \ell^2(\Z)$ into irreducible representations of $S^1$ acting on $L^2(S^1; \cc)$. Note that $\ell^2(\Z)$ is a completion of an infinite direct sum of the irreducible representations $\cc \cdot \exp(inx)$ of $S^1$; each appears with multiplicity $1$.
\end{example}
The Fourier transform was soon generalized to the theory of spherical harmonics, which studies the decomposition of $L^2(S^{n-1}; \cc)$ under the action of the group of rotations $\O_n \subseteq \GL_n(\RR)$. 
\begin{example}
    Let $\cH_j(\RR^n)$ denote the space of homogeneous harmonic polynomials $\RR^n \to \cc$ of degree $j$, and let $\cH_j(S^{n-1})$ denote the space of functions on $S^{n-1}$ obtained by restricting elements of $\cH_j(\RR^n)$ to $S^{n-1} \subseteq \RR^n$. Then, $L^2(S^{n-1}; \cc)$ is isomorphic to a completion of the direct sum $\bigoplus_{j\geq 0} \cH_j$, and each $\cH_j$ is an irreducible $\O_n$-representation appearing with multiplicity $1$.
\end{example}

Observing that $S^{n-1} \simeq \O_n/\O_{n-1}$, the theory of spherical harmonics can be generalized even further: if $G$ is a reductive algebraic group over $\cc$, $K$ is a maximal compact subgroup of $G(\cc)$, and $H\subseteq G$ is a closed subgroup, one can attempt to understand the decomposition of $L^2(G(\cc)/H(\cc); \cc)$ into irreducible $K$-representations. One can also state this goal in the $p$-adic setting, where the maximal compact subgroup $K\subseteq G(\cc)$ is replaced by the subgroup $G(\Z_p)\subseteq G(\QQ_p)$. Namely, if $H\subseteq G$ is a closed subgroup of a reductive algebraic group over $\Z_p$, one can study the decomposition of the space $C_c(G(\QQ_p)/H(\QQ_p); \cc)$ of compactly supported functions into irreducible $G(\Z_p)$-representations. This turns out to be especially understandable in the case of multiplicity $1$, in which case $H \subseteq G$ is called ``spherical''. Questions of this form have been placed into the context of the Langlands program by Sakellaridis and Venkatesh (among others) in \cite{sakellaridis-venkatesh}.

The archimedean and $p$-adic settings being too elaborate, it is often simpler to \textit{geometrize} such questions by studing their function field variants. 
\begin{recall}
    There is a standard analogy between $p$-adic number fields such as $\QQ_p$ (along with its ring of integers $\Z_p$) and function fields such as $\ol{\FF}_p\ls{t}$ (along with its ring of integers $\ol{\FF}_p\pw{t}$). There is a further analogy between $\ol{\FF}_p\ls{t}$ (along with its ring of integers $\ol{\FF}_p\pw{t}$) and $\cc\ls{t}$ (along with its ring of integers $\cc\pw{t}$). In this case, $G(\cc\ls{t})$ can be regarded as the $\cc$-points of the formal loop group of $G$, and hence acquires a natural topology. Therefore, instead of studying the decomposition of the space $C_c(G(\cc\ls{t})/H(\cc\ls{t}); \cc)$ of compactly supported functions into irreducible $G(\cc\pw{t})$-representations, we can further study a decomposition of the compactly supported cohomology ring $\H^\ast_c(G(\cc\ls{t})/H(\cc\ls{t}); \cc)$ into irreducible $G(\cc\pw{t})$-representations. (In the function field case, one can similarly consider the \'etale cohomology of $G(\ol{\FF}_p\ls{t})/H(\ol{\FF}_p\ls{t})$.)
\end{recall}

The ``spherical''\footnote{Apologies for the unfortunate, but standard, terminology! The terminology clash is even worse than one might expect: we will study both the sphere spectrum, as well as spectral decomposition for spherical varieties (and even propose the existence of a theory of spectral decomposition over the sphere spectrum for spherical varieties!).
In trying to keep with the title of \cite{grg-reg}, I recently learned that the phenomenon of ``chromatic aberration'' also goes by the name ``spherochromatism'' ... which is (unexpectedly) just perfect.
}
part of this cohomology ring is given by the \textit{$G(\cc\pw{t})$-equivariant} cohomology of $G(\cc\ls{t})/H(\cc\ls{t})$.
This cohomology ring can be itself be categorified: namely, one can consider the $\infty$-category of constructible $G(\cc\pw{t})$-equivariant sheaves of $\cc$-vector spaces on $G(\cc\ls{t})/H(\cc\ls{t})$, denoted $\Shv^c_{G(\cc\pw{t})}(G(\cc\ls{t})/H(\cc\ls{t}); \cc)$. This $\infty$-category will be the main topic of study in this article: in particular, when $G/H$ is an affine spherical variety of rank $1$ (the meaning of which will be explained below), we will study a spectral decomposition of this $\infty$-category.

In order to explain the precise sense in which this $\infty$-category admits a spectral decomposition, let us return to the de-categorified function field setting; in other words, consider the vector space of $G(\ol{\FF}_p\pw{t})$-invariant compactly supported functions $G(\ol{\FF}_p\ls{t})/H(\ol{\FF}_p\ls{t}) \to \cc$, denoted $C_{c, G(\ol{\FF}_p\pw{t})}(G(\ol{\FF}_p\ls{t})/H(\ol{\FF}_p\ls{t}); \cc)$. This vector space canonically admits an action of the algebra $C_{c, G(\ol{\FF}_p\pw{t}) \times G(\ol{\FF}_p\pw{t})}(G(\ol{\FF}_p\ls{t}); \cc)$ equipped with the product given by convolution. From now, let us assume (for simplicity) that $G$ is semisimple. One then has the following famous theorem:
\begin{theorem}[Satake isomorphism]
     There is an explicit isomorphism (defined by Macdonald)
     $$C_{c, G(\ol{\FF}_p\pw{t}) \times G(\ol{\FF}_p\pw{t})}(G(\ol{\FF}_p\ls{t}); \cc) \cong \cc[\bX_\ast(T)]^W,$$
     where $W$ is the Weyl group and $\bX_\ast(T)$ is the lattice of \textit{cocharacters} of $T$. 
\end{theorem}
The right-hand side is not quite the complexification of the representation ring of $G$, which would instead be isomorphic to $\cc[\bX^\ast(T)]^W$ by highest weight theory; instead, it is the complexification of the representation ring of the \textit{Langlands dual} group $\ld{G}$, which is defined so that the weights, coweights, roots, and coroots of $G$ are the coweights, weights, coroots, and roots of $\ld{G}$ (respectively). In other words, the Satake isomorphism gives an explicit isomorphism
$$C_{c, G(\ol{\FF}_p\pw{t}) \times G(\ol{\FF}_p\pw{t})}(G(\ol{\FF}_p\ls{t}); \cc) \cong K_0(\Rep_\cc(\ld{G})) \otimes \cc.$$
One therefore obtains an action
$$K_0(\Rep_\cc(\ld{G})) \otimes \cc \act C_{c, G(\ol{\FF}_p\pw{t})}(G(\ol{\FF}_p\ls{t})/H(\ol{\FF}_p\ls{t}); \cc),$$
and the task of providing a ``spectral decomposition'' of this vector space can be more precisely phrased as giving an explicit description of this action in terms of the Langlands dual group $\ld{G}$.

This interpretation of our task can be categorified, since both sides of the Satake isomorphism admit natural categorifications. A categorification of the Satake isomorphism itself is provided by the famous geometric Satake equivalence of Mirkovic-Vilonen \cite{mirkovic-vilonen}. To state it, let us switch back to the Laurent series ring $\cc\ls{t}$. Let $\Shv^c_{G(\cc\pw{t}) \times G(\cc\pw{t})}(G(\cc\ls{t}); \cc)$ denote the $\infty$-category of $G(\cc\pw{t}) \times G(\cc\pw{t})$-equivariant constructible sheaves of $\cc$-vector spaces on $G(\cc\ls{t})$ equipped with its perverse $t$-structure. Convolution defines a symmetric monoidal structure on the heart of the perverse $t$-structure, and then one has:
\begin{theorem}[Mirkovic-Vilonen]
    There is a symmetric monoidal equivalence
    $$\Shv^c_{G(\cc\pw{t}) \times G(\cc\pw{t})}(G(\cc\ls{t}); \cc)^\heartsuit \simeq \Rep(\ld{G})^\heartsuit$$
    of abelian categories. 
\end{theorem}
The na\"ive guess that this equivalence promotes to an equivalence of derived $\infty$-categories is \textit{false}; that is,
$$\Shv^c_{G(\cc\pw{t}) \times G(\cc\pw{t})}(G(\cc\ls{t}); \cc) \not \simeq \Rep(\ld{G})^\heartsuit.$$
We will return to this momentarily. For the moment, note that there is nevertheless a canonical action 
$$\Shv^c_{G(\cc\pw{t}) \times G(\cc\pw{t})}(G(\cc\ls{t}); \cc)^\heartsuit \act \Shv^c_{G(\cc\pw{t})}(G(\cc\ls{t})/H(\cc\ls{t}); \cc)^\heartsuit,$$
and this discussion suggests that the appropriate categorification of the task of providing a spectral decomposition of $C_{c, G(\cc\pw{t})}(G(\cc\ls{t})/H(\cc\ls{t}); \cc)$ would be: 
\begin{goal}\label{goal: prelim intro goal}
    Describe the action
    $$\Rep(\ld{G})^\heartsuit \act \Shv^c_{G(\cc\pw{t})}(G(\cc\ls{t})/H(\cc\ls{t}); \cc)^\heartsuit,$$
    perhaps by proving an equivalence between this category and the abelian category of quasicoherent sheaves on a quotient stack $\ld{M}/\ld{G}$ for some $\ld{G}$-space $\ld{M}$.
\end{goal}
Unfortunately, such an equivalence is generally not possible, since taking the heart of a $t$-structure is a rather severe process.
It is therefore natural to ask for a generalization of the Mirkovic-Vilonen equivalence describing the full $\infty$-category $\Shv^c_{G(\cc\pw{t}) \times G(\cc\pw{t})}(G(\cc\ls{t}); \cc)$. Such a description was provided by Bezrukavnikov-Finkelberg in \cite{bf-derived-satake} (following earlier work of Lusztig and Ginzburg; see \cite{lusztig-satake, ginzburg-langlands}), and the answer involves derived algebraic geometry. 
We will state the \textit{arithmetically sheared} (in the sense of \cite[Section 6.7]{bzsv}) version of the derived geometric Satake equivalence.
\begin{theorem}[Bezrukavnikov-Finkelberg; \cref{thm: derived satake}]\label{eq: intro bf derived}
    Let $\ld{\g}$ denote the Lie algebra of the Langlands dual group, viewed as a $\cc$-vector space, let $2\rho: \GG_m \to \ld{T}$ denote the sum of the positive coroots of $\ld{G}$, and let $\ld{\g}[2-2\rho]$ denote its $(2-2\rho)$-fold shift. Then, there is a monoidal equivalence\footnote{I will abusively write $\ld{\g}$ below, although the correct object on the coherent side is $\ld{\g}^\ast$. If $G$ is semisimple, these are isomorphic.}
    \begin{equation}
        \Shv^c_{G(\cc\pw{t}) \times G(\cc\pw{t})}(G(\cc\ls{t}); \cc) \simeq \Perf(\ld{\g}[2-2\rho]/\ld{G}[-2\rho]).
    \end{equation}
\end{theorem}
This is often known as the \textit{derived} geometric Satake equivalence, and it has some precedence in homotopy theory:
\begin{remark}
    Let $\fr{t}$ denote a Cartan subalgebra of $\g$, and let $W$ denote the Weyl group of $G$. Then, there is an isomorphism\footnote{Technically, this is isomorphism is not true; however, only for the purposes of this introduction, we will purposely conflate grading shifts with homological shifts.}
    $\spec \H^\ast_G(\ast; \cc) \cong \fr{t}[2]\mmod W$. There is also an isomorphism $\ld{\fr{t}}^\ast \cong \fr{t}$, so that the Chevalley restriction theorem gives an isomorphism $\fr{t}[2]\mmod W \cong \ld{\g}^\ast[2]\mmod \ld{G}$, and hence an isomorphism $\spec \H^\ast_G(\ast; \cc) \cong \ld{\g}^\ast[2]\mmod \ld{G}$. This implies an equivalence
    $$\Shv^c_{G(\cc\pw{t}) \times G(\cc\pw{t})}(G(\cc\pw{t}); \cc) \simeq \Shv^c_{G(\cc\pw{t})}(\ast; \cc) \simeq \Perf(\ld{\g}[2]\mmod \ld{G}),$$
    which is a restriction of the derived geometric Satake equivalence to equivariant sheaves on $G(\cc\pw{t}) \subseteq G(\cc\ls{t})$.
\end{remark}

There is an $\E{3}$-monoidal structure on the left-hand side of the derived geometric Satake equivalence, which comes from the geometry of the affine Grassmannian (and is spelled out in \cite{nocera-e3}).
Unfortunately, the fact that the derived geometric Satake equivalence is $\E{3}$-monoidal (in an appropriate sense) does not seem to be recorded anywhere in the literature (although it is closely related to work of Campbell-Raskin; see \cite{campbell-raskin-satake}), and we will also not address this point in our discussion.
We will give an argument for the above equivalence in \cref{thm: derived satake} which is slightly different from that of \cite{bf-derived-satake}; the key step (already accomplished in \cite{bfm}, and \cite{homology-langlands} in arbitrary characteristic, except for some carefully excluded primes) is the construction of a homomorphism
\begin{equation}\label{eq: map from Loops G to dual group}
    \spec \H_\ast^G(\Omega G; \cc) \to \ld{G} \times \ld{\g}^\ast[2]\mmod \ld{G}
\end{equation}
of group schemes over $\ld{\g}^\ast[2]\mmod \ld{G} \cong \spec \H^\ast_G(\ast; \cc)$.

We can now formulate the ``correct'' version of \cref{goal: prelim intro goal}.
\begin{goal}
    There is a canonical action 
    $$\Shv^c_{G(\cc\pw{t}) \times G(\cc\pw{t})}(G(\cc\ls{t}); \cc) \act \Shv^c_{G(\cc\pw{t})}(G(\cc\ls{t})/H(\cc\ls{t}); \cc),$$
    and we can state the task of providing a spectral decomposition of the latter $\infty$-category as explicitly describing the action
    $$\Perf(\ld{\g}[2-2\rho]/\ld{G}[-2\rho]) \act \Shv^c_{G(\cc\pw{t})}(G(\cc\ls{t})/H(\cc\ls{t}); \cc).$$
\end{goal}
The work of Ben-Zvi--Sakellaridis--Venkatesh \cite{bzsv} provides numerous conjectures about this description: namely, they conjecture in \cite[Conjecture 7.5.1]{bzsv} that if $G/H$ is a spherical $G$-variety (satisfying some other conditions), there is a graded \textit{Hamiltonian} $\ld{G}$-space $\ld{M}$ such that there is an equivalence\footnote{This is technically a slight lie: the left-hand side is replaced by a certain subcategory defined using the action of $\Shv^c_{G(\cc\pw{t}) \times G(\cc\pw{t})}(G(\cc\ls{t}); \cc)$.}
\begin{equation}\label{eq: intro bzsv display}
    \Shv^c_{G(\cc\pw{t})}(G(\cc\ls{t})/H(\cc\ls{t}); \cc) \simeq \Perf(\sh^{1/2} \ld{M}/\ld{G}(-2\rho)),
\end{equation}
and the action of $\Perf(\ld{\g}[2-2\rho]/\ld{G}[-2\rho])$ on the left-hand side is specified by the moment map $\mu: \ld{M}/\ld{G} \to \ld{\g}/\ld{G}$. Here, $\sh^{1/2}$ denotes a shearing, which converts gradings into homological shifts (more precisely, it sends a module in weight $2n$ to the same module shifted homologically by $2n$). Moreover, they give a precise construction of the predicted dual variety $\ld{M}$.
One of our main goals in this article is to show that \cite[Conjecture 7.5.1]{bzsv} is true for the simplest building blocks of spherical varieties:
\begin{theorem}[\cref{thm: rk 1 bzsv is true}]\label{intro: thm bzsv rk 1}
    Suppose $G/H$ is an affine spherical variety of rank $1$, and assume \cref{hypothesis: rank 1 weakly placid} holds\footnote{If the sheaf theory for the $G(\cc\pw{t})$-action on $G(\cc\ls{t})/H(\cc\ls{t})$ is sufficiently well-behaved, it should be possible to forego this hypothesis; but regardless of the sheaf-theoretic setup and the ultimate correctness of \cref{hypothesis: rank 1 weakly placid}, we believe that the calculations of \cref{sec: case by case} will be the key to proving any sort of Langlands duality.}. Then \cite[Conjecture 7.5.1]{bzsv} is true, i.e., there is an equivalence \cref{eq: intro bzsv display} for the dual variety $\ld{M}$ constructed in \cite[Section 4]{bzsv}.\footnote{In fact, we only prove a bare equivalence; namely, we do not check compatibility with the action of the spherical Hecke category $\Shv^c_{G(\cc\pw{t}) \times G(\cc\pw{t})}(G(\cc\ls{t}); \cc)$. We do not expect this to be an especially difficult task, but it is one we decided to omit.}

    More precisely, for each line of \cref{intro table: topology and dualities for rank 1 spherical varieties}, there is an isomorphism
    $$\ld{M}/\ld{G} \cong \ld{Y}/\ld{G}_X \times \text{``Normalization''},$$
    and an equivalence 
    $$\Shv^c_{G(\cc\pw{t})}(G(\cc\ls{t})/H(\cc\ls{t}); \cc) \simeq \Perf(\sh^{1/2} \ld{M}/\ld{G}),$$
    where the normalization term accounts for the parabolic subgroup stabilizing the open $B$-orbit in $G/H$, and $\ld{G}_X$ (isomorphic to $\SL_2$ in our case) denotes the dual group of \cite{sakellaridis-venkatesh} (see also \cite{gaitsgory-nadler, knop-schalke}).
\end{theorem}
\begin{table}[h]
\centering
{
\begin{tabular}{ |c|c|c|c|c|c|c| } 
 \hline
 Name & $X = G/H$ & Dual $\ld{Y}$ & Topological explanation \\
 \hline
 $A_n$ & $\PGL_{n+1}/\GL_n$ & $T^\ast(2n) \AA^2(2n,0)$ & Hopf fibration \\
 $B_n$ & $\SO_{2n+1}/\SO_{2n}$ & $T^\ast(2n) \AA^2(4n-2,0)$ & EHP sequence \\
 $C_n$ & $\Sp_{2n}/(\Sp_2 \times \Sp_{2n-2})$ & $T^\ast(4n-4) \AA^2(4n-2,0)$ & Hopf fibration \\
 $D_n$ & $\SO_{2n}/\mu_2 \cdot \SO_{2n-1}$ & $\sl_2((2n-2)(1 - \rho_{\SL_2}))$ & James splitting \\
 $\F_4$ & $\F_4/\Spin_9$ & $T^\ast(16) \AA^2(22, 0)$ & Exceptional Hopf fibration \\
 $\G_2$ & $\G_2/\SL_3$ & $T^\ast(6) \AA^2(10,0)$ & EHP sequence \\
 $B_3'$ & $\SO_7/\G_2$ & $\sl_2(6-6\rho_{\SL_2})$ & James splitting \\
 \hline
\end{tabular}
}
\vspace{1cm}
\caption{Table of dual varieties and topological phenomena corresponding to each of the rank one affine homogeneous spherical varieties with no ``roots of type N'' (such varieties are excluded by \cite{sakellaridis-venkatesh, bzsv}). For each of these varieties, the dual group is $\ld{G}_X = \SL_2$ (which is also equipped with a certain grading that we have omitted in this table). With varied columns, this table will appear again in the present article; see, in particular, \cref{table: topology and dualities for rank 1 spherical varieties}. This latter table also contains the ``normalization'' term. Here, the notation $\AA^2(i,j)$ denotes an affine $2$-space with coordinates in weights $-i$ and $-j$; and $T^\ast(j)(X)$ denotes the cotangent bundle with cotangent fibers placed in weight $j$.
\newline\newline
The reader should compare the numbers in this table to the points of evaluation of the L-functions appearing in the rightmost column of \cite[Table 1]{sakellaridis-rank-1}. Namely, the dual stack $T^\ast(2j)\AA^2(2i,0)$ in our table corresponds to $L(\std, \tfrac{i}{2}) L(\std, \tfrac{2j-i}{2})$ in \cite[Table 1]{sakellaridis-rank-1}, and similarly $\sl_2(2j-2j\rho_{\SL_2})$ in our table corresponds to $L(\mathrm{ad}, j)$ in \cite[Table 1]{sakellaridis-rank-1}. In fact, the numbers $i$ and $j$ can be read off entirely from the rational homotopy type of $X$.
}
\label{intro table: topology and dualities for rank 1 spherical varieties}
\end{table}
\begin{remark}
    In the real and $p$-adic settings, the analogue of \cref{intro: thm bzsv rk 1} was studied in \cite[Theorem 1]{gan-gomez} and \cite{sakellaridis-rank-1}. In proving \cref{intro: thm bzsv rk 1}, the work of Sakellaridis (especially \cite{sakellaridis-rank-1}) was very influential.
\end{remark}
\begin{remark}[Why rank $1$?]
    Most of this article does not restrict attention to affine spherical varieties of rank $1$; this assumption is only imposed in \cref{sec: case by case} for doing computations. The restriction to rank $1$ here is not for any particularly deep reason: these varieties have very simple equivariant cell structures, which makes them more amenable to calculations. These examples also capture many interesting phenomena expected in \cite{bzsv}, and for these examples, the resulting homotopy-theoretic explanations for these phenomena become easier to understand.
\end{remark}
\begin{remark}
    The reader might notice a conspicuous absence of loop-rotation equivariance in this article, which, under Langlands duality, conspires to a deformation quantization of the spectral/coherent side of \cref{intro: thm bzsv rk 1}. We have chosen to separate this topic into a different article \cite{Eodd-and-quantizations}, in order to give it the detailed treatment it deserves (as well as for the purposes of length).
\end{remark}
In the remainder of this introduction, we will: 
\begin{enumerate}
    \item explain the meaning of the terms in, and the proof of, the above theorem, and illustrate it in the example of the spherical $\PGL_2$-variety $\PGL_2/\GG_m$. When applied to the spherical $G \times G$-variety $G$, this discussion recovers the derived geometric Satake equivalence (even when $G$ is not of rank $1$).
    \item explain the generalization of the derived geometric Satake equivalence to coefficients in connective complex K-theory, and some limited analogues of the above theorem on relative rank $1$ spherical varieties.
    \item discuss some conjectures and expectations about a further generalization to coefficients in the sphere spectrum.
\end{enumerate}
The homotopically-minded reader is suggested to skip to (b) and return to (a) as needed, and the conjecturally-minded reader is suggested to skip to (c).

\subsection{The proof of \cref{intro: thm bzsv rk 1}}

The basic strategy to prove \cref{intro: thm bzsv rk 1} is discussed in \cref{thm: ordinary homology criterion satake}. Let us give a high-level summary of this argument. Assume throughout that $H\subseteq G$ is a connected closed subgroup of a connected reductive algebraic group (we will sometimes identify this with the inclusion of their maximal compact subgroups).
\begin{strategy}\label{strategy: intro}
    \begin{enumerate}
        \item First, \cref{intro: thm bzsv rk 1} reduces to a question in homotopy theory under \cref{hypothesis: rank 1 weakly placid}; this reduction relies on \cref{thm: ordinary homology criterion satake}, which relies heavily on \cref{hypothesis: rank 1 weakly placid}. This step does not require assuming that $G_\cc/H_\cc$ is rank $1$.
        \item On the Langlands dual side, recall (as mentioned before \cref{eq: intro bzsv display}) that Ben-Zvi--Sakellaridis--Venkatesh construct the dual Hamiltonian $\ld{G}$-space $\ld{M}$ using the spherical geometry of the quotient $G/H$. One important observation we make is that the conjectures of \cite{bzsv} in particular suggest that there is an isomorphism $\ld{M}\mmod \ld{G} \cong \ld{\fr{h}}^\ast[2]\mmod \ld{H}$ of invariant-theoretic quotients, compatible with the natural maps to $\ld{\g}^\ast[2]\mmod \ld{G}$, as well as a closed immersion $\kappa: \ld{\fr{h}}^\ast[2]\mmod \ld{H} \to \ld{M}$ called the \textit{Kostant section} whose image ``generates'' $\co_{\ld{M}}$ under the $\ld{G}$-action; see \cref{conj: generalized kostant slice}. We will implicitly assume the existence of the map $\kappa$ in the remainder of this introduction.
        
        In the case of affine homogeneous rank $1$ spherical varieties, we construct the Kostant section case-by-case. This idea, however, works in much more generality. For instance, in \cite{triple-product-cayley}, we study the dual of the spherical $\SO_3 \times \mathrm{PSO}_4$-variety $(\SO_3 \times \mathrm{PSO}_4)/\SO_3^\mathrm{diag}$, in which case the above expectations boil down essentially to the existence of the Cayley hyperdeterminant.
        \item Using a compact generation argument, \cref{intro: thm bzsv rk 1} is reduced to proving that there is an isomorphism
        $$\spec \H_\ast^H(\Omega (G/H); \cc) \cong \ld{\fr{h}}^\ast[2]\mmod \ld{H} \times_{\ld{M}/\ld{G}[-2\rho]} \ld{\fr{h}}^\ast[2]\mmod \ld{H}$$
        of graded group schemes over $\ld{\fr{h}}^\ast[2]\mmod \ld{H}$. In other words, $\spec \H_\ast^H(\Omega (G/H); \cc)$ is the stabilizer of the image of the Kostant section. It turns out that $\spec \H_\ast^H(\Omega (G/H); \cc)$ is a \textit{commutative} group scheme over $\ld{\fr{h}}^\ast[2]\mmod \ld{H}$. It is a bit of a miracle, therefore, that the above isomorphism predicts that the fiber product $\ld{\fr{h}}^\ast[2]\mmod \ld{H} \times_{\ld{M}/\ld{G}[-2\rho]} \ld{\fr{h}}^\ast[2]\mmod \ld{H}$ should be a commutative group scheme for every $\ld{M}$ as in \cite{bzsv} which is dual to an affine homogeneous spherical variety. Dually, it also might seem like a bit of a miracle that $\H_\ast^H(\Omega (G/H); \cc)$ is a commutative ring (after all, the Pontryagin product on the homology of a $1$-fold loop space such as $\Omega(S^i \vee S^j)$ will not be commutative!); but this, at least, has a nice explanation coming from the Deligne conjecture. See \cref{cor: loop homology E2}.
        
        The \textit{nonequivariant} homology of $\Omega(G/H)$ has been studied by many authors using Morse theory, at least in the case of symmetric spaces; see, e.g., \cite{bott-samelson-morse-symmetric, ziller-loops-globally-symmetric}.
    \end{enumerate}
    Together, the properties (b) and (c) of $\kappa$ imply that one can recover $\ld{M}$ from the $\spec \H_\ast^H(\Omega (G/H); \cc)$-action on $\ld{\fr{h}}^\ast[2]\mmod \ld{H}$, which ultimately leads to the proof of \cref{intro: thm bzsv rk 1}.
\end{strategy}
It is the isomorphism of (c) which we will establish in the rank $1$ case through case-by-case analysis, since the spaces $S^n$, $\CP^n$, $\HHP^n$, and $\OP^2$ form a finite list of such affine homogeneous spherical varieties up to finite covers (see \cite{ahiezer-rk-1}). 
Although most of the rank $1$ cases behave quite similarly to each other, each case showcases some interesting basic homotopy-theoretic facts (see \cref{table: topology and dualities for rank 1 spherical varieties}). In fact, we will establish the isomorphism of (c) even for homology with coefficients in $\Z[1/2]$, and in some cases with coefficients in $\ku$ (see the next section).
\begin{example}[(Spherical) geometrized spherical harmonics for $G/H = \SO_3/\SO_2$]
    Let us illustrate (c) in the case of the spherical $\PGL_2$-variety $\PGL_2/\GG_m$ (so $G = \PGL_2 \cong \SO_3$ and $H = \GG_m \cong \SO_2$). The Hopf fibration gives a homotopy equivalence $(\PGL_2/\GG_m)(\cc) \simeq S^2$, so (c) reduces to computing $\H^{S^1}_\ast(\Omega S^2; \cc)$. The Borel-equivariant analogue of this computation is quite simple: there is a homotopy fixed points spectral sequence
    $$E_2^{\ast,\ast} \cong \H_\ast(\Omega S^2; \cc) \otimes_\cc \H^\ast(BS^1; \cc) \Rightarrow \pi_\ast C_\ast(\Omega S^2; \cc)^{hS^1},$$
    with a single $d_2$-differential. This spectral sequence degenerates on the $E_3$-page, and gives an isomorphism
    $$\pi_\ast C_\ast(\Omega S^2; \cc)^{hS^1} \cong \cc\pw{x}[b]/bx,$$
    where $|b| = 2$ and $|x| = -2$. Replacing the left-hand side by $\H^{S^1}_\ast(\Omega S^2; \cc)$ simply has the effect of making $x$ into a polynomial (as opposed to power series) variable.

    Ignoring degrees for a moment, write $\AA^1 = \spec \H^\ast_{S^1}(\ast; \cc)$, and let $\kappa: \AA^1 \to T^\ast(\AA^2)$ denote the map sending $x$ to the point $(x,1)$ in the cotangent fiber over $(1,0)\in \AA^2$. If we equip $T^\ast(\AA^2)$ with its natural $\SL_2$-action (coming from the $\SL_2$-action on $\AA^2$), one can compute that there is an isomorphism
    $$\AA^1 \times_{T^\ast(\AA^2)/\SL_2} \AA^1 \cong \spec \cc[x,b]/bx,$$
    and hence an (ungraded) isomorphism 
    $$\spec \H^{S^1}_\ast(\Omega S^2; \cc) \cong \AA^1 \times_{T^\ast(\AA^2)/\SL_2} \AA^1.$$
    The right-hand side can be equipped with a grading such that the above isomorphism is one of graded schemes, which gives (c). In the case of $\PGL_2$, this leads to an equivalence
    $$\Shv^{c}_{\PGL_2\pw{t}}(\PGL_2\ls{t}/\GG_m\ls{t}; \cc) \simeq \Perf(T^\ast[2](\AA^2[2,0])/\SL_2[-2\rho]).$$
\end{example}

Before proceeding to the $\ku$-theoretic generalization, let us mention that many aspects of the Ben-Zvi--Sakellaridis--Venkatesh conjecture can be understood from the perspective of Hochschild (co)homology.
For the purposes of the introduction, it will be convenient to Borel-complete, i.e., to work with the ring $\pi_\ast \cc[\Omega (G/H)]^{hH}$. 
\begin{remark}
    The process of Borel-completion above amounts to replacing equivariant cochains $C^\ast_G(\ast; \cc)$ with the Borel-equivariant cochains $\cc^{hG} = C^\ast(BG; \cc)$ (given by the homotopy fixed points of $G$ acting trivially on $\cc$). Just as there is an isomorphism $\spec C^\ast_G(\ast; \cc) \cong {\ld{\g}}^\ast[2]\mmod \ld{G}$, there is an isomorphism $\spf \cc^{hG} \cong \widehat{\ld{\g}}^\ast[2]\mmod \ld{G}$. This might seem like a minor distinction, but it is in fact a crucial one (which becomes much more significant with $\ku$-coefficients).
\end{remark}
\begin{propo}[\cref{cor: thm 1 of BF for G mod H}, \cref{cor: loop homology E2}]
    The $\cc$-algebra $\cc[\Omega (G/H)]^{hH}$ can be identified with the relative Hochschild cohomology $\HC(\cc^{hH}/\cc^{hG})$ of the map $\cc^{hG} \to \cc^{hH}$ (this was already observed in \cite[Remark A.6]{grg-reg}). In particular, the Hochschild-Kostant-Rosenberg theorem implies that there is an isomorphism
    $$\spf \H^\ast(\cL (G/H)_{hG}; \cc) \cong T[-1]((\widehat{\ld{\fr{h}}}^\ast[2]\mmod \ld{H})/(\widehat{\ld{\g}}^\ast[2]\mmod \ld{G}),$$
    where the right-hand side is the $1$-shifted relative tangent bundle of the map $\widehat{\ld{\fr{h}}}^\ast[2]\mmod \ld{H} \to \widehat{\ld{\g}}^\ast[2]\mmod \ld{G}$, and the hats denote completion at the origin.
    Moreover, the Deligne conjecture equips $\cc[\Omega (G/H)]^{hH}$ with the structure of an $\E{2}$-$\cc^{hG}$-algebra.
\end{propo}
\begin{remark}
    In \cref{strategy: intro}, we said that \cref{conj: bzsv} often reduces to proving an isomorphism
    $$\spec \H^H_\ast(\Omega(G/H); \cc) \cong \ld{\fr{h}}^\ast\mmod \ld{H} \times_{\ld{M}/\ld{G}} \ld{\fr{h}}^\ast\mmod \ld{H},$$
    where we are ignoring gradings.
    The right-hand side is a group scheme over $\ld{\fr{h}}^\ast\mmod \ld{H}$, which we denote by $\ld{J}_X$; it is an analogue for $\ld{M}$ of the regular centralizer group scheme.\footnote{Perhaps it is better denoted by $\ld{J}_{\ld{M}/\ld{G}}$; but we are going to be viewing $\ld{J}_X$ as an object constructed from $X$, and then it is nontrivial that there is any relationship between $\ld{J}_X$ and $\ld{M}/\ld{G}$. Hence the notation $\ld{J}_X$.}
    If the left-hand side in the above isomorphism is replaced by $\spf \pi_\ast C^\ast_H(\Omega(G/H); \cc)^\ast$, i.e., the dual of equivariant \textit{co}homology, the preceding proposition says that the right-hand side must be replaced by $T^\ast[1](({\ld{\fr{h}}}^\ast\mmod \ld{H})/({\ld{\g}}^\ast\mmod \ld{G}))$; this is simply the \textit{Lie algebra} of $\ld{J}_X$. In other words, the dual of equivariant cohomology allows one to access the Lie algebra of $\ld{J}_X$; to understand $\ld{J}_X$ itself involves a ``decompletion'', which in homotopy theory is given by working with equivariant homology itself (since this is a predual of equivariant cohomology). 
    
    Summarizing, the various completions appearing in this article can, roughly speaking, be interpreted as follows:
    \begin{itemize}
        \item There is an isomorphism between $\H^H_\ast(\Omega(G/H); \cc)$ and the ring of functions on $\ld{J}_X = \ld{\fr{h}}^\ast\mmod \ld{H} \times_{\ld{M}/\ld{G}} \ld{\fr{h}}^\ast\mmod \ld{H}$.
        \item The $\H^\ast_H(\ast; \cc)$-linear dual of $\H^\ast_H(\Omega(G/H); \cc)$ can be identified with the ring of functions on the completion of $\ld{J}_X$ at the fibers of the map $\ld{J}_X \to \ld{\fr{h}}^\ast\mmod \ld{H}$; or equivalently, with the Lie algebra of $\ld{J}_X$.
        \item The Borel-equivariant homology $\pi_\ast \cc[\Omega(G/H)]^{hH}$ can be identified with the ring of functions on the completion of $\ld{J}_X \to \ld{\fr{h}}^\ast\mmod \ld{H}$ at the origin of $\ld{\fr{h}}^\ast\mmod \ld{H}$.
    \end{itemize}
\end{remark}
The above result can be viewed as a relative version of \cite[Theorem 1]{bf-derived-satake}. In a sense, most of this article can be regarded as an attempt to understand the decompletion of this Hochschild-Kostant-Rosenberg isomorphism.
As a perhaps helpful guide, a general principle about equivalences of the form conjectured in \cite{bzsv} lead to some analogies between topology and algebra, a limited collection of which we have recorded in \cref{table: analogies}.
\begin{table}[h]
\begin{tabular}{ |c|c|c|c|c|c| } 
 \hline
 A-side/topology & B-side/algebra \\
 \hline
 Spherical $G$-variety $X$ & $\ld{G}$-scheme $\ld{X}$ \\
 Free loop space $\cL X$ & (Twisted) cotangent bundle $T^\ast_\psi(\ld{X})$ \\
 Based loop space $\Omega X$ & Cotangent fiber \\
 Sheaves on $X\ls{t}/G\pw{t}$ with $\QQ$-coefficients & Perfect complexes on $T^\ast_\psi(\ld{X})/\ld{G}$ \\
 \hline
\end{tabular}
\vspace{1cm}
\caption{Analogies between topology and algebra. Passing from the left-hand to the right-hand column is roughly implemented by (rational) cohomology and the Hochschild-Kostant-Rosenberg theorem.}
\label{table: analogies}
\end{table}

The above relationship to Hochschild cohomology allows us to make a basic observation about the structure of a Hamiltonian $\ld{G}$-space on $\ld{M}$. First, we note the following, which is closely related to \cite[Section 5.2]{teleman-icm} and \cite[Section 5.1.5]{sakellaridis-icm}. We have ignored homological shifts/gradings for notational simplicity.
\begin{propo}[\cref{prop: lagrangian correspondence}]
    Let $\ld{J}_{\ld{G}}$ denote the group scheme of regular centralizers of $\ld{G}$, so that 
    $$\ld{J}_{\ld{G}} \cong \ld{\g}^\ast\mmod \ld{G} \times_{\ld{\g}^\ast/\ld{G}} \ld{\g}^\ast\mmod \ld{G},$$
    and similarly, let $\ld{J}_{\ld{H}}$ denote the group scheme of regular centralizers of $\ld{H}$.
    Then, there is a Lagrangian correspondence (interpreted in a derived sense):
    $$\xymatrix{
        & \ld{J}_{\ld{G}} \times_{\ld{\g}^\ast\mmod \ld{G}} \ld{\fr{h}}^\ast\mmod \ld{H} \ar[dl] \ar[dr] & \\
        \ld{J}_{\ld{H}} & & \ld{J}_{\ld{G}}.
    }$$
    Furthermore, there is an isomorphism
    $$(\ld{J}_{\ld{G}} \times_{\ld{\g}^\ast\mmod \ld{G}} \ld{\fr{h}}^\ast\mmod \ld{H}) \times_{J_{\ld{H}}} \ld{\fr{h}}^\ast\mmod \ld{H} \cong \ld{\fr{h}}^\ast\mmod \ld{H} \times_{\ld{M}/\ld{G}} \ld{\fr{h}}^\ast\mmod \ld{H}.$$
\end{propo}
This has an interesting consequence; let us state its ``morally correct'' form (but see \cref{cor: BJ and lag corr} for a precise statement):
\begin{corno}[\cref{cor: BJ and lag corr}]
    Let $\ld{M}$ denote the Hamiltonian $\ld{G}$-space dual to $G/H$, and assume that there is a ``Kostant slice'' $\kappa_{\ld{M}}: \ld{\fr{h}}^\ast\mmod \ld{H} \to \ld{M}$. Let $\ld{M}^\reg$ denote the $\ld{G}$-orbit of this Kostant slice. Similarly, let $\ld{\cM}$ denote the Hamiltonian $\ld{G} \times \ld{H}$-space dual to $(G \times H)/H^\mathrm{diag}$. Again assume that there is a Kostant slice $\kappa_{\ld{M}}: \ld{\fr{h}}^\ast\mmod \ld{H} \to \ld{\cM}$, and let $\ld{\cM}^\reg$ denote its $\ld{G} \times \ld{H}$-orbit. Then there is 
    an isomorphism 
    $$\ld{\cM}^\reg/\ld{H} \cong \ld{\g}^{\ast, \reg} \times_{\ld{\g}^\ast\mmod \ld{G}} \ld{\fr{h}}^\ast\mmod \ld{H},$$
    and 
    a Lagrangian correspondence
    \begin{equation*}
        \xymatrix{
        & \ld{\cM}^{\reg}/(\ld{G} \times \ld{H}) \ar[dl] \ar[dr] & \\
        \ld{\fr{h}}^{\ast, \reg}/\ld{H} & & \ld{\g}^{\ast, \reg}/\ld{G},
        }
    \end{equation*}
    as well as a Cartesian square
    $$\xymatrix{
    \ld{M}^\reg \ar[d] \ar[r] & \ld{\cM}^{\reg} \ar[d] \\
    \ld{\fr{h}}^{\ast}\mmod \ld{H} \ar[r]_-\kappa & \ld{\fr{h}}^{\ast, \reg}.
    }$$
\end{corno}
The Lagrangian correspondence of the preceding proposition, i.e., the Hamiltonian $\ld{G} \times \ld{H}$-space $\ld{\cM}$, ``controls'' transfer for $H \subseteq G$, in the sense that it can be thought of as a morphism from $\ld{\fr{h}}^{\ast}/\ld{H}$ to $\ld{\g}^{\ast}/\ld{G}$ in the algebraic (and $1$-shifted) avatar of Weinstein's category of symplectic varieties and Lagrangian correspondences between them. As explained in \cref{rmk: dirichlet and neumann}, the Cartesian square of the above result follows more generally from the perspective of Dirichlet and Neumann boundary conditions being swapped under Langlands duality. 
\begin{example}
    In \cref{ex: ggp}, we describe a curious calculation involving explicit descriptions of the regular centralizer group schemes for groups in the Gan-Gross-Prasad period. Namely, if $H\subseteq G$ is the inclusion $\SO_{2n} \subseteq \SO_{2n+1}$, so that $\ld{H} = \SO_{2n}$ and $\ld{G} = \Sp_{2n}$, then $\ld{\cM}^\ddag$ can be identified with $\Hom(\std_{2n}, \std_{2n})$.
\end{example}

A result of Safronov's from \cite{safronov-cs} can be used to translate the desired Hamiltonian $\ld{G}$-structure on $\ld{M}$ into the language of shifted symplectic geometry (\`a la \cite{ptvv}). This translation in turn predicts that equipping $\ld{M}$ with the structure of a Hamiltonian $\ld{G}$-space should, in particular, imply that the map
\begin{equation}\label{eq: intro coisotropic}
    \ld{\fr{h}}^\ast[2]\mmod \ld{H} \times_{\ld{M}/\ld{G}[-2\rho]} \ld{\fr{h}}^\ast[2]\mmod \ld{H} \to \ld{\g}^\ast[2]\mmod \ld{G} \times_{\ld{\g}^\ast[2-2\rho]/\ld{G}[-2\rho]} \ld{\g}^\ast[2]\mmod \ld{G} \cong \ld{J}_{\ld{G}}
\end{equation}
is Lagrangian, hence coisotropic (in an appropriate derived sense). It turns out that both \cref{prop: lagrangian correspondence} and the above coisotropicity are extremely general facts which have clean homotopy-theoretic explanations in terms of (higher) Hochschild cohomology.

\begin{example}
    In the special case when $G$ is replaced by $G \times G$ and $H$ is replaced by the diagonal embedding of $G$ (so that the associated homogeneous variety is just $(G\times G)/G^\mathrm{diag} \cong G$), the Hochschild cohomology $\HC(\cc^{hG}/\cc^{h(G\times G)})$ can be identified with the $\E{2}$-center $\dZ_\E{2}(\cc^{hG}/\cc)$. The Deligne conjecture equips $\cc[\Omega G]^{hG} \simeq \dZ_\E{2}(\cc^{hG}/\cc)$ with the structure of an $\E{3}$-$\cc$-algebra; this $\E{3}$-algebra structure is closely related to the $\E{3}$-monoidality of the derived geometric Satake equivalence.

    Note that $\pi_\ast \cc[\Omega G]^{hG}$ is equipped with the structure of a (graded) Poisson algebra whose Poisson bracket has weight $2$; this structure exists on the homotopy of \textit{any} $\E{3}$-algebra. It turns out that this Poisson structure in fact comes from a symplectic form (of weight $2$) on $\spf \pi_\ast \cc[\Omega G]^{hG}$. Moreover, \cref{strategy: intro} says, in particular, that there is an isomorphism
    $$\spf \pi_\ast \cc[\Omega G]^{hG} \cong \widehat{\ld{\g}}^\ast[2]\mmod \ld{G} \times_{\widehat{\ld{\g}}^\ast[2-2\rho]/\ld{G}[-2\rho]} \widehat{\ld{\g}}^\ast[2]\mmod \ld{G}.$$
\end{example}
Returning to relative Langlands, one can also show that if $H\subseteq G$ is a subgroup, $\cc[\Omega (G/H)]^{hH}$ admits the structure of an $\E{2}$-$\cc[\Omega G]^{hG}$-algebra.
The coisotropicity of \cref{eq: intro coisotropic} translates into the requirement that the natural map
$$\spec \H^H_\ast(\Omega (G/H); \cc) \to \spec \H^G_\ast(\Omega G; \cc)$$
is coisotropic. One can directly prove its Borel-completed analogue:
\begin{propo}[\cref{cor: desired map is coisotropic} and \cref{rmk: coisotropic correspondence}]
    The natural map
    $$\spf \H_H^\ast(\Omega (G/H); \cc)^\vee \to \spf \H_G^\ast(\Omega G; \cc)^\vee$$
    is coisotropic (in an appropriate derived sense).
\end{propo}
This result turns out to be a simple consequence of the fact that $\cc[\Omega (G/H)]^{hH}$ admits the structure of an $\E{2}$-$\cc[\Omega G]^{hG}$-algebra, and a general property of $\E{n}$-centers as established in \cite[Theorem 1.1]{francis}. We hope that further study of the relative Langlands program from the perspective of Hochschild (co)homology might shed more light into some of the structures predicted in \cite{bzsv}.

\subsection{$\ku$-theoretic aberrations}

In the course of proving \cref{intro: thm bzsv rk 1}, or even the derived geometric Satake equivalence, the reader will likely observe that many components of the proof do not depend very heavily on the particular choice of coefficient ring for the $\infty$-category of constructible sheaves. In particular, calculations such as that of $\H_\ast^G(\Omega G; \cc)$ (to construct the homomorphism \cref{eq: map from Loops G to dual group}) work equally well with $G$-equivariant $\cc$-(co)homology replaced by any well-behaved equivariant generalized cohomology theory. Motivated by this observation, our second goal in this article is to suggest that the (geometric) Langlands program should admit a generalization to sheaves with coefficients in more ``exotic'' rings, such as the sphere spectrum or complex cobordism. We will discuss the conceptual role of these coefficients in the next section.

Establishing this generalized form of the Langlands program is rather tricky, and so our focus in this article will be on the simpler example of \textit{connective complex K-theory}. Our focus is on this particular example for at least two reasons: first, it is a general principle in homotopy theory that statements about ordinary rational/integral cohomology which admit analogues for connective K-theory will likely admit generalizations to other complex-oriented spectra; second, it is mostly psychological, in that proving analogues for equivariant elliptic cohomology, etc., would requires further technical setup, and distracts from the main features of Langlands duality.\footnote{A third reason, in keeping with the epigraph of our previous article \cite{grg-reg}, is yet another quote of \sout{J. F. Adams} $\mathrm{E}_8$ from \cite{adams-e8-letter}: ``[To] consider the question [of torsion in the cohomology of $\mathrm{E}_8$] at all reveals a certain preoccupation with ordinary cohomology. Any impartial observer must marvel at your obsession with this obscure and unhelpful invariant. The author, like all respectable Lie groups, is much concerned to present a decorous and seemly appearance to the eyes of K-theory...'' (It should be said immediately that we do not study $\mathrm{E}_8$ in this article.) I do not have such strong feelings against ordinary cohomology, but the general thrust of this quote still applies: going from ordinary cohomology to K-theory should reveal deeper structures under Langlands duality.} An analogue of the derived geometric Satake equivalence with coefficients in periodic complex K-theory and elliptic cohomology was proved in \cite{grg-reg}.

If $G$ is a compact Lie group, Atiyah and Segal defined $G$-equivariant complex K-theory $\KU_G$ in \cite{segal-equiv-KU, atiyah-segal-original}: this is a generalized cohomology theory, viewed as a spectrum in the sense of homotopy theory, which classifies $G$-equivariant vector bundles on finite $G$-spaces. Direct sum and tensor products of $G$-equivariant vector bundles equips $\KU_G$ with the structure of a \textit{ring} spectrum; in fact, it is an $\Eoo$-ring, meaning (for instance) that the multiplication on cohomology can be refined by Adams operations.
Despite its definition, the geometric interpretation of cocycles for equivariant K-theory as equivariant vector bundles will play no role below. Two examples will play an important conceptual role:
\begin{example}
    When $G$ is the trivial group, $\KU_G$ is simply periodic complex K-theory $\KU$, and Bott periodicity gives a graded isomorphism $\pi_\ast \KU \cong \Z[\beta^{\pm 1}]$ with the Bott class $\beta$ in weight $2$. 
    On the other hand, when $G$ is a connected compact Lie group with complex representation ring $R_\cc(G)$, the coefficient ring $\pi_\ast \KU_G$ is the tensor product $R_\cc(G) \otimes_\Z \Z[\beta^{\pm 1}]$. In particular, if $G$ is a torus $T$, then $\spec \pi_\ast \KU_T$ is the corresponding algebraic torus $T_{\Z[\beta^{\pm 1}]}$ over $\Z[\beta^{\pm 1}]$.
\end{example}

Nonequivariant complex K-theory is in some sense the simplest generalized cohomology theory which is not just ordinary integral cohomology. In fact, the oft-cited analogies between them are more than coincidental: 
\begin{recall}
    There is an $\Eoo$-ring $\ku$ called \textit{connective} complex K-theory such that there is a graded isomorphism $\pi_\ast \ku \cong \Z[\beta]$. If we set $\beta = 0$, this $\Eoo$-ring simply degenerates to the Eilenberg-Maclane spectrum $\Z$ representing ordinary integral cohomology; and if we invert $\beta$, it recovers periodic complex K-theory. In other words, $\ku$ interpolates between $\Z$ and $\KU$, and can be viewed as a one-parameter deformation of the ring $\Z$ in a homotopy-theoretic direction (namely, along the Bott class $\beta$).
\end{recall}

If $G$ is a compact Lie group, one can also construct an $\Eoo$-ring $\ku_G$ called $G$-equivariant connective K-theory which interpolates between $G$-equivariant integral cohomology and $G$-equivariant (periodic) complex K-theory. In a precise sense (known in the homotopy-theoretic literature as \textit{complex-oriented/abelian descent}), the $\Eoo$-ring $\ku_G$ is determined by the $\Eoo$-rings $\ku_T$ for compact abelian Lie groups $T$, which are in turn determined by the $\Eoo$-ring $\ku_{S^1}$.
\begin{example}
    Since $\ku_{S^1}$ interpolates between $S^1$-equivariant integral cohomology and $\KU_{S^1}$, and there are isomorphisms $\spec \H^\ast_{S^1}(\ast; \Z) \cong \GG_a(2)$ (with coordinate in weight $-2$) and $\spec \pi_\ast \KU_{S^1} \cong \GG_{m,\Z[\beta^{\pm 1}]}$, one expects $\spec \pi_\ast \ku_{S^1}$ to interpolate between $\GG_m$ and $\GG_a$. In fact, equivariant connective K-theory is concocted so that $\spec \pi_\ast \ku_{S^1}$ is the canonical degeneration from $\GG_m$ to $\GG_a$:
    $$\spec \pi_\ast \ku_{S^1} \cong \spec \Z[\beta, x, \tfrac{1}{1+\beta x}],$$
    where $x$ is in weight $-2$. We will denote this group scheme (where the group structure makes $1 + \beta x$ into a grouplike element) by $\GG_\beta$. The case of a general compact abelian Lie group $T$ is a straightforward generalization:
    $$\spec \pi_\ast \ku_T \cong \Hom(\bX^\ast(T), \GG_\beta) =: T_\beta.$$
    Note that when $\beta = 0$, this group scheme is just $\fr{t}(2)$; and when $\beta$ is inverted, this group scheme is $T_{\Z[\beta^{\pm 1}]}$.
    This story is discussed further in \cref{sec: ku big section}.
\end{example}

Very simply, the effect of studying derived geometric Satake with coefficients in $\ku$ (instead of coefficients in $\Z$) is that the dual \textit{group} remains unchanged, and every appearance/consequence of the Cartan subalgebra $\fr{t}$ in the ``classical'' story is replaced by the group scheme $T_\beta$ over $\Z[\beta]$. In order to make this more precise, let us explain a $\ku$-theoretic analogue of the derived geometric Satake equivalence as proved in \cref{thm: ku derived satake}.
\begin{setup}\label{setup: intro Gbeta}
    Write $\sh^{1/2} \Z[\beta]$ to denote the polynomial $\Eoo$-$\Z$-algebra where $\beta$ lives in homological degree $2$, so that it is obtained as a shearing of the graded ring $\pi_\ast(\ku)$.
    \begin{itemize}
        \item There is a $\sh^{1/2} \Z[\beta]$-linear $\infty$-category $\Shv^{c,\Sat}_{(G \times G)\pw{t}}(G\ls{t}; \ku)^\faux$ of ``$G\pw{t} \times G\pw{t}$-equivariant sheaves of $\ku$-modules on $G\ls{t}$''. The definition of this $\infty$-category is given in \cref{cstr: Shv-Sat LG/H with ku}. However, this definition is extremely \textit{ad hoc}; hence the ``faux''. This $\infty$-category is essentially defined by playing around with the subcategory of locally constant sheaves. We hope it will agree with a ``correct'' definition (see \cref{rmk: true shvsat with ku is hard} for a little more on this point).
        \item Suppose $\ld{G}$ is a group scheme defined over $\Z$. Let $\GG_\beta^\vee$ denote the Cartier dual of $\GG_\beta$, and let $\ld{G}_\beta$ denote the group scheme over $\pi_\ast(\ku) \cong \Z[\beta]$ given by $\Hom(\GG_\beta^\vee, G_{\Z[\beta]})$. We will view $\ld{G}_\beta$ as a \textit{$\beta$-deformation} of $\ld{G}$. The quotient stack $\ld{G}_\beta/\ld{G}$ is related to the Hochschild-Kostant-Rosenberg filtration for the quotient stack $B\ld{G}$.
    \end{itemize}
\end{setup}
\begin{theorem}[Derived geometric Satake with $\ku$-theoretic coefficients; \cref{thm: ku derived satake}]\label{intro: derived ku satake}
    Let $G$ be a simply-laced and connected Lie group, and invert the order of the Weyl group $W$ (for simplicity). Then there is a $\sh^{1/2} \Z[1/|W|, \beta]$-linear equivalence
    $$\Shv^{c,\Sat}_{(G\times G)\pw{t}}(G\ls{t}; \ku)^\faux \simeq \Perf(\sh^{1/2} G(-2\rho)_\beta/\ld{G}(-2\rho)),$$
    where $G$ is the Chevalley split form of $G$ over $\Z[1/|W|]$, and $\ld{G}$ acts on $G$ by conjugation\footnote{That this makes sense is thanks to the simply-laced assumption on $G$.}.
\end{theorem}
\begin{remark}
    Upon setting $\beta = 0$, we have 
    $$G(-2\rho)_\beta|_{\beta = 0} = \Hom(\GG_a^\vee(2), G(-2\rho)) = \Hom(\hat{\GG}_a^\sharp(2), G(-2\rho)) = \g(2-2\rho);$$
    here, $\hat{\GG}_a^\sharp$ denotes the divided power hull of the origin in $\GG_a$, further completed at the divided power filtration. Because $G$ is simply-laced, there is an $\ld{G}$-equivariant isomorphism $\g(2-2\rho) \cong \ld{\g}^\ast(2-2\rho)$. It follows that upon setting $\beta = 0$, the left-hand side of \cref{intro: derived ku satake} becomes $\Shv^{c,\Sat}_{(G\times G)\pw{t}}(G\ls{t}; \Z[1/N])$, and the right-hand side becomes $\ld{\g}^\ast[2-2\rho]/\ld{G}[-2\rho]$. In other words, \cref{intro: derived ku satake} just reduces to the derived geometric Satake equivalence (in the simply-laced and connected case).
    
    On the other hand, upon inverting $\beta$, \cref{intro: derived ku satake} is related to the $\KU$-theoretic derived geometric Satake equivalence of \cite{grg-reg}, since
    $$G(-2\rho)_\beta|_{\beta^{-1}} = \Hom(\GG_m^\vee, G) = \Hom(\ul{\Z}, G) = G.$$
    All objects on the right hand sides of the above displayed isomorphisms are to be understood as base-changed from $\Z$ to $\Z[\beta^{\pm 1}]$ (we omitted this from the notation for readability): since $\beta$ lives in weight $2$ and is invertible, we may ignore the $-2\rho$-shift.
    In particular, \cref{intro: derived ku satake} implies that the $\beta$-adic filtration on $\Shv^{c,\Sat}_{(G\times G)\pw{t}/Z(G)\pw{t}}(G\ls{t}; \KU)^\faux$ corresponds to the Hochschild-Kostant-Rosenberg filtration on the free loop space of the quotient stack $B\ld{G}$.
\end{remark}
\begin{example}\label{ex: intro torus ku}
    Let us illustrate \cref{intro: derived ku satake} in the case when $G$ is a torus $T$. Identifying $\cL T$ with $T \times \Omega T$, we see that the $\infty$-category $\Shv^{c,\Sat}_{(G\times G)\pw{t}}(G\ls{t}; \ku)^\faux$ is just $\Shv^c_T(\Omega T; \ku)^\faux$. However, $\Omega T$ is simply the discrete set $\bX_\ast(T)$ of cocharacters of $T$, so that $\Shv^c_T(\Omega T; \ku)^\faux \simeq \bigoplus_{\bX_\ast(T)} \Shv^c_T(\ast; \ku)^\faux$. Almost by construction, there is an equivalence $\Shv^c_T(\ast; \ku)^\faux \simeq \Perf(\sh^{1/2} T_\beta)$. On the other hand, there is an equivalence $\bigoplus_{\bX_\ast(T)} \Mod_{\sh^{1/2} \Z[\beta]} \simeq \Perf(B\ld{T})$. Together, we obtain an equivalence 
    $$\Shv^c_T(\Omega T; \ku)^\faux \simeq \Perf(\sh^{1/2} T_\beta \times B\ld{T}),$$
    which is the right-hand side of \cref{intro: derived ku satake}.
\end{example}

Since the group scheme $G_\beta$ may seem somewhat mysterious, let us mention that it has an extremely concrete interpretation if, for instance, $\ld{G} = \PGL_n$. 
\begin{example}
    When $G = \SL_n$, and $G_\beta = \SL_{n,\beta}$ is the group scheme whose $R$-points (for $R$ being a graded $\Z[\beta]$-algebra) consists of those $n \times n$-matrices $A$ such that 
    $$\frac{\det(I + \beta A) - 1}{\beta} = 0.$$
    Since the derivative of the determinant is the trace, the specialization of this condition to $\beta = 0$ is simply the condition that $A$ is traceless. When $\beta$ is not required to be zero, the above equation should be thought of as a \textit{$\beta$-deformation} of the equation $\Tr(A) = 0$: for example, when $n=2$, it says that $\Tr(A) + \beta \det(A) = 0$. In particular, generically in $\beta$ (that is, when $\beta$ is invertible), $\SL_{n,\beta}$ is isomorphic to $\SL_n \times \Z[\beta^{\pm 1}]$ via the map $A \mapsto I+\beta A$.
    
    More conceptually, $G_\beta$ is a variant of $G$ whose Cartan subgroup is replaced by its $\beta$-deformation, but whose unipotent parts remain unchanged. For instance, there is an analogue of the Bruhat decomposition for $G_\beta$ where the big cell is $N^- \times T_\beta \times N$.
\end{example}

As mentioned before \cref{intro: derived ku satake}, the definition of the $\infty$-category $\Shv^{c,\Sat}_{(G\times G)\pw{t}}(G\ls{t}; \ku)^\faux$ is rather unsatisfying, and a better definition is warranted. The analogue of this desire on the spectral side would be the following:
\begin{conjecture}\label{conj: intro lifting to ku}
    The graded $\Z[\beta, 1/N]$-linear $\infty$-category $\QCoh(G(-2\rho)_\beta/\ld{G}(-2\rho))$ admits a lift to a $\ku[1/N]$-linear $\infty$-category along the degeneration $\ku[1/N] \rightsquigarrow \pi_\ast \ku[1/N] \cong \Z[\beta, 1/N]$.
\end{conjecture}
Note that \cref{conj: intro lifting to ku} does not ask for a lift of $G(-2\rho)$ to $\ku$, or even for a lift of the stack $G(-2\rho)_\beta/\ld{G}(-2\rho)$. A lift as in \cref{conj: intro lifting to ku} does exist if we restrict to quasicoherent sheaves over the regular locus $G_{\beta,\reg}(-2\rho)/\ld{G}(-2\rho)$.

Given \cref{intro: derived ku satake}, one is naturally led to wonder if there is an analogue of the relative Langlands program, and in particular of \cite[Conjecture 7.5.1]{bzsv}, in the context of $\ku$-theoretic coefficients. We do not have a conjecture as precise as that of \textit{loc. cit.} in this setting, but we do prove an analogue of \cref{intro: thm bzsv rk 1} for some affine homogeneous spherical varieties of rank $1$ with $\ku$-theoretic coefficients. For instance, we have the following result describing ``$\ku$-theoretic geometrized spherical harmonics'' for $\PGL_2/\GG_m = \SO_3/\SO_2$:
\begin{example}[\cref{cor: bzsv for CPn} and \cref{rmk: affine closure of SL2 mod Ga and Bbeta}]
    There is an equivalence
    $$\Shv^{c,\Sat}_{\PGL_2\pw{t}}(\PGL_2\ls{t}/\GG_m\ls{t}; \ku)^\faux \simeq \Perf(\sh^{1/2}\ld{V}_\beta/\SL_2(-2\rho)),$$
    where $\ld{V}_\beta$ is the affine closure of $\SL_2 \times^{\GG_a} (\GG_\beta \times \AA^1) \subseteq T^\ast(\AA^2_{\Z[\beta]} - \{0\})$, with $\GG_a$ acting on $\GG_\beta \times \AA^1 \subseteq \AA^2_{\Z[\beta]}$ via $b: (x,y) \mapsto (x,y+bx)$. The scheme $\ld{V}_\beta$ can be viewed as a $\beta$-deformation of $T^\ast(\AA^2)$; it can be explicitly identified with the open subscheme of $\AA^4_{\Z[\beta]} = \AA^4 \times_{\spec \Z} \spec \Z[\beta]$ given by the complement of the hypersurface
    $$1 + \beta (cB - aD) = 0,$$
    where $a$ is in weight $0$, $c$ is in weight $-2$, $B$ is in weight $0$, and $D$ is in weight $-2$.
\end{example}
The structure that seems to emerge out of these considerations is a $\beta$-deformation of the notion of a graded Hamiltonian $\ld{G}$-space, as explained in \cref{sec: everything deforms}; we call this notion a \textit{$\ku$-Hamiltonian space}.
To define this notion, it is convenient to use the language of shifted symplectic geometry as introduced in \cite{ptvv} (see \cref{recall: shifted sympl} for a brief review). 
\begin{recall}
    It was shown in \cite{safronov-cs} that the quotient stack $\ld{\g}(2)/\ld{G}$ admits a $1$-shifted symplectic structure, and that a graded Hamiltonian $\ld{G}$-space $\ld{M}$ is equivalent to the data of a \textit{Lagrangian morphism} $\ld{M}/\ld{G} \to \ld{\g}(2)/\ld{G}$ in the sense of \cite{ptvv}. In particular, the local geometric story of \cite{bzsv} can be restated as the expectation that for certain affine spherical $G_\cc$-varieties $G_\cc/H_\cc$, there is a dual Lagrangian morphism to $\ld{\g}(2)/\ld{G}$ such that \cite[Conjecture 7.5.1]{bzsv} holds.
\end{recall}

\begin{observe}
    Lagrangian morphisms to $G_\beta/\ld{G}$ are by definition \textit{$\ku$-Hamiltonian $\ld{G}$-varieties}. These are roughly graded $\ld{G}$-varieties $\ld{M}_\beta$ over $\Z[\beta]$ with some additional structure and a moment map $\ld{M}_\beta \to G_\beta$. See \cref{def: ku-hamiltonian}.
    Specializing $\beta = 0$, the aforementioned result from \cite{safronov-cs} says that $\ku$-Hamiltonian $\ld{G}$-spaces specialize to graded Hamiltonian $\ld{G}$-spaces. Upon inverting $\beta$, there is an isomorphism $G_\beta[\beta^{-1}]/\ld{G} \cong G/\ld{G} \times_{\spec \Z} \spec \Z[\beta^{\pm 1}]$, and another result from \cite{safronov-cs} then implies that $\ku$-Hamiltonian $\ld{G}$-spaces give rise to \textit{quasi-Hamiltonian} $\ld{G}$-spaces in the sense of \cite{amm-qham}.
\end{observe}
\begin{remark}
    There is some precedence for $\ku$-Hamiltonian spaces in the setting of integrable systems. Here are two such manifestations. Just as the phase space of a classical integrable system (such as the Calogero-Moser system and the Toda lattice) forms a Hamiltonian variety, the phase space of a ``relativistic'' integrable system and its degeneration to a ``nonrelativistic'' system (such as the Ruijsenaars-Schneider system degenerating into the Calogero-Moser system, and the relativistic Toda lattice degenerating into the classical Toda lattice) naturally forms a $\ku$-Hamiltonian space. In the case of the relativistic Toda lattice degenerating into the classical Toda lattice, this essentially follows from \cite{bfm} and the calculations of \cref{thm: ku homology LG and langlands mod center}. We will explain the example of the Ruijsenaars-Schneider system degenerating into the Calogero-Moser system in future work.
    
    Another manifestation of the notion of a $\ku$-Hamiltonian space is as follows. It is a theorem of \cite{amm-qham} that in the differential geometric setting, one can identify quasi-Hamiltonian $G$-spaces with Hamiltonian $LG$-spaces, where $LG$ is the loop group of the compact Lie group $G$. Roughly speaking, a $\ku$-Hamiltonian $G$-space in the differential geometric setting is a Hamiltonian $L_\beta G$-space equipped with a map to $\RR$, where $L_\beta G$ denotes the group over $\RR$ whose fiber over $\beta \in \RR$ is the group of maps from the circle of radius $|\beta|$ to $G$. That is, a $\ku$-Hamiltonian $G$-space is the data of a Hamiltonian $LG$-space along with the data of how it limits to a Hamiltonian $G$-space.
\end{remark}

In \cref{prop: Gbeta mod G is 1-shifted}, we show that $G_\beta/\ld{G}$ admits a $1$-shifted symplectic structure. Moreover, the $\ku$-theoretic calculations in the relative rank one cases of types $A_n$, $D_2$, and $\G_2$ show that there are equivalences of the form
$$\Shv^{c,\Sat}_{G\pw{t}}(G\ls{t}/H\ls{t}; \ku)^\faux \simeq \Perf(\sh^{1/2} \ld{M}_\beta/\ld{G}(-2\rho)),$$
where $\ld{M}_\beta$ is a graded $\ld{G}$-space over $\spec \Z[\beta]$. The action of $\Shv^{c,\Sat}_{(G\times G)\pw{t}}(G\ls{t}; \ku)^\faux$, equivalent to $\Perf(\sh^{1/2} G(-2\rho)_\beta/\ld{G}(-2\rho))$, on the left-hand side of the above equivalence suggests that there is a Lagrangian morphism $\ld{M}_\beta/\ld{G} \to G_\beta/\ld{G}$. One might therefore hope that the local geometric story of \cite{bzsv} admits a $\ku$-theoretic analogue.

\begin{slogan}\label{slogan: ku relative langlands}
    For certain affine spherical $G_\cc$-varieties $G_\cc/H_\cc$, there is a dual Lagrangian morphism $\ld{M}_\beta/\ld{G} \to G_\beta/\ld{G}$ (that is, a $\ku$-Hamiltonian space $\ld{M}_\beta$) such that an analogue of \cite[Conjecture 7.5.1]{bzsv} holds.
\end{slogan}

Unfortunately, the structure theory of quasi-Hamiltonian $\ld{G}$-spaces does not seem to be as well-developed as that of the theory of Hamiltonian $\ld{G}$-spaces, so it is hard at the current moment to make conjectures as precise as those in \cite{bzsv} regarding the nature of these $\ku$-Hamiltonian $\ld{G}$-spaces. We do, however, strongly believe that such a theory will play an important role in understanding $\ku$-theoretic deformations of the geometric Langlands program. 

In any case, we can at least explore consequences of \cref{slogan: ku relative langlands}. For instance, one of the consequences of \cref{slogan: ku relative langlands} is an analogue of \cref{conj: generalized kostant slice}, which says that there is an isomorphism $\ld{M}_\beta\mmod \ld{G} \cong T_{H,\beta}\mmod W_H$ of invariant-theoretic quotients, and that there is a Kostant section $T_{H,\beta}\mmod W_H \to \ld{M}_\beta$.
As in the discussion surrounding \cref{eq: intro coisotropic} above, equipping $\ld{M}_\beta$ with the structure of a $\ku$-Hamiltonian $\ld{G}$-space should furthermore imply that the map
\begin{equation}\label{eq: intro ku coisotropic}
    T_{H,\beta}\mmod W_H \times_{\ld{M}_\beta/\ld{G}} T_{H,\beta}\mmod W_H \to T_\beta\mmod W \times_{G_\beta/\ld{G}} T_\beta\mmod W
\end{equation}
is Lagrangian, hence coisotropic (in an appropriate derived sense). Here, we have ignored homological shifts, for simplicity. Under a putative $\ku$-theoretic version of the relative derived Satake equivalence, \cref{eq: intro ku coisotropic} identifies with the map
$$\spec \ku^H_\ast(\Omega(G/H)) \to \spec \ku^G_\ast(\Omega G),$$
so one wishes to see that this map is coisotropic.
Again, exactly as in the discussion surrounding \cref{eq: intro coisotropic}, this can be explained using homotopy-theoretic ideas regarding Hochschild (co)homology and $\E{n}$-centralizers. Namely, at least upon completion, one can show:
\begin{propo}[\cref{cor: desired map is coisotropic} and \cref{rmk: coisotropic correspondence}]
    The natural map
    $$\spf \ku_H^\ast(\Omega (G/H))^\vee \to \spf \ku_G^\ast(\Omega G)^\vee$$
    is coisotropic (in an appropriate derived sense).
\end{propo}

As the reader will likely observe (and as mentioned in the beginning of this subsection), the proofs of these statements for $\ku$ ultimately use very little on the specific structure of equivariant connective K-theory, and are almost axiomatic in nature: they only rely on certain basic properties of this equivariant cohomology theory.
While proving \cref{intro: thm bzsv rk 1} and its $\ku$-theoretic variant, we will also describe some calculations of independent interest along the way (such as \cref{rmk: Loops2 and BGa} and \cref{rmk: yun zhu centralizer}, which identifies the homology of $\Omega^2 \SU(n)$, even with $\Z$-coefficients, with the cohomology of the classifying stack of a shearing of the group scheme of length $n-1$ Witt vectors over $\Z$).

\subsection{Generalized coefficients}

One long-term goal, of which this project is a part, is to understand a version of the (relative) geometric Langlands program with coefficients in the sphere spectrum. We do not resolve this problem in this article (far from it!), but instead formulate some very na\"ive expectations in \cref{subsec: duality over MU expectations}; I hope that reporting my (meagre) partial progress will motivate further study into this topic. 

In order to motivate why this is a natural question, let us begin with some general remarks about the nature of the geometric Satake equivalence. A longstanding expectation has been that Langlands-type equivalences are of a ``motivic nature''. 
\begin{example}
    In the arithmetic incarnation of the Langlands program, results of this form are often very deep; for instance, V. Lafforgue has conjectured an independence of $\ell$ result for Langlands parametrization in the case of global function fields (see \cite[Conjecture 12.12]{v-lafforgue-18}).
\end{example}
As mentioned earlier in the introduction, it is often simpler to geometrize the archimedean and $p$-adic settings into the setting of complex curves. In this context, we will treat the word ``motivic'' as more a ``way of life'' instead of a precise mathematical word.
For instance, the motivic nature of the geometric Langlands equivalence could be interpreted as the expectation that spectral decompositions should exist for sheaves/automorphic forms valued in (modules) over an (essentially) arbitrary base ring.
\begin{example}
    One reflection of this motivic nature already appears in the geometric Satake equivalence from \cite{mirkovic-vilonen}, which (re)constructs the Chevalley split form of a reductive group scheme. It is also true of the derived geometric Satake equivalence (see, e.g., \cite{homology-langlands, ex-zhu-geom-sat-Z}); namely, \cref{eq: intro bf derived} still holds if the coefficients $\cc$ are replaced by some localization $\Z'$ of $\Z$, in which case $\ld{\g}[2-2\rho]/\ld{G}[-2\rho]$ is replaced by a lift to $\Z'$.
\end{example}

The field of motivic homotopy theory, as introduced by Morel-Voevodsky in \cite{morel-voevodsky}, suggests that the \textit{stable} motivic category is a more refined version of (integral) motives (although one which is perhaps less accessible by the general public). Taking this perspective into account suggests that one can generalize the motivic expectation of geometric Langlands equivalences to also include sheaves with coefficients in ring spectra. The discussion of the preceding section (e.g., \cref{intro: derived ku satake}) shows that this expectation is not implausible: namely, spectral decompositions \textit{do} exist, and their nature is modified according to the behaviour of Chern classes for the ring spectrum.

\begin{example}
    For instance, this relationship to Chern classes is the basic source of the difference between the case of ``ordinary'' derived geometric Satake (whose spectral side is $\ld{\g}[2-2\rho]/\ld{G}[-2\rho]$) and $\KU$-theoretic derived geometric Satake (whose spectral side is $G/\ld{G}$): indeed, Chern classes in integral cohomology and complex K-theory are very different from each other!\footnote{More precisely, the first Chern class is additive in integral cohomology, but is (essentially) multiplicative in complex K-theory. This distinction manifests itself in many ways in other (related) parts of mathematics; for instance, the Todd class appearing in the Grothendieck-Riemann-Roch theorem is just the ratio of the first Chern classes in rational cohomology and K-theory.} In the case of elliptic cohomology with associated elliptic curve $E$ (as studied in \cite{grg-reg}), the spectral side is in turn replaced by the moduli stack $\Bun_{\ld{G}}^{\ss,0}(E^\vee)$ of semistable degree zero $\ld{G}$-bundles on the dual elliptic curve $E^\vee$.
    The importance of Chern classes in geometric Langlands is reflected in our setting in \cref{ex: intro torus ku}, as well as in the classical setting of geometric Satake, where the the Chern class of the determinant line bundle on $\Gr_G$ for $G$ semisimple can be identified with a regular nilpotent element for the dual Lie algebra $\ld{\g}$ (see \cite{ginzburg-langlands} and \cref{prop: ginzburg coh}).
\end{example}

Motivated by this discussion, it is natural to wonder: 
\begin{question}\label{qn: chern and derived geom}
    Is there is an analogue of the derived geometric Satake equivalence with coefficients in an arbitrary ring spectrum $R$ which admits a theory of Chern classes?
\end{question}
Such ring spectra are called \textit{complex-oriented}. Associated to any complex-oriented ring spectrum $R$, one can define a graded ($1$-dimensional) formal group $\Ghat_R$ over $\pi_\ast(R)$ given by $\spf \pi_\ast(R^{hS^1}) = \spf \pi_\ast R^{\CP^\infty_+}$.
\begin{observe}
    Let $R$ be a complex-oriented ring spectrum. A generalization of the derived geometric Satake equivalence along the lines of \cref{intro: derived ku satake} should involve replacing the $1$-dimensional group scheme $\GG_\beta$ by a $1$-dimensional group scheme $\GG_R$ which is related to the complex-oriented structure of $R$. If $T$ is a maximal torus of a reductive algebraic group $G$, the group scheme $T_\GG := \Hom(\bX^\ast(T), \GG_R)$ would play the role of a Cartan subgroup of a ``$\GG_R$-deformation'' of the group scheme $G$.
    
    However, the existence of a theory of Chern classes only grants us access to the \textit{formal} group $\Ghat_R$, as opposed to an honest $1$-dimensional algebraic group. (For instance, in the case of periodic complex K-theory $\KU$, we have $\Ghat_\KU = \hat{\GG}_m$, as opposed to the multiplicative group $\GG_m$.) As discussed in \cite{survey} (motivated by the Atiyah-Segal completion theorem), the data of a decompletion $\GG_R$ of $\Ghat_R$ can be viewed as an algebraic incarnation of a \textit{genuine $S^1$-equivariant} analogue $R_{S^1}$ of $R$. Namely, $\GG_R$ can be understood as the graded group scheme $\spec \pi_\ast^{S^1}(R_{S^1})$ over $\pi_\ast(R)$ underlying $\spec R_{S^1}$. Note that the group structure on $\GG_R$ comes from the coproduct
    $$R_{S^1} \to R_{S^1 \times S^1} \xleftarrow{\sim} R_{S^1} \otimes_R R_{S^1};$$
    in particular, it is important that the assignment $T \mapsto R_T$ from tori satisfy the K\"unneth formula.
    Therefore, one should require the additional data of a genuine equivariant refinement of $R$ in order to answer \cref{qn: chern and derived geom}.
\end{observe}
\begin{remark}
    There is in fact a universal example of a complex-oriented ring spectrum, given by \textit{complex cobordism} $\MU$. This is an $\Eoo$-ring whose origin is geometric in nature (via cobordism classes of stably almost-complex manifolds), but nevertheless exerts strong control over the $\infty$-category of spectra (in a sense described below).
    In similar fashion, there is a universal ring $L_\ast$ carrying a $1$-dimensional formal group equipped with a coordinate; a theorem of Lazard's shows that $L_\ast$ -- called the \textit{Lazard ring} -- is isomorphic to a polynomial algebra on infinitely many generators (which encode the coefficients of the group law in the chosen coordinate).
    In \cite{quillen-formal-gps}, Quillen showed the following profound and deeply influential statement:
    \begin{thmno}[Quillen]
        The map $L_\ast \to \pi_\ast(\MU)$ classifying the $1$-dimensional formal group $\hat{\GG}_\MU$ is an isomorphism. In other words, the universal $1$-dimensional formal group with a coordinate can be identified with the homotopically-defined formal group $\hat{\GG}_\MU$.
    \end{thmno}
    If $G$ is a compact Lie group, there is also a notion of \textit{$G$-equivariant} complex cobordism $\MU_G$, defined using equivariant Thom spaces (and not geometrically via equivariant cobordism, thanks to the failure of equivariant transversality); see \cite{uribe-icm} for a survey.
    However, setting $\GG_\MU = \spec \pi_\ast^{S^1}(\MU_{S^1})$ does \textit{not} produce a $1$-dimensional group scheme! The problem is precisely the failure of the K\"unneth formula for the assignment $T \mapsto \MU_T$. Instead, as explained in \cite{hausmann-global}, the appropriate structure encoded by the assignment $T \mapsto \MU_T$ is that of a \textit{graded group law}.
    This is a functor $\GG$ from abelian compact Lie groups to graded commutative rings satisfying a certain condition which forces $\spec \GG(S^1)$ to behave like a $1$-dimensional group scheme over $\GG(\ast)$ --- but it is \textit{not} a group scheme in the usual sense. (See \cref{def: graded group laws} for further discussion.) In \cite[Theorem C]{hausmann-global}, Hausmann showed that the assignment $T \mapsto \pi_\ast^T \MU_T$ defines the universal graded group law; this can be regarded as an analogue of Quillen's theorem about $\MU$.
\end{remark}
A positive answer to \cref{qn: chern and derived geom} in the universal case of $\MU$ would therefore suggest that if $\ld{G}$ is a (split) reductive algebraic group (over $\Z$, say) with a chosen maximal torus $\ld{T}$, then \textit{every} graded group law $\GG$ defines a ``$\GG$-analogue'' $\ld{G}_\GG$ of $\ld{G}$, where the role of the Cartan subgroup is played by $\ld{T}_\GG := \spec \GG(\bX^\ast(\ld{T}))$. Unfortunately, I do not know how to define such a $\GG$-analogue. If $\GG$ comes from a $1$-dimensional algebraic group, \cref{setup: intro Gbeta} suggests defining $\ld{G}_\GG := \Hom(\GG^\vee, \ld{G})$, where $\GG^\vee$ is the Cartier dual of $\GG$ (but this definition is also somewhat lacking).\footnote{For instance, it cannot be correct if $\GG$ is an elliptic curve, since its Cartier dual is the zero group. However, one could instead define the quotient stack $\ld{G}_\GG/\ld{G}$ as $\Map(\Hom(\GG, B\GG_m), B\ld{G})$; this, too, is not quite correct in the case of an elliptic curve, but for more subtle reasons.} The definition of a graded group law, however, is so general that it is not clear how to define Cartier duals in this context (or even if it should be possible to do so!).
Nevertheless, we propose some expectations in \cref{subsec: duality over MU expectations} about a putative derived geometric Satake equivalence with coefficients in $\MU$, but (as the reader will see) we could not make it very far before getting stuck.
\begin{remark}
    It is quite easy to give a positive answer to \cref{qn: chern and derived geom} in the case of a torus: namely, for any reasonable definition of $T$-equivariant sheaves of $\MU$-modules on (ind-finite) $T$-spaces such that $\Shv_T^c(\ast; \MU) \simeq \Perf(T_{\GG_\MU})$, there will be an equivalence
    $$\Shv_{T \times T}^c(\cL T; \MU) \simeq \Perf(T_{\GG_\MU} \times B\ld{T}).$$
\end{remark}

\begin{remark}
    The sphere spectrum $S^0$ is \textit{not} complex-oriented, so it might not be clear that there should be an analogue of the derived geometric Satake equivalence with coefficients in $S^0$. Nevertheless, the sphere spectrum admits a ``local'' complex-orientation, in the sense that the unit map $S^0 \to \MU$ behaves as an fpqc cover. More precisely, work of Quillen and Landweber-Novikov suggests that rather than considering $\spec \pi_\ast(S^0)$, it is more appropriate to consider the stack $\Msfg$ of $1$-dimensional formal groups equipped with a square root of the dualizing line bundle; the stacky nature of $\Msfg$ corresponds to the failure of $S^0$ to be complex-oriented. This means that studying the dependence of an analogue of the derived geometric Satake equivalence with coefficients in $\MU$ on the choice of complex orientation should be the first step in understanding some analogue of the derived geometric Satake equivalence with coefficients in $S^0$.
\end{remark}
Most of the aforementioned difficulties come from attempting to view the effect of genuine equivariance on (spectral) algebraic geometry. However, some aspects of Langlands duality do not require working with equivariant cohomology. For instance, recall from \cite{ginzburg-langlands} (see \cref{prop: ginzburg coh}) that there is an isomorphism $\H^\ast(\Omega G; \cc) \cong U(\ld{\g}^e)$ of Hopf algebras, with $e$ being the principal nilpotent element corresponding to the Chern class of the determinant line bundle $\Omega G \to \CP^\infty$. In \cref{subsec: all the expectations}, we prove an analogue of this isomorphism for the sphere spectrum in the case $G = \SU(n)$. This calculation is rather simple, and it would be interesting to prove an analogue for arbitrary compact Lie groups.

%% file: intro/notation.tex
\subsection{Notation}

In writing this article, I discovered that it is extremely easy to fall into grading hell\footnote{Especially as a graduate student!}, and escaping it is a painful task; I hope the following list of conventions is helpful to the reader.
\begin{notation}
    We will {\textit{always}} use homological degrees. For instance, if $X$ is a space, a class $x\in \H^n(X; \QQ)$ in cohomology lives in homological degree $-n$.
\end{notation}
\begin{notation}
    Let $V$ be a finite-dimensional affine space over a (possibly graded) commutative ring $R$, so that $V = \spec \Sym^\ast_R(V^\ast)$. 
    We will denote $V(n)$ to denote the affine space $V$ with weight $n$. Then, we have
    $$V(n) \cong \spec \Sym^\ast_R(V(n)^\ast) \cong \spec \Sym^\ast_R(V^\ast(-n)).$$
    We will also write $\AA^n(i_1, \cdots, i_n)$ to denote the product $\prod_{j=1}^n \AA^1(i_j)$.
\end{notation}
\begin{notation}
    Let $V$ be a finite-dimensional affine space over a (possibly graded) commutative ring $R$, so that $V = \spec \Sym^\ast_R(V^\ast)$. 
    If $n\in \Z$, we will write $V[n]$ to denote the derived $R$-scheme which underlies the graded derived $R$-scheme $V[n](n) = \spec \sh^{1/2} \Sym^\ast_R(V(n)^\ast)$. Note that by \cref{lem: no e2 half shearing}, this definition may not be well-behaved unless $n$ is even (but this will be the case in all examples of interest).
    Note that
    $$\sh^{1/2} \Sym^\ast_R(V(n)^\ast) = \bigoplus_{j\geq 0} \sh^{1/2} V^\ast(-nj) = \bigoplus_{j \geq 0} V^\ast(-nj)[-nj],$$
    so that $V[n](n)$ is the graded derived $R$-scheme where the coordinate lives in degree $n$ and weight $n$. More generally, we will make it a point to distinguish between \textit{weights} and \textit{degrees}.
\end{notation}
\begin{notation}
    We will write $\GG_a^\sharp$ to denote the divided power hull of $\GG_a$ at the origin, and $\hat{\GG}_a^\sharp$ to denote its completion at the divided power filtration.
\end{notation}
\begin{warning}
    Note that $V[n]$ is generally \textit{not} equivalent to the derived $R$-scheme $\spec \Sym^\ast_R(V^\ast[-n])$! For example, suppose $n = -2$, and say that $V = \AA^1_R$ itself. Then the ring of functions $\co_{V[-2]}$ is the polynomial algebra $R[x]$ with $|x| = 2$, but $\Sym^\ast_R(V^\ast[-n]) = \Sym^\ast_R(R[2])$ is isomorphic to the sheared divided power algebra $\bigoplus_{j\geq 0} \Gamma^j_R(R)[2j]$ by d\'ecalage. Another isomorphism along similar lines is $\GG_a[2] \cong B^2 \GG_a^\sharp$. Of course, if $R$ is a $\QQ$-algebra, the two algebras $\co_{V[-2]}$ and $\Sym^\ast_R(R[2])$ are isomorphic to each other, but it is often (psychologically) safer to not make this assumption.
\end{warning}
\begin{notation}
    If $A$ is a ring spectrum with even homotopy groups, one obtains a graded affine scheme $\spec \pi_{\ast} A$. In particular, a class $x\in \pi_{n} A$ defines a map $\spec \pi_{\ast} A \to \AA^1(-n)$, i.e., lives in weight $n$. This is somewhat opposed to standard practice in homotopy theory, where a class in $\pi_{2n} A$ lives in weight $n$ (as opposed to our convention, where it has weight $2n$). 

    On the topic of ring spectra: we will often refer to \textit{$\E{n}$-rings} where $1\leq n \leq \infty$; the reader uninterested in homotopy-coherent subtleties can mostly treat $\E{1}$-rings as associative (differential graded) algebras, $\Eoo$-rings as commutative (differential graded) algebras, and $\E{n}$-rings as commutative (differential graded) algebras with a Poisson bracket of degree $n-1$. These replacements are correct for $\E{n}$-$\QQ$-algebras, but not so otherwise.
\end{notation}
\begin{notation}
    If $G$ is a topological group and $A$ is a ring spectrum, we will sometimes write $A^{hG}$ to denote the cochains $C^\ast(BG; A)$; this is the homotopy fixed points of the trivial $G$-action on $A$. This should be distinguished from the notation $C^\ast_G(\ast; A)$ or $A_G$, which will denote the \textit{genuine} equivariant version of $G$-equivariant $A$-cohomology. It differs from $A^{hG}$ by a completion (namely, $A^{hG}$ is the completion of $C^\ast_G(\ast; A)$ at the fiber of the map $C^\ast_G(\ast; A) \to A$).
\end{notation}
\begin{notation}
    If $X$ is a space and $A$ is a ring spectrum, we will on occasion write $A[X]$ to denote the $A$-chains on $X$; so $A[X] = A \otimes \Sigma^\infty_+ X$. If $G$ acts on $X$, the notation $A[X]^{hG}$ will denote the homotopy fixed points of $G$ acting on $A[X]$, so that $A[X]^{hG}$ can be identified with the Borel-equivariant homology $C_\ast(X; A)^{hG}$. Note, also, that what we call equivariant homology does not seem to be the same thing studied in the equivariant homotopy theory literature. For us, Borel-equivariant (say) homology will mean $C_\ast(X; A)^{hG}$. However, the equivariant homotopy theory literature might use the same term to denote $C_\ast(X_{hG}; A)$. This poses no problems in \textit{cohomology} (where $C^\ast(X; A)^{hG} = C^\ast(X_{hG}; A)$), but these notions crucially differ in the case of homology. For instance, if $X = \Omega Y$ with $G$ acting on $Y$ through based maps, the homotopy fixed points $C_\ast(X; A)^{hG}$ will admit the structure of an $\E{1}$-$A^{hG}$-algebra, but $C_\ast(X_{hG}; A)$ will not have any such multiplicative structure (it will generally be a highly divisible $A^{hG}$-module).
\end{notation}
\begin{notation}
    If $G$ is a compact Lie group acting on a finite space/CW-complex $X$, we will write $\cf_G(X)$ to denote the ``$G$-equivariant $\ku$-cochains of $X$'', i.e., $C^\ast_G(X; \ku)$. Its $\ku_G$-linear dual will be denoted $\cf_G(X)^\vee$; this is the ``$G$-equivariant $\ku$-chains of $X$'', i.e., $C^G_\ast(X; \ku)$. The reader uninterested in $\ku$-theoretic generalizations should simply replace $\cf_G(X)$ by $C^\ast_G(X; \QQ)$, and $\cf_G(X)^\vee$ by $C^G_\ast(X; \QQ)$.
\end{notation}
\begin{notation}
    We will often write $\Z'$ to denote a localization of the ring of integers. This will essentially always mean that the prime $2$ has been inverted. In some instances, it will denote the localization $\Z[1/|W|]$ obtained by inverting the order of a Weyl group.
\end{notation}
\begin{notation}
    The symbol $\HC$ will denote Hochschild cohomology, and $\HH$ will denote Hochschild homology.
\end{notation}
\begin{notation}
    We will always write a double-slash to mean GIT quotients, and not stacky quotients (which will be denoted by a single slash). For instance, $\g(2)\mmod G = \spec \Sym(\g^\ast(-2))^G$, while $\g(2)/G$ is a graded stack over $BG$.
\end{notation}
\begin{notation}
    If $X$ is a scheme over a base ring (often $\cc$), we will write $X\ls{t}$ to mean the functor $R \mapsto X(R\ls{t})$ (and similarly for $X\pw{t}$). The exception to this is that when $X = G/H$ is a homogeneous variety, we will sometimes write $X\ls{t}$ to mean $G\ls{t}/H\ls{t}$ (and similarly $X\pw{t}$ to mean $G\pw{t}/H\pw{t}$); this need not be the same as $(G/H)\ls{t}$, even if $G/H$ is affine. 
\end{notation}
\begin{notation}\label{notation: compact vs complexification}
    On the ``geometric''/topological side, we will often write $G$ to mean a compact Lie group, but also use the same symbol to denote its complexification $G_\cc$; the symbol $\g$ will denote the Lie algebra of the complexification $G_\cc$. Whether $G$ denotes the compact Lie group or its complexification should hopefully be clear from context and not cause confusion. For example, the equivariant cohomology $\H^\ast_G(\ast; k)$ is isomorphic to $\H^\ast_{G_\cc}(\ast; k)$; and similarly the based loop space $\Omega G$ is homotopy equivalent to $\Omega G_\cc$. Sometimes, we will say that $G$ is a compact Lie group and then write symbols such as $G\pw{t}$ or $G\ls{t}$; these should be understood to mean $G_\cc\pw{t}$ or $G_\cc\ls{t}$, respectively.
\end{notation}

%% file: intro/thanks.tex
\subsection{Acknowledgements}

First, I would like to thank David Gepner and Milton Lin for the opportunity to visit Johns Hopkins University in February 2023. 
I am deeply grateful to Yiannis Sakellaridis for extremely helpful conversations, his encouragement and patience, and showing me \cite[Table 1]{sakellaridis-rank-1} when I visited JHU. For me, this table made clear the analogy of Langlands dual L-functions for rank $1$ spherical varieties with the James splitting and the EHP + generalized Hopf fibrations; this was the starting point of the present article, so I owe its entire existence to Yiannis. I was also greatly inspired by Akshay Venkatesh's lectures at the Arizona Winter School in Spring 2022, and I am grateful to him and the organizers of the school for the opportunity to learn such beautiful mathematics. (I also enjoyed the chance to escape from Boston's chilly weather, if only for a week.)

I am also grateful to David Ben-Zvi, Justin Hilburn, Kevin Lin, Mark Macerato, Deven Manam, Andy Neitzke, Sam Raskin, Wyatt Reeves, Yiannis Sakellaridis (again!), Tomer Schlank, Akshay Venkatesh (again!), and Zhiwei Yun for helpful exchanges and enjoyable conversations; I'm especially grateful to Mark Macerato and Charles Fu for giving me helpful feedback on parts of this writeup.
Finally, as always, I am very thankful to my advisors Dennis Gaitsgory and Mike Hopkins for their support and encouragement.

%% file: connective-equiv-k-theory/shearing.tex
\subsection{Shearing and the Koszul sign rule}

The operation of \textit{shearing} will play a vital role in many of our constructions below. Outside of a few sources like \cite{rotinv} and \cite[Proposition 3.3.4]{raksit}, there does not seem to be much literature developing this notion from a homotopy-coherent perspective.
Let us recall the shearing functor $\tilde{\sh}: \Sp^\gr \xar{\sim} \Sp^\gr$ as constructed in \cite[Proposition 3.3.4]{raksit} (it is denoted there by ${\sh}$). Throughout this section, we will let $\Z^\ds$ denote the set of integers viewed as a discrete space.
\begin{construction}\label{cstr: shearing}
    Recall from the universal property of Day convolution that a lax $\E{2}$-monoidal functor $\Sp^\gr \to \Sp^\gr$ is the same data as a lax $\E{2}$-monoidal functor $\Z^\ds \times \Sp^\gr \to \Sp$. Let $f: \Z^\ds \to \Pic(\Sp)$ denote the $\E{2}$-composite
    $$\Z^\ds \xar{\Omega^2(\CP^\infty \to \BU)} \BU \times \Z \xar{J} \Pic(\Sp).$$
    This defines a lax $\E{2}$-monoidal functor via the composite
    $$\Z^\ds \times \Sp^\gr = \Z^\ds \times \Fun(\Z^\ds, \Sp) \xar{f, \ev} \Pic(\Sp) \times \Sp \xar{\otimes} \Sp.$$
    It is not difficult to see that the lax $\E{2}$-monoidal functor $\Sp^\gr \to \Sp^\gr$ constructed above is in fact a strictly $\E{2}$-monoidal equivalence. This functor will be denoted by $\tilde{\sh}$, and will be called shearing. Explicitly, it sends $M_\bull \mapsto M_\bull[2\bull]$, with $\bull$ denoting the weight. 
\end{construction}
\begin{remark}
    In \cite[Proposition 3.10]{framed-e2}, it is shown that the functor $\tilde{\sh}: \Sp^\gr \to \Sp^\gr$ is in fact a \textit{framed} $\E{2}$-monoidal functor. However, it \textit{cannot} be made into an $\E{3}$-monoidal functor (see \cite[Remark 3.11]{framed-e2}).
\end{remark}
Let $\MU$ denote the $\Eoo$-ring of complex cobordism. A first simple observation is the following (the reader uninterested in this generality can simply replace $\MU$ by $\Z$ in the statement).
\begin{lemma}\label{lem: Eoo shear}
    The shearing functor $\tilde{\sh}: \Mod_\MU^\gr \xar{\sim} \Mod_\MU^\gr$ admits a natural symmetric monoidal structure.
\end{lemma}
\begin{proof}
    The same argument as in \cref{cstr: shearing} will show that the shearing functor $\Mod_\MU^\gr \xar{\sim} \Mod_\MU^\gr$ is symmetric monoidal, as long as the map $\Z^\ds \xar{f} \Pic(\Sp) \to \Pic(\MU)$ admits an $\Eoo$-structure.
    Let $J: \BU \times \Z^\ds \to \Pic(\Sp)$ denote the $J$-homomorphism, so that $J$ is an $\Eoo$-map.
    Since $\MU$ is the Thom spectrum of the $\Eoo$-composite $\BU \to \BU\times \Z^\ds \xar{J} \Pic(\Sp)$, it can be understood as the initial $\Eoo$-ring $R$ equipped with a nullhomotopy of the $\Eoo$-map $\BU \xar{J} \Pic(\Sp) \to \Pic(R)$. In particular, there is a commutative diagram of $\Eoo$-maps:
    $$\xymatrix{
    \BU\times \Z = \Omega^\infty \ku \ar[r]^-J \ar[d] & \Pic(\Sp) \ar[d] \\
    \Z^\ds \simeq \Omega^\infty \tau_{\leq 0} \ku \ar[r] & \Pic(\MU),
    }$$
    which proves the desired claim.
\end{proof}
\begin{remark}
    One might wonder whether there is also an $\E{2}$-monoidal structure on the functor $\tilde{\sh}^{1/2}: \Sp^\gr \to \Sp^\gr$ which sends $M_\bull \mapsto M_\bull[\bull]$. In fact, one can already ask whether there is an $\E{2}$-monoidal structure on the functor $\tilde{\sh}^{1/2}: \Mod_\Z^\gr \to \Mod_\Z^\gr$ which sends $M_\bull \mapsto M_\bull[\bull]$.
    The essential difficulty is that of the Koszul sign rule. Namely, suppose that there was an $\E{2}$-monoidal structure on $\tilde{\sh}^{1/2}$. Applying $\tilde{\sh}^{1/2}$ to the graded $\Eoo$-algebra $\Z[x]$ with $x$ in degree zero and weight $1$ would produce a graded $\E{2}$-algebra $\Z[w]$ with $w$ in degree $1$ and weight $1$. The Koszul sign rule forces $w^2 = -w^2$, i.e., $2w^2 = 0$, which is a contradiction. This is one of the basic topological reasons for why we will work with evenly-graded objects throughout this article.
\end{remark}

\begin{notation}
    In the remainder of this article, we will write ${\sh}^{1/2}(M_\bull)$ to denote the underlying \textit{ungraded} spectrum of $\tilde{\sh}^{1/2}(M_\bull) = M_\bull[2\bull]$. 
\end{notation}
Let us note the following related result (which is just a fancy restatement of the Koszul sign rule):
\begin{lemma}\label{lem: no e2 half shearing}
    There is \textit{no} $\E{2}$-map $\Z^\ds \to \Pic(\Z)$ sending $1\mapsto \Z[1]$. 
    However, there is a unique $\E{1}$-map $\Z^\ds \to \Pic(\Sp)$ sending $1 \mapsto S^1$, and hence the method of \cref{cstr: shearing} produces an \textit{$\E{1}$-monoidal} structure on the functor $\tilde{\sh}^{1/2}: \Sp^\gr \to \Sp^\gr$.
\end{lemma}
\begin{proof}
    Let us first show that there is a unique $\E{1}$-map $\Z^\ds \to \Pic(\Sp)$ sending $1 \mapsto S^1$. This is easy, since $\Z^\ds$ is the group completion of the free $\E{1}$-space $\Z_{\geq 0}^\ds$ on a single class: the choice of $S^1\in \pi_0 \Pic(\Sp)$ defines an $\E{1}$-map $\Z_{\geq 0}^\ds \to \Pic(\Sp)$, which factors through $\Z^\ds$ since $\Pic(\Sp)$ is group-complete.

    To show that there is no $\E{2}$-map $\Z^\ds \to \Pic(\Z)$ sending $1\mapsto \Z[1]$, let us make the following observation. There is an fiber sequence of $\Eoo$-spaces
    \begin{equation}\label{eq: space PicZ}
        \tau_{\geq 1} \Pic(\Z) = \BGL_1(\Z) = \RP^\infty \to \Pic(\Z) \to \pi_0 \Pic(\Z) = \Z^\ds,
    \end{equation}
    where the final $\Z^\ds$ is generated by $S^1$. An $\E{2}$-map $\Z^\ds \to \Pic(\Z)$ would give an $\E{2}$-splitting of this fiber sequence, which we claim is impossible. Indeed, let $\picsp(\Z)$ denote the connective spectrum delooping $\Pic(\Z)$; then the above fiber sequence deloops to a cofiber sequence
    \begin{equation}\label{eq: picZ cofiber}
        \Sigma \gl_1(\Z) = \FF_2[1] \to \picsp(\Z) \to \Z.
    \end{equation}
    One can use the $J$-homomorphism and the $\Eoo$-map $\MSO \to \Z$ to show that $\picsp(\Z) \simeq \tau_{\leq 1} \ko$. (Applying $\Omega^\infty$, this amounts to the identification $\Pic(\Z) \simeq (\Z \times \BO)/\BSO$.) Therefore, the boundary map in \cref{eq: picZ cofiber} identifies with the first $k$-invariant of $\ko$, which is given by the composite
    \begin{equation}\label{eq: sq2}
        \Z \to \FF_2 \xar{\Sq^2} \FF_2[2].
    \end{equation}
    In particular, \cref{eq: picZ cofiber} is \textit{not} split as a cofiber sequence of spectra (this can be viewed as a manifestation of the Koszul sign rule). However, the above description of the boundary map also lets us show that \cref{eq: space PicZ} does not split as a fiber sequence of $\E{2}$-spaces. Namely, if \cref{eq: space PicZ} did split as a fiber sequence of $\E{2}$-spaces, then the twice-delooped fiber sequence
    $$K(\Z/2, 3) \to B^2 \Pic(\Z) \to B^2 \Z^\ds = \CP^\infty$$
    would also admit a splitting. But the boundary map $\CP^\infty \to K(\Z/2, 4)$ in this fiber sequence represents the generator of $\H^4(\CP^\infty; \Z/2)$, which is certainly nonzero.
\end{proof}
\begin{remark}
    There is no $\E{1} \otimes \mathbf{A}_2$-map $\Z^\ds \to \Pic(\Z)$ sending $1 \mapsto \Z[1]$, because the fiber sequence $\RP^\infty \to \Pic(\Z) \to \Z^\ds$ of \cref{eq: space PicZ} does not even split as a map of $\E{1} \otimes \mathbf{A}_2$-spaces. Indeed, the $\E{1} \otimes \mathbf{A}_2$-bar construction of the classifying map $\Z^\ds \to K(\Z/2, 2)$ of this fibration is already nontrivial (since it defines the nontrivial generator of $\H^4(\CP^2; \Z/2)$).
\end{remark}
\begin{remark}\label{rmk: Q does not save shearing}
    Lest \cref{lem: no e2 half shearing} seem like a problem specific to $\Pic(\Z)$, we note that the same problem persists for $\Pic(\QQ)$ (as well as for $\Pic(\FF_p)$ with $p>2$, but not for $\Pic(\FF_2) \cong \Z$). Indeed, the Postnikov fiber sequence for $\picsp(\QQ)$ is given by
    $$\Sigma \gl_1(\QQ) = \QQ^\times[1] \to \picsp(\QQ) \to \pi_0 \picsp(\QQ) = \Z^\ds.$$
    Recall that $\QQ^\times \cong \Z/2 \oplus \bigoplus_{\text{primes}} \Z$. Under this identification, the boundary map $\Z^\ds \to \QQ^\times[2]$ composes to give the map
    $$\Z \to \Z/2[2] \oplus \bigoplus_{\text{primes}} \Z[2] \xar{\pr} \Z/2[2],$$
    which identifies with \cref{eq: sq2}. Since this composite is not null as an $\E{2}$-map upon applying $\Omega^\infty$ by the argument of \cref{lem: no e2 half shearing}, the map $\Pic(\QQ) \to \Z^\ds$ does not admit an $\E{2}$-splitting.
\end{remark}
\begin{remark}\label{rmk: half shear on evenly graded}
    Consider the fully faithful functor $2: \Sp^\gr \hookrightarrow \Sp^\gr$ which doubles the weight. Then the composite
    $$\Sp^\gr \xar{2} \Sp^\gr \xar{\tilde{\sh}^{1/2}} \Sp^\gr$$
    identifies with the usual shearing functor, and hence admits an $\E{2}$-monoidal structure. Similarly, if we replace $\Sp^\gr$ by $\Mod_\MU^\gr$ (and in particular $\Mod_\Z^\gr$), the analogue of the above composite admits a symmetric monoidal structure. In particular, if $M_\bull$ is an $\E{n}$-algebra in graded spectra (resp. graded $\MU$- or $\Z$-module spectra) which is concentrated in even weights, its ``half-shear'' $\tilde{\sh}^{1/2}(M_\bull)$ admits an $\E{\min(n,2)}$-algebra structure in graded spectra (resp. $\E{n}$-algebra structure in $\MU$- or $\Z$-module spectra).
    
    Moreover, $\tilde{\sh}^{1/2}: \Mod_{\FF_2}^\gr \to \Mod_{\FF_2}^\gr$ admits a \textit{symmetric} monoidal structure. Indeed, there is an $\Eoo$-map $\Z^\ds \to \Pic(\FF_2)$ sending $1\mapsto \FF_2[1]$. There are many ways to see this; for instance, one can argue as in \cref{lem: Eoo shear}. Namely let $J_\RR: \BO \times \Z^\ds \to \Pic(\Sp)$ denote the real $J$-homomorphism, so that $J_\RR$ is an $\Eoo$-map, and $1 \in \Z^\ds \mapsto S^1 \in \Pic(\Sp)$.
    Since $\MO$ is the Thom spectrum of the $\Eoo$-composite $\BO \to \BO\times \Z^\ds \xar{J_\RR} \Pic(\Sp)$, it can be understood as the initial $\Eoo$-ring $R$ equipped with a nullhomotopy of the $\Eoo$-map $\BO \xar{J_\RR} \Pic(\Sp) \to \Pic(R)$. In particular, there is a commutative diagram of $\Eoo$-maps:
    $$\xymatrix{
    \BO\times \Z = \Omega^\infty \ko \ar[r]^-{J_\RR} \ar[d] & \Pic(\Sp) \ar[d] \\
    \Z^\ds \simeq \Omega^\infty \tau_{\leq 0} \ko \ar[r] & \Pic(\MO).
    }$$
    There is an $\Eoo$-orientation $\MO \to \FF_2$, so we obtain an $\Eoo$-map $\Z^\ds \to \Pic(\MO) \to \Pic(\FF_2)$ sending $1\mapsto \FF_2[1]$, as desired.
\end{remark}
The following observation will be useful below; in particular, the final sentence says that having polynomial homotopy generated by even-degree classes automatically forces ``formality'', even as an $\E{2}$-algebra. For instance, for any compact torus $T$, the $\Eoo$-$\Z$-algebra $C^\ast(BT; \Z)$ is ``formal'' as an $\E{2}$-$\Z$-algebra, i.e., $C^\ast(BT; \Z) \simeq \sh^{1/2} \pi_\bull C^\ast(BT; \Z)$ (More generally, the same is true of $C^\ast(BG; R)$ as long as $\H^\ast(BG; R)$ is a polynomial algebra on even classes.)
\begin{lemma}\label{lem: formality polynomial}
    If $R$ is an $\Eoo$-$\FF_2$-algebra, and $A$ is an $\E{2}$-$R$-algebra such that $\pi_\bull A$ is a finitely generated polynomial $R$-algebra (whose generators need not live in even degrees), there is an equivalence $A \simeq {\sh}^{1/2} \pi_\bull A$ of $\E{1}$-$R$-algebras.

    If $R$ is an $\Eoo$-ring with even homotopy, and $A$ is an $\E{2}$-$R$-algebra such that $\pi_\bull A$ is a finitely generated polynomial $R$-algebra generated by classes in even degrees, there is an equivalence $A \simeq {\sh}^{1/2} \pi_\bull A$ of $\E{1}$-$R$-algebras. 
    If $A$ furthermore admits an $\E{3}$-$R$-algebra structure, the equivalence $A \simeq {\sh}^{1/2} \pi_\bull A$ can be upgraded to one of $\E{2}$-$R$-algebras.
\end{lemma}
\begin{proof}
    Write $\pi_\bull A \cong R[x_1, \cdots, x_n]$ with $x_j$ in weight $2i_j$. Observe that the shearing ${\sh}^{1/2} R[x_j]$ is the free $\E{1}$-$R$-algebra on a class in degree $2i_j$. Similarly, if $R$ is an $\FF_2$-algebra and $y_j$ lives in weight $i_j$, the shearing ${\sh}^{1/2} R[y_j]$ is the free $\E{1}$-$R$-algebra on a class in degree $i_j$. This implies that there are $\E{1}$-$R$-algebra maps ${\sh}^{1/2} R[x_j] \to A$ for all $j$ such that $x_j$ is sent to the eponymous class on homotopy groups. Together, these define an $\E{1}$-$R$-algebra map $\bigotimes_{j=1}^n {\sh}^{1/2} R[x_j] \to A^{\otimes n} \to A$; the first map is evidently one of $\E{1}$-$R$-algebras, and the second map is one of $\E{1}$-$R$-algebras since $A$ is an $\E{2}$-$R$-algebra. By construction, this map induces an isomorphism on homotopy (and hence is an equivalence).

    Now assume that $A$ is an $\E{3}$-$R$-algebra such that $\pi_\bull A$ is a finitely generated polynomial $R$-algebra generated by classes in even degrees; we will first assume that $\pi_\bull A \cong \pi_\bull(R)[x_{2n}]$ with $x_{2n}$ in weight $2n$. 
    Let $R[x_{0,1}] = R[\Z_{\geq 0}]$ denote the flat \textit{graded} polynomial $R$-algebra on a class in weight $1$ and degree zero, so that $\tilde{\sh}^{\circ n}(R) = R[x_{2n,1}]$ is a graded polynomial $R$-algebra on a class in weight $1$ and degree $2n$. By \cite[Corollary 3.12]{framed-e2}, $R[x_{2n,1}]$ admits the structure of a framed $\E{2}$-algebra in $\Mod_R^\gr$. Let $\und(R[x_{2n,1}])$ denote the underlying ungraded $\E{2}$-$R$-algebra (in the body of this article, we will often simply omit ``$\und$'', but we keep it here for clarity), so that $\und(R[x_{2n,1}])$ is the free $\E{1}$-$R$-algebra on a class in degree $2n$. We will construct an equivalence $\und(R[x_{2n,1}]) \xar{\sim} A$ of $\E{2}$-$R$-algebras. For the case of multiple generators, note that if $\pi_\bull(A) = \pi_\bull(R)[x_1, \cdots, x_j]$, the above discussion produces an $\E{2}$-$R$-algebra map $\bigotimes_{j=1}^n {\sh}^{1/2} R[x_j] \to A^{\otimes n} \to A$; the first map is one of $\E{2}$-$R$-algebras by assumption, and the second map is one of $\E{2}$-$R$-algebras since $A$ is an $\E{3}$-$R$-algebra. The map $\bigotimes_{j=1}^n {\sh}^{1/2} R[x_j] \to A$ is an isomorphism on homotopy, hence an equivalence, as desired.

    Let us now return to the task of constructing an equivalence $\und(R[x_{2n,1}]) \xar{\sim} A$ of $\E{2}$-$R$-algebras. It suffices to show that $\und(R[x_{2n,1}])$ admits an \textit{even} cell structure as an $\E{2}$-$R$-algebra (this does \textit{not} require that $R$ have even homotopy!). Indeed, the class $x_{2n} \in \pi_{2n}(A)$ defines a map from the bottom $\E{2}$-cell into $A$; all obstructions to extending this map along the higher $\E{2}$-cells of $\und(R[x_{2n,1}])$ live in odd degrees, but the odd homotopy of $A$ vanishes, so such an extension $\und(R[x_{2n,1}]) \to A$ exists. By construction, this map is an isomorphism on homotopy, and hence is an equivalence. To construct an $\E{2}$-cell structure for $\und(R[x_{2n,1}])$, note that it in fact suffices to construct an $\E{2}$-cell structure for $R[x_{0,1}]$ in $\Alg_\E{2}(\Mod_R^\gr)$: indeed, the desired $\E{2}$-cell structure on $\und(R[x_{2n,1}])$ then follows from shearing and the fact that $\tilde{\sh}$ is $\E{2}$-monoidal by \cref{cstr: shearing}. Since $R[x_{0,1}]$ is an augmented $R$-algebra whose augmentation ideal is concentrated in positive weights, an $\E{2}$-cell structure for $R[x_{0,1}]$ is specified by a cell structure for the $2$-fold bar construction $\Bar^{(2)}(R[x_{0,1}])$. But this is a standard calculation: one finds that $\Bar^{(2)}(R[x_{0,1}]) \simeq \bigoplus_{n \geq 0} R[2n](n)$, at least as $R$-modules. For a reference in slightly different language, see \cite[Proposition 3.4.5]{rotinv}.
\end{proof}
\begin{remark}\label{rmk: shearing not e3}
    In \cite[Remark 3.11]{framed-e2}, we show that $\tilde{\sh}: \Sp^\gr \to \Sp^\gr$ \textit{cannot} be made into an $\E{3}$-monoidal functor. This implies that in the setting of \cref{lem: formality polynomial}, even if $A$ admits an $\E{4}$-$R$-algebra structure, the equivalence $A \simeq {\sh}^{1/2} \pi_\bull A$ need not upgrade to an equivalence of $\E{3}$-$R$-algebras. This is closely related to the subtlety of refining the derived geometric Satake equivalence (discussed in the present article as \cref{thm: derived satake}) into an $\E{3}$-monoidal equivalence. 
\end{remark}
\begin{remark}\label{rmk: modeling maps via shearing}
    Being concerned with the \textit{relative} Langlands program, it is important to understand whether \cref{lem: formality polynomial} can be used to model \textit{maps} between polynomial ring spectra, too. It is not too hard to show, for instance, that the map $C^\ast(BG; R) \to C^\ast(BH; R)$ induced by a homomorphism $H \to G$ can be modeled as the shearing of the induced map on cohomology when $R = \QQ$ (both $\H^\ast(BG; R)$ and $\H^\ast(BH; R)$ are polynomial on even degree classes), and in the cases of $\GL_n \to \GL_{n+1}$ and $R = \Z$, $\SO_n \to \SO_{n+1}$ when $R = \FF_2$ or $R = \Z[1/2]$, $\O_n \to \GL_n$ when $R = \FF_2$, etc.
    
    However, there are examples where the map $C^\ast(BG; R) \to C^\ast(BH; R)$ cannot be modeled so na\"ively; this should reflect subtle phenomena under Langlands duality at bad primes. For instance, consider the degree two map $\GG_m \to \GG_m^{(2)}$ with $R = \FF_2$. The homotopy fiber of the map $BS^1 \to B(S^1/\mu_2)$ is $B\mu_2 \simeq \RP^\infty$, and therefore 
    $$C^\ast(BS^1; \FF_2) \otimes_{C^\ast(B(S^1/\mu_2); \FF_2)} \FF_2 \simeq C^\ast(\RP^\infty; \FF_2).$$
    The induced map $\H^\ast(B(S^1/\mu_2); \FF_2) \to \H^\ast(BS^1; \FF_2)$ on cohomology rings is given by $\FF_2[x] \to \FF_2[y]$, which necessarily sends $x$ to zero (indeed, the map on cohomology is the derivative of the map $\GG_m \to \GG_m^{(2)}$ at the identity, and the derivative of a square vanishes modulo $2$).
    This implies that there is an isomorphism
    $$\H^\ast(BS^1; \FF_2) \otimes_{\H^\ast(B(S^1/\mu_2); \FF_2)} \FF_2 \cong \FF_2[y, \sigma(x)]/(\sigma(x)^2).$$
    The class $\sigma(x)$ lives in homological degree $-1$. The Serre spectral sequence starts with the above tensor product and converges to the $\FF_2$-cohomology of $\RP^\infty$; this spectral sequence degenerates, but there is an interesting multiplicative extension $\sigma(x)^2 = y$ on the $E_\infty$-page (and this is proved using topological input invisible to the map $\H^\ast(B(S^1/\mu_2); \FF_2) \to \H^\ast(BS^1; \FF_2)$). One finds in conclusion that $\H^\ast(B\mu_2; \FF_2) \cong \FF_2[\sigma(x)]$, as expected.
\end{remark}

%% file: connective-equiv-k-theory/equivariant-k-theory.tex
\subsection{Equivariant K-theory}

Let $G$ be a compact Lie group. Atiyah and Segal constructed \textit{$G$-equivariant K-theory} using the theory of $G$-equivariant vector bundles. We will review this theory here and describe the spectral algebro-geometric perspective on equivariant K-theory following \cite{survey}.
\begin{definition}
    A \textit{finite $G$-space} $X$ is a space with $G$-action which is constructed from finitely many $G$-cells of the form $G/H \times D^n$, where $H\subseteq G$ is a closed subgroup. Let $\cS(G)$ denote the $\infty$-category of finite $G$-spaces and $G$-equivariant maps between them.
\end{definition}
\begin{definition}
    Let $X$ be a finite $G$-space. A \textit{$G$-equivariant vector bundle} on $X$ is a vector bundle $\cV$ over $X$ equipped with a continuous $G$-action, such that the map $\cV \to X$ is $G$-equivariant. Let $\KU_G^0(X)$ denote the Grothendieck group of the monoid of $G$-equivariant vector bundles on $X$.
\end{definition}
Atiyah and Segal showed that the assignment $X \mapsto \KU_G^0(X)$ from the (opposite of the) homotopy category of finite $G$-spaces to groups extends to a cohomology theory which is represented in the homotopy category of $G$-spectra by a spectrum denoted by $\KU_G$.

In order to see all the structure on equivariant K-theory, it will be convenient to phrase the construction in terms of the $\infty$-category of orbispaces.
\begin{definition}
    Let $\Orb$ denote the \textit{global orbit $\infty$-category} as defined in \cite[Definition 2.7]{gepner-meier-equiv-tmf}. Heuristically, this is the full subcategory of the $\infty$-category of topological stacks spanned by objects of the form $\ast/G$. An \textit{orbispace} is a functor $\Orb^\op \to \Top$. Let $\Top_\Orb$ denote the $\infty$-category of orbispaces.

    Let $\Top_G$ denote the $\infty$-category of $G$-spaces, and let $\Orb_G$ denote the full subcategory of $\Top_G$ spanned by $G$-spaces of the form $G/H$ with $H\subseteq G$ being a closed subgroup.
    By \cite[Proposition 2.16]{gepner-meier-equiv-tmf}, there is a fully faithful functor $\Orb_G \to \Orb_{/\ast/G}$, whose essential image is spanned by those maps $\ast/H \to \ast/G$ which arise via an inclusion of subgroups $H\subseteq G$.
\end{definition}
\begin{remark}
    Note that $\cS(G)$ is the full subcategory of $\Top_G$ generated by $G$-spaces of the form $G/H$ (for closed subgroups $H\subseteq G$) under finite colimits.
\end{remark}
A more invariant construction of $\KU_G$, along with its $\Eoo$-ring structure, is as follows; see \cite[Section 4]{gepner-meier-equiv-tmf}.
\begin{construction}\label{cstr: KUG}
    Let $\Orb^\op \to \CAlg(\Prst)$ denote the functor sending $\ast/G \mapsto \Rep_\cc(G)$. Taking connective additive K-theory, we obtain a functor $K: \Orb^\op \to \CAlg$. The functor $K$ is a module over the constant functor $\Orb^\op \to \CAlg$ sending $\ast/G \mapsto K(\Vect_\cc) \simeq \tau_{\geq 0} (\KU)$. Therefore, inverting the Bott class $\beta \in \pi_2 \KU$ produces a functor $\Orb^\op \to \CAlg_\KU$ sending $\ast/G\mapsto \KU_G$.
    Right Kan extending along the functor $\Orb^\op \to \Top_\Orb^\op$ defines a lax symmetric monoidal functor $\Top_\Orb^\op \to \CAlg_\KU$ sending an orbispace $X/G \mapsto \KU_G(X)$.
\end{construction}
One important property of equivariant K-theory, which is also satisfied/posited to hold (depending on the construction) for equivariant analogues of other complex-oriented cohomology theories, is that it satisfies \textit{abelian descent}. Let us review this, following \cite[Section 4]{gepner-meier-equiv-tmf}.
\begin{definition}\label{def: orb for family of compact lie}
    Let $\cA$ denote a family of compact Lie groups (so that $\cA$ is closed under isomorphisms, subgroups, and quotients).
    Define $\Orb^\cA$ to be the full subcategory of $\Orb$ spanned by those $\ast/G$ with $G\in \cA$.
    For $\ast/G\in \Orb$, let $\Orb^\cA_G$ denote the full subcategory of $\Orb^\cA_{/\ast/G}$ spanned by those morphisms $\ast/H \to \ast/G$ which arise via an inclusion of subgroups $H\subseteq G$. Note that by \cite[Proposition 2.16]{gepner-meier-equiv-tmf}, one can identify $\Orb^\cA_G$ with the full subcategory of $\Orb_G$ spanned by those $G/H$ with $H\in \cA$.
\end{definition}
\begin{theorem}\label{thm: eq kthy ab extend}
    Let $\cA$ denote the family of abelian compact Lie groups.
    The functor $\Orb^\op \to \CAlg_\KU$ sending $\ast/G \mapsto \KU_G$ is right-Kan extended along the inclusion $\Orb^{\cA, \op} \hookrightarrow \Orb^\op$.
\end{theorem}
\begin{proof}
    For each $\ast/G$, one first notes that the inclusion $\Orb^\cA_G \to \Orb^\cA_{/\ast/G}$ is final, so we need to show that the canonical map $\KU_G \to \lim_{\ast/H\in \Orb^\cA_G} \KU_H$ is an equivalence. Let $E\cA \simeq \colim_{G/H\in \Orb^\cA_G} G/H$, so that $E\cA^K$ is empty if $K\not\in \cA$, and $E\cA^K \simeq \ast$ if $K\in \cA$. In fact, this property characterizes $E\cA$ up to weak equivalence. Then $\lim_{\ast/H\in \Orb^\cA_G} \KU_H \simeq \KU_G(E\cA)$, so we only need to show that the canonical map $\KU_G \to \KU_G(E\cA)$ is an equivalence. But this is \cite[Corollary 1.3]{adams-atiyah-segal}.
\end{proof}
\begin{remark}\label{rmk: KU ab ext}
    Instead of appealing to \cite[Corollary 1.3]{adams-atiyah-segal} in \cref{thm: eq kthy ab extend}, one can argue explicitly as follows in the case when $G$ is connected with torsion-free $\pi_1(G)$.
    Let $T$ be a maximal torus of $G$, and let $G/T$ be the flag variety. If $H\subseteq G$ is a closed subgroup, $H$ is abelian if and only if some conjugate $gHg^{-1}$ is contained in $T$, which in turn happens if and only if $(G/T)^H$ is nonempty. This implies that the $H$-invariants of the geometric realization $|(G/T)^{\times \bull+1}|$ is nonempty if and only if $H$ is abelian, in which case it is contractible. Therefore, by uniqueness of $E\cA$, there is a weak equivalence $E\cA \simeq |(G/T)^{\times \bull+1}|$. This implies that 
    $$\KU_G(E\cA) \simeq \Tot \KU_G((G/T)^{\times \bull+1}) \simeq \Tot \KU_G(G/T)^{\otimes_{\KU_G} \bull+1}.$$
    To conclude that this totalization is equivalent to $\KU_G$ by the unit map, it therefore suffices to show that the map $\KU_G \to \KU_G(G/T)$ induces a faithfully flat map on homotopy. But $\KU_G(G/T) = \KU_T$, so by $2$-periodicity, we only need to show that the map $R_\cc(G) \to R_\cc(T)$ on complex representation rings is faithfully flat. In fact, $R_\cc(T)$ is a free $R_\cc(G)$-module by the main theorem of \cite{pittie}, thanks to our assumption that $G$ is connected with torsion-free $\pi_1$.
\end{remark}

\begin{observe}
    The $\Eoo$-$\KU$-algebra $\KU_T$ is $2$-periodic, with $\pi_0$ given by the complex representation ring $R_\cc(T)$. In particular, $\spec \pi_0 \KU_T \cong \spec \Z[\bX^\ast(T)]$, where $\bX^\ast(T)$ is the lattice of characters. This is precisely the algebraic group $\Hom(\bX^\ast(T), \GG_m)$.
\end{observe}
\begin{prop}
    Let $\GG_{m,\KU}$ denote the $\Eoo$-$\KU$-scheme given by $\spec \KU[\Z]$. Let $T$ be an abelian compact Lie group. Then there is an equivalence $\spec \KU_T \simeq \Hom(\bX^\ast(T), \GG_{m,\KU})$.
\end{prop}
\begin{variant}\label{var: SAG KUG}
    Let $\cA$ denote the family of abelian compact Lie groups.
    Let $\Orb^\cA \to \Sch_{/\KU}$ denote the functor to spectral schemes over $\KU$ sending $\ast/T \mapsto \spec \KU_T$. That this is well-defined is essentially \cite[Proposition 4.4]{gepner-meier-equiv-tmf}.
    The left Kan extension of this functor along the inclusion $\Orb^\cA \hookrightarrow \Orb$ defines a functor $\Orb \to \Sch_{/\KU}$, which, by \cref{thm: eq kthy ab extend}, sends $\ast/G \mapsto \spec \KU_G$. Further left Kan extending along the inclusion $\Orb \to \Top_\Orb$ defines a functor $\Top_\Orb \to \Sch_{/\KU}$ sending $X/G \mapsto \spec \KU_G(X)$.
\end{variant}
\begin{remark}
    This construction can be extended further. Namely, consider the composite $\Orb^{\cA,\op} \to \Sch_{/\KU}^\op \to \CAlg(\LinCat_\KU)$, where the functor $\Sch_{/\KU}^\op \to \CAlg(\LinCat_\KU)$ is given by taking quasicoherent sheaves. Right Kan extending this functor along the inclusion $\Orb^{\cA, \op} \to \Top_\Orb^\op$ defines a functor $\Top_\Orb^\op \to \CAlg(\LinCat_\KU)$, which we will denote by $X/G \mapsto \Loc_G(X; \KU)$. The $\infty$-category $\Loc_G(X; \KU)$ could (somewhat abusively) be called the $\infty$-category of \textit{$G$-equivariant local systems of $\KU$-modules on $X$}. We will not use this notion below.
\end{remark}
\begin{remark} 
    The functor $\Top_\Orb \to \Sch_{/\KU}$ refines to a functor $(\Top_\Orb)_{/\ast/G} \to \Sch_{/\KU_G}$. In particular, if $X$ is a space with $G$-action, and $X/G$ denotes the associated orbispace (so that $X/G\in (\Top_\Orb)_{/\ast/G}$), there is a canonical map $\spec \KU_G(X) \to \spec \KU_G$. We will often write $\cf_G(X)$ to denote the associated $\Eoo$-$\KU_G$-algebra.
\end{remark}

%% file: connective-equiv-k-theory/equivariant-little-ku.tex
\subsection{Equivariant connective K-theory}\label{subsec: equiv little ku}

We will need a good theory of equivariant \textit{connective} K-theory. (This is \textit{not} the functor $K$ of \cref{cstr: KUG}.) It will be most convenient to adopt the spectral algebro-geometric perspective of \cref{var: SAG KUG}. To motivate the construction, let us briefly recall the definition of nonequivariant connective K-theory.
\begin{definition}
    Let $\ku$ denote the connective cover of complex K-theory $\KU$. Then $\pi_\ast \ku \cong \Z[\beta]$ with $|\beta| = 2$, so that $\ku[\beta^{-1}] = \KU$, and $\ku/\beta \simeq \Z$.
\end{definition}
Let us suggest some desiderata in the simple case of $S^1$-equivariance.
\begin{expect}\label{expect: kus1}
    By construction, there is an isomorphism $\spec \KU_{S^1} \simeq \GG_{m,\KU}$ of spectral schemes over $\KU$.
    Recall that $\spec \H^\ast(BS^1; \Z) \cong \hat{\GG}_a(2)$ as graded $\Z$-schemes, where the coordinate of $\hat{\GG}_a$ lives in weight $-2$. This can also be identified with the graded $\Z$-scheme $\spec \H^\ast_{S^1}(\ast; \Z) = \GG_a(2)$, since equipping the coordinate on $\GG_a$ with the nonzero weight $-2$ allows us to identify $\hat{\GG}_a(2) \cong \GG_a(2)$. Therefore, if $\Z_{S^1}$ denotes the $\Eoo$-$\Z$-algebra representing $S^1$-equivariant $\Z$-cohomology, one expects the appropriate notion of $S^1$-equivariant connective K-theory $\ku_{S^1}$ to be a sufficiently structured $\ku$-algebra such that there is a diagram where each square is Cartesian:
    $$\xymatrix{
    \GG_a(2) \ar@{=}[d] \ar[r] & \GG \ar@{=}[d] & (\GG_m)_{\Z[\beta^{\pm 1}]} \ar@{=}[d] \ar[l] \\
    \spec \pi_\ast \Z_{S^1} \ar[r] \ar[d] & \spec \pi_\ast \ku_{S^1} \ar[d] & \spec \pi_\ast \KU_{S^1} \ar[d] \ar[l] \\
    \spec \Z \ar[r] & \spec \Z[\beta] & \spec \Z[\beta^{\pm 1}] \ar[l].
    }$$
    In particular, one expects that there is an isomorphism of graded $\Z[\beta]$-group schemes
    $$\spec \pi_\ast \ku_{S^1} \cong \spec \Z[\beta][x, \tfrac{1}{1+\beta x}],$$
    where $x$ lives in weight $-2$ and the group structure is given by $x\mapsto x \otimes 1 + 1\otimes x + \beta x \otimes x$.
\end{expect}
Let us recall a construction of the group scheme $\spec \Z[\beta][x, \tfrac{1}{1+\beta x}]$.
\begin{recall}\label{recall: Rees}
    Let $\Z[y]$ denote the graded $\Z$-algebra where $y$ has weight $1$.
    The $(t-1)$-adically filtered ring $\Z[t^{\pm 1}]$ defines a commutative algebra object $\Mod_\Z^\fil$, which, by the Rees construction $\Mod_{\Z[y]}^\gr \simeq \Mod_\Z^\fil$, defines a commutative algebra object of $\Mod_{\Z[y]}^\gr$. This algebra is simply $\Z[y][t^{\pm 1}, \tfrac{t-1}{y}]$.
    The shearing autoequivalence of $\Mod_\Z^\gr$ sending $M_\bull\mapsto M_\bull[2\bull]$ sends $\Z[y]$ to a graded ring $\Z[\beta']$ with $\beta'$ in weight $1$ and degree $2$, so that shearing defines an equivalence 
    $$\Mod_\Z^\fil \simeq \Mod_{\Z[y]}^\gr \simeq \Mod_{\Z[\beta']}^\gr.$$
    Under this equivalence, the $(t-1)$-adically filtered ring $\Z[t^{\pm 1}]$ is sent to the graded $\Z[\beta']$-algebra $\Z[\beta', t^{\pm 1}, \tfrac{t-1}{\beta'}]$.
\end{recall}
We will adapt essentially the same argument to the spectral setting, with some minor homotopical issues complicating the story.
\begin{prop}[\cite{filtered-A1-Gm}]\label{prop: filt A1}
    The $\Eoo$-$\MU$-algebra $\MU[y]$ (with $y$ in degree zero) admits a grading where $y$ is in weight $1$.
    There is a symmetric monoidal equivalence $\Mod_{\MU[y]}^\gr \simeq \Mod_\MU^\fil$. In particular, if $\MU[\beta']$ denotes the graded $\Eoo$-$\MU$-algebra given by the shearing of $\MU[y]$ (so that $\beta'$ lives in degree $2$ and weight $1$), \cref{lem: Eoo shear} defines a symmetric monoidal equivalence $\Rees_{\beta'}: \Mod_\MU^\fil \xar{\sim} \Mod_{\MU[\beta']}^\gr$.
\end{prop}
\begin{lemma}\label{lem: oblv MUbeta'}
    Let $\oblv(\MU[\beta'])$ denote the ungraded $\Eoo$-$\MU$-algebra which underlies $\MU[\beta']$.
    There is an $\Eoo$-$\MU$-algebra map $\oblv(\MU[\beta']) \to \ku$ sending $\beta' \mapsto \beta$.
\end{lemma}
\begin{proof}
    Note that by construction, $\oblv(\MU[\beta'])$ can be identified with the Thom spectrum of the composite
    $$\BU \times \Z_{\geq 0} \to \BU \times \Z \xar{J} \Pic(\Sp).$$
    Let $\MU_\per$ denote the Thom spectrum of the map $J: \BU \times \Z \to \Pic(\Sp)$. Then there is a canonical map $\oblv(\MU[\beta']) \to \MU_\per$.
    There is an $\Eoo$-map $\MU_\per \to \KU$ which sends $\oblv(\beta')\in \pi_2 \MU_\per$ to the Bott class (see, e.g., \cite[Remark 6.3]{hoyois2021hilbert}), and hence an $\Eoo$-map $\oblv(\MU[\beta']) \to \KU$. which does the same. But $\oblv(\MU[\beta'])$ is connective, so this map factors through an $\Eoo$-map $\oblv(\MU[\beta']) \to \ku$, as desired.
\end{proof}
\begin{construction}\label{cstr: kuS1 and Gbeta}
    Let $I$ denote the fiber of the $\Eoo$-map $\MU[t^{\pm 1}] \to \MU$ given by $\MU$-chains of the map $\Z \to \ast$. Then the functor $\Z_{\geq 0} \to \Mod_\MU$ sending $n\mapsto I^{\otimes_{\MU[t^{\pm 1}]} n}$ defines an object $\F^\star_{(t-1)} \MU[t^{\pm 1}]\in \CAlg(\Mod_\MU^\fil)$.
    Applying the symmetric monoidal equivalence of \cref{prop: filt A1}, we obtain a graded $\Eoo$-$\MU[\beta']$-algebra $\Rees_{\beta'}(\F^\star_{(t-1)} \MU[t^{\pm 1}])$. Note that this graded $\Eoo$-$\MU[\beta']$-algebra admits a strictly cocommutative $\Eoo$-$\MU$-coalgebra structure by the map $t\mapsto t\otimes t$.

    Let $\ku_{S^1}$ denote the $\Eoo$-$\ku$-bialgebra given by
    $$\ku_{S^1} := \oblv(\Rees_{\beta'}(\F^\star_{(t-1)} \MU[t^{\pm 1}])) \otimes_{\oblv(\MU[\beta'])} \ku,$$
    where the $\Eoo$-map $\oblv(\MU[\beta']) \to \ku$ is given by \cref{lem: oblv MUbeta'}.
    We will write $\spec \ku_{S^1}$, or sometimes $\GG_{\ku,\beta}$, to denote the functor $\CAlg_\ku \to \Top$ which is corepresented by $\ku_{S^1}$. Since $\ku_{S^1}$ is a strictly cocommutative $\Eoo$-$\ku$-coalgebra, this functor in fact lands in the $\infty$-category $s\Ab$ of simplicial abelian groups. We will write $\GG_\beta$ to denote the underlying graded $\pi_\ast(\ku) = \Z[\beta]$-scheme.
\end{construction}
\begin{remark}
    By comparison to \cref{recall: Rees}, it is not difficult to see that
    $$\pi_\ast \Rees_{\beta'}(\F^\star_{(t-1)} \MU[t^{\pm 1}]) \cong \pi_\ast(\MU)[\beta', t^{\pm 1}, \tfrac{t-1}{\beta'}],$$
    where $\beta'$ is in weight $1$ and degree $2$, and $\tfrac{t-1}{\beta'}$ is in weight $-1$ and degree $-2$. Therefore,
    $$\pi_\ast \ku_{S^1} \cong \Z[\beta][t^{\pm 1}, \tfrac{t-1}{\beta}],$$
    as expressed in \cref{expect: kus1}.
\end{remark}
\begin{remark}
    Note that if we invert the Bott class in $\pi_2 \ku_{S^1}$, we obtain
    \begin{align*}
        \ku_{S^1}[\beta^{-1}] & = \oblv(\Rees_{\beta'}(\F^\star_{(t-1)} \MU[t^{\pm 1}])[\beta'^{-1}]) \otimes_{\oblv(\MU[\beta'^{\pm 1}])} \KU \\
        & \simeq \MU[\beta'^{\pm 1}][t^{\pm 1}] \otimes_{\MU[\beta'^{\pm 1}]} \KU \simeq \KU[t^{\pm 1}] = \KU_{S^1},
    \end{align*}
    as expressed in \cref{expect: kus1}.
\end{remark}
\begin{observe}\label{obs: kuT}
    It is not difficult to extend the above construction to arbitary abelian compact Lie groups $T$. Namely, given a functor $F: \CAlg_\ku \to s\Ab$ and an abelian compact Lie group $T$, one obtains a new functor $F_T: \CAlg_\ku \to s\Ab$ given by $\Hom(\bX^\ast(T), F)$. If $F$ is corepresentable, the same is true of $F_T$. Applied to the functor $\GG_{\ku,\beta}: \CAlg_\ku \to s\Ab$, we obtain an $\Eoo$-$\ku$-algebra $\ku_T$. This assignment evidently defines a functor from the $\infty$-category of abelian compact Lie groups to $\Eoo$-$\ku$-algebras. We will write the underlying graded group $\Z[\beta]$-scheme of $\spec \ku_T$ as $T_\beta$.
\end{observe}
In order to extend the functoriality of the assignment $T\mapsto \ku_T$, we need the following.
\begin{lemma}\label{lem: key completion}
    There is an $\Eoo$-$\ku$-algebra map $\ku_{S^1} \to \ku^{hS^1}$ which is given on homotopy by the map
    $$\Z[\beta][t^{\pm 1}, \tfrac{t-1}{\beta}] \to \Z[\beta]\pw{\hbar} \cong \Z\pw{t-1}[\beta]\pw{\hbar}/(\beta \hbar = t-1)$$
    sending $\tfrac{t-1}{\beta} \mapsto \hbar$.
\end{lemma}
\begin{proof}
    The orientation of the spectral formal multiplicative group $\hat{\GG}_m$ over $\KU$ defines an equivalence of $\Eoo$-rings $\KU^{hS^1} \xar{\sim} \KU\pw{t-1}$ (see \cite[Sections 4.4 and 6.5]{elliptic-ii}), where the map sends $\beta\hbar\mapsto t-1$ on $\pi_0$. This defines a canonical $\Eoo$-map $S\pw{t-1} \to \KU^{hS^1}$, which can be interpreted as a $BS^1$-family of $\Eoo$-map $S\pw{t-1} \to \KU$. Since $S\pw{t-1}$ is connective, this is the same as a $BS^1$-family of $\Eoo$-maps $S\pw{t-1} \to \ku$, i.e., an $\Eoo$-map $S\pw{t-1} \to \ku^{hS^1}$.

    This defines an $\Eoo$-map $\ku[t^{\pm 1}] \to \ku^{hS^1}$ via the $\ku$-linearization of the composite
    $$S[t^{\pm 1}] \to S\pw{t-1} \to \ku^{hS^1}.$$
    Note that the pushout $\ku^{hS^1} \otimes_{S[t^{\pm 1}]} S$ equips $\ku^{hS^1}/\beta\hbar$ with the structure of an $\Eoo$-$\ku^{hS^1}$-algebra.
    Let $J$ denote the fiber of the $\Eoo$-map $\ku^{hS^1} \to \ku^{hS^1}/\beta\hbar$, so that $\pi_\ast J$ is the ideal of $\pi_\ast \ku^{hS^1} = \Z[\beta]\pw{\hbar}$ generated by $\beta \hbar$. Let $\F^\star_{\beta \hbar} \ku$ denote the filtered $\Eoo$-ring $J^{\otimes_{\ku^{hS^1}} \ast}$.
    The above discussion produces a map $\F^\star_{(t-1)} \ku[t^{\pm 1}] \to \F^\star_{\beta\hbar} \ku$ of filtered $\Eoo$-$\ku$-algebras, and hence a map
    $$\Rees_{\beta'}(\F^\star_{(t-1)} \ku[t^{\pm 1}]) \to \Rees_{\beta'}(\F^\star_{\beta \hbar} \ku^{hS^1})$$
    of graded $\Eoo$-$\ku[\beta']$-algebras.  In particular, applying $\oblv$ and base-changing along the $\Eoo$-$\ku$-algebra map $\oblv(\ku[\beta']) \to \ku$ from \cref{lem: oblv MUbeta'} defines a map of ungraded $\Eoo$-$\ku$-algebras
    $$\oblv(\Rees_{\beta'}(\F^\star_{(t-1)} \ku[t^{\pm 1}])) \otimes_{\oblv(\ku[\beta'])} \ku \to \oblv(\Rees_{\beta'}(\F^\star_{\beta \hbar} \ku^{hS^1})) \otimes_{\oblv(\ku[\beta'])} \ku.$$
    But it is not difficult to see that the target is precisely $\ku^{hS^1}$.
\end{proof}
\begin{prop}\label{prop: Gbeta is preoriented}
    The group scheme $\GG_{\ku,\beta}$ is preoriented (compatibly with the orientation on $\GG_{m,\KU}$), and the construction from \cref{obs: kuT} extends to a functor $\Orb^{\cA} \to \Fun(\CAlg_\ku, \Top)$ sending $T\mapsto \spec \ku_T$.
\end{prop}
\begin{proof}
    Following \cite[Construction 3.13 and Proposition 4.4]{gepner-meier-equiv-tmf}, the desired functor can be defined as follows. First, note that the $\Eoo$-map of \cref{lem: key completion} defines a map 
    $$\Map_{\CAlg_\ku}(\ku^{hS^1}, \ku) \to \Map_{\CAlg_\ku}(\ku_{S^1}, \ku) = \GG_{\ku,\beta}(\ku).$$
    There is an obvious map
    $$\CP^\infty \to \Map_{\CAlg_\ku}(\ku^{hS^1}, \ku)$$
    of simplicial abelian groups. The resulting map $\CP^\infty \to \GG_{\ku,\beta}(\ku)$ defines a preorientation $\ast/S^1 \to \GG_{\ku,\beta}$, and hence a functor 
    $\Orb^{\cA} \times \CAlg_\ku \to \Top$
    sending
    $$(X, R) \mapsto \Map_{s\Ab_{\ast/S^1/}}(\Map(X, \ast/S^1), \GG_{\ku,\beta}(R)).$$
    This is adjoint to the desired functor $\Orb^{\cA} \to \Fun(\CAlg_\ku, \Top)$. It is not difficult to see that this functor sends $T \mapsto \ku_T$.
\end{proof}
Motivated by \cref{var: SAG KUG}, we are led to:
\begin{definition}\label{def: G-equiv ku}
    Let $\Top_\Orb \to \Fun(\CAlg_\ku, \Top)$ denote the functor given by left Kan extending the functor $\Orb^{\cA} \to \Fun(\CAlg_\ku, \Top)$ along the inclusion $\Orb^\cA \hookrightarrow \Top_\Orb$. It is not hard to see that this functor in fact lands in the full subcategory spanned by the representable functors, so we will denote this functor by $X/G \mapsto \spec \cf_G(X)$. The functor $\Top_\Orb \to \Fun(\CAlg_\ku, \Top)$ refines to a functor $(\Top_\Orb)_{/\ast/G} \to \Fun(\CAlg_\ku, \Top)$. In particular, if $X$ is a space with $G$-action, and $X/G$ denotes the associated orbispace (so that $X/G\in (\Top_\Orb)_{/\ast/G}$), there is a canonical map $\spec \cf_G(X) \to \spec \ku_G$. We will write $\ku^\ast_G(X)$ to denote $\pi_{-\ast} \cf_G(X)$.
\end{definition}
\begin{remark}
    This construction can be extended further. Namely, consider the composite $\Orb^{\cA,\op} \to \Fun(\CAlg_\ku, \Top)^\op \to \CAlg(\LinCat_\ku)$, where the functor $\Fun(\CAlg_\ku, \Top)^\op \to \CAlg(\LinCat_\ku)$ is given by taking quasicoherent sheaves. Right Kan extending this functor along the inclusion $\Orb^{\cA, \op} \to \Top_\Orb^\op$ defines a functor $\Top_\Orb^\op \to \CAlg(\LinCat_\ku)$, which we will denote by $X/G \mapsto \Loc_G(X; \ku)$. The $\infty$-category $\Loc_G(X; \ku)$ could (somewhat abusively) be called the $\infty$-category of \textit{$G$-equivariant local systems of $\ku$-modules on $X$}. We will not use this notion below.
\end{remark}
\begin{notation}\label{notn: MG defn}
    Let $\cM_G$ denote the underlying graded $\Z[\beta]$-scheme of $\spec \ku_G$, i.e., $\cM_G = \spec \pi_\ast \ku_G$.
\end{notation}
\begin{prop}\label{prop: Weyl invts}
    Let $G$ be a connected compact Lie group such that $\pi_1(G)$ is torsion-free, and let $T\subseteq G$ be a maximal torus with associated Weyl group $W$. Let $X/G$ be an orbispace over $\ast/G$, and let $X/T$ denote the associated orbispace of $\ast/T$.
    Upon inverting $|W|$, the natural map $\spec \cf_T(X) \to \spec \cf_G(X)$ exhibits the graded scheme $\spec \ku^\ast_G(X)$ as the GIT quotient $\spec \ku^\ast_T(X)\mmod W$. In particular, $\cM_G \cong T_\beta \mmod W$.
\end{prop}
\begin{proof}
    Following \cref{rmk: KU ab ext}, we can identify $\cf_G(X)$ with the totalization of the diagram
    $$\cf_T(X) \rightrightarrows \cf_T(X \times G/T) \threerightarrows \cf_T(X \times G/T \times G/T) \cdots$$
    There is an isomorphism $\cf_T(X \times G/T) \simeq \cf_T(X) \otimes_{\ku_T} \cf_T(G/T)$. We will argue the case $X = \ast$ below; the general case is no more difficult.
    
    To understand $\cf_T(G/T)$, we will use \cref{prop: gkm}. This is classical: $(G/T)^T = W$, so the set $V$ of vertices of the GKM diagram is given by $V = W$. The edges in the GKM diagram are indexed by positive roots $\alpha\in \Phi^+$: the edge labelled by $\alpha$ connects the vertices $w \mapsto s_\alpha w$. It follows that there is an equalizer diagram
    $$\pi_\ast \cf_T(G/T) \hookrightarrow \Map(W, \pi_\ast \ku_T) \rightrightarrows \prod_{(w,\alpha)} \pi_\ast \ku_{T_\alpha},$$
    where $T_\alpha$ is the kernel of the map $\alpha: T \to S^1$. Said differently, the map $\cf_T \otimes_\ku \cf_T \to \cf_T(G/T)$ defines a closed immersion $\spec \pi_\ast \cf_T(G/T) \hookrightarrow T_\beta \times_{\spec \Z[\beta]} T_\beta$ which exhibits $\spec \pi_\ast \cf_T(G/T)$ as the union of graphs of $W$ acting on $T_\beta$. Let us denote this union of graphs by $\cU_W$.
    
    Now, if $\pi_\ast \cf_T(G/T)$ is flat over $\pi_\ast \ku_T$, there is an isomorphism $\pi_\ast \cf_T((G/T)^{\times k+1}) \cong (\pi_\ast \cf_T(G/T))^{\otimes_{\pi_\ast \ku_T} k}$. The assumption that $\pi_1(G)$ is torsion-free implies (by \cite{pittie}) that $\pi_\ast \ku^{(G/T)_+}$ is flat (in fact free) over $\pi_\ast \ku$, from which the desired flatness of $\pi_\ast \cf_T(G/T)$ over $\pi_\ast \ku_T$ follows. We find that $\spec \pi_\ast \ku_G$ can be expressed as the geometric realization of the simplicial diagram 
    $$\spec \pi_\ast \cf_T((G/T)^{\times \bull+1}) \cong \spec (\pi_\ast \cf_T(G/T))^{\times_{\spec \pi_\ast \ku_T} \bull} \cong \cU_W^{\times_{T_\beta} \bull}.$$
    Let us now invert $|W|$.
    The map $T_\beta \to T_\beta \mmod W$ is faithfully flat (by \cite{pittie}). The argument of \cite[Proposition A.2]{gannon-coarse-quotient} shows that $T_\beta \times_{T_\beta\mmod W} T_\beta \cong \cU_W$: it is true after inverting $\beta$ (even if $|W|$ is not inverted) by \cite[Remark A.3]{gannon-coarse-quotient}, and when $\beta = 0$, it follows from \cite[Proposition A.2]{gannon-coarse-quotient}. (This result applies in the present situation because $|W|$ is inverted, and so the Chevalley-Shephard-Todd theorem continues to hold for the map $\fr{t} \to \fr{t}\mmod W$.) The geometric realization of the above simplicial diagram can therefore be identified with $T_\beta \mmod W$, as desired.
\end{proof}

%% file: connective-equiv-k-theory/k-homology.tex
\subsection{Equivariant $\ku$-homology}

The goal of this section is to set up the theory of equivariant $\ku$-homology. Fix a compact Lie group $G$ throughout.
\begin{definition}
    Let $\cf_G(-)^\vee: \Top(G) \to \Mod_{\ku_G}$ denote the functor given by sending $X/G \mapsto \cf_G(X)^\vee$, where $\cf_G(X)^\vee$ denotes the $\ku_G$-linear dual of $\cf_G(X)$. We will refer to $\cf_G(X)^\vee$ as the \textit{$G$-equivariant $\ku$-homology} of $X$, and often write the homotopy groups of this spectrum as $\ku^G_\ast(X)$.
\end{definition}
\begin{remark}
    Note that the functor $\cf_G(-)^\vee: \Top(G) \to \Mod_{\ku_G}$ is in fact symmetric monoidal. Since every $G$-space is naturally equipped with a diagonal map, this refines $\cf_G(-)^\vee$ to a functor $\Top(G) \to \coCAlg_{\ku_G}$.
\end{remark}
The above definition is badly behaved if $X$ is not a finite $G$-space. This is not special to the equivariant setting, as the following example shows.
\begin{example}
    The integral cohomology of the discrete space $\Z$ is given by the ring $\Map(\Z, \Z)$ of all functions $\Z \to \Z$. Given such a function $f$, one can define a new function $\Delta f: \Z \to \Z$ via $f(x+1) - f(x)$. Then, we formally have $f(x) = \sum_{k \geq 0} (\Delta^k f)(0) \binom{x}{k}$. This series converges in the completion of the ring $\Z[\binom{x}{k}]_{k\geq 0}$ of numerical polynomials. However, the $\Z$-linear dual of $\Map(\Z, \Z)$ is \textit{not} isomorphic to the group algebra $\Z[t^{\pm 1}] = \Z[\Z] = \H_\ast(\Z; \Z)$.

    The basic issue is the infinitude of $\Z$, which leads to a difference between $\Z^\Z$ and $\Z^{\oplus \Z}$.
    The simplest fix is to observe that $\Z$ admits a filtration by finite subsets $I_n = \{-n, \cdots, n\}$, and that $\H_\ast(I_n; \Z) = \Z\{t^{-n}, \cdots, t^n\}$ is indeed the $\Z$-linear dual of $\H^\ast(I_n; \Z) = \Map(I_n, \Z)$. Note that this filtration of $\Z$ equips it with the structure of a filtered group: namely, the addition on $\Z$ gives maps $I_n \times I_m \to I_{n+m}$ for each $n,m\in \Z$.
\end{example}
We will therefore rely on the following construction.
\begin{construction}\label{cst: ku-hmlgy filtered colimit}
    Let $\cf_G(-)^\vee: \Top_G \to \coCAlg_{\ku_G}$ denote the left Kan extension of the functor $\Top(G) \to \coCAlg_{\ku_G}$ along the inclusion $\Top(G) \to \Top_G$. Explicitly, if $X\in \Top_G$ is a $G$-space equipped with a presentation $X = \colim_{j\in \cJ} X_j$ as the filtered colimit of a filtered diagram $\cJ \to \Top(G)$ of finite $G$-spaces, then $\cf_G(X)^\vee$ is the filtered colimit $\colim_{j\in \cJ} \cf_G(X_j)^\vee$; we will refer to it as the \textit{$G$-equivariant $\ku$-homology of $X$}.
    Note that the forgetful functor $\coCAlg_{\ku_G} \to \Mod_{\ku_G}$ preserves colimits, so this filtered colimit can be computed in $\coCAlg_{\ku_G}$ or $\Mod_{\ku_G}$.
\end{construction}
In most examples of interest, there will be a geometrically defined presentation of $X$.
\begin{remark}
    Let $X$ be a \textit{finite} $G$-space equipped with an $\E{n}$-algebra structure in $G$-spaces. Then $\cf_G(X)^\vee$ admits an $\E{n}$-algebra structure in $\coCAlg_{\ku_G}$.
    If $X$ is not a finite $G$-space, but is equipped with an $\E{n}$-algebra structure in $\Top_G$, the definition of $G$-equivariant $\ku$-homology via \cref{cst: ku-hmlgy filtered colimit} does not guarantee the existence of an $\E{n}$-algebra structure on $\cf_G(X)^\vee \in \coCAlg_{\ku_G}$. Rather, if $\cJ$ is a filtered index category equipped with an $\E{n}$-monoidal structure and $\cJ \to \Top_G$ is an $\E{n}$-algebra object in the $\infty$-category $\Fun(\cJ, \Top_G)$ equipped with the Day convolution monoidal structure, \cref{cst: ku-hmlgy filtered colimit} will define an $\E{n}$-algebra structure on $\cf_G(X)^\vee \in \coCAlg_{\ku_G}$. We will refer to such a presentation of $X$ as a \textit{multiplicative} presentation.
\end{remark}
A basic fact about equivariant connective K-theory is the localization theorem. Although one can make statements about $\ku_G$ for arbitrary compact Lie groups $G$, we will restrict attention only to the case when $G = T$ is an abelian compact Lie group. In this case, we have the following simple observation.
\begin{lemma}\label{lem: atiyah localization}
    Let $T$ be an abelian compact Lie group, and let $X$ be a finite $T$-space. Let $T_0\subseteq T$ be a closed subgroup, and let $\cU_{T_0} \subseteq T_\beta$ denote the complement of the union of the closed subschemes $T'_\beta$ ranging over all closed subgroups $T' \subseteq T$ which do not contain $T_0$. Then the map $\cf_T(X) \to \cf_T(X^{T_0})$, and hence the map $\cf_T(X^{T_0})^\vee \to \cf_T(X)^\vee$, is an equivalence upon restriction to $\cU_{T_0}$.
\end{lemma}
\begin{proof}
    By induction on the orbit stratification on $X$, we are reduced to the case when $X = T/T_1$ for some closed subgroup $T_1 \subseteq T$. In this case, the fixed points $X^{T_0}$ is empty if $T_0 \not\subseteq T_1$, and $X^{T_0} = X$ if $T_0 \subseteq T_1$. It therefore suffices to show that $\cf_T(X) |_{T_\beta - T_{1,\beta}} = 0$ if $T_0 \not \subseteq T_1$; but this is clear, because $\cf_T(X) \cong \co_{T_{1,\beta}}$.
\end{proof}
\begin{remark}
    One special case of \cref{lem: atiyah localization} which is worth restating (coresponding to $T_0 = T$) is the following. Let $\punc{T}_\beta$ denote the complement of the union of the closed subscheme $T'_\beta$ ranging over all closed \textit{proper} subgroups $T' \subsetneq T$. Then the map $\cf_T(X) \to \cf_T(X^T)$, and hence the map $\cf_T(X^T)^\vee \to \cf_T(X)^\vee$, is an equivalence upon restriction to $\punc{T}_\beta$.
\end{remark}
\begin{lemma}\label{lem: injective coh}
    Let $T$ be a torus, and let $X$ be a finite $T$-space. If $\pi_\ast \cf_T(X)$ is a projective $\pi_\ast \ku_T$-module, the map $\pi_\ast \cf_T(X) \to \pi_\ast \cf_T(X^T)$ is an injection.
\end{lemma}
\begin{proof}
    Since the map $\cf_T(X) \to \cf_T(X^T) \to \cf_T(X^T)_{\punc{T}_\beta}$ factors as
    $\cf_T(X) \to \cf_T(X)|_{\punc{T}_\beta} \to \cf_T(X^T)|_{\punc{T}_\beta}$,
    and the map $\cf_T(X)|_{\punc{T}_\beta} \to \cf_T(X^T)|_{\punc{T}_\beta}$ is an equivalence by \cref{lem: atiyah localization}, it suffices to show that the map $\cf_T(X) \to \cf_T(X)|_{\punc{T}_\beta}$ induces an injection on homotopy groups. But $\pi_\ast \cf_T(X)$ was assumed to be a projective $\pi_\ast \ku_T$-module, so one is reduced to the case $X = \ast$, i.e., to showing that the map $\ku_T \to \ku_T|_{\punc{T}_\beta}$ induces an injection on homotopy groups. This, however, is clear, since the closed subscheme $T'_\beta \hookrightarrow T_\beta$ defined by each closed subgroup $T'\subseteq T$ is cut out by a regular sequence.
\end{proof}
\begin{definition}
    Let $X$ be a finite $T$-space equipped with a chosen presentation in terms of $T$-cells. Say that $X$ is a \textit{GKM space} if the following conditions are satisfied:
    \begin{itemize}
        \item $\pi_\ast \cf_T(X)$ is a projective $\pi_\ast \ku_T$-module;
        \item if $X^{(1)}$ denotes the equivariant $1$-skeleton of $X$, then $X^{(1)}$ consists of a finite number of spheres $S^\lambda$ meeting only at the fixed points, where $\lambda$ ranges over characters of $T$.
    \end{itemize}
    Let $V$ denote the set $X^T$ of fixed points, and let $E$ denote the set of characters $\lambda$ such that $S^\lambda \subseteq X^{(1)}$. There are two maps $E \rightrightarrows V$ sending $\lambda$ to the points $0,\infty\in S^\lambda \subseteq X^{(1)}$.
\end{definition}
\begin{prop}[Goresky-Kottwitz-MacPherson]\label{prop: gkm}
    Let $X$ be a finite GKM $T$-space equipped with a chosen presentation in terms of $T$-cells. For each character $\lambda: T \to S^1$, let $T_\lambda$ denote the kernel of $T$, and let $S(\lambda)$ denote the unit representation sphere, so that $\ku_{T_\lambda} \cong \cf_T(S(\lambda))$. Then there is an equalizer diagram
    $$\pi_\ast \cf_T(X) \hookrightarrow \pi_\ast \cf_T(X^T) \cong \Map(V, \pi_\ast \ku_T) \rightrightarrows \prod_{\lambda \in E} \pi_\ast \ku_{T_\lambda},$$ 
    where the two maps in the equalizer are defined in the evident manner.
\end{prop}
A general version of the Goresky-Kottwitz-MacPherson theorem is proved in \cite{generalized-gkm}.
\begin{proof}
    Let us first show that the maps $\cf_T(X) \to \cf_T(X^T)$ and $\cf_T(X^{(1)}) \to \cf_T(X^T)$ have the same images on homotopy. There is an evident map from the image of $\cf_T(X) \to \cf_T(X^T)$ on homotopy to the image of $\cf_T(X^{(1)}) \to \cf_T(X^T)$ on homotopy, which we will denote by $f$. The map $f$ is an injection by \cref{lem: injective coh}.
    Let $T'$ denote a proper closed subgroup of $T$ of codimension $1$, and let $U' \subseteq T'_\beta$ denote the complement of the union of the closed varieties $T''_\beta$ ranging over the proper closed subgroups $T''\subseteq T'$.
    By \cref{lem: atiyah localization}, the map $f$ is an isomorphism upon restriction to $U'\subseteq T'_\beta \subseteq T_\beta$ for each proper closed subgroup $T'\subseteq T$ of codimension $1$.
    
    Therefore, the locus $Z \subseteq T_\beta$ over which $f$ fails to be an isomorphism is contained in the union of closed subvarieties $T'_\beta$ for finitely many $T'\subseteq T$ of codimension at least $2$. However, the map $\cf_T(X) \to \cf_T(X)|_{T_\beta - Z}$ is an isomorphism (by Hartogs). Since the same is true of the map $\cf_T(X^T) \to \cf_T(X^T)|_{T_\beta - Z}$, and the map $\cf_T(X) \to \cf_T(X^T)$ factors through the map $\cf_T(X^{(1)}) \to \cf_T(X^T)$, the desired result follows.

    For the equalizer diagram, an easy induction on the cell structure of $X$ reduces us to the case $X = S^\lambda$ for a character $\lambda: T \to S^1$. In this case, the isomorphism $T/T_\lambda \cong S^\lambda$ defines an isomorphism $\ku_{T_\lambda} \cong \cf_T(S(\lambda))$. Since $S^\lambda \cong \Sigma S(\lambda)$, we obtain an equalizer diagram
    $$\pi_\ast \ku_T(S^\lambda) \to \pi_\ast \ku_T \oplus \pi_\ast \ku_T \cong \Map(\{0,\infty\}, \ku_T) \rightrightarrows \pi_\ast \ku_{T_\lambda}.$$
    This proves the desired claim.
\end{proof}
\begin{remark}\label{rmk: gkm minor ext}
    Note that the statement of \cref{prop: gkm} is natural in $X$, and in particular, one can use \cref{prop: gkm} to describe the $\pi_\ast \ku_T$-algebra structure on $\pi_\ast \cf_T(X)$. By dualizing \cref{prop: gkm}, one can also describe the $\pi_\ast \ku_T$-\textit{co}algebra structure on $\ku^T_\ast(X)$. Moreover, suppose that $X$ is a $T$-space equipped with a presentation $X = \colim_{j\in \cJ} X_j$ in terms of finite $T$-spaces, each of which is GKM and equipped with a chosen presentation in terms of $T$-cells, and such that the transition maps $X_j \to X_{j'}$ are maps of cellular $T$-spaces. Then \cref{prop: gkm} can be extended to compute $\pi_\ast \cf_T(X)$ and $\ku^T_\ast(X)$.
\end{remark}


%% file: equivalences/full-faithfulness.tex
\subsection{Full faithfulness of global sections}

In this section, we prove an analogue of a result of Ginzburg's from \cite{ginzburg-perverse}. We will closely follow \cite[Section 4.7]{quat-satake} and \cite[Section 8]{soergel-wendt}.
\begin{setup}
    Let $G$ be a compact Lie group, and fix a maximal torus $T\subseteq G$.
    Let $X$ be a finite $G$-space whose $G$-equivariant orbit stratification indexed by a poset $P$ (neccessarily finite). 
    Let $X_\lambda$ denote the stratum corresponding to $\lambda \in P$, and let $X_{\leq \lambda}$ denote its closure in $X$. 
    Suppose further that each $X_\lambda$ is a complex affine space of complex dimension $n_\lambda$ on which $G$ acts linearly.
    In particular, this implies that $\H^\ast_G(X; j_{\lambda, !} \ul{\QQ})$ is concentrated in even degrees for each $\lambda \in P$, where $j_\lambda: X_\lambda \hookrightarrow X_{\leq \lambda}$ denotes the inclusion.
    
    Let $X_{<\lambda} = X_{\leq \lambda} - X_\lambda$, and let $i_\lambda: X_{<\lambda} \hookrightarrow X_{\leq \lambda}$ denote the complementary closed embedding. We will also write $j_\lambda$ to denote the inclusion $X_\lambda \hookrightarrow X$.
    Let $\Shv_G^c(X; \QQ)$ denote the $\infty$-category of $G$-equivariant sheaves on $X$ which are constructible for the $G$-equivariant orbit stratification of $X$.
    Recall that the cohomology functor $\Gamma: \Shv_G^c(X; \QQ) \to \Mod_\QQ$ is given by $\ast$-pushforward to a point and then taking $G$-homotopy fixed points of the resulting $\QQ$-module with $G$-action.
\end{setup}
\begin{definition}\label{def: even sheaves}
    Let $\cf\in \Shv_G^c(X; \QQ)$. Say that $\cf$ is $\ast$-even if the $\ast$-pullback $j_\lambda^\ast \cf\in \Shv_G^c(X_\lambda; \QQ)$ is a direct sum of constant sheaves concentrated in even degrees for all $\lambda \in P$. Similarly, say that $\cf$ is $!$-even if the $!$-pullback $j_\lambda^! \cf\in \Shv_G^c(X_\lambda; \QQ)$ is a direct sum of constant sheaves concentrated in even degrees for all $\lambda \in P$. Say that $\cf$ is even if it is both $\ast$-even and $!$-even. Finally, say that $\cf$ is ($!$- or $\ast$-)odd if $\cf[1]$ is ($!$- or $\ast$-)even.
\end{definition}
The goal of this section is to prove the following result, by inducting on the stratification of $X$:
\begin{theorem}\label{thm: full faithful}
    Let $\cf$ and $\cg$ be even objects of $\Shv_G^c(X; \QQ)$. Then the map
    $$\Ext^\bull_{\Shv_G^c(X; \QQ)}(\cf, \cg) \to \Hom^\bull_{\H^\ast_G(X; \QQ)}(\H^\ast_G(X; \cf), \H^\ast_G(X; \cg))$$
    of graded $\QQ$-vector spaces is a graded isomorphism (the grading denoted by $\bull$), where the $\Hom$ on the right-hand side is taken in the $1$-category of graded $\H^\ast_G(X; \QQ)$-modules in $\Mod_\QQ^\heartsuit$.
\end{theorem}
\begin{remark}
    Although \cref{thm: full faithful} is stated only for $X$ being a \textit{finite} $G$-space, it can be extended easily to the situation when $X$ is not necessarily finite. Namely, suppose that $X$ is a $G$-space equipped with a presentation $X = \colim_{j\in \cJ} X_j$ in terms of finite $G$-spaces where each map $X_j \to X_{j'}$ is a closed embedding. In this case, we will write $\Shv_G^c(X; \QQ)$ to denote the inverse limit $\lim_{j\in \cJ} \Shv_G^c(X_j; \QQ)$ taken over $!$-pullbacks; and the meaning of evenness is exactly as in \cref{def: even sheaves}. With this definition of $\Shv_G^c(X; \QQ)$, \cref{thm: full faithful} continues to hold verbatim.
\end{remark}
\begin{remark}\label{rmk: full faithful is very general}
    The argument for \cref{thm: full faithful} below is sufficiently general that if $A$ is an $\Eoo$-ring with homotopy concentrated in even degrees, given a good theory of $G$-equivariant constructible sheaves of $A$-modules on a stratified (finite) $G$-space (including a six functor formalism), \cref{thm: full faithful} will continue to hold as long as $\pi_\ast A_G$ is concentrated in even degrees. In particular, it continues to hold if $A = \FF_2$. However, we will only describe the argument for \cref{thm: full faithful} with coefficients in $\QQ$.
\end{remark}

\begin{lemma}\label{lem: no exts}
    Let $R$ be an $\Eoo$-ring, and let $M_1 \to M_2 \to M_3$ be a cofiber sequence of $R$-modules such that each of $M_1$, $M_2$, and $M_3$ have homotopy concentrated in even degrees. Then there is a short exact sequence of graded $\pi_\ast R$-modules
    $$0 \to \pi_\ast M_1 \to \pi_\ast M_2 \to \pi_\ast M_3 \to 0.$$
\end{lemma}

\begin{lemma}\label{lem: evenness criterion}
    Let $F: \Shv_G^c(X; \QQ) \to \Mod_\QQ$ be an exact functor. Then $F$ sends $\ast$-even sheaves to $\QQ$-modules with even homotopy groups if and only if $F(j_{\lambda, !} \ul{\QQ})$ has even homotopy groups for each $\lambda \in P$.
\end{lemma}
\begin{proof}
    Suppose that $F$ sends $\ast$-even sheaves to a $\QQ$-module with even homotopy groups. We claim that $F(j_{\lambda, !} \ul{\QQ})$ has even homotopy groups for each $\lambda \in P$: for this, it suffices to show that for each $\lambda'\in \Lambda$, the pullback $j_{\lambda'}^\ast j_{\lambda, !} \ul{\QQ}$ is a direct sum of constant sheaves concentrated in even degrees. But this is clear, because this pullback is zero unless $\lambda' = \lambda$, in which case it is just $\ul{\QQ}$.

    Let us now show the other direction. Let $\cf\in \Shv_G^c(X; \QQ)$ be an $\ast$-even sheaf, and fix $\lambda \in P$ such that $X_\lambda$ is contained in the support of $\cf$.
    Then there is a recollement
    $$j_{\lambda, !} j_\lambda^! \cf \to \cf \to i_{\lambda, \ast} i_\lambda^\ast \cf.$$
    Since $j_\lambda: X_\lambda \hookrightarrow X$ is open, we can identify $j_\lambda^! = j_\lambda^\ast$, and so $j_\lambda^\ast \cf$ is a direct sum of constant sheaves concentrated in even degrees by assumption on $\cf$. This implies that $j_{\lambda, !} j_\lambda^! \cf$ is $\ast$-even (by the argument in the preceding paragraph), so that $F(j_{\lambda, !} j_\lambda^! \cf)$ has even homotopy groups by our assumption on $F$.
    Similarly, by induction on the strata contained in the support of $\cf$, we may assume that $F(i_{\lambda, \ast} i_\lambda^\ast \cf)$ has even homotopy groups. Since $F(\cf)$ is an extension of $F(i_{\lambda, \ast} i_\lambda^\ast \cf)$ by $F(j_{\lambda, !} j_\lambda^! \cf)$, this implies that $F(\cf)$ also has even homotopy groups.
\end{proof}
\begin{lemma}\label{lem: global sections even}
    The functor $\Gamma: \Shv_G^c(X; \QQ) \to \Mod_\QQ$ sends $\ast$-even sheaves to $\QQ$-modules with even homotopy groups.
\end{lemma}
\begin{proof}
    By \cref{lem: evenness criterion}, we need to show that if $\lambda \in \Lambda$, the global sections $\Gamma_G(X; j_{\lambda, !} \ul{\QQ})$ has homotopy concentrated in even degrees. This is true by our assumption on $X_\lambda$.
\end{proof}
\begin{lemma}\label{lem: homming is even}
    Let $\cg \in \Shv_G^c(X; \QQ)$ be $!$-even. Then the functor $\Shv_G^c(X; \QQ) \to \Mod_\QQ$ given by $\Map_{\Shv_G^c(X; \QQ)}(-, \cg)$ sends $\ast$-even sheaves to $\QQ$-modules with even homotopy groups.
\end{lemma}
\begin{proof}
    By \cref{lem: evenness criterion}, we need to show that if $\lambda \in \Lambda$, the $\QQ$-module $\Map_{\Shv_G^c(X; \QQ)}(j_{\lambda, !} \ul{\QQ}, \cg)$ has even homotopy. This $\QQ$-module can be identified with $\Map_{\Shv_G^c(X_\lambda; \QQ)}(\ul{\QQ}, j_\lambda^! \cg) = \Gamma_G(X_\lambda; j_\lambda^! \cg)$. Since $j_\lambda^! \cg$ is a direct sum of constant sheaves concentrated in even degrees (by assumption on $\cg$), the desired result again follows from the assumption that $\H^\ast_G(X_\lambda; \QQ)$ is concentrated in even degrees.
\end{proof}

\begin{lemma}\label{lem: Ext exact sequences}
    Let $\cf\in \Shv_G^c(X; \QQ)$ be $\ast$-even, and let $\cg\in \Shv_G^c(X; \QQ)$ be $!$-even. Then for each $\lambda \in P$ such that $X_\lambda$ is open in the support of $\cf$, there is an exact sequence
    $$0 \to \Ext^\bull_{\Shv_G^c(X_{<\lambda}; \QQ)}(i_\lambda^\ast \cf, i_\lambda^! \cf) \to \Ext^\bull_{\Shv_G^c(X; \QQ)}(\cf, \cg) \to \Ext^\bull_{\Shv_G^c(X_\lambda; \QQ)}(j_\lambda^! \cf, j_\lambda^\ast \cg) \to 0.$$
\end{lemma}
\begin{proof}
    Recall that there is a recollement cofiber sequence
    $$j_{\lambda, !} j_\lambda^! \cf \to \cf \to i_{\lambda, \ast} i_\lambda^\ast \cf.$$
    Applying $\Map_{\Shv_G^c(X; \QQ)}(-, \cg)$ produces a cofiber sequence of $\QQ$-modules
    \begin{equation}\label{eq: cofib homming}
        \Map_{\Shv_G^c(X; \QQ)}(i_{\lambda, \ast} i_\lambda^\ast \cf, \cg) \to \Map_{\Shv_G^c(X; \QQ)}(\cf, \cg) \to \Map_{\Shv_G^c(X; \QQ)}(j_{\lambda, !} j_\lambda^! \cf, \cg).
    \end{equation}
    Observe that $i_{\lambda, \ast} i_\lambda^\ast \cf$ and $j_{\lambda, !} j_\lambda^! \cf$ are both $\ast$-even, so that \cref{lem: homming is even} implies that each term in \cref{eq: cofib homming} has even homotopy. In particular, \cref{lem: no exts} implies that \cref{eq: cofib homming} induces a split exact sequence on homotopy groups.
    Note that by adjunction, we can rewrite the first term of \cref{eq: cofib homming} as $\Map_{\Shv_G^c(X_{<\lambda}; \QQ)}(i_\lambda^\ast \cf, i_\lambda^! \cg)$, and the final term of \cref{eq: cofib homming} as $\Map_{\Shv_G^c(X_\lambda; \QQ)}(j_\lambda^! \cf, j_\lambda^\ast \cg)$. Together with the above discussion, this proves the desired claim.
\end{proof}
\begin{lemma}\label{lem: coh exact sequence shriek and star}
    Let $\cf\in \Shv_G^c(X; \QQ)$ be $\ast$-even, and let $\cg\in \Shv_G^c(X; \QQ)$ be $!$-even. Then for each $\lambda \in P$ such that $X_\lambda$ is open in the support of $\cf$, there are exact sequences of graded $\QQ$-vector spaces
    \begin{align*}
        0 \to \H^\ast_G(X; j_{\lambda, !} j_\lambda^! \cf) & \to \H^\ast_G(X; \cf) \to \H^\ast_G(X_{<\lambda}; i_\lambda^\ast \cf) \to 0, \\
        0 \to \H^\ast_G(X_{<\lambda}; i_\lambda^! \cg) & \to \H^\ast_G(X; \cg) \to \H^\ast_G(X_{\lambda}; j_\lambda^\ast \cg) \to 0.
    \end{align*}
\end{lemma}
\begin{proof}
    We will only prove the first exact sequence; the second follows by an entirely analogous argument.
    Again, recall that there is a recollement cofiber sequence
    $$j_{\lambda, !} j_\lambda^! \cf \to \cf \to i_{\lambda, \ast} i_\lambda^\ast \cf.$$
    Applying $\Gamma_G(X; -)$ gives a cofiber sequence
    $$\Gamma_G(X; j_{\lambda, !} j_\lambda^! \cf) \to \Gamma_G(X; \cf) \to \Gamma_G(X_{<\lambda}; i_\lambda^\ast \cf).$$
    Observe that $i_{\lambda, \ast} i_\lambda^\ast \cf$ and $j_{\lambda, !} j_\lambda^! \cf$ are both $\ast$-even, so that \cref{lem: global sections even} implies that each term in this cofiber sequence has even homotopy. In particular, \cref{lem: no exts} implies that this cofiber sequence induces a split exact sequence on homotopy groups, as desired.
\end{proof}

\begin{lemma}\label{lem: fdtl class affine space}
    Let $V$ be a complex affine space on which $G$ acts linearly, and equip $V$ with the trivial stratification. Then:
    \begin{enumerate}
        \item The functor $\H^\ast_G(V; -): \Shv_G^c(\ast; \QQ) \xar{\sim} \Mod_{\H_G^\ast(\ast; \QQ)}(\Sp)$ is an equivalence, where the right-hand side denotes the $\infty$-category of $\H^\ast_G(V; \QQ)$-modules in spectra (i.e., the derived $\infty$-category of chain complexes of $\H^\ast_G(V; \QQ)$-modules).
        
        Moreover, if $\cf, \cg\in \Shv_G^c(V; \QQ)$ are sheaves such that $\H^\ast_G(V; \cf)$ is a projective $\H^\ast_G(V; \QQ)$-module, there is a graded isomorphism
        $$\Ext^\bull_{\Shv_G^c(V; \QQ)}(\cf, \cg) \xar{\sim} \Hom^\bull_{\H^\ast_G(V; \QQ)}(\H^\ast_G(V; \cf), \H^\ast_G(V; \cg)).$$
        Here, the $\Hom$ on the right-hand side is taken in the $1$-category of graded $\H^\ast_G(X; \QQ)$-modules in $\Mod_\QQ^\heartsuit$.
        \item The compactly supported equivariant cohomology $\H^\ast_{G,c}(V; \QQ)$ is isomorphic to a free $\H^\ast_G(\ast; \QQ)$-module generated by a single class $[V]$ in degree $\dim_\RR(V)$.
    \end{enumerate}
\end{lemma}
\begin{proof}
    Let us first show (a).
    Since $V$ is equipped with the trivial stratification, $\Shv_G^c(V; \QQ)$ is equivalent to the $\infty$-category of $G$-equivariant local systems on $V$. Because $V$ is a complex affine space, this is simply equivalent to the $\infty$-category $\Loc_G(\ast; \QQ)$ of $G$-equivariant local systems of $\QQ$-modules on a point. Almost by definition, there is an equivalence $\Loc_G(\ast; \QQ) \simeq \Mod_{C_G^\ast(\ast; \QQ)}(\Sp)$. Since $G$ is assumed to be a compact Lie group, the $\QQ$-algebra $\pi_\ast C_G^\ast(\ast; \QQ) = \H_G^{-\ast}(\ast; \QQ)$ is isomorphic to a graded polynomial $\QQ$-algebra on generators in even negative (homological) degrees. Since the free $\Eoo$-$\QQ$-algebra on a generator $x$ in even degree is isomorphic to the polynomial algebra $\QQ[x]$, choosing polynomial generators for $\pi_\ast C_G^\ast(\ast; \QQ)$ defines an equivalence $\H_G^\ast(\ast; \QQ) \cong C_G^\ast(\ast; \QQ)$ of $\Eoo$-$\QQ$-algebras. It follows that there is an equivalence $\Loc_G(\ast; \QQ) \cong \Mod_{\H_G^\ast(\ast; \QQ)}(\Sp)$.
    Finally, if $\cf, \cg\in \Shv_G^c(V; \QQ)$ are sheaves such that $\H^\ast_G(V; \cf)$ is a projective graded $\H^\ast_G(V; \QQ)$-module, the spectral sequence
    \begin{equation}\label{eq: sseq ext}
        E_2 = \Ext^\bull_{\H_G^\ast(\ast; \QQ)}(\H^\ast_G(\ast; \cf), \H^\ast_G(\ast; \cg)) \Rightarrow \pi_{-\ast} \Hom_{\Mod_{\H_G^\ast(\ast; \QQ)}(\Sp)}(\H^\ast_G(\ast; \cf), \H^\ast_G(\ast; \cg))
    \end{equation}
    degenerates at the $E_2$-page, where $\Ext^\bull$ denotes the \textit{graded} $\Ext$-groups taken internal to the $1$-category of graded $\H_G^\ast(\ast; \QQ)$-modules in $\Mod_\QQ^\heartsuit$. In fact, the $E_2$-page is concentrated entirely in the zero line, since $\H^\ast_G(V; \cf)$ is a projective graded $\H^\ast_G(V; \QQ)$-module (so there are no higher $\Ext$-groups). In particular, we have
    $$\Ext^\bull_{\Shv_G^c(V; \QQ)}(\cf, \cg) \xar{\sim} \Hom^\bull_{\H^\ast_G(V; \QQ)}(\H^\ast_G(V; \cf), \H^\ast_G(V; \cg)),$$
    as desired.
    
    Part (b) is simply the statement of (equivariant) Poincar\'e duality on an affine space.
\end{proof}
\begin{remark}
    It is natural to ask whether the statement of \cref{lem: fdtl class affine space}(a) is true for arbitrary $\cf$. Unfortunately, this need not be true. For example, suppose that $G = S^1$, and that $V$ is the trivial vector space (without loss of generality). Using the equivalence $\Loc_{S^1}(\ast; \QQ) \xar{\sim} \Mod_{\H^\ast_{S^1}(\ast; \QQ)}(\Sp)$ and the fact that $\H^\ast_{S^1}(\ast; \QQ) = \QQ[x]$ with $x$ in $\H^2_{S^1}(\ast; \QQ)$, one can define an $S^1$-local system on the point by the $\H^\ast_{S^1}(\ast; \QQ)$-module $\QQ[x]/x = \QQ$.
    The extension class
    $$\QQ[-2] \to \QQ[x]/x^2 \to \QQ$$
    defines a nontrivial element $\delta\in \pi_1 \Hom_{\Mod_{\H^\ast_{S^1}(\ast; \QQ)}(\Sp)}(\QQ, \QQ)$. However, this class cannot be seen from the graded $\Hom$ group: indeed,
    $$\Hom^\bull_{\H_{S^1}^\ast(\ast; \QQ)}(\QQ, \QQ) = \QQ$$
    in weight zero, generated by multiples of the identity map.
    The spectral sequence \cref{eq: sseq ext} still degenerates at the $E_2$-page in this case, but there is a nontrivial class in $\Ext^1_{\H_{S^1}^\ast(\ast; \QQ)}(\QQ, \QQ)$ which detects the class $\delta$.

    The above example is intended to illustrate the difference between the $\infty$-category of $\H^\ast_{S^1}(\ast; \QQ)$-modules in spectra, which can be identified with the derived $\infty$-category of chain complexes of $\H^\ast_{S^1}(\ast; \QQ)$-modules, and the category of graded $\H^\ast_{S^1}(\ast; \QQ)$-modules in $\Mod_\QQ^\heartsuit$.
\end{remark}

\begin{lemma}
    Let $R$ be a graded (discrete) $\QQ$-algebra.
    Fix two exact sequences
    \begin{align*}
        0 \to M_1 \to M_2 & \to M_3 \to 0,\\
        0 \to N_1 \to N_2 & \to N_3 \to 0
    \end{align*}
    of graded $R$-modules. Then there is a sequence
    $$0 \to \Hom^\bull_R(M_3, N_1) \to \Hom^\bull_R(M_2, N_2) \to \Hom^\bull_R(M_1, N_3)$$
    which is exact on the left (i.e., the second map is injective). Here, the $\Hom$s are taken in the $1$-category of graded $R$-modules in $\Mod_\QQ^\heartsuit$.

    If\footnote{In words: every graded $R$-linear map $M_2 \to N_1$ factors through the surjection $M_2 \twoheadrightarrow M_3$, and every graded $R$-linear map $M_1 \to N_3$ extends along $M_1 \hookrightarrow M_2$.} the maps
    \begin{align*}
        \Hom^\bull_R(M_3, N_1) & \hookrightarrow \Hom^\bull_R(M_2, N_1), \\
        \Hom^\bull_R(M_2, N_3) & \twoheadrightarrow \Hom^\bull_R(M_1, N_3)
    \end{align*}
    are isomorphisms, the above sequence is also exact in the middle.
\end{lemma}
\begin{proof}
    Exactness on the left is clear, since the map $\alpha$ factors as injections
    $$\Hom^\bull_R(M_3, N_1) \hookrightarrow \Hom^\bull_R(M_3, N_2) \hookrightarrow \Hom^\bull_R(M_2, N_2).$$
    Exactness in the middle given the assumptions follows from noting that the desired sequence can be written as the composite
    $$\xymatrix{
    & \Hom^\bull_R(M_3, N_1) \ar@{^(->}[d]^-\sim & & \Hom^\bull_R(M_1, N_3) \\
    0 \ar[r] & \Hom^\bull_R(M_2, N_1) \ar@{^(->}[r] & \Hom^\bull_R(M_2, N_2) \ar[r] & \Hom^\bull_R(M_2, N_3), \ar@{->>}[u]^-\sim
    }$$
    where the bottom row is exact in the middle (by left exactness of $\Hom$).
\end{proof}

\begin{lemma}\label{lem: exact sequences on hom}
    Let $\cf\in \Shv_G^c(X; \QQ)$ be $\ast$-even, and let $\cg\in \Shv_G^c(X; \QQ)$ be $!$-even. Suppose that the map $\H^\ast_G(X; \cf) \to \H^\ast_G(X_\lambda; j_\lambda^\ast \cf)$ is surjective for every $\lambda \in P$, and that the map $\H^\ast_G(X_\lambda; j_\lambda^! \cg) \to \H^\ast_G(X; \cg)$ is injective for every $\lambda \in P$. Then there is an exact sequence
    \begin{align*}
        \hspace*{-1cm} 0 \to \Hom^\bull_{\H^\ast_G(X_{<\lambda}; \QQ)}(\H^\ast_G(X_{<\lambda}; i_\lambda^\ast \cf), \H^\ast_G(X_{<\lambda}; i_\lambda^! \cg) & \to \Hom^\bull_{\H^\ast_G(X; \QQ)}(\H^\ast_G(X; \cf), \H^\ast_G(X; \cg))\\
        & \to \Hom^\bull_{\H^\ast_G(X_\lambda; \QQ)}(\H^\ast_G(X; j_{\lambda, !} j_\lambda^! \cf), \H^\ast_G(X_\lambda; j_\lambda^\ast \cg)).
    \end{align*}
    Here, the $\Hom$s are taken in the $1$-category of graded $\H^\ast_G(X_{<\lambda}; \QQ)$-modules (resp. $\H^\ast_G(X; \QQ)$- and $\H^\ast_G(X_\lambda; \QQ)$-modules) in $\Mod_\QQ^\heartsuit$.
\end{lemma}
\begin{proof}
    Applied to $R = \H^\ast(X; \QQ)$ and the exact sequences of \cref{lem: coh exact sequence shriek and star}, we see that the composite of \cref{lem: exact sequences on hom} is exact on the left.
    For exactness in the middle, \cref{lem: exact sequences on hom} says that we need to check:
    \begin{enumerate}
        \item Every graded $\H^\ast_G(X; \QQ)$-linear map $\H^\ast_G(X; \cf) \to \H^\ast_G(X_{<\lambda}; i_\lambda^! \cg)$ factors through the map $\H^\ast_G(X; \cf) \twoheadrightarrow \H^\ast_G(X_{<\lambda}; i_\lambda^\ast \cf)$.
        \item Every graded $\H^\ast_G(X; \QQ)$-linear map $\H^\ast_G(X; j_{\lambda, !} j_\lambda^! \cf) \to \H^\ast_G(X_\lambda; j_\lambda^\ast \cg)$ extends through the map $\H^\ast_G(X; j_{\lambda, !} j_\lambda^! \cf) \hookrightarrow \H^\ast_G(X; \cf)$.
    \end{enumerate}
    The proof of (b) is entirely analogous to that of (a), so we will only show (a). 
    Recall that $X_\lambda$ was assumed to be a complex affine space on which $G$ acts linearly. Therefore, the compactly supported equivariant cohomology $\H^\ast_{G,c}(X_\lambda; \QQ)$ is isomorphic to a free $\H^\ast_G(\ast; \QQ)$-module on a single generator (by \cref{lem: fdtl class affine space}(b)). Let us denote this generator by $[X_\lambda]$.
    
    Suppose we are given a graded $\H^\ast_G(X; \QQ)$-linear map $f: \H^\ast_G(X; \cf) \to \H^\ast_G(X_{<\lambda}; i_\lambda^! \cg)$. Then the following composite is zero
    $$\H^\ast_G(X; \cf) \xar{f} \H^\ast_G(X_{<\lambda}; i_\lambda^! \cg) \xar{\cdot [X_\lambda]} \H^\ast_G(X_{<\lambda}; i_\lambda^! \cg)[2n_\lambda],$$
    because $[X_\lambda]$ is zero in $\H^\ast_G(X_{<\lambda}; \QQ)$. Since $f$ is $\H^\ast_G(X; \QQ)$-linear, this implies that the following composite is also zero:
    $$\H^\ast_G(X; \cf)[-2n_\lambda] \xar{\cdot [X_\lambda]} \H^\ast_G(X; \cf) \xar{f} \H^\ast_G(X_{<\lambda}; i_\lambda^! \cg).$$
    In particular, the map $f: \H^\ast_G(X; \cf) \to \H^\ast_G(X_{<\lambda}; i_\lambda^! \cg)$ factors through the quotient $\H^\ast_G(X; \cf)/\im([X_\lambda])$. In order to show that the map $f$ factors through the map $\H^\ast_G(X; \cf) \twoheadrightarrow \H^\ast_G(X_{<\lambda}; i_\lambda^\ast \cf)$, it suffices to show that there is a dotted injection making the following diagram commute:
    $$\xymatrix{
    \H^\ast_G(X; \cf) \ar@{->>}[d] \ar@{->>}[dr] & \\
    \H^\ast_G(X_{<\lambda}; i_\lambda^\ast \cf) \ar@{^(-->}[r] & \H^\ast_G(X; \cf)/\im([X_\lambda]).
    }$$
    Equivalently, using \cref{lem: coh exact sequence shriek and star}, it suffices to show that $\H^\ast_G(X; j_{\lambda, !} j_\lambda^! \cf) \hookrightarrow \H^\ast_G(X; \cf)$ is contained in $\im([X_\lambda])$. 
    By Poincar\'e duality on the affine space $X_\lambda$, multiplication by $[X_\lambda]$ defines an isomorphism
    $$\H^\ast_G(X; j_{\lambda, \ast} j_\lambda^\ast \cf)[-2n_\lambda] \xar{\cdot [X_\lambda]} \H^\ast_G(X; j_{\lambda, !} j_\lambda^! \cf).$$
    Observe that there is a commutative diagram
    $$\xymatrix{
    & \H^\ast_G(X; \cf)[-2n_\lambda] \ar[r]^-{\cdot [X_\lambda]} \ar[d] & \H^\ast_G(X; \cf) \\
    \H^\ast_G(X_\lambda; j_\lambda^\ast \cf) \ar[r]^-\sim & \H^\ast_G(X; j_{\lambda, \ast} j_\lambda^\ast \cf)[-2n_\lambda] \ar[r]^-\sim_-{\cdot [X_\lambda]} & \H^\ast_G(X; j_{\lambda, !} j_\lambda^! \cf), \ar@{^(->}[u]
    }$$
    To show that $\H^\ast_G(X; j_{\lambda, !} j_\lambda^! \cf) \hookrightarrow \H^\ast_G(X; \cf)$ is contained in $\im([X_\lambda])$, it suffices that the left vertical map be surjective; but this is precisely our assumption on $\cf$.
\end{proof}

\begin{lemma}\label{lem: hypotheses-descend}
    Let $\cf\in \Shv_G^c(X; \QQ)$ be $\ast$-even, and suppose that the map $\H^\ast_G(X; \cf) \to \H^\ast_G(X_\lambda; j_\lambda^\ast \cf)$ is surjective for every $\lambda \in P$. Then the same is true of $i_\mu^\ast \cf$ for any $\mu\in P$.

    Similarly, let $\cg\in \Shv_G^c(X; \QQ)$ be $!$-even, and suppose that the map $\H^\ast_G(X_\lambda; j_\lambda^! \cg) \to \H^\ast_G(X; \cg)$ is injective for every $\lambda \in P$. Then the same is true of $i_\mu^! \cg$ for any $\mu\in P$.
\end{lemma}
\begin{proof}
    The proof for $\cg$ is analogous to the proof for $\cf$, so we will only prove the latter. First, it is clear that $i_\mu^\ast \cf$ is $\ast$-even for any $\mu\in P$. To prove the surjectivity claim, it evidently suffices to assume that $X_\lambda$ is contained in the support of $i_\mu^\ast \cf$ (otherwise $\H^\ast_G(X_\lambda; j_\lambda^\ast \cf) = 0$). Let $\lambda$ be such that $X_\lambda$ is contained in $X_{<\mu} = X_{\leq \mu} - X_\mu$, so that $i_\mu: X_{<\mu} \hookrightarrow \supp(\cf)$ is the inclusion. Then, we have maps
    $$\H^\ast_G(X; \cf) \to \H^\ast_G(X_{<\mu}; i_\mu^\ast \cf) \to \H^\ast_G(X_\lambda; j_\lambda^\ast \cf).$$
    The composite is surjective, and hence the map $\H^\ast_G(X_{<\mu}; i_\mu^\ast \cf) \to \H^\ast_G(X_\lambda; j_\lambda^\ast \cf)$ is surjective. This gives the desired claim, since $j_\lambda^\ast \cf \cong j_\lambda^\ast i_\mu^\ast \cf$.
\end{proof}

\begin{proof}[Proof of \cref{thm: full faithful}]
    Let us begin by showing that if $\cf\in \Shv_G^c(X; \QQ)$ is even, the map $\H^\ast_G(X; \cf) \to \H^\ast_G(X_\lambda; j_\lambda^\ast \cf)$ is surjective for every $\lambda \in P$; and that if $\cg\in \Shv_G^c(X; \QQ)$ is even, the map $\H^\ast_G(X_\lambda; j_\lambda^! \cg) \to \H^\ast_G(X; \cg)$ is injective for every $\lambda \in P$. 
    The claim for $\cf$ follows from the fact that it is even, and hence $!$-even, and the second exact sequence of \cref{lem: coh exact sequence shriek and star}. Similarly, the claim for $\cg$ follows from the fact that it is even, and hence $\ast$-even, and the first exact sequence of \cref{lem: coh exact sequence shriek and star}.
    
    Let us now prove \cref{thm: full faithful}. Assume that $\cf\in \Shv_G^c(X; \QQ)$ is 
    even, and that $\cg\in \Shv_G^c(X; \QQ)$ is 
    even.
    The preceding paragraph implies that the assumptions of \cref{lem: exact sequences on hom} are satisfied.
    We will show by induction that the canonical map $\Ext^\bull_{\Shv_G^c(X; \QQ)}(\cf, \cg) \to \Hom^\bull_{\H^\ast_G(X; \QQ)}(\H^\ast_G(X; \cf), \H^\ast_G(X; \cg))$ is an isomorphism. Let $\lambda \in P$ be such that $X_\lambda$ is open in the union of the supports of $\cf$ and $\cg$.
    There is a map of sequences
    $$\xymatrix{
    \Ext^\bull_{\Shv_G^c(X_{<\lambda}; \QQ)}(i_\lambda^\ast \cf, i_\lambda^! \cf) \ar[r]^-\alpha \ar[d] & \Hom^\bull_{\H^\ast_G(X_{<\lambda}; \QQ)}(\H^\ast_G(X_{<\lambda}; i_\lambda^\ast \cf), \H^\ast_G(X_{<\lambda}; i_\lambda^! \cg) \ar[d] \\
    \Ext^\bull_{\Shv_G^c(X; \QQ)}(\cf, \cg) \ar[r]^-\beta \ar[d] & \Hom^\bull_{\H^\ast_G(X; \QQ)}(\H^\ast_G(X; \cf), \H^\ast_G(X; \cg)) \ar[d] \\
    \Ext^\bull_{\Shv_G^c(X_\lambda; \QQ)}(j_\lambda^! \cf, j_\lambda^\ast \cg) \ar[r]_-\gamma & \Hom^\bull_{\H^\ast_G(X_\lambda; \QQ)}(\H^\ast_G(X; j_{\lambda, !} j_\lambda^! \cf), \H^\ast_G(X_\lambda; j_\lambda^\ast \cg)).
    }$$
    By \cref{lem: Ext exact sequences}, the left vertical composite is a short exact sequence. By \cref{lem: exact sequences on hom}, the right vertical composite is left exact (i.e., the first map is injective, and it is exact in the middle). By \cref{lem: hypotheses-descend} and the inductive hypothesis, the map denoted $\alpha$ is a graded isomorphism. The map $\gamma$ is also a graded isomorphism: indeed, since $\cf$ is $!$-even, $j_\lambda^! \cf$ is a direct sum of constant sheaves (in even degrees); similarly, since $\cg$ is $\ast$-even, $j_\lambda^\ast \cg$ is also a direct sum of constant sheaves (in even degrees). In particular, $\H^\ast_G(X_\lambda; j_\lambda^! \cf)$ is a projective $\H^\ast_G(X_\lambda; \QQ)$-module, so \cref{lem: fdtl class affine space}(a) implies that the map denoted $\gamma$ is also a graded isomorphism. This implies that the map denoted $\beta$ is also a graded isomorphism, as desired.
\end{proof}

%% file: equivalences/geometric-satake-review.tex
\subsection{Review of derived geometric Satake}\label{subsec: review satake}

Using the results of the preceding section, we will now review an argument for the derived geometric Satake equivalence of \cite{bf-derived-satake}, because it will serve as a model for the arguments appearing later in this article. We will use this as an opportunity to review some facts from geometric representation theory, and to point out some simplifications coming from homotopy theory (see \cref{rmk: e2-formality in satake}).
\begin{setup}
    Fix a connected semisimple compact Lie group $G$ throughout this section, and let $G_\cc$ denote the associated complex algebraic group over $\cc$ (so that $G_\cc$ is a connected semisimple group). Let $\ld{G}$ denote the Chevalley form over $\Z$ of the split semisimple algebraic group (over $\Z$) whose root datum is Langlands dual to that of $G_\cc$. We will simply write $\ld{G}$ to denote the base-change of $\ld{G}$ to $\QQ$; the base over which $\ld{G}$ is defined will be clear from context.

    Later in this section, we will also fix a Borel subgroup $\ld{B} \subseteq \ld{G}$, and write $\ld{N}$ to denote its unipotent radical. Let $\Phi$ denote the set of roots of $\ld{G}$ (i.e., coroots of $G_\cc$), and $\Lambda$ will denote the character lattice of $\ld{G}$ (i.e., cocharacter lattice of $G_\cc$). The choice of $\ld{B}$ defines a subset $\Phi^+$ of positive roots, and $\Delta\subseteq \Phi^+$ will denote a base of simple roots. Let $\Lambda^+$ denote the subset of dominant weights of $\ld{G}$, so that $\Lambda^+$ is in (order-preserving) bijection with the set of orbits of the $G$-action on $\Omega G$ by \cite[Theorem 1.6.1 and Equation 2.1.1]{zhu-grass}. We will add checks above each of these symbols to denote coroots, positive coroots, and simple coroots, respectively. 
\end{setup}
\begin{definition}
    Let $\Shv_{G\pw{t}}^c(\Gr_G; \QQ)$ denote the $\infty$-category of $G$-equivariant sheaves of $\QQ$-modules on $\Gr_G(\cc)$ which are constructible for the orbit stratification. More precisely, if we write $\Gr_G(\cc)$ as the direct limit $\colim_{\lambda \in \Lambda^+} \Gr_G^{\leq \lambda}(\cc)$ of the finite-dimensional $G$-equivariant Schubert strata $\Gr_G^{\leq \lambda}(\cc)$ for $\lambda \in \Lambda^+$, the $\infty$-category $\Shv_{G\pw{t}}^c(\Gr_G; \QQ)$ is defined to be the inverse limit $\lim_{\lambda\in \Lambda^+} \Shv_G^c(\Gr_G^{\leq \lambda}(\cc); \QQ)$ along $!$-pullbacks.
\end{definition}
In addition to \cref{thm: full faithful}, there are two key results needed to prove the derived Satake theorem. The first of these is the following.
\begin{theorem}[Abelian geometric Satake, \cite{mirkovic-vilonen}]\label{thm: abelian satake}
    Let $\Perv_{G(\co)}(\Gr_G; \QQ)$ denote the abelian $1$-category of $G(\co)$-equivariant perverse sheaves on $\Gr_G$, so that $\Perv_{G(\co)}(\Gr_G; \QQ)$ admits a symmetric monoidal structure arising from convolution on the affine Grassmannian.
    Then there is a symmetric monoidal equivalence $\Perv_{G(\co)}(\Gr_G; \QQ) \simeq \Rep(\ld{G})$.
\end{theorem}
\begin{remark}
    We will not need to spell out the definition of $\Perv_{G(\co)}(\Gr_G; \QQ)$ in the remainder of this article; in fact, all that we will need is the consequence that there is a fully faithful functor $\Rep(\ld{G}) \hookrightarrow \Shv_{G\pw{t}}^c(\Gr_G; \QQ)$.
\end{remark}

\begin{definition}\label{def: Shv-Sat GrG}
    The action of $\Perv_{G(\co)}(\Gr_G; \QQ)$ on $\Shv_{G\pw{t}}^c(\Gr_G; \QQ)$ via convolution defines, via \cref{thm: abelian satake}, an action of $\Rep(\ld{G})$ on $\Shv_{G\pw{t}}^c(\Gr_G; \QQ)$.
    Let $\IC_0\in \Shv_{G\pw{t}}^c(\Gr_G; \QQ)$ denote the pushforward $i_! \ul{\QQ}$ of the constant sheaf along the inclusion $i: \{\ast\} \hookrightarrow \Omega G$ of the basepoint.
    Let $\Shv_{G\pw{t}}^{c,\Sat}(\Gr_G; \QQ)$ denote the full subcategory of $\Shv_{G\pw{t}}^c(\Gr_G; \QQ)$ generated by $\IC_0$ under the action of $\Rep(\ld{G})$.
\end{definition}
\begin{notation}
    Let $2\rho = \sum_{\ld{\alpha}\in \ld{\Phi}^+} \ld{\alpha}$ denote the sum of the positive coroots, so that it defines a homomorphism $2\rho: \GG_m \to \ld{T}$. The adjoint action defines an action of $\GG_m$ on $\ld{\g}$ via $(2 - 2\rho)$, which fixes the element $e$. The adjoint action of $\ld{G}$ on $\ld{\g}$ refines to a \textit{graded} action if $\ld{G}$ is equipped with the grading coming from $2\rho$. We will \textit{only} view $\ld{\g}$ (resp. $\ld{G}$) as a graded scheme via the $(2 - 2\rho)$-action (resp. $-2\rho$-action). To emphasize this, we will denote these graded schemes as $\ld{\g}(2 - 2\rho)$ and $\ld{G}(-2\rho)$.

    For instance, if $\ld{G} = \SL_2$, the grading equips the entries of an element $\begin{psmallmatrix}
        a & b\\
        c & -a
    \end{psmallmatrix} \in \sl_2$ with the gradings where $a$ lives in weight $-2$, $b$ lives in weight $0$, and $c$ lives in weight $-4$; similarly, the entries an element $\begin{psmallmatrix}
        a & b\\
        c & d
    \end{psmallmatrix} \in \SL_2$ have the gradings here $a$ and $d$ live in weight $0$, $b$ lives in weight $2$, and $c$ lives in weight $-2$.
\end{notation}
We then have (see \cite{bf-derived-satake}):
\begin{theorem}[Derived geometric Satake]\label{thm: derived satake}
    There is a monoidal equivalence 
    $$\Shv_{G\pw{t}}^{c, \Sat}(\Gr_G; \QQ) \simeq \Perf(\sh^{1/2}\ld{\g}(2 - 2\rho)/\ld{G}(-2\rho))$$
    of $\QQ$-linear $\infty$-categories.
\end{theorem}
\begin{remark}\label{rmk: remarks about derived satake}
    Let us make two comments regarding \cref{thm: derived satake}. 
    \begin{enumerate}
        \item The shifts appearing in \cref{thm: derived satake} are different than those which appear in \cite{bf-derived-satake}; this is because \cref{thm: derived satake} is stated using the \textit{arithmetic shearing} of \cite[Section 6.7]{bzsv}.
        \item One can wonder about an analogue of \cref{thm: derived satake} for sheaves with coefficients in other commutative rings. The discussion in \cite{campbell-raskin-satake} suggests that there should be an equivalence 
        $$\Shv_{G\pw{t}}^{c, \Sat}(\Gr_G; \Z) \simeq \Coh((\{1\} \times_{\ld{G}_\Z} \{1\})/\ld{G}_\Z),$$
        with $\ld{G}_\Z$ being the Chevalley split form of the Langlands dual. Presumably this can be proved using the ideas of \cite{campbell-raskin-satake}. The $\infty$-category on the right-hand side of this equivalence is generally \textit{not} equivalent to $\Coh((\{0\} \times_{\ld{\g}_\Z} \{0\})/\ld{G}_\Z) \simeq \Perf(\ld{\g}_\Z^\ast[2]/\ld{G}_\Z)$; but such an equivalence will exist after inverting some integer $n$.
        \footnote{Indeed, there is an isomorphism $\{1\} \times_{\ld{G}_\Z} \{1\} \cong \{1\} \times_{\widehat{\ld{G}}_\Z} \{1\}$, where $\widehat{\ld{G}}_\Z$ denotes the completion of $\ld{G}_\Z$ at $1\in \ld{G}_\Z$. It therefore suffices to show that there is an isomorphism $\widehat{\ld{G}}_{\Z[1/n]} \cong \widehat{\ld{\g}}_{\Z[1/n]}$ after inverting some integer $n$, where $\widehat{\ld{\g}}_{\Z[1/n]}$ denotes the completion of $\ld{\g}_{\Z[1/n]}$ at the origin.
        For this, it suffices that $\ld{G}_{\Z[1/n]}$ admit the datum of an $\ld{G}_{\Z[1/n]}$-equivariant splitting of the map $\cI \to \cI/\cI^2 = \ld{\g}_{\Z[1/n]}^\ast$, where $\cI$ is the ideal sheaf cutting out $\{1\} \subseteq \ld{G}_{\Z[1/n]}$. (Since $\ld{G}_{\Z[1/n]}$ is smooth, this is exactly the datum of a \textit{quasi-logarithm} $\ld{G}_{\Z[1/n]} \to \ld{\g}_{\Z[1/n]}$ in the sense of \cite[Section 1.8]{kazhdan-varshavsky-quasilog}.) However, such a splitting exists over $\QQ$ since $\Rep(\ld{G}_\QQ)$ is semisimple, and hence is defined over $\Z[1/n]$ for some $n\gg 0$, as desired. Sometimes such a splitting exists without inverting any primes: for example, such a quasi-logarithm exists over $\Z$ for $\ld{G} = \GL_n$ (see \cite[Proposition 6.3]{friedlander-negron}).}
    \end{enumerate}
\end{remark}
To prove \cref{thm: derived satake}, we need a few ingredients. First, we need to explain the theory of the Kostant slice.
\begin{recall}
    Since $\ld{\g}$ is semisimple, the map $\ld{\g} \to \ld{\g}^\ast$ from the Killing form is a $\ld{G}$-equivariant isomorphism. This implies that $\ld{\g} \mmod \ld{G} \cong \ld{\g}^\ast \mmod \ld{G}$, and the Chevalley restriction theorem gives an isomorphism $\ld{\fr{g}}^\ast \mmod \ld{G} \cong \fr{t}\mmod W$. In particular, there is an isomorphism $\fr{t}\mmod W \cong \ld{\g}\mmod \ld{G}$.
    There is also an \textit{ungraded} isomorphism $\spec \H^\ast_G(\ast; \QQ) \cong \fr{t} \mmod W$.
    
    This can be upgraded to a graded isomorphism as follows. Let $\ld{\g}(2) = \spec \Sym^\ast_\QQ(\ld{\g}^\ast(-2))$ denote the $\QQ$-vector space where the coordinate has weight $-2$. Similarly, let $\fr{t}(2)$ denote $\spec \Sym^\ast_\QQ(\fr{t}^\ast(-2))$ denote the $\QQ$-vector space where the coordinate has weight $-2$. Then, there is a graded isomorphism
    $$\spec \H^\ast_G(\ast; \QQ) \cong \fr{t}(2) \mmod W \cong \ld{\g}(2)\mmod \ld{G}.$$
\end{recall}
\begin{definition}\label{def: slice kostant def}
    Fix a Borel subgroup $\ld{B} \subseteq \ld{G}$, and let $\ld{N}$ denote its unipotent radical.
    Let $e\in \ld{\g}$ be a principal nilpotent element, i.e., an element $e\in \ld{\fr{n}}$ such that for each simple root $\alpha\in \Delta$, the image of $e$ under the following composite is nonzero:
    $$\ld{\fr{n}} \twoheadrightarrow \ld{\fr{n}}/[\ld{\fr{n}}, \ld{\fr{n}}] \cong \prod_{\alpha\in \Delta} \GG_a \xar{\pr_\alpha} \GG_a.$$
    Let $i: \sl_2 \to \ld{\g}$ denote the Lie algebra homomorphism arising from the Jacobson-Morozov theorem, so that $i$ sends $\begin{psmallmatrix}
        0 & 1\\
        0 & 0
    \end{psmallmatrix}\in \sl_2$ to $e\in \ld{\g}$. Let $h \in \ld{\g}$ denote the image of $\begin{psmallmatrix}
        1 & 0\\
        0 & -1
    \end{psmallmatrix}\in \sl_2$, and let $f\in \ld{\g}$ denote the image of $\begin{psmallmatrix}
        0 & 0\\
        -1 & 0
    \end{psmallmatrix}\in \sl_2$, so that $f\in \fr{b}^-$. The adjoint action of $h$ equips $e$ with weight $2$ and $f$ with weight $-2$. In fact, the adjoint action of $h$ on $\ld{\g}$ equips it with a grading where all the weights are \textit{even} integers.

    Let $\ld{\g}^f$ denote the centralizer of $f$ in $\ld{\g}$. 
    The \textit{Kostant slice} is defined to be the affine subspace $\kappa: e + \ld{\g}^f \hookrightarrow \ld{\g}$, so that $\kappa$ is a closed immersion. The grading on $\ld{\g}$ via $(2 - 2\rho)$ restricts to a grading on $e + \ld{\g}^f \subseteq \ld{\g}$.
\end{definition}
\begin{theorem}[{Kostant, \cite[Theorem 1.2]{kostant-whittaker}}]\label{thm: kostant}
    Fix notation as in \cref{def: slice kostant def}, and equip $e+\ld{\fr{b}}^-$ with the grading coming from the adjoint action of $2 - 2\rho$. Then the adjoint action of $\ld{N}^-$ on $e+\ld{\fr{b}}^-$ is one of graded schemes.
    The natural maps
    $$e + \ld{\g}^f \to (e + \ld{\fr{b}}^-)/\ld{N}^- \to \ld{\g}(2)\mmod \ld{G} \leftarrow \fr{t}(2)\mmod W$$
    are isomorphisms of graded schemes. Moreover, the closed subscheme $e + \ld{\g}^f \subseteq \ld{\g}(2 - 2\rho)$ meets every regular orbit exactly once, and transversally so.
\end{theorem}
\begin{remark}
    Thanks to \cref{thm: kostant}, we will often abusively view the Kostant slice as a graded map $\kappa: \fr{t}\mmod W \cong e + \ld{\g}^f \subseteq \ld{\g}$, which provides a section to the quotient map $\ld{\g} \to \ld{\g}\mmod \ld{G} \cong \fr{t}\mmod W$.
\end{remark}
Motivated by \cref{thm: kostant}, we are led to the following definition.
\begin{definition}
    Let $\ld{C}$ denote the (graded) subgroup scheme of $\ld{G}(-2\rho) \times \ld{\g}(2 - 2\rho)$ consisting of those pairs $(g, x)$ such that $x$ is regular and $\Ad_g(x) = x$.
    Write $\ld{J}$ to denote the fiber product $\ld{C} \times_{\ld{\g}} (e + \ld{\g}^f)$. The projection $\ld{J} \to (e + \ld{\g}^f)$ equips $\ld{J}$ with the structure of a graded group scheme over $(e + \ld{\g}^f) \cong \ld{\g}(2)\mmod \ld{G}$; we will call $\ld{J}$ the \textit{group scheme of regular centralizers}. The fibers of the map $\ld{J} \to \ld{\g}(2)\mmod \ld{G}$ are commutative group schemes, so $\ld{J}$ is a commutative group scheme over $\ld{\g}(2)\mmod \ld{G}$.
\end{definition}
\begin{example}\label{ex: classical groups reg centralizer}
    In classical cases, the group scheme of regular centralizers can be described very explicitly. (The descriptions below are closely related to well-known descriptions of spectral curves in the Hitchin fibration for classical groups.)
    \begin{itemize}
        \item When $\ld{G} = \GL_n$, we may identify $\ld{\g}\mmod \ld{G} \cong \AA^n$ with coordinates $c_1, \cdots, c_n$. Since a matrix commutes with a regular matrix $x\in \gl_n$ if and only if it can be expressed as a polynomial in $x$, the Cayley-Hamilton theorem implies that $\ld{J}$ is the group scheme whose fiber over $(c_1, \cdots, c_n)$ is given by the group of units in $k[t]/(t^n + c_1 t^{n-1} + \cdots + c_n)$. Specifically, any such unit $f(t)$ defines an element of $\GL_n$ by evaluation on $\kappa(c_1, \cdots, c_n)$. Note that the fiber of $\ld{J}$ over $(c_1, \cdots, c_n) = (0, \cdots, 0)$ is the group of units in $k[t]/t^n$; this is isomorphic to $\GG_m \times \WW_{n-1}$, where $\WW_{n-1}$ is the (integral) Witt vectors of length $n-1$.
        The case of $\GL_n$ leads to explicit descriptions of $\ld{J}$ for other classical groups.
        \item If $\ld{G} = \Sp_{2n}$, we may identify $\ld{\g}\mmod \ld{G} \cong \AA^n$ with coordinates $p_1, \cdots, p_n$. Observe that $k[t]/(t^{2n} + p_1 t^{2n-2} + \cdots + p_n)$ admits the structure of a symplectic vector space: the symplectic pairing sends 
        $$(f,g) \mapsto \text{coefficient of }t^{2n-1}\text{ in }f(t)g(-t).$$
        Using this, one can view $\ld{J}$ as the group scheme whose fiber over $(p_1, \cdots, p_n)$ is the subgroup of those units $f(t)$ in $k[t]/(t^{2n} + p_1 t^{2n-2} + \cdots + p_n)$ such that $f(t)^{-1} = f(-t)$. Note that multiplication by $f(t)$ leaves the above symplectic pairing invariant.
        \item If $\ld{G} = \SO_{2n+1}$, we may identify $\ld{\g}\mmod \ld{G} \cong \AA^n$ with coordinates $p_1, \cdots, p_n$. Observe that $k[t]/(t^{2n+1} + p_1 t^{2n-1} + \cdots + p_n t)$ admits the structure of a quadratic vector space: the associated symmetric bilinear form pairing sends 
        $$(f,g) \mapsto \text{coefficient of }t^{2n}\text{ in }f(t)g(-t).$$
        Just as above, one can describe $\ld{J}$ as the group scheme whose fiber over $(p_1, \cdots, p_n)$ is the subgroup of those units $f(t)$ in $k[t]/(t^{2n+1} + p_1 t^{2n-1} + \cdots + p_n t)$ such that $f(t)^{-1} = f(-t)$.
        \item If $\ld{G} = \SO_{2n}$, we may identify $\ld{\g}\mmod \ld{G} \cong \AA^n$ with coordinates $p_1, \cdots, p_{n-1}, c_n$.
        Observe that $k[t, v]/(tv-c_n, t^{2n-2} + p_1 t^{2n-4} + \cdots + p_{n-1} + v^2)$ admits the structure of a quadratic vector space: the associated symmetric bilinear form sends
        $$(f,g) \mapsto \text{coefficient of }t^{2n-2}\text{ in }f(t,v)g(-t,-v).$$
        Using this, one can view $\ld{J}$ as the group scheme whose fiber over $(p_1, \cdots, p_{n-1}, c_n)$ is the subgroup of those units $f(t,v)$ in $k[t,v]/(tv - c_n, t^{2n-2} + p_1 t^{2n-4} + \cdots + p_{n-1} + v^2)$ such that $f(t,v)^{-1} = f(-t,-v)$.
    \end{itemize}
    It is also possible to use the description of the regular centralizer group scheme for $\GL_n$ to get at $\ld{J}_{\ld{G}}$ for exceptional $\ld{G}$. Here is a fun example:
    \begin{itemize}
        \item Suppose $\ld{G} = \G_2$, so that it is the subgroup of $\GL_7$ stabilizing a certain $3$-form $\phi$ on $\AA^7$ (this $3$-form was discovered by Engel and Reichel). If an element of $\G_2$ has eigenvalues $\lambda_1, \lambda_2$, we will set $\lambda_3 = -(\lambda_1 + \lambda_2)$; the corresponding matrix in $\GL_7$ has eigenvalues $0, \lambda_1, \lambda_2, \lambda_3, -\lambda_1, -\lambda_2, -\lambda_3$. Moreover, the scheme $\g_2^\ast\mmod \G_2$ can be described as $\spec k[c_2, c_6]$, where
        \begin{align*}
            c_2 & = \lambda_1^2 + \lambda_2^2 + \lambda_1 \lambda_2, \\
            c_6 & = \lambda_1^2 \lambda_2^2 (\lambda_1 + \lambda_2)^2.
        \end{align*}
        If $e_j$ denotes the $j$th elementary symmetric polynomial in the variables $0, \lambda_1, \lambda_2, \lambda_3, -\lambda_1, -\lambda_2, -\lambda_3$, then one directly calculates that
        $$e_j = \begin{cases}
            1 & j=0,\\
            -2c_2 & j=2,\\
            c_2^2 & j=4, \\
            -c_6 & j=6, \\
            0 & j\text{ odd, or }j\geq 8.
        \end{cases}$$
        Therefore, the characteristic polynomial of an element of $\G_2$ with eigenvalues $\lambda_1, \lambda_2$ is given by
        $$\chi(t) = t^7 - 2c_2 t^5 + c_2^2 t^3 - c_6 t.$$
        The ring $k[c_2, c_6, t]/\chi(t)$ is a free $k[c_2, c_6]$-module of rank $7$, and the Engel-Reichel $3$-form $\phi$ defines a $k[c_2, c_6]$-linear $3$-form on $k[c_2, c_6, t]/\chi(t)$. This $3$-form can be described explicitly, but its description is not very enlightening (at least, to me).
        One can describe $\ld{J}_{\G_2}$ as the group scheme over $\g_2^\ast\mmod \G_2$ whose fiber over $(c_2, c_6)$ is given by a certain subgroup of those units $f(t)\in (k[t]/\chi(t))^\times \subseteq \GL_7$ such that the automorphism of $k[t]/\chi(t)$ given by multiplication by $f(t)$ stabilizes the aforementioned $3$-form.
    \end{itemize}
\end{example}
\begin{theorem}[\cite{bfm, homology-langlands}]\label{thm: classical homology loops G}
    There is an isomorphism
    $$\spec \H^G_\ast(\Omega G; \QQ) \cong \ld{J}$$
    of commutative graded group schemes over $\ld{\g}(2)\mmod \ld{G}$.
\end{theorem}
\begin{remark}
    Suppose $G$ is simply-connected.
    The Borel-completion of $\H^G_\ast(\Omega G; \QQ)$ is simply $\pi_\ast \QQ[\Omega G]^{hG}$, which is Koszul dual to $\QQ[G]$ acting on $\QQ[\Omega G]$. According to rational homotopy theory, this can in turn be viewed as a deformation of $\pi_\ast \QQ[\Omega G]$ along $\{0\} \hookrightarrow \spec \pi_\ast \QQ^{hG} \cong \widehat{\ld{\fr{t}}(2)}\mmod W$. But $\pi_\ast \QQ[\Omega G]$ can be identified with the universal enveloping algebra $U(\pi_\ast(\Omega G)_\QQ)$, which is isomorphic to the dual Lie algebra $(\ld{\g}^e)^\ast$ by \cref{thm: classical homology loops G} (or \cref{prop: ginzburg coh}). \cref{thm: classical homology loops G} implies that the relevant deformation of $\pi_\ast \QQ[\Omega G]$ along $\{0\} \hookrightarrow \widehat{\ld{\fr{t}}(2)}\mmod W$ is precisely the relative enveloping algebra of Lie algebra of (the completion of) $\ld{J}$.
\end{remark}
\begin{corollary}[\cite{ngo-invent}]\label{cor: classifying stack of J}
    The classifying stack $B_{\ld{\g}(2)\mmod \ld{G}} \ld{J}$ is isomorphic to the stacky quotient $\ld{\g}(2 - 2\rho)^\reg/\ld{G}(-2\rho)$. In particular, there is a graded isomorphism
    $$(\ld{\g}(2)\mmod \ld{G} \times \ld{G}(-2\rho))/\ld{J} \cong \ld{\g}(2 - 2\rho)^\reg.$$
\end{corollary}
\begin{proof}
    For notational simplicity, let us drop the gradings.
    By definition of $\ld{J}$, the classifying stack $B_{\ld{\g}(2)\mmod \ld{G}} \ld{J}$ is a $1$-stack whose objects given by the $\ld{G}$-orbit of the Kostant slice $e + \ld{\g}^f \subseteq \ld{\g}$, and such that an isomorphism $x \xar{\sim} y$ is given by an element $g\in \ld{G}$ such that $\Ad_g(x) = y$. By \cref{thm: kostant}, the $\ld{G}$-orbit of the Kostant slice is precisely $\ld{\g}^\reg$, so we obtain the desired result.
\end{proof}

We will also need to use the following general result, discussed in \cite[Section 6.5]{bf-derived-satake}; this result relies on several results from \cite[Section 6.1]{bbd-perverse}.
\begin{setup}\label{setup: changing coefficients and bases}
    Let $F$ be a number field, and let $R$ be a localization of the ring of integers $\co_F$. Let $\ell$ be a prime which is invertible in $R$, and fix an isomorphism $\cc \cong \ol{\QQ_\ell}$. Let $X$ be an $R$-scheme of finite type equipped with an action of an affine group scheme $H$ over $R$, such that the set of orbits is finite\footnote{In \cite{bf-derived-satake}, the group scheme $H$ is assumed to be smooth, but this is not necessary.}.
    If $R'$ is a ring equipped with a ring map $R \to R'$, let $X_{R'}$ denote the base-change of $X$ to $R'$ (and similarly for $H_{R'}$).
\end{setup}
\begin{theorem}[{\cite[Proposition 5]{bf-derived-satake}}]\label{thm: finite type changing coefficients and bases}
    In \cref{setup: changing coefficients and bases}, there is a localization $R \subseteq R'$ such that for any $k$-point $R' \twoheadrightarrow k$ with $k$ being a finite field, there are equivalences (which are functorial in $X$)
    \begin{align*}
        \Shv^{c,\et}_{H_{\ol{\FF}_q}}(X_{\ol{\FF}_q}; \ol{\QQ_\ell}) & \xar{\sim} \Shv^{c,\et}_{H_{\ol{F}}}(X_{\ol{F}}; \ol{\QQ_\ell}) \xar{\sim} \Shv^{c,\et}_{H_\cc}(X_\cc; \ol{\QQ_\ell}) \\
        & \xar{\sim} \Shv^c_{H(\cc)}(X(\cc); \ol{\QQ_\ell}) \xar{\sim} \Shv^{c}_{H(\cc)}(X(\cc); \cc).
    \end{align*}
\end{theorem}
\begin{remark}\label{rmk: coeff change infty cat}
    We are implicitly stating in \cref{thm: finite type changing coefficients and bases} that the displayed equivalences are of presentable symmetric monoidal stable $\infty$-categories. 
    However, we will not prove this here.
\end{remark}
We will apply \cref{thm: finite type changing coefficients and bases} to the case when $X$ is an ind-finite $R$-scheme (namely, the affine Grassmannian of $G$): see \cref{thm: changing coefficients and bases}. In order to use \cref{thm: finite type changing coefficients and bases} in this manner, we need the following result.
\begin{theorem}[Quillen, Garland-Raghunathan, \cite{garland-raghunathan, mitchell-buildings}]\label{thm: quillen}
    Let $G_c$ be a compact Lie group, and let $G$ denote the associated reductive algebraic group over $\cc$. Then there is a homotopy equivalence $G(\cc\ls{t})/G(\cc\pw{t}) \simeq \Omega G_c$ which is equivariant for the left-action of $G_c \subseteq G(\cc) \subseteq G(\cc\pw{t})$ on the left-hand side and the action of $G_c$ on the right-hand side given by conjugation.
\end{theorem}

The final input we need to prove \cref{thm: derived satake} is the following elementary observation.
\begin{observe}\label{obs: formal defn}
    Let $A$ be an $\E{n}$-$\QQ$-algebra spectrum with $n\geq 1$. View $\pi_\bull A$ as a \textit{graded} $\QQ$-algebra spectrum (where $\pi_j A$ is placed in degree zero and weight $j$). Then, its half-shearing $\sh^{1/2}(\pi_\bull A)$ has underlying $\QQ$-module spectrum $\bigoplus_{j\in \Z} \pi_j(A)[j]$. Note that \cref{lem: no e2 half shearing} and \cref{rmk: Q does not save shearing} say that $\sh^{1/2}(\pi_\bull A)$ will \textit{a priori} only admit an $\E{1}$-algebra structure in graded $\QQ$-modules; but if $A$ only has \textit{even} homotopy groups, \cref{rmk: half shear on evenly graded} says that $\sh^{1/2} \pi_\bull A$ will indeed admit an $\Eoo$-$\QQ$-algebra structure in graded $\QQ$-modules.
    If $A$ is formal, then there is an equivalence $A \simeq \sh^{1/2}(\pi_\bull A)$ of $\E{1}$-$\QQ$-algebra spectra. (Here, we have implicitly applied the forgetful functor $\Mod_\QQ^\gr \to \Mod_\QQ$ to the right-hand side.)
\end{observe}

Let us now prove \cref{thm: derived satake}, ignoring any questions of $\E{3}$-monoidality.
\begin{proof}[Proof of \cref{thm: derived satake}]
    For notational simplicity, let us write $\cC$ to denote $\Shv_{G\pw{t}}^{c,\Sat}(\Gr_G; \QQ)$. Let $\tilde{\cC}$ denote the base-change $\cC \otimes_{\Rep(\ld{G})} \Mod_\QQ$, so that $\IC_0$ is a compact generator of $\tilde{\cC}$ (by definition of $\cC$). It follows from the Barr-Beck theorem \cite[Theorem 4.7.3.5]{HA} that there is an equivalence $\Phi: \tilde{\cC} \xar{\sim} \Perf_{\End_{\tilde{\cC}}(\IC_0)}$, implemented by the functor $\Hom_{\tilde{\cC}}(\IC_0, -)$.
    Recall that $\Shv_{G\pw{t}}^c(\Gr_G; \QQ)$ is equipped with an $\E{3}$-monoidal structure (see, e.g., \cite{nocera-e3}) where $\IC_0$ is the monoidal unit, which equips $\End_{\tilde{\cC}}(\IC_0)$ with an $\E{3}$-monoidal structure. Note that the $\E{3}$-monoidal structure on $\Shv_{G\pw{t}}^c(\Gr_G; \QQ)$ induces one on $\cC$. (For simplicity, we will ignore questions of $\E{3}$-monoidality in the argument below.)
    Write $\ld{\cR}\in \Shv_{G\pw{t}}^{c,\Sat}(\Gr_G; \QQ)$ denote the (perverse) sheaf obtained by the action of $\co_{\ld{G}}\in \Rep(\ld{G})$ on $\IC_0$. By definition of $\tilde{\cC}$, we can identify $\End_{\tilde{\cC}}(\IC_0) \simeq \Hom_{\cC}(\IC_0, \IC_0 \star \ld{\cR})$.

    A key claim proved in \cite{bf-derived-satake} is that $\End_{\tilde{\cC}}(\IC_0)$ is in fact formal. In other words, by \cref{obs: formal defn}, there is an isomorphism 
    \begin{equation}\label{eq: displayed ext}
        \End_{\tilde{\cC}}(\IC_0) \simeq \Hom_{\cC}(\IC_0, \IC_0 \star \ld{\cR}) \cong \sh^{1/2}(\Ext^\bull_{\Shv_{G\pw{t}}^c(\Gr_G; \QQ)}(\IC_0, \IC_0 \star \ld{\cR})).
    \end{equation}
    We will only give a sketch of this below (but see \cref{rmk: e2-formality in satake} for an alternative argument).
    Note that the claimed formality of the $\QQ$-algebra $\End_{\tilde{\cC}}(\IC_0)$ can be proved after base-changing to $\cc$.
    In this case, we can identify
    $$\End_{\tilde{\cC}}(\IC_0) \otimes_\QQ \cc \cong \End_{\tilde{\cC} \otimes_\QQ \cc}(\IC_0).$$
    Let us write $\cd$ to denote $\Shv_{G_{\ol{\FF}_q}}^{c,\et}(\Gr_{G,_{\ol{\FF}_q}}; \ol{\QQ_\ell})$, and $\tilde{\cd}$ to denote the base-change of the $\Rep(\ld{G}_{\ol{\QQ_\ell}})$-module $\cd$ to $\Mod_{\ol{\QQ_\ell}}$.
    Applying \cref{thm: changing coefficients and bases}, we can identify
    $$\End_{\tilde{\cC} \otimes_\QQ \cc}(\IC_0) \simeq \End_{\tilde{\cd}}(\IC_0) \simeq \Hom_{\Shv_{G_{\ol{\FF}_q}\pw{t}}^{c, \et}(\Gr_{G,_{\ol{\FF}_q}}; \ol{\QQ_\ell})}(\IC_0, \IC_0 \star \ld{\cR}).$$
    As discussed in the next paragraph, $\IC_0 \star \ld{\cR}$ is a direct sum over $\lambda \in \Lambda^+$ of a finite number of copies of $\IC_\lambda[\pdb{2\rho, \lambda}](\pdb{\ld{\rho}, \lambda})$. Using \cite[Example 3.1.4]{bezr-yun-koszul}, one finds that both $\IC_0$ and $\IC_0 \star \ld{\cR}$ are pure of weight zero. Therefore, \cite[Lemma 3.1.5]{bezr-yun-koszul} implies that Frobenius acts on $\pi_{-j} \Hom_{\Shv_{G_{\ol{\FF}_q}\pw{t}}^{c,\et}(\Gr_{G,_{\ol{\FF}_q}}; \ol{\QQ_\ell})}(\IC_0, \IC_0 \star \ld{\cR})$ by multiplication\footnote{\textit{A priori}, one needs to fix a square root $q^{1/2}$ of $q$ to state this result. However, as shown below, $\End_{\tilde{\cC}}(\IC_0)$, and hence $\End_{\tilde{\cd}}(\IC_0)$, is concentrated entirely in even degrees. Therefore, the integer $j$ can be assumed to be even without any loss of generality, and no choice of a square root of $q$ is needed.} by $q^{j/2}$. Since the action of Frobenius respects the ring structure on $\End_{\tilde{\cd}}(\IC_0)$, it follows that the action of Frobenius splits the Postnikov filtration on $\End_{\tilde{\cd}}(\IC_0)$ as a $\ol{\QQ_\ell}$-algebra. In particular, $\End_{\tilde{\cd}}(\IC_0)$ is formal as a $\ol{\QQ_\ell}$-algebra, as desired.

    In order to compute the $\Ext$-algebra on the right-hand side of \cref{eq: displayed ext}, we will apply \cref{thm: full faithful}. It is easy to see that $\IC_0$ is even (it is supported only on the basepoint of $\Omega G$). The fact that $\IC_0 \star \ld{\cR}$ is even is a consequence of the Peter-Weyl theorem and the proof of \cref{thm: abelian satake}. Namely, since $\co_{\ld{G}} = \bigoplus_{\lambda\in \Lambda^+} \End(V_\lambda) \cong \bigoplus_{\lambda \in \Lambda^+} V_\lambda \otimes V_\lambda^\ast$, one can identify $\IC_0 \star \ld{\cR}$ with the direct sum $\bigoplus_{\lambda \in \Lambda^+} \IC_\lambda[\pdb{2\ld{\rho}, \lambda}] \otimes V_\lambda^\ast$. 
    However, if $\mu \leq \lambda$ and $j_\mu: \Gr_G^\mu \hookrightarrow \Gr_G^{\leq \lambda}$ is the inclusion, $j_\mu^\ast \IC_\lambda \cong j_\mu^! \IC_\lambda \cong \ul{\QQ}_{\Gr_G^\mu}[2\dim \Gr_G^\mu]$.
    In particular, each $\IC_\lambda$ is even in the sense of \cref{thm: full faithful}.

    Applying \cref{thm: full faithful}, we conclude that there is a graded isomorphism
    \begin{multline*}
        \Ext^\bull_{\Shv_{G\pw{t}}^c(\Gr_G; \QQ)}(\IC_0, \IC_0 \star \ld{\cR}) \cong \Hom^\bull_{\H^\ast_G(\Omega G; \QQ)}(\H^\ast_G(\Omega G; \IC_0), \H^\ast_G(\Omega G; \IC_0 \star \ld{\cR})) \\
        \cong \Hom^\bull_{\H^\ast_G(\ast; \QQ)}(\H^\ast_G(\Omega G; \IC_0), \H^\ast_G(\Omega G; \IC_0 \star \ld{\cR}))^{\spec \H^G_\ast(\Omega G; \QQ)}.
    \end{multline*}
    But $\H^\ast_G(\Omega G; \IC_0) \cong \H^\ast_G(\ast; \QQ) \cong \co_{\fr{t}(2)\mmod W}$, and $\H^\ast_G(\Omega G; \IC_0 \star \ld{\cR}) \cong \co_{\fr{t}(2)\mmod W \times \ld{G}}$. Note that here and below, the symbol $\co$ denotes the \textit{classical} (and not derived) ring of functions; in the two cases above, this distinction does not matter, but it will momentarily.
    This, along with \cref{thm: classical homology loops G}, implies that 
    \begin{multline*}
        \Ext^\bull_{\Shv_{G\pw{t}}^c(\Gr_G; \QQ)}(\IC_0, \IC_0 \star \ld{\cR}) \cong \Hom^\bull_{\H^\ast_G(\ast; \QQ)}(\H^\ast_G(\Omega G; \IC_0), \H^\ast_G(\Omega G; \IC_0 \star \ld{\cR}))^{\spec \H^G_\ast(\Omega G; \QQ)} \\
        \cong \Hom^\bull_{\co_{\fr{t}(2)\mmod W}}(\co_{\fr{t}(2)\mmod W}, \co_{\fr{t}(2)\mmod W \times \ld{G}(-2\rho) })^{\ld{J}} 
        \cong \co_{\fr{t}(2)\mmod W \times \ld{G}}^{\ld{J}(2)} \cong \co_{(\fr{t}(2)\mmod W \times \ld{G}(-2\rho) )/\ld{J}}.
    \end{multline*}
    \cref{cor: classifying stack of J} and the Chevalley restriction theorem precisely identifies this with $\co_{\ld{\g}^\ast(2 - 2\rho)^\reg}$. Since this is the classical (and not derived) ring of functions\footnote{We emphasize that the desired result would be \textit{false} if one was to instead take {derived} rings of functions. Indeed, it is not true that the inclusion $\AA^2 - \{0\} \subseteq \AA^2$ induces an isomorphism on derived global sections of the structure sheaf (namely, $\pi_{-1} \co_{\AA^2 - \{0\}} \cong \QQ[x,y]/(x^\infty, y^\infty)$), although it is certainly true at the level of classical rings of functions.}, and the complement of $\ld{\g}^\ast(2 - 2\rho)^\reg \subseteq \ld{\g}^\ast(2 - 2\rho)$ has codimension $\geq 2$, there is an isomorphism $\co_{\ld{\g}^\ast(2 - 2\rho)^\reg} \cong \co_{\ld{\g}^\ast(2 - 2\rho)}$ by the algebraic Hartogs lemma. We conclude that $\Ext^\bull_{\Shv_{G\pw{t}}^c(\Gr_G; \QQ)}(\IC_0, \IC_0 \star \ld{\cR}) \cong \co_{\ld{\g}^\ast(2 - 2\rho)}$. This implies that there is an isomorphism
    $$\End_{\tilde{\cC}}(\IC_0) \cong \sh^{1/2}(\Ext^\bull_{\Shv_{G\pw{t}}^c(\Gr_G; \QQ)}(\IC_0, \IC_0 \star \ld{\cR})) \cong \sh^{1/2}(\co_{\ld{\g}^\ast(2 - 2\rho)}) \cong \co_{\ld{\g}^\ast[2 - 2\rho]}.$$
    This implies that $\tilde{\cC} \simeq \Perf(\ld{\g}^\ast[2 - 2\rho])$, and so $\cC \simeq \Perf(\ld{\g}^\ast[2 - 2\rho]/\ld{G}[-2\rho])$, as desired.
\end{proof}
\begin{remark}\label{rmk: e2-formality in satake}
    The proof of \cref{eq: displayed ext} using \cref{thm: changing coefficients and bases} can be circumvented using \cref{lem: formality polynomial}: namely, $\End^\bull_{\tilde{\cC}}(\IC_0) \simeq \Ext^\bull_{\Shv_{G\pw{t}}^c(\Gr_G; \QQ)}(\IC_0, \IC_0 \star \ld{\cR})$ is isomorphic to the finitely generated polynomial algebra $\co_{\ld{\g}^\ast(2)}$. One knows that $\End^\bull_{\tilde{\cC}}(\IC_0)$ can be made into an $\E{3}$-algebra (see e.g., the discussion in \cite[Proposition 16.2.8]{bzsv}), and hence \cref{lem: formality polynomial} shows that $\End_{\tilde{\cC}}(\IC_0)$ is formal as an $\E{2}$-$\QQ$-algebra. 
    However, lifting this to formality as an $\E{3}$-algebra is rather nontrivial. Indeed, ring spectra with polynomial homotopy need \textit{not} be formal as $\E{3}$-algebras (as mentioned in \cref{rmk: shearing not e3}).
    
    However (at least for $\QQ$-coefficients), the argument using \cref{thm: changing coefficients and bases} is more general, since it does not assume that the relevant $\Ext$-algebra is polynomial. This argument will be useful as a model for \cref{thm: ordinary homology criterion satake} below.
\end{remark}
\begin{remark}
    It is a consequence of the above proof that under the equivalence of \cref{thm: derived satake}, the following diagram commutes:
    \begin{equation}\label{eq: compatibility with kostant slice}
        \xymatrix{
        \Shv_{G\pw{t}}^{c, \Sat}(\Gr_G; \QQ) \ar[rr]^-\sim_-{\text{\cref{thm: derived satake}}} \ar[d]_-{\Gamma_{G\pw{t}}(\Gr_G; -)} & & \Perf(\ld{\g}[2-2\rho]/\ld{G}[-2\rho]) \ar[d]^-{\kappa^\ast} \\
        \Shv_G^c(\ast; \QQ) \ar[rr]^-\sim_-{\text{\cref{lem: fdtl class affine space}(a)}} & & \Perf(\fr{t}[2] \mmod W).
        }
    \end{equation}
    Here, the map $\kappa$ is the map of ungraded derived schemes underlying the half-shearing of the Kostant slice $\kappa: \fr{t}(2) \mmod W \to \ld{\g}(2 - 2\rho)^\reg/\ld{G}(-2\rho)$. Note that since these two schemes have even weights, \cref{rmk: half shear on evenly graded} ensures that the half-shearing of $\kappa$ will indeed be a map of $\Eoo$-$\QQ$-schemes.
\end{remark}

The proof of the main result of this article will follow the same outline. We will still need \cref{thm: abelian satake} as input; but the heart of our work lies in proving an analogue of \cref{thm: classical homology loops G}.

%% file: equivalences/spherical-varieties.tex
\subsection{Spherical varieties}

In this section, we will review some of the theory of spherical varieties.
Since the examples we will study in this article are rather simple (from the perspective of representation theory), we do not, strictly speaking, need the general theory. However, the recollections of this section will nevertheless be useful in placing basic phenomena that we will observe later into a broader context (see \cref{subsec: stating the bzsv conjecture}). 

We will not give any proofs in this section, but instead refer to \cite{brion-luna-vust, luna-vust, timashev-homogeneous-spaces, bzsv} for details; in particular, this section is \textit{not} intended to be an introduction to the theory of spherical varieties or to the theory of their Hamiltonian duals. (Instead, the reader should see \cite{perrin-survey-spherical} for a very readable introduction to spherical varieties.) The base field in this section will always be the complex numbers, $G$ will always be a connected reductive algebraic group over $\cc$, $B\subseteq G$ will denote a chosen Borel subgroup, and $N$ will be its unipotent radical.
\begin{definition}
    A subgroup $H\subseteq G$ is called \textit{spherical} if any of the following equivalent conditions are satisfied:
    \begin{enumerate}
        \item For any $G$-variety $X$ and any $H$-fixed point $x\in X$, the closure $\ol{G\cdot x}$ contains finitely many $G$-orbits.
        \item There are finitely many $H$-orbits in the flag variety $G/B$ of $G$.
        \item There is an open $H$-orbit in $G/B$.
        \item The action of $B$ on $G/H$ has an open dense orbit.
    \end{enumerate}
    An irreducible $G$-variety $X$ is called \textit{spherical} if it is normal and admits a dense open $B$-orbit $\punc{X}\subseteq X$. In this case, $X$ also contains an open $G$-orbit given by $G\cdot \punc{X}$. If $x\in \punc{X}$ and $H$ is its stabilizer, there is an isomorphism $\punc{X} = G/H$, and $H$ is a spherical subgroup of $G$.
\end{definition}
Before delving into examples, let us mention that the condition of being a spherical $G$-variety is relevant for our purposes because of the following result:
\begin{theorem}[{\cite[Theorem 3.2.1]{gaitsgory-nadler}}]\label{thm: spherical gaitsgory nadler}
    Let $H\subseteq G$ be a subgroup. Then the following conditions are equivalent:
    \begin{enumerate}
        \item $G/H$ is a spherical $G$-variety.
        \item The group $H(\cc\ls{t})$ acts on $\Gr_G(\cc) = G(\cc\ls{t})/G(\cc\pw{t})$ with countably many orbits.
        \item The group $G(\cc\pw{t})$ acts on $(G/H)(\cc\ls{t})$ with countably many orbits.
    \end{enumerate}
\end{theorem}
\begin{remark}
    We refer the reader to \cite{gaitsgory-nadler} for a proof of \cref{thm: spherical gaitsgory nadler}, but since the argument is so short, let us recall why (b) implies (a). Suppose $\lambda: \GG_m \to G$ is a subgroup, so that we obtain a point $x_\lambda\in \Gr_G(\cc)$. Then the $G$-orbit $X_\lambda = G \cdot x_\lambda\subseteq \Gr_G$ is a flag variety of $G$, and by (b), the number of $H(\cc\ls{t})$-orbits intersecting $X_\lambda$ is countable. This implies that there is an $H(\cc\ls{t})$-orbit which intersects $X_\lambda$ in an open set. If we choose a point $y\in X_\lambda$ in this open set, this implies that there is a surjection $\fr{h} \twoheadrightarrow T_y X_\lambda$. If $\fr{p}_y$ is the Lie algebra of the parabolic subgroup of $G$ stabilizing $y$, the tangent space $T_y X_\lambda$ can be identified with $\g/\fr{p}_y$. In particular, if we choose $\lambda$ to be regular, $\fr{p}_y$ is isomorphic to a Borel subalgebra $\fr{b}\subseteq \g$, and hence there is a surjection $\fr{h} \twoheadrightarrow \g/\fr{b}$. But this implies that $H$ has an open orbit in $G/B$, so $H$ is spherical.
\end{remark}

There are a lot of examples of spherical varieties: it includes the class of flag varieties, symmetric spaces (essentially by the Iwasawa decomposition), and toric varieties.
\begin{example}
    The quotient $\GL_n/\GL_{n-1}$ is an affine spherical $\GL_n$-variety; it is isomorphic to the variety $\{(x,V) \in \cc^{n+1} \times \Gr_n(\cc^{n+1}) | x\not\in V\}$. The fact that the $\cc$-points of $\GL_n/\GL_{n-1}$ is homotopy equivalent to $S^{2n-1}$ motivates the terminology ``spherical''.
\end{example}
\begin{example}\label{ex: gp case}
    As mentioned above, any symmetric space is a spherical variety. In particular, since $G$ is the fixed points of the involution on $G\times G$ which swaps the two factors, we see that $G \cong (G \times G)/G^\Delta$ is a spherical $G\times G$-variety. This will often be called the \textit{group case}.
\end{example}
\begin{example}\label{ex: spherical subgroups of PGL2}
    Suppose $G = \PGL_2$. Since the flag variety of $G$ is isomorphic to $\PP^1$, a subgroup $H\subseteq \PGL_2$ is spherical if and only if it has an open orbit in $\PP^1$. This is equivalent to saying that it is a subgroup of positive dimension. It is not difficult to see that all positive-dimensional subgroups of $\PGL_2$ can be conjugated either to $\PGL_2$ itself, the diagonal torus $\GG_m\subseteq \PGL_2$, its normalizer $\N_{\PGL_2}(\GG_m) \cong \mathrm{PO}_2 \subseteq \PGL_2$, or $S\cdot N\subseteq \PGL_2$, where $N$ is the strictly upper-triangular matrices and $S\subseteq \GG_m$. In general, a spherical subgroup $H\subseteq G$ is called \textit{horospherical} if $H$ contains the unipotent radical of the Borel $B \subseteq G$; the motivation for this term being, of course, that horocycles in $\SL_2(\RR)/\SO_2$ are orbits of the subgroup of strictly upper-triangular matrices in $\SL_2(\RR)$. These kinds of spherical varieties are \textit{not} considered in the present article.
\end{example}
\begin{warning}\label{warning: G/T not spherical}
    If $G$ is a semisimple algebraic group and $T\subseteq G$ is a maximal torus, the quotient $G/T$ is generally \textit{not} a spherical $G$-variety. Indeed, there generally will not be an open dense $T$-orbit in $G/B$, since $|\Phi^-|$ is often larger than $\mathrm{rank}(T)$, where $\Phi^-$ is the set of negative roots of $G$. For instance, although the quotient $\SL_2/\GG_m$ is a spherical $\SL_2$-variety, the quotient $\SL_3/T$ is not a spherical $\SL_3$-variety.
\end{warning}
\begin{remark}
    There is a finite list of closed connected spherical subgroups of simple algebraic groups: see \cite{knop-rohrle, kramer-spherical-subgroups}.
\end{remark}
\begin{example}\label{ex: toric}
    Let $G$ be a torus $T$. Then a $T$-variety $X$ is spherical if it is normal and contains a dense orbit, and hence is precisely an affine toric variety. Let $\Lambda$ denote the monoid of weights of $T$. Note that $\co_X$ is a $T$-submodule of $\co_T$, and so $\co_X = \bigoplus_{\lambda \in S_X} \cc_\lambda$ for some subset $S_X \subseteq \Lambda$. A standard fact from the theory of affine toric varieties is that a subset $S_X \subseteq \Lambda$ arises from an affine toric variety if and only if $S_X = C \cap \Lambda$ for some convex cone $C \subseteq \Lambda_\RR$ generated by finitely many elements of $\Lambda$ which span $\Lambda_\RR$. Equivalently, if $\ld{C} \subseteq \ld{\Lambda}$ denotes the dual cone, one observes that $C$ spans $\Lambda_\RR$ if and only if $\ld{C}$ is strictly convex (i.e., contains no line). Therefore, affine toric varieties are classified by strictly convex rational polyhedral cones of $\ld{\Lambda}_\RR$.
\end{example}
\cref{ex: toric} is the first indication that certain spherical varieties admit interesting combinatorial data. In particular, this combinatorial data will be useful in defining the \textit{Langlands dual group} to a spherical variety. We will recall some generalities on defining this dual group below, and then explain its manifestation in examples.

To define this dual group following \cite{sakellaridis-venkatesh}, let us now suppose that $X$ is a homogeneous quasi-affine spherical $G$-variety. In this case, if $\punc{X}\subseteq X$ is the open $B$-orbit, we will write $H$ to be the stabilizer of a point $\punc{X}(\cc)$, so that $X = G/H$ and $B\cdot H\subseteq G$ is open.

\begin{construction}\label{cstr: weight lattice for spherical}
    Let $\Frac(\co_X)$ denote the fraction field of $\co_X$, and let $\Frac(\co_X)^{(B)}$ denote the subset of $\Frac(\co_X) - \{0\}$ consisting of the nonzero rational $B$-eigenfunctions. Then the lattice $\dX_X$ is simply the group of $B$-eigencharacters, and there is an exact sequence
    $$1 \to \cc^\times \to \Frac(\co_X)^{(B)} \to \dX_X \to 1;$$
    in other words, for a fixed $\lambda \in \dX_X$, the functions $f\in \Frac(\co_X)^{(B)}$ which are $\chi$-eigenvectors are all proportional by a scalar in $\cc^\times$ (this follows from $X$ being spherical).
    Let $\Lambda_X$ denote the dual lattice to $\dX_X$. Then $\Lambda_X$ defines a torus $T_X$, and we will write $\fr{t}_X$ to denote $\Lambda_X \otimes \QQ$.
    The rank of the lattice $\Lambda_X$ (which is also the rank of $\dX_X$) is called the \textit{rank} of $X$.
\end{construction}
\begin{remark}
    Suppose $X = G/H$ is a homogeneous quasi-affine $G$-variety, and let $\dX_X = \Frac(\co_X)^{(B)}/\cc^\times$ as above. It is not difficult to see that $X$ is spherical if and only if $\dX_X$ is a lattice of finite rank.
    If $K$ is a maximal compact subgroup of $G(\cc)$, \cite{akhiezer-rank-topology} shows that 
    $$\rank(X) = \dim(K\backslash X(\cc)).$$
    This is a purely topological description of the rank of $X$.
\end{remark}

\begin{construction}\label{cstr: associated parabolic}
    The stabilizer of the open $B$-orbit $\punc{X}\subseteq X$ is a parabolic subgroup $P(X)$. We will write $L(X)$ to denote the Levi quotient of $P(X)$; it will often be viewed as a subgroup of $P(X)$ when convenient.
    Let $T$ be a maximal torus of $B\cap L(X)$; then the torus $T_X$ from above can be identified with $T/(T\cap B)$. The $T_X$-orbit of a point in the open $B$-orbit $\punc{X}(\cc)$ defines an embedding $T_X\hookrightarrow \punc{X}(\cc)$. In other words, the $B$-action on $\punc{X}$ defines a $T$-action on $\punc{X}\mmod N = \spec \co_X^N$, and this $T$-action factors through the quotient $T \twoheadrightarrow T_X$.
\end{construction}
\begin{remark}
    In \cite[Lemma 3.1]{knop-asymptotic-invariants}, Knop showed that if $X$ is quasi-affine, the set of coroots in the span of $\Delta_{L(X)}$ in $\Lambda$ is precisely the set of coroots $\ld{\alpha} \in \Lambda$ which are perpendicular to $\Lambda_X$.
\end{remark}

\begin{construction}\label{cstr: g-invt valuations and little weyl}
    Suppose $v: \Frac(\co_X)^\times \to \QQ$ is a discrete valuation which is trivial on $\cc^\times$. Then the restriction of $v$ to $\Frac(\co_X)^{(B)}$ defines a homomorphism $\Lambda_X \to \QQ$, i.e., a point of $\fr{t}_X$. It is known that the map from $G$-invariant valuations to $\fr{t}_X$ is an injection, and so we will write $\cV\subseteq \fr{t}_X$ to denote the subspace of $G$-invariant valuations. Let $\ld{\Lambda}_X^+$ denote the intersection $\Lambda_X \cap \cV$ of $G$-invariant $\Z$-valued valuations.

    It turns out that the subset $\cV \subseteq \fr{t}_X$ is a fundamental domain for the Weyl group $W_X$ of a root system in $\Lambda$ (where the weight lattice is $\Lambda_X$). In other words, the reflections over faces of $\cV$ of codimension $1$ generate a finite reflection subgroup $W_X \subseteq \GL(\fr{t}_X)$, and this Weyl group $W_X$ is called the \textit{little Weyl group} of $X$. One can canonically identify $W_X$ with a subgroup of $W$ which normalizes the Weyl group $W_{L(X)}$ of $L(X)$ (with respect to the chosen torus $T$).
\end{construction}
\begin{remark}\label{rmk: knop microlocal weyl group}
    The definition of the little Weyl group given above does not immediately relate to the microlocal nature of $X$. In \cite{knop-asymptotic-invariants, knop-set-of-orbits}, Knop gave an alternative construction of $W_X$ using the Hamiltonian $G$-action on $T^\ast X$. Very briefly, let us review this construction. 
    The quotient map $\g \to \g/\fr{h}$ defines an inclusion $(\g/\fr{h})^\ast \subseteq \g^\ast$, and we will denote this by $\fr{h}^\perp$ (it can be viewed as a subspace of $\g$ via the isomorphism $\g^\ast \cong \g$ given by the Killing form). Consider the moment map $\mu: T^\ast X \cong (G \times \fr{h}^\perp)/H \to \g^\ast$ of the Hamiltonian $G$-action on $T^\ast X$. Composing with the characteristic polynomial map $\g^\ast \to \g^\ast \mmod G \cong \fr{t}^\ast\mmod W$ defines a map $T^\ast X \to \fr{t}^\ast \mmod W$. Observe also that the quotient map $T \twoheadrightarrow T_X$ induces an inclusion $\fr{t}_X^\ast \hookrightarrow \fr{t}^\ast$.

    Fix a character $\chi: T_X \to \GG_m$. Then, there is a $(P(X), \chi)$-eigenfunction $f_\chi \in \co_{\punc{X}}$ (unique up to scalar multiplication) defines a section $d\log(f_\chi): \punc{X} \to T^\ast \punc{X}$. This section is independent of the choice of $f_\chi$, since $f_\chi$ is unique up to scalar multiplication. Ranging over all characters $\chi$, one obtains a map $\fr{t}_X^\ast \times \punc{X} \to T^\ast \punc{X}$. If $\cP$ denotes the set of conjugates of the parabolic subgroup $P(X)$, we further obtain a map $\fr{t}_X^\ast \times (\cP \times \punc{X}) \to T^\ast X$. 
    Knop showed that the image of this map is dense, and that there is an isomorphism $(T^\ast X)\mmod G \cong \fr{t}_X^\ast \mmod W_X$.
    Said slightly differently, the fiber product $T^\ast X \times_{\fr{t}^\ast \mmod W} \fr{t}^\ast$ generally has multiple irreducible components. If $C$ is an irreducible component which dominates $T^\ast X$, we obtain a covering $C \to T^\ast X$, and $W_X$ is the Galois group of this covering. In particular, note that this construction describes $W_X$ as a subquotient of $W$. (However, there is in fact a canonical embedding $W_X \hookrightarrow W$.)
    
    In \cite{knop-set-of-orbits}, Knop reinterpreted the above construction as follows: if $\co_B(X)$ is the set of $B$-orbits in $X$, Knop constructed an action of $W$ on $\co_B(X)$. There is a canonical bijection between $\co_B(X)$ and the set of irreducible components of $T^\ast X \times_{\g^\ast} \tilde{\g}$ (given by taking the conormal bundle), where $\tilde{\g}$ is the Grothendieck-Springer resolution. The action of $W$ on $\co_B(X)$ can be understood as arising from the action of the Steinberg scheme $\tilde{\g} \times_{\g^\ast} \tilde{\g}$ by convolution and the isomorphism of \cite[Theorem 3.4.1]{chriss-ginzburg}. In any case, the stabilizer of the open $B$-orbit is isomorphic to $W_X \ltimes W_{L(X)}$. A related result was proved in \cite{ressayre-knop}: namely, if $H\subseteq G$ is a reductive spherical subgroup and $X = G/H$, the Weyl group of $H$ can be recovered as the stabilizer inside $W$ of a(ny) minimal rank $B$-orbit on $X$ viewed as an element of $\co_B(X)$.
\end{remark}
\begin{remark}\label{rmk: gaitsgory nadler orbits}
    Continuing \cref{thm: spherical gaitsgory nadler}, one can show (see \cite[Proposition 4.10]{luna-vust} or \cite[Theorem 3.2.1]{gaitsgory-nadler}) that the $G(\cc\pw{t})$-orbits on $(G/H)(\cc\ls{t})$ are in bijection with $H(\cc\ls{t})$-orbits on $\Gr_G(\cc)$, which in turn are in bijection with $\ld{\Lambda}_X/W_X \cong \ld{\Lambda}_X^+$. This generalizes the Cartan decomposition, in the sense that when applied to the group case of \cref{ex: gp case}, it recovers the standard parametrization of the $G(\cc\pw{t})$-orbits on $\Gr_G$. The bijection between $G(\cc\pw{t})$-orbits on $(G/H)(\cc\ls{t})$ and $\ld{\Lambda}_X^+$ sends a map $\lambda: \co_{G/H} \to \cc\ls{t}$ to the valuation given by the composite
    $$\co_{G/H} \to \co_{G/H} \otimes_\cc \co_G \xar{\lambda} \co_G\ls{t} \xar{v_t} \Z.$$
    This is a $G$-invariant discrete valuation of $\co_{G/H}$, 
\end{remark}

\begin{construction}
    Let $\cV^\perp$ denote the cone $\{\chi\in \fr{t}_X^\ast | \pdb{\chi, v}\leq 0 \text{ for each }v\in \cV\}$. Let $\Sigma_X$ denote the set of generators of intersections of extremal rays of $\cV^\perp$ with $\Lambda_X$. It turns out that the elements of $\Sigma_X$ are linearly independent; they are known as the \textit{spherical roots} of $X$. In fact, they form the set of simple roots of the based root system mentioned in \cref{cstr: g-invt valuations and little weyl}.
\end{construction}

\begin{remark}
    It turns out that for each spherical root $\gamma\in \Sigma_X$, there is some element $n\in \{\tfrac{1}{2}, 1, 2\}$ such that $\gamma' = n\gamma$ is either a positive root of $G$, or is the sum $\alpha + \beta$ of two positive roots which are orthogonal to each other and $\alpha$ and $\beta$ are elements of some system of simple roots. These simple roots need not correspond to the choice of $B$! Let $\Delta_X$ denote the set $\{\gamma'| \gamma\in \Sigma_X\}$; then $\Delta_X$ is called the set of \textit{normalized spherical roots}. Moreover, if $\Phi_X$ denotes the set of $W_X$-translates of $\Delta_X$, it is shown in \cite[Proposition 2.2.1]{sakellaridis-venkatesh} that the pair $(\Phi_X, W_X)$ defines a root system (called the \textit{normalized spherical root system} of $X$) where $\Delta_X$ forms a set of simple roots. Let $(\ld{\Phi}_X, W_X)$ denote the dual root system, and $\ld{\Delta}_X$ the set of simple coroots.
\end{remark}
\begin{theorem}[{\cite[Proposition 2.2.2]{sakellaridis-venkatesh}, \cite{knop-schalke}}]\label{thm: sv forms root datum}
    Suppose that $\Sigma_X$ does not contain any elements of the form $2\alpha$ for $\alpha$ being a root of $G$. Then, $(\Lambda_X, \Phi_X, \ld{\Lambda}_X, \ld{\Phi}_X)$ forms a root datum, with associated split complex reductive group $G_X$.
\end{theorem}
\begin{definition}\label{def: SV dual group}
    Let $\ld{G}_X$ denote the complex reductive group with maximal torus $\ld{T}_X$ with root datum given by the dual of that of \cref{thm: sv forms root datum}. We will refer to $\ld{G}_X$ as the \textit{(Langlands) dual group} of $X$. It admits a morphism to $\ld{G}$. Also see \cite{gaitsgory-nadler, knop-schalke}.
\end{definition}
\begin{example}
    As in \cref{ex: gp case}, if $X = G$ is viewed as a spherical $G\times G$-variety, the group $\ld{G}_X$ is simply the Langlands dual $\ld{G}$ of $G$ itself.
\end{example}
\begin{example}[Spherical $\PGL_2$-varieties]\label{ex: spherical pgl2-varieties}
    Recall the classification of spherical subgroups $H\subseteq \PGL_2$ from \cref{ex: spherical subgroups of PGL2}. Let us describe the root datum of $X = \PGL_2/H$ from \cref{thm: sv forms root datum} in each case.
    \begin{enumerate}
        \item If $H = \PGL_2$, the quotient $X$ is a point, and everything is trivial.
        \item If $H = \GG_m$, the orbits of $B$ on $X$ are the same as orbits of $\GG_m$ on $\PP^1$. There are therefore three orbits, given by $\GG_m$ (the open orbit) and the points $0$ and $\infty$.
        To describe the spherical roots, let us instead consider $\SL_2/\GG_m \cong (\PP^1 \times \PP^1) - \PP^1_\mathrm{diag}$. Note that $\co_{\SL_2/\GG_m} = \co_{\SL_2}^{\GG_m} \cong \bigoplus_{n\geq 0} V_{n\alpha}$, where $\alpha$ is the positive root of $\SL_2$ and $V_{n\alpha}$ is the representation with highest weight $n$. It follows that $\Lambda_X \cong \Z$, generated by $\alpha$. A little calculation implies that $\cV\subseteq \fr{t}_X$ identifies with $\{v\in \fr{t}_X | \pdb{v, \alpha} \leq 0\}$. This implies that $\Sigma_X = \Delta_X = \{\alpha\}$, and so $\ld{G}_X = \PGL_2$. If we worked with $\PGL_2/\GG_m$ instead, we would find that $\ld{G}_X = \SL_2$.
        \item If $H = \N_{\PGL_2}(\GG_m)$, the sublattice $\Lambda_X \subseteq \Lambda_{\PGL_2/\GG_m}$ has index two. In particular, by (b) above, we see that $\Lambda_X = \Z\cdot 2\alpha$, and $\Sigma_X = \{2\alpha\}$. In particular, \cref{thm: sv forms root datum} does not apply to this particular case.
        \item If $H = S\cdot N\subseteq \PGL_2$, the orbits of $B$ on $X$ are the same as orbits of $H$ on $\PP^1$. There are therefore two orbits, given by $\AA^1$ (the open orbit) and the point $\infty$. Let us assume for simplicity that $S = \{1\}$. Again, $\Lambda_X \cong \Z$, and one now calculates that $\Sigma_X$ is empty. One therefore finds that $\ld{G}_X = \ld{T}$. In general, the dual group of horospherical varieties is the Cartan subgroup.
    \end{enumerate}
    The cases (b), (c), and (d) above are known as types $T$, $N$, and $U$. The spherical $\PGL_2 \times \PGL_2$-variety $\PGL_2$ (i.e., the group case of \cref{ex: gp case}) is known as type $G$.
\end{example}
\begin{remark}
    If $\alpha$ is a simple root of $G$ (or $\alpha$ and $\beta$ are two orthogonal simple roots of $G$) and $P_\alpha$ (or $P_{\alpha\beta}$) is the associated parabolic subgroup, then the spherical variety $\punc{X} P_\alpha/U_{P_\alpha}$ is isomorphic to one of $\PGL_2/\PGL_2$, $\PGL_2/T$ for $T$ being a torus, $\PGL_2/\N_{\PGL_2}(T)$, or $(\PGL_2 \times \PGL_2)/\PGL_2^\mathrm{diag}$. Correspondingly, the unique element of $\Sigma_X$ is a normalized spherical root, and its type is as defined in \cref{ex: spherical pgl2-varieties}. In particular, the condition of \cref{thm: sv forms root datum} asks that $X$ have no normalized spherical root of type $N$.
\end{remark}
\begin{remark}\label{rmk: sv arthur sl2}
    Assume from now on that $X$ does not have any spherical roots of type $N$.
    As in \cite[Section 3.6]{sakellaridis-venkatesh}, the embedding $\ld{G}_X \hookrightarrow \ld{G}$ commutes with the image of a principal $\SL_2 \to \ld{L}(X)$. In particular, there is a map $\iota: \ld{G}_X \times \SL_2 \to \ld{G}$ such that upon restriction to the diagonal torus $\GG_m \subseteq \SL_2$, the map $\GG_m \to \ld{L}(X)$ is given by $2\rho_{L(X)} = \sum_{\alpha\in \Phi^+_{L(X)}} \alpha$ (regarded as a coweight of $\ld{G}$).
    Since we will mainly deal with spherical varieties of rank $1$ below, where $\ld{G}_X$ itself will sometimes be $\SL_2$, we will distinguish the $\SL_2$ above with a superscript: namely, we will write it as $\SL_2^\arth$.
\end{remark}

%% file: equivalences/whittaker-induction.tex
\subsection{Whittaker induction and \cite[Conjecture 7.5.1]{bzsv}}\label{subsec: stating the bzsv conjecture}

In this section, we will review the notion of Whittaker induction (following \cite[Section 3.4]{bzsv}), and the statement of \cite[Conjecture 7.5.1]{bzsv}. This construction takes as input a map $H \times \SL_2^\arth \to G$ and produces a functor from Hamiltonian $H$-spaces to Hamiltonian $G$-spaces. We warn the reader that our notation will differ slightly from that of \cite[Section 3.4]{bzsv}.
\begin{recall}
    A \textit{Hamiltonian $G$-space} is a smooth symplectic variety $M$ (with symplectic form $\omega$) equipped with a Hamiltonian $G$-action (i.e., the map $i: \g \to T_M$ given by the derivative of the $G$-action lands in the subspace of Hamiltonian vector fields on $M$). The moment map $\mu: M \to \g^\ast$ is characterized by the property that for each $x\in \g$, we have $d\pdb{\mu, x} = \pdb{i(x), \omega}$. We will often simply specify a Hamiltonian $G$-space as the pair $(M, \omega)$ along with its moment map. There will frequently be a grading present, which we encode by an action of $\GG_{m,\rot}$ on $M$, $G$, and $\omega$. We will say that $(M, \omega, \mu: M \to \g^\ast)$ is a \textit{graded} Hamiltonian $G$-space (for a given $\GG_{m,\gr}$-action on $G$) if $M$ has a $\GG_{m,\gr}$-action which acts on $\omega$ with weight $2$, and the moment map $\mu$ is $\GG_{m,\gr}$-equivariant.
\end{recall}

Let us review the basic example of Whittaker induction.
\begin{example}
    Let $G$ be a connected reductive group (over $\cc$), and let $e\in \g$ be a principal nilpotent element, so that the Jacobson-Morozov theorem produces a map $\SL_2^\arth \to G$. Let $H$ be the trivial group, and let $M$ denote the trivial Hamiltonian $H$-space. Then the Whittaker induction of $M$ along the map $\iota: \{1\} \times \SL_2^\arth \to G$ is given by $\WInd_{\iota}^G(M) = (e + \fr{b}^-) \times^{N^-} G$, where $N^-$ is the unipotent radical of the opposite Borel subgroup $B^-$ corresponding to $e$, and $\fr{b}^-$ is the Lie algebra of $B^-$. Note that \cref{thm: kostant} gives isomorphisms
    $$\WInd_\iota^G(M)/G \cong (e + \fr{b}^-)/N^- \cong \g\mmod G.$$
\end{example}
Let us now describe the construction in general.
\begin{construction}\label{cstr: n mod n- hamiltonian HN}
    Suppose we are given a map $H \times \SL_2^\arth \to G$ of reductive algebraic groups over $\cc$ such that $H$ centralizes the map $\SL_2^\arth \to G$. Let $f\in \g$ be the image of $\begin{psmallmatrix}
        0 & 0\\
        1 & 0
    \end{psmallmatrix}\in \sl_2^\arth$ inside $\g$. The action of $\GG_m^\arth \subseteq \SL_2^\arth$ on $\g$ defines a decomposition
    $$\g = \fr{z}^\arth \oplus \ol{\fr{n}} \oplus \fr{n}_0 \oplus \fr{n},$$
    where $\fr{z}^\arth$ is the centralizer of $\sl_2^\arth \to \g$, and $\ol{\fr{n}}$, $\fr{n}_0$, and $\fr{n}$ are the negative, zero, and positive weight spaces. Let $N$ denote the associated unipotent subgroup of $G$. Note that all the weights of the $\GG_m^\arth$-action on $\g$ are integers, and that $e\in \fr{n}$. Note that the orthogonal complement to $\fr{z}^\arth \subseteq \g$ is a Levi subalgebra $\fr{l} \subseteq \g$. Let $L\subseteq G$ denote the associated subgroup.
    
    Let $\fr{n}_+$ denote the subspace of ${\fr{n}}$ of elements with weight $\geq 2$, and let $N_+$ denote the associated unipotent subgroup. One can then equip ${\fr{n}}/\fr{n}_+$ with the structure of a Hamiltonian $H{N}$-space. There is an $H$-invariant symplectic form $\omega$ on ${\fr{n}}/\fr{n}_+$, given by $\omega(x,y) = \pdb{f, [x,y]}$.\footnote{Note that this symbol is well-defined: if $x\in \fr{n}_+$, then $[x,y]$ lives in weight $\geq 3$, so $\pdb{f, [x,y]} = 0$ since $f$ has weight $-2$. Moreover, this form is indeed nondegenerate: if $x \in {\fr{n}}$ is nonzero of weight $1$, then $[f,x]$ is a nonzero element of weight $-1$. This implies that there is some $y\in {\fr{n}}$ of weight $1$ such that $\pdb{[f,x], y} = \pdb{f, [x,y]}$ is nonzero, as desired.}
    Since $H$ preserves $\omega$, we obtain a homomorphism $H \to \Sp({\fr{n}}/\fr{n}_+)$, and hence a map $\fr{h} \to \fr{sp}_{{\fr{n}}/\fr{n}_+}$.
    The group $H$ acts on ${\fr{n}}/\fr{n}_+$ by the adjoint action. Moreover, the group ${N}$ acts on ${\fr{n}}/\fr{n}_+ \cong {N}/N_+$ via translation. The moment map $\mu: {\fr{n}}/\fr{n}_+ \to \fr{h}^\ast \oplus {\fr{n}}^\ast$ is defined as follows:
    \begin{itemize}
        \item The map ${\fr{n}}/\fr{n}_+ \to \fr{h}^\ast$ is adjoint to the map
        $${\fr{n}}/\fr{n}_+ \oplus \fr{h} \to {\fr{n}}/\fr{n}_+ \oplus \fr{sp}_{{\fr{n}}/\fr{n}_+} \xar{(x,g) \mapsto \frac{1}{2} \omega(gx, x)} \g_a.$$
        \item The map ${\fr{n}}/\fr{n}_+ \to \fr{n}^\ast$ is given by the composite
        $${\fr{n}}/\fr{n}_+ \xar{\omega} ({\fr{n}}/\fr{n}_+)^\ast \xar{x \mapsto f + x} {\fr{n}}^\ast.$$
        Here, $f$ is viewed as an element of ${\fr{n}}^\ast$ via the identification ${\fr{n}}^\ast \cong \fr{n}$. Under this isomorphism, the image of ${\fr{n}}/\fr{n}_+$ inside $\fr{n}$ is simply $f + \fr{n}_1$, where $\fr{n}_1$ is the weight $1$ eigenspace.
    \end{itemize}
\end{construction}
\begin{remark}\label{rmk: grading on n/n-}
    There is a natural grading defined on ${\fr{n}}/\fr{n}_+$, as well as a natural $\GG_{m,\gr}$-action on ${N}$ via the conjugation action of $\GG_m^\arth$. If $H$ is equipped with the trivial $\GG_{m,\gr}$-action, the Hamiltonian $H{N}$-space ${\fr{n}}/\fr{n}_+$ from \cref{cstr: n mod n- hamiltonian HN} can be viewed as a {graded} Hamiltonian $H{N}$-space.
\end{remark}

\begin{definition}\label{def: whittaker induction}
    Fix a map $\iota: H \times \SL_2^\arth \to G$ of reductive algebraic groups over $\cc$ such that $H$ centralizes the map $\SL_2^\arth \to G$.
    The conjugation action of $\GG_m^\arth$ on $G$ composed with the square character equips $G$ with a grading (which we will think of as a $\GG_{m,\gr}$-action).
    Let $M$ be a graded Hamiltonian $H$-space. Then the \textit{Whittaker induction} $\WInd_\iota^G(M)$ is defined as
    $$\WInd_\iota^G(M) = (M \times {\fr{n}}/\fr{n}_+) \times^{H{N}}_{\fr{h}^\ast \oplus {\fr{n}}^\ast} (T^\ast G),$$
    where $T^\ast G$ is regarded as a Hamiltonian $H{N}$-space via restriction along $H{N} \subseteq G$. There is a natural grading on $\WInd_\iota^G(M)$, coming from the grading on $M$, the grading on ${\fr{n}}/\fr{n}_+$ from \cref{rmk: grading on n/n-}, and the grading on $T^\ast G$ coming from the $\GG_{m,\gr}$-action on $G$.
    In particular, note that there is an isomorphism of stacks
    $$\WInd_\iota^G(M)/G \cong ( (M \times {\fr{n}}/\fr{n}_+) \times_{\fr{h}^\ast \oplus {\fr{n}}^\ast} \fr{g}^\ast )/H{N}.$$
\end{definition}
The simplest way to understand Whittaker induction in the case when $M$ is a symplectic $H$-representation is as follows.
\begin{lemma}[{\cite[Section 3.4.8]{bzsv}}]\label{lem: wind and ind}
    Suppose $M$ is a symplectic $H$-representation, and fix an isomorphism $\g^\ast \cong \g$. Then there is an isomorphism of stacks
    $$\WInd_\iota^G(M)/G \cong (M \oplus (\fr{h}^\perp \cap \g^e))/H$$
    over $BG$.
\end{lemma}
\begin{proof}
    Using \cite[Lemma 2.1]{gan-ginzburg}, one obtains an inclusion $f + \g^e \subseteq f + \fr{n}_+^\perp$ which is a slice of the $N$-action on $f + \fr{n}_+^\perp$. Therefore, there is an isomorphism
    \begin{align*}
        N \times (M \times_{\fr{h}^\ast} \g^e) & \to (M \times {\fr{n}}/\fr{n}_+) \times_{\fr{h}^\ast \oplus {\fr{n}}^\ast} \g^\ast\\
        & \cong \{(v,x)\in M \times (f + \fr{n}_+^\perp) \text{ such that } \mu(v) = x|_{\fr{h}}\},
    \end{align*}
    sending $(n, v, y) \mapsto (v, n \cdot (f + y))$. This isomorphism is $H$-equivariant, so it follows that $\WInd_\iota^G(M)/G$ is isomorphic to $(M \times_{\fr{h}^\ast} \g^e)/H$ as stacks over $BG$. This implies the desired claim, since $M \times_{\fr{h}^\ast} \g^e \cong M \oplus (\fr{h}^\perp \cap \g^e)$.
\end{proof}
\begin{remark}
    An alternative way to describe Whittaker induction using the language of shifted symplectic geometry (reviewed briefly in \cref{sec: ku-hamiltonian}) is as follows. Recall (see \cite{safronov-cs}) that a Lagrangian morphism $L \to \fr{h}^\ast/H$ is equivalent to the data of a Hamiltonian $H$-space $M$; the correspondence sets $L = M/H$. Moreover, intersecting Lagrangian correspondences produces another Lagrangian correspondence (see \cref{recall: shifted sympl}).
    From this perspective, one can describe Whittaker induction as follows. Let $\iota: H \times \SL_2^\arth \to G$ be a map of reductive algebraic groups over $\cc$ such that $H$ centralizes the map $\SL_2^\arth \to G$. Let $f + \g^e$ denote the Slodowy slice of the $\fr{sl}_2$-triple; then, there is a Lagrangian correspondence
    $$\xymatrix{
    & (f + \g^e)/H \ar[dl] \ar[dr] & \\
    \fr{h}^\ast/H & & \g^\ast/G,
    }$$
    and Whittaker induction amounts to intersecting the above Lagrangian correspondence with the Lagrangian $M/H \to \fr{h}^\ast/H$. (This will produce a Lagrangian morphism $(M \times_{\fr{h}^\ast} (f + \g^e)/H \to \g^\ast/G$, which is identified with the Hamiltonian $G$-variety of \cref{def: whittaker induction}.)
\end{remark}


Let us now recall a statement of \cite[Conjecture 7.5.1]{bzsv}; our presentation will follow \cite[Section 4.3]{bzsv}. Assume for now that $X$ is an affine spherical $G$-variety over $\cc$ which is the affine closure of its open $G$-orbit (for instance, this holds if $X$ is affine and homogeneous).
\begin{definition}\label{def: colors}
    A \textit{color} of $X$ is an irreducible $B$-stable divisor which is not $G$-stable (if $X$ is homogeneous, this is simply an irreducible $B$-stable divisor). 
    Following \cite[Definition 4.3.4]{bzsv}, a standard parabolic $P\subseteq G$ is said to be of \textit{even spherical type} if the spherical $P/U_P$-variety $X^\circ P/U_P$ is isomorphic to either the spherical $\SO_{2n+1}$-variety $\SO_{2n+1}/\SO_{2n}$ or the spherical $\G_2$-variety $\G_2/\SL_3$. (Note that there are diffeomorphisms $\SO_{2n+1}/\SO_{2n} \cong S^{2n}$ and $\G_2/\SL_3 \cong S^6$.) A color $D$ is said to be of \textit{even spherical type} if it meets $X^\circ P$ for a standard parabolic $P$ of even spherical type. Let $\cC_X$ denote the set of colors of $X$ of even spherical type.

    Suppose that the elements of $\cC_X$ freely generate a direct summand of $\ld{\Lambda}_X$. Let $\cd_X$ denote the set of dominant $W_X$-translates of $\cC_X\subseteq \ld{\Lambda}_X$, and let $\cd_X^{\max}$ denote the subset of maximal elements of $\cd_X$ (with respect to the ordering via coroots of $\ld{G}_X$). Let $S_X$ denote the $\ld{G}_X$-representation with highest weights $\cd_X^{\max}$. It is expected (see \cite[Conjecture 4.3.16]{bzsv}) that $S_X$ admits an $\ld{G}_X$-invariant symplectic form.
\end{definition}
\begin{example}[{\cite[Example 4.3.9]{bzsv}}]\label{ex: SX for GL2 mod Gm}
    Consider the example of the spherical $\GL_2$-variety $X = \GL_2/\GG_m$ (in which case $\ld{G}_X = \ld{G} = \GL_2$). Then $U\backslash X^\circ \cong \GG_m^2$ via the map $\begin{psmallmatrix}
        a & b\\
        c & d
    \end{psmallmatrix} \mapsto (b,d^{-1} \det)$. The colors of $X$ are given by the vanishing loci of $b$ and $d$, and are both of even spherical type. As explained in \cite[Example 4.3.9]{bzsv}, this implies that $\cC_X$ is the subset $\{\ld{\alpha}_1, -\ld{\alpha}_2\}$ of $\ld{\Lambda}_X = \ld{\Lambda}$, which in turn implies that $S_X = \AA^2 \oplus (\AA^2)^\ast \cong T^\ast(\AA^2)$ as an $\ld{G}_X$-representation.
    However, as remarked in \cite{bzsv}, the condition that the elements of $\cC_X$ freely generate a direct summand of $\ld{\Lambda}_X$ is \textit{not} true in the example of $\PGL_2/\GG_m$ (whose dual group is $\ld{G}_X = \ld{G} = \SL_2$). Nevertheless, the variant of \cref{def: colors} discussed in \cite[Section 4.4]{bzsv} shows that $S_X$ is the $\ld{G}_X$-representation $T^\ast(\AA^2)$.
\end{example}
\begin{example}\label{ex: GLn mod GLn1 for n at least 2}
    For $n>2$, the spherical $\GL_n$-variety $X = \GL_n/\GL_{n-1}$ still has $\ld{G}_X = \GL_2$, but the representation $S_X$ is zero. (I am very grateful to Justin Hilburn and Yiannis Sakellaridis for this point.) 
    For instance, when $n=3$, the Whittaker induction $\WInd_\iota^{\GL_3} S_X$ along the map $\iota: \GL_2 \times \SL_2^\arth \to \GL_3$ of \cref{rmk: sv arthur sl2} can be identified with $T^\ast(\GL_3/\GL_2)$ using \cref{lem: wind and ind}.
\end{example}
\begin{example}\label{ex: SO4 mod SO3}
    Consider the example of the spherical $\SO_4/\mu_2$-variety $\SO_4/\mu_2\cdot \SO_3$ (in which case $\ld{G}_X = \SL_2$). Since $\Spin_4 \cong \SL_2 \times \SL_2$, there is an isomorphism $\SO_4/\mu_2 \cong \SO_3 \times \SO_3$, under which the embedding of $\SO_3$ into $\SO_4/\mu_2$ is given by the diagonal. Therefore, there is an isomorphism $\SO_4/\mu_2 \SO_3 \cong \SO_3$, and this spherical $\SO_4/\mu_2$-variety can be understood as the group case for $\SO_3$. Using this, one can show that $\ld{G}\backslash \WInd_\iota^{\ld{G}} S_X \cong \sl_2/\ld{G}_X$. 
\end{example}

The following is a slight variant of \cite[Conjecture 7.5.1]{bzsv}.
\begin{conjecture}\label{conj: bzsv}
    Suppose $X$ is a smooth affine spherical $G$-variety over $\cc$ which is the affine closure of its open $G$-orbit, and with no roots of type N. Let $\iota: \ld{G}_X \times \SL_2^\arth \to \ld{G}$ denote the map of \cref{rmk: sv arthur sl2}. Suppose that $S_X$ admits an $\ld{G}_X$-invariant symplectic form, and let $\ld{M}$ denote $\WInd_\iota^{\ld{G}} S_X$. Then:
    \begin{itemize}
        \item There is an equivalence\footnote{The $\infty$-category on the left-hand side is defined in \cref{def: Shv-Sat LG/H} using the action of $\Shv^{c,\Sat}_{(G\times G)\pw{t}}(G\ls{t}; \QQ)$ on $\Shv^c_{G\pw{t}}(X\ls{t}; \QQ)$.}
        $$\Shv^{c,\Sat}_{G(\cc\pw{t})}(X(\cc\ls{t}); \QQ) \cong \Perf(\sh^{1/2} \ld{M}/\ld{G}(-2\rho)).$$
        \item This equivalence is equivariant for the actions of $\Shv^{c,\Sat}_{G(\cc\pw{t})\times G(\cc\pw{t})}(G(\cc\ls{t}); \QQ)$ and $\Perf(\ld{\g}^\ast[2-2\rho]/\ld{G}[-2\rho])$ under the equivalence of \cref{thm: derived satake}.
    \end{itemize}
\end{conjecture}

\begin{remark}\label{rmk: pointings and bzsv}
    One of the requirements for the equivalence of \cref{conj: bzsv} is the ``pointing'' of \cite[Section 7.5.2]{bzsv}. Namely, the pushforward of the constant sheaf along $i: X(\cc\pw{t}) \to X(\cc\ls{t})$ must be sent under the equivalence of \cref{conj: bzsv} to the structure sheaf of $\sh^{1/2} \ld{M}/\ld{G}$. This implies, in particular, that
    $$\End_{\Shv^{c,\Sat}_{G(\cc\pw{t})}(X(\cc\ls{t}); \QQ)}(i_\ast \ul{\QQ}_{X(\cc\pw{t})}) \simeq \co_{\sh^{1/2} \ld{M}}^{\ld{G}}.$$
    The left-hand side is simply $C^\ast_{G(\cc\pw{t})}(X(\cc\pw{t}); \QQ) \simeq C^\ast_G(X; \QQ)$, while the right-hand side is $\co_{\sh^{1/2} \ld{M}\mmod \ld{G}}$. Therefore, the ``pointing'' requirement can be restated as the existence of an equivalence of $\E{1}$-$\QQ$-algebras $C^\ast_G(X; \QQ) \simeq \co_{\sh^{1/2} \ld{M}\mmod \ld{G}}$.
    If $X = G/H$, the left-hand side is exactly $C^\ast_H(\ast; \QQ) \simeq \sh^{1/2} \H^\ast_H(\ast; \QQ)$, so this equivalence can be rephrased as a graded isomorphism
    \begin{equation}\label{eq: displayed pointing and bzsv}
        \ld{M}\mmod \ld{G} \cong \spec \H^\ast_H(\ast; \QQ) \cong \ld{\fr{h}}^\ast(2)\mmod \ld{H}.
    \end{equation}
    Using \cref{lem: wind and ind}, one can identify $\ld{M}\mmod \ld{G} \cong (S_X \oplus (\ld{\g}_X^\perp \cap \ld{\g}^e))\mmod \ld{G}_X$; it might be possible to prove the resulting identification with $\ld{\fr{h}}^\ast(2)\mmod \ld{H}$ in a direct manner (without having first established \cref{conj: bzsv}).
    One approach to proving \cref{eq: displayed pointing and bzsv} is to construct the Cartan $\fr{t}_H$ and the Weyl group $W_H$ of $H$ from $\ld{M}$. 
\end{remark}
There are several variants of \cref{conj: bzsv}, e.g., where one allows for some ramification. For instance, in the case of tame ramification, local geometric Langlands suggests the following (which is closely related to \cite[Conjecture 1.1.3]{finkelberg-ginzburg-travkin}):
\begin{conjecture}\label{conj: iwahori-ramification}
    Let $I\subseteq G\pw{t}$ be an Iwahori subgroup associated to a Borel $B\subseteq G$. Suppose $X$ is a smooth affine spherical $G$-variety over $\cc$ which is the affine closure of its open $G$-orbit, and with no roots of type N. Let $\ld{M}$ denote its dual Hamiltonian $\ld{G}$-space \`a la \cite{bzsv}. 
    Then:
    \begin{itemize}
        \item There is an equivalence
        $$\Shv^{c,\Sat}_I(X(\cc\ls{t}); \QQ) \simeq \Perf(\sh^{1/2} (\ld{M} \times_{\ld{\g}} \tilde{\ld{\g}})/\ld{G}),$$
        and the image of $\IC_0 = i_! \ul{\QQ}$ under the above equivalence should be the structure sheaf of $\sh^{1/2} (\ld{M} \times_{\ld{\g}} \tilde{\ld{\g}})/\ld{G}$. Here, $\tilde{\ld{\g}} \cong T^\ast(\ld{G}/\ld{N})/\ld{T}$ denotes the Grothendieck-Springer resolution over $\ld{\g}^\ast$; and $\Shv^{c,\Sat}_I(X(\cc\ls{t}); \QQ)$ denotes the full subcategory of $\Shv^{c}_I(X(\cc\ls{t}); \QQ)$ generated by $\IC_0$ under the action of $\Shv^{c,\Sat}_{I\times I}(G(\cc\ls{t}); \QQ)$ via convolution.
        \item This equivalence should be equivariant for the natural action of
        $$\Shv^{c,\Sat}_{I\times I}(G(\cc\ls{t}); \QQ) 
        \simeq \Perf(\sh^{1/2} (\tilde{\ld{\g}} \times_{\ld{\g}^\ast} \tilde{\ld{\g}})/\ld{G}).$$
        on both sides.
        This equivalence is provided by \cite{bezrukavnikov-two-geometric}.
    \end{itemize}
\end{conjecture}
There is an obvious variant of \cref{conj: iwahori-ramification} for standard parahorics, where the relevant replacement of the equivalence of \cite{bezrukavnikov-two-geometric} is proved in \cite{chen-dhillon-parabolic}.
\begin{remark}
    Suppose, for instance, that $X = G/H$. As in \cref{rmk: pointings and bzsv}, the first part of \cref{conj: iwahori-ramification} then implies that there should be an isomorphism
    $$(\ld{M} \times_{\ld{\g}} \tilde{\ld{\g}})\mmod \ld{G}
    \cong \spec \H^\ast_T(G/H; \QQ) \cong \ld{\fr{t}}^\ast(2)\times_{\ld{\g}^\ast(2)\mmod \ld{G}} \ld{\fr{h}}^\ast(2)\mmod \ld{H}.$$
    This might be easier to prove than the isomorphism described in \cref{eq: displayed pointing and bzsv}. 
\end{remark}

Let us end this section by stating our main result regarding affine homogeneous rank one spherical varieties, as classified by \cite{ahiezer-rk-1}, and their conjectured dual Hamiltonian spaces. This relies on \cref{table: topology and dualities for rank 1 spherical varieties}, part of which is essentially lifted from \cite[Section A.3.6]{sakellaridis-venkatesh} and \cite[Table 1]{sakellaridis-rank-1}.
The main result of this article is the following, proved in \cref{sec: case by case}:
\begin{theorem}\label{thm: rk 1 bzsv is true}
    Assuming \cref{hypothesis: rank 1 weakly placid}, the first part of 
    \cref{conj: bzsv} is true for each of the affine homogeneous rank $1$ spherical varieties in \cref{intro table: topology and dualities for rank 1 spherical varieties}. 
\end{theorem}
\begin{remark}\label{rmk: dual splits as main term and normalization}
    For each row in \cref{intro table: topology and dualities for rank 1 spherical varieties}, there is a graded isomorphism
    $$\ld{M}/\ld{G} \cong \ld{Y}/\ld{G}_X \times ``\text{Normalization}"$$
    of stacks, where the normalization term is described in \cref{table: topology and dualities for rank 1 spherical varieties}. 
    In fact, this is true more generally: this essentially follows from \cref{lem: wind and ind}, which identifies $S_X \oplus (\ld{\g}_X^\perp \cap \ld{\g}^e)$ with $\ld{Y} \times ``\text{Normalization}"$. This is at least somewhat surprising: for instance, when $X = \GL_{n+1}/\GL_n$ (so $\ld{G}_X = \GL_2$), we have $S_X = T^\ast(\AA^2)$ when $n=1$ by \cref{ex: SX for GL2 mod Gm}; but $S_X = 0$ for $n\geq 2$. Nevertheless, $\ld{Y}$ always identifies with $T^\ast(\AA^2)$ as $\GL_2$-schemes.

    Empirical evidence suggests that in many cases, the normalization term can be identified with $\fr{l}_X^\wedge\mmod L^\wedge_X$, where $L^\wedge_X$ is the subgroup of $\ld{G}$ from \cite{knop-schalke}. (This is \textit{not} quite the Langlands dual of the Levi subgroup $L(X)$.) However, we do not know how to prove this in general; any general statement would be very interesting (it is closely related to \cref{rmk: suggestion kostant for duals}).
    
    In any case, \cref{thm: rk 1 bzsv is true} therefore states that there is a grading on $\ld{Y}$ and an equivalence 
    $$\Shv^{c,\Sat}_{G\pw{t}}(X\ls{t}; \QQ) \cong \Perf(\sh^{1/2}\ld{Y}/\ld{G}_X \times ``\text{Normalization}");$$
    it is this form of \cref{thm: rk 1 bzsv is true} which we will prove below. See \cref{table: topology and dualities for rank 1 spherical varieties}.
\end{remark}
\begin{remark}
    As mentioned in the introduction, the analogue of \cref{thm: rk 1 bzsv is true} in the real and $p$-adic cases was proved by Gan and Gomez as \cite[Theorem 1]{gan-gomez}.
\end{remark}
\begin{remark}
    Since most of the examples in \cref{intro table: topology and dualities for rank 1 spherical varieties} are topologically just spheres, one might wonder about the  case of $\SO_9/\Spin_7$ (indeed, the Hopf fibration gives a homotopy equivalence $\SO_9/\Spin_7 \simeq \RP^{15}$). However, despite being homotopy equivalent to a sphere, the spherical $\SO_9$-variety $\SO_9/\Spin_7$ is \textit{not} of rank $1$. Instead, one can compute that it is of rank $2$: its dual group is $\ld{G}_X = \SL_2 \times \SL_2$ (see \cite[Line 9 of final table]{sakellaridis-spherical-functions}). In particular, it does not fall under the purview of \cref{thm: rk 1 bzsv is true}. One can study this example using \cref{thm: ordinary homology criterion satake} below, but we will not do so in the present article.
\end{remark}
As mentioned at the end of \cref{subsec: review satake}, the proof of \cref{thm: rk 1 bzsv is true} is modeled after the proof of \cref{thm: derived satake}.
In fact, in the course of proving \cref{thm: rk 1 bzsv is true}, we will accomplish more than what is stated: namely, our argument also suggests that there is a ``chromatic deformation'' of \cref{conj: bzsv}, which we will explain in greater detail after the proof of \cref{thm: rk 1 bzsv is true} (which will in turn occupy \cref{sec: case by case}).

%% file: equivalences/G-to-H-equiv-homology.tex
\subsection{The $G$-equivariant $\ku$-(co)homology of $\cL(G/H)$}\label{subsec: G of LGmodH}

Fix a compact Lie group $G$, and let $H\subseteq G$ be a closed subgroup.
Throughout this section, we will always assume that $H$ and $G$ are connected, and also (for simplicity) that $G/H$ has finite fundamental group (so that $\Omega(G/H)$ has finitely many connected components). Recall \cref{notation: compact vs complexification}: we will write $G\ls{t}$ or $G\pw{t}$ below to mean $G_\cc\ls{t}$ or $G_\cc\pw{t}$, respectively.

Before proceeding to the proof of \cref{thm: rk 1 bzsv is true}, we need some general statements about the $G$-equivariant $\ku$-(co)homology of the free loop space $\cL(G/H) = \Map(S^1, G/H)$. Let us remind the reader uninterested in $\ku$, but interested in $\Z$, that one should feel free to replace $\ku_G$ everywhere below by the equivariant cochains $C^\ast_G(\ast; \Z)$ (and similarly for other symbols involving $\ku$); and moreover, any instance of the symbol $\beta$ can be set to zero to obtain a corresponding statement over $\Z$.

The following basic result is an analogue of the algebro-geometric fact that $G\backslash T^\ast(G/H) \cong H\backslash (\g/\fr{h})^\ast$, or its homotopic analogue that $(G/H)_+^{hG} \simeq (S^{\g/\fr{h}})^{hH}$.
\begin{prop}\label{prop: going from G of LG/H to H of Loops G/H}
    Let $H$ act on $G/H$, and hence on $\Omega(G/H)$, by conjugation (equivalently, left-translation). Then the $G$-space $\cL(G/H)$ is $G$-equivariantly homotopy equivalent to $\Ind_H^G \Omega(G/H)$. In particular, there is an equivalence of orbispaces
    $$G \backslash \cL(G/H) \simeq H\backslash \Omega(G/H).$$
\end{prop}
\begin{proof}
    It is a classical fact that the map $m: G \times \Omega G \to \cL G$ sending $(g, \gamma)$ to the loop $\gamma_g: t \mapsto g \gamma(t)$ is a homotopy equivalence. The left action of $G$ on $\cL G$ is simply given by
    $$G\ni g': \gamma_g(t) \mapsto \gamma_{g' g}(t),$$
    which allows us to identify $\cL G \simeq \Ind_{\{1\}}^G \Omega G$. Recall that there is a principal fibration
    $$H \to G \to G/H,$$
    which gives equivalences $\Omega (G/H) \simeq \Omega G/\Omega H$ and $\cL(G/H) \simeq \cL G/\cL H$. Since and the diagram
    $$\xymatrix{
    H \times \Omega H \ar[r] \ar[d]_-m & G \times \Omega G \ar[d]^-m \\
    \cL H \ar[r] & \cL G
    }$$
    commutes, we find that there is an equivalence of $G$-spaces
    $$\cL G/\cL H \simeq (\Ind_{\{1\}}^G \Omega G)/(\Ind_{\{1\}}^H \Omega H) \simeq \Ind_H^G \Omega(G/H),$$
    as desired. Alternatively, this also follows from \cref{eq: hochschild of orbispace} using that $\ast \times_{\ast \times_{\ast/G} \ast/H} \ast \simeq \Omega(G/H)$.
\end{proof}

Let us now state our key assumption; only the first part will be used in this subsection. The other assumptions will be used and explained in \cref{subsec: using the regular centralizer}.
\begin{hypothesis}\label{hypothesis: rank 1 weakly placid}
    Say that a connected reductive subgroup $H_\cc\subseteq G_\cc$ is \textit{optimal} if:
    \begin{enumerate}
        \item If $G$ (resp. $H$) denotes the maximal compact subgroup of $G(\cc)$ (resp. $H(\cc)$), there is a homotopy equivalence of orbifolds
        $$G_\cc(\cc\pw{t}) \backslash G_\cc(\cc\ls{t})/H_\cc(\cc\ls{t}) \simeq G\backslash \cL(G/H).$$
        \item The $G_\cc$-action on $G_\cc/H_\cc$ is defined over a finite type $\Z$-algebra $R\subseteq \cc$, and that the $G_R(R\pw{t})$-action on $(G_R/H_R)(R\ls{t})$ is placid in the sense of \cref{def: collected action}.
        \item Let $\ld{\cR}\star \IC_0 \in \Shv^{c,\Sat}_{G\pw{t}}(G\ls{t}/H\ls{t}; \QQ)$ denote the sheaf obtained by the action of $\co_{\ld{G}}\in \Rep(\ld{G})$ on $\IC_0$. Then the sheaves $\IC_0, \IC_0 \star \ld{\cR} \in \Shv^{c,\et}_{G(\ol{\FF_q}\pw{t})}(G(\ol{\FF_q}\ls{t})/H(\ol{\FF_q}\ls{t}))$ are pure of weight zero for $q\gg 0$.
        \item For the placid presentation of the $G_R(R\pw{t})$-action on $(G_R/H_R)(R\ls{t})$ from \cref{def: collected action}, the evenness assumption of \cref{thm: full faithful} are satisfied for $\IC_0$ and $\IC_0 \star \ld{\cR}$.
    \end{enumerate}
\end{hypothesis}

\begin{remark}
    We will assume \cref{hypothesis: rank 1 weakly placid}(a) for the remainder of this section. In the case of symmetric varieties --- in particular, for the cases of \cref{thm: rk 1 bzsv is true} \textit{not} of types $\G_2$ or $B'_3$ --- this was proved in \cite{mitchell-buildings}.
\end{remark}

\begin{remark}\label{rmk: cohomology going from G of LG/H to H of Loops G/H}
    There is a multiplicative presentation of $\Omega(G/H)$ as $\colim X_\lambda$ via finite $H$-spaces $X_\lambda$, and the induced $G$-spaces $\Ind_H^G X_\lambda$ defines a presentation of $\cL(G/H)$ by finite $G$-spaces.
    It follows that there is an equivalence $\cf_{G\pw{t}}(G\ls{t}/H\ls{t}) \cong \cf_H(\Omega(G/H))$ of $\Eoo$-$\ku_G$-algebras, where the right-hand side is viewed as an $\Eoo$-$\ku_G$-algebra via its natural $\ku_H$-algebra structure and the canonical map $\ku_G \to \ku_H$.
\end{remark}
\begin{warning}\label{warning: G to H only on cohomology}
    Although there is an equivalence $\cf_{G\pw{t}}(G\ls{t}/H\ls{t}) \cong \cf_H(\Omega(G/H))$ of $\Eoo$-$\ku_G$-algebras, there is \textit{not} an equivalence $\cf_{G\pw{t}}(G\ls{t}/H\ls{t})^\vee \cong \cf_H(\Omega(G/H))^\vee$ of $\ku_G$-modules. Indeed, $\cf_{G\pw{t}}(G\ls{t}/H\ls{t})^\vee$ denotes the \textit{$\ku_G$-linear} dual of $\cf_{G\pw{t}}(G\ls{t}/H\ls{t})$, while $\cf_H(\Omega(G/H))^\vee$ denotes the \textit{$\ku_H$-linear} dual of $\cf_H(\Omega(G/H))$. (To clarify, this dual is not taken in the naive sense: rather, if $\Omega(G/H) = \colim X_\lambda$ as in \cref{rmk: cohomology going from G of LG/H to H of Loops G/H}, $\cf_H(\Omega(G/H))^\vee$ means $\colim_\lambda \cf_H(X_\lambda)^\vee$.)
\end{warning}
\begin{remark}\label{rmk: splitting is not S1-equivariant}
    \cref{prop: going from G of LG/H to H of Loops G/H} breaks the natural symmetry on $G \backslash \cL(G/H)$. Namely, since the action of $G$ on $\cL(G/H)$ is defined via the $G$-action on $G/H$, the orbispace $G \backslash \cL(G/H)$ has an action of the circle $S^1_\rot$ given by rotating loops. However, this structure is not naturally visible on the orbispace $\Omega(G/H)/H$. Indeed, the proof of \cref{prop: going from G of LG/H to H of Loops G/H} used the splitting $G \times \Omega G \xar{\sim} \cL G$; but this splitting is \textit{not} $S^1_\rot$-equivariant.
\end{remark}

A slight variant of \cref{prop: going from G of LG/H to H of Loops G/H} lets us describe the $G$-equivariant $\ku$-\textit{co}homology of $\cL(G/H)$. The following result is proved nonequivariantly in \cite[Remark A.6]{grg-reg}.
\begin{prop}\label{prop: equiv coh as hochschild}
    There is an $S^1_\rot$-equivariant equivalence of $\Eoo$-$\ku_G$-algebras
    $$\cf_{G\pw{t}}(G\ls{t}/H\ls{t}) \simeq \HH(\ku_H/\ku_G),$$
    where the right-hand side denotes the relative Hochschild homology of the $\Eoo$-map $\ku_H \to \ku_G$ (equipped with its natural $S^1$-action).
\end{prop}
\begin{proof}
    Since $G/H$ is itself the fiber product $\ast \times_{\ast/G} \ast/H$ in orbispaces, there is an equivalence 
    $$G\backslash \cL(G/H) \simeq G\backslash (\ast \times_{\cL(\ast/G)} \cL(\ast/H)) \simeq \ast/G \times_{\cL(\ast/G)} \cL(\ast/H).$$
    But $\cL(\ast/G) \simeq \ast/G \times_{\ast/G \times \ast/G} \ast/G$, where the two maps $\ast/G \to \ast/G \times \ast/G$ are both given by the diagonal. Therefore, we can identify
    \begin{align}
        G\backslash \cL(G/H) & \simeq \ast/G \times_{\ast/G \times_{\ast/G \times \ast/G} \ast/G} (\ast/H \times_{\ast/H \times \ast/H} \ast/H) \nonumber \\
        & \simeq \ast/H \times_{\ast/H \times_{\ast/G} \ast/H} \ast/H. \label{eq: hochschild of orbispace}
    \end{align}
    By construction of equivariant $\ku$, it follows that there is an equivalence of $\Eoo$-$\ku_G$-algebras
    $$\cf_{G\pw{t}}(G\ls{t}/H\ls{t}) \simeq \ku_H \otimes_{\ku_H \otimes_{\ku_G} \ku_H} \ku_H = \HH(\ku_H/\ku_G).$$
    Moreover, the equivalence of \cref{eq: hochschild of orbispace} is manifestly $S^1$-equivariant, so we obtain the desired claim.
\end{proof}
Base-changing \cref{prop: equiv coh as hochschild} along the map $\ku \to \QQ$ and using the Hochschild-Kostant-Rosenberg theorem (in the form proved in \cite{raksit, toen-hkr}), one finds:
\begin{corollary}\label{cor: thm 1 of BF for G mod H}
    Let $\fr{c}_{\ld{H}}^\ast = \spec \sh^{1/2} \H^\ast_H(\ast; \QQ)$ and $\fr{c}_{\ld{G}}^\ast = \spec \sh^{1/2} \H^\ast_G(\ast; \QQ)$ are (the shearings of) the Chevalley bases for $\ld{H}$ and $\ld{G}$, respectively, so that $\fr{c}_{\ld{H}}^\ast \cong \ld{\fr{h}}^\ast[2]\mmod \ld{H}$ (and similarly for $\fr{c}_{\ld{G}}^\ast$).
    There is an equivalence
    $$C^\ast_{G\pw{t}}(G\ls{t}/H\ls{t}; \QQ) \simeq L\Omega^\ast_{\fr{c}_{\ld{H}}^\ast/\fr{c}_{\ld{G}}^\ast},$$
    where the right-hand side is derived Hodge cohomology. In other words, there is an isomorphism
    $$\spec C^\ast_{G\pw{t}}(G\ls{t}/H\ls{t}; \QQ) \cong T[-1](\fr{c}_{\ld{H}}^\ast/\fr{c}_{\ld{G}}^\ast).$$
    This equivalence identifies the loop rotation action of $\GG_m^\rot$ on the left-hand side with the de Rham differential on the right-hand side; this implies that there is an isomorphism
    $$\spec C^\ast_{G\pw{t} \rtimes \GG_m^\rot}(G\ls{t}/H\ls{t}; \QQ) \cong \Def_\hbar(\fr{c}_{\ld{H}}^\ast/\fr{c}_{\ld{G}}^\ast),$$
    where the right-hand side denotes the deformation to the normal cone of the morphism $\fr{c}_{\ld{H}}^\ast \to \fr{c}_{\ld{G}}^\ast$ in the sense of \cite[Section 9.2]{gr-ii}, living over the base $\spec C^\ast_{\GG_m^\rot}(\ast; \QQ) \cong \AA^1_\hbar[2]$.
\end{corollary}
\begin{example}
    For instance, in the group case (so $G = H \times H$ with $H$ embedded diagonally), the first isomorphism of \cref{cor: thm 1 of BF for G mod H} says that there is an isomorphism
    $$\spec C^\ast_{H\pw{t}}(\Gr_H; \QQ) \cong T[-1](\fr{c}_{\ld{H}}^\ast/\fr{c}_{\ld{H}\times \ld{H}}^\ast) \cong T[-2](\fr{c}_{\ld{H}}^\ast/\spec(\QQ)).$$
    Similarly, there is an isomorphism
    $$\spec C^\ast_{H\pw{t} \rtimes \GG_m^\rot}(\Gr_H; \QQ) \cong \mathrm{Def}_\hbar(\fr{c}_{\ld{H}}^\ast/\fr{c}_{\ld{H}\times \ld{H}}^\ast);$$
    since $\fr{c}_{\ld{H} \times \ld{H}}^\ast \cong \fr{c}_{\ld{H}}^\ast \times_{\spec(\QQ)} \fr{c}_{\ld{H}}^\ast$, this recovers the isomorphism of \cite[Theorem 1]{bf-derived-satake}. We have already discussed this perspective before in \cite[Example A.8]{grg-reg}.
\end{example}

Let us note that away from characteristic zero, or especially well-behaved cases like the map $\GL_n \to \GL_{n+1}$, one needs to be very careful in modeling the map $C^\ast_G(\ast; R) \to C^\ast_H(\ast; R)$ via shearing. See \cref{rmk: modeling maps via shearing} for an example of the sort of subtlety that can arise.
\begin{remark}\label{rmk: suggestion kostant for duals}
    Suppose $G/H$ is an affine spherical $G$-variety, and assume that \cref{eq: displayed pointing and bzsv} of \cref{rmk: pointings and bzsv} holds for $G/H$ (which would follows from \cref{conj: bzsv}). Then \cref{cor: thm 1 of BF for G mod H} implies that if $\mu: \ld{M} \to \ld{\g}^\ast$ denotes the moment map, there is an isomorphism (ignoring gradings for simplicity)
    $$\H^\ast_{G\pw{t}}(G\ls{t}/H\ls{t}; \QQ) \simeq L\Omega^\ast_{\ld{M}\mmod \ld{G} /\ld{\g}^\ast\mmod \ld{G}}.$$
    In fact, more is true: taking cohomology defines a functor 
    $$\Shv^{c,\Sat}_{G\pw{t}}(G\ls{t}/H\ls{t}; \QQ) \xar{\H^\ast_G(-; \QQ)} \Mod(\H^\ast_G(\cL G/\cL H; \QQ)).$$
    By \cref{cor: thm 1 of BF for G mod H}, the right-hand side can be identified with the $\infty$-category of modules over the shearing of $L\Omega^\ast_{\ld{\fr{h}}^\ast(2)\mmod \ld{H} / \ld{\g}^\ast(2)\mmod \ld{G}}$, i.e., the $\infty$-category of perfect complexes over the shearing of the $(-1)$-shifted tangent bundle $T[-1](\ld{\fr{h}}^\ast(2)\mmod \ld{H} / \ld{\g}^\ast(2)\mmod \ld{G})$. Under the isomorphism $\H^\ast_G(\cL G/\cL H; \QQ) \cong \H^\ast_H(\Omega(G/H); \QQ)$, the $\H^\ast_G(\ast;\QQ)$-module structure on $\H^\ast_G(\cL G/\cL H; \QQ)$ factors through the canonical map $\H^\ast_G(\ast; \QQ) \to \H^\ast_H(\ast; \QQ)$. This defines a factorization
    \begin{equation}\label{eq: factorization of cohomology}
        \xymatrix{
        \Shv^{c,\Sat}_{G\pw{t}}(G\ls{t}/H\ls{t}; \QQ) \ar@{-->}[d] \ar[rr]^-{\H^\ast_G(-; \QQ)} & & \Mod(\H^\ast_G(\cL G/\cL H; \QQ)) \ar@{-->}[dll] \ar[d]^-{\mathrm{forget}} \\
        \Mod(\H^\ast_H(\ast; \QQ)) \ar[rr]_-{\mathrm{restriction}} & & \Mod(\H^\ast_G(\ast; \QQ)),
        }
    \end{equation}
    which makes the triangles commute.
    
    By \cref{conj: bzsv}, $\Shv^{c,\Sat}_{G\pw{t}}(G\ls{t}/H\ls{t}; \QQ)$ can be identified with $\Perf(\sh^{1/2} \ld{M}/\ld{G})$. In particular, there is a natural map $\H^\ast_H(\ast; \QQ) \to \H^\ast_G(\cL G/\cL H; \QQ)$, which gives a functor
    $$\Perf(\sh^{1/2} \ld{M}/\ld{G}) \to \Perf(\sh^{1/2} T[-1](\ld{\fr{h}}^\ast(2)\mmod \ld{H} / \ld{\g}^\ast(2)\mmod \ld{G})) \to \Perf(\sh^{1/2} \ld{\fr{h}}^\ast(2)\mmod \ld{H}).$$
    When $G = H \times H$, so that $\ld{M} = T^\ast \ld{H}$, this is precisely the Kostant functor of \cite[Section 2.6]{bf-derived-satake}. This functor is compatible with the commutative diagram \cref{eq: factorization of cohomology}, in that the following diagram is its analogue on the spectral side:
    $$\xymatrix{
    \Perf(\sh^{1/2} \ld{M}/\ld{G}) \ar@{-->}[d]_-{\kappa_{\ld{M}}} \ar[rr]^-{\text{``Kostant functor''}} & & \Perf(\sh^{1/2} T[-1](\ld{\fr{h}}^\ast(2)\mmod \ld{H} / \ld{\g}^\ast(2)\mmod \ld{G})) \ar@{-->}[dll]^-{\text{zero section}^\ast} \ar[d] \\
    \Perf(\sh^{1/2} \ld{\fr{h}}^\ast(2)\mmod \ld{H}) \ar[rr]_-{\mathrm{restriction}} & & \Perf(\sh^{1/2} \ld{\g}^\ast(2)\mmod \ld{G}),
    }$$
    where we will now describe the dotted map denoted $\kappa_{\ld{M}}$.
    
    Let us ignore gradings in the following discussion.
    It is natural to expect that the above analogue of the Kostant functor is induced by pullback along a certain map
    $$\kappa_{\ld{M}}: \ld{\fr{h}}^\ast\mmod \ld{H} \to \ld{M}/\ld{G}.$$
    For instance, when $G = H \times H$, so that $\ld{M} = T^\ast \ld{H}$, the map $\kappa$ is simply the Kostant slice for $\ld{H}$. Moreover, in the general case, the compatibility of the equivalence of \cref{conj: bzsv} with the action of the Satake category implies that there is a commutative square
    $$\xymatrix{
        \ld{\fr{h}}^\ast\mmod \ld{H} \ar[r]^-{\kappa_{\ld{M}}} \ar[d] & \ld{M}/\ld{G} \ar[d]^-\mu \\
        \ld{\g}^\ast\mmod \ld{G} \ar[r]^-{\kappa} & \ld{\g}^\ast/\ld{G}.
    }$$
    Therefore, \cref{eq: displayed pointing and bzsv} and \cref{cor: thm 1 of BF for G mod H} together make the following concrete prediction (see \cref{rmk: extension of generalized kostant} for an additional component of our prediction): 
    \begin{conjecture}\label{conj: generalized kostant slice}
        On the spectral side of \cref{conj: bzsv}, if $\ld{M}$ is the Hamiltonian $\ld{G}$-space dual to $G/H$, there is an isomorphism 
        $$\ld{M} \mmod \ld{G} \cong \ld{\fr{h}}^\ast\mmod \ld{H}.$$
        In general, if $\ld{M}$ is the Hamiltonian $\ld{G}$-space which is dual to an affine (not necessarily homogeneous) spherical $G$-variety, there is a ``Kostant section'' $\kappa_{\ld{M}}: \ld{M} \mmod \ld{G} \to \ld{M}/\ld{G}$ which makes the following square commute:
        $$\xymatrix{
        \ld{M} \mmod \ld{G} \ar[r]^-{\kappa_{\ld{M}}} \ar[d] & \ld{M}/\ld{G} \ar[d]^-\mu \\
        \ld{\g}^\ast\mmod \ld{G} \ar[r]^-{\kappa} & \ld{\g}^\ast/\ld{G}.
        }$$
        Furthermore, $\kappa_{\ld{M}}$ can be refined (non-uniquely) to a map $\ld{M} \mmod \ld{G} \to \ld{M}$ such that the algebra of regular functions on its $\ld{G}$-orbit is isomorphic to $\co_{\ld{M}}$. 

        Motivated by \cite{knop-asymptotic-invariants, knop-set-of-orbits} (see \cref{rmk: knop microlocal weyl group} for a brief summary), we further expect the following. Suppose that $\ld{M}$ can be written as $T^\ast \ld{X}$ for a spherical $\ld{G}$-variety $\ld{X}$, and let $\co_{\ld{B}}(\ld{X})$ (resp. $\co_B(G/H)$) denote the poset of $\ld{B}$-orbit closures in $\ld{X}$ (resp. $B$-orbit closures in $G/H$) equipped with the Bruhat order. Then there is a bijection $\co_{\ld{B}}(\ld{X}) \leftrightarrow \co_B(G/H)$ which is equivariant for the action of the Weyl group $W_{\ld{G}} \cong W_G$ on either side described in \cite{knop-set-of-orbits}.\footnote{Said differently, there is a $W$-equivariant bijection between the sets of irreducible components of $T^\ast(G/H) \times_{\fr{b}^\ast} \{0\}$ and $\ld{M} \times_{\ld{\fr{b}}^\ast} \{0\}$, where $T^\ast(G/H) \to \fr{b}^\ast$ (resp. $\ld{M} \to \ld{\fr{b}}^\ast$) is the moment map for the $B$- (resp. $\ld{B}$-)action.} 
        
        Under this bijection, the results of \cite{knop-asymptotic-invariants, knop-set-of-orbits, ressayre-knop} (see \cref{rmk: knop microlocal weyl group} for a brief summary) should give an isomorphism between the stabilizers of a minimal rank $B$-orbit on $G/H$ and of the closure $\ld{X}$ of the open $\ld{B}$-orbit in $\ld{X}$; in other words, there is an isomorphism
        $$W_{\ld{X}} \ltimes W_{L(X)} \cong W_H.$$
        If $\fr{t}_H$ denotes a Cartan of $H$ and $\fr{t}_{\ld{X}}$ denotes the Cartan of $\ld{X}$, there should be a canonical $W_{L(X)}$-Galois covering $\fr{t}_H \to \fr{t}_{\ld{X}}^\ast$ which is equivariant for the surjection $W_H \twoheadrightarrow W_{\ld{X}}$.
    \end{conjecture}
    Some brief comments regarding \cref{conj: generalized kostant slice}:
    \begin{itemize}
        \item Just as the Kostant slice plays a crucial role in the geometric Langlands program, we expect the Kostant section $\kappa_{\ld{M}}$ to play a central role in the story of relative geometric Langlands.
        \item Since $\ld{M} = \Ind_{\ld{G}_X}^{\ld{G}} (S_X \oplus (\ld{\g}_X^\perp \cap \ld{\g}^e))$ by \cref{lem: wind and ind}, the first part of \cref{conj: generalized kostant slice} is equivalent to the statement that 
        $$(S_X \oplus (\ld{\g}_X^\perp \cap \ld{\g}^e))\mmod \ld{G}_X \cong \ld{\fr{h}}^\ast\mmod \ld{H}.$$
        The latter being a polynomial ring (by Chevalley restriction and Chevalley-Shephard-Todd), we see that \cref{conj: generalized kostant slice} forces in particular that the action of $\ld{G}_X$ on $S_X \oplus (\ld{\g}_X^\perp \cap \ld{\g}^e)$ be \textit{coregular}.
        \item The penultimate paragraph of \cref{conj: generalized kostant slice} clearly generalizes to the case when $T^\ast(G/H)$ is replaced by a more general Hamiltonian $G$-space. In fact, this generalization of the penultimate part of \cref{conj: generalized kostant slice} appears as \cite[Conjecture 1.1.1]{finkelberg-ginzburg-travkin}; I am grateful to Akshay Venkatesh and Zhiwei Yun for informing me of this paper. 
        \item The middle part of \cref{conj: generalized kostant slice} should be a starting point for proving \cref{conj: iwahori-ramification}.
        \item The final paragraph of \cref{conj: generalized kostant slice} attempts to provide a ``geometric reason'' for the expected isomorphism $(T^\ast \ld{X})\mmod \ld{G} \cong \fr{t}_H\mmod W_H$: indeed, as described in \cref{rmk: knop microlocal weyl group}, Knop has shown that $(T^\ast \ld{X})\mmod \ld{G} \cong \fr{t}_{\ld{X}}^\ast\mmod W_{\ld{X}}$, and the final part of \cref{conj: generalized kostant slice} identifies this with $\fr{t}_H \mmod W_H$.
    \end{itemize}
\end{remark}

\begin{corollary}\label{cor: loop homology E2}
    The $\E{1}$-$\ku_G$-algebra structure obtained via the $\Eoo$-map $\ku_G \to \ku_H$ on the $\ku_H$-linear dual of $\cf_H(\Omega(G/H))$ -- which is \emph{not} $\cf_H(\Omega(G/H))^\vee$ in the notation of \cref{cst: ku-hmlgy filtered colimit} -- refines to an $\E{2}$-$\ku_G$-algebra structure.
\end{corollary}
\begin{proof}
    Taking the $\ku_H$-linear dual of the right-hand side of \cref{prop: equiv coh as hochschild} produces the Hochschild \textit{co}homology $\HC(\ku_H/\ku_G)$. By the Deligne conjecture (in the form proved in \cite[Section 5.3]{HA}), this admits the structure of an $\E{2}$-$\ku_G$-algebra. On the other hand, by \cref{prop: going from G of LG/H to H of Loops G/H}, the right-hand side of \cref{prop: equiv coh as hochschild} can be identified with the equivariant cohomology $\cf_H(\Omega(G/H))$. The desired result follows.
\end{proof}

\begin{remark}
    One can also identify the $\ku_H$-linear dual of $\cf_H(\Omega G)$ with the $\E{2}$-centralizer of the map $\ku_G \to \ku_H$; see \cref{recall: centers}. The $\E{2}$-structure on the $\ku_H$-linear dual of $\cf_H(\Omega(G/H))$ is essentially the reason for the $\E{2}$-monoidal structure on the relative Langlands category $\Shv^{c,\Sat}_{G\pw{t}}(G\ls{t}/H\ls{t}; \QQ)$ (to be defined below) from \cite[Remark 7.5.12 and Section 16]{bzsv}.
    
    In the special case when $G = H \times H$ and $H$ is embedded diagonally, one can identify $\HC(\ku_H/\ku_G)$ with the \textit{$\E{2}$-Hochschild cohomology} $\HC_\E{2}(\ku_H/\ku)$. The Deligne conjecture therefore equips the $\ku_H$-linear dual of $\cf_H(\Omega H)$ with an $\E{3}$-algebra structure, and again this is essentially the source of the folklore $\E{3}$-monoidal structure on the spherical Hecke category $\Shv^{c,\Sat}_{H\pw{t} \times H\pw{t}}(H\ls{t}; \QQ)$. We will use the perspective of $\E{n}$-centers to establish some coisotropicity results in \cref{subsec: ku and lagrangians}.
\end{remark}
\begin{warning}\label{warning: completion hochschild}
    The reader should keep \cref{warning: G to H only on cohomology} in mind: the $\ku_H$-linear dual of $\cf_H(\Omega(G/H))$ is \textit{not} equivalent to the $\ku_G$-linear dual of $\cf_{G\pw{t}}(G\ls{t}/H\ls{t})$. In fact, as mentioned in \cref{cor: loop homology E2}, the $\ku_H$-linear dual of $\cf_H(\Omega(G/H))$ is also not equivalent to the equivariant homology $\cf_H(\Omega(G/H))^\vee$; the former is only a \textit{completion} of the latter.
\end{warning}
\begin{remark}\label{rmk: K-equiv-homology}
    There are, of course, many mild variants of \cref{cor: loop homology E2}. For instance, suppose $K \subseteq H$ is a closed subgroup. Then $\ku_K$ is an $\Eoo$-$\ku_H$-algebra; in particular, it is a $\ku_H$-bimodule in $\ku_G$-modules. Therefore, one can consider the Hochschild cohomology $\HC(\ku_H/\ku_G; \ku_K)$ with coefficients in the bimodule $\ku_K$. Just as in \cref{cor: loop homology E2}, one can identify $\HC(\ku_H/\ku_G; \ku_K)$ with the $\ku_K$-linear dual of $\cf_K(\Omega(G/H))$. It follows, for instance, that $\pi_\ast \HC(C^\ast_H(\ast; \Z)/C^\ast_G(\ast; \Z); C^\ast_K(\ast; \Z))$ is a completion of $\H^K_\ast(\Omega(G/H); \Z)$.
\end{remark}
\begin{example}\label{ex: hochschild and homology of loops Sn}
    Let us illustrate \cref{cor: loop homology E2}, or rather, the identification of the $\ku_H$-linear dual of $\cf_H(\Omega(G/H))$ with Hochschild cohomology, after base-changing along the map $\ku_H \to \ku \to \Z$ (to $\Z[1/2]$ in the second example) in two simple cases:
    \begin{enumerate}
        \item Let $H = \SU(n-1)\subseteq \SU(n) = G$. Then $G/H \simeq S^{2n-1}$, and so there is an isomorphism
        $$\pi_\ast \cf_H(\Omega(G/H))^\vee \otimes_{\ku_H} \Z \cong \H_\ast(\Omega S^{2n-1}; \Z) \cong \Z[y],$$
        where $y$ lives in weight $2n-2$. On the other hand, the map $\H^\ast_G(\ast; \Z) \to \H^\ast_H(\ast; \Z)$ identifies with the map 
        $$\Z[c_1, \cdots, c_n] \to \Z[c_1, \cdots, c_{n-1}]$$
        sending $c_n \mapsto 0$, where the $i$th Chern class $c_i$ lives in weight $-2i$. Taking Hochschild homology along this map identifies
        $$\HH(\sh^{1/2}\H^\ast_H(\ast; \Z)/\sh^{1/2}\H^\ast_G(\ast; \Z)) \simeq \sh^{1/2}\H^\ast_H(\ast; \Z) \otimes_\Z \HH(\Z/\sh^{1/2}\Z[c_n]).$$
        But $\pi_\ast \HH(\Z/\sh^{1/2} \Z[c_n])$ is isomorphic to the divided power algebra $\Z\pdb{\sigma^2(c_n)}$, where $\sigma$ denotes ``suspension'', so that $\sigma^2(c_n)$ lives in degree $-2n + 2$; it follows that there is an isomorphism
        $$\pi_\ast \HH(\sh^{1/2}\H^\ast_H(\ast; \Z)/\sh^{1/2}\H^\ast_G(\ast; \Z)) \cong \Z[c_1, \cdots, c_{n-1}]\pdb{\sigma^2(c_n)}.$$
        This in turn implies that there is an isomorphism
        $$\pi_\ast \HC(\sh^{1/2}\H^\ast_H(\ast; \Z)/\sh^{1/2}\H^\ast_G(\ast; \Z)) \cong \Z[c_1, \cdots, c_{n-1}]\pw{y}$$
        where the class $y$ in weight $2n-2$ is dual to $\sigma^2(c_n)$. Killing $c_1, \cdots, c_n$ (i.e., base-changing along $\sh^{1/2}\H^\ast_H(\ast; \Z) \to \Z$) precisely recovers a completion of $\H_\ast(\Omega S^{2n-1}; \Z)$. We will discuss a generalization (and decompletion) of this calculation in \cref{thm: homology of loops SV}.
        
        \item Let $H = \SO_{2n} \subseteq \SO_{2n+1} = G$ with $n>0$, and recall that we are replacing $\Z$ by $\Z' = \Z[1/2]$. Then $G/H \simeq S^{2n}$, and so a standard argument with the Serre spectral sequence shows that there is an isomorphism
        $$\pi_\ast \cf_H(\Omega(G/H))^\vee \otimes_{\ku_H} \Z' \cong \H_\ast(\Omega S^{2n}; \Z') \cong \Z'[y, z]/z^2,$$
        where $z$ lives in weight $2n-1$ and $y$ lives in weight $4n-2$. On the other hand, the map $\H^\ast_G(\ast; \Z') \to \H^\ast_H(\ast; \Z')$ identifies with the map 
        $$\Z'[p_1, \cdots, p_{n-1}, p_n] \to \Z'[p_1, \cdots, p_{n-1}, p_n^{1/2}]$$
        sending $p_n \mapsto (p_n^{1/2})^2$, where the $i$th Pontryagin class $p_i$ lives in weight $-4i$ and the Euler class $p_n^{1/2}$ lives in weight $-2n$. Taking Hochschild homology along this map identifies
        $$\HH(\sh^{1/2}\H^\ast_H(\ast; \Z')/\sh^{1/2}\H^\ast_G(\ast; \Z')) \simeq \sh^{1/2}\Z'[p_1, \cdots, p_{n-1}] \otimes_{\Z'} \HH(\sh^{1/2} \Z'[p_n^{1/2}]/\sh^{1/2} \Z'[p_n]),$$
        and so computing the Hochschild cohomology from \cref{cor: loop homology E2} amounts to computing the Hochschild cohomology $\HC(\Z'[p_n^{1/2}]/\Z'[p_n])$. \cref{lem: HC of verschiebung} implies that there is an isomorphism
        $$\pi_\ast \HC(\sh^{1/2}\H^\ast_H(\ast; \Z')/\sh^{1/2}\H^\ast_G(\ast; \Z')) \cong \Z'[p_1, \cdots, p_{n-1}, p_n^{1/2}]\pw{w}/p_n^{1/2} w,$$
        with $w$ in in weight $4n-2$. Upon killing $p_1, \cdots, p_{n-1}, p_n^{1/2}$ (i.e., base-changing along $\sh^{1/2}\H^\ast_H(\ast; \Z) \to \Z$), one precisely recovers a completion of $\H_\ast(\Omega S^{2n}; \Z')$. We will discuss a generalization (and decompletion) of this calculation in \cref{thm: homology of loops SV}.
    \end{enumerate}
    \begin{lemma}\label{lem: HC of verschiebung}
        Let $x$ be a class in homological degree $2n$, and let $j\geq 1$. Then there is an isomorphism
        $$\pi_\ast \HC(\Z[x]/\Z[x^j]) \cong \Z[x]\pw{w}/jx^{j-1}w,$$
        where $w$ lives in weight $-2nj-2$. 
    \end{lemma}
    \begin{proof}
        Since $\Z[x] = \sh \Z[x_n]$ where $x_n$ lives in weight $n$ and degree $0$, it suffices to work in the graded setting and assume that $x$ lives in weight $n$ and degree $0$.
        Let us first work in the ungraded setting; fix a nonconstant polynomial $g(x)\in \Z[x]$, and consider $\HC(\Z[x]/\Z[g])$. There is an isomorphism
        $$\Z[x] \otimes_{\Z[g]} \Z[x] \cong \Z[x, x']/(g(x) - g(x')) \cong \Z[x, z]/z f,$$
        where $z = x' - x$ and $f = \frac{g(x) - g(x + z)}{z}$. (If $x$ has weight $n$ and $g$ is homogeneous of degree $j$, the class $z$ lives in degree $0$ and weight $n$, and $f$ lives in degree $0$ and weight $n(j-1)$.)
        Our goal is to compute $\pi_\ast \End_{\Z[x,z]/zf}(\Z[x])$, where the map $\Z[x,z]/zf \to \Z[x]$ sends $z\mapsto 0$. There are several ways to compute this: one is to note that there is a presentation
        $$\Z[x] \simeq (\Z[x,z,u]\pdb{v}/(zf, u^2), d(u) = z, d(v) = uf)$$
        of $\Z[x]$ as a $\Z[x,z]/zf$-algebra. If $x$ has weight $n$ and $g$ is homogeneous of degree $j$, the class $u$ is in degree $1$ and weight $n$, and $v$ is a divided power class in degree $2$ and weight $nj$.
        This implies that there is an equivalence
        $$\End_{\Z[x,z]/zf}(\Z[x]) \simeq (\Z[x,u']\pw{w}/{u'}^2, d(u') = f(z=0) w),$$
        where $u'$ is dual to $u$ and $w$ is dual to $v$.
        If $x$ has weight $n$ and $g$ is homogeneous of degree $j$, the class $u'$ is in degree $-1$ and weight $-n$, and $w$ is in degree $-2$ and weight $-nj$.
        It follows that there is a class $w\in \pi_{-2} \End_{\Z[x,z]/zf}(\Z[x])$ such that $f(z=0)w = 0 \in \pi_{-2} \End_{\Z[x,z]/zf}(\Z[x])$, which gives an isomorphism
        $$\pi_\ast \HC(\Z[x]/\Z[g]) = \pi_\ast \End_{\Z[x,z]/zf}(\Z[x]) \cong \Z[x]\pw{w}/g'(0) w.$$
        If $x$ has weight $n$ and $g$ is homogeneous of degree $j$, the class $w$ lives in $\pi_{-2,-nj} \End_{\Z[x,z]/zf}(\Z[x])$, and we obtain a graded isomorphism
        $$\pi_\ast \HC(\Z[x]/\Z[x^j]) \cong \Z[x]\pw{w}/jx^{j-1}w,$$
        which gives the desired calculation by shearing.
    \end{proof}
\end{example}
\begin{remark}
    As in \cref{lem: HC of verschiebung}, one can also compute $\pi_\ast \HC(\pi_\ast \ku_{S^1}/\pi_\ast \ku_{\SU(2)})$ to obtain the following:
    $$\pi_\ast \HC(\pi_\ast \ku_{S^1}/\pi_\ast \ku_{\SU(2)}) \cong \Z[\beta, x, \tfrac{1}{1+\beta x}]\pw{w}/w(x - \ol{x}).$$
    Here, $w$ lives in degree $0$ and weight $2$, and $\ol{x} = -\frac{x}{1 + \beta x}$ is the negative of $x$ under the group law on $\GG_\beta = \spec \pi_\ast \ku_{S^1}$. When $\beta = 0$, this recovers \cref{lem: HC of verschiebung} for $j=2$ and $n = -1$.
    There is a spectral sequence whose $E_1$-page is $\pi_\ast \HC(\pi_\ast \ku_{S^1}/\pi_\ast \ku_{\SU(2)})$ which converges to $\pi_\ast \HC(\ku_{S^1}/\ku_{\SU(2)})$; this spectral sequence degenerates. Since the $2$-series of $x$ is $[2](x) = (1+\beta x)(x - \ol{x})$, we find that
    $$\pi_\ast \HC(\ku_{S^1}/\ku_{\SU(2)}) \cong \Z[\beta, x, \tfrac{1}{1+\beta x}]\pw{w}/w[2](x),$$
    where $w$ lives in degree $2$.
\end{remark}
\begin{remark}\label{rmk: general map complex reflection}
    The reader might observe that one can analyze $\HH(\pi_\ast \ku_H/\pi_\ast \ku_G)$ essentially using the combinatorics of the weight lattices and Weyl groups of $H$ and $G$. More generally, therefore, let $W_1 \to W_2$ be a homomorphism of finite groups acting on vector spaces $V_1 \to V_2$ over a field $k$ (possibly of nonzero characteristic). Then there is a map $V_1\mmod W_1 \to V_2\mmod W_2$, and hence one can consider the Hochschild homology $\HH(V_1\mmod W_1 / V_2\mmod W_2)$. This should be an interesting invariant associated to homomorphisms of finite groups, but it is likely only well-behaved if the map $V_1\mmod W_1 \to V_2\mmod W_2$ is an affine bundle.
    
    For instance, \cref{prop: equiv coh as hochschild} and \cref{rmk: expected G mod T} below imply that the Hochschild cohomology $\HC(\fr{t}/\fr{t}\mmod W)$ is closely related to a completion of the ``regular centralizer'' $\ld{\fr{t}}^\ast \times_{\ld{G}\backslash \ol{T^\ast(\ld{G}/\ld{N})}} \ld{\fr{t}}^\ast$ for $\ol{T^\ast(\ld{G}/\ld{N})}$. 
    More generally, if $\g$ is a simple Lie algebra (over $\cc$) with Weyl group $W$ and Cartan subalgebra $\fr{t}$, $V$ is a miniscule representation of $\g$, and $W_H$ is the stabilizer of a weight of $V$, one obtains a Gelfand pair $(W, W_H)$. The combinatorial properties of this Gelfand pair ``control'' the Hochschild cohomology $\HC(\fr{t}\mmod W_H/\fr{t}\mmod W)$. One therefore expects $\HC(\fr{t}\mmod W_H/\fr{t}\mmod W)$ to be an interesting combinatorial invariant of $V$.

    One example which does not come from Lie-theoretic data is the dihedral group $D_{2n}$ acting on $\AA^2 = \spec \cc[x_1, x_2]$ by the matrices $s = \begin{psmallmatrix}
        1 & 0\\
        0 & -1
    \end{psmallmatrix}$ and $r = \begin{psmallmatrix}
        \zeta_n & 0 \\
        0 & \zeta_n^{-1}
    \end{psmallmatrix}$. The algebra of invariants $\cc[x_1, x_2]^{D_{2n}}$ is simply $\cc[x_1 x_2, x_1^n + x_2^n]$. As in \cref{lem: HC of verschiebung}, one finds that
    $$\pi_\ast \HC(\AA^2 / \AA^2\mmod D_{2n}) \cong \cc[x_1, x_2]\pw{w_1, w_2}/\begin{psmallmatrix}
        x_1 & x_2 \\
        x_2^{n-1} & x_1^{n-1}
    \end{psmallmatrix}\vec{w}.$$
\end{remark}

\begin{example}
    Consider the subgroup $G^\mathrm{diag} \subseteq G \times G$, so that $(G \times G)/G^\mathrm{diag} \simeq G$ (this is the ``group case'' of \cref{ex: gp case}). Then \cref{prop: equiv coh as hochschild} says that there is an $S^1$-equivariant equivalence of $\Eoo$-$\ku_G$-algebras
    $$\cf_{G\times G}(\cL G) \simeq \HH(\ku_G/\ku_{G \times G}).$$
    By construction of equivariant $\ku$, there is an equivalence $\ku_{G \times G} \simeq \ku_G \otimes_{\ku} \ku_G$, so that the right-hand side can be identified with the factorization homology
    $$\HH(\ku_G/\ku_G \otimes_{\ku} \ku_G) \simeq \int_{S^2} \ku_G/\ku.$$
    Note that by \cref{prop: going from G of LG/H to H of Loops G/H}, the left-hand side can be identified with $\cf_G(\Omega G)$, so \cref{prop: equiv coh as hochschild} describes the $G$-equivariant $\ku$-\textit{co}homology of the affine Grassmannian:
    \begin{equation}\label{eq: ku-coh of GrG}
        \ku^\ast_G(\Omega G) \simeq \pi_\ast \int_{S^2} \ku_G/\ku.
    \end{equation}
    Note that there is an equivalence
    $$\HH(\ku_G/\ku_G \otimes_{\ku} \ku_G) \simeq \ku_G \otimes_\ku \HH(\ku/\ku_G).$$
    Upon killing the Bott class $\beta$, \cref{eq: ku-coh of GrG} implies that
    $$C^\ast_G(\Omega G; \Z) \simeq \int_{S^2} C^\ast_G(\ast; \Z)/\Z.$$
    As argued in \cite[Example A.8]{grg-reg}, this recovers \cite[Theorem 1]{bf-derived-satake} and \cite[Section 1.7]{ginzburg-langlands} upon rationalization.
\end{example}
\begin{remark}\label{rmk: S1rot and torus}
    Unlike \cref{prop: going from G of LG/H to H of Loops G/H}, \cref{prop: equiv coh as hochschild} gives an $S^1_\rot$-equivariant equivalence. In particular, it allows us to calculate the $S^1_\rot$-equivariant cohomology $\ku^\ast_{G \times G \times S^1_\rot}(\cL G) \simeq \ku^\ast_{G \times S^1_\rot}(\Omega G)$. We will discuss this in a future article, since addressing loop rotation in the detail it deserves will take us too far afield.

    However, since it is not very difficult to make explicit, let us explicate \cref{prop: equiv coh as hochschild} (or rather, its variant for Hochschild cohomology describing $\ku^{G \times S^1_\rot}_\ast(\Omega G)$) in the case when $G = T$ is a torus. (At the beginning of this section, we asked that $G$ be simply-connected; this is obviously not true for a torus, but that assumption was necessary only when doing computations with Hochschild (co)homology. We will \textit{not} use this perspective below.)
    
    As in \cite[Proposition 3.3.4]{grg-reg}, the associative graded ring $\ku^{T \times S^1_\rot}_\ast(\Omega T)$ can be identified with the algebra of $\GG_\beta$-differential operators on the dual torus $\ld{T}$. This is an analogue of the algebra of (asymptotic) differential operators. Let us assume for simplicity that $T$ is of rank $1$; then the algebra $\ku^{T \times S^1_\rot}_\ast(\Omega T)$ is the $F$-Weyl algebra $F\cd_{\square, \GG_m}$ of \cite[Definition 4.4.1]{generalized-n-series} for $F(x,y) = x + y + \beta xy$, at least up to completion. Explicitly, when $T = S^1$, we have
    $$\ku^{T \times S^1_\rot}_\ast(\Omega T) \cong \Z[\beta, \hbar, \tfrac{1}{1+\beta \hbar}]\{x, a^{\pm 1}\}[\tfrac{1}{1+\beta x}]/([x,a] = a \hbar(1 + \beta x)).$$
    Here, the curly brackets denotes the free associative algebra generated by the elements enclosed within. The classes $\hbar$ and $x$ live in weight $-2$ (they are the $S^1$-equivariant Chern classes for $\ku$), $\beta$ lives in weight $2$, and $a$ lives in weight zero. Let us note two specializations of this associative algebra:
    \begin{enumerate}
        \item If $\beta = 0$, the right-hand side above simply becomes $\Z[\hbar]\{x, a^{\pm 1}\}/([x,a] = \hbar a)$, which is precisely the algebra of asymptotic differential operators on $\ld{T} = \spec \Z[a^{\pm 1}]$ over $\Z$. Namely, $x = \hbar a\partial_a$; see \cite[Example 4.4.2]{generalized-n-series}.
        \item If $\beta$ is inverted, all elements can be pushed to degree zero. Namely, let $q = 1 + \beta \hbar$ and $\Theta = 1 + \beta x$. Then there is an isomorphism
        \begin{multline*}
            \Z[\beta^{\pm 1}, \hbar, \tfrac{1}{1+\beta \hbar}]\{x, a^{\pm 1}\}[\tfrac{1}{1+\beta x}]/([x,a] - a \hbar(1 + \beta x))
            \cong \Z[\beta^{\pm 1}, q^{\pm 1}]\{\Theta^{\pm 1}, a^{\pm 1}\}/(\Theta a - q a\Theta),
        \end{multline*}
        since
        \begin{align*}
            \Theta a & = (1 + \beta x) a = a + \beta x a = a + \beta a (x + \hbar + \beta \hbar x)\\
            & = a (1+\beta x)(1 + \beta \hbar) = qa\Theta.
        \end{align*}
        In particular, $\ku^{T \times S^1_\rot}_\ast(\Omega T)[1/\beta] = \KU^{T \times S^1_\rot}_\ast(\Omega T)$ can be identified with the $q$-Weyl algebra of $\ld{T} = \spec \Z[a^{\pm 1}]$. Namely, $\Theta = q^{a\partial_a}$; see \cite[Example 4.4.3]{generalized-n-series}.
    \end{enumerate}
    In general, $\ku^{T \times S^1_\rot}_\ast(\Omega T)$ interpolates between the algebra of asymptotic differential and $q$-difference operators on $\ld{T}$.
\end{remark}

%% file: equivalences/sheaves-equivalences.tex
\subsection{Using the regular centralizer}\label{subsec: using the regular centralizer}

Again, recall \cref{notation: compact vs complexification}: if $G$ is a compact Lie group, we will write $G\ls{t}$ or $G\pw{t}$ below to mean $G_\cc\ls{t}$ or $G_\cc\pw{t}$, respectively.
\begin{definition}\label{def: Shv-Sat LG/H}
    Let $\Shv^c_{G\pw{t}}(G\ls{t}/H\ls{t}; \QQ)$ denote the $\infty$-category of $G$-equivariant sheaves of $\QQ$-modules on $G\ls{t}/H\ls{t}$ which are constructible for the orbit stratification on $G\ls{t}/H\ls{t}$. 
    There is a natural left-action of the $\E{3}$-monoidal $\infty$-category $\Shv^c_{(G\times G)\pw{t}}(G\ls{t}; \QQ)$ on $\Shv^c_{G\pw{t}}(G\ls{t}/H\ls{t}; \QQ)$, and in particular, a left-action of $\Rep(\ld{G})$ by \cref{thm: abelian satake}.
    Let $\IC_0\in \Shv^c_{G\pw{t}}(G\ls{t}/H\ls{t}; \QQ)$ denote the pushforward $i_! \ul{\QQ}$ of the constant sheaf along the inclusion $i: G\pw{t}/H\pw{t} \hookrightarrow G\ls{t}/H\ls{t}$ of the constant loops.
    Let $\Shv^{c,\Sat}_{G\pw{t}}(G\ls{t}/H\ls{t}; \QQ)$ denote the full subcategory of $\Shv^c_{G\pw{t}}(G\ls{t}/H\ls{t}; \QQ)$ generated by $\IC_0$ under the action of $\Rep(\ld{G})$.
\end{definition}
\begin{example}
    In the group case of \cref{ex: gp case}, the $\infty$-category $\Shv^{c,\Sat}_{(G\times G)\pw{t}}((G\times G)\ls{t}/G^\diag\ls{t}; \QQ)$ from \cref{def: Shv-Sat LG/H} agrees with \cref{def: Shv-Sat GrG}.
\end{example}
We can now state our main criterion for proving equivalences of the form \cref{conj: bzsv}.
We will now use \textit{all} the hypotheses of \cref{hypothesis: rank 1 weakly placid}.
\begin{remark}\label{rmk: meaning of the hypothesis}
    Before proceeding to the argument, let us make one comment about \cref{hypothesis: rank 1 weakly placid}: as the terminology suggests, a subgroup being optimal is a rather idealized situation; it is essentially the smallest set of hypotheses needed to make the argument of \cref{thm: derived satake} go through.
    It is, of course, possible that examples of interest (even the ones considered in this article) do not satisfy this condition. Verifying \cref{hypothesis: rank 1 weakly placid} is likely no easy task, and working out the appropriate subtleties of the sheaf theories involved in \cref{conj: bzsv} will be very important to understanding microlocal aspects of relative geometric Langlands.
    However, the point of \cref{hypothesis: rank 1 weakly placid} is to isolate some hard sheaf-theoretic subtleties, contingent upon which we may state \cref{thm: ordinary homology criterion satake}, whose criteria (we believe) extract the key ways in which (relative) Langlands duality is born. This has the effect of giving psychologically more manageable conjectures at ``category level $0$''.

    The first two assumptions of \cref{hypothesis: rank 1 weakly placid} go into proving the formality of the graded (derived) algebra $\Hom_{\Shv^c_{G\pw{t}}(G\ls{t}/H\ls{t}; \QQ)}(\IC_0, \IC_0 \star \ld{\cR})$. In the case when $G/H$ is $\GL_n/\O_n$ or $\GL_{2n}/\Sp_{2n}$, this formality was proved as \cite[Theorem 23]{chen-yi-formality}. The remaining assumption in \cref{hypothesis: rank 1 weakly placid} is concerned with the applicability of \cref{thm: full faithful} to the present situation.
\end{remark}

\begin{theorem}\label{thm: ordinary homology criterion satake}
    Let $H_\cc\subseteq G_\cc$ be a closed connected reductive subgroup which is optimal in the sense of \cref{hypothesis: rank 1 weakly placid}.
    Let $\ld{M}$ denote a ``dual'' affine graded $\ld{G}$-space (as prescribed by \cite{bzsv} if $H_\cc$ is a spherical subgroup).

    Suppose that there is a ``Kostant section'' $\kappa_{\ld{M}}: \ld{\fr{h}}^\ast(2)\mmod \ld{H} \hookrightarrow \ld{M}$ (see, e.g., \cref{conj: generalized kostant slice}) such that:
    \begin{enumerate}
        \item Let $\ld{J}_{X}'$ denote the (possibly non-flat) group scheme $\ld{\fr{h}}^\ast(2)\mmod \ld{H} \times_{\ld{M}/\ld{G}(-2\rho)} \ld{\fr{h}}^\ast(2)\mmod \ld{H}$ over $\ld{\fr{h}}^\ast(2)\mmod \ld{H}$. Then the algebra of regular functions on $(\ld{\fr{h}}^\ast(2)\mmod \ld{H} \times \ld{G}(-2\rho))/\ld{J}_{X}'$ is isomorphic to $\co_{\ld{M}}$.
        (For instance, this holds by the algebraic Hartogs lemma if $\ld{M}$ is normal and the $\ld{G}$-orbit $\ld{M}^\reg$ of the image of $\kappa_{\ld{M}}$ is open with complement of codimension $\geq 2$.)
        \item Define
        \begin{equation}\label{eq: JX def}
            \ld{J}_X := \spec \H^H_\ast(\Omega (G/H); \QQ).
        \end{equation} 
        There is an isomorphism of graded group schemes over $\ld{\fr{h}}^\ast(2)\mmod \ld{H}$:
        $$\ld{J}_X \cong \ld{\fr{h}}^\ast(2)\mmod \ld{H} \times_{\ld{M}/\ld{G}(-2\rho)} \ld{\fr{h}}^\ast(2)\mmod \ld{H} = \ld{J}_X'.$$
    \end{enumerate}
    Then there is an equivalence of $\QQ$-linear $\infty$-categories
    \begin{equation}\label{eq: display rel sat criterion}
        \Shv^{c,\Sat}_{G\pw{t}}(G\ls{t}/H\ls{t}; \QQ) \simeq \Perf(\sh^{1/2}\ld{M}/\ld{G}(-2\rho)).
    \end{equation}
    Moreover, this equivalence fits into a commutative diagram
    $$\xymatrix{
    \Shv^{c,\Sat}_{G\pw{t}}(G\ls{t}/H\ls{t}; \QQ) \ar[r]^\sim \ar[d]^-{\mathrm{cohomology}} & \Perf(\sh^{1/2}\ld{M}/\ld{G}(-2\rho)) \ar[d]^-{\kappa^\ast} \\
    \Shv_H(\ast; \QQ) \ar[r]_-\sim & \Perf(\sh^{1/2} \ld{\fr{h}}^\ast(2)\mmod \ld{H}),
    }$$
    where the cohomology functor $\Shv^{c,\Sat}_{G\pw{t}}(G\ls{t}/H\ls{t}; \QQ) \to \Shv_G(\ast; \QQ)$ factors through the canonical functor $\Shv_H(\ast; \QQ) \to \Shv_G(\ast; \QQ)$ by \cref{eq: factorization of cohomology}.

    Let $\mu: \ld{M}/\ld{G}(-2\rho) \to \ld{\g}^\ast(2\rho - 2)/\ld{G}(-2\rho)$ denote the moment map, and assume that there is a commutative diagram
    $$\xymatrix{
    \ld{\fr{h}}^\ast(2)\mmod \ld{H} \ar[r]^-{\kappa_{\ld{M}}} \ar[d] & \ld{M}/\ld{G}(-2\rho) \ar[d]^-\mu \\
    \ld{\g}^\ast(2)\mmod \ld{G} \ar[r]^-{\kappa} & \ld{\g}^\ast(2 - 2\rho)/\ld{G}(-2\rho),
    }$$
    so that there is an induced map
    $$\ld{\fr{h}}^\ast(2)\mmod \ld{H} \times_{\ld{M}/\ld{G}(-2\rho)} \ld{\fr{h}}^\ast(2)\mmod \ld{H} \cong \ld{J}_X \to \ld{\g}^\ast(2)\mmod \ld{G} \times_{\ld{\g}^\ast(2 - 2\rho)/\ld{G}(-2\rho)} \ld{\g}^\ast(2)\mmod \ld{G} \cong \ld{J}.$$
    If the isomorphism of (b) fits into a commutative diagram
    $$\xymatrix{
    \spec \H_\ast^H(\Omega(G/H); \QQ) \ar[rr]^-\sim_-{\textrm{(b)}} \ar[d] & & \ld{\fr{h}}^\ast(2)\mmod \ld{H} \times_{\ld{M}/\ld{G}(-2\rho)} \ld{\fr{h}}^\ast(2)\mmod \ld{H} \ar[d] \\
    \spec \H_\ast^G(\Omega G; \QQ) \ar[rr]^-\sim_-{\textrm{\cref{thm: derived satake}}} & & \ld{\g}^\ast(2)\mmod \ld{G} \times_{\ld{\g}^\ast(2 - 2\rho)/\ld{G}(-2\rho)} \ld{\g}^\ast(2)\mmod \ld{G},
    }$$
    then the equivalence \cref{eq: display rel sat criterion} is equivariant for the left-action of $\Shv_{(G \times G)\pw{t}}^{c,\Sat}(G\ls{t}; \QQ) \simeq \Perf(\sh^{1/2} \ld{\g}^\ast(2 - 2\rho)/\ld{G}(-2\rho))$ via \cref{thm: derived satake}.
\end{theorem}
\begin{proof}
    The proof will follow the first half of the proof of \cref{thm: derived satake}.
    Let $\cC$ denote $\Shv^{c,\Sat}_{G\pw{t}}(G\ls{t}/H\ls{t}; \QQ)$, so that $\cC$ admits a left-action of $\Shv^{c,\Sat}_{G\pw{t}}(\Gr_G; \QQ) \simeq \Shv^{c,\Sat}_{(G\times G)\pw{t}}(G\ls{t}; \QQ)$. In particular, \cref{thm: abelian satake} implies that $\cC$ admits a left-action of $\Rep(\ld{G})$.
    Let $\tilde{\cC}$ denote the base-change $\cC \otimes_{\Rep(\ld{G})} \Mod_\QQ$, so that $\IC_0$ is a compact generator of $\tilde{\cC}$ (by definition of $\cC$).
    It follows from the Barr-Beck theorem \cite[Theorem 4.7.3.5]{HA} that there is an equivalence $\Phi: \tilde{\cC} \xar{\sim} \Perf_{\End_{\tilde{\cC}}(\IC_0)}$, implemented by the functor $\Hom_{\tilde{\cC}}(\IC_0, -)$.
    Recall that $\ld{\cR}\star \IC_0 \in \Shv^{c,\Sat}_{G\pw{t}}(G\ls{t}/H\ls{t}; \QQ)$ denotes the sheaf obtained by the action of $\co_{\ld{G}}\in \Rep(\ld{G})$ on $\IC_0$. By definition of $\tilde{\cC}$, we can identify $\End_{\tilde{\cC}}(\IC_0) \simeq \Hom_{\cC}(\IC_0, \IC_0 \star \ld{\cR})$.

    The same argument as in the proof of \cref{thm: derived satake} shows that $\End_{\tilde{\cC}}(\IC_0)$ is formal. The key input needed is the hypothesis that if $q$ is a sufficiently large prime number, the objects 
    $$\IC_0, \IC_0 \star \ld{\cR} \in \Shv^{c,\et}_{G(\ol{\FF_q}\pw{t})}(G(\ol{\FF_q}\ls{t})/H(\ol{\FF_q}\ls{t}))$$
    are pure of weight zero. This allows us to use \cite[Lemma 3.1.5]{bezr-yun-koszul} and \cref{thm: changing coefficients and bases} to obtain the formality of $\End_{\tilde{\cC}}(\IC_0)$. It follows that
    \begin{equation}\label{eq: Endtilde C and sh1/2}
        \End_{\tilde{\cC}}(\IC_0) \simeq \sh^{1/2}(\Ext^\bull_{\Shv^c_{G\pw{t}}(G\ls{t}/H\ls{t}; \QQ)}(\IC_0, \IC_0 \star \ld{\cR})).
    \end{equation}

    To compute this $\Ext$-algebra, we will use \cref{thm: full faithful}. One can show that the assumptions of \cref{thm: full faithful} are satisfied for $\IC_0$ and $\IC_0 \star \ld{\cR}$, so the cited result gives a graded isomorphism
    \begin{multline*}
        \Ext^\bull_{\Shv^c_{G\pw{t}}(G\ls{t}/H\ls{t}; \QQ)}(\IC_0, \IC_0 \star \ld{\cR}) \\
        \cong \Hom^\bull_{\H^\ast_{G\pw{t}}(G\ls{t}/H\ls{t}; \QQ)}(\H^\ast_{G\pw{t}}(G\ls{t}/H\ls{t}; \IC_0), \H^\ast_{G\pw{t}}(G\ls{t}/H\ls{t}; \IC_0 \star \ld{\cR})).
    \end{multline*}
    There is an isomorphism
    $$\H^\ast_{G\pw{t}}(G\ls{t}/H\ls{t}; \IC_0) \cong \H^\ast_G(G/H; \QQ) \cong \H^\ast_H(\ast; \QQ) \cong \co_{\ld{\fr{h}}^\ast(2)\mmod \ld{H}},$$
    and hence an isomorphism
    \begin{align*}
        \H^\ast_{G\pw{t}}(G\ls{t}/H\ls{t}; \IC_0 \star \ld{\cR}) & \cong \H^\ast_{G\pw{t}}(G\ls{t}/H\ls{t}; \IC_0) \otimes_{\H^\ast_{(G \times G)\pw{t}}(G\ls{t}; \IC_0)} \H^\ast_{(G\times G)\pw{t}}(G\ls{t}; \IC_0 \star \ld{\cR}) \\
        & \cong \co_{\ld{\fr{h}}^\ast(2)\mmod \ld{H}} \otimes_{\co_{\ld{\g}^\ast(2)\mmod \ld{G}}} \co_{\ld{\g}^\ast(2)\mmod \ld{G} \times \ld{G}(-2\rho)} \\
        & \cong \co_{\ld{\fr{h}}^\ast(2)\mmod \ld{H} \times \ld{G}(-2\rho)}.
    \end{align*}
    Moreover, there is an isomorphism $\H^\ast_{G}(\cL (G/H); \QQ) \cong \H^\ast_H(\Omega(G/H); \QQ)$ by \cref{rmk: cohomology going from G of LG/H to H of Loops G/H}, so we find that there is a graded isomorphism
    \begin{multline*}
        \Ext^\bull_{\Shv^c_{G\pw{t}}(G\ls{t}/H\ls{t}; \QQ)}(\IC_0, \IC_0 \star \ld{\cR}) 
        \cong \Hom^\bull_{\H^\ast_{H}(\Omega (G/H); \QQ)}(\co_{\ld{\fr{h}}^\ast(2)\mmod \ld{H}}, \co_{\ld{\fr{h}}^\ast(2)\mmod \ld{H} \times \ld{G}(-2\rho)})\\
        \cong \Hom^\bull_{\H^\ast_{H}(\ast; \QQ)}(\co_{\ld{\fr{h}}^\ast(2)\mmod \ld{H}}, \co_{\ld{\fr{h}}^\ast(2)\mmod \ld{H} \times \ld{G}(-2\rho)})^{\spec \H^H_\ast(\Omega (G/H); \QQ)} 
        \cong \co_{\ld{\fr{h}}^\ast(2)\mmod \ld{H} \times \ld{G}(-2\rho)}^{\spec \H^H_\ast(\Omega (G/H); \QQ)}.
    \end{multline*}
    By (b), there is an isomorphism $\spec \H^H_\ast(\Omega (G/H); \QQ) = \ld{J}_X \cong \ld{J}_X'$, and hence there is an isomorphism
    $$\Ext^\bull_{\Shv^c_{G\pw{t}}(G\ls{t}/H\ls{t}; \QQ)}(\IC_0, \IC_0 \star \ld{\cR}) \cong \co_{(\ld{\fr{h}}^\ast(2)\mmod \ld{H} \times \ld{G}(-2\rho))/\ld{J}_X'}.$$
    By (a), this is further isomorphic to $\co_{\ld{M}}$. Using \cref{eq: Endtilde C and sh1/2}, it follows that $\tilde{\cC}$ is equivalent to the $\infty$-category $\Perf_{\sh^{1/2}\co_{\ld{M}}} \simeq \Perf(\sh^{1/2} (\ld{M}))$. This in turn implies that $\cC$ is equivalent to the $\infty$-category $\Perf(\sh^{1/2} \ld{M}/\ld{G}(-2\rho))$, as desired.
\end{proof}

\begin{remark}
    In this article, we will only focus on applying \cref{thm: ordinary homology criterion satake} to the case when $G_\cc/H_\cc$ is a spherical $G_\cc$-variety of rank $1$. However, in \cite{triple-product-cayley}, we use \cref{thm: ordinary homology criterion satake} to study \cref{conj: bzsv} for the spherical subgroups $\PGL_2^\mathrm{diag} \subseteq \PGL_2^{\times 3}$ and $\G_2 \subseteq \SO_8/\mu_2$ of relative rank $3$. In \cite{quat-satake}, \cref{thm: ordinary homology criterion satake} was used in the case when $G_\cc/H_\cc = \GL_{2n}/\Sp_{2n}$ to show that the dual quotient stack $\ld{M}/\GL_{2n}(-2\rho)$ identifies with $\gl_n(4-4\rho)/\GL_n(-4\rho)$. This relies on the following isomorphism of group schemes over $\spec \H_{\Sp_{2n}}^\ast(\ast; \QQ) \cong \fr{t}^n(4)\mmod \Sigma_n$:
    \begin{equation}\label{eq: sp2n-equiv homology of symplectic grassmannian}
        \spec \H^{\Sp_{2n}}_\ast(\Omega (\GL_{2n}/\Sp_{2n}); \QQ) \cong \fr{t}^n(4)\mmod \Sigma_n \times_{\gl_n^\ast(4-4\rho)/\GL_n(-4\rho)} \fr{t}^n(4)\mmod \Sigma_n.
    \end{equation}
    Here, the map $\fr{t}^n(4)\mmod \Sigma_n \to \gl_n^\ast(4-4\rho)/\GL_n(-4\rho)$ is given by the Kostant slice. 
    However, more is true: \cref{eq: sp2n-equiv homology of symplectic grassmannian} can be refined to an isomorphism
    $$\spec \H^{\Sp_2^n}_\ast(\Omega (\GL_{2n}/\Sp_{2n}); \QQ) \cong \fr{t}^n(4) \times_{\tilde{\gl_n}(4-4\rho)/\GL_n(-4\rho)} \fr{t}^n(4)$$
    of graded group schemes over $\spec \H_{\Sp_2^n}^\ast(\ast; \QQ) \cong \fr{t}^n(4)$. Here, the map $\fr{t}^n(4) \to \tilde{\gl_n}$ is the Kostant slice for the Grothendieck-Springer resolution of $\gl_n$. (See also \cref{rmk: K-equiv-homology} for an interpretation via Hochschild cohomology.)
    
    In future work, we will study some more exotic examples. For instance, we will see that if \cref{hypothesis: rank 1 weakly placid} is satisfied in the case of the spherical $\mathrm{E}_6$-variety $\mathrm{E}_6/\F_4$ (where both $\mathrm{E}_6$ and $\F_4$ denote the simply-connected forms), then there is an equivalence
    $$\Shv^{c,\Sat}_{\mathrm{E}_6\pw{t}}(\mathrm{E}_6\ls{t}/\F_4\ls{t}; \QQ) \simeq \Perf(\fr{pgl}_3^\ast[8-8\rho]/\PGL_3[-8\rho] \times \g_2^\ast[2]\mmod \G_2).$$
    In other words, if $\ld{M}$ is the Hamiltonian $\mathrm{E}_6$-variety which is dual to $\mathrm{E}_6/\F_4$, then there is an isomorphism of stacks
    \begin{equation}\label{eq: M for E6 mod F4}
        \ld{M}/\mathrm{E}_6 \cong \fr{pgl}_3^\ast/\PGL_3 \times \g_2^\ast\mmod \G_2.
    \end{equation}
    This is shown by computing an isomorphism of graded group schemes over $\spec \H^\ast_{\F_4}(\ast; \QQ) \cong \fr{pgl}_3^\ast(8) \mmod \PGL_3 \times \g_2^\ast(2)\mmod \G_2$:
    $$\spec \H^{\F_4}_\ast(\Omega(\mathrm{E}_6/\F_4); \QQ) \cong \left( \fr{pgl}_3^\ast(8)\mmod \PGL_3 \times_{\fr{pgl}_3^\ast(8-8\rho)/\PGL_3(-8\rho)} \fr{pgl}_3^\ast(8) \mmod \PGL_3\right) \times \g_2^\ast(2)\mmod \G_2.$$
    In fact, more is true: there is an isomorphism
    $$\spec \H^{\Spin_8}_\ast(\Omega(\mathrm{E}_6/\F_4); \QQ) \cong \left(\AA^2(8,8) \times_{\tilde{\fr{pgl}_3}(8-8\rho)/\PGL_3(-8\rho)} \AA^2(8,8)\right) \times \g_2^\ast(2)\mmod \G_2$$
    of graded group schemes over $\spec \H^\ast_{\Spin_8}(\ast; \QQ) \cong \AA^2(8,8) \times \g_2^\ast(2)\mmod \G_2$ (see also \cref{rmk: K-equiv-homology} for an interpretation via Hochschild cohomology).
\end{remark}

Note the following consequence of the argument of \cref{thm: ordinary homology criterion satake}.
\begin{lemma}\label{lem: hom to G^ x h mod H}
    There is a homomorphism of graded group schemes
    $$\spec \H^H_\ast(\Omega (G/H); \QQ) \to \ld{G}(-2\rho) \times \ld{\fr{h}}^\ast(2)\mmod \ld{H}$$
    over $\ld{\fr{h}}^\ast(2)\mmod \ld{H} \cong \spec \H^\ast_H(\ast; \QQ)$.
\end{lemma}
\begin{remark}\label{rmk: a composite map in homology}
    One can see \cref{lem: hom to G^ x h mod H} more directly as follows.
    \cref{thm: classical homology loops G} and \cref{cor: classifying stack of J} together imply that there is a homomorphism
    $$\spec \H^G_\ast(\Omega G; \QQ) \cong \ld{J} \to \ld{G}(-2\rho) \times \ld{\g}^\ast(2)\mmod \ld{G}$$
    of graded group schemes over $\ld{\g}^\ast(2)\mmod \ld{G} \cong \spec \H^\ast_G(\ast; \QQ)$. Base-changing along the map $\H^\ast_G(\ast; \QQ) \to \H^\ast_H(\ast; \QQ)$, we obtain a homomorphism
    $$\spec \H^H_\ast(\Omega G; \QQ) \to \ld{G}(-2\rho) \times \ld{\fr{h}}^\ast(2)\mmod \ld{H}.$$
    Composition with the map $\spec \H^H_\ast(\Omega (G/H); \QQ) \to \spec \H^H_\ast(\Omega G; \QQ)$ induced by the map $G \to G/H$ produces the desired homomorphism.

    Note that this leads to the following description of $\ld{M}$ (ignoring gradings for simplicity), which amounts to turning \cref{thm: ordinary homology criterion satake} on its head. Namely, recall that we defined $\ld{J}_X = \spec \H^H_\ast(\Omega (G/H); \QQ)$.
    Then, \cref{lem: hom to G^ x h mod H} gives a map $B_{\ld{\fr{h}}^\ast\mmod \ld{H}} \ld{J}_X \to B_{\ld{\fr{h}}^\ast\mmod \ld{H}} (\ld{J} \times_{\ld{\g}^\ast\mmod \ld{G}} \ld{\fr{h}}^\ast\mmod \ld{H})$ of stacks over $\ld{\fr{h}}^\ast\mmod \ld{H}$. Moreover, the target identifies with $\ld{\g}^{\ast,\reg}/\ld{G} \times_{\ld{\g}^\ast\mmod \ld{G}} \ld{\fr{h}}^\ast\mmod \ld{H}$ by \cref{cor: classifying stack of J}. In particular, there is a canonical composite
    $$B_{\ld{\fr{h}}^\ast\mmod \ld{H}} \ld{J}_X \to B_{\ld{\g}^\ast\mmod \ld{G}} \ld{J} \cong \ld{\g}^\reg/\ld{G} \to B\ld{G}.$$
    One can then identify $\ld{M}$ with the affine closure of $B_{\ld{\fr{h}}^\ast\mmod \ld{H}} \ld{J}_X \times_{B\ld{G}} \spec(\QQ)$.
    Said differently, $\ld{M}/\ld{G}$ can be thought of as the ``affine closure of $B_{\ld{\fr{h}}^\ast\mmod \ld{H}} \ld{J}_X$ relative to $B\ld{G}$''.

    Observe that one can \textit{define} some scheme $\ld{M}^\ddag$ in this way even if $G/H$ is not spherical (one just needs $H$ to be reductive; in fact, the below definition can be made for \textit{any} homomorphism $H \to G$) 
    \begin{equation}\label{eq: M as G mod JX}
        {\ld{M}}^\ddag = \ol{(\ld{G} \times \ld{\fr{h}}^\ast\mmod \ld{H})/\ld{J}_X}.
    \end{equation}
    One generally has very little control of $\ld{M}^\ddag$ if it is constructed via \cref{eq: M as G mod JX}. However, one can at the very least use \cref{eq: BJ span} to obtain a Lagrangian morphism $\ld{M}^\ddag/\ld{G} \to \ld{\g}^\ast/\ld{G}$, which by \cref{prop: safronov lag maps} can be understood as equipping $\ld{M}^\ddag$ with the structure of a Hamiltonian $\ld{G}$-space.
    
    The definition of $\ld{M}^\ddag$ is rigged so that if \cref{hypothesis: rank 1 weakly placid} is satisfied, there is an equivalence of categories
    $$\Shv^{c,\Sat}_{G\pw{t}}(G\ls{t}/H\ls{t}; \QQ) \simeq \Perf(\sh^{1/2}\ld{M}^\ddag/\ld{G}(-2\rho)).$$
    \begin{notno}
        We will write $\ld{M}^\ddag$ to denote the scheme of \cref{eq: M as G mod JX}, and when $H\subseteq G$ is spherical, we will write $\ld{M}$ (without the $\ddag$) to denote the scheme defined in \cite[Section 4]{bzsv}. The question of proving \cref{conj: bzsv} (when $H\subseteq G$ is spherical) becomes about identifying $\ld{M}^\ddag$ with $\ld{M}$. In \cref{sec: case by case}, we verify this for rank $1$ homogeneous spherical varieties by explicit calculation. In the general case, understanding \cref{conj: JX and GcheckX} should be an important step in this process.
    \end{notno}

    The perspective on $\ld{M}$ as being $\ld{M}^\ddag$ (defined by \cref{eq: M as G mod JX}) leads to several interesting and nontrivial structures. For instance, $G/H$ has an action of its $G$-equivariant automorphism group $\N_G(H)/H$, and hence $\N_G(H)/H$ acts on $\ld{J}_X$. The above construction of $\ld{M}^\ddag$ therefore shows that there is a natural $\N_G(H)/H$-action on $\ld{M}^\ddag$, and hence an expected $\N_G(H)/H$-action on $\ld{M}$, which commutes with its Hamiltonian $\ld{G}$-action. This action is highly interesting; for instance, when $H = T \subseteq G$, we check in \cref{rmk: expected G mod T} that $\ld{M}^\ddag$ is the affine closure of $T^\ast(\ld{G}/\ld{N})$; the above action of $\N_G(T)/T \cong W$ turns out to be Gelfand-Graev action of the Weyl group on $\ol{T^\ast(\ld{G}/\ld{N})}$ (as described by Ginzburg-Kazhdan in \cite{ginzburg-kazhdan}).
\end{remark}

The following asks for a refinement of \cref{lem: hom to G^ x h mod H}.
\begin{conjecture}\label{conj: JX and GcheckX}
    Suppose $G/H$ is a spherical $G$-variety.
    The homomorphism of \cref{lem: hom to G^ x h mod H} fits into a commutative diagram
    $$\xymatrix{
    \spec \H^H_\ast(\Omega (G/H); \QQ) \ar[r] \ar[d] & \spec \H^H_\ast(\Omega G; \QQ) \ar[d] \ar[r]^-\sim & \ld{J}_{\ld{G}} \times_{\ld{\g}^\ast(2)\mmod \ld{G}} \ld{\fr{h}}^\ast(2)\mmod \ld{H}\\
    \ld{G}_X(-2\rho_{\ld{G}}) \times \ld{\fr{h}}^\ast(2)\mmod \ld{H} \ar[r] & \ld{G}(-2\rho_{\ld{G}}) \times \ld{\fr{h}}^\ast(2)\mmod \ld{H} &
    }$$
    of graded group schemes over $\ld{\fr{h}}^\ast(2)\mmod \ld{H}$, where the homomorphism $\ld{G}_X \to \ld{G}$ is that of \cref{def: SV dual group}, and the vertical maps are closed immersions.
\end{conjecture}
\begin{remark}
    \cref{conj: JX and GcheckX} should in some sense follow from the \textit{abelian} Satake equivalence of \cite{gaitsgory-nadler} via the Tannakian formalism (as in \cite{homology-langlands}). Namely, let $\QQ(Z)$ denote the tensor category studied in \cite{gaitsgory-nadler}. If \cref{hypothesis: rank 1 weakly placid}(a) is satisfied (for example, $G/H$ is a symmetric variety for $G$), \cref{prop: going from G of LG/H to H of Loops G/H} shows that equivariant homology defines a functor from $\QQ(Z)$ to the abelian $1$-category $\coMod_{\H^H_\ast(\Omega (G/H); \QQ)}(\QCoh(\ld{\fr{h}}^\ast(2)\mmod \ld{H}))$. 
    There is a symmetric monoidal equivalence $\QQ(Z) \simeq \Rep(\ld{G}_{X,GN})$ by \cite{gaitsgory-nadler}, where $\ld{G}_{X,GN}$ is the Gaitsgory-Nadler dual group. If $\ld{G}_{X,GN} \cong \ld{G}_X$, and there is an analogue of \cite[Lemma 2.2]{homology-langlands} in this context, the Tannakian formalism would give a homomorphism $\spec \H^H_\ast(\Omega (G/H); \QQ) \to \ld{G}_X \times \ld{\fr{h}}^\ast(2)\mmod \ld{H}$. The desired analogue of \cite[Lemma 2.2]{homology-langlands} is closely related to the other items in \cref{hypothesis: rank 1 weakly placid}.
\end{remark}
\begin{remark}
    Note that \cref{conj: JX and GcheckX} allows one to consider the affine $k$-scheme $V^\ddag_X$ defined by 
    $$\co_{V^\ddag_X} := \co_{\ld{\fr{h}}^\ast(2)\mmod \ld{H} \times \ld{G}_X(-2\rho_{\ld{G}})}^{\spec \H^H_\ast(\Omega (G/H); \QQ)};$$
    this satisfies the property that $\ld{M}^\ddag \cong \Ind^{\ld{G}}_{\ld{G}_X} V^\ddag_X$.
    As in the proof of \cref{thm: ordinary homology criterion satake}, $\co_{V^\ddag_X}$ can be identified with the graded $\Ext$-algebra $\Ext^\bull_{\Shv^c_{G\pw{t}}(G\ls{t}/H\ls{t}; \QQ)}(\IC_0, \IC_0 \star \ld{\cR}_X)$, where $\ld{\cR}_X$ is the ``regular sheaf'' on $\Shv^c_{G\pw{t}}(G\ls{t}/H\ls{t}; \QQ)$ corresponding to $\co_{\ld{G}_X}$. Let $R_X := \Map_{\Shv^c_{G\pw{t}}(G\ls{t}/H\ls{t}; \QQ)}(\IC_0, \IC_0 \star \ld{\cR}_X)$. If $R_X$ is formal, one would be able to identify the full subcategory of $\Shv^c_{G\pw{t}}(G\ls{t}/H\ls{t}; \QQ)$ generated by $\IC_0$ under the action of $\co_{\ld{G}_X}$ with $\Perf(\sh^{1/2} V^\ddag_X/\ld{G}_X)$.
    
    In particular, were \cref{conj: bzsv} to hold, \cref{lem: wind and ind} would tell us that $V^\ddag_X$ is a graded vector space with $\ld{G}_X$-action. If one somehow knew this \textit{a priori}, and also that $R_X$ admits the structure of an $\E{2}$-$k$-algebra, then \cref{lem: formality polynomial} would \textit{automatically} imply that $R_X$ is formal as an $\E{1}$-$k$-algebra (since $\pi_\ast R_X = \co_{V^\ddag_X}$ would be a polynomial algebra with generators in even weights). This would sidestep needing to use parts of \cref{hypothesis: rank 1 weakly placid} to prove placidity. Morally, this is how one proves \cref{thm: rk 1 bzsv is true}; but since we do not actually directly show (even in the rank $1$ cases) that $R_X$ admits an $\E{2}$-$k$-algebra structure, the actual approach taken here is to assume \cref{hypothesis: rank 1 weakly placid} for these examples in order to deduce formality.
\end{remark}

\begin{remark}\label{rmk: reg centr and endoscopy}
    The group scheme $\spec \H^H_\ast(\Omega(G/H); \QQ)$ can be described in terms of the regular centralizer group schemes $\ld{J}_{\ld{G}}$ and $\ld{J}_{\ld{H}}$ for $\ld{G}$ and $\ld{H}$. There is a fiber sequence of $\E{1}$-spaces
    $$\Omega H \to \Omega G \to \Omega (G/H),$$
    which gives an equivalence
    $$C^H_\ast(\Omega(G/H); \QQ) \simeq C^H_\ast(\Omega G; \QQ) \otimes_{C^H_\ast(\Omega H; \QQ)} C_H^\ast(\ast; \QQ)$$
    of $\E{1}$-$\QQ$-algebras. It follows from \cref{thm: classical homology loops G}, for instance, that if the map $C^H_\ast(\Omega H; \QQ) \to C^H_\ast(\Omega G; \QQ)$ is flat, there is an isomorphism
    $$\spec \H^H_\ast(\Omega(G/H); \QQ) \cong (\ld{J}_{\ld{G}} \times_{\ld{\g}^\ast\mmod \ld{G}} \ld{\fr{h}}^\ast\mmod \ld{H}) \times_{\ld{J}_{\ld{H}}} \ld{\fr{h}}^\ast\mmod \ld{H}$$
    of group schemes over $\ld{\fr{h}}^\ast\mmod \ld{H}$. Therefore, the study of the $H$-action on $\Omega(G/H)$ is closely related to understanding the map $\ld{J}_{\ld{G}} \times_{\ld{\g}^\ast\mmod \ld{G}} \ld{\fr{h}}^\ast\mmod \ld{H} \to \ld{J}_{\ld{H}}$ (which plays an important role in Langlands transfer).
    
    For instance, let $G = \SL_2$ and $H = \GG_m$; then, $\ld{J}_{\ld{H}} \cong T^\ast \GG_m$, while $\ld{J}_{\ld{G}} \times_{\ld{\g}^\ast\mmod \ld{G}} \ld{\fr{h}}^\ast\mmod \ld{H}$ is isomorphic to the affine blowup $(T^\ast \GG_m)[\tfrac{e^x-1}{x}]$ of $T^\ast \GG_m$. It follows from the preceding discussion that 
    \begin{align*}
        \spec \H^H_\ast(\Omega(G/H); \QQ) & \cong (T^\ast \GG_m)[\tfrac{e^x-1}{x}] \times_{T^\ast \GG_m} \g_m \\
        & \cong (T^\ast \GG_m)[\tfrac{e^x-1}{x}] \times_{\GG_m} \{1\} \cong \spec \QQ[x,\tfrac{e^x-1}{x}]/(x\cdot \tfrac{e^x-1}{x}).
    \end{align*}
    One can verify that this isomorphism holds by computing $\H^H_\ast(\Omega(G/H); \QQ)$ independently; see \cref{ex: borel homology loops S2}.
\end{remark}
\begin{remark}[Singular support]
    It is natural to ask a criterion for when an object of $\Shv^{c,\Sat}_{G\pw{t}}(G\ls{t}/H\ls{t}; \QQ)$ is compact in terms of the equivalence \cref{eq: display rel sat criterion}. If $\cf$ is a compact object of $\Shv^{c,\Sat}_{G\pw{t}}(G\ls{t}/H\ls{t}; \QQ)$, the equivariant cohomology $\H^\ast_{G\pw{t}}(G\ls{t}/H\ls{t}; \cf)$ is finite-dimensional as a vector space. If $\Phi(\cf) \in \Perf(\sh^{1/2}\ld{M}/\ld{G}(-2\rho))$ denotes the corresponding object under the equivalence of \cref{eq: display rel sat criterion}, the first commutative diagram of \cref{thm: ordinary homology criterion satake} implies that the set-theoretic support $\supp(\Phi(\cf))$ intersects the image of $\kappa$ in a zero-dimensional scheme. I expect that $\cf$ is compact if and only if it is set-theoretically supported on the nullcone $\cN_{\ld{M}} := \ld{M} \times_{\ld{\fr{h}}\mmod \ld{H}} \{0\}$ of $\ld{M}$; in other words, that
    \begin{equation}\label{eq: nilpotent singular support}
        \Shv^{c,\Sat}_{G\pw{t}}(G\ls{t}/H\ls{t}; \QQ)^\omega \simeq \Perf_{\cN_{\ld{M}}}(\sh^{1/2} \ld{M}/\ld{G}).
    \end{equation}
    Let $\mu: \ld{M} \to \ld{\g}^\ast$ denote the moment map, and let $\cN_{\ld{\g}}\subseteq \ld{\g}^\ast$ denote the nullcone of $\ld{\g}^\ast$. Then there is a canonical map $\cN_{\ld{M}} \to \mu^{-1}(\cN_{\ld{\g}})$, where $\mu^{-1}(\cN_{\ld{\g}})$ is the \textit{derived} preimage of $\cN_{\ld{\g}}$ under the moment map. It turns out that the map $\cN_{\ld{M}} \to \mu^{-1}(\cN_{\ld{\g}})$ \textit{nearly} induces an isomorphism on reduced schemes (it is an ``Artinian'' thickening\footnote{An easy way to see this is as follows. There is a Cartesian square
    $$\xymatrix{
    \cN_{\ld{M}} \ar[d] \ar[r] \ar[d] & \{0\} \ar[d] \\
    \mu^{-1}(\cN_{\ld{\g}}) \ar[r] & \ld{M}\mmod \ld{G} \times_{\ld{\g}^\ast\mmod \ld{G}} \{0\},
    }$$
    which follows from writing $\cN_{\ld{M}} = \ld{M} \times_{\ld{M}\mmod \ld{G}} \{0\}$, $\mu^{-1}(\cN_{\ld{\g}}) = \ld{M} \times_{\ld{\g}^\ast\mmod \ld{G}} \{0\}$, and $\{0\} = \ld{M} \mmod \ld{G} \times_{\ld{M}\mmod G} \{0\}$. The fiber product $\ld{M}\mmod \ld{G} \times_{\ld{\g}^\ast\mmod \ld{G}} \{0\}$ is an Artinian thickening of $\{0\}$, which implies the desired claim. 
    In fact, since \cref{conj: generalized kostant slice} says that $\ld{M}\mmod \ld{G} \cong \ld{\fr{h}}^\ast \mmod \ld{H}$, one expects an isomorphism $\ld{M}\mmod \ld{G} \times_{\ld{\g}^\ast\mmod \ld{G}} \{0\} \cong \spec C^\ast(G/H; \QQ)$. One again sees that this fiber product is necessarily an Artinian thickening of a point, this time because $G/H$ is a finite CW-complex, so $C^\ast(G/H; \QQ)$ is a finite $\QQ$-module.}); so the singular support in \cref{eq: nilpotent singular support} cannot quite be replaced with singular support contained in $\mu^{-1}(\cN_{\ld{\g}})$.
    For the cases covered by \cref{thm: rk 1 bzsv is true}, the expectation \cref{eq: nilpotent singular support} does indeed hold. In the type T cases, this follows from the argument of \cite[Section 4.5 and 4.6]{braverman-finkelberg}; and in the type G cases, this follows from the argument of \cite[Theorem 12.5.3]{arinkin-gaitsgory-singsupp} applied to the simplest case $G = \PGL_2$.
\end{remark}

There are some expected equivalences of the form \cref{eq: display rel sat criterion} which do not fit into the parameters of \cref{thm: ordinary homology criterion satake}, and are difficult to make ``combinatorial'' since they are not spherical varieties. 
\begin{example}\label{rmk: expected G mod T}
    For simplicity, we will ignore gradings in the following discussion. Let $G_\cc$ be an almost simple algebraic group over $\cc$, and let $T_\cc \subseteq G_\cc$ be a maximal torus. Note that $(G_\cc \times T_\cc)/T_\cc^\mathrm{diag} \cong G_\cc$ is generally not a spherical $G_\cc \times T_\cc$-variety (just for dimension reasons). Nevertheless, one expects that there is a $\QQ$-linear equivalence
    \begin{equation}\label{eq: expected satake G/T}
        \Shv^{c,\Sat}_{(G \times T)\pw{t}}(G\ls{t}; \QQ) \simeq \Perf(\sh^{1/2} \ol{T^\ast(\ld{G}/\ld{N})}/(\ld{G} \times \ld{T}))
    \end{equation}
    which is equivariant for the left-action of $\Shv_{(G \times T)\pw{t} \times (G\times T)\pw{t}}^{c,\Sat}((G\times T)\ls{t}; \QQ) \simeq \Perf(\sh^{1/2} \ld{\g}^\ast(2 - 2\rho)/\ld{G}(-2\rho) \times \ld{\fr{t}}^\ast[2]/\ld{T})$ via \cref{thm: derived satake}. Note that \cref{eq: expected satake G/T} implies that there is a $\QQ$-linear equivalence
    $$\Shv^{c,\Sat}_{G\pw{t}}(G\ls{t}/T\ls{t}; \QQ) \simeq \Perf(\sh^{1/2} \ol{T^\ast(\ld{G}/\ld{N})}/\ld{G}).$$
    Moreover, were the equivalence \cref{eq: expected satake G/T} true, the canonical action of the Weyl group $W = N_G(T)/T$ on the left-hand side of \cref{eq: expected satake G/T} should correspond to the $\ld{G} \times \ld{T}$-equivariant (semi-classical) Gelfand-Graev action on $\ol{T^\ast(\ld{G}/\ld{N})}$ as studied in \cite{ginzburg-kazhdan}. For the case when $G$ has semisimple rank $1$, see \cref{cor: bzsv for CPn} and \cref{rmk: type N PGL2 mod PO2}. When $G = \PGL_3$, the affine closure $\ol{T^\ast(\SL_3/\ld{N})}$ was explicitly identified in \cite{jia-affine-closure} as the closure of the minimal nilpotent orbit in $\fr{so}_8$ (this was also studied earlier by Levasseur-Stafford in \cite{levasseur-stafford}, as well as Kazhdan in \cite{kazhdan-minimal-rep-SO8} in the present context of restriction along $T \subseteq \PGL_3$).  The action of $W = \Sigma_3$ on $\ol{T^\ast(\SL_3/\ld{N})}$ can then be identified as the $\Sigma_3$-action coming from triality on $\fr{so}_8$.
    
    Let us observe now that the criteria (a) and (b) of \cref{thm: ordinary homology criterion satake} can be checked to hold (see below), if we set $\ld{M} = \ol{T^\ast(\ld{G}/\ld{N})}$; this suggests that \cref{eq: expected satake G/T} might indeed hold, if the left-hand side were sufficiently well-behaved. (See also \cite[Theorem 2.3.1]{ginzburg-riche}.)
    
    First, we will show \cref{thm: ordinary homology criterion satake}(b). The relevant Kostant section can be defined as follows. Fix a nondegenerate character $\psi: \ld{\fr{n}} \to \GG_a$, and define the map 
    $$\kappa: \ld{\fr{t}}^\ast \to \ld{\fr{b}}^\ast \hookrightarrow \Ind^{\ld{G}}_{\ld{N}} \ld{\fr{b}}^\ast \hookrightarrow \ol{T^\ast(\ld{G}/\ld{N})},$$
    where the first map is given by the inclusion $\ld{\fr{t}}^\ast \subseteq \ld{\fr{t}}^\ast \oplus \ld{\fr{n}}^\ast$ sending $x \mapsto (x, \psi)$. Then, there is an isomorphism
    \begin{align*}
        \ld{\fr{t}}^\ast \times_{\ld{M}/(\ld{G} \times \ld{T})} \ld{\fr{t}}^\ast & \cong \ld{\fr{t}}^\ast \times_{T^\ast(\ld{G}/\ld{N})/(\ld{G} \times \ld{T})} \ld{\fr{t}}^\ast \\
        & \cong \ld{\fr{t}}^\ast \times_{\ld{\fr{b}}^\ast/\ld{B}} \ld{\fr{t}}^\ast.
    \end{align*}
    It follows from the main result of \cite{abg-iwahori-satake}, or equivalently \cite[Theorem 6.1]{homology-langlands} (see also \cite[Section 4.1]{grg-reg} for a $2$-periodified analogue), that there is an isomorphism
    $$\ld{\fr{t}}^\ast \times_{\ld{\fr{b}}^\ast/\ld{B}} \ld{\fr{t}}^\ast \cong \spec \H^T_\ast(\Omega G; \QQ)$$
    of graded group schemes over $\ld{\fr{t}}^\ast$. This implies \cref{thm: ordinary homology criterion satake}(b).
    
    It remains to check \cref{thm: ordinary homology criterion satake}(a). For this, we need to check that there is an isomorphism 
    $$\co_{(\ld{G} \times \ld{T} \times \ld{\fr{t}}^\ast)/(\ld{\fr{t}}^\ast \times_{\ld{\fr{b}}^\ast/\ld{B}} \ld{\fr{t}}^\ast)} \cong \co_{\ol{T^\ast(\ld{G}/\ld{N})}}.$$
    The $\ld{G} \times \ld{T}$-orbit of the image of $\kappa: \ld{\fr{t}}^\ast \to T^\ast(\ld{G}/\ld{N})$ is the regular locus $T^\ast(\ld{G}/\ld{N})^\reg$, so that there is an isomorphism
    $$\co_{(\ld{G} \times \ld{T} \times \ld{\fr{t}}^\ast)/(\ld{\fr{t}}^\ast \times_{\ld{\fr{b}}^\ast/\ld{B}} \ld{\fr{t}}^\ast)} \cong \co_{T^\ast(\ld{G}/\ld{N})^\reg}.$$
    The inclusion $T^\ast(\ld{G}/\ld{N})^\reg \subseteq T^\ast(\ld{G}/\ld{N})$ has complement of codimension $\geq 2$, since it can be identified with the inclusion $\Ind_{\ld{N}}^{\ld{G}} \ld{\fr{b}}^{\ast,\reg} \subseteq \Ind_{\ld{N}}^{\ld{G}} \ld{\fr{b}}^{\ast}$; the claim follows from \cref{lem: some codim stuff}(a) along with the observation that the inclusion $\ld{\fr{b}}^{\ast,\reg} \subseteq \ld{\fr{b}}^{\ast}$ has complement of codimension $\geq 2$.
    The inclusion $T^\ast(\ld{G}/\ld{N}) \subseteq \ol{T^\ast(\ld{G}/\ld{N})}$ also has complement of codimension $\geq 2$ (e.g., by \cref{lem: some codim stuff}(b)), so that there are isomorphisms
    $$\co_{T^\ast(\ld{G}/\ld{N})^\reg} \cong \co_{T^\ast(\ld{G}/\ld{N})} \cong \co_{\ol{T^\ast(\ld{G}/\ld{N})}}$$
    by the algebraic Hartogs lemma.
    This verifies \cref{thm: ordinary homology criterion satake}(a), as desired. 
\end{example}
\begin{remark}
    One can extend \cref{rmk: expected G mod T} to more general Levi subgroups as follows. Let $I\subseteq \ld{\Delta}$ denote a subset of the simple roots of $G$, let $P_I$ denote the associated parabolic subgroup, and let $L_I \subseteq P_I$ denote a fixed Levi factor. The subset $I$ defines a subset of simple roots of $\ld{G}$, which we will also denote by $I$ (for simplicity). Following \cite{macerato-levi-restriction}, let $\psi_I: \ld{N} \to \GG_a$ denote the additive character given by the composite
    $$\ld{N} \to \ld{N}/[\ld{N}, \ld{N}] \cong \prod_\Delta \GG_a \xar{\proj_I} \prod_I \GG_a \xar{\sum} \GG_a,$$
    so that $\psi_I$ defines an element of $\ld{\fr{n}}^\ast$. This allows one to define the Whittaker reduction $T^\ast(\ld{G}/_{\psi_I} \ld{N})$, which admits a natural $\ld{G}$-action. Let $\ol{T^\ast(\ld{G}/_{\psi_I} \ld{N})}$ denote its affine closure. We then expect the following:
    \begin{conjecture}\label{conj: spherical levi}
        Let $L_I$ denote the above Levi subgroup of $G$. Then there is an isomorphism
        $$\ol{T^\ast(\ld{G}/_{\psi_I} \ld{N})} \cong \ol{(\ld{G} \times \ld{\fr{l}}_I^\ast\mmod \ld{L}_I)/\ld{J}_X},$$
        where $\ld{J}_X = \spec \H^{L_I}_\ast(\Omega(G/L_I); k)$ (so that $\ld{J}_X$ can be identified with the kernel of the homomorphism $\ld{J}_{\ld{G}} \times_{\ld{\g}^\ast\mmod \ld{G}} \ld{\fr{l}}_I^\ast\mmod \ld{L}_I \to \ld{J}_{\ld{L}_I}$ as in \cref{eq: J and lag corr}).
        
        More strongly, there is an equivalence
        $$\Shv^{c,\Sat}_{G\pw{t}}(G\ls{t}/L_I\ls{t}; \QQ) \simeq \Perf(\sh^{1/2} \ol{T^\ast(\ld{G}/_{\psi_I} \ld{N})}/\ld{G})$$
        which is equivariant for the left-action of $\Shv_{(G \times G)\pw{t}}^{c,\Sat}(G\ls{t}; \QQ) \simeq \Perf(\sh^{1/2} \ld{\g}^\ast/\ld{G})$ via \cref{thm: derived satake}; we have again omitted the grading on the spectral side for simplicity.
    \end{conjecture}
    If the hypotheses of \cref{thm: ordinary homology criterion satake} are satisfied for $L_I \subseteq G$, the equivalence of categories in \cref{conj: spherical levi} follows from the first isomorphism therein. In future work \cite{gannon-me} with Tom Gannon, we will address the first part of \cref{conj: spherical levi} using ``Whittaker descent''.
    
    When $L_I$ is spherical\footnote{Note that there is a classification of spherical Levi subgroups of simple linear algebraic groups using Kr\"amer's classification \cite{kramer-spherical-subgroups}; see \cite[Theorem 4.1]{brundan-spherical-levi}. For instance, if $G$ is a classical group, the only possibilities are $\GL_j \times \GL_{n-j} \subseteq \GL_j$, $\SO_2 \times \SO_{2n-1} \subseteq \SO_{2n+1}$, $\GL_n \subseteq \SO_{2n+1}$, $\GG_m \times \Sp_{2n-2} \subseteq \Sp_{2n}$, $\GL_n \subseteq \Sp_{2n}$, $\SO_2 \times \SO_{2n-2} \subseteq \SO_{2n}$, and $\GL_n \subseteq \SO_{2n}$.}, the relationship between the spectral side of \cref{conj: spherical levi} and the predicted dual variety of \cref{conj: bzsv} seems to be very nontrivial. For instance, when $\ld{G} = \GL_{2n+1}$ and $I$ corresponds to the partition $[n,n+1]$, the first part of \cref{conj: spherical levi} (which will be proved in future work) along with \cref{ex: mirabolic satake and variants} shows that there is an isomorphism
    $$\ol{T^\ast(\GL_{2n+1}/_{\psi_I} \ld{N})} \cong T^\ast(\GL_{2n+1}/(\GL_n \times \GL_{n+1})).$$
    However, this isomorphism does not seem so easy to see directly.
    \cref{conj: bzsv} in the more general case of the spherical Levi $\GL_j \times \GL_{n-j} \subseteq \GL_n$ is work-in-progress of Chen-Macerato-Nadler-O'Brien.
    
    Let us remark on one interesting consequence of \cref{conj: spherical levi}. Write $\N_G(L_I)$ to denote the normalizer of $L_I \subseteq G$. There is a natural action of the relative Weyl group $W_I = \N_G(L_I)/L_I$ on the left-hand side of \cref{conj: spherical levi}, which defines an action of $W_I$ on the right-hand side. Based on \cref{rmk: expected G mod T}, it is natural to hope that there is in fact an action of $W_I$ on $\ol{T^\ast(\ld{G}/_{\psi_I} \ld{N})}$ which commutes with its natural $\ld{G}$-action. This would be a parabolic variant of the (semi-classical) Gelfand-Graev action; see \cref{rmk: a composite map in homology}. It implies the following extension of \cref{conj: spherical levi}:
    $$\Shv^{c,\Sat}_{G\pw{t}}(G\ls{t}/\N_G(L_I)\ls{t}; \QQ) \simeq \Perf(\sh^{1/2} (\ol{T^\ast(\ld{G}/_{\psi_I} \ld{N})}/W_I)/\ld{G}).$$
    In particular, the Hamiltonian $\ld{G}$-``space'' which is dual to the spherical $G$-variety $G/\N_G(L_I)$ (which often has roots of type N) would be the \textit{stack} $\ol{T^\ast(\ld{G}/_{\psi_I} \ld{N})}/W_I$.
\end{remark}

Let us now shift gears somewhat. The following result is related to the discussion in \cite[Section 5.1.5]{sakellaridis-icm} and to \cite[Section 5.2]{teleman-icm}.
\begin{prop}\label{prop: lagrangian correspondence}
    Let $H\subseteq G$ be a closed subgroup. Then there is a Lagrangian correspondence (interpreted in a derived sense)
    $$\xymatrix{
    & \ld{J}_{\ld{G}} \times_{\ld{\g}^\ast\mmod \ld{G}} \ld{\fr{h}}^\ast\mmod \ld{H} \ar[dl] \ar[dr] & \\
    \ld{J}_{\ld{H}} & & \ld{J}_{\ld{G}},
    }$$
    where the left map restricts to the zero section of $\ld{J}_{\ld{H}}$ when pulled back to the identity section of $\ld{J}_{\ld{G}}$.
\end{prop}
\begin{proof}
    The desired claim follows from the analogous statement at the level of Lie algebras. It is a classical fact that the Lie algebra of $\ld{J}_{\ld{G}}$ can be identified with the cotangent bundle $T^\ast(\ld{\g}\mmod \ld{G}) \cong T^\ast(\ld{\g}^\ast\mmod \ld{G})$, and similarly for $\ld{J}_{\ld{H}}$. We therefore need to see that there is a Lagrangian correspondence
    $$\xymatrix{
    & T^\ast(\ld{\g}^\ast\mmod \ld{G}) \times_{\ld{\g}^\ast\mmod \ld{G}} \ld{\fr{h}}^\ast\mmod \ld{H} \ar[dl] \ar[dr] & \\
    T^\ast(\ld{\fr{h}}^\ast\mmod \ld{H}) & & T^\ast(\ld{\g}^\ast\mmod \ld{G}).
    }$$
    More generally, if $Y \to Z$ is a map between schemes, there is a Lagrangian correspondence
    $$\xymatrix{
    & T^\ast(Z) \times_Z Y \ar[dl] \ar[dr] & \\
    T^\ast Y & & T^\ast Z.
    }$$
    This is of course well-known if $Y \to Z$ is a smooth map of smooth schemes, but the same continues to hold in general (see, e.g., \cite[Theorem 2.8]{calaque-shifted-cotangent}).
    Taking $Y \to Z$ to be the map $\ld{\fr{h}}^\ast\mmod \ld{H} \to \ld{\g}^\ast\mmod \ld{G}$, we win.
\end{proof}

\begin{remark}\label{rmk: lag corr in homogeneous case}
    The left map in \cref{prop: lagrangian correspondence} is precisely the one of \cref{rmk: reg centr and endoscopy}. Note that the map $\ld{J}_X \to \ld{J}_{\ld{G}}$ of \cref{rmk: a composite map in homology} is simply obtained by intersecting this Lagrangian correspondence with the identity section of $\ld{J}_{\ld{H}}$; in other words, $\ld{J}_X$ is the kernel of the homomorphism $\ld{J}_{\ld{G}} \times_{\ld{\g}^\ast\mmod \ld{G}} \ld{\fr{h}}^\ast\mmod \ld{H} \to \ld{J}_{\ld{H}}$. Concretely, there is a commutative diagram
    \begin{equation}\label{eq: J and lag corr}
        \xymatrix{
        & \ld{J}_X \ar[dr] \ar[dl] & & \\
        \ld{\fr{h}}^\ast\mmod \ld{H} \ar[dr] & & \ld{J}_{\ld{G}} \times_{\ld{\g}^\ast\mmod \ld{G}} \ld{\fr{h}}^\ast\mmod \ld{H} \ar[dl] \ar[dr] & \\
        & \ld{J}_{\ld{H}} & & \ld{J}_{\ld{G}},
        }
    \end{equation}
    where the square is Cartesian. This implies that the map $\ld{J}_X \to \ld{J}_{\ld{G}}$ is Lagrangian (in a derived sense). Moreover, it implies that there is an isomorphism
    $$\mathrm{Lie}(\ld{J}_X) \cong T^\ast[1](\ld{\fr{h}}^\ast\mmod \ld{H} / \ld{\g}^\ast\mmod \ld{G}),$$
    where the right-hand side denotes the $1$-shifted cotangent bundle. The formula \cref{eq: M as G mod JX} also shows that
    \begin{align*}
        \ld{M}^\ddag & \cong \ol{(\ld{J}_{\ld{H}} \times \ld{G})/(\ld{J}_{\ld{G}} \times_{\ld{\g}^\ast\mmod \ld{G}} \ld{\fr{h}}^\ast\mmod \ld{H})} \\
        & \cong \ol{(\ld{J}_{\ld{H}} \times_{\ld{\g}^\ast\mmod \ld{G}} T^\ast(\ld{G}/_\psi \ld{N}))/(\ld{J}_{\ld{G}} \times_{\ld{\g}^\ast\mmod \ld{G}} \ld{\fr{h}}^\ast\mmod \ld{H})}.
    \end{align*}
    The final isomorphism comes from the identification $T^\ast(\ld{G}/_\psi \ld{N}) \cong \ld{G} \times \ld{\g}^\ast\mmod \ld{G}$ via \cref{thm: kostant}. If $\ld{G}$ has trivial center, for instance, the group scheme $\ld{J}_{\ld{G}}$ is connected, and so we find that \cref{eq: M as G mod JX} can be rewritten to describe $\co_{\ld{M}^\ddag}$ as the Poisson centralizer
    $$\co(\ld{M}^\ddag) \cong \co(\ld{J}_{\ld{H}} \times_{\ld{\g}^\ast\mmod \ld{G}} T^\ast(\ld{G}/_\psi \ld{N}))^{\co(\ld{\g}^\ast\mmod \ld{G})}.$$
    This is a formula analogous to \cite[Theorem 1.3.3]{ginzburg-kazhdan}. Again, if one defines $\ld{M}^\ddag$ in this way, the question of proving \cref{conj: bzsv} (when $H\subseteq G$ is spherical) now becomes about identifying $\ld{M}^\ddag$ with the prescription of \cite{bzsv}.
    
    Furthermore, using the main result of \cite{bfm} and \cref{thm: ordinary homology criterion satake} (all of which is related to \cite[Theorem 5.3]{teleman-icm}), the above diagram \cref{eq: J and lag corr} can be identified with
    $$\xymatrix@=.75em{
        & \spec \H^H_\ast(\Omega(G/H); \QQ) \ar[dr] \ar[dl] & & \\
        \spec \H_H^\ast(\ast; \QQ) \ar[dr] & & \spec \H^H_\ast(\Omega G; \QQ) \ar[dl] \ar[dr] & \\
        & \spec \H^H_\ast(\Omega H; \QQ) & & \spec \H^G_\ast(\Omega G; \QQ).
    }$$
    The long composite on the right-hand side of the above diagram will be Lagrangian, hence coisotropic (but this has to be interpreted in a derived sense; for example, it need not be a closed immersion!). We will study this phenomenon of coisotropicity (upon completion) in much greater detail in \cref{subsec: ku and lagrangians} from the perspective of Hochschild cohomology and centralizers \`a la \cite[Section 5.3]{HA} and \cite{francis}. See \cref{rmk: coisotropic correspondence}, for instance.
\end{remark}
\begin{remark}\label{rmk: extension of generalized kostant}
    \cref{prop: lagrangian correspondence} should admit the following generalization to the non-homogeneous case. Suppose $\ld{M}$ is the Hamiltonian $\ld{G}$-space dual to a spherical $G$-variety $X$ satisfying the hypotheses of \cref{conj: bzsv}. If one assumes \cref{conj: generalized kostant slice}, i.e., that there is a Kostant section $\kappa_{\ld{M}}: \ld{M}\mmod \ld{G} \to \ld{M}$ making the relevant diagram commute, it is not hard to see that the map $\kappa_{\ld{M}}$ factors as a map 
    $$\tilde{\kappa}_{\ld{M}}: \ld{M}\mmod \ld{G} \to \ld{M}/_\psi \ld{N}.$$
    Indeed, \cref{thm: kostant} gives an isomorphism $\ld{\g}^\ast\mmod \ld{G} \cong \ld{G}\backslash T^\ast(\ld{G}/_\psi \ld{N})$, so that the Kostant slice $\ld{\g}^\ast\mmod \ld{G} \to \ld{\g}^\ast/\ld{G}$ is the ``kernel'' of Whittaker reduction. This implies that there is an isomorphism
    $$\ld{M}/_\psi \ld{N} \cong \ld{\g}^\ast\mmod \ld{G} \times_{\ld{\g}^\ast/\ld{G}} \ld{M}/\ld{G}.$$
    The commutative diagram of \cref{conj: generalized kostant slice} (and the universal property of fiber products) now gives the desired map $\tilde{\kappa}_{\ld{M}}$.
    We then expect: 
    \begin{conjectureno}
        In addition to \cref{conj: generalized kostant slice}, there is a Lagrangian correspondence
        $$\xymatrix{
        & \ld{J}_{\ld{G}} \times_{\ld{\g}^\ast\mmod \ld{G}} \ld{M}\mmod \ld{G} \ar[dl] \ar[dr] & \\
        \ld{M}/_\psi \ld{N} & & \ld{J}_{\ld{G}},
        }$$
        where the left map restricts to $\tilde{\kappa}_{\ld{M}}$ when pulled back to the identity section of $\ld{J}_{\ld{G}}$. 
    \end{conjectureno}
    When $\ld{G}$ is a torus, the above conjecture is essentially the same as asking that the closed immersion $\ld{M}\mmod \ld{G} \to \ld{M}$ exhibits $\ld{M}\mmod \ld{G}$ as a Lagrangian in $\ld{M}$.
    
    Note that \cite[Example 8.4.5]{bzsv} says that when $X = G/H$ is homogeneous, \cref{conj: bzsv} implies that the Whittaker reduction $\ld{M}/_\psi \ld{N}$ is isomorphic to $\ld{J}_{\ld{H}}$. (Working with $\ld{M}^\ddag$ instead of $\ld{M}$, the identification of $\ld{M}^\ddag/_\psi \ld{N}$ with $\ld{J}_{\ld{H}}$ is almost tautological given that $\ld{J}_X$ is the kernel of the homomorphism $\ld{J}_{\ld{G}} \times_{\ld{\g}^\ast\mmod \ld{G}} \ld{\fr{h}}^\ast\mmod \ld{H} \to \ld{J}_{\ld{H}}$.) Therefore, \cref{conj: generalized kostant slice} and \cref{conj: bzsv} imply that the above Lagrangian correspondence specializes to that of \cref{prop: lagrangian correspondence}.
\end{remark}

\cref{prop: lagrangian correspondence} has an interesting consequence.
\begin{construction}\label{cstr: theta}
    The homomorphisms
    \begin{align*}
        \ld{J}_{\ld{G}} \times_{\ld{\g}^\ast\mmod \ld{G}} \ld{\fr{h}}^\ast\mmod \ld{H} & \to \ld{G} \times \ld{\fr{h}}^\ast\mmod \ld{H}, \\
        \ld{J}_{\ld{G}} \times_{\ld{\g}^\ast\mmod \ld{G}} \ld{\fr{h}}^\ast\mmod \ld{H} & \to \ld{J}_{\ld{H}} \to \ld{H} \times \ld{\fr{h}}^\ast\mmod \ld{H}
    \end{align*}
    define a closed immersion 
    $$\ld{J}_{\ld{G}} \times_{\ld{\g}^\ast\mmod \ld{G}} \ld{\fr{h}}^\ast\mmod \ld{H} \to \ld{G} \times \ld{H} \times \ld{\fr{h}}^\ast\mmod \ld{H}$$
    of group schemes over $\ld{\fr{h}}^\ast\mmod \ld{H}$. 
    Let $\ld{\cM}^\ddag$ denote the affine closure 
    $$\ld{\cM}^\ddag = \ol{(\ld{G} \times \ld{H} \times \ld{\fr{h}}^\ast\mmod \ld{H})/(\ld{J}_{\ld{G}} \times_{\ld{\g}^\ast\mmod \ld{G}} \ld{\fr{h}}^\ast\mmod \ld{H})},$$
    so that
    \begin{align*}
        \dim(\ld{\cM}^\ddag) & = \dim(\ld{G}) + \dim(\ld{H}) + \mathrm{rank}(\ld{H}) - \mathrm{rank}(\ld{G})\\
        & = 2 \left( \dim(\ld{H}/N_{\ld{H}}) + \dim(\ld{G}/B_{\ld{G}}) \right),
    \end{align*}
    where $N_{\ld{H}}$ is the unipotent radical of a Borel subgroup of $\ld{H}$, and $B_{\ld{G}}$ is a Borel subgroup of $\ld{G}$.
    It can be shown that $\ld{\cM}^\ddag$ admits the structure of a Hamiltonian $\ld{G} \times \ld{H}$-space (in fact, this is a consequence of the second part of \cref{cor: BJ and lag corr} below and \cref{prop: safronov lag maps}).
\end{construction}
\begin{corollary}\label{cor: BJ and lag corr}
    {Define} $\ld{M}^\ddag$ as in \cref{eq: M as G mod JX}, and let $\ld{M}^{\ddag,\reg} = (\ld{G} \times \ld{\fr{h}}^\ast\mmod \ld{H})/\ld{J}_X$ denote the $\ld{G}$-orbit of the map $\kappa_{\ld{M}^\ddag}: \ld{\fr{h}}^\ast\mmod \ld{H} \to \ld{M}^\ddag$.
    Let $\ld{\cM}^\ddag$ denote the Hamiltonian $\ld{G} \times \ld{H}$-space of \cref{cstr: theta}, and define $\ld{\cM}^{\ddag, \reg}$ similarly. Then there is an isomorphism 
    $$\ld{\cM}^{\ddag,\reg}/\ld{H} \cong \ld{\g}^{\ast, \reg} \times_{\ld{\g}^\ast\mmod \ld{G}} \ld{\fr{h}}^\ast\mmod \ld{H},$$
    and a diagram
    \begin{equation}\label{eq: BJ span}
        \xymatrix{
        & \ld{M}^{\ddag, \reg}/\ld{G} \ar[dr] \ar[dl] & & \\
        \ld{\fr{h}}^\ast\mmod \ld{H} \ar[dr]_-\kappa & & \ld{\cM}^{\ddag, \reg}/(\ld{H} \times \ld{G}) \ar[dl] \ar[dr] & \\
        & \ld{\fr{h}}^{\ast, \reg}/\ld{H} & & \ld{\g}^{\ast, \reg}/\ld{G},
        }
    \end{equation}
    where the square is Cartesian, the long composite on the right-hand side is a Lagrangian morphism, and the span at the bottom of the diagram is a ($1$-shifted) Lagrangian correspondence. In particular, there is a Cartesian square
    \begin{equation}\label{eq: theta and whit}
        \xymatrix{
        \ld{M}^{\ddag, \reg} \ar[d] \ar[r] & \ld{\cM}^{\ddag, reg} \ar[d] \\
        \ld{\fr{h}}^{\ast}\mmod \ld{H} \ar[r]_-\kappa & \ld{\fr{h}}^{\ast, \reg}.
        }
    \end{equation}
\end{corollary}
\begin{proof}
    By \cref{cstr: theta},
    $$\ld{\cM}^{\ddag, \reg} \cong (\ld{G} \times \ld{H} \times \ld{\fr{h}}^\ast\mmod \ld{H})/(\ld{J}_{\ld{G}} \times_{\ld{\g}^\ast\mmod \ld{G}} \ld{\fr{h}}^\ast\mmod \ld{H}),$$
    so that
    \begin{align*}
        \ld{\cM}^{\ddag, \reg}/\ld{H} & \cong (\ld{G} \times \ld{\fr{h}}^\ast\mmod \ld{H})/(\ld{J}_{\ld{G}} \times_{\ld{\g}^\ast\mmod \ld{G}} \ld{\fr{h}}^\ast\mmod \ld{H}) \\
        & \cong ((\ld{G} \times \ld{\g}^\ast\mmod \ld{G})/\ld{J}_{\ld{G}}) \times_{\ld{\g}^\ast\mmod \ld{G}} \ld{\fr{h}}^\ast\mmod \ld{H} \\
        & \cong \ld{\g}^{\ast,\reg} \times_{\ld{\g}^\ast\mmod \ld{G}} \ld{\fr{h}}^\ast\mmod \ld{H},
    \end{align*}
    as desired. It follows from this identification that the desired diagram \cref{eq: BJ span} then becomes 
    $$\xymatrix{
        & \ld{M}^{\ddag, \reg}/\ld{G} \ar[dr] \ar[dl] & & \\
        \ld{\fr{h}}^\ast\mmod \ld{H} \ar[dr]_-\kappa & & (\ld{\g}^{\ast, \reg} \times_{\ld{\g}^\ast\mmod \ld{G}} \ld{\fr{h}}^\ast\mmod \ld{H})/\ld{G} \ar[dl] \ar[dr] & \\
        & \ld{\fr{h}}^{\ast, \reg}/\ld{H} & & \ld{\g}^{\ast, \reg}/\ld{G},
    }$$
    which satisfies the desired properties since it is obtained by taking classifying stacks of the diagram in \cref{prop: lagrangian correspondence} via the identifications
    $$B_{\ld{\g}^\ast\mmod \ld{G}} \ld{J}_{\ld{G}} \cong \ld{\g}^{\ast, \reg}/\ld{G}, \ B_{\ld{\fr{h}}^\ast\mmod \ld{H}} \ld{J}_{\ld{H}} \cong \ld{\fr{h}}^{\ast, \reg}/\ld{H}, \ B_{\ld{\fr{h}}^\ast\mmod \ld{H}} \ld{J}_X \cong \ld{M}^{\ddag, \reg}/\ld{G}$$
    coming from \cref{cor: classifying stack of J}.
\end{proof}

\begin{remark}
    One can use \cref{thm: classical homology loops G} to identify 
    $$\spec \H^H_\ast(\Omega ((G\times H)/H); \QQ) \cong \ld{J}_{\ld{G}} \times_{\ld{\g}^\ast\mmod \ld{G}} \ld{\fr{h}}^\ast\mmod \ld{H}.$$
    Suppose that the subgroup $H^\mathrm{diag} \subseteq G \times H$ satisfies the hypotheses of \cref{hypothesis: rank 1 weakly placid}. As in \cref{rmk: a composite map in homology}, we then expect that if $H^\mathrm{diag} \subseteq G \times H$ is a spherical subgroup, $\ld{\cM}^\ddag$ is isomorphic to the Hamiltonian $\ld{G} \times \ld{H}$-space $\ld{\cM}$ dual to $H^\mathrm{diag} \subseteq G \times H$. Note that similarly to the Cartesian square \cref{eq: theta and whit}, there is a Cartesian square
    \begin{equation}\label{eq: symplectic reduction G/H}
        \xymatrix{
        T^\ast(G/H) \ar[r] \ar[d] & T^\ast((G \times H)/H^\mathrm{diag}) \ar[d]^-\mu \\
        \{0\} \ar[r] & \fr{h}^\ast.
        }
    \end{equation}
    In other words, the diagram analogous to \cref{eq: BJ span} in this case is the restriction to regular loci of
    $$\xymatrix@=.75em{
        & T^\ast(G/H)/G \cong (\g/\fr{h})^\ast/H \ar[dr] \ar[dl] & & \\
        BH \ar[dr]_-{\{0\}} & & T^\ast((G \times H)/H^\mathrm{diag})/(G\times H) \cong \g^\ast/H \ar[dl] \ar[dr] & \\
        & \fr{h}^{\ast}/H & & \g^\ast/G,
    }$$
    where again the square is Cartesian, the long composite on the right-hand side is the moment map for $T^\ast(G/H)$, and the span at the bottom of the diagram is a ($1$-shifted) Lagrangian correspondence.
\end{remark}
\begin{remark}\label{rmk: dirichlet and neumann}
    Assume now that $\ld{\cM}^\ddag \cong \ld{\cM}$, and similarly $\ld{M}^\ddag \cong \ld{M}$.
    The square of \cref{eq: theta and whit} then says that the Whittaker reduction of the $\ld{H}$-action on $\ld{\cM}$ identifies with $\ld{M}$. Since the dual to $T^\ast((G \times H)/H^\mathrm{diag})$ is $\ld{\cM}$, and the dual to $T^\ast(G/H)$ is $\ld{M}$, the squares \cref{eq: theta and whit} and \cref{eq: symplectic reduction G/H} showcase the Langlands duality between ``symplectic reduction at $0$'' and ``Whittaker reduction''. In the language of quantum field theories, this is the duality between the Dirichlet and Neumann boundary conditions.
\end{remark}
\begin{example}
    If $H$ is a Levi subgroup (spherical or not!) of $G$ with associated parabolic $P$ and unipotent radical $N_P$, for instance, it turns out that one can identify $\ld{\cM}^\ddag$ with the affine closure of $T^\ast(\ld{G}/N_{\ld{P}}^-)$; this will follow from the proof in \cite{gannon-me} of the first part of \cref{conj: spherical levi}. The span at the bottom of the diagram \cref{eq: BJ span} identifies with the restriction to regular loci of the Lagrangian correspondence
    $$\xymatrix{
    & \tilde{\ld{\g}}_{\ld{P}}/\ld{G} \cong T^\ast(\ld{G}/N_{\ld{P}}^-)/(\ld{G} \times \ld{H}) \ar[dl] \ar[dr] & \\
    \ld{\fr{h}}^{\ast}/\ld{H} & & \ld{\g}^{\ast}/\ld{G}
    }$$
    coming from the parabolic Grothendieck-Springer resolution (see \cite{safronov-symplectic-implosion}). This span extends to the affine closure $\ol{T^\ast(\ld{G}/N_{\ld{P}}^-)}$, i.e., there is a span
    $$\xymatrix{
    & \ol{T^\ast(\ld{G}/N_{\ld{P}}^-)}/(\ld{G} \times \ld{H}) \ar[dl] \ar[dr] & \\
    \ld{\fr{h}}^{\ast}/\ld{H} & & \ld{\g}^{\ast}/\ld{G}.
    }$$
    Let us make the following pleasant observation: all constructions on the topological side depend only on the choice of Levi $H\subseteq G$, and \textit{not} on the parabolic $P$. Although the first span \textit{does} rely on the choice of parabolic to even define $\tilde{\ld{\g}}_{\ld{P}}$, the formula for $\ld{\cM}^\ddag$ shows that $\ol{T^\ast(\ld{G}/N_{\ld{P}}^-)}$ does {not} depend on the choice of parabolic.
\end{example}
\begin{example}[Gan-Gross-Prasad]\label{ex: ggp}
    For instance, suppose $H\subseteq G$ is the inclusion $\SO_{2n} \subseteq \SO_{2n+1}$, so that $\ld{H} = \SO_{2n}$ and $\ld{G} = \Sp_{2n}$. In this case, one can show that
    \begin{equation}\label{eq: std tensor std}
        \ld{\cM}^\ddag \cong \Hom(\std_{2n}, \std_{2n})
    \end{equation}
    as a $\SO_{2n} \times \Sp_{2n}$-variety (this is the Gan-Gross-Prasad period), so that \cref{thm: ordinary homology criterion satake} implies that if \cref{hypothesis: rank 1 weakly placid} is true in this example, the $\infty$-category $\Shv^{c,\Sat}_{(G\times H)\pw{t}}((G\times H)\ls{t}/H^\mathrm{diag}\ls{t}; \QQ)$ is equivalent to $\Perf(\sh^{1/2} \Hom(\std_{2n}, \std_{2n})/(\SO_{2n} \times \Sp_{2n})(-2\rho))$ for a certain grading on $\Hom(\std_{2n}, \std_{2n})$. When $n=1$, this gives an alternative perspective on \cref{cor: bzsv for CPn} in the special case $\PGL_2/\GG_m = \SO_3/\SO_2$.
    
    To prove \cref{eq: std tensor std}, one must describe a morphism
    \begin{equation}\label{eq: Mddag for theta}
        (\ld{G} \times \ld{H} \times \ld{\fr{h}}^\ast\mmod \ld{H})/(\ld{J}_{\ld{G}} \times_{\ld{\g}^\ast\mmod \ld{G}} \ld{\fr{h}}^\ast\mmod \ld{H}) \hookrightarrow \Hom(\std_{2n}, \std_{2n}).
    \end{equation}
    We will construct this morphism below; the key step is to define a Kostant slice 
    $$\kappa: \ld{\fr{h}}^\ast\mmod \ld{H} \to \Hom(\std_{2n}, \std_{2n})$$
    as in \cref{conj: generalized kostant slice}. With a bit more work, one can check that it is an open immersion with complement of codimension $2$, from which \cref{eq: std tensor std} follows.

    Let us now describe the Kostant slice $\kappa$.
    Recall from \cref{ex: classical groups reg centralizer} that:
    \begin{itemize}
        \item We may identify 
        $$\ld{\g}^\ast\mmod \ld{G} \cong \spec k[p_1, \cdots, p_n] \cong \spec \H^\ast_{\SO_{2n+1}}(\ast; k),$$
        and $\ld{J}_{\ld{G}}$ is the group scheme whose fiber over $\vec{p} := (p_1, \cdots, p_n)$ is the subgroup of those units $f(t) \in k[t]/(t^{2n} + p_1 t^{2n-2} + \cdots + p_n)$ such that $f(t)^{-1} = f(-t)$.
        Recall that $k[t]/(t^{2n} + p_1 t^{2n-2} + \cdots + p_n)$ admits the structure of a symplectic vector space.
        \item We may identify 
        $$\ld{\fr{h}}^\ast\mmod \ld{H} \cong \spec k[p_1, \cdots, p_{n-1}, c_n] \cong \spec \H^\ast_{\SO_{2n}}(\ast; k),$$
        and $\ld{J}_{\ld{H}}$ is the group scheme whose fiber over $(\vec{p}, c_n) := (p_1, \cdots, p_{n-1}, c_n)$ is the subgroup of those units $f(t,v) \in k[t, v]/(tv-c_n, t^{2n-2} + p_1 t^{2n-4} + \cdots + p_{n-1} + v^2)$ such that $f(t,v)^{-1} = f(-t,-v)$.
        Recall that $k[t, v]/(tv-c_n, t^{2n-2} + p_1 t^{2n-4} + \cdots + p_{n-1} + v^2)$ admits the structure of a quadratic vector space.
    \end{itemize}
    The map $\pi: \ld{\fr{h}}^\ast \mmod \ld{H} \to \ld{\g}^\ast\mmod \ld{G}$ is induced by the inclusion
    $$\pi: k[p_1, \cdots, p_n] \hookrightarrow k[p_1, \cdots, p_{n-1}, c_n], \ p_n \mapsto c_n^2.$$
    The map $\pi$ induces a map of $2n$-dimensional $k$-vector spaces
    \begin{multline*}
        \varphi_{\vec{p}, c_n}: k[t]/(t^{2n} + p_1 t^{2n-2} + \cdots + p_{n-1} t^2 + c_n^2) \\
        \cong k[t]/(t^{2n} + p_1 t^{2n-2} + \cdots + p_n) \otimes_{\co_{\ld{\g}^\ast\mmod \ld{G}}} \co_{\ld{\fr{h}}^\ast \mmod \ld{H}} \\
        \to k[t, v]/(tv-c_n, t^{2n-2} + p_1 t^{2n-4} + \cdots + p_{n-1} + v^2);
    \end{multline*}
    in other words, this is a linear map $\varphi_{\vec{p}, c_n}: \std_{2n} \to \std_{2n}$. That is, $\pi$ induces a map
    $$\kappa: \ld{\fr{h}}^\ast \mmod \ld{H} \cong \spec k[p_1, \cdots, p_{n-1}, c_n] \to \Hom(\std_{2n}, \std_{2n}), \ (\vec{p}, c_n) \mapsto \varphi_{\vec{p}, c_n}.$$
    This is the relevant Kostant slice, and so we obtain a map
    $$\ld{G} \times \ld{H} \times \ld{\fr{h}}^\ast\mmod \ld{H} \to \Hom(\std_{2n}, \std_{2n})$$
    by taking the $\ld{G} \times \ld{H}$-orbit of $\kappa$. This map factors through the quotient by $\ld{J}_{\ld{G}} \times_{\ld{\g}^\ast\mmod \ld{G}} \ld{\fr{h}}^\ast\mmod \ld{H}$, and therefore produces the desired map \cref{eq: Mddag for theta}.
\end{example}
\begin{remark}\label{rmk: minimal nilpotent orbits}
    One rather curious phenomenon which can be observed empirically, but we do not have a general explanation for, is a relationship between $\ld{\cM}^\ddag$ and closures of minimal nilpotent orbits. Namely, suppose $H\subseteq G$ is a closed subgroup (not necessarily spherical), so that $H^\mathrm{diag} \subseteq G \times H$; this need not be spherical even if $H \subseteq G$ is spherical. Assume that there is a reductive group scheme $\ld{K}$ equipped with an embedding $\ld{G} \times \ld{H} \subseteq \ld{K}$ of group schemes such that $\ld{G}$ is the centralizer of $\ld{H}$, and conversely $\ld{H}$ is the centralizer of $\ld{G}$. In other words, $(\ld{G}, \ld{H})$ form a \textit{(reductive) dual pair in $\ld{K}$}. It turns out that $\ld{\cM}^\ddag$, as defined in \cref{cstr: theta}, can sometimes be identified with the closure of the minimal nilpotent orbit of $\ld{K}$.
    
    Let us give three examples to illustrate this; a general explanation of this observation would be very interesting!
    \begin{itemize}
        \item Let $\SL_3 \subseteq \SO_8 = \ld{K}$ via the adjoint action of $\SL_3$ on $\sl_3$ equipped with its Killing form. From \cite[Table 1]{deligne-gross-exceptional-series}, one finds that the centralizer of $\SL_3$ is $\GG_m^2 \subseteq \SO_8$, and $(\SL_3, \GG_m^2)$ forms a dual pair in $\SO_8$. (It corresponds to the inclusion $H = \GG_m^2/\mu_3 \subseteq \PGL_3 = G$ of the maximal torus; note that this is \textit{not} a spherical subgroup!) Then $\ld{\cM}^\ddag$ can be identified with $\ol{T^\ast(\SL_3/\ld{N})}$, and this is isomorphic to the closure of the minimal nilpotent orbit in $\fr{so}_8$. This was explicitly studied in \cite{jia-affine-closure, kazhdan-minimal-rep-SO8, levasseur-stafford}.
        \item Let $\SL_2 \subseteq \SL_3 \subseteq \SO_8 = \ld{K}$. Again from \cite[Table 1]{deligne-gross-exceptional-series}, one finds that the centralizer of $\SL_2$ is $\SL_2^{\times 3}$, and $(\SL_2, \SL_2^{\times 3})$ forms a dual pair in $\SO_8$. (It corresponds to $H = \PGL_2^\mathrm{diag} \subseteq \PGL_2^{\times 3} = G$.) Interpreting \cite[Equation 4.1]{moore-tachikawa} correctly, one is led to the prediction that $\ld{\cM}^\ddag$ can again be identified with the closure of the minimal nilpotent orbit in $\fr{so}_8$. See \cite[Remark 3.10]{triple-product-cayley}.
        \item Let $\SL_3 \subseteq \mathrm{E}_6 = \ld{K}$ as in \cite[Equation 1]{deligne-gross-exceptional-series}. Again from \cite[Table 1]{deligne-gross-exceptional-series}, one finds that the centralizer of $\SL_3$ is $\SL_3^{\times 2}$, and $(\SL_3, \SL_3^{\times 2})$ forms a dual pair in $\mathrm{E}_6$. (It corresponds to $H = \PGL_3^\mathrm{diag} \subseteq \PGL_3^{\times 2} = G$.) Using \cite[Section 4.2]{moore-tachikawa}, proved mathematically in \cite[Section 5(iv)]{bfn-ring-objects}, one finds that $\ld{\cM}^\ddag$ can be identified with the closure of the minimal nilpotent orbit in $\mathrm{E}_6$.
    \end{itemize}
    The preceding \cref{ex: ggp} also fits into this paradigm, up to some ``metaplectic'' correction. Namely, $(\SO_{2n}, \Sp_{2n})$ forms a dual pair in $\Sp_{4n^2}$, and the closure of the minimal nilpotent orbit of $\Sp_{4n^2}$ can be identified with $\AA^{4n^2}\mmod \{\pm 1\}$ (i.e., the image of the moment map $\AA^{4n^2} \to \fr{sp}_{4n^2}$). Up to this order two quotient, one can therefore identify the closure of the minimal nilpotent orbit of $\Sp_{4n^2}$ with the tensor product of the standard representations of $\SO_{2n}$ and $\Sp_{2n}$, which is what appeared in \cref{ex: ggp}.
\end{remark}

To generalize the discussion of this section to $\ku$-theoretic coefficients, we need some preliminary results.

%% file: equivalences/ku-satake.tex
\subsection{$\ku$-theoretic derived geometric Satake}

Recall from the proof of \cref{thm: derived satake} that the key step in the argument, once given \cref{thm: abelian satake}, is \cref{thm: classical homology loops G}. Our goal in this section is to prove a $\ku$-theoretic analogue of \cref{thm: classical homology loops G} in the case when $G$ is assumed to be simply-laced and connected. Throughout, $\GG_\beta^\vee$ will denote the Cartier dual of $\GG_\beta$.
\begin{definition}\label{def: ku loop space}
    Let $X$ be a (possibly graded) scheme over a commutative ring $R$. The \textit{$\ku$-loop space} $\cL_\beta X$ of $X$ is defined to be the graded $R[\beta]$-scheme $\Map(B\GG_\beta^\vee, X_{R[\beta]})$.
\end{definition}
The following is a slight variant of the main result of \cite{toen-hkr} (see also \cite[Corollary 6.1]{moulinos-formal-group}).
\begin{lemma}\label{lem: ku loops and HKR}
    Let $X$ be a derived scheme over a commutative ring $R$. Then the pushforward of the structure sheaf along the canonical map $\cL_\beta(X)/\GG_m \to \spec(R[\beta])/\GG_m$ corresponds (under the equivalence between quasicoherent sheaves on $\spec(R[\beta])/\GG_m$ and filtered $R$-modules) to the Hochschild-Kostant-Rosenberg filtration on the Hochschild homology $\HH(X/R)$.
\end{lemma}
\begin{proof}
    In \cite{toen-hkr}, it is shown that if $R$ is a $\Z_p$-algebra, $W$ is the ring scheme of $p$-typical Witt vectors, and $F$ is its Frobenius, the pushforward of the structure sheaf along the canonical map $\Map(W[F = \beta^{p-1}], X)/\GG_m \to \spec(R[\beta])/\GG_m$ corresponds to the Hochschild-Kostant-Rosenberg filtration on $\HH(X/R)$. It therefore suffices to identify $\GG_\beta^\vee$ with a completion of $W[F = \beta^{p-1}]$. As shown in \cite[Proposition C.6]{thh-xn}, there is an isomorphism
    \begin{equation}\label{eq: G-beta-dual}
        \GG_\beta^\vee \cong \spf \Z_p\left[\beta, \tfrac{y(y-\beta)\cdots(y-(n-1)\beta)}{n!}\right]^\wedge
    \end{equation}
    where the element $y$ is primitive (i.e., the coproduct sends $y \mapsto y \otimes 1 + 1 \otimes y$) and lives in weight $2$. Here, the completion is taken with respect to the $\beta$-deformed divided power filtration (i.e., with respect to $\frac{1}{n!} \prod_{j=0}^{n-1} (y-j\beta)$ for $n\geq 1$). The desired identification with the completion of $W[F = \beta^{p-1}]$ is now given by \cite[Remark C.7]{thh-xn}.
    Using the arithmetic fracture square, it only remains to prove the lemma when $R$ is a $\QQ$-algebra. In this case, $\GG_\beta^\vee$ is isomorphic to $\spf \QQ[\beta]\pw{y} \cong \hat{\GG}_a$, from which the desired result follows since the the Hochschild-Kostant-Rosenberg filtration on $\HH(X/R)$ splits, and the Hodge complex of $X$ over $R$ can be identified with the global sections of the mapping stack $\Map(B\hat{\GG}_a, X)$.
\end{proof}
\begin{definition}\label{def: Hbeta defn}
    Let $H$ be a graded algebraic group over a commutative ring $R$. Let $H_\beta$ denote the graded group scheme over $R[\beta]$ given by $\Hom(\GG_\beta^\vee, H_{R[\beta]})$. Note that there is a canonical action of $H_{R[\beta]}$ on $H_\beta$ by conjugation, and the quotient stack $H_\beta/H_{R[\beta]}$ is isomorphic to $\cL_\beta(BH) = \Map(B\GG_\beta^\vee, BH_{R[\beta]})$.
\end{definition}
\begin{lemma}\label{lem: Hbeta invert and mod beta}
    If $H$ is a graded algebraic group, there is an isomorphism $H_\beta[\beta^{-1}]/\GG_m \cong H$, and a graded isomorphism $H_\beta/\beta \cong \fr{h}(2)$.
\end{lemma}
\begin{proof}
    Since $\Map(B\widehat{\GG_a^\sharp(-2)}, BH) \cong \fr{h}(2)/H$ by \cref{ex: BGa sharp}, and $\Map(B\Z, BH) \cong H/H$, it suffices to show that $\GG_\beta^\vee[\beta^{-1}]/\GG_m \cong \Z$, while $\GG_\beta^\vee/\beta$ is isomorphic to the completion $\widehat{\GG_a^\sharp(-2)}$ of the PD-hull of the origin in $\GG_a(-2)$ at the divided power filtration. This in turn follows from the fact that $\GG_\beta[\beta^{-1}]/\GG_m \cong \GG_m$ and $\GG_\beta/\beta \cong \GG_a(2)$, and that $\Z$ (resp. $\GG_a^\sharp(-2)$) is the Cartier dual of $\GG_m$ (resp. $\GG_a(2)$).
\end{proof}
\begin{remark}\label{rmk: morava k-theory Gbeta}
    In \cref{def: ku loop space}, there was no reason to restrict to considering maps out of $B\GG_\beta^\vee$: we could have considered \textit{any} $1$-dimensional group scheme over $\Z[\beta]$ in place of $\GG_\beta$. (This sort of philosophy fits very naturally into \cref{subsec: all the expectations}, which more generally suggests that it should be very interesting to study the \textit{universal} case, where $\Z[\beta]$ is replaced with the Lazard ring carrying the universal formal group law itself.) We can even input $1$-dimensional \textit{formal} group schemes into the above cosntruction. For instance, a particularly important example which arises naturally in chromatic homotopy theory is the following. Fix a prime $p$, and consider the formal group over $\QQ[\beta]$ whose logarithm is given by the invertible ``$p$-typical polylogarithmic'' power series
    $$\ell_F(x) = \sum_{j\geq 0} \beta^{p^{nj}-1} \tfrac{x^{p^{nj}}}{p^j}.$$
    Here, the class $x$ lives in weight $-2$.
    That the power series $F(x,y) = \ell_F^{-1}(\ell_F(x) + \ell_F(y))$ has coefficients in $\Z_{(p)}[\beta]$ is a consequence of Hazewinkel's functional equation lemma \cite[Section I.2]{hazewinkel-book}; write $\widehat{\GG}_{k_\Z(n)}$ to denote the associated formal group law over $\Z_{(p)}[\beta]$. When base-changed to $\FF_p[\beta]$, we will denote it by $\widehat{\GG}_{k(n)}$. The Cartier dual of $\widehat{\GG}_{k(n)}$ was computed in \cite[Example 4.5.14]{generalized-n-series}, where it was shown that 
    $$\widehat{\GG}_{k(n)}^\vee \cong \spec \FF_p[\beta][y_0, y_1, \cdots]/(y_{n+j-1}^p - \beta^{p^j(p^n-1)} y_j).$$
    Here, the classes $y_j$ live in weight $2p^j$. Observe that for $\beta = 0$, one recovers the Cartier dual $\GG_a^\sharp$ of $\widehat{\GG}_a$. For $n = 1$, the formal group $\widehat{\GG}_{k(1)}$ is isomorphic to the $p$-typification of $\widehat{\GG}_\beta$.

    To connect this to chromatic homotopy theory, note that $\ell_F(x)$ only depends on $\beta$ through $\beta^{p^n - 1}$; so $\widehat{\GG}_{k_\Z(n)}$ is in fact defined over $\Z_{(p)}[\beta^{p^n - 1}]$. The class $\beta^{p^n - 1}$ is often denoted $v_n$, and the resulting formal group law over $\Z_{(p)}[v_n]$ is the one associated to the complex orientation of (a form of) integral Morava K-theory of height $n$. When $H$ is a group scheme over $\FF_p$, the group scheme $\Hom(\widehat{\GG}_{k(n)}^\vee, H)$ is closely related to the combinatorial constructions of \cite{generalized-n-series}, and we expect it to capture a lot of interesting aspects of the modular representation theory of $H$. For instance, when $n=1$, it essentially reduces to the group scheme $H_\beta$ studied in the present section (the only difference is the completion of $\GG_\beta$). In general, it is an interesting deformation of the Lie algebra $\fr{h}(2)$ to $\FF_p[\beta]$. We will not discuss this generalization of $H_\beta$ further here, but plan to in \cite{Eodd-and-quantizations}.
\end{remark}
\begin{lemma}\label{lem: Lie and hom to Ga}
    Let $R$ be a commutative ring, and let $K$ be a commutative group scheme over $R$ with Cartier dual $K^\vee$. Then there is an isomorphism $\Hom(K^\vee, \GG_a) \cong \Lie(K)$.
\end{lemma}
\begin{lemma}\label{lem: B beta}
    Let $B\subseteq \GL_2$ denote the Borel subgroup of upper-triangular matrices, graded by the action of $2n\rho: \GG_m \to \GL_2$. If $R_\ast$ is a graded $\Z[\beta]$-algebra, the group scheme $B_\beta(R_\ast)$ is isomorphic to the subgroup of $B(R_\ast)$ consisting of matrices of the form $\begin{psmallmatrix}
        1+\beta x & \beta y\\
        0 & 1+\beta w
    \end{psmallmatrix}$, where $x,w\in \GG_\beta(R) \subseteq R_{-2}$ and $y\in R_{2n-2}$.
\end{lemma}
\begin{proof}
    There is a graded extension
    $$\GG_a(-2n) \to B \to \GG_m^2,$$
    which implies that there is a graded extension
    $$\GG_a(-2n)_\beta \to B_\beta \to (\GG_m^2)_\beta.$$
    By construction, $\GG_{m,\beta}$ is the Cartier dual of $\GG_\beta^\vee$, i.e., $\GG_{m,\beta} \cong \GG_\beta$. Moreover, \cref{lem: Lie and hom to Ga} gives an isomorphism $\GG_a(-2n)_\beta \cong \Lie(\GG_\beta)(-2n) \cong \GG_a(2-2n)$. It follows that there is an extension
    $$\GG_a(2-2n) \to B_\beta \to \GG_\beta^2.$$
    This extension precisely classifies matrices of the form $\begin{psmallmatrix}
        1+\beta x & \beta y\\
        0 & 1+\beta w
    \end{psmallmatrix}$ with $x,w\in \GG_\beta(R) \subseteq R_{-2}$ and $y\in R_{2n-2}$.
\end{proof}
In the remainder of this section, we will always invert the order $N = |W|$ of the Weyl group $W$, and write $\Z' = \Z[1/N]$; in particular, $\pi_\ast \ku \cong \Z'[\beta]$. 
\begin{lemma}\label{lem: ramified W-cover}
    Let $\ld{G}$ be a reductive group over a commutative ring $R$. Let $\ld{G}_\beta^\reg$ denote the open subscheme consisting of those elements $x\in \ld{G}_\beta$ such that the centralizer $Z_{\ld{G}}(x) \subseteq \ld{G}$ has minimal dimension (i.e., the rank of $\ld{G}$). Similarly, let $\ld{B}_\beta^\reg$ denote $\ld{B}_\beta \cap \ld{G}_\beta^\reg\subseteq \ld{G}_\beta$. Then the morphism $\ld{B}_\beta^\reg/\ld{B} \to \ld{G}_\beta^\reg/\ld{G}$ is a ramified $W$-Galois cover.
\end{lemma}
\begin{definition}\label{def: beta deformed kostant}
    Let $e\in \ld{G}$ be a principal unipotent element, so that $e$ defines a homomorphism $\SL_2 \to \ld{G}$ by the Jacobson-Morozov theorem. This homomorphism is $\GG_m$-equivariant if $\ld{G}$ (resp. $\SL_2$) is graded by $2\rho: \GG_m \to \ld{G}$ (resp. the restriction of the $2\rho$-grading on $\ld{G}$ to $\SL_2$). There is an induced homomorphism $B_{\SL_2,\beta} \to \ld{G}_\beta$, where $B_{\SL_2}\subseteq \SL_2$ is the Borel subgroup of upper-triangular matrices. Let $\begin{psmallmatrix}
        1 & \beta\\
        0 & 1
    \end{psmallmatrix}\in B_{\SL_2,\beta}$ denote the element defined by \cref{lem: B beta}; the image of this element in $\ld{G}_\beta$ will be denoted by $e_\beta$.

    Let $\tilde{\kappa}: \ld{T}_\beta \to \ld{B}_\beta$ denote the map sending $x \mapsto e_\beta x$, so that $\tilde{\kappa}$ induces a map $\ld{T}_\beta \to \ld{B}_\beta/\ld{B}$ (which will also be denoted by $\tilde{\kappa}$). We will refer to this as the \textit{$\beta$-deformed Kostant slice}. It is not difficult to see that $\tilde{\kappa}: \ld{T}_\beta \to \ld{B}_\beta/\ld{B}$ is $W$-equivariant, so \cref{lem: ramified W-cover} implies that the composite
    $$\ld{T}_\beta \xar{\tilde{\kappa}} \ld{B}_\beta/\ld{B} \to \ld{G}_\beta/\ld{G}$$
    descends to a morphism $\ld{T}_\beta\mmod W \to \ld{G}_\beta/\ld{G}$, which we will denote by $\kappa$. The map $\kappa: \ld{T}_\beta\mmod W \to \ld{G}_\beta/\ld{G}$ will also be called the \textit{$\beta$-deformed Kostant slice}. It defines a graded map $\ld{T}_\beta\mmod W \to \ld{G}(-2\rho)_\beta/\ld{G}(-2\rho)$. (See \cite[Section 2.1]{bfm} for a related construction. Note that if $\ld{G}$ is simply-connected, \cite{steinberg-slice} gives an identification $\ld{G}\mmod \ld{G} \cong \ld{T}\mmod W$, as well as the construction of a slice $\ld{T}\mmod W \to \ld{G}$. The above construction is weaker than this, in the sense that the slice lands in the quotient stack $\ld{G}/\ld{G}$, and not $\ld{G}$ itself.)
\end{definition}
Recall that if $G$ is simply-laced and connected, its Langlands dual group $\ld{G}$ is isogenous to $G$ itself. In particular, the action of $G$ on itself by conjugation induces an action of $\ld{G}$ on $G$. For instance, if $G$ is simply-connected, $\ld{G}$ is the quotient of (the Chevalley split form of) $G$ by its center; and the action of $G$ on itself by conjugation descends to an action of $\ld{G} = G/Z(G)$. Similarly, if $G$ is adjoint, $\ld{G}$ is a $\pi_1(G)_\mathrm{tors}$-cover of $G$, and so the action of $G$ on itself by conjugation restricts to an action of $\ld{G}$ on $G$.
\begin{theorem}\label{thm: ku homology LG and langlands mod center}
    Let $G$ be a simply-laced and connected compact Lie group with associated reductive group $G_\cc$ over $\cc$, so that $(G \times G)/Z(G)^\mathrm{diag}$ acts on $\cL G$ by left and right translation. 
    Then there is a graded $\Z'[\beta]$-linear isomorphism
    $$\spec \ku^{G/Z(G)}_\ast(\Omega G) \cong \ld{T}^\mathrm{ad}_\beta \mmod W \times_{\ld{G}^\mathrm{ad}(-2\rho)_\beta/\ld{G}(-2\rho)} \ld{T}^\mathrm{ad}_\beta \mmod W$$
    of group schemes over $\ld{T}^\mathrm{ad}_\beta\mmod W = \spec \ku_{G/Z(G)}^\ast(\ast)$.
\end{theorem}
\begin{corollary}\label{cor: sc ku homology LG and langlands}
    In the setup of \cref{thm: ku homology LG and langlands mod center}, there is a graded $\Z'[\beta]$-linear isomorphism
    $$\spec \ku^{G}_\ast(\Omega G) \cong T_\beta \mmod W \times_{G(-2\rho)_\beta/\ld{G}(-2\rho)} T_\beta \mmod W$$
    of group schemes over $T_\beta\mmod W = \spec \ku_{G}^\ast(\ast)$, where $G$ on the right-hand side denotes the split form of the compact Lie group over $\Z'$, and $\ld{G}$ (which is a central quotient of $G$) acts on $G$ by conjugation.
\end{corollary}
\begin{remark}
    \cref{cor: sc ku homology LG and langlands} is a simultaneous generalization of \cite[Proposition 4.1.5 and Theorem 4.2.5]{grg-reg} (which is in turn related to \cite[Theorem 2.12 and Theorem 2.15]{bfm}).
\end{remark}
The following result is essentially \cite[Step II of Theorem 6.1]{homology-langlands}.
\begin{lemma}\label{lem: centralizer borel flat}
    The scheme $\ld{T}_\beta \times_{\ld{B}_\beta/\ld{B}} \ld{T}_\beta$ is flat over $\ld{T}_\beta$ after inverting $|W|$.
\end{lemma}
\begin{proof}
    There is a closed immersion
    $$\ld{T}_\beta \times_{\ld{B}_\beta/\ld{B}} \ld{T}_\beta \subseteq \ld{T}_\beta \times \ld{B},$$
    which exhibits the left-hand side as the subgroup of those $(x, g)$ such that $g$ stabilizes $\kappa(x)$; in particular, it is cut out by $\dim \ld{N}$ equations, so that the fibers of the projection $\ld{T}_\beta \times_{\ld{B}_\beta/\ld{B}} \ld{T}_\beta \to \ld{T}_\beta$ are at most $\dim\ld{B} - \dim\ld{N} = \dim\ld{T}$-dimensional. To prove that $\ld{T}_\beta \times_{\ld{B}_\beta/\ld{B}} \ld{T}_\beta$ is flat over $\ld{T}_\beta$, it suffices to show that all the fibers of the projection
    $$\ld{T}_\beta \times_{\ld{B}_\beta/\ld{B}} \ld{T}_\beta \to \ld{T}_\beta$$
    are $\dim\ld{T}$-dimensional. 

    Let $g\in \ld{B}$ and $x\in \ld{T}_\beta$. Since $\kappa(x) = e_\beta x$, we have $\Ad_g \kappa(x) = \Ad_g(e_\beta) \Ad_g(x) \in \ld{B}_\beta \subseteq \ld{B}_{\Z[\beta]}$. Since $x$ is semisimple, the same is true of $\Ad_g(x)$. It therefore suffices to show that the subgroup $Z_{\ld{B}}(e_\beta) = \{g\in \ld{B} | \Ad_g(e_\beta) = e_\beta\} \subseteq \ld{B}$ is $\dim\ld{T}$-dimensional. If $\Ad_g(e_\beta) = e_\beta$, then $\Ad_g(e) = e$, where $e\in \ld{\g}$ is the associated nilpotent element; this implies that $\dim Z_{\ld{B}}(e_\beta) \leq \dim Z_{\ld{B}}(e)$. Therefore, it suffices to show that the centralizer $Z_{\ld{B}}(e)$ is $\dim\ld{T}$-dimensional, which is even true with a smaller set of primes inverted (see \cite{keny-nilp}).
\end{proof}
\begin{proof}[Proof of \cref{thm: ku homology LG and langlands mod center}]
    It suffices to show that there is a graded isomorphism
    $$T_\beta\mmod W \times_{G(-2\rho)_\beta/\ld{G}(-2\rho)} T_\beta \mmod W \cong \spec \ku^{G}_\ast(\Omega G).$$
    By \cref{prop: Weyl invts}, $\ku^{G}_\ast(\Omega G) \cong \ku^T_\ast(\Omega G)^W$, and so by definition of $\kappa$, it suffices to show that there is a $W$-equivariant graded isomorphism
    $${T}_\beta \times_{B_\beta/\ld{B}} {T}_\beta \cong \spec \ku^T_\ast(\Omega G),$$
    for which it in turn suffices to show that there is a $W$-equivariant graded isomorphism
    \begin{equation}\label{eq: Tad and ku homology}
        \ld{T}^\mathrm{ad}_\beta \times_{\ld{B}^\mathrm{ad}_\beta/\ld{B}} \ld{T}^\mathrm{ad}_\beta \cong \spec \ku^{T/Z(G)}_\ast(\Omega G)
    \end{equation}
    over $\ld{T}^\mathrm{ad}_\beta$.
    
    For simplicity, we will ignore gradings in the following discussion. By equivariant formality for $\Omega G$, the scheme $\spec \ku^{T/Z(G)}_\ast(\Omega G)$ is flat over $\spec \pi_\ast \ku_{T/Z(G)} \cong \ld{T}_\beta$. 
    The scheme $\ld{T}_\beta \times_{\ld{B}_\beta/\ld{B}} \ld{T}_\beta$ is also flat over $\ld{T}_\beta$ by \cref{lem: centralizer borel flat}. Therefore, the argument of \cite[Section 4.3]{bfm} reduces us to checking the isomorphism \cref{eq: Tad and ku homology} in the case when $G$ (equivalently $\ld{G}$) has semisimple rank $1$ (this, as usual, is through a Hartogs/codimension $2$ argument). For this, it in turn suffices to prove \cref{eq: Tad and ku homology} when $\ld{G} = \GL_2$, $\SL_2$ and $\PGL_2$. (Below, we will calculate both sides of \cref{eq: Tad and ku homology} as schemes over $\ld{T}_\beta$. But checking that the isomorphism \cref{eq: Tad and ku homology} is one of group schemes is not difficult: using flatness over $\ld{T}_\beta$, one observes that the coproduct/group structure is determined by the coproduct/group structure over the complement $\punc{\ld{T}}_\beta$ of all root hypersurfaces. But then this reduces to the case when $G$ is a torus itself, where the isomorphism is evidently one of group schemes.)

    Let $\ld{G} = \PGL_2$, so $G = \SL_2$, and $Z(G)$ is trivial. Recall that there is a closed immersion
    $$\ld{T}_\beta \times_{\ld{B}_\beta/\ld{B}} \ld{T}_\beta\subseteq \ld{T}_\beta \times \ld{B},$$
    which exhibits the left-hand side as the subgroup of those $(x, g)$ such that $g$ stabilizes $\kappa(x)$. 
    Since $\kappa(x)\in \ld{B}_\beta$ is the matrix
    $$\kappa(x) = \begin{psmallmatrix}
        1 & \beta \\
        0 & 1
    \end{psmallmatrix} \begin{psmallmatrix}
        1 & 0 \\
        0 & 1+\beta x
    \end{psmallmatrix} = \begin{psmallmatrix}
        1 & \beta(1+\beta x)\\
        0 & 1+\beta x
    \end{psmallmatrix},$$
    it follows that that if $g = \begin{psmallmatrix}
        a & b\\
        0 & 1
    \end{psmallmatrix} \in \ld{G} = \PGL_2$, we have 
    $$\Ad_g \kappa(x) = \begin{psmallmatrix}
        1 & a \beta (1 + \beta {x}) - b\beta x\\
        0 & 1+\beta x
    \end{psmallmatrix}.$$
    It follows that $g$ fixes $\kappa(x)$ if and only if
    $$a \beta (1 + \beta {x}) - b\beta x = \beta (1 + \beta x),$$
    i.e., if and only if
    $$b = \tfrac{a-1}{x} (1+\beta x).$$
    It follows that there is an isomorphism
    $$\ld{T}_\beta \times_{\ld{B}_\beta/\ld{B}} \ld{T}_\beta \cong \spec \Z'[\beta, x, \tfrac{1}{1+\beta x}, a^{\pm 1}, \tfrac{a - 1}{x}].$$
    Now using \cite[Theorem 3.2.12]{grg-reg}\footnote{Note that \cite[Theorem 3.2.12]{grg-reg} asks that $\GG$ be an \textit{oriented} group scheme over $A$, but this is in fact not necessary: it suffices that $\GG$ be preoriented, so that \cref{prop: gkm} continues to hold. In this case, the desired preorientation of $\GG_{\ku,\beta}$ is given by \cref{prop: Gbeta is preoriented}.} or \cref{thm: homology of loops SV} with $V$ being the weight $2$ representation of $S^1$, one sees that 
    $$\ku^{S^1}_\ast(\Omega G) \cong \Z'[\beta, x, \tfrac{1}{1 + \beta x}, a^{\pm 1}, \tfrac{a-1}{[2](x)}],$$
    which implies that
    $$\ku^{S^1/(\Z/2)}_\ast(\Omega G) \cong \Z'[\beta, [2](x), \tfrac{1}{1 + \beta [2](x)}, a^{\pm 1}, \tfrac{a-1}{[2](x)}].$$
    Therefore, $\spec \ku^{T/Z(G)}_\ast(\Omega G)$ is isomorphic to $\ld{T}_\beta \times_{\ld{B}_\beta/\ld{B}} \ld{T}_\beta$, as desired. 
    A similar calculation proves \cref{eq: Tad and ku homology} when ${G} = \GL_2, \PGL_2$ (note that in this case $\pi_1(G) \cong \Z/2$ is not zero; but the result is still true by direct calculation).
\end{proof}
\begin{example}
    It is important that the order of $W$ is inverted. For instance, \cref{thm: ku homology LG and langlands mod center} implies that $\spec \ku_\ast(\Omega G)$ is a subgroup of $\ld{G}$, given by the centralizer of $e_\beta\in \ld{G}_\beta$. Rationally (in fact, after inverting $|W|$), this is isomorphic to the centralizer of $e \in \ld{\g}$, which identifies with $\spec \H_\ast(\Omega G; \QQ)[\beta]$ by \cref{thm: classical homology loops G}. Note, however, that $\ku_\ast(\Omega G)$ is (unsurprisingly) \textit{not} isomorphic to $\H_\ast(\Omega G; \Z)[\beta]$, even as algebras. For instance, let $G = \G_2$. Then \cite{bott-space-of-loops} shows that
    $$\H_\ast(\Omega \G_2; \Z) \cong \Z[a, b, c]/(a^2 - 2b),$$
    where $a$ is in weight $2$, $b$ is in weight $4$, and $c$ is in weight $10$. 
    On the other hand, a slight refinement of \cite[Proposition 7.1]{clarke-k-thy-loops-lie} shows that
    $$\ku_\ast(\Omega \G_2) \cong \Z[\beta, a, b, c]/(a^2 - 2b - \beta a).$$
    Of course, when $2$ is inverted, we may express $b = \tfrac{a^2 - \beta a}{2}$, and $\ku_\ast(\Omega \G_2)[1/2] \cong \H_\ast(\Omega \G_2; \Z[1/2])[\beta]$.
\end{example}


We now make the following rather contrived definition.
\begin{construction}\label{cstr: Shv-Sat LG/H with ku}
    Let $G$ be a simply-laced connected compact Lie group, so that \cref{thm: ku homology LG and langlands mod center} gives a homomorphism
    $$\spec \ku^G_\ast(\Omega G) \hookrightarrow \ld{G}(-2\rho) \times T_\beta \mmod W$$
    of group schemes over $T_\beta\mmod W$.
    Let $H_\cc\subseteq G_\cc$ be a connected {reductive} subgroup, and let $\cM_H$ be the graded $\Z[\beta]$-scheme defined in \cref{notn: MG defn}. Then there is a homomorphism (in fact, closed imersion) of graded group schemes
    $$\spec \ku^H_\ast(\Omega (G/H)) \to \ld{G}(-2\rho) \times \cM_H$$
    over $\cM_H \cong \spec \ku^\ast_H(\ast)$ constructed as the following composite:
    \begin{align*}
        \spec \ku^H_\ast(\Omega (G/H)) & \to \spec \ku^H_\ast(\Omega G) \cong \spec \ku^G_\ast(\Omega G) \times_{\spec \ku_G^\ast(\ast)} \spec \ku_H^\ast(\ast) \\
        & \hookrightarrow (\ld{G}(-2\rho) \times T_\beta \mmod W) \times_{T_\beta \mmod W} \cM_H \cong \ld{G}(-2\rho) \times \cM_H.
    \end{align*}
    Let $A_\beta$ denote the $\co_{\cM_H}$-algebra of regular functions on the quotient $(\ld{G} \times \cM_H)/\spec \ku^H_\ast(\Omega (G/H))$ over $\cM_H$, so that $A_\beta$ admits a canonical grading, as well as a canonical action of $\ld{G}(-2\rho)$. Define $\Shv^{c,\Sat}_{G\pw{t}}(G\ls{t}/H\ls{t}; \ku)^\faux$ to be the $\sh^{1/2} \Z'[\beta]$-linear $\infty$-category $\Perf(\sh^{1/2}\spec A_\beta/\ld{G}(-2\rho))$.
\end{construction}
The reason that $H_\cc$ is assumed to be reductive is precisely thanks to the proof of \cref{thm: ordinary homology criterion satake}, which implies \cref{prop: defns ShvSat ku and Q} below; this result ensures consistency in notation (and showing that \cref{thm: ku derived satake} implies \cref{thm: derived satake} at least in the simply-laced and connected case):
\begin{prop}\label{prop: defns ShvSat ku and Q}
    Let $G$ be a simply-laced and connected compact Lie group, let $G_\cc$ be the associated algebraic group over $\cc$, and let $H_\cc\subseteq G_\cc$ be a connected {reductive} subgroup which is optimal in the sense of \cref{hypothesis: rank 1 weakly placid}.
    Let $\sh^{1/2} \Z'[\beta] \to \QQ$ denote the $\Eoo$-map given by rationalization and sending $\beta \mapsto 0$. Then there is an equivalence 
    $$\Shv^{c,\Sat}_{G\pw{t}}(G\ls{t}/H\ls{t}; \ku)^\faux \otimes_{\sh^{1/2} \Z'[\beta]} \QQ \simeq \Shv^{c,\Sat}_{G\pw{t}}(G\ls{t}/H\ls{t}; \QQ).$$
\end{prop}
\begin{remark}\label{rmk: true shvsat with ku is hard}
    \cref{cstr: Shv-Sat LG/H with ku} is \textit{not} a good definition (hence the ``faux''). For instance, it is not defined intrinsically to the $G$-action on $G/H$, and instead uses the algebra of functions on ${\ld{G}}$; in other words, as we said in the introduction, it is defined by ``playing games'' with the subcategory of locally constant sheaves. Its only saving grace is \cref{prop: defns ShvSat ku and Q}.

    Instead, \cref{cstr: Shv-Sat LG/H with ku} is intended to be a replacement of \cref{def: Shv-Sat LG/H} in the case where one does not have an analogue of \cref{thm: abelian satake}. It might be the case that there is a well-behaved $\ku_G$-linear $\infty$-category $\tilde{\Shv}^{c}_{G\pw{t}}(G\ls{t}/H\ls{t}; \ku)^\faux$ of constructible sheaves of $G\pw{t}$-equivariant $\ku$-modules defined on ind-finite $G\pw{t}$-spaces such as $G\ls{t}/H\ls{t}$. In the nonequivariant case, there \textit{is} a well-behaved notion of constructible sheaves of $\ku$-modules defined using exit path categories \`a la \cite[Appendix A]{HA}. Roughly, then, $\Shv^{c,\Sat}_{G\pw{t}}(G\ls{t}/H\ls{t}; \ku)^\faux$ would be an ``associated graded for the Postnikov filtration'' on $\tilde{\Shv}^{c}_{G\pw{t}}(G\ls{t}/H\ls{t};\ku)$.
    
    Namely, suppose that there was a $\ku_G$-linear $\infty$-category ${\Shv}^{c}_{G\pw{t}}(G\ls{t}/H\ls{t};\ku)$, and an $\Eoo$-coalgebra $\tilde{\ld{\cR}}$ in $\Alg_\E{2}({\Shv}^{c}_{(G\times G)\pw{t}}(\Gr_G;\ku))$ which lifted the regular sheaf (corresponding to $\co_{\ld{G}} \in \Perv_{G\pw{t} \times G\pw{t}}(G\ls{t}; \Z)$). Let $\tilde{\cC}$ denote the de-equivariantization of ${\Shv}^{c}_{G\pw{t}}(G\ls{t}/H\ls{t};\ku)$ with respect to this $\E{2}$-Hopf algebra, and let $\IC_0 \in {\Shv}^{c}_{G\pw{t}}(G\ls{t}/H\ls{t};\ku)$ denote the $!$-pushforward of the constant sheaf in ${\Shv}^{c}_{G\pw{t}}(G\pw{t}/H\pw{t};\ku) \simeq {\Shv}^{c}_{H\pw{t}}(\ast;\ku)$.
    Then the full subcategory $\cC$ of $\tilde{\cC}$ generated by $\IC_0 \star \tilde{\ld{\cR}}$ can be identified with left modules over the $\ku$-algebra $\tilde{A} := \End_{\tilde{\cC}}(\IC_0 \star \tilde{\ld{\cR}})$. Finally, let ${\Shv}^{c,\Sat}_{G\pw{t}}(G\ls{t}/H\ls{t};\ku)$ denote the equivariantization of $\cC$ with respect to the $\E{2}$-Hopf algebra $\tilde{\ld{\cR}}$. If the assumptions of \cref{thm: full faithful} applied (see also \cref{rmk: full faithful is very general}), we could identify $\pi_\ast \tilde{A} = A_\beta$, and the coaction of $\tilde{\ld{\cR}}$ on $\tilde{A}$ induces the $\co_{\ld{G}}$-coaction on $A_\beta$ from \cref{cstr: Shv-Sat LG/H with ku}. In particular, $\Shv^{c,\Sat}_{G\pw{t}}(G\ls{t}/H\ls{t}; \ku)^\faux$ is an ``associated graded'' of ${\Shv}^{c,\Sat}_{G\pw{t}}(G\ls{t}/H\ls{t};\ku)$.

    Of course, it would be ideal to work with ${\Shv}^{c,\Sat}_{G\pw{t}}(G\ls{t}/H\ls{t};\ku)$ itself instead of the \textit{ad hoc} $\infty$-category of \cref{cstr: Shv-Sat LG/H with ku}. However, carrying out the above program has proven to be very challenging (for both technical and conceptually interesting reasons), and \cref{cstr: Shv-Sat LG/H with ku} is an attempt to salvage the situation somewhat.
\end{remark}

\cref{cor: sc ku homology LG and langlands} implies the following generalization of \cref{thm: derived satake}:
\begin{theorem}\label{thm: ku derived satake}
    Let $G$ be a simply-laced and connected compact Lie group.
    There is an $\E{2}$-monoidal equivalence of $\sh^{1/2} \Z'[\beta]$-linear $\infty$-categories
    $$\Shv^{c,\Sat}_{(G\times G)\pw{t}}(G\ls{t}; \ku)^\faux \simeq \Perf(\sh^{1/2} G(-2\rho)_\beta/\ld{G}(-2\rho)),$$
    where again the symbol $G$ on the right-hand side denotes the split form of the compact Lie group over $\Z'$.
\end{theorem}
\begin{proof}
    By \cref{cstr: Shv-Sat LG/H with ku}, we need to check that there is a graded $\Z'[\beta]$-linear isomorphism
    $$\co_{{G}(-2\rho)_\beta} \cong \co_{(\ld{G}(-2\rho) \times {T}_\beta \mmod W)/({T}_\beta \mmod W \times_{G(-2\rho)_\beta/\ld{G}(-2\rho)} T_\beta \mmod W)},$$
    for which it in turn suffices to check that there is a graded $\Z'[\beta]$-linear isomorphism
    $$\co_{\ld{G}(-2\rho)_\beta} \cong \co_{(\ld{G}(-2\rho) \times \ld{T}_\beta \mmod W)/(\ld{T}_\beta \mmod W \times_{\ld{G}(-2\rho)_\beta/\ld{G}(-2\rho)} \ld{T}_\beta \mmod W)}.$$
    For simplicity, let us momentarily ignore gradings.
    We claim that the $\ld{G}$-orbit of the image of $\kappa: \ld{T}_\beta \mmod W \to \ld{G}_\beta$ is isomorphic to $\ld{G}_\beta^\reg$, so that the right-hand side above is $\co_{\ld{G}_\beta^\reg}$. For this, it suffices to show that the $\ld{B}$-orbit of the image of $\tilde{\kappa}: \ld{T}_\beta \to \ld{B}_\beta$ is isomorphic to $\ld{B}_\beta^\reg$ where $\tilde{\kappa}: \ld{T}_\beta \to \ld{B}_\beta/\ld{B}$ is as in \cref{def: beta deformed kostant}. First, we claim that $\ld{B} \cdot \im(\tilde{\kappa}) \subseteq \ld{B}_\beta^\reg$. Indeed, if $x\in \ld{T}_\beta$, then for a fixed $b\in \ld{B}$, there is an isomorphism $Z_{\ld{B}}(b \cdot \tilde{\kappa}(x)) \to Z_{\ld{B}}(\tilde{\kappa}(x))$ sending $g\mapsto bgb^{-1}$.
    It is not difficult to show that both $\ld{B} \cdot \im(\tilde{\kappa})$ and $\ld{B}_\beta^\reg$ are flat over $\Z'[\beta]$ (the latter because it is an open subscheme of $\ld{B}_\beta$, which is smooth over $\Z'[\beta]$), so it suffices to show that the inclusion $\ld{B} \cdot \im(\tilde{\kappa}) \subseteq \ld{B}_\beta^\reg$ is an isomorphism after inverting $\beta$ and setting $\beta = 0$. By \cref{lem: Hbeta invert and mod beta}, this is equivalent to the following pair of well-known facts: the inclusions $\ld{B} \cdot \im(\tilde{\kappa}: \ld{T} \to \ld{B}) \subseteq \ld{B}^\reg$ and $\ld{B} \cdot \im(\tilde{\kappa}: \ld{\fr{t}} \to \ld{\fr{b}}) \subseteq \ld{\fr{b}}^\reg$ are isomorphisms. (This reduces to the Jordan decomposition and the fact that $e_{\beta=1}$ and $e_{\beta = 0}$ are representatives for the open $\ld{B}$-orbits in the $\ld{N}$ and $\fr{n}$, respectively.)
    
    Since $\ld{G}_\beta$ is normal and irreducible, the desired isomorphism $\co_{\ld{G}(-2\rho)^\reg_\beta} \cong \co_{\ld{G}(-2\rho)_\beta}$ is therefore a consequence of the following claim (and the algebraic Hartogs lemma): the complement of $\ld{G}(-2\rho)^\reg_\beta\subseteq \ld{G}(-2\rho)_\beta$ is of codimension $\geq 2$. This complement is flat over $\Z'[\beta]$, and so it suffices to check the claim after inverting $\beta$ and setting $\beta = 0$. Again, this reduces to the well-known facts (see e.g., \cite[Theorem 4.13]{humphreys-conjugacy-classes} for the group case) that the closed subschemes in $\ld{G}$ and $\ld{\g}$ of irregular elements has complement of codimension $\geq 2$.
\end{proof}
\begin{remark}\label{rmk: HKR and ku satake}
    Since $\ku$ interpolates between $\Z$ and $\KU$, it can be understood as encoding the $\beta$-adic filtration on $\KU$. In the same way, the right-hand side of \cref{thm: ku derived satake} interpolates between $\Perf(\ld{\g}[2-2\rho]/\ld{G}[-2\rho])$ and $\Perf(G/\ld{G})$. Therefore, \cref{thm: ku derived satake} along with \cref{lem: ku loops and HKR} say that the $\beta$-adic filtration on the topological/A-side corresponds to the Hochschild-Kostant-Rosenberg filtration on the free loop space of $B\ld{G}$ (up to the issue of replacing $\ld{G}$ with its simply-connected cover, and changing gradings by $-2\rho$).
\end{remark}

The next result is true by definition, but we are restating it as such to draw the analogy to \cref{thm: ordinary homology criterion satake}.
\begin{prop}\label{prop: ku rk 1 reg centr --> thm}
    Let $G$ be a simply-laced and connected compact Lie group, and let $H_\cc \subseteq G_\cc$ be a closed connected reductive subgroup.
    Let $\cM_H$ be the graded $\Z[\beta]$-scheme defined in \cref{notn: MG defn}.
    Let $\ld{M}^\ddag_\beta$ denote the scheme $\spec A_\beta$ from \cref{cstr: Shv-Sat LG/H with ku}, so that there is a $\beta$-deformed ``Kostant section'' $\cM_H \to \ld{M}^\ddag_\beta/\ld{G}(-2\rho)$, and an isomorphism
    $$\spec \ku_\ast^{H}(\Omega (G/H)) \cong \cM_H \times_{\ld{M}^\ddag_\beta/\ld{G}(-2\rho)} \cM_H$$
    over $\cM_H$.
    Then there is an equivalence of $\sh^{1/2} \Z'[\beta]$-linear $\infty$-categories
    \begin{equation}\label{eq: display ku rel sat criterion}
        \Shv^{c,\Sat}_{G\pw{t}}(G\ls{t}/H\ls{t}; \ku)^\faux \simeq \Perf(\sh^{1/2}\ld{M}^\ddag_\beta/\ld{G}(-2\rho)).
    \end{equation}

    Suppose, further, that there is a morphism $\mu: \ld{M}^\ddag_\beta/\ld{G}(-2\rho) \to G(-2\rho)_\beta/\ld{G}(-2\rho)$ over $\Z'[\beta]$ such that there is a commutative diagram
    $$\xymatrix{
    \cM_H \ar[r]^-{\kappa_{\ld{M}^\ddag_\beta}} \ar[d] & \ld{M}^\ddag_\beta/\ld{G}(-2\rho) \ar[d]^-\mu \\
    T_{\beta}\mmod W \ar[r]^-{\kappa} & G(-2\rho)_\beta/\ld{G}(-2\rho),
    }$$
    so that there is an induced map
    $$\cM_H \times_{\ld{M}^\ddag_\beta/\ld{G}(-2\rho)} \cM_H \cong \ld{J}_{X,\beta} \to T_\beta\mmod W \times_{G(-2\rho)_\beta/\ld{G}(-2\rho)} T_\beta\mmod W.$$
    If the isomorphism of (b) fits into a commutative diagram
    $$\xymatrix{
    \spec \ku_\ast^H(\Omega(G/H)) \ar[rr]^-\sim_-{\textrm{(b)}} \ar[d] & & \cM_H \times_{\ld{M}^\ddag_\beta/\ld{G}(-2\rho)} \cM_H \ar[d] \\
    \spec \ku_\ast^G(\Omega G) \ar[rr]^-\sim_-{\textrm{\cref{thm: derived satake}}} & & T_\beta\mmod W \times_{G(-2\rho)_\beta/\ld{G}(-2\rho)} T_\beta\mmod W,
    }$$
    then the equivalence \cref{eq: display ku rel sat criterion} is equivariant for the left-action of $\Shv_{(G \times G)\pw{t}}^{c,\Sat}(\Gr_G; \ku)^\faux$, identified with $\Perf(\sh^{1/2}G(-2\rho)_\beta/\ld{G}(-2\rho))$ via \cref{thm: ku derived satake}.
\end{prop}
\begin{remark}
    Most of the results in \cref{subsec: using the regular centralizer} continue to hold in the $\ku$-theoretic context, with essentially the same arguments. Assume that $G$ is a simply-laced and connected compact Lie group, and let $H\subseteq G$ be a closed subgroup satisfying the same hypotheses. Let us briefly summarize the resulting picture:
    \begin{itemize}
        \item It will not be clear that $\ld{M}^\ddag_\beta$ is well-behaved (for example, is it of finite type?).
        \item Just as with \cref{conj: JX and GcheckX}, one expects that there is a commutative diagram
        $$\xymatrix{
        \spec \ku^H_\ast(\Omega (G/H)) \ar[r] \ar[d] & \spec \ku^H_\ast(\Omega G; \QQ) \ar[d] \ar[r]^-\sim & \ld{J}_{\ld{G},\beta} \times_{T_\beta\mmod W} \cM_H\\
        \ld{G}_X(-2\rho_{\ld{G}}) \times \cM_H \ar[r] & \ld{G}(-2\rho_{\ld{G}}) \times \cM_H &
        }$$
        of graded group schemes over $\cM_H$, where the homomorphism $\ld{G}_X \to \ld{G}$ is that of \cref{def: SV dual group}, and the vertical maps are closed immersions.
        \item There is a diagram 
        \begin{equation}\label{eq: ku variant J and lag corr}
        \xymatrix{
            & \ld{J}_{X,\beta} \ar[dr] \ar[dl] & & \\
            \cM_H \ar[dr] & & \ld{J}_{\ld{G},\beta} \times_{T_\beta\mmod W} \cM_H \ar[dl] \ar[dr] & \\
            & \ld{J}_{\ld{H},\beta} & & \ld{J}_{\ld{G},\beta},
            }
        \end{equation}
        where the square is Cartesian, and the span at the bottom of the diagram is a Lagrangian correspondence. This identifies with the diagram
        $$\xymatrix@=.75em{
        & \spec \ku^H_\ast(\Omega(G/H)) \ar[dr] \ar[dl] & & \\
        \spec \pi_\ast \ku_H \ar[dr] & & \spec \ku^H_\ast(\Omega G) \ar[dl] \ar[dr] & \\
        & \spec \ku^H_\ast(\Omega H) & & \spec \ku^G_\ast(\Omega G).
        }$$
        The long composite on the right-hand side of the above diagram will be Lagrangian, hence coisotropic (but this has to be interpreted in a derived sense; for example, it need not be a closed immersion!). We will study this phenomenon of coisotropicity (upon completion) in much greater detail in \cref{subsec: ku and lagrangians} from the perspective of Hochschild cohomology and centralizers \`a la \cite[Section 5.3]{HA} and \cite{francis}. See \cref{rmk: coisotropic correspondence}, for instance.
        \item Let $\ld{\cM}_\beta$ denote the affine closure of the quotient of $\ld{G} \times \ld{H} \times \cM_H$ by $\ld{J}_{\ld{G},\beta} \times_{T_\beta\mmod W} \cM_H$. Then \cref{eq: ku variant J and lag corr} implies that there is a diagram
        \begin{equation}\label{eq: ku variant BJ span}
        \xymatrix{
            & \ld{M}^{\ddag,\reg}_\beta/\ld{G} \ar[dr] \ar[dl] & & \\
            \cM_H \ar[dr]_-\kappa & & \ld{\cM}_\beta^{\ddag, \reg}/(\ld{H} \times \ld{G}) \ar[dl] \ar[dr] & \\
            & H_\beta^{\reg}/\ld{H} & & G_\beta^{\reg}/\ld{G},
        }
        \end{equation}
        where the square is Cartesian, and the span at the bottom of the diagram is a ($1$-shifted) Lagrangian correspondence. (The $1$-shifted symplectic structures on ${H}_\beta/\ld{H}$ and ${G}_\beta/\ld{G}$ are discussed in \cref{def: ku-hamiltonian}.)
    \end{itemize}
\end{remark}
\begin{remark}
    It would be interesting to use \cref{prop: ku rk 1 reg centr --> thm} to prove the $\ku$-theoretic analogue of \cref{conj: spherical levi} for spherical Levi subgroups of $G$.
\end{remark}

%% file: each-type/big-table.tex
\begin{landscape}
\begin{table}
\centering
\resizebox{18cm}{!}{
\begin{tabular}{ |c|c|c|c|c|c|c|c| } 
 \hline
 Name & Citation & $X = G/H$ & Type & $\ld{Y}$ & $\ld{G}_X$ grading & ``Normalization'' & Topological phenomenon \\
 \hline
 $A_n$ & \cref{cor: bzsv for CPn} & $\PGL_{n+1}/\GL_n$ & T & $T^\ast(2n) \AA^2(2n,0)$ & $2n\rho_{\ld{G}_X}$ & $\gl_{n-1}[2]\mmod \GL_{n-1}$ & Hopf fibration $S^1 \to \Omega\CP^n \to \Omega S^{2n+1}$ \\
 $B_n$ & \cref{thm: bzsv for Bn} & $\SO_{2n+1}/\SO_{2n}$ & T & $T^\ast(2n) \AA^2(4n-2,0)$ & $(4n-2)\rho_{\ld{G}_X}$ & $\fr{sp}_{2n-2}[2]\mmod \Sp_{2n-2}$ & EHP sequence $S^{2n-1} \to \Omega S^{2n} \to \Omega S^{4n-1}$ \\
 $C_n$ & \cref{thm: bzsv for Cn} & $\Sp_{2n}/(\Sp_2 \times \Sp_{2n-2})$ & T & $T^\ast(4n-4) \AA^2(4n-2,0)$ & $(4n-2)\rho_{\ld{G}_X}$ & $(\fr{sp}_2 \times \fr{sp}_{2n-4})[2]\mmod (\Sp_2 \times \Sp_{2n-4})$ & Hopf fibration $S^3 \to \Omega\HHP^{n-1} \to \Omega S^{4n-1}$ \\
 $D_n$ & \cref{thm: bzsv for Dn} & $\SO_{2n}/\mu_2 \cdot \SO_{2n-1}$ & G & $\sl_2(2n-2-(2n-2)\rho_{\ld{G}_X})$ & $(2n-2)\rho_{\ld{G}_X}$ & $\fr{spin}_{2n-3}[2]\mmod \Spin_{2n-3}$ & James splitting for $\Omega S^{2n-1}$ \\
 $\F_4$ & \cref{thm: bzsv for F4} & $\F_4/\Spin_9$ & T & $T^\ast(16) \AA^2(22, 0)$ & $22\rho_{\ld{G}_X}$ & $\fr{sp}_6[2]\mmod \Sp_6$ & Exceptional Hopf fibration $S^7 \to \Omega \OP^2 \to \Omega S^{23}$ \\
 $\G_2$ & \cref{thm: bzsv for G2} & $\G_2/\SL_3$ & T & $T^\ast(6) \AA^2(10,0)$ & $10\rho_{\ld{G}_X}$ & $\sl_2[2]\mmod \SL_2$ & EHP sequence $S^5 \to \Omega S^6 \to \Omega S^{11}$ \\
 $B_3'$ & \cref{thm: bzsv for B3'} & $\SO_7/\G_2$ & G & $\sl_2(6-6\rho_{\ld{G}_X})$ & $6\rho_{\ld{G}_X}$ & $\fr{sp}_2[2]\mmod \Sp_2$ & James splitting for $\Omega S^7$ \\
 $N(A_1)$ & \cref{rmk: type N PGL2 mod PO2} & $\PGL_2/\mathrm{PO}_2$ & N & $(T^\ast(2) \AA^2(2,0))/(\Z/2)$ & $2\rho_{\ld{G}_X}$ & $0$ & Antipodal action on $S^2$ \\
 $N(B_n)$ & \cref{rmk: SO2n+1 mod normalizer of SO2n} & $\SO_{2n+1}/\N_{\SO_{2n+1}}(\SO_{2n})$ & N & $(T^\ast(2n)\AA^2(4n-2,0))/(\Z/2)$ & $(4n-2)\rho_{\ld{G}_X}$ & $\fr{sp}_{2n-2}[2]\mmod \Sp_{2n-2}$ & Antipodal action on $S^{2n}$ \\
 $N(\G_2)$ & \cref{rmk: G2 mod normalizer of SL3} & $\G_2/\N_{\G_2}(\SL_3)$ & N & $(T^\ast(6) \AA^2(10, 0))/(\Z/2)$ & $10\rho_{\ld{G}_X}$ & $\sl_2[2]\mmod \SL_2$ & Antipodal action on $S^6$ \\
 \hline
\end{tabular}
}
\vspace{1cm}
\caption{Table of dualities for affine homogeneous rank one spherical varieties; only the final three varieties have roots of type N. For each of these varieties, $\ld{G}_X = \SL_2$. The column labeled by $\ld{Y}$ should be thought of encoding the cotangent bundle to the dual Hamiltonian $\ld{G}$-space to $X$. The dual group $\ld{G}_X$ is equipped with the grading via the cocharacter specified above.
The normalization term is of the form $\ld{\fr{l}}^\ast[2]\mmod \ld{L}$,  where $\ld{\fr{l}}$ is the Lie algebra of the subgroup $L^\wedge_X$ of $\ld{G}$ from \cite{knop-schalke} (which is \textit{not} the Langlands dual of the Levi $L(X)$!); in particular, see \cite[Final column of Table 3]{knop-schalke} for the groups $L^\wedge_X$ in the present case. 
\newline
\newline
We have also included the examples of type N (in the last three rows) to illustrate the phenomenon that their dual is really a \textit{stack}, as opposed to a scheme; the action of $\Z/2$ on $T^\ast(\AA^2)$ is via the symplectic form. Notice that this phenomenon only shows up with \textit{even-dimensional} spheres; topologically, this is just because the antipodal map on $S^n$ is not homotopic to the identity precisely in even dimensions.
\newline
\newline
The table also includes some topological phenomena corresponding to each of the rank one affine homogeneous spherical varieties. Most of these are standard, except perhaps for the ``exceptional Hopf fibration'' for the octonionic projective plane, which is proved as \cite[Theorem 1.2]{davis-mahowald}.
The point is that these homotopy-theoretic aspects of $\Omega X$ control essentially all of the properties of the Langlands dual of $X$. For instance, if $S^{i} \to \Omega X \to \Omega S^{j}$ is a fibration as in the table with both $i$ and $j$ being odd, the graded scheme $\ld{Y}$ can be identified with $T^\ast(j-i) \AA^2(j-1,0)$; and, if $\Omega X = \Omega S^{2j+1}$, the graded scheme $\ld{Y}$ can be identified with $\sl_2(2j-2j\rho_{\ld{G}_X})$.
\newline
\newline
The reader should compare the numbers in this table to the points of evaluation of the L-functions appearing in the rightmost column of \cite[Table 1]{sakellaridis-rank-1}. Namely, the dual stack $T^\ast(2j)\AA^2(2i,0)$ in our table corresponds to $L(\std, \tfrac{i}{2}) L(\std, \tfrac{2j-i}{2})$ in \cite[Table 1]{sakellaridis-rank-1}, and similarly $\sl_2(2j-2j\rho_{\SL_2})$ in our table corresponds to $L(\mathrm{ad}, j)$ in \cite[Table 1]{sakellaridis-rank-1}. Note, in particular, that these numbers can be read purely off of the (rational) homotopy type of $G/H$!
}
\label{table: topology and dualities for rank 1 spherical varieties}
\end{table}
\end{landscape}

%% file: each-type/case-by-case-summary.tex
\subsection{Summary}

The main results of this section can be summarized below in \cref{table: topology and dualities for rank 1 spherical varieties}, refining \cref{intro table: topology and dualities for rank 1 spherical varieties}.
This table should be read as follows: there is an equivalence of $\QQ$-linear $\infty$-categories
\begin{equation}\label{eq: mold rk 1 equivalences}
    \Shv^{c,\Sat}_{G\pw{t}}(G\ls{t}/H\ls{t}; \QQ) \simeq \Perf(\sh^{1/2}(\ld{Y}/\ld{G}_X) \times \text{``Normalization''}),
\end{equation}
where $\ld{Y}$ is a $\ld{G}_X$-space, and the normalization term can be identified with $\fr{l}_X^\wedge\mmod L^\wedge_X$, where $L^\wedge_X$ is the subgroup of $\ld{G}$ from \cite{knop-schalke}. This is \textit{not} the Lie algebra of the dual Levi $\ld{L}(X)$.
Here, $\ld{G}_X = \SL_2$ for each of the varieties in \cref{table: topology and dualities for rank 1 spherical varieties}. 
\begin{warning}
    \textit{Throughout this section, we will assume \cref{hypothesis: rank 1 weakly placid} holds for these spherical varieties.} As mentioned in \cref{rmk: meaning of the hypothesis}, it is possible that this hypothesis simply fails; so it should perhaps be stated at the outset that our actual goal in this section is to explicitly verify the conditions of \cref{thm: ordinary homology criterion satake}, and that \cref{hypothesis: rank 1 weakly placid} \textit{only} comes in when using \cref{thm: ordinary homology criterion satake} to prove an equivalence of $\infty$-categories.
\end{warning}
\begin{remark}
    Below, we will only prove the bare equivalence \cref{eq: mold rk 1 equivalences}. We have not proved compatibility with the equivalence of \cref{thm: derived satake}, but we expect it to be possible using the second part of \cref{thm: ordinary homology criterion satake}. One can check the desired compatibility in type $A_n$. In fact, \cref{conj: spherical levi} (essentially) implies \cref{cor: bzsv for CPn}.
\end{remark}

During the course of proving the equivalence \cref{eq: mold rk 1 equivalences} in types $A_n$, $D_2$, and $\G_2$, we will in fact prove an analogue of the above equivalence for $\Shv^{c,\Sat}_{G\pw{t}}(G\ls{t}/H\ls{t}; \ku)^\faux$; the conceptual intepretation of the results thus obtained will be discussed in a later section. The reader only interested in \cref{thm: rk 1 bzsv is true} should simply set $\beta = 0$ everywhere.
\begin{remark}
    Whenever we work with integral (i.e., non-rational) coefficients, we will always invert the order of the Weyl group of $G$.
\end{remark}

The proof of the equivalence \cref{eq: mold rk 1 equivalences} for the spherical varieties of \cref{table: topology and dualities for rank 1 spherical varieties} relies on the criterion of \cref{thm: ordinary homology criterion satake} (or, in the $\ku$-theoretic case, \cref{prop: ku rk 1 reg centr --> thm}). Namely, in each case, we will:
\begin{itemize}
    \item compute the $H$-equivariant homology $\H^H_\ast(\Omega (G/H); \QQ)$ of $\Omega(G/H)$; and
    \item show that $\spec \H^H_\ast(\Omega (G/H); \QQ)$ can be identified with the group scheme $\ld{J}_X'$ from \cref{thm: ordinary homology criterion satake}(a).
\end{itemize}
The key input into the first part is the calculation of the $H$-equivariant homology of the based loop space $\Omega S^V$ of the representation spheres (i.e., one-point compactifications) of (unitary) $H$-representations $V$. This is accomplished in \cref{cor: G-equiv homology of Loops SV}. The underlying analytic spaces of each of the rank one spherical varieties in \cref{table: topology and dualities for rank 1 spherical varieties} are either representation spheres themselves, or can be built as a quotient of a representation sphere, so \cref{cor: G-equiv homology of Loops SV} lets us describe the $H$-equivariant ($\ku$-)homology of $\Omega(G/H)$.
Once $\H^H_\ast(\Omega (G/H); \QQ)$ has been computed, the second part is rather straightforward. Indeed, the group scheme $\ld{J}_X = \spec \H^H_\ast(\Omega (G/H); \QQ)$ is quite simple in each example, and the main difficulty in checking \cref{thm: ordinary homology criterion satake}(a) is in bookkeeping weights.
\begin{remark}
    In applying \cref{thm: ordinary homology criterion satake}, we will see that the only piece of $\spec \H_H^\ast(\ast; \QQ) \cong \ld{\fr{t}}_H^\ast(2)\mmod W_H \cong \ld{\fr{h}}^\ast(2)\mmod \ld{H}$ which ``interacts'' with $\ld{Y}/\ld{G}_X$ via the Kostant section of \cref{thm: ordinary homology criterion satake} is the coordinate corresponding to the highest degree fundamental invariant of $W_H$ acting on $\ld{\fr{t}}_H^\ast(2)$. The remainder of $\ld{\fr{t}}_H^\ast(2)\mmod W_H$ makes up the normalization term $\fr{l}_X^\wedge\mmod L^\wedge_X$.
\end{remark}

Since this section is somewhat technical, let us make an observation about the main qualitative difference between root types T and G which appears in the course of the proof. Let us focus on the prototypical cases of types $A_1$ (which is of root type T) and $D_2$ (which is of root type G). In these cases, there is no normalization term; in general, this normalization term comes from cohomology classes in $\H^\ast(BG; \QQ)$ which do not ``interact'' with $\H_\ast(\Omega(G/H); \QQ)$.
\begin{enumerate}
    \item In type $A_1$, there are isomorphisms of derived $\QQ$-schemes
    \begin{multline*}
        \spec \sh^{1/2} \H^{\SO_2}_\ast(\Omega (\SO_3/\SO_2); \QQ) \cong \spec \QQ[b, x]/bx 
        \cong \AA^1[2] \times_{T^\ast[2] \AA^2[2,0]/\SL_2[-2\rho_{\SL_2}]} \AA^1[2],
    \end{multline*}
    where $b$ lives in degree $2$, and $x$ lives in degree $-2$. Let $\ld{J}_X$ denote the above group scheme. Then $\ld{J}_X$ is \textit{not} flat over $\AA^1[2] = \spec \QQ[x]$.

    The nonflatness of $\ld{J}_X$ over $\ld{\fr{h}}^\ast[2] \mmod \ld{H}$ is characteristic of the case of roots of type T. Topologically, this corresponds to the observation that the $H$-invariant subspace $(G/H)^H$ is just $S^0$. (In general, one can use Atiyah-Bott localization and \cref{thm: ordinary homology criterion satake} to see that the phenomenon of $\Omega (G/H)^H$ being rationally contractible is Langlands dual to the triviality of the stabilizer in $\ld{G}$ of a generic point in $\ld{M}$.)
    \item In type $D_2$, there is an isomorphism $\SO_4/\mu_2 \cong \SO_3 \times \SO_3$, so that $\SO_4/\mu_2 \SO_3 \simeq \SO_3$. Then, there are isomorphisms of derived $\QQ$-schemes
    \begin{multline*}
        \spec \sh^{1/2} \H^{\SO_3}_\ast(\Omega (\SO_4/\mu_2 \SO_3); \QQ) \cong \spec \sh^{1/2} \H^{\SO_3}_\ast(\Omega \SO_3; \QQ)\\
        \cong \spec \QQ[x, a^{\pm 1}, \tfrac{a - a^{-1}}{2x}]^{\Z/2}
        \cong \AA^1[2]\mmod (\Z/2) \times_{\sl_2[2-2\rho_{\SL_2}]/\SL_2[-2\rho_{\SL_2}]} \AA^1[2] \mmod (\Z/2),
    \end{multline*}
    where the map $\AA^1[2]\mmod (\Z/2) \to \sl_2[2-2\rho_{\SL_2}]/\SL_2[-2\rho_{\SL_2}]$ is the Kostant slice sending $x^2 \mapsto \begin{psmallmatrix}
        0 & 1\\
        x^2 & 0
    \end{psmallmatrix}$.
    Let $\ld{J}_X$ denote the above group scheme. Then $\ld{J}_X$ is flat over $\AA^1[2]\mmod (\Z/2) = \spec \QQ[x^2]$, and can be identified with the group scheme of regular centralizers of $\ld{G}_X = \SL_2$.

    The identification of $\ld{J}_X$ over $\ld{\fr{h}}^\ast[2] \mmod \ld{H}$ with (the product of the normalization factor with) the group scheme of regular centralizers of $\SL_2$ is characteristic of the case of roots of type G. Topologically, this corresponds to the observation that the $T$-invariant subspace $(G/H)^T$ is just $S^1$ (where $T$ is a maximal torus of $H$). 
\end{enumerate}

\begin{remark}\label{rmk: similarities between spheres}
    As the reader will see, the calculations of this section are quite repetitive; it is possible to handle all the type T and type G cases simultaneously, but at the risk of confusing oneself with various gradings (in other words, the repetitive nature of this section is mostly for my purposes, and one can conglomerate these calculations into a more uniform argument if desired).
    However, the similarity between the examples in \cref{table: topology and dualities for rank 1 spherical varieties} seems to be related to much deeper phenomena in relative Langlands duality; we hope to explore this in future work. Let us briefly mention the sort of phenomena we hope to explain.
    
    If $X$ is an affine (spherical) $G$-variety, it is often the case that there is some subgroup ${G}' \subseteq {G}$ and a $G'$-variety $X'$ such that \textit{up to grading}, the dual stack $\ld{M}/\ld{G}$ associated to the $G$-variety $X$ via \cref{conj: bzsv} can be identified with the dual stack $\ld{M}'/\hat{G}'$ associated to $X'$. In particular, the dual groups of $X$ and $X'$ are isomorphic.
    
    For instance, suppose $X = \SO_{2n+1}/\SO_{2n}$. By \cref{thm: rk 1 bzsv is true}, the underlying \textit{ungraded} dual stack associated to the $\SO_{2n+1}$-variety $X$ is $\ld{M}/\Sp_{2n} \cong T^\ast(\AA^2)/\SL_2 \times \AA^{n-1}$. Let $G' = \SO_3 \times \SO_{2n-2} \subseteq {G} = \SO_{2n+1}$. However, by \cref{thm: rk 1 bzsv is true}, this is \textit{also} the underlying ungraded dual stack associated to the ${G}'$-variety ${G}'/(\SO_2 \times \SO_{2n-2})$. The underlying homotopy type of $X$ is the sphere $S^{2n}$, and the underlying homotopy type of ${G}'/(\SO_2 \times \SO_{2n-2})$ is $S^2$; there is no direct relation between these homotopy types, but nevertheless, the dual stacks are isomorphic up to grading. (In this sense, the grading plays an absolutely vital role, and getting it right is one of the reasons for the repetitive/tedious nature of this section.) Spencer Leslie has told us that he is working on describing the relationship between such ${G}'$ and the ``endoscopy of $X$''; see \cite{leslie-endoscopic-symmetric-varieties}.
\end{remark}

%% file: each-type/homology-of-loops-sphere.tex
\subsection{Homology of loops on a sphere}\label{subsec: ku Loops rep sphere}

The proof of \cref{thm: rk 1 bzsv is true} rests on a key topological calculation, namely that of the equivariant homology $\ku^T_\ast(\Omega S^V)$ for $S^V$ being the one-point compactification of a (unitary) representation $V$ of a (compact) torus $T$. To illustrate this calculation, let us begin with two simple (but exemplifying) examples.
\begin{notation}
    Recall the group scheme $\GG_\beta$ from \cref{cstr: kuS1 and Gbeta}, i.e., the graded group scheme over $\Z'[\beta]$ given by $\spec \Z'[\beta, t^{\pm 1}, \tfrac{t-1}{\beta}]$ where $\beta$ lives in weight $2$, and with coproduct determined by the formula $t\mapsto t \otimes t$. The invertible class $t$ defines a homomorphism $\GG_\beta \to (\GG_m)_{\Z'[\beta]}$ which exhibits $\GG_\beta$ as an affine blowup of $(\GG_m)_{\Z'[\beta]}$ at the identity of $\GG_m$. The kernel of this map is given by a group scheme $\GG_\beta^0$, whose underlying graded $\Z'[\beta]$-scheme is isomorphic to $\AA^1(-2) = \spec \Z'[\beta, \tfrac{t-1}{\beta}]$, and whose group structure is given by
    $$\tfrac{t-1}{\beta}\mapsto \tfrac{t-1}{\beta} \otimes 1 + 1 \otimes \tfrac{t-1}{\beta} + \beta \tfrac{t-1}{\beta} \otimes \tfrac{t-1}{\beta}.$$
    A compact torus $T$ defines a group scheme $T_\beta$ given by $\Hom(\bX^\ast(T), \GG_\beta)$, and hence a subgroup scheme $T_\beta^0 = \ker(T_\beta \to T)$. Let $\lambda: T \to S^1$ be a character, let $T_\lambda$ denote its kernel, and let $n\geq 0$. Define $T_{\lambda, \beta}^{[n]}$ to be the closed subscheme of $T_\beta$ given by the $n$th infinitesimal neighborhood of $T_{\lambda, \beta}$; similarly for $T_{\lambda, \beta}^{0, [n]}$. Let $x_\lambda \in \co_{T_\beta}$ denote the function which cuts out $T_{\lambda, \beta}$.

    Since $\GG_\beta = \spec \Z'[\beta, t^{\pm 1}, \tfrac{t-1}{\beta}]$, if we call $x = \tfrac{t-1}{\beta}$ (so that $x$ lives in weight $-2$), this ring can be identified with $\spec \Z'[\beta, x, \tfrac{1}{1+\beta x}]$. It will often be more convenient to consider this presentation of $\GG_\beta$.
    For $n\in \Z$, we will write $[n](x)$ to denote the $n$-fold sum of $x$ in the group structure on $\GG_\beta$, so that $[n](x) = \frac{(1+\beta x)^n - 1}{\beta}$.
\end{notation}

\begin{example}\label{ex: borel homology loops S2}
    Let $\std$ denote the standard $1$-dimensional complex representation of $S^1$, and consider the one-point compactification $S^\std$ (so that its underlying nonequivariant space is $S^2$). We will be interested in computing the (for the moment) Borel-equivariant homology $\H^{S^1}_\ast(\Omega S^\std; \Z) = \pi_\ast \Z[\Omega S^\std]^{hS^1}$. There is a homotopy fixed points spectral sequence
    $$E_2^{\ast,\ast} \cong \H_\ast(\Omega S^\std; \Z) \otimes_\Z \H^\ast(\CP^\infty; \Z) \Rightarrow \pi_\ast \Z[\Omega S^\std]^{hS^1}.$$
    The $E_2$-page can be computed easily to be 
    $$E_2^{\ast,\ast} \cong \Z[a,b]\pw{x}/a^2,$$
    where $x \in \H^2(\CP^\infty; \Z)$ is the first Chern class, $a\in \H_1(\Omega S^\std; \Z)$ coming from the inclusion $S^1 \subseteq \Omega S^2$, and $b\in \H_2(\Omega S^\std; \Z)$ coming from the map $S^2 \to \Omega S^2$ adjoint to the Hopf fibration $S^3 \to S^2$. There is a single differential $d_2(a) = bx$ (if $\std$ is replaced by $\std^{\otimes_\cc n}$ for some $n\geq 1$, this differential is simply replaced by $nbx$). After running this differential, the spectral sequence is concentrated in even degrees, and we find that
    $$\pi_\ast \Z[\Omega S^\std]^{hS^1} \cong \Z[b]\pw{x}/bx.$$
    Compare to \cref{ex: hochschild and homology of loops Sn}(b) with $n=1$.
    Exactly the same calculation holds with $\Z$ replaced by $\ku$:
    $$\pi_\ast \ku[\Omega S^{\std^{\otimes_\cc n}}]^{hS^1} \cong \Z[\beta, b]\pw{x}/b[n](x).$$
\end{example}


\begin{remark}
    Note that the quotient $\pi_\ast \Z[\Omega S^\std]^{hS^1}/x$ is precisely the nonequivariant homology $\H_\ast(\Omega S^\std; \Z)$. Indeed, \cref{ex: borel homology loops S2} says that the class $b$ is $x$-torsion in $\pi_\ast \Z[\Omega S^\std]^{hS^1}$; therefore, if we kill $x$, the class $bx$ in degree zero bumps up to a class $\sigma(bx)$ in degree $1$. 

    One interesting observation is that the homotopy {quotient} $\pi_\ast \Z[\Omega S^\std]_{hS^1}$ (which would compute what is traditionally called equivariant homology) is \textit{not} concentrated in even degrees: namely, the fact that $b$ is $x$-torsion implies that the $\Z\pw{x}$-linear dual of $\pi_\ast \Z[\Omega S^\std]^{hS^1}$ will have odd homotopy groups.
    
    Let us mention that the fact that $\pi_\ast \Z[\Omega S^\std]^{hS^1}$ is concentrated in even degrees is an absolutely crucial fact (related to the subtleties of \cref{lem: no e2 half shearing}), which has important implications in the Langlands duality of \cref{thm: rk 1 bzsv is true}, and emphasizes the role of equivariance in our discussion. 
\end{remark}
\begin{example}\label{ex: borel homology loops S3}
    Again, let $\std$ denote the standard $1$-dimensional complex representation of $S^1$, and consider the one-point compactification $S^{\std\oplus \RR}$ (so that its underlying nonequivariant space is $S^3$). We will be interested in computing the (for the moment) Borel-equivariant homology $\H^{S^1}_\ast(\Omega S^{\std \oplus \RR}; \Z) = \pi_\ast \Z[\Omega S^{\std \oplus \RR}]^{hS^1}$. There is a homotopy fixed points spectral sequence
    $$E_2^{\ast,\ast} \cong \H_\ast(\Omega S^{\std \oplus \RR}; \Z) \otimes_\Z \H^\ast(\CP^\infty; \Z) \Rightarrow \pi_\ast \Z[\Omega S^{\std \oplus \RR}]^{hS^1}.$$
    Now, the $E_2$-page is simply $\Z[b]\pw{x}$, where again $x \in \H^2(\CP^\infty; \Z)$ is the first Chern class, and $b\in \H_2(\Omega S^3; \Z)$ is the generator. The entire spectral sequence is concentrated in even degrees, so there can be no differentials, and the spectral sequence degenerates. This implies that
    $$\pi_\ast \Z[\Omega S^{\std \oplus \RR}]^{hS^1} \cong \Z[b]\pw{x}.$$
    Compare to \cref{ex: hochschild and homology of loops Sn}(a) with $n=1$.
\end{example}
\begin{remark}
    Note that in \cref{ex: borel homology loops S3}, the class $bx$ is topologically nilpotent, so that $1 + bx$ is invertible. The inclusion $S^\RR \subseteq S^{\std \oplus \RR}$ induces a map $\Z \cong \Omega S^\RR \to \Omega S^{\std \oplus \RR}$; this map is in fact just the inclusion of the $S^1$-fixed points of $\Omega S^{\std \oplus \RR}$. Therefore, there is a map $\Z[\Omega S^\RR] \cong \Z[a^{\pm 1}] \to \Z[\Omega S^{\std \oplus \RR}]^{hS^1}$, and it is not difficult to check that this map sends $a \mapsto 1 + bx$.
\end{remark}

In order to describe the main calculation, we need to introduce some notation.
\begin{setup}\label{setup: unitary rep of torus}
    Let $T$ be a \textit{compact} torus, and let $V$ be a (unitary) representation of $T$ with no nonzero fixed vectors. Let $\Lambda(V)$ denote the set of weights of $V$, and let $\chi_V: T \to S^1$ to denote the character of $V$. Moreover, if $\lambda: T \to S^1$ is a character, we will write $T_\lambda$ to denote the kernel of $\lambda$. Note that $\pdb{\chi_V, \lambda}$ is the dimension of the $\lambda$-weight space of $V$ (as a complex vector space).
\end{setup}
\begin{setup}\label{setup: inverting 2}
    Throughout, when we talk about coefficients in a homology theory, we will \textit{invert the prime $2$}. This will be implicit in the notation, and we will write $\Z'$ to denote $\Z[1/2]$.
\end{setup}

\begin{definition}\label{def: coord axes}
    In \cref{setup: unitary rep of torus}, fix an integer $j \in \Z$, and define $\cC_V$ to be the graded $\Z'[\beta]$-scheme given by the union
    $$\cC_V(-j) = (T_\beta^0 \times \{0\}) \cup \bigcap_{\lambda \in \Lambda(V)} (T_{\lambda, \beta}^{0, [\pdb{\chi_V, \lambda}]} \times \AA^1(-j)).$$
    We will call $\cC_V(-j)$ the \textit{$V$-coordinate axes with weight $j$}.
\end{definition}
\begin{definition}\label{def: affine blowup}
    In \cref{setup: unitary rep of torus}, consider the blowup
    $$X := \mathrm{Bl}\left( \bigcap_{\lambda \in \Lambda(V)} (T_{\lambda, \beta}^{[\pdb{\chi_V, \lambda}]} \times \{1\}) \subseteq T_\beta \times \GG_m \right).$$
    Let $\Bl_V$ denote the complement of the proper preimage of $\{0_{T_\beta}\} \times \GG_m$ from $X$, so that $\Bl_V$ is an affine blowup of $T_\beta \times \GG_m$.
\end{definition}
\begin{example}
    Suppose $T = S^1$, and let $V$ denote the weight $n$ representation of $T$. Then there is an isomorphism
    $$\cC_V(-j) \cong \spec \Z'[\beta, x,b]/bx,$$
    where $x$ lives in weight $-2$ and $b$ lives in weight $j$. This is the reason for the terminology of \cref{def: coord axes}: we are more concerned with the weight of the function $b$ on $\cC_V(-j)$.
    Similarly, there is an isomorphism
    $$\Bl_V \cong \spec \Z'[\beta, t^{\pm 1}, \tfrac{t-1}{\beta}, a^{\pm 1}, \tfrac{(a-1)\beta}{t-1}] \cong \spec \Z'[\beta, x, \tfrac{1}{1+\beta x}, a^{\pm 1}, \tfrac{a-1}{x}],$$
    where $t$ lives in weight $0$ (so $\frac{t-1}{\beta}$ lives in weight $2$), and $a$ lives in weight $0$ (so $\tfrac{(a-1)\beta}{t-1}$ lives in weight $2$).
\end{example}
\begin{example}\label{ex: CV and BlV examples}
    More generally, let $T = (S^1)^m$, and let $V$ denote the representation $\bigoplus_{i=1}^m d_i \std_i^{\otimes_\cc n_i}$. Then 
    $$\cC_V(-j) \cong \spec \Z'[\beta, x_1, \cdots, x_m, \prod_{i=1}^m \tfrac{1}{1+\beta x_i}, b]/b\prod_{i=1}^m [n_i](x_i)^{d_i},$$
    with each $x_i$ in weight $-1$ and $b$ in weight $j$.
    Similarly, we have 
    $$\Bl_V \cong \spec \Z'[\beta, x_1, \cdots, x_m, \prod_{i=1}^m \tfrac{1}{1+\beta x_i}, a^{\pm 1}, \prod_{i=1}^m \tfrac{a-1}{[n_j](x_j)^{d_j}}],$$
    where each $x_i$ lives in weight $-1$, and $a$ lives in weight $0$ (so $\prod_{i=1}^m \tfrac{a-1}{[n_j](x_j)^{d_j}}$ lives in weight $2\sum_{i=1}^m d_j = \dim_\RR(V)$). This immediately implies:
\end{example}
\begin{lemma}\label{lem: relation BV and CV}
    In \cref{setup: unitary rep of torus}, the fiber of the projection map $\cB_V \to \GG_m$ over $\{1\}\in \GG_m$ is isomorphic to $\cC_V(-\dim_\RR(V))$.
\end{lemma}
\begin{theorem}\label{thm: homology of loops SV}
    In \cref{setup: unitary rep of torus}, let $S^V$ denote the one-point compactification of $V$.
    Then there are graded isomorphisms of $\pi_\ast \ku_T$-algebras
    \begin{align*}
        \ku^T_\ast(\Omega S^V) & \cong \co_{\cC_V(2-2\dim_\RR(V))}, \\
        \ku^T_\ast(\Omega S^{V\oplus \RR}) & \cong \co_{\Bl_V}.
    \end{align*}
    In particular, both $\ku^T_\ast(\Omega S^{V})$ and $\ku^T_\ast(\Omega S^{V\oplus \RR})$ are concentrated in even weights, and are graded commutative $\pi_\ast \ku_T$-algebras.
\end{theorem}
\begin{remark}\label{rmk: happens to be commutative}
    If $T$ acts on a pointed space $X$ (and we are given some multiplicative presentation of $\Omega X$ as a $T$-space), the equivariant homology $\ku^T_\ast(\Omega X)$ need not be a commutative $\Z'[\beta]$-algebra in general: \textit{a priori}, it is only an associative $\Z'[\beta]$-algebra, since $\Omega X$ generally only admits the structure of an $\E{1}$-space. Although $\Omega S^V$ is still generally only an $\E{1}$-space (unless $V$ is isomorphic to $\emptyset$, $\RR$, or $\RR^3$), \cref{thm: homology of loops SV} implies that $\ku^T_\ast(\Omega S^{V})$ and $\ku^T_\ast(\Omega S^{V\oplus \RR})$ are in fact concentrated in even weights, and generated (as a $\pi_\ast \ku_T$-algebra) by a single class. In particular, it is necessarily a commutative $\Z'[\beta]$-algebra\footnote{This is a manifestation of the fact that in classical algebra, an associative algebra being commutative is a \textit{property}, whereas in homotopy theory, it is more \textit{structure}.}. We will implicitly use this observation throughout this article, by rewriting \cref{thm: homology of loops SV} as a pair of graded isomorphisms of $\Z'[\beta]$-schemes
    \begin{align*}
        \spec \ku^T_\ast(\Omega S^V) & \cong {\cC_V(2-2\dim_\RR(V))}, \\
        \spec \ku^T_\ast(\Omega S^{V\oplus \RR}) & \cong {\Bl_V}.
    \end{align*}
\end{remark}
\begin{proof}[Proof of \cref{thm: homology of loops SV}]
    Let us first compute the \textit{Borel} $T$-equivariant $\ku$-homology of $\Omega S^V$. Write $V = \cc^n$; then there is a homotopy fixed points spectral sequence
    $$E_2^{\ast,\ast} \cong \pi_\ast \ku^{hT} \otimes_{\Z'[\beta]} \ku_\ast(\Omega S^{\cc^n}) \Rightarrow \pi_\ast \ku[\Omega S^V]^{hT}.$$
    The EHP sequence for $\Omega S^V$ is the fibration given by
    $$S^{2n-1} \to \Omega S^{2n} \to \Omega S^{4n-1}.$$
    Note that the map $S^{2n-1} \to \Omega S^{2n} \simeq \Omega \Sigma S^{2n-1}$ is \textit{not} defined equivariantly: the unit sphere $S(V)$ does not have $T$-fixed points (by assumption on $V$), so there is no $T$-equivariant basepoint of $S(V)$ with respect to which the reduced suspension can be constructed.
    In any case, the EHP sequence splits after inverting the prime $2$ (this is the reason for \cref{setup: inverting 2}), which implies that there is an isomorphism $$\ku_\ast(\Omega S^{2n}) \cong \Z'[\beta, a, b]/a^2,$$
    where $|a| = 2n-1$ and $|b| = 4n-2$. (See, e.g., \cref{ex: hochschild and homology of loops Sn} for a massively overblown derivation of this isomorphism.)
    \begin{remark}
        That the EHP sequence splits after inverting the prime $2$ goes all the way back to Serre: on \cite[p. 281]{serre-loop-splitting}, he showed that there is a equivalence
	$$\can \times \Omega [\iota_{2n}, \iota_{2n}]: S^{2n-1} \times \Omega
	S^{4n-1} \to \Omega S^{2n}$$
        after inverting $2$. Here, $[\iota_{2n}, \iota_{2n}] \in \pi_{4n-1}(S^{2n})$ denotes the Whitehead product of $\iota_{2n} \in \pi_{2n}(S^{2n})$ with itself.
    \end{remark}
    
    It follows that the $E_2$-page can be identified with
    $$E_2^{\ast,\ast} \cong \co_{T_\beta}^\wedge[a, b]/a^2.$$
    Each weight $\lambda \in \Lambda(V)$ defines a function $x_\lambda\in \co_{T_\beta}$, and there is a single differential in this spectral sequence, given by
    \begin{equation}\label{eq: diffl in Loops SV sseq}
        d_{\dim_\RR(V)}(a) = b\prod_{\lambda \in \Lambda(V)} x_\lambda.
    \end{equation}
    One can see this by reducing to the case when $V$ is one-dimensional, in which case \cref{eq: diffl in Loops SV sseq} follows from \cref{ex: borel homology loops S2}. 
    The $E_{\dim_\RR(V)+1}$-page of the spectral sequence is then concentrated entirely in even degrees, and therefore degenerates.
    This implies that 
    \begin{equation}\label{eq: borel equiv homology Loops SV}
        \pi_\ast \ku[\Omega S^V]^{hT} \cong \co_{T_\beta}^\wedge[b]/b\prod_{\lambda \in \Lambda(V)} x_\lambda.
    \end{equation}
    Note that this is nearly $\co_{\cC_V(2-2\dim_\RR(V))}$, except for the completion.
    
    The above calculation of $\pi_\ast \ku[\Omega S^V]^{hT}$ is enough to imply the desired calculation of $\ku^T_\ast(\Omega S^V)$. Indeed, let $e_j$ be a given basis vector of $\bX^\ast(T)$, let $T_j$ denote the kernel of $e_j: T \to S^1$, and let $V_j$ denote the fixed locus $V^{T_j}$.
    Let $\cU_j$ denote the complement of the union of the closed subschemes $T'_\beta$ ranging over all closed subgroups $T' \subseteq T$ which do not contain $T_j$. Then \cref{lem: atiyah localization} gives an isomorphism 
    \begin{equation}\label{eq: localization for Loops SV}
        \ku^T_\ast(\Omega S^V) |_{\cU_j} \simeq \ku^T_\ast(\Omega S^{V_j})|_{\cU_j}.
    \end{equation}
    Indeed, since $(\Omega S^V)^{T_j} \simeq \Omega (S^V)^{T_j}$, it suffices to note that $(S^V)^{T_j} \simeq S^{V_j}$.
    Using the fracture square and induction on the dimension of $V$, one finds that there is a Cartesian square
    $$\xymatrix{
    \ku^T_\ast(\Omega S^V) \ar[r] \ar[d] & \pi_\ast \ku[\Omega S^V]^{hT} \ar[d] \\
    \pi_\ast \ku_T |_{T_\beta - \bigcup_{1\leq j \leq m} T_{j,\beta}} \ar[r] & \pi_\ast \ku^{tT}.
    }$$
    This precisely has the effect of correcting the completion in \cref{eq: borel equiv homology Loops SV}, which recovers $\co_{\cC_V(2-2\dim_\RR(V))}$.

    Turning to $\Omega S^{V\oplus \RR}$, let us first compute the \textit{Borel} $T$-equivariant $\ku$-homology of $\Omega S^{V\oplus \RR}$. Write $V = \cc^n$; then there is a homotopy fixed points spectral sequence
    $$E_2^{\ast,\ast} \cong \pi_\ast \ku^{hT} \otimes_{\Z'[\beta]} \ku_\ast(\Omega S^{\cc^n \oplus \RR}) \Rightarrow \pi_\ast \ku[\Omega S^{V\oplus \RR}]^{hT}.$$
    There is an isomorphism $\ku_\ast(\Omega S^{\cc^n \oplus \RR}) \cong \ku_\ast[b]$, where $|b| = 2n$. This implies that the entire spectral sequence is concentrated in even degrees, so there are no differentials, and we find that
    \begin{equation}\label{eq: borel equiv homology Loops SV plus R}
        \pi_\ast \ku[\Omega S^{V\oplus \RR}]^{hT} \cong \Z'[\beta, b]\pw{x_1, \cdots, x_n}.
    \end{equation}
    To compute $\ku^T_\ast(\Omega S^{V\oplus \RR})$, we will again use the fracture square. Again, it is not difficult to reduce to the case when $T = S^1$ and $V = {\std^{\otimes_\cc n}}$. Then there is a Cartesian square
    $$\xymatrix{
    \ku^{S^1}_\ast(\Omega S^{\std^{\otimes_\cc n} \oplus \RR}) \ar[r] \ar[d] & \ku^{S^1}_\ast(\Omega S^{\std^{\otimes_\cc n} \oplus \RR})^\wedge_{x} \ar[d] \\
    \ku^{S^1}_\ast(\Omega S^{\std^{\otimes_\cc n} \oplus \RR})[x^{-1}] \ar[r] & \ku^{S^1}_\ast(\Omega S^{\std^{\otimes_\cc n} \oplus \RR})^\wedge_x [x^{-1}].
    }$$
    Note that the $S^1$-fixed point set $(\Omega S^{\std^{\otimes_\cc n} \oplus \RR})^{S^1}$ is simply $\Omega S^1 \cong \Z$, so that \cref{lem: atiyah localization} lets us identify the bottom-left corner with $\pi_\ast \ku_T[a^{\pm 1}]$.
    By \cref{eq: borel equiv homology Loops SV plus R}, we can identify the above Cartesian square with
    $$\xymatrix{
    \ku^{S^1}_\ast(\Omega S^{\std^{\otimes_\cc n} \oplus \RR}) \ar[r] \ar[d] & \Z'[\beta, b]\pw{x} \ar[d] \\
    \Z'[\beta, a^{\pm 1}, x^{\pm 1}, \tfrac{1}{1+\beta x}] \ar[r] & \Z'[\beta, b]\ls{x},
    }$$
    where the bottom map sends $a \mapsto 1 + bx$. It follows that
    $$\ku^{S^1}_\ast(\Omega S^{\std^{\otimes_\cc n} \oplus \RR}) \cong \Z'[\beta, x, \tfrac{1}{1+\beta x}, a^{\pm 1}, \tfrac{a-1}{x}] \cong \co_{\Bl_V},$$
    where $\tfrac{a-1}{x} \mapsto b$.
\end{proof}
\begin{remark}
    The evenness of \cref{thm: homology of loops SV} is always true for $\ku^T_\ast(\Omega S^{V\oplus \RR})$ (in fact, more generally for $\ku^T_\ast(\Omega S^{V\oplus \RR^{2n+1}})$ for any $n\geq 0$), but it is \textit{not} true for $\ku^T_\ast(\Omega S^{V\oplus \RR^{2n}})$ if $n>0$.
\end{remark}
The following result is an immediate consequence of \cref{thm: homology of loops SV} and \cref{prop: Weyl invts}.
\begin{corollary}\label{cor: G-equiv homology of Loops SV}
    Let $G$ be a connected compact Lie group whose $\pi_1$ is torsion-free, and let $T\subseteq G$ be a maximal torus with associated Weyl group $W$. If $V$ is a (unitary) $G$-representation with no nonzero $T$-fixed vectors, there are graded isomorphisms of $\pi_\ast \ku_G$-schemes
    \begin{align*}
        \spec \ku^G_\ast(\Omega S^V) & \cong \cC_V(2-2\dim_\RR(V)) \mmod W, \\
        \spec \ku^G_\ast(\Omega S^{V\oplus \RR}) & \cong \Bl_V \mmod W.
    \end{align*}
\end{corollary}
\begin{remark}
    If $H$ is a compact Lie group with maximal torus $T$ and $V$ is a (unitary) $H$-representation with no nonzero $T$-fixed vectors, then we implicitly view $\ku^H_\ast(\Omega S^V)$ as a commutative algebra over $\pi_\ast \ku_H$ as in \cref{rmk: happens to be commutative}. In the case-by-case analysis below, this is in fact not as abusive as it might seem: namely, if $G$ is a compact Lie group and $H\subseteq G$ is a closed inclusion of subgroups such that $S^V \cong G/H$ as $H$-equivariant spaces, the natural $\E{1}$-algebra structure on the $\ku_H$-linear dual of $\cf_H(\Omega S^V)$ in fact upgrades to an $\E{2}$-algebra structure via \cref{cor: loop homology E2}. In particular, the commutative ring structure on the completion of $\ku^H_\ast(\Omega S^V)$ does indeed have a homotopical origin.
\end{remark}

%% file: each-type/type-An.tex
\subsection{Type $A_n$: $\GL_{n+1}/\GL_n$}

Our goal in this section is to prove \cref{thm: rk 1 bzsv is true} in type $A_n$, i.e., for the spherical $\GL_{n+1}$-variety $\GL_{n+1}/\GL_n$. 
We will write $G = \GL_{n+1}$ and $H = \GL_n$, so $\ld{G} = \GL_{n+1}$. Recall from \cref{table: topology and dualities for rank 1 spherical varieties} that $\ld{G}_X = \GL_2$. Equip $\GL_2$ with the grading where the entries of a matrix $\begin{psmallmatrix}
    a & b\\
    c & d
\end{psmallmatrix}$ have the following weights: $a$ and $d$ live in weight zero, $b$ lives in weight $2n$, and $c$ lives in weight $-2n$; we will write $\GL_2(-2n\rho_{\ld{G}_X})$ to denote this graded group.
\begin{theorem}[\cref{thm: rk 1 bzsv is true} in type $A_n$]\label{thm: bzsv for An}
    There is an equivalence of $\QQ$-linear $\infty$-categories
    $$\Shv_{G\pw{t}}^{c,\Sat}(G\ls{t}/H\ls{t}; \QQ) \simeq \Perf(T^\ast[2n](\AA^2[2n,0])/\GL_2[-2n\rho_{\ld{G}_X}] \times \gl_{n-1}[2]\mmod \GL_{n-1}).$$
\end{theorem}

\begin{example}
    For instance, if $n=1$, we have $\ld{G}_X = \ld{G}$, and so \cref{thm: bzsv for An} states that there is an equivalence of $\QQ$-linear $\infty$-categories
    $$\Shv_{\GL_2\pw{t}}^{c,\Sat}(\GL_2\ls{t}/\GG_m\ls{t}; \QQ) \simeq \Perf(T^\ast[2](\AA^2[2,0])/\GL_2[-2\rho_{\ld{G}_X}]).$$
    The Koszul dual of this statement is also proved as \cite[Theorem 1.8(2)]{braverman-finkelberg}. A variant of this equivalence is a special case of \cref{cor: bzsv for CPn}:
    \begin{equation}\label{eq: PGL2 mod Gm remark}
        \Shv_{\PGL_2\pw{t}}^{c,\Sat}(\PGL_2\ls{t}/\GG_m\ls{t}; \QQ) \simeq \Perf(T^\ast[2](\AA^2[2,0])/\SL_2[-2\rho_{\ld{G}_X}]).
    \end{equation}
    In other words, the Hamiltonian $\SL_2$-variety dual to $\PGL_2/\GG_m$ is $T^\ast(\AA^2)$. Let us quickly check part of \cref{conj: generalized kostant slice} in this case. The poset of $B$-orbit closures in $\PGL_2/\GG_m$ is the same as the poset of $\GG_m$-orbit closures in $\PP^1 \cong \PGL_2/B$; this is simply
    $$\{0\} \subseteq \PP^1 \supseteq \{\infty\}.$$
    Similarly, the poset of $\ld{B}$-orbit closures in $\AA^2$ is given by 
    $$\{(0,0)\} \subseteq \{y=0\} \subseteq \AA^2,$$
    which is indeed in bijection with the above set. (One can also see this poset of $\ld{B}$-orbit closures by computing the set of irreducible components of $T^\ast(\AA^2) \times_{\ld{\fr{b}}^\ast} \{0\}$. This is cut out inside $T^\ast(\AA^2) \cong \AA^4$ by the ideal $(v_1 w_1 - v_2 w_2, v_2 w_1)$, whose primary decomposition is
    $$(v_1 w_1 - v_2 w_2, v_2 w_1) = (v_1 v_2) \cap (w_1, w_2) \cap (w_1^2, v_2 w_1, v_1 w_1 - v_2 w_2, v_2^2);$$
    so the fiber product $T^\ast(\AA^2) \times_{\ld{\fr{b}}^\ast} \{0\}$ has three irreducible components.)

\end{example}
\begin{remark}
    There is a relationship between \cref{thm: bzsv for An} and the mirabolic Satake equivalence of \cite{mirabolic-satake}, which studies the spectral decomposition of the spherical $\GL_n \times \GL_{n-1}$-variety $(\GL_n \times \GL_{n-1})/\GL_{n-1}$. In our language, their main result states that there is an equivalence
    \begin{multline*}
        \Shv^{c,\Sat}_{(\GL_n \times \GL_{n-1})\pw{t}}((\GL_n \times \GL_{n-1})\ls{t}/\GL_{n-1}\ls{t}; \QQ) \\
        \simeq \Perf^\sh(\GL_n\backslash T^\ast(\Hom(\AA^n, \AA^{n-1}))/\GL_{n-1}),
    \end{multline*}
    where $\Perf^\sh$ denotes the $\infty$-category of perfect complexes on a shearing. (We have omitted the precise gradings for brevity.)
    If $n=2$, the above equivalence specializes to
    $$\Shv^{c,\Sat}_{(\GL_2 \times \GL_1)\pw{t}}((\GL_2 \times \GL_1)/\GL_1\ls{t}); \QQ) \simeq \Perf^\sh(\GL_2\backslash T^\ast \AA^2/\GG_m),$$
    and forgetting the $\GG_m$-quotient on the coherent side is equivalent to extending $\GL_2\pw{t} \times \GL_1\pw{t}$-equivariance to $\GL_2\pw{t} \times \GL_1\ls{t}$-equivariance on the left-hand side; this in turn recovers \cref{thm: bzsv for An} for the spherical $\GL_2$-variety $\GL_2/\GG_m$.

    If $n>2$, then using the equivalence 
    $$\Shv^{c,\Sat}_{\GL_{n-1}\pw{t}}(\ast; \QQ) \simeq \Perf^\sh(\gl_{n-1}^\ast/_\psi N_{n-1})$$
    with $N_{n-1}$ being the unipotent radical of a Borel subgroup of $\gl_{n-1}$ and $\psi$ being a nondegenerate character on its Lie algebra, the mirabolic Satake equivalence implies that
    $$\Shv^{c,\Sat}_{\GL_n\pw{t}}(\GL_n\ls{t}/\GL_{n-1}\ls{t}; \QQ) \simeq \Perf^\sh(\GL_n\backslash T^\ast(\Hom(\AA^n, \AA^{n-1}))/_\psi N_{n-1}).$$
    Justin Hilburn informed me that using the results of \cite{nakajima-takayama}, one obtains an isomorphism of stacks
    $$\GL_n\backslash T^\ast(\Hom(\AA^n, \AA^{n-1}))/_\psi N_{n-1} \cong T^\ast(\AA^2)/\GL_2 \times \gl_{n-2}\mmod \GL_{n-2},$$
    which shows that \cref{thm: bzsv for An} (up to the question of matching up gradings) is in fact a consequence of the mirabolic Satake equivalence of \cite{mirabolic-satake}. This also follows from the Cartesian square \cref{eq: theta and whit} applied to the example of $H^\diag = \GL_{n-1}^\diag \subseteq G \times H = \GL_n \times \GL_{n-1}$ (see \cref{rmk: dirichlet and neumann}).
\end{remark}
\begin{remark}\label{rmk: Loops2 and BGa}
    One can use the analogue of \cref{thm: bzsv for An} for sheaves with coefficients in $\Z$ to show that there is an equivalence
    \begin{multline*}
        \Loc(\Omega(\U(n+1)/\U(n)); \Z) \simeq \\
        \Perf(\{0\} \times_{\gl_{n+1}[2]\mmod \GL_{n+1}} (T^\ast[2n](\AA^2[2n,0])/\GL_2[-2n\rho_{\ld{G}_X}] \times \gl_{n-1}[2]\mmod \GL_{n-1}))
    \end{multline*}
    describing (derived) \textit{local systems} on $\Omega(\U(n+1)/\U(n)) \simeq \Omega S^{2n+1}$. Here, the map to $\gl_{n+1}[2]\mmod \GL_{n+1}$ is via the moment map. However, one can compute that there is an isomorphism
    $$(T^\ast(\AA^2)/\GL_2 \times \gl_{n-1}\mmod \GL_{n-1}) \times_{\gl_{n+1}\mmod \GL_{n+1}} \{0\} \cong B\GG_a,$$
    where we have ignored shifts for notational simplicity.
    Putting shifts back in, we conclude that there is a Fourier equivalence
    \begin{equation}\label{eq: loc on Omega S2n+1}
        \Loc(\Omega(\U(n+1)/\U(n)); \Z) \simeq \Perf(B\GG_a[-2n]);
    \end{equation}
    this equivalence is quite simple to see using Koszul duality, but it is satisfying to see the right-hand side fall out of \cref{thm: bzsv for An}. This equivalence sends the skyscraper sheaf at the basepoint of $\Omega(\U(n+1)/\U(n)) \simeq \Omega S^{2n+1}$ to the structure sheaf of $B\GG_a[-2n]$.
    
    Using the equivalence \cref{eq: loc on Omega S2n+1} to compute endomorphisms of the skyscraper sheaf at the basepoint of $\Omega(\U(n+1)/\U(n))$, we find that there is an isomorphism
    $$\pi_\ast \Z[\Omega^2 S^{2n+1}] \cong \pi_\ast \sh^{1/2} \Gamma(B\GG_a(-2n); \co).$$
    More generally, one can show that if $j>n$ is even, there is an equivalence of $\E{n+1}$-$\Z$-algebras
    \begin{equation}
        \Z[\Omega^{n+1} S^{j+1}]
        \simeq \sh^{1/2} \Gamma(B^n \GG_a(-j); \co);
    \end{equation}
    this follows from the fact that $\Z[\Omega^{n+1} S^{j+1}]$ is the free $\E{n+1}$-$\Z$-algebra on a class in degree $j-n$, hence is the shearing $\sh^{1/2} \Free_\E{n+1}(\Sigma^{-n} \Z(j))$. But for any $j$ and $n$, there is an equivalence 
    $$\Free_\E{n+1}(\Sigma^{-n} \Z(j)) \simeq \LSym_\Z(\Z[-n](j)) \simeq \Gamma(B^n \GG_a(-j); \co)$$
    of $\E{n+1}$-$\Z$-algebras, coming from Koszul duality.
    

    Let $j\geq 1$ be an integer, fix a prime $p$, and let $J_{p^j - 1}(S^{2n}) \subseteq \Omega S^{2n+1}$ denote the $(p^j - 1)$st partial James construction. 
    Using the EHP fiber sequence
    $$J_{p^j - 1}(S^{2n}) \to \Omega S^{2n+1} \to \Omega S^{2np^j+1},$$
    one can similarly show that there is a Fourier equivalence
    \begin{equation}\label{eq: loc on Omega of james}
        \Loc(J_{p^j - 1}(S^{2n}); \FF_p) \simeq \Perf(B\alpha_{p^j}[-2n]),
    \end{equation}
    which sends the skyscraper sheaf at the basepoint of $J_{p^j - 1}(S^{2n})$ to the structure sheaf of $B\alpha_{p^j}[-2n]$.
    Computing endomorphisms of this skyscraper sheaf, we find that there is an isomorphism\footnote{One conceptual way to compute the cohomology of $B\GG_a(-2)$ and $B\alpha_{p^j}(-2)$ is as follows (assume $p>2$ for simplicity). Let $W$ denote the $p$-typical Witt ring scheme, let $V: F_\ast W \to W$ denote the Verschiebung, let $W_n$ denote the quotient ring scheme of $p$-typical Witt vectors of length $n$, and let $W_n[F^j]$ denote the kernel of $j$-fold Frobenius on $W_n$ (so that $W_1 = \GG_a$ and $W_1[F^j] = \alpha_{p^j}$). All of these group schemes admit a natural action of $\GG_m$ where the $j$th ghost coordinate lives in weight $2p^j$ (to compute the cohomology of $B\GG_a(-2m)$, say, one simply replaces $2p^j$ by $2mp^j$). Then, there is a graded isomorphism
    $$\H^\ast(BW_n[F^j]; \co) \cong \FF_p[\zeta_n, \cdots, \zeta_{n+j-1}] \otimes_{\FF_p} \FF_p[\tau_0, \cdots, \tau_{j-1}]/(\tau_i^2 | 0\leq i \leq j-1),$$
    where $\zeta_i\in \H^2(BW_n[F^j]; \co)$ and $\tau_i \in \H^1(BW_n[F^j]; \co)$ both live in weight $2p^j$. Note the similarity with the dual Steenrod algebra (this is in fact not a coincidence, but explaining this is outside the scope of the present article)!
    
    To see this, first observe that $\H^\ast(BW; \co) \cong \FF_p[\tau_0, \cdots]/(\tau_i^2 | i \geq 0)$ and that $\H^\ast(B^2 F^n_\ast W; \co) \cong \FF_p[\zeta_n, \cdots]$; this follows, for instance, either from the existence of ghost coordinates, or from the fact that $W$ is an extension of the group schemes $F^j_\ast \GG_a^\sharp$, and that there are isomorphisms 
    $$\H^\ast(BF^j_\ast \GG_a^\sharp; \co) \cong \FF_p[\tau_j]/(\tau_j^2), \ \H^\ast(B^2 F^j_\ast \GG_a^\sharp; \co) \cong \FF_p[\zeta_j].$$
    Next, note that there is an exact sequence
    $$F^n_\ast W[F^j] \xar{V^n} W[F^j] \to W_n[F^j]$$
    of flat group schemes. This induces an exact sequence
    $$BW[F^j] \to BW_n[F^j] \to B^2 F^n_\ast W[F^j]$$
    of commutative group stacks. However, since $W[F^j]$ is the kernel of $F^j: W \to F^j_\ast W$, the above calculation of $\H^\ast(BW; \co)$ and $\H^\ast(B^2 F^n_\ast W; \co)$ implies that $\H^\ast(BW[F^j]; \co) \cong \FF_p[\tau_0, \cdots, \tau_{j-1}]/(\tau_i^2 | 0\leq i \leq j-1)$ and that $\H^\ast(B^2 F^n_\ast W[F^j]; \co) \cong \FF_p[\zeta_n, \cdots, \zeta_{n+j-1}]$. The desired calculation of $\H^\ast(BW_n[F^j]; \co)$ follows from this.
    }
    $$\pi_\ast \FF_p[\Omega J_{p^j - 1}(S^{2n})] \cong \pi_\ast \sh^{1/2} \Gamma(B\alpha_{p^j}(-2n); \co).$$
\end{remark}
The proof of \cref{thm: bzsv for An} will occupy the remainder of this section; in fact, we will prove a $\ku$-theoretic deformation of \cref{thm: bzsv for An} below in \cref{cor: bzsv for ku in An}.
\begin{lemma}\label{lem: PGLn+1/GLn and CPn}
    There is a homotopy equivalence $(\GL_{n+1}/\GL_n)(\cc) \simeq S^{2n+1}$. Furthermore, if $B\subseteq G$ is the subgroup of upper-triangular matrices, the Levi quotient $L(X)$ of the parabolic subgroup stabilizing the open $B$-orbit in $\GL_{n+1}/\GL_n$ is given by $\GL_{n-1} \times \GG_m^{\times 2}$.
\end{lemma}
\begin{lemma}\label{lem: pi* ku Un}
    There is an isomorphism of graded $\pi_\ast \ku$-algebras
    $$\pi_\ast \ku_{\U(n)} \cong \Z'[\beta, c_1, \cdots, c_n, \tfrac{1}{1 + \beta c_1 + \cdots + \beta^n c_n}] \cong \co_{T^n_\beta\mmod \Sigma_n},$$
    where $c_j$ lives in weight $-2j$.
\end{lemma}
\begin{proof}
    Let $T^n$ denote the standard diagonal torus of $\U(n)$, so that
    $$\pi_\ast \ku_{T^n} \cong \Z'[\beta, x_1, \cdots, x_n, \prod_{i=1}^n \tfrac{1}{1+\beta x_i}].$$
    Since the Weyl group of $T^n$ inside $\U(n)$ is the symmetric group $\Sigma_n$, \cref{prop: Weyl invts} says that $\pi_\ast \ku_{\U(n)} \cong (\pi_\ast \ku_{T^n})^{\Sigma_n}$. The action of $\Sigma_n$ on $\pi_\ast \ku_{T^n}$ is simply given by permuting the $x_j$. If we set $c_j$ to denote the $j$th elementary symmetric polynomial in the variables $x_1, \cdots, x_n$, the lemma follows immediately.
\end{proof}
\begin{prop}\label{prop: U(n) homology Loops S2n+1}
    There is an isomorphism of graded $\pi_\ast \ku_{\U(n)}$-algebras
    $$\ku_\ast^{\U(n)}(\Omega S^{2n+1}) \cong \Z'[\beta, c_1, \cdots, c_n, \tfrac{1}{1 + \beta c_1 + \cdots + \beta^n c_n}, a^{\pm 1}, \tfrac{a-1}{c_n}],$$
    where $c_j$ lives in weight $-2j$ and $a$ lives in weight $0$.
\end{prop}
\begin{proof}
    By \cref{cor: G-equiv homology of Loops SV}, $\ku^{\U(n)}_\ast(\Omega S^{2n+1}) \cong \ku^{T^n}_\ast(\Omega S^{2n+1})^{\Sigma_n}$.
    As a $T^n$-representation, $S^{2n+1}$ is the one-point compactification of $V = \RR \oplus \bigoplus_{j=1}^n \std_j$. \cref{thm: homology of loops SV} says that
    $$\ku^{T^n}_\ast(\Omega S^{2n+1}) \cong \Z'[\beta, x_1, \cdots, x_n, \prod_{i=1}^n \tfrac{1}{1+\beta x_i}, a^{\pm 1}, \tfrac{a-1}{x_1 \cdots x_n}],$$
    so since the action of $\Sigma_n$ simply permutes the $x_j$ and leaves $y$ invariant, we see from \cref{lem: pi* ku Un} (and $c_n = x_1 \cdots x_n$) that
    $$\ku^{\U(n)}_\ast(\Omega S^{2n+1}) \cong \Z'[\beta, c_1, \cdots, c_n, \tfrac{1}{1 + \beta c_1 + \cdots + \beta^n c_n}, a^{\pm 1}, \tfrac{a-1}{c_n}],$$
    as desired.
\end{proof}

\begin{corollary}\label{cor: U(n) homology Loops CPn}
    There is an isomorphism of graded $\pi_\ast \ku_{\U(n)}$-algebras
    $$\ku_\ast^{\U(n)}(\Omega \CP^n) \cong \Z'[\beta, c_1, \cdots, c_n, \tfrac{1}{1 + \beta c_1 + \cdots + \beta^n c_n}, b]/bc_n,$$
    where $c_j$ lives in weight $-2j$ and $b$ lives in weight $2n$.
\end{corollary}
\begin{proof}
    Let $V$ denote the $\U(n)$-representation $\cc^n$, so that the $\U(n)$ action on $\CP^n$ is obtained by viewing it as $\CP(V \oplus \RR)$. We then have the generalized $\U(n)$-equivariant Hopf fibration
    $$S^1 \to S^{V\oplus \RR} \to \CP(V \oplus \RR),$$
    which induces a $\U(n)$-equivariant fibration of $\E{1}$-spaces
    $$\Omega S^1 \cong \Z \to \Omega S^{V\oplus \RR} \to \Omega \CP(V \oplus \RR).$$
    This implies that there is an equivalence of $\E{1}$-$\ku_{\U(n)}$-algebras
    $$\cf_{\U(n)}(\Omega \CP^n)^\vee \simeq \cf_{\U(n)}(\Omega S^{2n+1})^\vee \otimes_{\ku[\Omega S^1]} \ku.$$
    By \cref{prop: U(n) homology Loops S2n+1}, $\ku_\ast^{\U(n)}(\Omega S^{2n+1})$ is a flat $\ku_\ast(\Omega S^1) \cong \Z'[\beta, a^{\pm 1}]$-module, so we obtain an isomorphism of graded $\pi_\ast \ku_{\U(n)}$-algebras
    $$\ku_\ast^{\U(n)}(\Omega \CP^n) \cong \ku_\ast^{\U(n)}(\Omega S^{2n+1}) \otimes_{\Z'[\beta, a^{\pm 1}]} \Z'[\beta] \cong \ku_\ast^{\U(n)}(\Omega S^{2n+1})/(a-1).$$
    The desired result follows from the calculation of \cref{prop: U(n) homology Loops S2n+1}: the class $\tfrac{a-1}{c_n}\in \ku_\ast^{\U(n)}(\Omega S^{2n+1})$ is sent to the class denoted $b\in \ku_{2n}^{\U(n)}(\Omega \CP^n)$ under the generalized Hopf fibration.
\end{proof}
\begin{remark}
    The generalized Hopf fibration above also shows that
    $$\ku_\ast^{T^n}(\Omega \CP^n) \cong \ku_\ast^{T^n}(\Omega S^{2n+1})/(a-1) \cong \co_{\Bl_V \times_{\GG_m} \{1\}}.$$
    In particular, \cref{lem: relation BV and CV} implies that there is an isomorphism $\spec \ku_\ast^{T^n}(\Omega \CP^n) \cong \cC_V(-\dim_\RR(V))$, and hence an isomorphism
    $$\spec \ku_\ast^{\U(n)}(\Omega \CP^n) \cong \cC_V(-\dim_\RR(V))\mmod \Sigma_n.$$
\end{remark}



Before proceeding, we need the following lemma.
\begin{lemma}\label{lem: some codim stuff}
    The following statements hold:
    \begin{enumerate}
        \item Let $G_2 \subseteq G_1$ be a closed subgroup scheme. Let $Y \to Z$ be an open immersion of schemes with $G_2$-action whose complement has codimension $\geq d$. Then the induced map $\Ind_{G_2}^{G_1} Y \to \Ind_{G_2}^{G_1} Z$ is an open immersion of schemes with $G_1$-action whose complement has codimension $\geq d$.
        \item Let $Y$ be an integral quasi-affine variety such that $\pi_0 \Gamma(Y; \co_Y)$ is Noetherian. Then the map $Y \to \ol{Y}$ to its affine closure is an open immersion whose complement has codimension $\geq 2$.
    \end{enumerate}
\end{lemma}
\begin{proof}
    Part (a) is clear.
    For part (b), let $Z\subseteq \ol{Y}$ denote a closed subscheme of $Y$ associated to a minimal prime in the complement of $Y\subseteq \ol{Y}$, so that $Y\subseteq \ol{Y} - Z$. Then there are maps
    $$\H^0(\co_{\ol{Y}}) \xar{f} \H^0(\co_{\ol{Y} - Z}) \xar{g} \H^0(\co_Y).$$
    The map $g$ is an isomorphism, and the composite is also an isomorphism (by assumption). Therefore, $f$ is also an isomorphism. We claim that this forces $\fr{q}$ is necessarily of height $\geq 2$. Indeed, let $A = \H^0(\co_{\ol{Y}})$, so that $A$ is Noetherian. The desired claim then follows from the more general observation: if $A$ is Noetherian and $\fr{p}\subseteq A$ is a height one prime ideal corresponding to a closed subscheme $Z\subseteq \spec(A)$ with complement $U\subseteq \spec(A)$, the inclusion $A\subseteq \Gamma(\co_{U})$ is strict.
    To see this, standard arguments reduce us to the case when $A$ is local. Since $A$ is Noetherian, $\fr{p}$ is the radical of any nonzero $a\in \fr{p}$. Therefore, the inclusion $A\subseteq \Gamma(\co_{U})$ corresponds to the inclusion $A\subseteq A[1/a]$, which is evidently strict.
\end{proof}

\begin{construction}\label{cstr: CPn kappa ku}
    Let $\ld{V}_\beta = \AA^1 \times T^n_\beta \mmod \Sigma_n$; we will denote a point of $\ld{X}$ by $(c_0, c_1, \cdots, c_n)$ where $c_j$ lives in weight $-2j$. Write $\ld{V}$ to denote $\ld{V}_\beta/\beta$.

    Let $\kappa: T^n_\beta \mmod \Sigma_n \to \ld{V}_\beta$ denote the map $\vec{c} = (c_1, \cdots, c_n) \mapsto (1, \vec{c})$. 
    Equip $\GL_2$ with the grading where the entries of a matrix $\begin{psmallmatrix}
        a & b\\
        c & d
    \end{psmallmatrix}$ have the following weights: $a$ and $d$ live in weight zero, $b$ lives in weight $2n$, and $c$ lives in weight $-2n$.
    Let $\Mir_2(-2n)$ denote the mirabolic subgroup of $\GL_2$ of matrices of the form $\begin{psmallmatrix}
        a & b\\
        0 & 1
    \end{psmallmatrix}$, so that $\Mir_2(-2n)$ is an extension of $\GG_m$ by $\GG_a(-2n)$.
    There is an action of $\Mir_2(-2n)$ on $\ld{V}_\beta$, where $\begin{psmallmatrix}
        a & b\\
        0 & 1
    \end{psmallmatrix}\in \Mir_2(-2n)$ acts via
    $$\ld{V}_\beta \ni (c_0, \cdots, c_n) \to (ac_0 - bc_n, c_1, \cdots, c_n).$$
\end{construction}
\begin{prop}\label{prop: ku reg centr for S2n+1}
    There is an isomorphism of graded group schemes over $T^n_\beta \mmod \Sigma_n$:
    \begin{align*}
        \spec \ku^{\U(n)}_\ast(\Omega S^{2n+1}) & \cong T^n_\beta \mmod \Sigma_n \times_{\ld{V}_\beta/\Mir_2(-2n)} T^n_\beta \mmod \Sigma_n.
    \end{align*}
    Moreover, the $\Mir_2(-2n)$-orbit of $\kappa(T^n_\beta \mmod \Sigma_n) \subseteq \ld{V}_\beta$ has complement of codimension $\geq 2$.
\end{prop}
\begin{proof}
    By \cref{cor: U(n) homology Loops CPn}, it suffices to show that there is an isomorphism of graded schemes over $T^n_\beta \mmod \Sigma_n$:
    \begin{align*}
        T^n_\beta \mmod \Sigma_n \times_{\ld{V}_\beta/\Mir_2(-2n)} T^n_\beta \mmod \Sigma_n & \cong \spec \Z'[\beta, c_1, \cdots, c_n, \tfrac{1}{1 + \beta c_1 + \cdots + \beta^n c_n}, a^{\pm 1}, \tfrac{a-1}{c_n}].
    \end{align*}
    There is a closed immersion
    $$T^n_\beta \mmod \Sigma_n \times_{\ld{V}_\beta/\Mir_2(-2n)} T^n_\beta \mmod \Sigma_n \hookrightarrow T^n_\beta \mmod \Sigma_n \times \Mir_2(-2n)$$
    which exhibits $T^n_\beta \mmod \Sigma_n \times_{\ld{V}_\beta/\Mir_2(-2n)} T^n_\beta \mmod \Sigma_n$ as the subscheme of pairs $(\vec{c}, \begin{psmallmatrix}
        a & b\\
        0 & 1
    \end{psmallmatrix})$ such that $b$ stabilizes $\kappa(\vec{c})$. But $\begin{psmallmatrix}
        a & b\\
        0 & 1
    \end{psmallmatrix}$ sends $\kappa(\vec{c}) \mapsto (a - bc_n, \vec{c})$, so the necessary condition is that $b = \tfrac{a-1}{c_n}$, as desired.
    This also shows that the $\Mir_2(-2n)$-orbit of $\kappa(T^n_\beta \mmod \Sigma_n)$ is the complement of the closed subscheme $(0, \ast, \cdots, \ast, 0) \subseteq \ld{V}_\beta$. This closed subscheme has codimension $\geq 2$, as desired.
\end{proof}
\begin{observe}
    Equip $\ld{G}_X = \GL_2$ with the grading via the action of $2n\rho_{\ld{G}_X}$, so that if $\begin{psmallmatrix}
        a & b\\
        c & d
    \end{psmallmatrix}\in \GL_2$, the elements $a$ and $d$ have weight $0$, $b$ has weight $2n$, and $c$ has weight $-2n$. Let $V$ denote the affine space $\AA^2(2n,0)$, so that there is an action of $\ld{G}_X = \GL_2$ on $V$ via
    $$\begin{psmallmatrix}
        a & b\\
        c & d
    \end{psmallmatrix} \cdot (x,y) = (ax + cy, bx + dy);$$
    here, $x$ lives in degree $-2n$ and $y$ lives in degree zero. There is an isomorphism
    $$\Mir_2(-2n)\backslash \ld{G}_X \cong V - \{0\}, \ \begin{psmallmatrix}
        a & b\\
        c & d
    \end{psmallmatrix} \mapsto (c,d),$$
    and the above action of $\ld{G}_X$ on $V$ restricts on $V - \{0\}$ to the right-action of $\ld{G}_X$ on $\Mir_2(-2n)\backslash \ld{G}_X$.

    There is a $\ld{G}_X$-equivariant fibration
    \begin{equation}\label{eq: Vbeta non affine closure fibration}
        \Ind_{\Mir_2(-2n)}^{\ld{G}_X} \ld{V}_\beta \to \Mir_2(-2n)\backslash \ld{G}_X \cong V - \{0\}
    \end{equation}
    whose fibers are isomorphic to $\ld{V}_\beta$. 
    Let $\ol{\Ind_{\Mir_2(-2n)}^{\ld{G}_X} \ld{V}_\beta}$ denote the affine closure of $\Ind_{\Mir_2(-2n)}^{\ld{G}_X} \ld{V}_\beta$, so that there is a $\ld{G}_X$-equivariant fibration
    $$\ol{\Ind_{\Mir_2(-2n)}^{\ld{G}_X} \ld{V}_\beta} \to \ol{\Mir_2(-2n)\backslash \ld{G}_X} \cong V$$
    whose fibers are isomorphic to $\ld{V}_\beta$. 
    Finally, let $\ld{M}^\ddag_\beta$ denote the induction
    $$\ld{M}^\ddag_\beta := \Ind^{\ld{G}}_{\ld{G}_X} \ol{\Ind_{\Mir_2(-2n)}^{\ld{G}_X} \ld{V}_\beta}.$$
    \cref{lem: some codim stuff} implies:
    \begin{lemma}\label{lem: ind type An codim 2}
        There is an open immersion
        $$\Ind_{\Mir_2(-2n)}^{\ld{G}} \ld{V}_\beta \hookrightarrow \Ind^{\ld{G}}_{\ld{G}_X} \ol{\Ind_{\Mir_2(-2n)}^{\ld{G}_X} \ld{V}_\beta} = \ld{M}^\ddag_\beta$$
        which exhibits the target as the affine closure of the source, and whose complement is of codimension $\geq 2$.
    \end{lemma}
    \begin{proof}
        By \cref{lem: some codim stuff}, it suffices to show that there is an open immersion
        $$\Ind_{\Mir_2(-2n)}^{\ld{G}_X} \ld{V}_\beta \hookrightarrow \ol{\Ind_{\Mir_2(-2n)}^{\ld{G}_X} \ld{V}_\beta}$$
        which exhibits the target as the affine closure of the source, and whose complement is of codimension $\geq 2$. The statement about being the affine closure is true by definition, and the fact that the complement is of codimension $\geq 2$ is a consequence of the fact that the open subscheme $\Mir_2(-2n)\backslash \ld{G}_X \hookrightarrow \ol{\Mir_2(-2n)\backslash \ld{G}_X}$ has complement of codimension $2$.
    \end{proof}
    The map $\kappa: T^n_\beta \mmod \Sigma_n \to \ld{V}_\beta$ defines a locally closed immersion
    $$T^n_\beta \mmod \Sigma_n \xar{\kappa} \ld{V}_\beta \hookrightarrow \Ind_{\Mir_2(-2n)}^{\ld{G}} \ld{V}_\beta \hookrightarrow \Ind^{\ld{G}}_{\ld{G}_X} \ol{\Ind_{\Mir_2(-2n)}^{\ld{G}_X} \ld{V}_\beta} \cong \ld{M}^\ddag_\beta.$$
    We will denote this map by $\kappa_{\ld{M}^\ddag_\beta}$.
\end{observe}
\begin{remark}
    Here (and in the remaining sections), the stack $\ld{M}^\ddag_\beta$ and its fiber $\ld{M}^\ddag$ over $\beta = 0$ only act as crutches. \cref{lem: wind and ind} implies that $\ld{M}^\ddag$ is isomorphic to the dual variety $\ld{M}$ of \cref{thm: rk 1 bzsv is true}, but we have opted to use different notation (as in \cref{eq: M as G mod JX}.
\end{remark}
\begin{lemma}
    The $\ld{G}$-orbit of the image of $\kappa_{\ld{M}^\ddag_\beta}$ has complement of codimension $\geq 2$.
    Moreover, there is an isomorphism of graded group schemes over $T^n_\beta \mmod \Sigma_n$:
    \begin{align*}
        \spec \ku^{\U(n)}_\ast(\Omega S^{2n+1}) & \cong T^n_\beta \mmod \Sigma_n \times_{\ld{M}^\ddag_\beta/\ld{G}(-2\rho)} T^n_\beta \mmod \Sigma_n.
    \end{align*}
\end{lemma}
\begin{proof}
    For the first statement, \cref{lem: ind type An codim 2} implies that it suffices to show that the $\ld{G}$-orbit of the image of the composite
    $$T^n_\beta \mmod \Sigma_n \xar{\kappa} \ld{V}_\beta \hookrightarrow \Ind_{\Mir_2(-2n)}^{\ld{G}} \ld{V}_\beta$$
    has complement of codimension $\geq 2$. 
    Let $\ld{V}_\beta^\reg$ denote the $\Mir_2(-2n)$-orbit of $\kappa(T^n_\beta \mmod \Sigma_n)$ inside $\ld{V}_\beta$. Applying \cref{lem: some codim stuff} to the inclusion $\ld{V}_\beta^\reg \hookrightarrow \ld{V}_\beta$, it suffices to show that $\ld{V}_\beta^\reg \subseteq \ld{V}_\beta$ has complement of codimension $\geq 2$; but this is precisely \cref{prop: ku reg centr for S2n+1}.
\end{proof}
The above lemma combined with \cref{prop: ku rk 1 reg centr --> thm} and the isomorphism 
$$\ld{M}^\ddag_\beta/\ld{G}(-2\rho) \cong (\ol{\Ind_{\Mir_2(-2n)}^{\ld{G}_X} \ld{V}_\beta})/\ld{G}_X(-2n\rho_{\ld{G}_X})$$
implies:
\begin{corollary}\label{cor: bzsv for ku in An}
    Recall that $G = \GL_{n+1}$ and $H = \GL_n$, so $\ld{G} = \GL_{n+1}$ and $\ld{G}_X = \GL_2$.
    There is an equivalence of $\sh^{1/2}(\Z'[\beta])$-linear $\infty$-categories
    $$\Shv_{G\pw{t}}^{c,\Sat}(G\ls{t}/H\ls{t}; \ku)^\faux \simeq \Perf(\sh^{1/2} \ol{\Ind_{\Mir_2(-2n)}^{\ld{G}_X} \ld{V}_\beta}/\ld{G}_X(-2n\rho_{\ld{G}_X})).$$
\end{corollary}

The following simple observation is helpful for bookkeeping weights.
\begin{lemma}\label{lem: weighted cotangent}
    There is a graded isomorphism $T^\ast(j) \AA^2(m,n) = \AA^2(m,n) \times \AA^2(j-m, j-n)$.
\end{lemma}

\begin{proof}[Proof of \cref{thm: bzsv for An}]
    There is an isomorphism $\ld{V}_\beta/\beta \cong \fr{t}^{n-1}(2)\mmod \Sigma_{n-1} \times \AA^2(0,2n)$, so \cref{lem: weighted cotangent} and the fibration \cref{eq: Vbeta non affine closure fibration} defines an isomorphism
    $$\Ind_{\Mir_2(-2n)}^{\ld{G}_X} \ld{V} \cong \fr{t}^{n-1}(2)\mmod \Sigma_{n-1} \times T^\ast(2n)(V - \{0\}),$$
    and hence 
    $$\ol{\Ind_{\Mir_2(-2n)}^{\ld{G}_X} \ld{V}} \cong \fr{t}^{n-1}(2)\mmod \Sigma_{n-1} \times T^\ast(2n)(V).$$
    \cref{cor: bzsv for ku in An} now implies \cref{thm: bzsv for An}.
\end{proof}
\begin{corollary}\label{cor: bzsv for CPn}
    Let $G = \PGL_{n+1}$ and $H = \GL_n$, so $\ld{G} = \SL_{n+1}$ and $\ld{G}_X = \SL_2$.
    There is an equivalence of $\sh^{1/2}(\Z'[\beta])$-linear $\infty$-categories
    $$\Shv_{G\pw{t}}^{c,\Sat}(G\ls{t}/H\ls{t}; \ku)^\faux \simeq \Perf(\sh^{1/2} \ld{M}^\ddag_\beta/\SL_{n+1}[-2\rho]).$$
    When $\beta = 0$, this specializes to an equivalence of $\QQ$-linear $\infty$-categories
    $$\Shv_{G\pw{t}}^{c,\Sat}(G\ls{t}/H\ls{t}; \QQ) \simeq \Perf(T^\ast[2n](\AA^2[2n,0])/\SL_2[-2n\rho_{\ld{G}_X}] \times \fr{gl}_{n-1}[2]\mmod \GL_{n-1}).$$
\end{corollary}
\begin{proof}
    Recall that $\PGL_{n+1} = \GL_{n+1}/\GL_1^\mathrm{diag}$, so that there are equivalences (where the second line comes from \cref{cor: bzsv for ku in An})
    \begin{multline*}
        \Shv_{\PGL_{n+1}\pw{t}}^{c,\Sat}(\PGL_{n+1}\ls{t}/\GL_n\ls{t}; \ku)^\faux \\
        \simeq \Shv_{\GL_{n+1}\pw{t}}^{c,\Sat}(\GL_{n+1}\ls{t}/\GL_n\ls{t}; \ku)^\faux \otimes_{\Shv_{\GL_1\pw{t}}(\GL_1\ls{t}; \ku)^\faux} \LMod_{\sh^{1/2} \Z'[\beta]} \\
        \simeq \Perf(\sh^{1/2} \ld{M}^\ddag_\beta/\GL_{n+1}[-2\rho]) \otimes_{\Perf(B\GG_m)} \LMod_{\sh^{1/2} \Z'[\beta]} \\
        \simeq \Perf\left(\sh^{1/2} \ld{M}^\ddag_\beta/\GL_{n+1}[-2\rho] \times_{B\GG_m} \spec \sh^{1/2} \Z'[\beta]\right)
    \end{multline*}
    However, the displayed fiber product is precisely $\sh^{1/2} \ld{M}^\ddag_\beta/\SL_{n+1}[-2\rho]$. The claim about identifying its reduction modulo $\beta$ with $T^\ast[2n](\AA^2[2n,0])/\SL_2[-2\rho_{\SL_2}] \times \fr{gl}_{n-1}[2]\mmod \GL_{n-1}$ follows from the construction of $\ld{M}^\ddag_\beta$.
\end{proof}

\begin{remark}\label{rmk: type N PGL2 mod PO2}
    In the case of the spherical variety $\PGL_2/\GG_m$, \cref{cor: bzsv for CPn} states that there is an equivalence
    $$\Shv_{\PGL_2\pw{t}}^{c,\Sat}(\PGL_2\ls{t}/\GG_m\ls{t}; \QQ) \simeq \Perf(T^\ast[2n](\AA^2[2n,0])/\SL_2[-2n\rho_{\ld{G}_X}]).$$
    The $\PGL_2$-variety $\PGL_2/\GG_m$ has a natural action of the Weyl group $\Z/2 = \N_{\PGL_2}(\GG_m)/\GG_m$ (under the homotopy equivalence $(\PGL_2/\GG_m)(\cc) \simeq S^2$, this is the antipodal action). This equips the left-hand side of the above equivalence with a natural $\Z/2$-action. One can show that under this equivalence, the resulting $\Z/2$-action on the right-hand side identifies with the natural $\Z/2$-action on $T^\ast(\AA^2)$ via the symplectic form. 
    That is, the action sends
    $$(v_1, v_2), (w_1, w_2) \mapsto (-w_2, w_1), (v_2, -v_1).$$
    See \cref{conj: spherical levi} and the surrounding discussion for an expected generalization to arbitrary semisimple algebraic groups. The normalizer $\N_{\PGL_2}(\GG_m)$ can be identified with $\mathrm{PO}_2$, which implies that there is an equivalence
    $$\Shv_{\PGL_2\pw{t}}^{c,\Sat}(\PGL_2\ls{t}/\mathrm{PO}_2\ls{t}; \QQ) \simeq \Perf(T^\ast[2](\AA^2[2,0])/(\SL_2[-2\rho_{\ld{G}_X}] \times \Z/2)).$$
    Note that the spherical root of $\PGL_2/\mathrm{PO}_2$ is (by definition) of type N, and so  the spherical $\PGL_2$-variety $\PGL_2/\mathrm{PO}_2$ is excluded by \cite{sakellaridis-venkatesh, bzsv}.
    Nevertheless, the preceding equivalence shows that the Hamiltonian $\SL_2$-``space'' which should be dual (in the sense of \cite{bzsv}) to $\PGL_2/\mathrm{PO}_2$ is the \textit{stack} $T^\ast(\AA^2)/(\Z/2)$.
\end{remark}

%% file: each-type/type-Bn.tex
\subsection{Type $B_n$: $\SO_{2n+1}/\SO_{2n}$}

Our goal in this section is to prove \cref{thm: rk 1 bzsv is true} in type $B_n$, i.e., for the spherical $\SO_{2n+1}$-variety $\SO_{2n+1}/\SO_{2n}$. Note that if $(V,q)$ is a quadratic space and $v\in V$ with $q(v) = 1$, then $\SO_V/\SO_{v^\perp}$ can be identified with the hyperboloid $\{w\in V | q(w) = 1\}$.
Write $G = \SO_{2n+1}$ and $H = \SO_{2n}$, so that $\ld{G} = \Sp_{2n}$. Recall from \cref{table: topology and dualities for rank 1 spherical varieties} that $\ld{G}_X = \SL_2$. In this section, we will only consider coefficients in $\Z'$ (instead of $\ku$).
\begin{theorem}[\cref{thm: rk 1 bzsv is true} in type $B_n$]\label{thm: bzsv for Bn}
    There is an equivalence of $\QQ$-linear $\infty$-categories
    $$\Shv_{G\pw{t}}^{c,\Sat}(G\ls{t}/H\ls{t}; \QQ) \simeq \Perf(T^\ast[2n](\AA^2[4n-2,0])/\SL_2[-(4n-2)\rho_{\ld{G}_X}] \times \fr{sp}_{2n-2}[2]\mmod \Sp_{2n-2}).$$
\end{theorem}

The proof of \cref{thm: bzsv for Bn} will take up the remainder of this section.

\begin{lemma}
    There is a homotopy equivalence $(\SO_{2n+1}/\SO_{2n})(\cc) \simeq S^{2n}$. Moreover, if $B\subseteq G$ is the Borel subgroup of upper-triangular matrices, the Levi quotient $L(X)$ of the parabolic subgroup stabilizing the open $B$-orbit in $\SO_{2n+1}/\SO_{2n}$ is given by $\SO_{2n-1} \times \GG_m$.
\end{lemma}
\begin{lemma}\label{lem: SO2n equiv coh of point}
    Let $W = (\Z/2)^{n-1} \rtimes \Sigma_n$ denote the Weyl group of $\SO_{2n}$. Then there is an isomorphism
    $$\H^\ast_{\SO_{2n}}(\ast; \Z') \cong \Z'[p_1, \cdots, p_{n-1}, c_n],$$
    where the injective map $\H^\ast_{\SO_{2n}}(\ast; \Z') \to \H^\ast_{T^{n}}(\ast; \Z')$ sends $p_j$ to the $j$th elementary symmetric polynomial in the variables $x_1^2, \cdots, x_{n-1}^2$ (so $p_j$ lives in weight $-4j$), and $c_n \mapsto x_1 \cdots x_n$.
\end{lemma}
\begin{prop}\label{prop: SO2n homology Loops S2n}
    There is a graded isomorphism of $\H_{\SO_{2n}}^\ast(\ast; \Z')$-algebras
    $$\H^{\SO_{2n}}_\ast(\Omega S^{2n}; \Z') \cong \Z'[p_1, \cdots, p_{n-1}, c_n, b]/bc_n,$$
    where $b$ lives in weight $4n-2$. (In particular, $\H^{\SO_{2n}}_\ast(\Omega S^{2n}; \Z')$ is \textit{not} flat over $\H_{\SO_{2n}}^\ast(\ast; \Z')$.)
\end{prop}
\begin{proof}
    The restriction of the $\SO_{2n}$-action on $S^{2n}$ to the maximal torus $T^n \subseteq \SO_{2n}$ exhibits $S^{2n}$ as the one-point compactification of the standard $n$-dimensional complex representation $\std$. \cref{cor: G-equiv homology of Loops SV} implies that there is a graded isomorphism of $\H_{\SO_{2n}}^\ast(\ast; \Z')$-algebras
    $$\H^{\SO_{2n}}_\ast(\Omega S^{2n}; \Z') \cong \left(\Z'[x_1, \cdots, x_n, b]/bx_1 \cdots x_n\right)^W,$$
    where $b$ lives in weight $4n-2$. The $W$-invariants on the right-hand side can be computed using \cref{lem: SO2n equiv coh of point} (note that the action of $W$ on $b$ is trivial), and gives the desired calculation.
\end{proof}
\begin{definition}\label{def: kostant Bn}
    Let $\ld{V}$ denote the graded affine space $\AA^1(2-2n) \times \fr{t}^n(2)\mmod W$. There is an action of $\GG_a(2-4n)$ on $\ld{V}$, where $\begin{psmallmatrix}
        1 & b\\
        0 & 1
    \end{psmallmatrix}\in \GG_a(2-4n)$ sends $(z,\vec{p}, c_n) \mapsto (z - bc_n, \vec{p}, c_n)$. Note that $b$ lives in weight $4n-2$.
    
    Equip $\ld{G}_X = \SL_2$ with the grading coming from $(4n-2)\rho_{\SL_2}$, so that the entries of a matrix $\begin{psmallmatrix}
        a & b\\
        c & d
    \end{psmallmatrix}$ are equipped with the following weights: $a$ and $d$ have weight $0$, $b$ has weight $4n-2$, and $c$ has weight $2-4n$. Let $V$ denote the affine space $\AA^2(4n-2, 0)$, so that there is an action of $\ld{G}_X = \SL_2$ on $\ld{V}$ via
    $$\begin{psmallmatrix}
        a & b\\
        c & d
    \end{psmallmatrix} \cdot (x,y) = (ax + cy, bx + dy);$$
    here, $x$ lives in degree $2-4n$ and $y$ lives in degree $0$.
    There is an isomorphism
    $$\GG_a(2-4n)\backslash \ld{G}_X \cong \ld{V} - \{0\}, \ \begin{psmallmatrix}
        a & b\\
        c & d
    \end{psmallmatrix} \mapsto (c,d),$$
    and the above action of $\ld{G}_X$ on $\ld{V}$ restricts on $\ld{V} - \{0\}$ to the right-action of $\ld{G}_X$ on $\GG_a(2-4n)\backslash \ld{G}_X$.
    
    Let $\kappa: \fr{t}^n(2)\mmod W \to \ld{V}$ denote the map sending $(\vec{p}, c_n) \mapsto (0, \vec{p}, c_n)$.
    Let $W'$ denote the Weyl group of $\SO_{2n-1}$, so \cref{lem: weighted cotangent} implies that there is an isomorphism
    $$\Ind_{\GG_a(2-4n)}^{\ld{G}_X} \ld{V} \cong \fr{t}^{n-1}(2)\mmod W' \times T^\ast(2n)(V - \{0\}).$$
    In particular, there is an open immersion
    $$\Ind_{\GG_a(2-4n)}^{\ld{G}_X} \ld{V} \hookrightarrow \fr{t}^{n-1}(2)\mmod W' \times T^\ast(2n)(V)$$
    which exhibits the target as the affine closure of $\Ind_{\GG_a(2-4n)}^{\ld{G}_X} \ld{V}$. Inducing along the map $\ld{G}_X \to \ld{G}$ produces an open immersion
    $$\Ind_{\GG_a(2-4n)}^{\ld{G}} \ld{V} \hookrightarrow \fr{t}^{n-1}(2)\mmod W' \times \Ind_{\ld{G}_X}^{\ld{G}} T^\ast(2n)(V),$$
    and we will write $\ld{M}^\ddag = \fr{t}^{n-1}(2)\mmod W' \times \Ind_{\ld{G}_X}^{\ld{G}} T^\ast(2n)(V)$.

    The map $\kappa$ defines a locally closed immersion
    $$\fr{t}^n(2)\mmod W \xar{\kappa} \ld{V} \hookrightarrow \Ind_{\GG_a(2-4n)}^{\ld{G}} \ld{V} \hookrightarrow \fr{t}^{n-1}(2)\mmod W' \times \Ind_{\ld{G}_X}^{\ld{G}} T^\ast(2n)(V) = \ld{M}^\ddag,$$
    which we will denote by $\kappa_{\ld{M}^\ddag}$.
\end{definition}
\begin{remark}
    In \cref{def: kostant Bn}, it does not make sense to ask that $\kappa: \fr{t}^n(2)\mmod W \to \ld{V}$ instead send $(\vec{p}, c_n) \mapsto (1, \vec{p}, c_n)$. Indeed, the point $1\in \AA^1(2-2n)$ is not well-defined, since it would have to be cut out by the ideal $(z-1)$, which is not homogeneous (i.e., is not a graded ideal).
\end{remark}

\begin{prop}\label{lem: constructible set quotient}
    Let $C = \AA^1 = \spec \Z'[c]$, and let $\AA^2 = \spec \Z'[z, c]$. Let $\ul{\GG_a}$ denote the constant group scheme over $C$ acting on $\AA^2$ by $b\cdot (z,c) = (z+bc, c)$. Let $f: C \to Y$ denote the map $c \mapsto (0,c)$, and let $J$ denote the stabilizer group scheme of the image of $f$ (over $C$). Then $J$ is not flat over $C$, and the quotient $\ul{\GG_a}/J$ is isomorphic to $\AA^2 - \{0\}$. Moreover, if $\ul{\SL_2}$ is the constant group scheme over $C$ with $J \subseteq \ul{\GG_a}$ embedded via the upper-triangular matrices, the affine closure of $\ul{\SL_2}/J$ can be identified with $T^\ast \AA^2$.
\end{prop}
\begin{proof}
    It is immediate that $J \cong \spec \Z'[b,c]/bc$, so is not flat over $C$. (This implies that the $\ul{\GG_a}$-orbit of the image of $f$ is only a constructible subset of $\AA^2$: indeed, this $\ul{\GG_a}$-orbit is the standard example $(\AA^2 - \{z=0\}) \cup \{(0,0)\}$ of a constructible subset.)
    For the claim about $\ul{\GG_a}/J$, let $\Mir_2 \subseteq \GL_2$ denote (as usual) the mirabolic subgroup of matrices of the form $\begin{psmallmatrix}
        a & b\\
        0 & 1
    \end{psmallmatrix}$, let $\ul{\Mir_2}$ denote the constant group scheme over $C$, and let $\tilde{J} \subseteq \Mir_2$ denote the subgroup of matrices cut out by the equation $a-1 = bc$. Then $\ul{\GG_a}$ (resp. $J$) is the kernel of the determinant $\ul{\Mir_2} \xar{\det} \ul{\GG_m}$ (resp. $\tilde{J} \to \ul{\Mir_2} \xar{\det} \ul{\GG_m}$). This implies that $\ul{\GG_a}/J \cong \ul{\Mir_2}/\tilde{J}$, and the latter is isomorphic to $\AA^2 - \{0\}$ via the map $\ul{\Mir_2} \to \AA^2 - \{0\}$ sending $\begin{psmallmatrix}
        a & b \\
        0 & 1
    \end{psmallmatrix}, c \mapsto (a+bc,c)$, as desired. (In the $\beta$-deformed case, this is \cref{prop: ku reg centr for S2n+1}.) Next, there is an isomorphism ${\SL_2}/{\GG_a} \cong \AA^2 - \{0\}$, which implies that $\ul{\SL_2}/J \cong \ul{\GL_2}/\tilde{J}$ is isomorphic to the complement of the zero section in $T^\ast(\AA^2 - \{0\})$. (It can alternatively be described as $\SL_2 \times^{\GG_a} (\AA^2 - \{0\})$.) Its affine closure is $T^\ast(\AA^2)$, as claimed.
\end{proof}
\begin{prop}\label{prop: ordinary reg centr type Bn}
    There is an isomorphism of graded group schemes over $\fr{t}^n(2)\mmod W$:
    $$\spec \H^{\SO_{2n}}_\ast(\Omega S^{2n}; \Z') \cong \fr{t}^n(2)\mmod W \times_{\ld{M}^\ddag/\ld{G}(-2\rho)} \fr{t}^n(2)\mmod W.$$
    Moreover, if $\ld{J}_{X}$ denotes the above group scheme over $\fr{t}^n(2)\mmod W$, the algebra of regular functions on $(\fr{t}^n(2)\mmod W \times \ld{G})/\ld{J}_{X}$ is isomorphic to $\co_{\ld{M}^\ddag}$.
\end{prop}
\begin{proof}
    There is an isomorphism
    $$\fr{t}^n(2)\mmod W \times_{\ld{M}^\ddag/\ld{G}(-2\rho)} \fr{t}^n(2)\mmod W \cong \fr{t}^n(2)\mmod W \times_{\ld{V}/\GG_a(2-4n)} \fr{t}^n(2)\mmod W,$$
    as well as a closed immersion
    $$\fr{t}^n(2)\mmod W \times_{\ld{V}/\GG_a(2-4n)} \fr{t}^n(2)\mmod W \hookrightarrow \fr{t}^n(2)\mmod W \times \GG_a(2-4n),$$
    which exhibits $\fr{t}^n(2)\mmod W \times_{\ld{V}/\GG_a(2-4n)} \fr{t}^n(2)\mmod W$ as the subscheme of pairs $(\vec{p}, c_n, b)$ such that $b$ stabilizes $\kappa(\vec{p}, c_n)$. But by definition of $\kappa$, this happens if and only if $bc_n = 0$. In other words, there is an isomorphism of graded schemes over $\fr{t}^n(2)\mmod W$:
    $$\fr{t}^n(2)\mmod W \times_{\ld{V}/\GG_a(2-4n)} \fr{t}^n(2)\mmod W \cong \spec \Z'[p_1, \cdots, p_{n-1}, c_n, b]/bc_n.$$
    The first part of the proposition therefore follows from \cref{prop: SO2n homology Loops S2n}.
    The second part of the proposition follows from \cref{lem: constructible set quotient} (rather, its obvious variant for graded affine spaces).
\end{proof}
\begin{proof}[Proof of \cref{thm: bzsv for Bn}]
    This follows from \cref{thm: ordinary homology criterion satake} and \cref{prop: ordinary reg centr type Bn}, along with the isomorphism between $\ld{M}^\ddag/\ld{G}(-2\rho)$ and $T^\ast(2n)(V)/\ld{G}_X(-(4n-2)\rho_{\ld{G}_X}) \times \fr{t}^{n-1}(2)\mmod W'$.
\end{proof}

\begin{remark}\label{rmk: SO2n+1 mod normalizer of SO2n}
    \cref{thm: bzsv for Bn} can be used to prove a variant for the rank one spherical $\SO_{2n+1}$-variety $\SO_{2n+1}/\N_{\SO_{2n+1}}(\SO_{2n})$, where $\N_{\SO_{2n+1}}(\SO_{2n})$ is the normalizer of $\SO_{2n} \subseteq \SO_{2n+1}$. The dual group $\ld{G}_X$ in this case is again $\SL_2$. The quotient $\N_{\SO_{2n+1}}(\SO_{2n})/\SO_{2n}$ is isomorphic to $\Z/2$, and it acts on $(\SO_{2n+1}/\SO_{2n})(\cc) \simeq S^{2n}$ via the antipodal action. As in \cref{rmk: type N PGL2 mod PO2}, this equips the left-hand side of the equivalence of \cref{thm: bzsv for Bn} with a natural $\Z/2$-action. One can show that under this equivalence, the resulting $\Z/2$-action on the right-hand side of \cref{thm: bzsv for Bn} identifies with the natural $\Z/2$-action on $T^\ast(\AA^2)$ via the symplectic form. Using this, one finds that there is an equivalence
    \begin{multline*}
        \Shv_{\SO_{2n+1}\pw{t}}^{c,\Sat}(\SO_{2n+1}\ls{t}/\N_{\SO_{2n+1}}(\SO_{2n})\ls{t}; \QQ) \simeq \\
        \Perf(T^\ast[2n](\AA^2[4n-2,0])/(\SL_2[-(4n-2)\rho_{\ld{G}_X}] \times \Z/2) \times \fr{sp}_{2n-2}[2]\mmod \Sp_{2n-2}).
    \end{multline*}
    The spherical root of $\SO_{2n+1}/\N_{\SO_{2n+1}}(\SO_{2n})$ is of type N, and so it is excluded by \cite{sakellaridis-venkatesh, bzsv}; nevertheless, the preceding equivalence shows that it does admit a dual, given by the \textit{stack} $\Ind_{\SL_2}^{\Sp_{2n}} (T^\ast(\AA^2)/(\Z/2) \times \AA^{n-1})$. 
\end{remark}

%% file: each-type/type-Cn.tex
\subsection{Type $C_n$: $\Sp_{2n}/(\Sp_2 \times \Sp_{2n-2})$}

Our goal in this section is to prove \cref{thm: rk 1 bzsv is true} in type $C_n$, i.e., for the spherical $\Sp_{2n}$-variety $\Sp_{2n}/(\Sp_2 \times \Sp_{2n-2})$. Let $G = \Sp_{2n}$ and $H = \Sp_2 \times \Sp_{2n-2}$, so that $\ld{G} = \SO_{2n+1}$. Recall from \cref{table: topology and dualities for rank 1 spherical varieties} that $\ld{G}_X = \SL_2$.
\begin{theorem}[\cref{thm: rk 1 bzsv is true} in type $C_n$]\label{thm: bzsv for Cn}
    There is an equivalence of $\QQ$-linear $\infty$-categories
    $$\Shv_{G\pw{t}}^{c,\Sat}(G\ls{t}/H\ls{t}; \QQ) \simeq \Perf(T^\ast[4n-4](\AA^2[4n-2,0])/\SL_2[-(4n-2)\rho_{\ld{G}_X}] \times \ld{\fr{h}}^\ast[2]\mmod \ld{H}),$$
    where $H = \Sp_2 \times \Sp_{2(n-2)}$.
\end{theorem}

The proof of \cref{thm: bzsv for Cn} will take up the remainder of this section.

\begin{lemma}
    There is a homotopy equivalence $(\Sp_{2n}/(\Sp_2 \times \Sp_{2n-2}))(\cc) \simeq \HHP^{n-1}$. Moreover, if $B\subseteq G$ is the Borel subgroup of upper-triangular matrices, the Levi quotient $L(X)$ of the parabolic subgroup stabilizing the open $B$-orbit in $\Sp_{2n}/(\Sp_2 \times \Sp_{2n-2})$ is given by $\GL_2 \times \Sp_{2n-4}$.
\end{lemma}
\begin{lemma}\label{lem: sp2n equiv coh of point}
    Let $W$ denote the Weyl group of $\Sp_{2n}$. Then there is an isomorphism
    $$\H_{\Sp_{2n}}^\ast(\ast; \QQ) \cong \QQ[p_1, \cdots, p_n],$$
    where the map $\H_{\Sp_{2n}}^\ast(\ast; \QQ) \to \H_{T^n}^\ast(\ast; \QQ)$ sends $p_j$ to the $j$th elementary symmetric polynomial in the variables $x_1\ol{x}_1, \cdots, x_n \ol{x}_n$. Here, $\ol{x} = -x$.
\end{lemma}
Although one can give an argument for the following result using \cref{thm: homology of loops SV}, it is simpler to give an argument ``from scratch''.
\begin{prop}\label{prop: ku Sp2n-2 homology Loops HPn-1}
    There is an isomorphism of graded $\H_{\Sp_{2n-2}}^\ast(\ast; \QQ)$-algebras
    $$\H^{\Sp_{2n-2}}_\ast(\Omega \HHP^{n-1}; \QQ) \cong \QQ[p_1, \cdots, p_{n-1}, b]/bp_{n-1},$$
    where $b$ lives in weight $4n-2$.
\end{prop}
\begin{proof}
    Let $T^{n-1} \subseteq \Sp_{2n-2}$ denote the maximal torus. Then the homotopy fixed points spectral sequence for $\pi_\ast \QQ[\Omega \HHP^{n-1}]^{hT^{n-1}}$ is given by
    $$E_2^{\ast,\ast} \cong \H_\ast(\Omega \HHP^{n-1}; \QQ) \otimes_{\QQ} \pi_\ast \QQ^{hT^{n-1}} \Rightarrow \pi_\ast \QQ[\Omega \HHP^{n-1}]^{hT^{n-1}}.$$
    To compute the $E_2$-page, observe that the Hopf fibration $S^3 \to S^{4n-1} \to \HHP^{n-1}$ implies that there is an equivalence $\Omega \HHP^{n-1} \simeq S^3 \times \Omega S^{4n-1}$. This gives an isomorphism
    $$\H_\ast(\Omega \HHP^{n-1}; \QQ) \cong \QQ[a,b]/a^2,$$
    where $a$ lives in weight $3$ and $b$ lives in weight $4n-2$. Therefore, 
    $$E_2^{\ast,\ast} \cong \QQ[a, b]\pw{x_1, \cdots, x_{n-1}}/a^2,$$
    where each $x_j$ lives in weight $-2$. Recall that the action of $T^{n-1}$ on $\HHP^{n-1}$ is induced by the inclusion $T^{n-1} \subseteq \Sp_{2n-2} \subseteq \U(4n-4)$ given by the representation $\bigoplus_{j=1}^{n-1} \std \oplus \std^{-1}$. This forces a single differential in the above spectral sequence, given by
    $$d_2(a) = bx_1 \ol{x}_1 \cdots x_{n-1} \ol{x}_{n-1}.$$
    After running this differential, the spectral sequence is concentrated in even degrees, and we find that there is an isomorphism
    $$\pi_\ast \QQ[\Omega \HHP^{n-1}]^{hT^{n-1}} \cong \QQ[b]\pw{x_1, \cdots, x_{n-1}}/bx_1 \ol{x}_1 \cdots x_{n-1} \ol{x}_{n-1}.$$
    To calculate $\H^{T^{n-1}}_\ast(\Omega \HHP^{n-1}; \QQ)$ itself (and not just its completion $\pi_\ast \QQ[\Omega \HHP^{n-1}]^{hT^{n-1}}$), the strategy of \cref{thm: homology of loops SV} reduces us to showing that the restriction $\H^{T^{n-1}}_\ast(\Omega \HHP^{n-1}; \QQ)|_{\punc{\fr{t}}} = \H_{T^{n-1}}^\ast(\ast; \QQ)$. By \cref{lem: atiyah localization}, there is an isomorphism $\H^{T^{n-1}}_\ast(\Omega \HHP^{n-1}; \QQ)|_{\punc{\fr{t}}} \cong \H^{T^{n-1}}_\ast(\Omega (\HHP^{n-1}; \QQ)^{T^{n-1}})|_{\punc{\fr{t}}}$. It therefore suffices to show that $\Omega (\HHP^{n-1})^{T^{n-1}}$ is contractible, but this is a consequence of the simple observation that $(\HHP^{n-1})^{T^{n-1}} \cong S^0$. This discussion gives an isomorphism of graded $\H_{T^{n-1}}^\ast(\ast; \QQ)$-algebras
    $$\H^{T^{n-1}}_\ast(\Omega \HHP^{n-1}; \QQ) \cong \QQ[x_1, \cdots, x_{n-1}, b]/bx_1 \ol{x}_1 \cdots x_{n-1} \ol{x}_{n-1}.$$
    This isomorphism is $W$-equivariant (where $W$ is the Weyl group of $\Sp_{2n-2}$), so \cref{prop: Weyl invts} implies that there is an isomorphism of graded $\H_{\Sp_{2n-2}}^\ast(\ast; \QQ)$-algebras
    $$\H^{\Sp_{2n-2}}_\ast(\Omega \HHP^{n-1}; \QQ) \cong \left(\QQ[x_1, \cdots, x_{n-1}, b]/bx_1 \ol{x}_1 \cdots x_{n-1} \ol{x}_{n-1}\right)^W.$$
    Noting that the action of $W$ leaves $b$ invariant, \cref{lem: sp2n equiv coh of point} computes the right-hand side; the resulting answer is precisely the right-hand side of the proposition.
\end{proof}
\begin{notation}
    Let $W$ denote the Weyl group of $\Sp_2 \times \Sp_{2n-2}$, so that it is the product of $\Z/2$ with the Weyl group of $\Sp_{2n-2}$. There is a natural action of $W$ on the torus $\GG_m \times T^{n-1} \subseteq \Sp_2 \times \Sp_{2n-2}$, and hence an action of $W$ on $\fr{t}^n$. It is an easy consequence of \cref{lem: sp2n equiv coh of point} that there is an isomorphism
    $$\fr{t}^n \mmod W \cong \spec \QQ[p_1', p_1, \cdots, p_{n-1}],$$
    where $p_1'$ lives in weight $-4$, and $p_j$ lives in weight $-4j$.
\end{notation}
\begin{construction}
    Let $\ld{V}$ denote the graded affine scheme $\AA^1(-2) \times \fr{t}^n \mmod W$, and let $\kappa: \fr{t}^n \mmod W \to \ld{V}$ denote the map sending $(p_1', \vec{p}) \mapsto (0, p_1', \vec{p})$. There is an action of $\GG_a(2-4n)$ on $\ld{V}$, where $b\in \GG_a(2-4n)$ sends
    $$(z, p_1', \vec{p}) \mapsto (z - bp_{n-1}, p_1', \vec{p}).$$
    Equip $\ld{G}_X = \SL_2$ with the grading coming from $(4n-2)\rho_{\SL_2}$, so that the entries of a matrix $\begin{psmallmatrix}
        a & b\\
        c & d
    \end{psmallmatrix}$ are equipped with the following weights: $a$ and $d$ have weight $0$, $b$ has weight $4n-2$, and $c$ has weight $2-4n$. Let $\ld{V}$ denote the affine space $\AA^2(4n-2, 0)$, so that there is an action of $\ld{G}_X = \SL_2$ on $\ld{V}$ via
    $$\begin{psmallmatrix}
        a & b\\
        c & d
    \end{psmallmatrix} \cdot (x,y) = (ax + cy, bx + dy);$$
    here, $x$ lives in degree $2-4n$ and $y$ lives in degree $0$.
    There is an isomorphism
    $$\GG_a(2-4n)\backslash \ld{G}_X \cong V - \{0\}, \ \begin{psmallmatrix}
        a & b\\
        c & d
    \end{psmallmatrix} \mapsto (c,d),$$
    and the above action of $\ld{G}_X$ on $\ld{V}$ restricts on $\ld{V} - \{0\}$ to the right-action of $\ld{G}_X$ on $\GG_a(2-4n)\backslash \ld{G}_X$.
    There is a fibration
    $$\Ind_{\GG_a(2-4n)}^{\ld{G}_X} \ld{V} \to \GG_a(2-4n)\backslash \ld{G}_X \cong \ld{V} - \{0\}$$
    whose fibers are isomorphic to $\ld{V}$. Let $\ol{\Ind_{\GG_a(2-4n)}^{\ld{G}_X} \ld{V}}$ denote the affine closure of $\Ind_{\GG_a(2-4n)}^{\ld{G}_X} \ld{V}$, so that there is a $\ld{G}_X$-equivariant fibration
    $$\ol{\Ind_{\GG_a(2-4n)}^{\ld{G}_X} \ld{V}} \to \ol{\GG_a(2-4n)\backslash \ld{G}_X} \cong V$$
    whose fibers are isomorphic to $\ld{V}$. Let $\ld{M}^\ddag$ denote the induction
    $$\ld{M}^\ddag = \Ind_{\ld{G}_X}^{\ld{G}} \ol{\Ind_{\GG_a(2-4n)}^{\ld{G}_X} \ld{V}}.$$
    \cref{lem: some codim stuff} implies:
    \begin{lemma}
        There is an open immersion
        $$\Ind_{\GG_a(2-4n)}^{\ld{G}} \ld{V} \to \Ind_{\ld{G}_X}^{\ld{G}} \ol{\Ind_{\GG_a(2-4n)}^{\ld{G}_X} \ld{V}} = \ld{M}^\ddag$$
        which exhibits the target as the affine closure of the source, and whose complement is of codimension $\geq 2$.
    \end{lemma}
    The map $\kappa: \fr{t}^n \mmod W \to \ld{V}$ defines a locally closed immersion
    $$\fr{t}^n \mmod W \xar{\kappa} \ld{V} \hookrightarrow \Ind_{\GG_a(2-4n)}^{\ld{G}} \ld{V} \hookrightarrow \Ind_{\ld{G}_X}^{\ld{G}} \ol{\Ind_{\GG_a(2-4n)}^{\ld{G}_X} \ld{V}},$$
    which we will denote by $\kappa_{\ld{M}^\ddag}$.
\end{construction}
\begin{prop}\label{prop: ku reg centr type Cn}
    There is an isomorphism of graded group schemes over $\fr{t}^n \mmod W$:
    $$\spec \H^{\Sp_2 \times \Sp_{2n-2}}_\ast(\Omega \HHP^{n-1}; \QQ) \cong \fr{t}^n \mmod W \times_{\ld{M}^\ddag/\ld{G}(-2\rho)} \fr{t}^n \mmod W.$$
    Moreover, if $\ld{J}_X$ denotes the above group scheme over $\fr{t}^n \mmod W$, the algebra of regular functions on $(\fr{t}^n \mmod W \times \ld{G})/\ld{J}_X$ is isomorphic to $\co_{\ld{M}^\ddag}$.
\end{prop}
\begin{proof}
    There is an isomorphism 
    $$\fr{t}^n \mmod W \times_{\ld{M}^\ddag/\ld{G}(-2\rho)} \fr{t}^n \mmod W \cong \fr{t}^n \mmod W \times_{\ld{V}/\GG_a(2-4n)} \fr{t}^n \mmod W,$$
    as well as a closed immersion
    $$\fr{t}^n \mmod W \times_{\ld{V}/\GG_a(2-4n)} \fr{t}^n \mmod W \hookrightarrow \fr{t}^n \mmod W \times \GG_a(2-4n),$$
    which exhibits $\fr{t}^n \mmod W \times_{\ld{V}/\GG_a(2-4n)} \fr{t}^n \mmod W$ as the subscheme of tuples $(p_1', \vec{p}, b)$ such that $b$ stabilizes $\kappa(p_1', \vec{p})$. By definition of $\kappa$, this happens if and only if $bp_{n-1} = 0$, which gives an isomorphism of graded schemes over $\fr{t}^n \mmod W$:
    $$\fr{t}^n \mmod W \times_{\ld{V}/\GG_a(2-4n)} \fr{t}^n \mmod W \cong \spec \QQ[p_1', p_1, \cdots, p_{n-1}, b]/bp_{n-1}.$$
    On the other hand, \cref{prop: ku Sp2n-2 homology Loops HPn-1} gives an isomorphism of graded $\pi_\ast \H_{\Sp_{2n-2} \times \Sp_2}^\ast(\ast; \QQ)$-algebras
    $$\H^{\Sp_2 \times \Sp_{2n-2}}_\ast(\Omega \HHP^{n-1}; \QQ) \cong \QQ[p_1', p_1, \cdots, p_{n-1}, b]/bp_{n-1},$$
    which implies the first part of the proposition. 
    The second part of the proposition follows from \cref{lem: constructible set quotient} (rather, its obvious variant for graded affine spaces).
\end{proof}
\cref{prop: ku reg centr type Cn} and \cref{thm: ordinary homology criterion satake} imply:
\begin{corollary}\label{cor: ku bzsv for Cn}
    Let $G = \Sp_{2n}$ and $H = \Sp_2 \times \Sp_{2n-2}$, so $\ld{G} = \SO_{2n+1}$ and $\ld{G}_X = \SL_2$.
    There is an equivalence of $\QQ$-linear $\infty$-categories
    $$\Shv_{G\pw{t}}^{c,\Sat}(G\ls{t}/H\ls{t}; \QQ) \simeq \Perf(\sh^{1/2} \ld{M}^\ddag/\SO_{2n+1}).$$
\end{corollary}
\begin{proof}[Proof of \cref{thm: bzsv for Cn}]
    Combining \cref{cor: ku bzsv for Cn} with \cref{prop: defns ShvSat ku and Q}, we see that there is an equivalence of $\QQ$-linear $\infty$-categories
    $$\Shv_{G\pw{t}}^{c,\Sat}(G\ls{t}/H\ls{t}; \QQ) \simeq \Perf(\sh^{1/2} \ld{M}^\ddag/\SO_{2n+1}).$$
    It suffices to describe $\ld{M}^\ddag/\SO_{2n+1}$. Let $W'$ denote the Weyl group of $\Sp_{2n-4}$, so that it acts on $\fr{t}^{n-1}(2)$ such that $\fr{t}^{n-1}(2)\mmod W' \cong \spec \QQ[p_1', p_1, \cdots, p_{n-2}]$.
    Recall that there is a $\ld{G}_X$-equivariant fibration
    $$\ol{\Ind_{\GG_a(2-4n)}^{\ld{G}_X} \ld{V}} \to \ol{\GG_a(2-4n)\backslash \ld{G}_X} \cong V$$
    whose fibers are isomorphic to $\ld{V}$.
    By \cref{lem: weighted cotangent}, this implies that $\ol{\Ind_{\GG_a(2-4n)}^{\ld{G}_X} \ld{V}} \cong T^\ast(4n-4)(V)$, and so
    $$\ld{M}^\ddag \cong \fr{t}^{n-1}(2)\mmod W' \times \Ind_{\ld{G}_X}^{\ld{G}} T^\ast(4n-4)(V);$$
    this implies that $\ld{M}^\ddag/\SO_{2n+1}$ is isomorphic to $\fr{t}^{n-1}(2)\mmod W' \times T^\ast(4n-4)(V)/\ld{G}_X(-(4n-2)\rho_{\ld{G}_X})$, which implies the desired claim.
\end{proof}

%% file: each-type/type-Dn.tex
\subsection{Type $D_n$: $\mathrm{PSO}_{2n}/\SO_{2n-1}$}

Our goal in this section is to prove \cref{thm: rk 1 bzsv is true} in type $D_n$, i.e., for the spherical $\SO_{2n}/\mu_2$-variety $\SO_{2n}/\mu_2\cdot \SO_{2n-1}$.
Write $G = \SO_{2n}/\mu_2$ and $H = \SO_{2n-1}$, so $\ld{G} = \Spin_{2n}$ and $\ld{G}_X = \Spin_3 = \SL_2$. 
As a shorthand, we will write $N = 2n-2$.
\begin{theorem}[\cref{thm: rk 1 bzsv is true} in type $D_n$]\label{thm: bzsv for Dn}
    There is an equivalence of $\QQ$-linear $\infty$-categories
    $$\Shv_{G\pw{t}}^{c,\Sat}(G\ls{t}/H\ls{t}; \QQ) \cong \Perf(\sl_2[N-N\rho_{\ld{G}_X}]/\SL_2[-N\rho_{\ld{G}_X}] \times \fr{spin}_{N-1}[2]\mmod \Spin_{N-1}).$$
\end{theorem}
\begin{example}\label{ex: bzsv D2}
    If $n=2$, \cref{thm: bzsv for Dn} says that there is an equivalence of $\QQ$-linear $\infty$-categories
    $$\Shv_{(\SO_4/\mu_2)\pw{t}}^{c,\Sat}((\SO_4/\mu_2)\ls{t}/\SO_3\ls{t}; \QQ) \cong \Perf(\sl_2[2-2\rho]/\SL_2[-2\rho]).$$
    By \cref{ex: SO4 mod SO3}, the left-hand side can therefore be identified with $\Shv_{(\SO_3 \times \SO_3)\pw{t}}^{c,\Sat}(\SO_3\ls{t}; \QQ)$, and the above equivalence is simply the derived Satake equivalence of \cref{thm: derived satake} for $\SO_3$.
\end{example}
The proof of \cref{thm: bzsv for Dn} will occupy the remainder of this section.
\begin{lemma}\label{lem: SO2n/SO2n-1 and S2n-1}
    There is a homotopy equivalence $(\SO_{2n}/\mu_2\cdot \SO_{2n-1})(\cc) \simeq \RP^{2n-1}$.
    Furthermore, if $B\subseteq G$ is the Borel subgroup of upper-triangular matrices, the Levi quotient $L(X)$ of the parabolic subgroup stabilizing the open $B$-orbit in $\SO_{2n}/\mu_2\cdot \SO_{2n-1}$ is given by $\SO_{2n-2}/\mu_2 \times \GG_m$.
\end{lemma}
\begin{lemma}\label{lem: cohomology Spin2n-1}
    Let $W = (\Z/2)^{n-1} \rtimes \Sigma_{n-1}$ denote the Weyl group of $\SO_{2n-1}$. Then there is an isomorphism of graded $\Z'$-algebras
    $$\H^\ast_{\SO_{2n-1}}(\ast; \Z') \cong \Z'[p_1, \cdots, p_{n-1}] \cong \co_{\fr{t}^{n-1}\mmod W},$$
    where $p_j$ lives in weight $-4j$. The map $\H^\ast_{\SO_{2n-1}}(\ast; \Z') \to \H^\ast_{T^{n-1}}(\ast; \Z')$ sends $p_j$ to the $j$th elementary symmetric polynomial in the variables $x_1^2,\cdots,x_{n-1}^2$.
\end{lemma}
\begin{prop}\label{prop: homology RP2n-1}
    There is an isomorphism of graded $\H^\ast_{\SO_{2n-1}}(\ast; \Z')$-algebras
    $$\H^{\SO_{2n-1}}_\ast(\Omega \RP^{2n-1}; \Z') \cong \Z'[x_1, \cdots, x_{n-1}, a^{\pm 1}, \tfrac{a-a^{-1}}{\prod_{i=1}^{n-1} x_i}]^{W}.$$
    Here, the action of the $j$th $\Z/2 \subseteq (\Z/2)^{n-1} \subseteq W$ sends $x_j \mapsto -{x}_j$ and $a \mapsto a^{-1}$, and the symmetric group acts by permuting the variables $x_1^2, \cdots, x_{n-1}^2$ (and leaves $a$ invariant).
\end{prop}
\begin{proof}
    The restriction of the $\Spin_{2n-1}$-action on $\Spin_{2n}/\Spin_{2n-1}$ (which is homotopy equivalent to $S^{2n-1}$) to $T^{n-1}\subseteq \Spin_{2n-1}$ exhibits $S^{2n-1}$ as the one-point compactification of $\std \oplus \RR$, where $\std$ is the standard $(n-1)$-dimensional complex representation of $T^{n-1}$.
    \cref{cor: G-equiv homology of Loops SV} implies that there is a graded isomorphism of $\H^\ast_{\Spin_{2n-1}}(\ast; \Z')$-algebras
    $$\H^{\Spin_{2n-1}}_\ast(\Omega S^{2n-1}; \Z') \cong \Z'[x_1, \cdots, x_{n-1}, a^{\pm 1}, \tfrac{a-1}{x_1 \cdots x_{n-1}}]^{W}.$$
    The calculation of $\H^{\SO_{2n-1}}_\ast(\Omega \RP^{2n-1}; \Z')$ is a consequence of the above description of $\H^{\Spin_{2n-1}}_\ast(\Omega S^{2n-1}; \Z')$, the fact that $\Spin_{2n-1}/\mu_2 \cong \SO_{2n-1}$, and the fact that $S^{2n-1}/(\Z/2) \cong \RP^{2n-1}$.
\end{proof}
\begin{construction}
    Let $W'$ denote the Weyl group of $\Spin_{N-1}$. Equip $\SL_2$ with the grading coming from $-N\rho$, and consider the $\ld{G}_X = \SL_2$-scheme $\fr{t}^{n-2}(2)\mmod W' \times \sl_2(N-N\rho_{\ld{G}_X})$ (where $\SL_2$ acts only on the factor $\sl_2(N-N\rho_{\ld{G}_X})$, via the adjoint action). Define
    $$\ld{M}^\ddag = \fr{t}^{n-2}(2)\mmod W' \times \Ind_{\ld{G}_X}^{\ld{G}} \sl_2(N-N\rho_{\ld{G}_X}).$$
    Let $\kappa: \fr{t}^{n-1}(2)\mmod W \to \fr{t}^{n-2}(2)\mmod W' \times \sl_2(N-N\rho_{\ld{G}_X})$ denote the closed immersion sending
    $$(p_1, \cdots, p_{n-1}) \mapsto (p_1, \cdots, p_{n-2}), \begin{psmallmatrix}
        0 & 1\\
        p_{n-1} & 0
    \end{psmallmatrix}.$$
    There is a closed immersion 
    \begin{align*}
        \fr{t}^{n-2}(2)\mmod W' \times \sl_2(N-N\rho_{\ld{G}_X}) & \cong \fr{t}^{n-2}(2)\mmod W' \times \Ind_{\ld{G}_X}^{\ld{G}_X} \sl_2(N-N\rho_{\ld{G}_X}) \\
        & \hookrightarrow \fr{t}^{n-2}(2)\mmod W' \times \Ind_{\ld{G}_X}^{\ld{G}} \sl_2(N-N\rho_{\ld{G}_X}) \cong \ld{M}^\ddag,
    \end{align*}
    and hence $\kappa$ defines a closed immersion
    $$\fr{t}^{n-1}(2)\mmod W \xar{\kappa} \fr{t}^{n-2}(2)\mmod W' \times \sl_2(N-N\rho_{\ld{G}_X}) \to \ld{M}^\ddag.$$
    We will denote the above map by $\kappa_{\ld{M}^\ddag}$.
\end{construction}

\begin{prop}\label{prop: ordinary reg centr for Dn}
    The $\ld{G}$-orbit of the image of $\kappa_{\ld{M}^\ddag}$ has complement of codimension $\geq 2$. Moreover, there is an isomorphism of graded group schemes over $\fr{t}^{n-1}(2)\mmod W$:
    $$\spec \H^{\SO_{2n-1}}_\ast(\Omega \RP^{2n-1}; \Z') \cong \fr{t}^{n-1}(2)\mmod W \times_{\ld{M}^\ddag/\ld{G}(-2\rho)} \fr{t}^{n-1}(2)\mmod W.$$
    In particular, the conditions of \cref{thm: ordinary homology criterion satake} hold for the spherical $\SO_{2n}/\mu_2$-variety $\SO_{2n}/\mu_2 \SO_{2n-1}$.
\end{prop}
\begin{proof}
    Let us denote by $Y$ the $\SL_2$-orbit of the image of $\kappa: \fr{t}^{n-1}(2)\mmod W \to \fr{t}^{n-2}(2)\mmod W' \times \sl_2(N-N\rho_{\ld{G}_X})$.
    Then $Y \cong \fr{t}^{n-2}(2)\mmod W' \times \sl_2^\reg(N-N\rho_{\ld{G}_X})$, and it is well-known that the complement of $\sl_2^\reg \subseteq \sl_2$ has complement of codimension $\geq 2$. Applying \cref{lem: some codim stuff} to the $\SL_2$-equivariant inclusion $Y \hookrightarrow \fr{t}^{n-2}(2)\mmod W' \times \sl_2(N-N\rho_{\ld{G}_X})$ and the map $\SL_2 \to \ld{G} = \Spin_{2n}$, we conclude that the $\ld{G}$-orbit of the image of $\kappa_{\ld{M}^\ddag}$ has complement of codimension $\geq 2$.

    To prove the second part of the proposition, \cref{prop: homology RP2n-1} reduces us to showing that there is a graded isomorphism
    \begin{equation}\label{eq: n-shifted reg centr sl2}
        \spec \Z'[x_1, \cdots, x_{n-1}, a^{\pm 1}, \tfrac{a-a^{-1}}{\prod_{i=1}^{n-1} x_i}]^{W} \cong \fr{t}^{n-1}(2)\mmod W \times_{\ld{M}^\ddag/\ld{G}(-2\rho)} \fr{t}^{n-1}(2)\mmod W.
    \end{equation}
    The $j$th copy of $\Z/2 \subseteq (\Z/2)^{n-1} \subseteq W$ sends $x_j \mapsto -x_j$ and $a \mapsto a^{-1}$, so there is an isomorphism
    \begin{align}
        \Z'[x_1, \cdots, x_{n-1}, a^{\pm 1}, \tfrac{a-a^{-1}}{\prod_{i=1}^{n-1} x_i}]^{W} & \cong \Z'[p_1, \cdots, p_{n-1}, a + a^{-1}, \tfrac{a-a^{-1}}{x_1 \cdots x_{n-1}}] \nonumber \\
        & \cong \co_{\fr{t}^{n-2}(2)\mmod W'} \otimes_{\Z'} \Z'[p_{n-1}, a + a^{-1}, \tfrac{a-a^{-1}}{x_1 \cdots x_{n-1}}], \label{eq: Weyl invts reg centr SL2}
    \end{align}
    where we recall that $p_j$ is the $j$th elementary symmetric polynomial in the variables $x_1^2, \cdots, x_{n-1}^2$. In particular, $p_{n-1} = (x_1 \cdots x_{n-1})^2$.

    On the other hand, there is a graded isomorphism
    \begin{align*}
        \hspace*{-1cm}
        \fr{t}^{n-1}(2)\mmod W \times_{\ld{M}^\ddag/\ld{G}(-2\rho)} \fr{t}^{n-1}(2)\mmod W & \cong  \fr{t}^{n-1}(2)\mmod W \times_{(\fr{t}^{n-2}(2)\mmod W' \times \sl_2(N-N\rho_{\ld{G}_X}))/\SL_2(-N\rho)} \fr{t}^{n-1}(2)\mmod W\\
        & \cong \fr{t}^{n-2}(2)\mmod W' \times (\AA^1(N)\mmod (\Z/2) \times_{\sl_2(N-N\rho_{\ld{G}_X})/\SL_2(-N\rho)} \AA^1(N)\mmod (\Z/2)).
    \end{align*}
    By construction, the map $\AA^1(N)\mmod (\Z/2) \to \sl_2(N-N\rho_{\ld{G}_X})/\SL_2(-N\rho)$ is precisely a shifted version of the Kostant slice, so the discussion in \cite[Remark B.4]{grg-reg} implies that if we write $p_{n-1}$ to denote the coordinate on $\AA^1(N)\mmod (\Z/2)$, there is an isomorphism
    \begin{align*}
        \AA^1(N)\mmod (\Z/2) \times_{\sl_2(N-N\rho_{\ld{G}_X})/\SL_2(-N\rho)} \AA^1(N)\mmod (\Z/2) & \cong \spec \Z'[c_{n-1}, a^{\pm 1}, \tfrac{a - a^{-1}}{c_{n-1}}]^{\Z/2} \\
        & \cong \spec \Z'[p_{n-1}, a + a^{-1}, \tfrac{a-a^{-1}}{c_{n-1}}],
    \end{align*}
    where $c_{n-1}^2 = p_{n-1}$. This, along with \cref{eq: Weyl invts reg centr SL2}, implies \cref{eq: n-shifted reg centr sl2}; it is not difficult to observe that \cref{eq: n-shifted reg centr sl2} is in fact an isomorphism of group schemes.
\end{proof}
\begin{remark}\label{ex: SO2n+1 and GLn+1 good}
    Using \cref{prop: homology RP2n-1}, one finds that there are isomorphisms
    \begin{multline*}
        \spec \H^{\GL_n}_\ast(\Omega (\SO_{2n+2}/\SO_{2n+1}); \QQ) \cong \spec \H^{\SO_{2n+1}}_\ast(\Omega (\SO_{2n+2}/\SO_{2n+1}); \QQ) \otimes_{\H^\ast_{\SO_{2n+1}}(\ast; \QQ)} \H^\ast_{\GL_n}(\ast; \QQ) \\
        \cong \AA^1(2n)\mmod (\Z/2) \times_{\sl_2^\ast(2n-2n\rho)/\PGL_2(-2n\rho)} \AA^1(2n) \times \fr{gl}_{n-1}^\ast(2)\mmod \GL_{n-1}\\
        \cong \left( \AA^1(2n) \times_{\tilde{\sl_2}(2n-2n\rho)/\PGL_2(-2n\rho)} \AA^1(2n)\right) \times \fr{gl}_{n-1}^\ast(2)\mmod \GL_{n-1}
    \end{multline*}
    of group schemes over $\spec \H^\ast_{\GL_n}(\ast; \QQ) \cong \AA^1(2n) \times \fr{gl}_{n-1}^\ast(2)\mmod \GL_{n-1}$. Here, the map $\AA^1 \to \tilde{\sl_2}$ is obtained by pulling back the Kostant section $\AA^1 \to \sl_2$ along the Grothendieck-Springer resolution. (See also \cref{rmk: K-equiv-homology} for an interpretation via Hochschild cohomology.)
\end{remark}
\begin{proof}[Proof of \cref{thm: bzsv for Dn}]
    This follows from \cref{prop: ordinary reg centr for Dn} and \cref{thm: ordinary homology criterion satake}, along with the identification between $\ld{M}^\ddag/\ld{G}(-2\rho)$ and $\fr{t}^{n-2}(2)\mmod W' \times \sl_2(N-N\rho_{\ld{G}_X})/\ld{G}_X(-N\rho_{\ld{G}_X})$.
\end{proof}

\begin{remark}\label{rmk: levi sqrt Dn}
    Note that the normalization term $\fr{t}^{n-2}(2)\mmod W'$ identifies with $\fr{so}_{N-1}(2)\mmod \SO_{N-1}$, which is the invariant quotient for the group $L^\wedge_X$ from \cite{knop-schalke}.
\end{remark}

%% file: each-type/type-F4.tex
\subsection{Type $\F_4$: $\F_4/\Spin_9$}

Our goal in this section is to prove \cref{thm: rk 1 bzsv is true} in type $\F_4$, i.e., for the spherical $\F_4$-variety $\F_4/\Spin_9$. Let $G = \F_4$ and $H = \Spin_9$, so that $\ld{G} = \F_4$. Recall from \cref{table: topology and dualities for rank 1 spherical varieties} that $\ld{G}_X = \SL_2$.
\begin{theorem}[\cref{thm: rk 1 bzsv is true} in type $\F_4$]\label{thm: bzsv for F4}
    There is an equivalence of $\QQ$-linear $\infty$-categories
    $$\Shv_{G\pw{t}}^{c,\Sat}(G\ls{t}/H\ls{t}; \QQ) \simeq \Perf(T^\ast[16](\AA^2[22,0])/\SL_2[-22\rho_{\ld{G}_X}] \times \fr{sp}_6^\ast[2]\mmod \Sp_6).$$
\end{theorem}

The proof of \cref{thm: bzsv for F4} will take up the remainder of this section.

\begin{lemma}
    There is a homotopy equivalence $(\F_4/\Spin_9)(\cc) \simeq \OP^2$. Moreover, if $B\subseteq G$ is the Borel subgroup of upper-triangular matrices, the Levi quotient $L(X)$ of the parabolic subgroup stabilizing the open $B$-orbit in $\F_4/\Spin_9$ is given by $\Sp_6 \times \GG_m$.
\end{lemma}
\begin{prop}\label{prop: ku F4 homology OP2}
    Let $W = (\Z/2)^4 \rtimes \Sigma_4$ denote the Weyl group of $\Spin_9$.
    There is an isomorphism of graded $\H_{\Spin_9}^\ast(\ast; \QQ)$-algebras
    $$\H^{\Spin_9}_\ast(\Omega \OP^2; \QQ) \cong \QQ[p_1, \cdots, p_4, b]/b p_4,$$
    where $b$ lives in weight $22$.
\end{prop}
\begin{proof}
    The argument is essentially the same as that of the preceding subsections. Let $T^4 \subseteq \Spin_9$ denote the maximal torus; we will begin by describing $\H^{T^4}_\ast(\Omega \OP^2; \QQ)$. The homotopy fixed points spectral sequence for $\pi_\ast \QQ[\Omega \OP^2]^{hT^4}$ is given by
    $$E_2^{\ast,\ast} \cong \H_\ast(\Omega \OP^2; \QQ) \otimes_{\QQ} \pi_\ast \H_{T^4}^\ast(\ast; \QQ) \Rightarrow \pi_\ast \QQ[\Omega \OP^2]^{hT^4}.$$
    To compute the $E_2$-page, we need to compute $\H_\ast(\Omega \OP^2; \QQ)$. Using the Atiyah-Hirzebruch spectral sequence, we will first calculate $\H_\ast(\Omega \OP^2; \QQ)$. Although there is no Hopf fibration $S^7 \to S^{23} \to \OP^2$ (otherwise, the cofiber of the map $S^{23} \to \OP^2$ would provide a contradiction to the Hopf invariant one problem), we can instead compute $\H_\ast(\Omega \OP^2; \QQ)$ using the Serre spectral sequence for the fibration
    $$\Omega \OP^2 \to \ast \to \OP^2$$
    and the fact that $\H_\ast(\OP^2; \QQ)$ is isomorphic to a free graded $\QQ$-module on classes $\{1, x_1, x_2\}$ in weights $0$, $8$, and $16$.
    This is a standard argument: one finds that $\H_\ast(\Omega \OP^2; \QQ) \cong \QQ[a, b]/a^2$ where $a$ lives in weight $7$ and $b$ lives in weight $22$; the differentials in the Serre spectral sequence are given by
    $$d^8(b^j x_1) = ab^j, \ d^8(b^j x_2) = ab^j x_1, \ d^{23}(a b^j x_2) = b^{j+1}.$$
    The Atiyah-Hirzebruch spectral sequence for $\H_\ast(\Omega \OP^2; \QQ)$ degenerates at the $E_1$-page (with no multiplicative extensions), and we obtain an isomorphism $\H_\ast(\Omega \OP^2; \QQ) \cong \QQ[a, b]/a^2$.
    Returning to the homotopy fixed points spectral sequence, the above discussion implies that
    $$E_2^{\ast,\ast} \cong \QQ[a, b]\pw{x_1, \cdots, x_4}/a^2.$$
    There is a single differential
    $$d_2(a) = b x_1 \ol{x_1} \cdots x_4 \ol{x_4},$$
    and the spectral sequence is concentrated in even degrees after running this differential. It therefore collapses on the $E_3$-page, and we find that there is an isomorphism 
    $$\pi_\ast \QQ[\Omega \OP^2]^{hT^4} \cong \QQ[b]\pw{x_1, \cdots, x_4}/b x_1 \ol{x_1} \cdots x_4 \ol{x_4}.$$
    To calculate $\H^{T^4}_\ast(\Omega \OP^2; \QQ)$ itself (and not just its completion $\pi_\ast \QQ[\Omega \OP^2]^{hT^4}$), the strategy of \cref{thm: homology of loops SV} reduces us to showing that the restriction $\H^{T^4}_\ast(\Omega \OP^2; \QQ)|_{\punc{\fr{t}}} = \pi_\ast \H_{T^4}^\ast(\ast; \QQ)$. By \cref{lem: atiyah localization}, there is an isomorphism $\H^{T^4}_\ast(\Omega \OP^2; \QQ)|_{\punc{\fr{t}}} \cong \H^{T^4}_\ast(\Omega (\OP^2)^{T^4}; \QQ)|_{\punc{\fr{t}}}$. It therefore suffices to show that $\Omega (\OP^2)^{T^4}$ is contractible, but this is a consequence of the simple observation that $(\OP^2)^{T^4} \cong S^0$. This discussion gives an isomorphism of graded $\pi_\ast \H_{T^4}^\ast(\ast; \QQ)$-algebras
    $$\H^{T^4}_\ast(\Omega \OP^2; \QQ) \cong \QQ[x_1, \cdots, x_4, b]/bx_1 \ol{x_1} \cdots x_{4} \ol{x_4}.$$
    This isomorphism is $W$-equivariant (where $W$ is the Weyl group of $\Spin_9$), so \cref{prop: Weyl invts} implies that there is an isomorphism of graded $\pi_\ast \H_{\Spin_9}^\ast(\ast; \QQ)$-algebras
    $$\H^{\Spin_9}_\ast(\Omega \OP^2; \QQ) \cong \left(\QQ[x_1, \cdots, x_{4}, b]/bx_1 \ol{x}_1 \cdots x_4 \ol{x_4}\right)^W.$$
    Noting that the action of $W$ leaves $b$ invariant, \cref{lem: cohomology Spin2n-1} computes the right-hand side; the resulting answer is precisely the right-hand side of the proposition.
\end{proof}
\begin{construction}
    Let $\ld{V}$ denote the graded affine scheme $\AA^1(-6) \times \fr{t}^4\mmod W$, and let $\kappa: \fr{t}^4\mmod W \to \ld{V}$ denote the map sending $\vec{p} \mapsto (0, \vec{p})$. There is an action of $\GG_a(-22)$ on $\ld{V}$, where $b\in \GG_a(-22)$ sends
    $$(x, \vec{p}) \mapsto (x - bp_4, \vec{p}).$$
    Equip $\ld{G}_X = \SL_2$ with the grading coming from $22\rho_{\SL_2}$, so that the entries of a matrix $\begin{psmallmatrix}
        a & b\\
        c & d
    \end{psmallmatrix}$ are equipped with the following weights: $a$ and $d$ have weight $0$, $b$ has weight $22$, and $c$ has weight $-22$. Let $\ld{V}$ denote the affine space $\AA^2(22, 0)$, so that so that there is an action of $\ld{G}_X = \SL_2$ on $\ld{V}$ via
    $$\begin{psmallmatrix}
        a & b\\
        c & d
    \end{psmallmatrix} \cdot (x,y) = (ax + cy, bx + dy);$$
    here, $x$ lives in degree $-22$ and $y$ lives in degree $0$.
    There is an isomorphism
    $$\GG_a(-22)\backslash \ld{G}_X \cong \ld{V} - \{0\}, \ \begin{psmallmatrix}
        a & b\\
        c & d
    \end{psmallmatrix} \mapsto (c,d),$$
    and the above action of $\ld{G}_X$ on $\ld{V}$ restricts on $\ld{V} - \{0\}$ to the right-action of $\ld{G}_X$ on $\GG_a(-22)\backslash \ld{G}_X$.
    There is a fibration
    $$\Ind_{\GG_a(-22)}^{\ld{G}_X} \ld{V} \to \GG_a(-22)\backslash \ld{G}_X \cong \ld{V} - \{0\}$$
    whose fibers are isomorphic to $\ld{V}$.
    Let $\ol{\Ind_{\GG_a(-22)}^{\ld{G}_X} \ld{V}}$ denote the affine closure of $\Ind_{\GG_a(-22)}^{\ld{G}_X} \ld{V}$, so that there is a $\ld{G}_X$-equivariant fibration
    $$\ol{\Ind_{\GG_a(-22)}^{\ld{G}_X} \ld{V}} \to \ol{\GG_a(-22)\backslash \ld{G}_X} \cong V$$
    whose fibers are isomorphic to $\ld{V}$. Let $\ld{M}^\ddag$ denote the induction
    $$\ld{M}^\ddag = \Ind_{\ld{G}_X}^{\ld{G}} \ol{\Ind_{\GG_a(-22)}^{\ld{G}_X} \ld{V}}.$$
    \cref{lem: some codim stuff} implies:
    \begin{lemma}
        There is an open immersion
        $$\Ind_{\GG_a(-22)}^{\ld{G}} \ld{V} \to \Ind_{\ld{G}_X}^{\ld{G}} \ol{\Ind_{\GG_a(-22)}^{\ld{G}_X} \ld{V}} = \ld{M}^\ddag$$
        which exhibits the target as the affine closure of the source, and whose complement is of codimension $\geq 2$.
    \end{lemma}
    The map $\kappa: \fr{t}^4\mmod W \to \ld{V}$ defines a locally closed immersion
    $$\fr{t}^4\mmod W \xar{\kappa} \ld{V} \hookrightarrow \Ind_{\GG_a(-22)}^{\ld{G}} \ld{V} \hookrightarrow \Ind_{\ld{G}_X}^{\ld{G}} \ol{\Ind_{\GG_a(-22)}^{\ld{G}_X} \ld{V}},$$
    which we will denote by $\kappa_{\ld{M}^\ddag}$.
\end{construction}
\begin{prop}\label{prop: ku reg centr type F4}
    There is an isomorphism of graded group schemes over $\fr{t}^4\mmod W$:
    $$\spec \H^{\Spin_9}_\ast(\Omega \OP^2; \QQ) \cong \fr{t}^4\mmod W \times_{\ld{M}^\ddag/\ld{G}(-2\rho)} \fr{t}^4\mmod W.$$
    Moreover, if $\ld{J}_X$ denotes the above group scheme over $\fr{t}^4\mmod W$, the algebra of regular functions on $(\fr{t}^4\mmod W \times \ld{G})/\ld{J}_X$ is isomorphic to $\co_{\ld{M}^\ddag}$.
\end{prop}
\begin{proof}
    There is an isomorphism 
    $$\fr{t}^4\mmod W \times_{\ld{M}^\ddag/\ld{G}(-2\rho)} \fr{t}^4\mmod W \cong \fr{t}^4\mmod W \times_{\ld{V}/\GG_a(-22)} \fr{t}^4\mmod W,$$
    as well as a closed immersion
    $$\fr{t}^4\mmod W \times_{\ld{V}/\GG_a(-22)} \fr{t}^4\mmod W \hookrightarrow \fr{t}^4\mmod W \times \GG_a(-22),$$
    which exhibits $\fr{t}^4\mmod W \times_{\ld{V}/\GG_a(-22)} \fr{t}^4\mmod W$ as the subscheme of tuples $(\vec{p}, b)$ such that $b$ stabilizes $\kappa(\vec{p})$. By definition of $\kappa$, this happens if and only if $bp_4 = 0$, which gives an isomorphism of graded schemes over $\fr{t}^4\mmod W$:
    $$\fr{t}^4\mmod W \times_{\ld{V}/\GG_a(-22)} \fr{t}^4\mmod W \cong \spec \QQ[p_1, \cdots, p_4, b]/bp_4.$$
    On the other hand, \cref{prop: ku F4 homology OP2} gives an isomorphism of graded $\pi_\ast \ku_{\Spin_9}$-algebras
    $$\ku^{\Spin_9}_\ast(\Omega \OP^2) \cong \QQ[p_1, \cdots, p_4, b]/bp_4,$$
    which implies the first part of the proposition. 
    The second part of the proposition follows from \cref{lem: constructible set quotient} (rather, its obvious variant for graded affine spaces).
\end{proof}
\cref{prop: ku reg centr type F4} and \cref{prop: ku rk 1 reg centr --> thm} imply:
\begin{corollary}\label{cor: ku bzsv for F4}
    Let $G = \F_4$ and $H = \Spin_9$, so $\ld{G} = \F_4$ and $\ld{G}_X = \SL_2$.
    There is an equivalence of $\QQ$-linear $\infty$-categories
    $$\Shv_{G\pw{t}}^{c,\Sat}(G\ls{t}/H\ls{t}; \QQ) \simeq \Perf(\sh^{1/2} \ld{M}^\ddag/\F_4).$$
\end{corollary}
\begin{proof}[Proof of \cref{thm: bzsv for F4}]
    Combining \cref{cor: ku bzsv for F4} with \cref{prop: defns ShvSat ku and Q}, we see that if $\ld{M}^\ddag = \ld{M}^\ddag/\beta$, there is an equivalence of $\QQ$-linear $\infty$-categories
    $$\Shv_{G\pw{t}}^{c,\Sat}(G\ls{t}/H\ls{t}; \QQ) \simeq \Perf(\sh^{1/2} \ld{M}^\ddag/\F_4).$$
    It suffices to compute $\ld{M}^\ddag/\F_4$. 
    Let $W'$ denote the Weyl group of $\Spin_7$, so that $W' \cong (\Z/2)^3 \rtimes \Sigma_3$ acts on $\fr{t}^3(2)$ such that $\fr{t}^3(2)\mmod W' \cong \spec \QQ[p_1, p_2, p_3]$. Recall that there is a $\ld{G}_X$-equivariant fibration
    $$\ol{\Ind_{\GG_a(-22)}^{\ld{G}_X} \ld{V}} \to \ol{\GG_a(-22)\backslash \ld{G}_X} \cong V$$
    whose fibers are isomorphic to $\ld{V}$.
    By \cref{lem: weighted cotangent}, this implies that $\ol{\Ind_{\GG_a(-22)}^{\ld{G}_X} \ld{V}} \cong T^\ast(16)(V)$, and so the isomorphism between $\ld{M}^\ddag/\F_4$ and $\fr{t}^3(2)\mmod W' \times T^\ast(16)(V)/\ld{G}_X(-22\rho_{\ld{G}_X})$ gives the desired claim.
\end{proof}

%% file: each-type/type-G2.tex
\subsection{Type $\G_2$: $\G_2/\SL_3$}

Our goal in this section is to prove \cref{thm: rk 1 bzsv is true} in type $\G_2$, i.e., for the spherical $\G_2$-variety $\G_2/\SL_3$. Let $G = \G_2$ and $H = \SL_3$, so that $\ld{G} = \G_2$. Recall from \cref{table: topology and dualities for rank 1 spherical varieties} that $\ld{G}_X = \SL_2$.
\begin{theorem}[\cref{thm: rk 1 bzsv is true} in type $\G_2$]\label{thm: bzsv for G2}
    There is an equivalence of $\QQ$-linear $\infty$-categories
    $$\Shv_{G\pw{t}}^{c,\Sat}(G\ls{t}/H\ls{t}; \QQ) \simeq \Perf(T^\ast[6](\AA^2[10,0])/\SL_2(-10\rho_{\ld{G}_X}) \times \fr{sl}_2[2]\mmod \SL_2).$$
\end{theorem}

The proof of \cref{thm: bzsv for G2} will take up the remainder of this section; in fact, we will prove a $\ku$-theoretic deformation.

\begin{lemma}
    Because $\g_2/\sl_3 \cong T^\ast(\std_3)$ as an $\SL_3$-representation, there is a homotopy equivalence $(\G_2/\SL_3)(\cc) \simeq S^6$. Moreover, if $B\subseteq G$ is the Borel subgroup of upper-triangular matrices, the Levi quotient $L(X)$ of the parabolic subgroup stabilizing the open $B$-orbit in $\G_2/\SL_3$ is given by $\GL_2$.
\end{lemma}
\begin{lemma}\label{lem: SUn equiv coh of point}
    Let $W = \Sigma_n$ denote the Weyl group of $\SL_n$. Then there is an isomorphism
    $$\pi_\ast \ku_{\SU(n)} \cong \Z'[\beta, c_2, \cdots, c_n],$$
    where the injective map $\pi_\ast \ku_{\SU(n)} \to \pi_\ast \ku_{T^{n-1}}$ sends $c_j$ to the $j$th elementary symmetric polynomial in the variables $x_1, \cdots, x_{n-1}, x_n$, where $x_n$ is the sum of $x_1, \cdots, x_{n-1}$ in the group structure on $\GG_\beta$.
\end{lemma}
\begin{prop}\label{prop: ku homology Loops S6}
    There is an isomorphism of graded $\pi_\ast \ku_{\SU(3)}$-algebras
    $$\ku^{\SU(3)}_\ast(\Omega S^6) \cong \Z'[\beta, c_2, c_3, b]/bc_3,$$
    where $b$ is in weight $10$.
\end{prop}
\begin{proof}
    The restriction of the $\SU(3)$-action on $S^6$ to $T^2\subseteq \SU(3)$ exhibits $S^6$ as the one-point compactification of the $T^2$-representation with weights $\lambda_1$, $\lambda_2$, and $\lambda_1 + \lambda_2$. Therefore, \cref{cor: G-equiv homology of Loops SV} implies that
    $$\ku^{\SU(3)}_\ast(\Omega S^6) \cong \left(\Z'[\beta, x_1, x_2, \tfrac{1}{(1+\beta x_1)(1+\beta x_2)}, b]/b x_1 x_2 x_3\right)^{\Sigma_3}.$$
    The action of $\Sigma_3$ permutes $x_1$, $x_2$, and $x_3$, and leaves $b$ invariant. Therefore, \cref{lem: SUn equiv coh of point} implies that this ring of invariants can be identified with $\Z'[\beta, c_2, c_3, b]/bc_3$, as desired.
\end{proof}
\begin{definition}
    Let $\ld{V}_\beta$ denote the graded affine scheme $\AA^1(-4) \times T^2_\beta\mmod \Sigma_3$, and let $\kappa: T^2_\beta\mmod \Sigma_3 \to \ld{V}_\beta$ denote the map sending $\vec{c} \mapsto (0, \vec{c})$. There is an action of $\GG_a(-10)$ on $\ld{V}_\beta$, where $b\in \GG_a(-10)$ sends
    $$(z, \vec{c}) \mapsto (z - bc_3, \vec{c}).$$
    Equip $\ld{G}_X = \SL_2$ with the grading coming from $10\rho_{\SL_2}$, so that the entries of a matrix $\begin{psmallmatrix}
        a & b\\
        c & d
    \end{psmallmatrix}$ are equipped with the following weights: $a$ and $d$ have weight $0$, $b$ has weight $10$, and $c$ has weight $-10$. Let $\ld{V}$ denote the affine space $\AA^2(10, 0)$, so that there is an action of $\ld{G}_X = \SL_2$ on $\ld{V}$ via
    $$\begin{psmallmatrix}
        a & b\\
        c & d
    \end{psmallmatrix} \cdot (x,y) = (ax + cy, bx + dy);$$
    here, $x$ lives in degree $-10$ and $y$ lives in degree $0$.
    There is an isomorphism
    $$\GG_a(-10)\backslash \ld{G}_X \cong V - \{0\}, \ \begin{psmallmatrix}
        a & b\\
        c & d
    \end{psmallmatrix} \mapsto (c,d),$$
    and the above action of $\ld{G}_X$ on $\ld{V}$ restricts on $\ld{V} - \{0\}$ to the right-action of $\ld{G}_X$ on $\GG_a(-10)\backslash \ld{G}_X$.
    There is a fibration
    $$\Ind_{\GG_a(-10)}^{\ld{G}_X} \ld{V}_\beta \to \GG_a(-10)\backslash \ld{G}_X \cong \ld{V} - \{0\}$$
    whose fibers are isomorphic to $\ld{V}_\beta$. Let $\ol{\Ind_{\GG_a(-10)}^{\ld{G}_X} \ld{V}_\beta}$ denote the affine closure of $\Ind_{\GG_a(-10)}^{\ld{G}_X} \ld{V}_\beta$, so that there is a $\ld{G}_X$-equivariant fibration
    $$\ol{\Ind_{\GG_a(-10)}^{\ld{G}_X} \ld{V}_\beta} \to \ol{\GG_a(-10)\backslash \ld{G}_X} \cong \ld{V}$$
    whose fibers are isomorphic to $\ld{V}_\beta$. Let $\ld{M}^\ddag_\beta$ denote the induction
    $$\ld{M}^\ddag_\beta = \Ind_{\ld{G}_X}^{\ld{G}} \ol{\Ind_{\GG_a(-10)}^{\ld{G}_X} \ld{V}_\beta}.$$
    \cref{lem: some codim stuff} implies:
    \begin{lemma}
        There is an open immersion
        $$\Ind_{\GG_a(-10)}^{\ld{G}} \ld{V}_\beta \to \Ind_{\ld{G}_X}^{\ld{G}} \ol{\Ind_{\GG_a(-10)}^{\ld{G}_X} \ld{V}_\beta} = \ld{M}^\ddag_\beta$$
        which exhibits the target as the affine closure of the source, and whose complement is of codimension $\geq 2$.
    \end{lemma}
    The map $\kappa: T^2_\beta\mmod \Sigma_3 \to \ld{V}_\beta$ defines a locally closed immersion
    $$T^2_\beta\mmod \Sigma_3 \xar{\kappa} \ld{V}_\beta \hookrightarrow \Ind_{\GG_a(-10)}^{\ld{G}} \ld{V}_\beta \hookrightarrow \Ind_{\ld{G}_X}^{\ld{G}} \ol{\Ind_{\GG_a(-10)}^{\ld{G}_X} \ld{V}_\beta},$$
    which we will denote by $\kappa_{\ld{M}^\ddag_\beta}$.
\end{definition}
\begin{prop}\label{prop: ku reg centr type G2}
    There is an isomorphism of graded group schemes over $T^2_\beta\mmod \Sigma_3$:
    $$\spec \ku^{\SU(3)}_\ast(\Omega S^6) \cong T^2_\beta\mmod \Sigma_3 \times_{\ld{M}^\ddag_\beta/\ld{G}(-2\rho)} T^2_\beta\mmod \Sigma_3.$$
    Moreover, if $\ld{J}_{X, \beta}$ denotes the above group scheme over $T^2_\beta\mmod \Sigma_3$, the algebra of regular functions on $(T^2_\beta\mmod \Sigma_3 \times \ld{G})/\ld{J}_{X, \beta}$ is isomorphic to $\co_{\ld{M}^\ddag_\beta}$.
\end{prop}
\begin{proof}
    There is an isomorphism 
    $$T^2_\beta\mmod \Sigma_3 \times_{\ld{M}^\ddag_\beta/\ld{G}(-2\rho)} T^2_\beta\mmod \Sigma_3 \cong T^2_\beta\mmod \Sigma_3 \times_{\ld{V}_\beta/\GG_a(-10)} T^2_\beta\mmod \Sigma_3,$$
    as well as a closed immersion
    $$T^2_\beta\mmod \Sigma_3 \times_{\ld{V}_\beta/\GG_a(-10)} T^2_\beta\mmod \Sigma_3 \hookrightarrow T^2_\beta\mmod \Sigma_3 \times \GG_a(-10),$$
    which exhibits $T^2_\beta\mmod \Sigma_3 \times_{\ld{V}_\beta/\GG_a(-10)} T^2_\beta\mmod \Sigma_3$ as the subscheme of pairs $(\vec{c}, b)$ such that $b$ stabilizes $\kappa(\vec{c})$. By definition of $\kappa$, this happens if and only if $bc_3 = 0$, which gives an isomorphism of graded schemes over $T^2_\beta\mmod \Sigma_3$:
    $$T^2_\beta\mmod \Sigma_3 \times_{\ld{V}_\beta/\GG_a(-10)} T^2_\beta\mmod \Sigma_3 \cong \spec \Z'[\beta, c_2, c_3, b]/bc_3.$$
    \cref{prop: ku homology Loops S6} therefore implies the first part of the proposition. 
    The second part of the proposition follows from \cref{lem: constructible set quotient} (rather, its obvious variant for graded affine spaces).
\end{proof}
\cref{prop: ku reg centr type G2} and \cref{prop: ku rk 1 reg centr --> thm} imply:
\begin{corollary}\label{cor: ku bzsv for G2}
    Let $G = \G_2$ and $H = \SL_3$, so $\ld{G} = \G_2$ and $\ld{G}_X = \SL_2$.
    There is an equivalence of $\sh^{1/2}(\Z'[\beta])$-linear $\infty$-categories
    $$\Shv_{G\pw{t}}^{c,\Sat}(G\ls{t}/H\ls{t}; \ku)^\faux \simeq \Perf(\sh^{1/2} \ld{M}^\ddag_\beta/\G_2(-2\rho)).$$
\end{corollary}
\begin{proof}[Proof of \cref{thm: bzsv for G2}]
    Combining \cref{cor: ku bzsv for G2} with \cref{prop: defns ShvSat ku and Q}, we see that if $\ld{M}^\ddag = \ld{M}^\ddag_\beta/\beta$, there is an equivalence of $\QQ$-linear $\infty$-categories
    $$\Shv_{G\pw{t}}^{c,\Sat}(G\ls{t}/H\ls{t}; \QQ) \simeq \Perf(\sh^{1/2} \ld{M}^\ddag/\G_2(-2\rho)).$$
    It suffices to compute $\ld{M}^\ddag/\G_2(-2\rho)$. Let $W' = \Z/2$ denote the Weyl group of $\SL_2$, so that it acts on $\fr{t}^1(2)$ such that $\fr{t}^1(2)\mmod W' \cong \AA^1(4)$.
    Recall that there is a $\ld{G}_X$-equivariant fibration
    $$\ol{\Ind_{\GG_a(-10)}^{\ld{G}_X} \ld{V}_\beta} \to \ol{\GG_a(-10)\backslash \ld{G}_X} \cong V$$
    whose fibers are isomorphic to $\ld{V}_\beta$, which implies (by setting $\beta = 0$) that there is a $\ld{G}_X$-equivariant fibration
    $$\ol{\Ind_{\GG_a(-10)}^{\ld{G}_X} \ld{V}} \to \ol{\GG_a(-10)\backslash \ld{G}_X} \cong V$$
    whose fibers are isomorphic to $\ld{V}$.
    By \cref{lem: weighted cotangent}, this implies that $\ol{\Ind_{\GG_a(-10)}^{\ld{G}_X} \ld{V}} \cong T^\ast(6)(V)$, and so the isomorphism between $\ld{M}^\ddag/\G_2(-2\rho)$ and $\fr{t}^1(2)\mmod W' \times T^\ast(6)(V)/\ld{G}_X(-10\rho_{\ld{G}_X})$ implies the desired result.
\end{proof}
\begin{remark}\label{rmk: G2 mod normalizer of SL3}
    \cref{thm: bzsv for G2} can be used to prove a variant for the rank one spherical $\G_2$-variety $\G_2/\N_{\G_2}(\SL_3)$, where $\N_{\G_2}(\SL_3)$ is the normalizer of $\SL_3 \subseteq \G_2$. The dual group $\ld{G}_X$ in this case is again $\SL_2$. The quotient $\N_{\G_2}(\SL_3)/\SL_3$ is isomorphic to $\Z/2$, with a generator being given by the longest element of the Weyl group of $\G_2$. Moreover, $\N_{\G_2}(\SL_3)/\SL_3 \cong \Z/2$ acts on $(\G_2/\SL_3)(\cc) \simeq S^6$ via the antipodal action. As in \cref{rmk: type N PGL2 mod PO2}, this equips the left-hand side of the equivalence f \cref{thm: bzsv for G2} with a natural $\Z/2$-action. One can show that under this equivalence, the resulting $\Z/2$-action on the right-hand side of \cref{thm: bzsv for G2} identifies with the natural $\Z/2$-action on $T^\ast(\AA^2)$ via the symplectic form. Using this, one finds that there is an equivalence
    $$\Shv_{\G_2\pw{t}}^{c,\Sat}(\G_2\ls{t}/\N_{\G_2}(\SL_3)\ls{t}; \QQ) \simeq \Perf(T^\ast[6](\AA^2[10])/(\SL_2[-10\rho_{\ld{G}_X}] \times \Z/2) \times \sl_2[2]\mmod \SL_2).$$
    The spherical root of $\G_2/\N_{\G_2}(\SL_3)$ is of type N, and so it is excluded by \cite{sakellaridis-venkatesh, bzsv}; nevertheless, the preceding equivalence shows that it does admit a dual, given by the \textit{stack} $\Ind_{\SL_2}^{\G_2} (T^\ast(\AA^2)/(\Z/2) \times \AA^1)$. 
\end{remark}

%% file: each-type/type-B3prime.tex
\subsection{Type $B_3'$: $\SO_7/\G_2$}

We will only work with coefficients in $\Z'$ below. 
Our goal in this section is to prove \cref{thm: rk 1 bzsv is true} in type $B_3'$, i.e., for the spherical $\SO_{7}$-variety $\SO_{7}/\G_2$. Write $G = \SO_7$ and $H = \G_2$, so that $\ld{G} = \Sp_{6}$. Recall from \cref{table: topology and dualities for rank 1 spherical varieties} that $\ld{G}_X = \SL_2$.
\begin{theorem}[\cref{thm: rk 1 bzsv is true} in type $B_3'$]\label{thm: bzsv for B3'}
    There is an equivalence of $\QQ$-linear $\infty$-categories
    $$\Shv_{G\pw{t}}^{c,\Sat}(G\ls{t}/H\ls{t}; \QQ) \cong \Perf(\sl_2[6-6\rho_{\ld{G}_X}]/\SL_2[-6\rho_{\ld{G}_X}] \times \fr{sp}_2[2]\mmod \Sp_2).$$
\end{theorem}

\begin{lemma}\label{lem: levi stab B3'}
    There is a homotopy equivalence $(\Spin_7/\G_2)(\cc) \cong S^7$, which implies that $(\SO_7/\G_2)(\cc) \cong \RP^7$. Moreover, if $B\subseteq G$ is the Borel subgroup of upper-triangular matrices, the Levi quotient $L(X)$ of the parabolic subgroup stabilizing the open $B$-orbit in $\SO_{7}/\G_2$ is given by $\GL_3$.
\end{lemma}
\begin{lemma}\label{lem: G2 equiv coh of point}
    Let $W = \Sigma_3 \times \Z/2$ denote the Weyl group of $\G_2$. Then there is an isomorphism
    $$\H^\ast_{\G_2}(\ast; \Z') \cong \Z'[c_2, c_6],$$
    where the injective map $\H^\ast_{\G_2}(\ast; \Z') \to \H^\ast_{T^{2}}(\ast; \Z')$ sends
    \begin{align*}
        c_2 & \mapsto x_1^2 + x_2^2 + x_1 x_2, \\
        c_6 & \mapsto x_1^2 x_2^2 (x_1 + x_2)^2.
    \end{align*}
\end{lemma}
\begin{proof}
The action of $\Sigma_3 \subseteq W$ on $\H^\ast_{T^{2}}(\ast; \Z') = \Z'[x_1, x_2]$ is given by the reduced standard representation (i.e., $x_1$, $x_2$, and $-(x_1 + x_2)$ are permuted), and the action of $\Z/2\subseteq W$ is given by negating the $x_i$ simultaneously. It follows that 
$$\H^\ast_{\G_2}(\ast; \Z') = \H^\ast_{T^{2}}(\ast; \Z')^W \cong (\Z'[x_1, x_2]^{\Sigma_3})^{\Z/2} \cong \Z'[x_1^2 + x_2^2 + x_1 x_2, x_1 x_2 (x_1 + x_2)]^{\Z/2},$$
which is precisely $\Z'[c_2, c_6]$, as desired.
\end{proof}

\begin{prop}\label{prop: homology RP7}
    There is an isomorphism of graded $\pi_\ast \H^\ast_{\G_2}(\ast; \Z')$-algebras
    $$\H^{\G_2}_\ast(\Omega \RP^{7}; \Z') \cong \Z'[x_1, x_2, a^{\pm 1}, \tfrac{a-a^{-1}}{x_1 x_2 x_3}]^{W},$$
    where $x_3 = x_1 + x_2$ is the sum of the two weights $x_1$ and $x_2$. 
    Here, the action of $\Z/2 \subseteq W$ sends $x_j \mapsto -{x}_j$ and $a \mapsto a^{-1}$, and the symmetric group acts by permuting the variables $x_1, x_2, x_3$ (and leaves $a$ invariant).
\end{prop}
\begin{proof}
    The restriction of the $\G_2$-action on $(\Spin_7/\G_2)(\cc) \simeq S^7$ exhibits it as the one-point compactification of $\std \oplus \RR$, where $\std$ is the standard $3$-dimensional complex representation of $\G_2$. In particular, the restriction of the $\G_2$-action on $S^7$ to $T^2\subseteq \G_2$ exhibits it as the one-point compactification of the $3$-dimensional representation of $T^2$ with weights $\lambda_1$, $\lambda_2$, and $\lambda_1 + \lambda_2$. Given this observation, the isomorphism of the proposition is a consequence of \cref{cor: G-equiv homology of Loops SV} (just as with \cref{prop: homology RP2n-1}).
\end{proof}

\begin{construction}
    Let $W'$ denote the Weyl group of $\Sp_2$, so that $\fr{t}^1(2)\mmod W' \cong \AA^1(4)$. Equip $\ld{G}_X = \SL_2$ with the grading by $-6\rho$, and consider the $\ld{G}_X = \SL_2$-scheme $\fr{t}^1(2)\mmod W' \times \sl_2(6-6\rho_{\ld{G}_X})$ (where $\SL_2$ acts only on the factor $\sl_2(6-6\rho_{\ld{G}_X})$ by the adjoint action). Define
    $$\ld{M}^\ddag = \fr{t}^1(2)\mmod W' \times \Ind_{\ld{G}_X}^{\ld{G}} \sl_2(6-6\rho_{\ld{G}_X}).$$
    Let $\kappa: \fr{t}^2(2)\mmod W \to \fr{t}^1(2)\mmod W' \times \sl_2(6-6\rho_{\ld{G}_X})$ denote the closed immersion sending
    $$(c_2, c_6) \mapsto c_2, \begin{psmallmatrix}
        0 & 1\\
        c_6 & 0
    \end{psmallmatrix}.$$
    There is a closed immersion
    \begin{align*}
        \fr{t}^1(2)\mmod W' \times \sl_2(6-6\rho_{\ld{G}_X}) & \cong \fr{t}^1(2)\mmod W' \times \Ind_{\ld{G}_X}^{\ld{G}_X} \sl_2(6-6\rho_{\ld{G}_X}) \\
        & \hookrightarrow \fr{t}^1(2)\mmod W' \times \Ind_{\ld{G}_X}^{\ld{G}} \sl_2(6-6\rho_{\ld{G}_X}) \cong \ld{M}^\ddag,
    \end{align*}
    and hence $\kappa$ defines a closed immersion
    $$\fr{t}^{2}(2)\mmod W \xar{\kappa} \fr{t}^1(2)\mmod W' \times \sl_2(6-6\rho_{\ld{G}_X}) \to \ld{M}^\ddag.$$
    We will denote the above map by $\kappa_{\ld{M}^\ddag}$.
\end{construction}

\begin{prop}\label{prop: ordinary reg centr for B3'}
    The $\ld{G}$-orbit of the image of $\kappa_{\ld{M}^\ddag}$ has complement of codimension $\geq 2$. Moreover, there is an isomorphism of graded group schemes over $\fr{t}^{2}(2)\mmod W$:
    $$\spec \H^{\G_2}_\ast(\Omega \RP^7; \Z') \cong \fr{t}^{2}(2)\mmod W \times_{\ld{M}^\ddag/\ld{G}(-2\rho)} \fr{t}^{2}(2)\mmod W.$$
    In particular, the conditions of \cref{thm: ordinary homology criterion satake} hold for the spherical $\SO_7$-variety $\SO_7/\G_2$.
\end{prop}
\begin{proof}
    Let us denote by $Y$ the $\SL_2$-orbit of the image of $\kappa: \fr{t}^{2}(2)\mmod W \to \fr{t}^1(2)\mmod W' \times \sl_2(6-6\rho_{\ld{G}_X})$.
    Then $Y \cong \fr{t}^1(2)\mmod W' \times \sl_2^\reg(6)$, and it is well-known that the complement of $\sl_2^\reg \subseteq \sl_2$ has complement of codimension $\geq 2$. Applying \cref{lem: some codim stuff} to the $\SL_2$-equivariant inclusion $Y \hookrightarrow \fr{t}^1(2)\mmod W' \times \sl_2(6-6\rho_{\ld{G}_X})$ and the map $\SL_2 \to \ld{G} = \Sp_6$, we conclude that the $\ld{G}$-orbit of the image of $\kappa_{\ld{M}^\ddag}$ has complement of codimension $\geq 2$.

    By \cref{prop: homology RP7}, there is an isomorphism of graded schemes over $\fr{t}^{2}(2)\mmod W$:
    $$\spec \H^{\G_2}_\ast(\Omega \RP^{7}; \Z') \cong \Z'[x_1, x_2, a^{\pm 1}, \tfrac{a-a^{-1}}{x_1 x_2 (x_1 + x_2)}]^{W}.$$
    The action of $\Sigma_3 \subseteq W$ on $\H^\ast_{T^{2}}(\ast; \Z') = \Z'[x_1, x_2]$ is given by the reduced standard representation (i.e., $x_1$, $x_2$, and $-(x_1 + x_2)$ are permuted, and $a$ is fixed), and the action of $\Z/2\subseteq W$ is given by negating the $x_i$ simultaneously and sending $a \mapsto a^{-1}$. It follows that
    $$\Z'[x_1, x_2, a^{\pm 1}, \tfrac{a-a^{-1}}{x_1 x_2 (x_1 + x_2)}]^{W} \cong \Z'[c_2, c_6, a + a^{-1}, \tfrac{a - a^{-1}}{x_1 x_2 (x_1 + x_2)}].$$
    On the other hand, there is a graded isomorphism
    \begin{align*}
        \fr{t}^{2}(2)\mmod W \times_{\ld{M}^\ddag/\ld{G}(-2\rho)} \fr{t}^{2}(2)\mmod W & \cong  \fr{t}^{2}(2)\mmod W \times_{(\fr{t}^1(2)\mmod W' \times \sl_2(6-6\rho_{\ld{G}_X}))/\SL_2(-6\rho)} \fr{t}^{2}(2)\mmod W\\
        & \cong \fr{t}^1(2)\mmod W' \times (\AA^1(6)\mmod (\Z/2) \times_{\sl_2(6-6\rho_{\ld{G}_X})/\SL_2(-6\rho)} \AA^1(6)\mmod (\Z/2)).
    \end{align*}
    By construction, the map $\AA^1(6)\mmod (\Z/2) \cong \AA^1(12) \to \sl_2(6-6\rho_{\ld{G}_X})/\SL_2(-6\rho)$ is precisely a shifted version of the Kostant slice, so the discussion in \cite[Remark B.4]{grg-reg} implies that if we write $c_6$ to denote the coordinate on $\AA^1(12)$, there is an isomorphism
    \begin{align*}
        \AA^1(6)\mmod (\Z/2) \times_{\sl_2(6-6\rho_{\ld{G}_X})/\SL_2(-6\rho)} \AA^1(6)\mmod (\Z/2) & \cong \spec \Z'[c_3, a^{\pm 1}, \tfrac{a - a^{-1}}{x_1 x_2 (x_1 + x_2)}]^{\Z/2} \\
        & \cong \spec \Z'[c_6, a + a^{-1}, \tfrac{a-a^{-1}}{c_3}],
    \end{align*}
    where $c_3^2 = c_6$. Therefore, there is an isomorphism 
    $$\fr{t}^{2}(2)\mmod W \times_{\ld{M}^\ddag/\ld{G}(-2\rho)} \fr{t}^{2}(2)\mmod W \cong \spec \Z'[c_2, c_6, a + a^{-1}, \tfrac{a - a^{-1}}{x_1 x_2 (x_1 + x_2)}],$$
    which gives the desired claim.
\end{proof}
\begin{remark}
    Using \cref{prop: homology RP7}, one finds (by arguing as in \cref{ex: SO2n+1 and GLn+1 good}) that there is an isomorphism
    $$\spec \H^{\SL_3}_\ast(\Omega (\SO_7/\G_2); \QQ) \cong \left( \AA^1(6) \times_{\tilde{\sl_2}(6-6\rho)/\SL_2(-6\rho)} \AA^1(6)\right) \times \fr{sp}_2^\ast(2)\mmod \Sp_2$$
    of graded group schemes over $\spec \H^\ast_{\SL_3}(\ast; \QQ) \cong \AA^1(6) \times \fr{sp}_2^\ast(2)\mmod \Sp_2$. Here, again, the map $\AA^1 \to \tilde{\sl_2}$ is obtained by pulling back the Kostant section $\AA^1 \to \sl_2$ along the Grothendieck-Springer resolution. (See also \cref{rmk: K-equiv-homology} for an interpretation via Hochschild cohomology.)
\end{remark}

\begin{proof}[Proof of \cref{thm: bzsv for B3'}]
    This follows from \cref{prop: ordinary reg centr for B3'} and \cref{thm: ordinary homology criterion satake}, along with the isomorphism between $\ld{M}^\ddag/\ld{G}(-2\rho)$ and $\fr{t}^1(2)\mmod W' \times \sl_2(6-6\rho_{\ld{G}_X})/\ld{G}_X(-6\rho_{\ld{G}_X})$.
\end{proof}
\begin{remark}\label{rmk: levi sqrt B3'}
    Note that the normalization term $\fr{t}^{n-2}(2)\mmod W'$ identifies with $\fr{so}_{3}(2)\mmod \SO_{3}$, which is the invariant quotient for the group $L^\wedge_X$ from \cite{knop-schalke}.
\end{remark}

%% file: ku-hamiltonian/beta-hamiltonian.tex
\subsection{$\ku$-Hamiltonian spaces}\label{sec: ku-hamiltonian}

Our goal in this section is to place the calculations of the preceding section into a broader context. The basic topic of study in this section is the $\beta$-deformation of a group scheme introduced in \cref{def: Hbeta defn}. We will soon focus on the graded quotient stack $\ld{G}(-2\rho)_\beta/\ld{G}(-2\rho)$ appearing in \cref{thm: ku derived satake}.
\begin{remark}
    Since we will work in the setting of graded schemes, and both $\ld{G}(-2\rho)_\beta$ and $\ld{G}(-2\rho)$ are shifted by the same cocharacter of $\ld{T}$, we can (and will) simply ignore this cocharacter. In other words, we will focus only on the quotient stack $\ld{G}_\beta/\ld{G}$, instead of $\ld{G}(-2\rho)_\beta/\ld{G}(-2\rho)$.
\end{remark}
We will begin with a brief review of the theory of shifted symplectic stacks; for optimal generality, we will work over a discrete commutative ring $R$.  Many of the results below were proved in \cite{ptvv, calaque-lagr, safronov-cs} under the assumption that $R$ is a $\QQ$-algebra, but this is often superfluous.
\begin{recall}[\cite{ptvv, calaque-lagr}]\label{recall: shifted sympl}
    Let $X$ be a derived $R$-stack which admits a (co)tangent complex which is a perfect $\co_X$-module, and let $\F^\star_\H \dR_{X/R}$ denote the Hodge-filtered de Rham complex of $X$ relative to $R$, so that $\gr^n_\H \dR_{X/R} = (\wedge^n L_{X/R})[-n]$. A \textit{closed $j$-form of degree $n$} on $X$ is a global section $\omega$ of $\F^{\geq j}_\H \dR_{X/R}[n+j]$; let $\Omega^{j,\cl}_{BG/R, n}$ denote $\H^0(X; \F^{\geq j}_\H \dR_{X/R}[n+j])$. A closed $2$-form $\omega$ of degree $n$ defines an \textit{$n$-shifted symplectic structure} on $X$ if the section of $(\wedge^2 L_{X/R})[n] \simeq \Sym^2(L_{X/R}[1])[n-2]$ defined by the image of $\omega$ under the map 
    $$\F^{\geq 2}_\H \dR_{X/R}[n+2] \to \gr^2_\H \dR_{X/R}[n+2] \cong (\wedge^2 L_{X/R})[n]$$
    defines an equivalence $T_{X/R} \xar{\sim} L_{X/R}[n]$. If $X$ is an $n$-shifted symplectic stack (the closed $2$-form will be left implicit in the notation), let $\ol{X}$ denote $X$ equipped with the opposite symplectic structure.

    Let $X$ be an $n$-shifted symplectic stack, and let $f: L \to X$ be a morphism of derived $R$-stacks, where $L$ and $f$ each admit perfect (co)tangent complexes. An \textit{isotropic structure} on $f$ is a nullhomotopy of the composite
    $$T_L \to f^\ast T_X \stackrel{\omega}{\simeq} f^\ast L_X[n] \to L_L[n].$$
    An isotropic structure is called \textit{Lagrangian} if the above composite is a cofiber sequence. If $X$ and $Y$ are $n$-shifted symplectic stacks, a \textit{Lagrangian correspondence} is a Lagrangian morphism $L \to X \times \ol{Y}$. By \cite[Theorem 1.2]{safronov-cs}, if $L_1 \to X \times \ol{Y}$ and $L_2 \to Y \times \ol{Z}$ are Lagrangian correspondences, the fiber product $L_1 \times_Y L_2 \to X \times \ol{Z}$ is also a Lagrangian correspondence. As a special case, if $L_1, L_2 \to X$ are Lagrangian morphisms to an $n$-shifted symplectic stack, the fiber product $L_1 \times_X L_2$ admits the structure of an $(n-1)$-shifted symplectic stack.
\end{recall}
\begin{prop}\label{prop: BG 2-shifted}
    Let $G$ be a split reductive group over $R$, and assume that there is a nondegenerate $G$-invariant quadratic form $q$ on $\g$. If its torsion primes are inverted in $R$, any such $q$ on $\g$ (more precisely, a choice of lift of $q$ along $\H^4_\dR(BG/R) \to \Sym^2(\g^\ast)^G$, which does indeed exist) defines a $2$-shifted symplectic structure on $BG$.
\end{prop}
\begin{proof}
    Any nondegenerate $G$-invariant quadratic form $q$ on $\g$ defines a section of $\H^2(BG; \wedge^2 L_{BG/R})$. The object $L_{BG/R}\in \Perf(BG)$ can be identified with the coadjoint representation $\g^\ast[-1]$. The underlying graded $R$-algebra of the de Rham complex $\dR_{BG/R}$ can be identified with $\Gamma^\ast_R(L_{BG/R}[-1]) \cong \Gamma^\ast_R(\g^\ast[-2])$. There is a d\'ecalage isomorphism $\Gamma^n_R(M[-2]) \simeq \Sym^n_R(M)[-2n]$ for any bounded-below $R$-module $M$, and hence an isomorphism $\Gamma^\ast_R(\g^\ast[-2]) \cong \sh \Sym^\ast_R(\g^\ast(-1))$. 
    
    Now, the Hodge-de Rham spectral sequence runs
    $$E_1^{i,j} \cong \H^j(BG; \wedge^i L_{BG/R}) \cong \H^{i-j}(G; \Sym^i(\g^\ast)) \Rightarrow \H^{i+j}_\dR(BG/R).$$
    By \cite[Theorem 10.2]{totaro-bg}, the Hodge-de Rham spectral sequence for $BG$ degenerates at the $E_1$-page since the torsion primes for $G$ are inverted in $R$. Therefore, the map $\H^4(BG; \F^{\geq 2}_\H \dR_{BG/R}) \to \H^2(BG; \wedge^2 L_{BG/R})$ is surjective, and any choice of a lift of $q$ along the map defines a closed $2$-form of degree $2$ on $BG$. It is easy to see that this closed $2$-form defines a $2$-shifted symplectic structure if and only if the associated $G$-invariant bilinear form on $\g$ is nondegenerate.
\end{proof}
\begin{remark}\label{rmk: when nondeg quad}
    Suppose $G$ is semisimple and simply-connected over $R = \FF_p$. Then $\g$ admits a nondegenerate invariant quadratic form if and only if:
    \begin{itemize}
        \item $p\nmid (n+1)$ if $G$ is of type $A_n$ (but $p$ can be arbitrary if $G = \GL_n$);
        \item $p\neq 2$ if $G$ is of type $B_n$ or $C_n$ with $n\geq 2$, $D_n$ with $n\geq 4$, $\F_4$, or $\mathrm{E}_7$;
        \item $p\neq 3$ if $G$ is of type $\G_2$ or $\mathrm{E}_6$;
        \item $p$ arbitrary if $G$ is of type $\mathrm{E}_8$.
    \end{itemize}
    Indeed, $p$ satisfies the above conditions if and only if $\g$ admits a \textit{nondegenerate} invariant symmetric \textit{bilinear} form over $\FF_p$.
    See \cite[Section 6.4(b)]{jantzen-reps-lie-prime-char} and \cite[Proposition 4]{gross-nebe} for a reference; in particular, the table following \cite[Proposition 4]{gross-nebe} determines the order of the cokernel of the resulting map $\g_\Z \to \g_\Z^\ast$. Since a quadratic form is nondegenerate when the same is true of its associated symmetric bilinear form (by definition), we only need to check that in each of these cases, the bilinear form admits a quadratic refinement. The relevant $G$-invariant bilinear form on $\g$ is induced from the bilinear form $\pdb{-,-} = \frac{1}{2h^\vee} \kappa(-,-)$ defined on $\g_\Z$, where $\kappa$ is the Killing form and $h^\vee$ is the dual Coxeter number. By \cite[Proposition 4]{gross-nebe}, the bilinear form $\pdb{-,-}: \g_\Z \times \g_\Z \to \Z$ is \textit{even}, and so it automatically admits a quadratic refinement as desired.

    Looking at the torsion primes for $G$, we see from \cref{prop: BG 2-shifted} that $BG$ admits a $2$-shifted symplectic structure over $\FF_p$ if:
    \begin{itemize}
        \item $2(n+1)\neq 0\pmod{p}$ if $G$ is of type $A_n$ (but $p$ can be arbitrary if $G = \GL_n$);
        \item $p\neq 2$ if $G$ is of type $B_n$ or $C_n$ with $n\geq 2$, or $D_n$ with $n\geq 4$;
        \item $p\neq 2,3$ if $G$ is of type $\G_2$, $\F_4$, $\mathrm{E}_6$, or $\mathrm{E}_7$;
        \item $p\neq 2,3,5$ if $G$ is of type $\mathrm{E}_8$.
    \end{itemize}
\end{remark}
From now, all reductive groups $G$ over $R$ will be split, and we will assume that enough primes are inverted in $R$ so that the hypotheses of \cref{prop: BG 2-shifted} are satisfied.
\begin{prop}[{\cite[Theorem 2.5]{ptvv}, \cite[Theorem 3.5]{safronov-cs}}]\label{prop: mapping stack shifted}
    Let $X$ be an \textit{$d$-oriented stack}, i.e., a stack equipped with a map $\Gamma(X; \co_X) \to R[-d]$ such that for any animated $R$-algebra $A$ and any $\cf\in \Perf(X \otimes_R A)$, the induced map $\Gamma(X \otimes_R A; \cf)^\vee \to \Gamma(X \otimes_R A; \cf[-d]^\vee)$ is an isomorphism.
    The data of an $n$-shifted symplectic structure on $Y$ equips the mapping stack $\Map(X, Y)$ with an $(n-d)$-shifted symplectic structure.
    More generally, if $L \to Y$ is a Lagrangian morphism, the induced map $\Map(X, L) \to \Map(X, Y)$ acquires a natural Lagrangian structure.
\end{prop}
\begin{example}\label{ex: BGa sharp}
    Let $\GG_a^\sharp(-2)$ denote the divided power hull of the origin in $\GG_a$ equipped with a $\GG_m$-action placing its coordinate in weight $2$, and let $\widehat{\GG_a^\sharp(-2)}$ denote the completion of $\GG_a^\sharp(-2)$ at the divided power filtration. Then $X = B\widehat{\GG_a^\sharp(-2)}$ is a $1$-oriented stack: the cohomology of its structure sheaf is isomorphic to $R[\epsilon]/\epsilon^2$ with $\epsilon$ in homological degree $-1$ (and weight $-2$). Moreover, if $Y$ is any derived $R$-stack which admits a cotangent complex, and ${T_Y^\sharp(2)}$ denotes the PD-hull of the zero section of the tangent bundle of $Y$, the mapping stack $\Map(B\widehat{\GG_a^\sharp(-2)}, Y)$ can be identified with the stack ${BT_Y^\sharp(2)}$.
    It follows from \cref{prop: mapping stack shifted} that the data of an $n$-shifted symplectic structure on $Y$ equips ${BT_Y^\sharp(2)}$ with an $(n-1)$-shifted symplectic structure.

    Let $G$ be a reductive group over $R$ equipped with a nondegenerate $G$-invariant quadratic form on $\g$. When $Y = BG$, we may identify $T_{BG} \cong \g[1] \in \Perf(BG)$, which implies that ${BT_{BG}^\sharp(2)} \cong \g(2)/G$. It follows from \cref{prop: BG 2-shifted} that $\g(2)/G$ admits a $1$-shifted symplectic structure of weight $2$; forgetting the grading, we see that $\g/G$ admits a $1$-shifted symplectic structure.
\end{example}
\begin{example}
    Let $S^1$ denote the constant stack $B\Z$. Then $S^1$ is a $1$-oriented stack, since $\Gamma(S^1; \co) \cong C^\ast(S^1; R)$, and the circle admits a canonical orientation. It follows from \cref{prop: mapping stack shifted} that the data of an $n$-shifted symplectic structure on $Y$ equips $\Map(S^1, Y)$ with an $(n-1)$-shifted symplectic structure. In particular, if $G$ is a reductive group over $R$, then applying \cref{prop: BG 2-shifted} implies that $\Map(S^1, BG) \cong G/G$ admits a $1$-shifted symplectic structure.
\end{example}
\begin{prop}[{\cite{safronov-cs}}]\label{prop: safronov lag maps}   
    Let $R = \cc$, and let $G$ be a complex reductive group.
    A Lagrangian morphism $L \to \g/G$ is equivalent to the data of a Hamiltonian $G$-space. Similarly, a Lagrangian morphism $L \to G/G$ is equivalent to the data of a quasi-Hamiltonian $G$-space in the sense of \cite{amm-qham}.
\end{prop}
The reader not familiar with the definition of a quasi-Hamiltonian $G$-space can take the second part of \cref{prop: safronov lag maps} to be a definition. Motivated by \cref{lem: Hbeta invert and mod beta}, we make the following observation.
\begin{lemma}\label{lem: Gbeta 1-oriented}
    The stack $B\GG_\beta^\vee$ admits a canonical $1$-orientation of weight $2$.
\end{lemma}
\begin{proof}
    We will just construct the orientation on $\GG_\beta^\vee$, and leave verifying the properties of \cref{prop: mapping stack shifted} to the reader. To compute the cohomology $\H^\ast(B\GG_\beta^\vee; \co)$, let us first compute $\GG_\beta^\vee$. As shown in \cite[Proposition C.6]{thh-xn}, the Cartier dual $\GG_\beta^\vee$ is isomorphic to $\spf \Z[\beta, \frac{1}{n!} \prod_{j=0}^{n-1} (y-j\beta)]^\wedge$ where the element $y$ is primitive (i.e., the coproduct sends $y \mapsto y \otimes 1 + 1 \otimes y$) and lives in weight $2$. Here, the completion is taken with respect to the $\beta$-deformed divided power filtration (i.e., with respect to $\frac{1}{n!} \prod_{j=0}^{n-1} (y-j\beta)$ for $n\geq 1$). Computing the cohomology of the trivial $\GG_\beta^\vee$-representation in the standard manner shows that $\H^0(B\GG_\beta^\vee; \co) \cong R$, $\H^1(B\GG_\beta^\vee; \co)$ is isomorphic to the submodule of primitive elements in $\co_{\GG_\beta^\vee}$, and $\H^j(B\GG_\beta^\vee; \co)$ is zero for $j>1$. It is not difficult to see that the only primitive elements in $\co_{\GG_\beta^\vee}$ are scalar multiples of $y$, and so the cohomology ring $\H^\ast(B\GG_\beta^\vee; \co)$ is exterior on a single class in cohomological degree $1$ and weight $2$. This generator of $\H^1(B\GG_\beta^\vee; \co)$ gives the desired $1$-orientation of $B\GG_\beta^\vee$.
\end{proof}
\begin{remark}
    The $1$-orientation on $B\GG_\beta^\vee$ is closely connected to the preorientation on $\GG_{\beta, \ku}$ from \cref{prop: Gbeta is preoriented}. Indeed, the $1$-orientation on $B\GG_\beta^\vee$ can be viewed as a map $B\GG_\beta^\vee \to B\GG_a$. Since there is an isomorphism 
    $$\Map(B\GG_\beta^\vee, B\GG_a) \cong \Hom(\GG_\beta^\vee, \GG_a)/\GG_a \cong \Hom(\GG_\beta^\vee, \GG_a) \times B\GG_a,$$
    we can identify
    $$\H^1(B\GG_\beta^\vee; \co) \cong \pi_0 \Map(B\GG_\beta^\vee, B\GG_a) \cong \pi_0 \Hom(\GG_\beta^\vee, \GG_a),$$
    and the $1$-orientation on $B\GG_\beta^\vee$ can be viewed as a homomorphism $\GG_\beta^\vee \to \GG_a$. However, $\Hom(\GG_\beta^\vee, \GG_a)$ is isomorphic to the Lie algebra of $\GG_\beta$, and a section of this Lie algebra is precisely the datum of a preorientation on $\GG_\beta$.
\end{remark}
\begin{prop}\label{prop: Gbeta mod G is 1-shifted}
    Let $G$ be a reductive group over $R$.
    The choice of a nondegenerate $G$-invariant quadratic form on $\g$ equips $\Map(B\GG_\beta^\vee, BG) \cong G_\beta/G$ with a $1$-shifted symplectic structure of weight $2$.
    If $G$ is semisimple and $G'$ is centrally isogenous of $G$\footnote{For instance, $G$ is simply-connected and $G' = G/Z(G)$; or $G$ is adjoint and $G'$ is the natural $\pi_1(G)_\mathrm{tors}$-cover.}, the quotient stack $G'_\beta/G$ admits a $1$-shifted symplectic structure of weight $2$.
\end{prop}
\begin{proof}
    This is a consequence of \cref{prop: BG 2-shifted}, \cref{prop: mapping stack shifted}, and \cref{lem: Gbeta 1-oriented}.
\end{proof}
\begin{definition}\label{def: ku-hamiltonian}
    Let $G$ be a reductive group over $R$.
    Let $M$ be a graded $R[\beta]$-stack which admits a cotangent complex, and suppose $M$ is equipped with a $G$-action. A \textit{$\ku$-Hamiltonian structure} on $M$ is a Lagrangian morphism $M/G \to G_\beta/G$. The resulting $G$-equivariant map $M \to G_\beta$ will be called the \textit{$\ku$-moment map}. If $G$ is semisimple, a slight variant of this definition would be to ask for a Lagrangian morphism $M/G \to G_\beta/G$; we will sometimes refer to such data also as a $\ku$-Hamiltonian structure on $M$.
\end{definition}
\begin{remark}
    Motivated by \cref{lem: ku loops and HKR}, one alternative name for $\ku$-Hamiltonian structures might be ``HKR-filtered quasi-Hamiltonian structures'' (where HKR stands for Hochschild-Kostant-Rosenberg).
\end{remark}
\begin{remark}
    As in \cref{rmk: morava k-theory Gbeta}, it is not really crucial to treat $\GG_\beta$ as the fundamental object here: we could have considered {any} $1$-dimensional group scheme over $\Z[\beta]$ in place of $\GG_\beta$. Suppose we permit $1$-dimensional \textit{formal} group schemes. Fix a prime $p$, let $\widehat{\GG}_{k_\Z(n)}$ be the $1$-dimensional formal group over $\Z_{(p)}[\beta]$ from \cref{rmk: morava k-theory Gbeta}, and let $G$ be a reductive group over $\Z_{(p)}[\beta]$. One can then prove that $\Map(B\widehat{\GG}_{k_\Z(n)}^\vee, BG)$ admits a $1$-shifted symplectic structure over $\Z_{(p)}[\beta]$; Lagrangian morphisms to this stack provide a(n integral) Morava K-theoretic analogue of $\ku$-Hamiltonian spaces. Again, we will not study these objects further in the present article, but we hope to in \cite{Eodd-and-quantizations}.
\end{remark}
\begin{remark}\label{rmk: ku-ham and radius}
    In the setting of differential geometry (as studied in \cite{amm-qham}), where $G$ is replaced by a compact Lie group, \cite[Theorem 8.3]{amm-qham} proves that there is an equivalence between quasi-Hamiltonian $G$-spaces and Hamiltonian $LG$-spaces (the latter needing some care to define because of infinite-dimensional analytic issues). It seems likely that the notion of a $\ku$-Hamiltonian structure in the setting of differential geometry might be equivalent to the theory of Hamiltonian spaces for $\Map(S^1_\beta, G)$, with $S^1_\beta$ being defined as
    $$S^1_\beta = \{(\beta, z) | \|z\| = \beta\}\subseteq \RR \times \cc,$$
    and the $\Map(S^1_\beta, G)$-Hamiltonian space is equipped with a compatible fibration to $\RR$. In other words, the radii of the fibers of $S^1_\beta \to \RR$ should be related to the parameter $\beta$.
\end{remark}
\begin{prop}\label{prop: beta cotangent of G/UP}
    Let $P\subseteq G$ be a parabolic subgroup, let $U_P\subseteq P$ denote its unipotent radical, and let $L = P/U_P$ denote the Levi quotient. Then $\Ind_{U_P}^G P_\beta$ admits the structure of a $\ku$-Hamiltonian $G \times L$-space where the $\ku$-moment map $\Ind_{U_P}^G P_\beta \to G_\beta$ is given by conjugation.
\end{prop}
\begin{proof}
    As shown in \cite[Lemma 3.4]{safronov-cs}, the maps $BP \to BG$ and $BP \to BL$ define a Lagrangian correspondence $BP \to BG \times BL$, essentially because there is an exact sequence
    $$0 \to \fr{p} \to \g \oplus \fr{l} \to \fr{p}^\ast \to 0.$$
    It follows from \cref{prop: mapping stack shifted} and \cref{lem: Gbeta 1-oriented} that there is a Lagrangian correspondence 
    $$P_\beta/P \cong \Map(B\GG_\beta^\vee, BP) \to \Map(B\GG_\beta^\vee, BG \times BL) \cong G_\beta/G \times L_\beta/L.$$
    Since $P_\beta/P \times_{BG \times BL} \spec \Z[\beta] \cong \Ind_{U_P}^G P_\beta$, this produces the desired $\ku$-Hamiltonian $G \times L$-structure on $\Ind_{U_P}^G P_\beta$.
\end{proof}
\begin{remark}
    Upon inverting $\beta$ and quotienting by $\GG_m$, \cref{prop: beta cotangent of G/UP} reduces to \cite[Theorem 9]{boalch}.
\end{remark}
\begin{prop}\label{prop: B mod N is ku hamiltonian}
    Let $B\subseteq G$ be a Borel subgroup with unipotent radical $N$. Then $\Ind_N^G B_\beta$ admits the structure of a $\ku$-Hamiltonian $G$-space, where the $\ku$-moment map $\Ind_N^G B_\beta \to G_\beta$ is given by conjugation.
\end{prop}
\begin{proof}
    Let us first show show that the map $(\Ind_N^G B_\beta)/G \cong B_\beta/N \to G_\beta/G$ admits a Lagrangian structure. Using \cref{prop: beta cotangent of G/UP}, \cref{recall: shifted sympl}, and the isomorphism $B_\beta/B \times_{G_\beta/G \times T_\beta/T} (G_\beta/G \times T_\beta) \cong B_\beta/N$, it suffices to show that the map $T_\beta \to T_\beta/T$ is a Lagrangian morphism. For this, note that the tangent complex to $H_\beta/H$ is given by the complex $\ul{\fr{h}} \to T_{H_\beta}$, where $\ul{\fr{h}}$ denotes $\fr{h} \otimes \co_{H_\beta}$, and the differential $\fr{h} \to \Gamma(H_\beta; \co)$ is given by the adjoint action $\xi \mapsto \xi^R - \xi^L$. Here, $\xi^L$ and $\xi^R$ denote the vector fields generating the left and right actions of $H$ on $H_\beta$. (Note that there is an isomorphism $T_{H_\beta} \cong \ul{\fr{h}}$.) When $H$ is commutative (such as $H = T$), $\xi^R = \xi^L$, and so the tangent complex is split. There is an obvious cofiber sequence
    $$\ul{\fr{t}} \to \ul{\fr{t}}\oplus \ul{\fr{t}}[1] \to \ul{\fr{t}}[1],$$
    which identifies with the cofiber sequence
    $$T_{T_\beta} \to T_{T_\beta/T} \cong L_{T_\beta/T} \to L_{T_\beta}[1].$$
    This gives the desired Lagrangian structure.
\end{proof}
\begin{remark}\label{rmk: affine closure of SL2 mod Ga and Bbeta}
    Suppose $G = \SL_2$, so that the affine closure of $G/N \cong \SL_2/\GG_a \cong \AA^2 - \{0\}$ is smooth. Using the algebraic Hartogs lemma, one can show that the $\ku$-Hamiltonian structure on the $G$-space $\Ind_N^G B_\beta \cong \SL_2 \times^{\GG_a} B_\beta$ from \cref{prop: B mod N is ku hamiltonian} extends to a $\ku$-Hamiltonian structure on its affine closure $\ol{\SL_2 \times^{\GG_a} B_\beta}$.
    It seems reasonable to expect that this is true in general, i.e., that the $\ku$-Hamiltonian structure on the $G$-space $\Ind_N^G B_\beta$ from \cref{prop: B mod N is ku hamiltonian} extends to $\ol{\Ind_N^G B_\beta}$. This is not immediately clear from the perspective of derived algebraic geometry, since neither of the affine closures $\ol{G/N}$ or $\ol{\Ind_N^G B_\beta}$ are smooth outside of the rank one case.

    For the sake of concreteness, let us describe the affine closure $\ol{\SL_2 \times^{\GG_a} B_\beta}$ explicitly. Fix $\begin{psmallmatrix}
        a & b\\
        c & d
    \end{psmallmatrix} \in \SL_2$ and $(x,y)\in B_\beta$ with $a$, $d$, and $y$ in weight $0$, $c$ and $x$ in weight $-2$, and $b$ in weight $2$. Then the action of $\begin{psmallmatrix}
        1 & z \\
        0 & 1
    \end{psmallmatrix} \in \GG_a$ (with $z$ in weight $2$) sends
    $$\begin{psmallmatrix}
        a & b\\
        c & d
    \end{psmallmatrix}\mapsto \begin{psmallmatrix}
        a & az + b\\
        c & cz + d
    \end{psmallmatrix}, \ (x,y) \mapsto (x, y + z[-2](x)),$$
    where $[-2](x) = -\frac{\beta x^2 + 2x}{1 + \beta x}$ is the $(-2)$-series of $x$ in the group law on $\GG_\beta$.
    In particular, the $\GG_a$-action fixes $a$, $c$, and $x$, as well as $ay - [-2](x) b$ and $cy - [-2](x) d$. If we write $B = ay - [-2](x) b$ (in weight $0$) and $D = cy - [-2](x) d$ (in weight $-2$), the only relation is
    \begin{equation}\label{eq: reln closure Tbeta A2}
        cB - aD = [-2](x).
    \end{equation}
    In other words, $\ol{\SL_2 \times^{\GG_a} B_\beta}$ is cut out inside $\AA^4_{\Z[\beta]} \times_{\spec \Z[\beta]} \GG_\beta$ by the above equation, where the affine space has coordinates $a,c,B,D$.
    
    If we had instead replaced the $\GG_a$-action on $B_\beta$ by conjugation with the $\GG_a$-action on $\ld{V}_\beta = \AA^1 \times \GG_\beta$ (following the notation of \cref{cstr: CPn kappa ku} with $n = 1$), the above equations would continue to hold if $[-2](x)$ was replaced by $\ol{x}$.
    The analogue of the equation \cref{eq: reln closure Tbeta A2} implies that $\ol{x} = cB - aD$, so $\ol{\SL_2 \times^{\GG_a} V_\beta}$ is the open subscheme of $\AA^4_{\Z[\beta]}$ given by the complement of the hypersurface
    $$1 + \beta (cB - aD) = 0.$$
    Note that when $\beta = 0$, this is the entirety of $\AA^4$. Just as with $\ol{\SL_2 \times^{\GG_a} B_\beta}$, the scheme $\ol{\SL_2 \times^{\GG_a} V_\beta}$ also admits a $\ku$-Hamiltonian structure for its natural $\SL_2$-action, and perhaps deserves to be called $T^\ast_\beta \AA^2$.
\end{remark}

\begin{remark}\label{rmk: subgroup coisotropic}
    If $H \subseteq G$ is a closed subgroup, it is natural to ask whether there is a ``$\ku$-theoretic'' cotangent bundle $T^\ast_\beta(G/H)$ whose fiber at $\beta = 0$ is $T^\ast(G/H)$, and which admits the structure of a $\ku$-Hamiltonian $G$-space? It seems rather difficult to define such an object for arbitrary subgroups $H$. However, in the case that the annhilator of ${\fr{h}} \subseteq \g$ (under an invariant bilinear form on $\g$) is itself a Lie subalgebra, we expect this to be possible, but we will not study this topic here. (For instance, \cref{prop: B mod N is ku hamiltonian} fits into this general class of examples.) Upon inverting $\beta$, i.e., working with quasi-Hamiltonian $G$-spaces, this was shown in \cite{balibanu-mayrand}. 
\end{remark}

%% file: ku-hamiltonian/duals-are-lagrangian.tex
\subsection{Functoriality of Hochschild cohomology}\label{subsec: ku and lagrangians}

Following \cite{bzsv}, it is natural to hope:
\begin{expect}\label{expect: lag to Gbeta}
    Suppose that $G$ is a simply-laced and connected compact Lie group, $H_\cc \subseteq G_\cc$ is a closed reductive spherical subgroup, and $\ld{M}_\beta$ is a ``dual'' affine graded $\ld{G}$-variety over $\Z[\beta]$ satisfying the hypotheses of \cref{prop: ku rk 1 reg centr --> thm} equipped with a morphism $\ld{M}_\beta \to G_\beta$. Then, there is an equivalence $\Shv^{c,\Sat}_{G\pw{t}}(G\ls{t}/H\ls{t}; \ku)^\faux \simeq \Perf(\sh^{1/2}\ld{M}_\beta/\ld{G}(-2\rho))$ which is compatible with the action of 
    $$\Shv_{(G \times G)\pw{t}}^{c,\Sat}(G\ls{t}; \ku)^\faux \simeq \Perf(\sh^{1/2} G(-2\rho)_\beta/\ld{G}(-2\rho)),$$
    the equivalence coming frm \cref{thm: ku derived satake}.
    Based on \cref{conj: bzsv}, we expect that the map $\ld{M}_\beta/\ld{G} \to G/\ld{G}$ admits a Lagrangian structure, i.e., that $\ld{M}_\beta$ admits the structure of a $\ku$-Hamiltonian $\ld{G}$-space.
\end{expect}
In the generality of \cref{prop: ku rk 1 reg centr --> thm}, it is not clear how one might prove \cref{expect: lag to Gbeta}. This can be shown in the case of $\PGL_2/\GG_m$, but contributions from the Levi factor/Whittaker induction present difficulties in proving \cref{expect: lag to Gbeta} for types $A_n$ with $n>1$, $C_n$, $D_2$, and $\G_2$ as studied in \cref{sec: case by case}. However, we can use the discussion in \cref{subsec: G of LGmodH} to prove some partial results along these lines. The starting point of this discussion is the following.
\begin{observe}\label{obs: self int and E3 poisson}
    Recall that if $R$ is an $\Eoo$-ring and $A$ is a (nonunital) $\E{n}$-$R$-algebra with $n\geq 2$, the homotopy groups $\pi_\ast(A)$ admit the structure of a graded (nonunital) Poisson algebra over $\pi_\ast(R)$, where the graded Poisson bracket has weight $n-1$ (i.e., if $f,g$ are functions in weights $i,j$ respectively, the graded Poisson bracket $\{f,g\}$ is in weight $i+j+n-1$) and is $\pi_\ast(R)$-linear. This graded Poisson bracket comes from action of the generator of $\pi_{n-1} \Conf_2(\RR^n) \cong \pi_{n-1} S^{n-1}$ on $\pi_\ast(R)$. (See, e.g., \cite[Example 4.5]{lawson-dl}.) 

    Let $\widehat{\cf_{G}(\Omega G)}^\vee$ denote the $\ku_G$-linear dual of $\cf_{G\times G}(\cL G)$.
    There is a graded Poisson structure of weight $2$ on $\spf \pi_\ast \widehat{\cf_{G}(\Omega G)}^\vee$ arising from the $\E{3}$-algebra structure on $\widehat{\cf_{G}(\Omega G)}^\vee$ (viewed as the $\E{2}$-center $\dZ_\E{2}(\ku_G/\ku)$ via \cref{cor: loop homology E2}). Using \cref{thm: ku homology LG and langlands mod center}, this implies that (a completion of) the fiber product $T_\beta \mmod W \times_{G_\beta/\ld{G}} T_\beta \mmod W$ admits a graded Poisson structure where the graded Poisson bracket has weight $2$.
\end{observe}
\begin{remark}
    The Kostant slice $\kappa: \fr{t} \mmod W \to \ld{\g}/\ld{G}$ is Lagrangian for the $1$-shifted symplectic structure on $\ld{\g}/\ld{G}$ from \cref{ex: BGa sharp}; see, e.g., \cite[Proposition 4.18]{safronov-geom-quant}. \cref{recall: shifted sympl} implies that the self-intersection $\fr{t} \mmod W \times_{\ld{\g}/\ld{G}} \fr{t} \mmod W$ admits a ($0$-shifted) symplectic structure; the underlying Poisson structure can be identified with the $\beta = 0$ degeneration of \cref{obs: self int and E3 poisson}.
\end{remark}

If \cref{expect: lag to Gbeta} holds, the map
\begin{equation}\label{eq: langlands dual map coisotropic}
    T_{H,\beta}\mmod W_H \times_{\ld{M}_\beta/\ld{G}} T_{H,\beta}\mmod W_H \to T_\beta \mmod W \times_{G_\beta/\ld{G}} T_\beta \mmod W
\end{equation}
from \cref{prop: ku rk 1 reg centr --> thm} will, in particular, be coisotropic (in an appropriate derived sense). See \cref{prop: lagrangian correspondence} and \cref{rmk: lag corr in homogeneous case}. Our goal in this section is to show that this consequence of \cref{expect: lag to Gbeta} is always true; hopefully some variant of our discussion below could imply \cref{expect: lag to Gbeta} itself.
In order to explain this, we need to rephrase the graded Poisson bracket on $T_\beta \mmod W \times_{G_\beta/\ld{G}} T_\beta \mmod W$ in homotopy-theoretic terms. As the reader will observe, the coisotropicity of \cref{eq: langlands dual map coisotropic} is a rather general phenomenon.
Let us begin by reviewing the notion of an $\E{n}$-center.
\begin{recall}[{\cite[Section 5.3]{HA} and \cite{francis}}]\label{recall: centers}
    Let $\cC$ be a presentably symmetric monoidal stable $\infty$-category with unit $\unit$, and let $f: A \to B$ be a morphism in $\Alg_{\E{n}}(\cC)$. The centralizer $\dZ_\E{n}(f)$ is the universal $\E{n}$-algebra object of $\cC$ equipped with the data of commutative diagram
    $$\xymatrix{
    A \ar[dr]_-f \ar[r] & \dZ_\E{n}(f) \otimes A \ar[d] \\
    & B
    }$$
    in $\Alg_\E{n}(\cC)$.
    The existence of centralizers is proved in \cite[Theorem 5.3.1.14]{HA}. One can explicitly identify $\dZ_\E{n}(f) = \Map_{\Mod_A^\E{n}}(A, B) = \Map_{\int_{S^{n-1}} A/\cC}(A, B)$.
    If $f$ is the identity map on $A$, we will simply write $\dZ_\E{n}(A/\cC)$ to denote $\dZ_\E{n}(\id_A)$; moreover, if $R$ is an $\Eoo$-ring, we will write $\dZ_\E{n}(A/R)$ to denote $\dZ_\E{n}(A/\Mod_R)$. 
\end{recall}
\begin{remark}\label{rmk: center to centralizer}
    In the setup of \cref{recall: centers}, there is a canonical $\E{n}$-algebra map $\dZ_\E{n}(A/\cC) \to \dZ_\E{n}(f)$ defined using the universal property of $\dZ_\E{n}(f)$ and the commutative diagram
    $$\xymatrix{
    A \ar[dr]_-{\id_A} \ar[r] & \dZ_\E{n}(A/\cC) \otimes A \ar[d] & \\
    & A \ar[r]_-f & B
    }$$
    in $\Alg_\E{n}(\cC)$.
\end{remark}
\begin{remark}\label{rmk: abstracted scenario}
    The scenario of \cref{prop: ku rk 1 reg centr --> thm} can be modeled as follows. Under the hypotheses of \cref{prop: ku rk 1 reg centr --> thm}, the map \cref{eq: langlands dual map coisotropic} can be identified with the composite map
    \begin{equation}\label{eq: want map coisotropic}
        \spec \ku^H_\ast(\Omega(G/H)) \to \spec \ku^H_\ast(\Omega G) \to \spec \ku^G_\ast(\Omega G).
    \end{equation}
    Let $R = \ku$, $A = \ku_G$, and $B = \ku_H$, so that there is a map $f: A \to B$ of $\Eoo$-$R$-algebras. Following \cref{cor: loop homology E2} and \cref{warning: completion hochschild}, we can identify a completion of $\cf_H(\Omega(G/H))^\vee$ with the Hochschild cohomology $\dZ_\E{1}(B/A)$, and a completion of $\cf_G(\Omega G)^\vee$ with the $\E{2}$-center $\dZ_\E{2}(A/R)$. It is also easy to see that a completion of $\cf_H(\Omega G)^\vee$ can be identified with the centralizer $\dZ_\E{2}(f)$. In particular, the above composite can be identified, at least upon completion, with a map 
    $$\spec \pi_\ast \dZ_\E{1}(B/A) \to \spec \pi_\ast \dZ_\E{2}(f) \to \spec \pi_\ast \dZ_\E{2}(A/R).$$
\end{remark}
\begin{lemma}\label{lem: map on centers}
    Let $\cC$ be a presentably symmetric monoidal stable $\infty$-category, let $A\in \Alg_\E{n+1}(\cC)$, and let $B\in \Alg_\E{n-1}(\Mod_A(\cC))$ with unit map $f: A \to B$. Then there is a canonical $\E{n}$-$\dZ_\E{n}(A/\cC)$-algebra structure on the $\E{n}$-$A$-algebra $\dZ_\E{n-1}(B/\Mod_A(\cC))$ such that the unit map factors as a composite
    $$\dZ_\E{n}(A/\cC) \to \dZ_\E{n}(f) \to \dZ_\E{n-1}(B/A).$$
\end{lemma}
\begin{proof}
    Let $g: B \to B'$ be a map of $\E{n-1}$-$A$-algebras. Then there is a map $\mu_g: \dZ_\E{n}(A/\cC) \otimes \dZ_\E{n-1}(g) \to \dZ_\E{n-1}(g)$ of $\E{n-1}$-$A$-algebras defined using the universal property of $\dZ_\E{n-1}(g)$ as follows. Recall that there are commutative diagrams
    $$\xymatrix{
    A \ar[r] \ar[dr]_-{\id_A} & \dZ_\E{n}(A/\cC) \otimes A \ar[d] \\
    & A
    }, \ \xymatrix{
    B \ar[r] \ar[dr]_-{g} & \dZ_\E{n-1}(g) \otimes_A B \ar[d] \\
    & B',
    }$$
    where the maps in the first diagram are of $\E{n}$-algebras in $\cC$, and the maps in the second diagram are of $\E{n-1}$-$A$-algebras. Since $A$ is an $\E{n+1}$-algebra, the first diagram can be upgraded to a commutative diagram of $\E{n}$-$A$-algebras. Tensoring these two diagrams over $A$ produces a commutative diagram
    $$\xymatrix{
    B \ar[r] \ar[dr]_-{g} & (\dZ_\E{n}(A/\cC) \otimes \dZ_\E{n-1}(g)) \otimes_A B \ar[d] \\
    & B'
    }$$
    of $\E{n-1}$-$A$-algebras, which gives the desired map of $\E{n-1}$-$A$-algebras
    $$\mu_g: \dZ_\E{n}(A/\cC) \otimes \dZ_\E{n-1}(g) \to \dZ_\E{n-1}(g).$$
    It is not difficult to see that this map is compatible with composition in $g$, in the sense that if $g': B' \to B''$ is another morphism and $c: \dZ_\E{n-1}(g) \otimes_A \dZ_\E{n-1}(g') \to \dZ_\E{n-1}(g' \circ g)$ is the composition coming from functoriality of $\E{n-1}$-centers, there is a commutative digram of $\E{n-1}$-$A$-algebras
    $$\xymatrix{
    (\dZ_\E{n}(A/\cC) \otimes \dZ_\E{n-1}(g)) \otimes_A (\dZ_\E{n}(A/\cC) \otimes \dZ_\E{n-1}(g')) \ar[r]^-{\mathrm{mult}_{\dZ_\E{n}(A/\cC)} \otimes c} \ar[d]_-{\mu_g \otimes \mu_{g'}} & \dZ_\E{n}(A/\cC) \otimes \dZ_\E{n-1}(g' \circ g) \ar[d]^-{\mu_{g' \circ g}}\\
    \dZ_\E{n-1}(g) \otimes_A \dZ_\E{n-1}(g') \ar[r] & \dZ_\E{n-1}(g' \circ g).
    }$$
    Since $A$ is an $\E{n+1}$-algebra, the $\infty$-category $\Alg_\E{n-1}(\LMod_A(\cC))$ admits an $\E{1}$-monoidal structure.
    Taking $g = g' = \id_B$, we find that the map $\mutil := \mu_{\id_B}$ can be upgraded to a map $\dZ_\E{n}(A/\cC) \otimes \dZ_\E{n-1}(B/A) \to \dZ_\E{n-1}(B/A)$ of $\E{1}$-algebras in $\E{n-1}$-$A$-algebras, i.e., of $\E{n}$-$A$-algebras. In particular, this equips $\dZ_\E{n-1}(B/A)$ with the structure of an $\E{n}$-$\dZ_\E{n}(A/\cC)$-algebra.
    
    The factorization through the centralizer $\dZ_\E{n}(f)$ is a consequence of the construction of $\mutil$. Namely, the unit map $\dZ_\E{n}(A/\cC) \to \dZ_\E{n-1}(B/A)$ can be regarded as giving a commutative diagram
    $$\xymatrix{
    B \ar[r] \ar[dr]_-{\id_B} & (\dZ_\E{n}(A/\cC) \otimes A) \otimes_A B \ar[d] \\
    & B
    }$$
    of $\E{n-1}$-$A$-algebra maps, which is in turn obtained via a commutative diagram
    $$\xymatrix{
    A \ar[d]_-f \ar[r] & \dZ_\E{n}(A/\cC) \otimes A \ar[d] \\
    B \ar[r] \ar[dr]_-{\id_B} & (\dZ_\E{n}(A/\cC) \otimes A) \otimes_A B \ar[d] \\
    & B,
    }$$
    where the square is given by the tensor product in $\cC$. The universal property of the centralizer $\dZ_\E{n}(f)$ implies that there is a factorization
    $$\xymatrix{
    B \ar[r] \ar[d]_-{\id_B} & (\dZ_\E{n}(A/\cC) \otimes A) \otimes_A B \ar[d] \\
    B \ar[r] \ar[dr]_-{\id_B} & (\dZ_\E{n}(f) \otimes A) \otimes_A B \ar[d] \\
    & B,
    }$$
    which in turn gives the desired $A$-linear map $\dZ_\E{n}(f) \otimes A \to \dZ_\E{n-1}(B/A)$ factoring the unit $\dZ_\E{n}(A/\cC) \otimes \dZ_\E{n-1}(B/A)$.
\end{proof}
\begin{prop}\label{prop: nonunital}
    Let $\cC$ be a presentably symmetric monoidal stable $\infty$-category, let $A\in \Alg_\E{n+1}(\cC)$, and let $B\in \Alg_\E{n}(\Mod_A(\cC))$ with unit map $f: A \to B$. Then the fiber $\fib(f)$ admits the structure of a nonunital $\E{n+1}$-algebra in $\cC$.
\end{prop}
\begin{proof}
    The map $f$ factors as a composite
    $$A \xar{g} \dZ_\E{n}(B/\cC) \xar{h} B,$$
    where $g$ is a map of $\E{n+1}$-algebras in $\cC$, and the $\E{n}$-map $h: \dZ_\E{n}(B/\cC) \to B$ is the unit. Let $T_{B/\cC}^\E{n}$ denote the $\E{n}$-cotangent complex of $B$ (viewed as an object of $\Alg_\E{n}(\cC)$), so that $T_{B/\cC}^\E{n}[-n]$ admits the structure of a nonunital $\E{n+1}$-algebra in $\cC$, and there is a map $T_{B/\cC}^\E{n}[-n] \to \dZ_\E{n}(B/\cC)$ of nonunital $\E{n+1}$-algebras. We claim that there is a Cartesian square
    $$\xymatrix{
    \fib(f) \ar[r] \ar[d] & T_{B/\cC}^\E{n}[-n] \ar[d] \\
    A \ar[r]_-g & \dZ_\E{n}(B/\cC),
    }$$
    so that $\fib(f)$ is canonically equipped with the structure of a nonunital $\E{n+1}$-algebra in $\cC$. To see this, note that the factorization of $f$ as $h\circ g$ implies a fiber sequence
    $$\fib(g) \to \fib(f) \to \fib(h).$$
    It therefore remains to identify $\fib(h)$ with $T_{B/\cC}^\E{n}[-n]$; but this follows from \cite[Theorem 1.1]{francis} (or equivalently \cite[Theorem 7.3.5.1]{HA}).
\end{proof}
\begin{remark}
    \cref{prop: nonunital} immediately implies \cite[Corollary 8.8]{hill-lawson-Ek}, which states that if $A$ is an $\Eoo$-algebra in $\Sp$, $B$ is an $\E{n}$-$B$-algebra, and $\alpha\in \pi_\ast(A)$ maps to zero in $\pi_\ast(B)$, then all the $\E{n+1}$-algebra Dyer-Lashof operations on $\alpha$ also map to zero in $B$. Indeed, if $f:A \to B$ is the unit map, $\alpha$ lifts to $\pi_\ast \fib(f)$ by assumption; but $\fib(f)$ is a nonunital $\E{n+1}$-algebra by \cref{prop: nonunital}, so all $\E{n+1}$-power operations on $\alpha$ must also map to zero in $\pi_\ast(B)$, as desired.
\end{remark}
We are now in a position to explain the ``coisotropic'' property of the map \cref{eq: langlands dual map coisotropic} (or more generally of \cref{eq: want map coisotropic}). To motivate the use of this term, let us make the following observation.
\begin{observe}\label{obs: what does it mean to be coisotropic}
    Let $P_\ast$ be a graded adic algebra with Poisson bracket of weight $n$, and let $Q_\ast$ be a graded commutative adic $P_\ast$-algebra. The map $\spf Q_\ast \to \spf P_\ast$ is coisotropic if the unit map $P_\ast \to Q_\ast$ is surjective, and its kernel is closed under the Poisson bracket. Suppose that the Poisson structure on $P_\ast$ arises from an $\E{n+1}$-ring $P$ with $\pi_\ast P \cong P_\ast$, and $Q_\ast$ arises as the homotopy groups of an $\E{n}$-$P$-algebra $Q$ with unit map $f: P \to Q$. Then, \cref{prop: nonunital} implies that $\fib(f)$ admits a nonunital $\E{n+1}$-algebra structure. We therefore find that the map $\spf Q_\ast \to \spf P_\ast$ is the inclusion of a coisotropic \textit{subvariety} precisely when the map $\pi_\ast P \to \pi_\ast Q$ is surjective (since its kernel is then $\pi_\ast \fib(f) \subseteq \pi_\ast P$, which is closed under the Poisson bracket).
\end{observe}
Using \cref{obs: what does it mean to be coisotropic} with $n=2$, one can show that the map of \cref{eq: langlands dual map coisotropic} is very nearly coisotropic: the only obstructions are given by a completion issue, and that the map \cref{eq: langlands dual map coisotropic} need not be a closed immersion (i.e., the unit map is $\ku^G_\ast(\Omega G) \to \ku^H_\ast(\Omega(G/H))$ need not be surjective). Keeping these in mind, we are led to the following.
\begin{observe}\label{cor: desired map is coisotropic}
    Let $H\subseteq G$ be a closed subgroup of a compact Lie group such that the $\ku_G$-linear dual $\widehat{\cf_G(\Omega G)^\vee}$ of $\cf_G(\Omega G)$ and the $\ku_H$-linear dual $\widehat{\cf_H(\Omega(G/H))^\vee}$ of $\cf_H(\Omega(G/H))$ are concentrated in even degrees\footnote{This is a very mild condition, which can be checked to hold in all examples in this article. However, it is somewhat subtle, in the sense that working equivariantly is crucial. See, e.g., \cref{ex: borel homology loops S2}.}.
    \cref{rmk: abstracted scenario} and \cref{lem: map on centers} together imply that the unit map $\widehat{\cf_G(\Omega G)^\vee} \to \widehat{\cf_H(\Omega(G/H))^\vee}$ exhibits $\widehat{\cf_H(\Omega(G/H))^\vee}$ as an $\E{2}$-$\widehat{\cf_G(\Omega G)^\vee}$-algebra, and so \cref{prop: nonunital} implies that the fiber of the unit map admits a nonunital $\E{3}$-algebra structure. 
    For instance, this means that if the unit map induces a surjection on homotopy (which is not common!), the map $\spf \pi_\ast \widehat{\cf_H(\Omega(G/H))^\vee} \to \spf \pi_\ast \widehat{\cf_G(\Omega G)^\vee}$ will be the inclusion of a coisotropic subvariety. 
\end{observe}
\begin{remark}
    Applying \cref{lem: map on centers} to \cref{prop: nonunital} is of course quite ``lossy'', in the sense that the $\E{n}$-algebra structure from \cref{lem: map on centers} is of a very specific kind. The setup of \cref{lem: map on centers} hopefully exhibits further special features which allows us to understand \cref{expect: lag to Gbeta} further.
\end{remark}
\begin{remark}\label{rmk: coisotropic correspondence}
    More generally, let $A\in \Alg_\E{n+1}(\cC)$, and let $B\in \Alg_\E{n}(\Mod_A(\cC))$ with unit map $f: A \to B$. Then one has the following diagram:
    \begin{equation}\label{eq: hochschild square lag corr}
        \xymatrix{
           & \dZ_\E{n-1}(B/A) & & \\
            B \ar[ur] & & \dZ_\E{n}(f) \ar[ul] & \\
            & \dZ_\E{n}(B) \ar[ul] \ar[ur] & & \dZ_\E{n}(A), \ar[ul]
        }
    \end{equation}
    where, under some mild finiteness conditions, the square exhibits $\dZ_\E{n-1}(B/A)$ as the tensor product $B \otimes_{\dZ_\E{n}(B)} \dZ_\E{n}(f)$. The span 
    \begin{equation}\label{eq: centers lag corr}
        \xymatrix{
           & \dZ_\E{n}(f) & \\
            \dZ_\E{n}(B) \ar[ur] & & \dZ_\E{n}(A) \ar[ul] 
        }
    \end{equation}
    exhibits $\dZ_\E{n}(f)$ as an $\E{n}$-$\dZ_\E{n}(B) \otimes \dZ_\E{n}(A)$-algebra, so that as in \cref{obs: what does it mean to be coisotropic}, the induced diagram
    $$\xymatrix{
    & \spf \pi_\ast \dZ_\E{n}(f) \ar[dr] \ar[dl] & \\
    \spf \pi_\ast \dZ_\E{n}(B) & & \spf \pi_\ast \dZ_\E{n}(A)
    }$$
    should be viewed as a coisotropic correspondence. (It would be a coisotropic correspondence if the map $\pi_\ast \dZ_\E{n}(B) \otimes \pi_\ast \dZ_\E{n}(A) \to \pi_\ast \dZ_\E{n}(f)$ were surjective.) It would be interesting to know conditions on $A$ and $B$ which guarantee that $\spf \pi_\ast \dZ_\E{n}(A)$ and $\spf \pi_\ast \dZ_\E{n}(B)$ are (formal) \textit{symplectic} schemes with symplectic form of weight $n$, and then whether the above span is a \textit{Lagrangian} correspondence.
\end{remark}
\begin{example}
    When $A = \ku_G$ and $B = \ku_H$, \cref{eq: hochschild square lag corr} becomes the diagram
    $$\xymatrix{
        & \widehat{\cf_H(\Omega(G/H))^\vee} & & \\
        \ku_H \ar[ur] & & \widehat{\cf_H(\Omega G)^\vee} \ar[ul] & \\
        & \widehat{\cf_H(\Omega H)^\vee} \ar[ul] \ar[ur] & & \widehat{\cf_G(\Omega G)^\vee} \ar[ul]
    }$$
    which, upon base-changing along $\ku \to \Z$ (i.e., replacing $\ku_G$ everywhere by $C^\ast_G(\ast; \Z)$, etc.), is a completion of the ring of functions on the diagram in \cref{rmk: lag corr in homogeneous case}. 
    
    Similarly, suppose $A \to B$ is a map of commutative $\QQ$-algebras, and let $X = \spec(A)$ and $Y = \spec(B)$ denote the associated affine schemes, so that there is a map $Y \to X$. Suppose, also, that $n\geq 2$. Then, the Hochschild-Kostant-Rosenberg theorem implies that upon applying $\spf$ to \cref{eq: hochschild square lag corr}, we obtain the completion of the diagram
    $$\xymatrix{
    & T^\ast[n-1](Y/X) \ar[dr] \ar[dl] & & \\
    Y \ar[dr] & & T^\ast[n](X) \times_X Y \ar[dr] \ar[dl] & \\
    & T^\ast[n](Y) & & T^\ast[n](X)
    }$$
    at the respective zero sections. 
    Here, as usual, $T^\ast[n](X) = \spec \sh^{1/2} \Sym_{\co_X}(T_X(-n))$. 
    The subdiagram
    $$\xymatrix{
    & T^\ast[n](X) \times_X Y \ar[dr] \ar[dl] & \\
    T^\ast[n](Y) & & T^\ast[n](X)
    }$$
    is a decompletion of \cref{eq: centers lag corr}. As discussed in \cref{rmk: coisotropic correspondence}, this span is a coisotropic correspondence for essentially homotopy-theoretic reasons; if $X$ and $Y$ are smooth over $\QQ$, it is even a ($n$-shifted) Lagrangian correspondence (see, e.g., \cite{calaque-shifted-cotangent}).
\end{example}
\begin{remark}
    The coisotropicity of \cref{cor: desired map is coisotropic} can be understood from the perspective of boundary theories for topological quantum field theories (TQFTs for short). Namely, consider an $n$-dimensional TQFT $Z: \mathrm{Bord}_n \to \Cat_{(\infty,n-1)}$ (valued in, e.g., some $(\infty,n)$-category of ``presentable'' $\cc$-linear $(\infty,n-1)$-categories). Then $Z(S^j)$ is naturally with the structure of an $\E{j+1}$-algebra object in the $(\infty,n-j)$-category of $\cc$-linear $(\infty,n-j-1)$-categories, since $S^j$ is an $\E{j+1}$-algebra in $\mathrm{Bord}_n$. Indeed, if $I$ is a set of $n$ points on $\RR^{j+1}$, the manifold $S^{j+1} - I$ can be viewed as a cobordism $\sqcup_{i=0}^n S^j \rightsquigarrow S^j$.
    
    A boundary theory $\cB$ for $Z$ is a natural transformation from the trivial $n$-dimensional TQFT (i.e., whose value on a point is the unit object) to $Z$, so that $\cB$ is an $(n-1)$-dimensional TQFT. In particular, $\cB(S^{n-2})$ is an $\E{n-1}$-algebra object of the $\E{n-1}$-monoidal $\infty$-category $Z(S^{n-2})$. The cell decomposition $S^{n-1} \cong D^{n-1} \sqcup_{S^{n-2}} D^{n-1}$ implies that $Z(S^{n-1})$ is the $\E{n}$-algebra (in $\cc$-modules) given by the $\E{n-1}$-center of $Z(S^{n-2})$. This implies that $\cB(S^{n-2})$ is an $\E{n-1}$-$Z(S^{n-1})$-algebra.
    
    A special case of the above situation is when $Z$ is determined by an $\E{n}$-algebra $A$ in $\Mod_\cc(\Sp)$, in which case $\cB$ is determined by an $\E{n}$-$A$-algebra $B$. The statement that $\cB(S^{n-2})$ is an $\E{n-1}$-$Z(S^{n-1})$-algebra then translates to the statement that $\dZ_\E{n-2}(B/A)$ is an $\E{n-1}$-$\dZ_\E{n-1}(A/\cc)$-algebra, which is precisely \cref{lem: map on centers}.

    Returning to the general case, the structure of an $\E{n-1}$-$Z(S^{n-1})$-algebra on $\cB(S^{n-2})$ implies that there is a map $\spec \pi_\ast \cB(S^{n-2}) \to \spec \pi_\ast Z(S^{n-1})$, at least if $n \geq 3$. \cref{prop: nonunital} (and the general discussion in \cref{cor: desired map is coisotropic}) implies that this map should be coisotropic (in an appropriate sense) for the Poisson bracket of weight $n-1$ on $\pi_\ast Z(S^{n-1})$. It will be too much of a digression to discuss this here, but this coisotropicity can in fact be deduced from the secondary product structure (as discussed in \cite{ben-zvi-susy}) on local operators in extended TQFTs arising via topological descent.
\end{remark}

%% file: ku-hamiltonian/S1-equiv-MU.tex
\subsection{Equivariant spectra, graded groups, and duality over $\MU$}\label{subsec: duality over MU expectations}

Motivated by the discussion of the previous sections, we suggest that there is a generalization of the Langlands duality equivalences of \cref{thm: ku derived satake} and \cref{thm: rk 1 bzsv is true} to ``coefficients in complex cobordism''. It might be that the story suggested below is not the ``correct'' one, but nevertheless it leads to some interesting questions that we expect to play an important role in the future. In order to motivate our presentation, we will need to recall some results of Hausmann regarding equivariant complex cobordism.

In order to state results analogous to \cref{thm: ku derived satake} and \cref{thm: rk 1 bzsv is true}, we need to incorporate genuine equivariance into this picture. The primary motivation for this discussion is the work of Hausmann (see \cite{hausmann-global}).
\begin{notation}
    In what follows, we will write $\Sp_\glob$ to denote the $\infty$-category of global spectra as defined in \cite{schwede-global}, so that each compact Lie group $G$ defines a symmetric monoidal restriction functor $\Sp_\glob \to \Sp_G$ to the $\infty$-category of genuine $G$-equivariant global spectra. These functors are jointly conservative over all $G$. Note that \cite{schwede-global} works in the model-categorical setting, but we will work with the corresponding $\infty$-categories; the reader is referred to \cite{global-lax} for a discussion of global spectra in $\infty$-categorical language. We will not really need this theory below, but it is useful to review the context in which global complex cobordism arose.
\end{notation}
\begin{construction}
    Let $\cA$ denote the family of abelian compact Lie groups, let $\Orb^\cA$ denote the $\infty$-category of \cref{def: orb for family of compact lie}, and let $\Orb^\cA_\ast$ denote the (non-full) subcategory of pointed objects in $\Orb^\cA$. Then $\Orb^\cA_\ast$ is equivalent to the (nerve of the) topological category of abelian compact Lie groups. 

    Let $\ul{R}$ be an $\Eoo$-algebra in $\Sp_\glob$. Then there is a natural lax symmetric monoidal functor $\ul{\tau}_{\geq \star}: \Mod_{\ul{R}}(\Sp_\glob) \to \Fun(\Orb^{\cA,\op}_\ast, \Mod_{\tau_{\geq \star} R}^\fil)$ sending $\ul{A} \in \Sp_\glob$ to the functor $T \mapsto {\tau}_{\geq \star}^T(A_T)$.
    This constrution can be slightly modified as follows: let $\Lat$ denote the $1$-category of lattices, so that Pontryagin duality naturally gives a fully faithful functor $\Lat \subseteq \Orb^{\cA,\op}_\ast$. Restricting $\ul{\tau}_{\geq \star}$ along this inclusion defines a functor 
    $$\Mod_{\ul{R}}(\Sp_\glob) \to \Fun(\Lat, \Mod_{\tau_{\geq \star} R}^\fil).$$
    We will write $\ul{\pi}_\ast$ to denote the composite functor
    $$\Mod_{\ul{R}}(\Sp_\glob) \xar{\ul{\tau}_{\geq \star}} \Fun(\Orb^{\cA,\op}_\ast, \Mod_{\tau_{\geq \star} R}^\fil) \to \Fun(\Orb^{\cA,\op}_\ast, \Mod_{\pi_\ast R}^\gr).$$
\end{construction}
Let $\ul{\MU}$ denote the global complex cobordism spectrum, so that $\ul{\tau}_{\geq \star}(\ul{\MU})$ defines a functor $\Orb^{\cA,\op}_\ast \to \Mod_{\tau_{\geq \star} \MU}^\fil$. In fact, this refines to a functor $\Orb^{\cA,\op}_\ast \to \CAlg_{\tau_{\geq \star} \MU}^\fil$, which in particular defines a functor $\Orb^{\cA,\op}_\ast \to \CAlg_{\pi_{\ast} \MU}^\gr$.
Recall that the graded ring $\pi_\ast(\MU)$ classifies the universal graded ($1$-dimensional) formal group law by Quillen's \cite{quillen-formal-gps}. Hausmann proved a similar characterization of $\ul{\pi}_\ast(\ul{\MU})$, too.
\begin{definition}\label{def: graded group laws}
    A \textit{graded ($1$-dimensional) group law} is a functor $\co_\GG: \Lat \to \CAlg^{\heartsuit,\gr}$ equipped with an element $x\in \co_\GG(\Z)_{-2}$ in weight $-2$ such that for every lattice $\Lambda$ and split injective homomorphism $\chi: \Z \to \Lambda$, the sequence of graded abelian groups
    $$0 \to \co_\GG(\Lambda)_{\ast+2} \xar{\cdot \chi^\ast x} \co_\GG(\Lambda)_\ast \to \co_\GG(\Lambda/\Z\chi)_\ast \to 0$$
    is exact. Such an element $x$ will be called a \textit{coordinate}.
    Say that an $n$-tuple $x^{(1)}, \cdots, x^{(n)}$ of coordinates on a graded group law $\co_\GG$ is \textit{strict} if each $x^{(j)}$ is a multiple of $x^{(1)}$ by a unit $\lambda_j \in \co_\GG(\Z)_0$ whose restriction along the map $\co_\GG(\Z) \to \co_\GG(\{0\})$ is $1$.

\end{definition}
\begin{example}\label{ex: graded group coproduct}
    Graded group laws $\co_\GG: \Lat \to \CAlg^{\heartsuit,\gr}$ which preserve coproducts are simply specified by $\co_\GG(\{0\})$, $\co_\GG(\Z)$ viewed as a Hopf algebra over $\co_\GG(\{0\})$, and a regular element of $\co_\GG(\Z)_{-2}$ which generates the augmentation ideal. In particular, graded generalized groups which preserve coproducts are the same as $1$-dimensional linear graded algebraic groups over a graded commutative ring $R$ (namely $\co_\GG(\{0\})$) which are weight-connected over $R$.
\end{example}
The following result then gives a universal property for $\ul{\pi}_\ast(\ul{\MU})$, analogous to Quillen's theorem \cite{quillen-formal-gps} establishing $\pi_\ast(\MU)$ as the universal ring carrying a ($1$-dimensional) \textit{formal} group law. Note that the ``global'' perspective is crucial in even formulating the universal property below.
\begin{theorem}[{Hausmann, \cite[Theorem A]{hausmann-global}; Comeza\~na, \cite[XXVIII Theorem 5.3]{comezana}; L\"offler, \cite{loffler-MUG}}]\label{thm: hausmann}
    The ring $\ul{\pi}_\ast \ul{\MU}$ is concentrated in even weights, and $(\ul{\pi}_\ast \ul{\MU})(\Lambda)$ is a free $\pi_\ast(\MU)$-module for each lattice $\Lambda$.
    
    Moreover, for any $n \geq 0$, let $x_\tau^{(1)}, \cdots, x_\tau^{(n)}\in \ul{\pi}_\ast(\ul{\MU}^{\otimes n})$ denote the $n$ different tautological complex orientations of $\ul{\MU}$. If $(\co_\GG, x^{(1)}, \cdots, x^{(n)})$ is a strict $n$-tuple of coordinates on a graded group law, there is a unique homomorphism of graded group laws $(\ul{\pi}_\ast(\ul{\MU}^{\otimes n}), x_\tau^{(1)}, \cdots, x_\tau^{(n)}) \to (\co_\GG, x^{(1)}, \cdots, x^{(n)})$ which sends $x_\tau^{(j)} \mapsto x^{(j)}$.
\end{theorem}
\begin{construction}\label{cstr: Guniv}
    If $\Lambda$ is a lattice, let $\GG_\univ(\Lambda)$ denote the graded scheme $\spec \ul{\pi}_\ast(\ul{\MU})(\Lambda)$. 
    If $M$ is a finite $T$-equivariant spectrum, let $\cf_{T,M; \MU}$ denote the graded quasicoherent sheaf over $\GG_\univ(\bX^\ast(T))$ specified by $\pi_\ast^T(M^\vee \otimes \MU_T)$. This defines a functor $(\Sp_T^\fin)^\op \to \Perf^\gr(\GG_\univ(\bX^\ast(T)))$, which should be viewed as the analogue of the functor $\MU^\ast(-): (\Sp^\fin)^\op \to \Perf^\gr_{\pi_\ast \MU}$. If $M$ is the suspension spectrum of a finite $T$-space $X$, let $\cf_T(X; \MU) = \cf_{T,M; \MU}$. If $X$ is a $T$-space, let $\cf_T(X; \MU)^\vee$ denote the $\co_{\GG_\univ(\bX^\ast(T))}$-linear dual of $\cf_T(X; \MU)$, and if $X = \colim_{i\in \cI} X_i$ is an ind-finite $T$-space, let $\cf_T(X; \MU)^\vee = \colim_{i\in \cI} \cf_T(X_i; \MU)^\vee$.
\end{construction}

Since the arguments establishing \cref{thm: classical homology loops G} and \cref{thm: ku homology LG and langlands mod center} ultimately reduce to the case of \textit{torus}-equivariant (co)homology, if one is to generalize these results to equivariant $\MU$, it is first natural to ask whether equivariant $\MU$ admits abelian descent. Perhaps the most na\"ive formulation of this question is the following: if $G$ is a connected Lie group whose $\pi_1$ is torsion-free, $T\subseteq G$ is a maximal torus, and $W$ is its Weyl group, does restriction induce an isomorphism $\pi_\ast^G(\MU_G) \xar{\sim} \pi_\ast^T(\MU_T)^W$? This turns out to be \textit{false}:
\begin{lemma}[{\cite[Remark 1.2]{schwede-chern-MUG}}]\label{lem: equiv MU res not injective}
    The restriction map $\pi_\ast^{\U(n)}(\MU_{\U(n)}) \to \pi_\ast^T(\MU_T)$ is not injective if $n>1$.
\end{lemma}
\begin{remark}
    It is also not known whether $\pi_\ast^G(\MU_G)$ is concentrated in even degrees for a general compact Lie group $G$. This is a famous open problem in equivariant algebraic topology; see \cite{uribe-icm}.
\end{remark}
\begin{remark}\label{rmk: what is MUG really}
    \cref{lem: equiv MU res not injective} leads to the issue of whether an appropriate analogue of \cref{thm: classical homology loops G} and \cref{thm: ku homology LG and langlands mod center} should use $G$-equivariant $\MU$, or instead a more limited analogue which is built from $A$-equivariant $\MU$ over all compact abelian subgroups $A$ of $G$. 
    The resolution of this issue in the context of Langlands duality is not entirely clear to me, although it seems that the latter analogue of equivariant $\MU$ should be more relevant.
\end{remark}

Naturally, one is interested in proving a ``universal'' analogue of \cref{thm: derived satake} and \cref{thm: ku derived satake}. We will make some speculations about the form of such an equivalence, and discuss an actual mathematical statement relating to these speculations in \cref{subsec: all the expectations}. Far too many components of this discussion do not have well-behaved foundations at the moment, and so we will label everything below as a series of expectations. It is quite likely that these expectations are too na\"ive, and that more refinement is needed to make them ``correct''. I apologize in advance for the speculative nature of this discussion!
\begin{expect}\label{expect: shvsyn G equiv}
    Let $G$ be a topological group, and suppose that $X$ is a stratified finite space with $G$-action respecting the stratification. Then, there should be a $\F^\star_\ev (\MU_G)$-linear $\infty$-category $\Shv_G^c(X; \MU)^\Syn$ of equivariant ``synthetic'' constructible sheaves of $\MU$-modules on $X$. This theory should admit a well-behaved six functor formalism, and should also extend to ind-finite stratified $G$-spaces. (Note that the question of constructing such an $\infty$-category is already interesting when $X$ is a point, since we are then asking for a synthetic analogue of the $\infty$-category $\Perf_{\MU_G}(\Sp_G)$!)
    
    Changing coefficients of the underlying $\infty$-category $\Shv_G^c(X; \MU)$ along the unit map $\MU \to \Z$ should produce the $\infty$-category $\Shv_G^c(X; \Z)$. Furthermore, when $G = T$ is a torus, the corresponding graded categories $\Shv_G^c(X; \MU)^\gr$ should satisfy the property that there is an equivalence $\Shv_T^c(\ast; \MU)^\gr \simeq \Perf^\gr(\GG_\univ(\bX^\ast(T)))$, and the functor
    $$\Shv_T^c(X; \MU)^\gr \to \Shv_T^c(\ast; \MU)^\gr \simeq \Perf^\gr(\GG_\univ(\bX^\ast(T)))$$
    given by pushing forward along $X \to \ast$ should send the constant sheaf to the graded quasicoherent sheaf $\cf_T(X; \MU)$ on $\GG_\univ(\bX^\ast(T))$.
\end{expect}
For a sufficiently robust theory as in \cref{expect: shvsyn G equiv}, it should be possible to define the $\infty$-category $\Shv_{G\pw{t}/Z(G)\pw{t}}^{c,\gr}(\Gr_G; \MU)$. In order to obtain a spectral decomposition of this $\infty$-category analogous to \cref{thm: derived satake} and \cref{thm: ku derived satake}, we need an analogue of the quotient stack $G_\beta/G$ from \cref{def: Hbeta defn}.
\begin{expect}\label{expect: univ probing}
    Let $\GG: \Lat^\op \to \Aff^{\heartsuit, \gr}$ be a graded generalized group law, and let $H$ be a graded group scheme over $\GG(\{0\})$. Then there should be a graded stack $H_\GG/H$ over the classifying stack $B_{\GG(\{0\})} H$ satisfying the following properties. First, when $H$ is a torus $T$ (viewed as a constant group scheme over $\GG(\{0\})$), one should have $H_\GG = \GG(\bX^\ast(T))$. Second, if $\GG$ preserves products, so that it can be identified with the data of the weight-connected $1$-dimensional linear graded algebraic group $\GG(\Z)$ over $\co_\GG(\{0\})$ by \cref{ex: graded group coproduct}, there should be an isomorphism
    \begin{equation}\label{eq: HGmod H}
        H_\GG/H \cong \Map(\ul{\Hom}(\GG(\Z), B\GG_m), BH).
    \end{equation}
    Note that $\ul{\Hom}(\GG(\Z), B\GG_m)$ is the shifted Cartier dual of $\GG(\Z)$. Let us remark that it is not obvious how one might define the Cartier dual of a graded group $\GG$ which does not necessarily preserve products: the $\co_{\GG(\{0\})}$-linear dual of $\co_{\GG(\Z)}$ does \textit{not} necessarily admit a ring structure.
\end{expect}
\begin{remark}
    Suppose $\GG: \Lat^\op \to \Aff^{\heartsuit, \gr}$ is a graded generalized group law which preserves products, and let $H$ be a graded group scheme over $\GG(\{0\})$. The condition on $\GG$ means that it is fully determined by $\GG(\{0\})$, $\GG(\Z)$, and a coordinate on $\GG(\Z)$ by \cref{ex: graded group coproduct}. Let $H_\GG$ denote the fiber product $H_\GG/H \times_{B_{\GG(\{0\})} H} \GG(\{0\})$, where $H_\GG/H$ is defined as in \cref{eq: HGmod H}. Then $H_\GG$ can be viewed as a ``probing'' of $H$ by $\GG$. If $H$ is semisimple and base-changed from $\Z$, the construction $H \rightsquigarrow H_\GG$ ``replaces'' the Cartan $T$ of $H$ by $T_\GG$, and leaves the unipotent part alone. For instance, if $\GG = \GG_a(2)$, we can identify $H_\GG$ with $\fr{h}(2)$; if $\GG = \GG_m$, we can identify $H_\GG$ with $H$; and if $\GG = \GG_\beta$, we can identify $H_\GG$ with $H_\beta$.
\end{remark}
\begin{expect}[Derived Satake over $\MU$]\label{expect: sphere satake}
    Let $G$ be a simply-laced simply-connected semisimple algebraic group, and invert the integer $|W|$ in  $\MU$.
    Let $\GG_\univ$ denote the graded group law of \cref{cstr: Guniv}, let $\ld{G}_{\GG_\univ}/\ld{G}$ denote the graded stack of \cref{expect: univ probing}, and let $\ld{G}_{\GG_\univ}\mmod \ld{G}$ denote its coarse moduli space. Then:
    \begin{itemize}
        \item $\pi^{G/Z(G)}_\ast (\MU_{G/Z(G)})$ is concentrated in even degrees, and there is an isomorphism
        $$\spec \pi^{G/Z(G)}_\ast (\MU_{G/Z(G)}) \cong \ld{G}_{\GG_\univ}\mmod \ld{G},$$
        as well as a ``Kostant slice''
        $$\kappa_\univ: \ld{G}_{\GG_\univ}\mmod \ld{G} \to \ld{G}_{\GG_\univ}/\ld{G};$$
        \item there is a (ramified) $W$-cover 
        $$\spec \pi^{T/Z(G)}_\ast(\MU_{T/Z(G)}) \cong \ld{T}_{\GG_\univ} \to \ld{G}_{\GG_\univ}\mmod \ld{G},$$
        and if $\ld{T}_{\GG_\univ} \to \ld{G}_{\GG_\univ}/\ld{G}$ is the composite of this cover with $\kappa_\univ$, there is an isomorphism
        $$\spec \cf_{T/Z(G)}(\Omega G; \MU)^\vee \cong \ld{G}_{\GG_\univ}\mmod \ld{G} \times_{\ld{G}_{\GG_\univ}/\ld{G}} \ld{T}_{\GG_\univ};$$
        \item there is a $\pi^{G/Z(G)}_\ast (\MU_{G/Z(G)})$-linear equivalence
        $$\Shv_{G\pw{t}/Z(G)\pw{t}}^{c,\gr}(\Gr_G; \MU) \simeq \Perf^\gr(\ld{G}_{\GG_\univ}/\ld{G})$$
        such that changing coefficients along the unit map $\MU \to \QQ$ produces the derived geometric Satake equivalence.
        \item the above equivalence fits into a commutative diagram
        $$\xymatrix{
        \Shv_{G\pw{t}/Z(G)\pw{t}}^{c,\gr}(\Gr_G; \MU) \ar[r]^-\sim \ar[d]_-{\text{pushforward}} & \Perf^\gr(\ld{G}_{\GG_\univ}/\ld{G}) \ar[d]^-{\kappa_\univ^\ast} \\
        \Shv_{G\pw{t}/Z(G)\pw{t}}^{c,\gr}(\ast; \MU) \ar[r]_-\sim & \Perf^\gr(\ld{G}_{\GG_\univ}\mmod \ld{G}).
        }$$
    \end{itemize}
    The expected equivalence above can be viewed as describing (a categorification of) the $G/Z(G)$-equivariant stable homotopy type of $\Omega G$ in terms of the Langlands dual group, i.e., via the combinatorics of $G$.
\end{expect}
For instance, if a sheaf theory as in \cref{expect: shvsyn G equiv} exists, and $T$ is a torus, there will be an equivalence
$$\Shv_{T\pw{t}}^c(\Gr_T; \MU)^\gr \simeq \Shv_T^c(\Omega T; \MU)^\gr \simeq \Perf^\gr(\GG_\univ(\bX^\ast(T))/\ld{T}),$$
where the action on $\ld{T}$ on $\GG_\univ(\bX^\ast(T))$ is trivial. In other words, $\GG_\univ(\bX^\ast(T)) =: \ld{T}_{\GG_\univ}$ plays the role that the Lie algebra $\ld{\fr{t}}^\ast(2)$ played in \cref{thm: derived satake}.

\begin{remark}
    Suppose that the Kostant slice lifts to a map $\ld{G}_{\GG_\univ}\mmod \ld{G} \to \ld{G}_{\GG_\univ}$. When restricted to the subcategory of local systems, \cref{expect: sphere satake} should then give an equivalence between $\Loc_{G\pw{t}/Z(G)\pw{t}}^{\gr}(\Gr_G; \MU)$ and $\Perf^\gr(\ld{G}^\reg_{\GG_\univ}/\ld{G})$, where $\ld{G}^\reg_{\GG_\univ}$ is the $\ld{G}$-orbit of the slice $\ld{G}_{\GG_\univ}\mmod \ld{G} \to \ld{G}_{\GG_\univ}$. Forgetting about equivariance, one therefore expects an equivalence between $\Loc^{\gr}(\Gr_G; \MU)$ and the $\infty$-category of graded perfect complexes on $\ld{G}^\reg_{\GG_\univ}/\ld{G} \times_{\ld{G}_{\GG_\univ}\mmod \ld{G}} \{0\}$. 

    Since $\Loc(\Gr_G; \MU) \simeq \coMod_{\MU[\Omega G]}$, understanding $\Loc^{\gr}(\Gr_G; \MU)$ amounts to understanding $\pi_\ast \MU[\Omega G]$ as a bialgebra over $\pi_\ast \MU$. At least upon inverting $|W|$, one can identify $\pi_\ast \MU[\Omega G]$ with $\H_\ast(\Omega G; \Z') \otimes_{\Z'} \pi_\ast(\MU)$, at least as \textit{algebras}. This need not be an isomorphism of coalgebras. 
    
    Since $\spec \H_\ast(\Omega G; \Z') \cong Z_e(\ld{G}_{\Z'})$ is the centralizer of a regular nilpotent element in $\ld{\g}_{\Z'}$, one finds that if $\ld{\cN}^\reg$ denotes the regular locus in the nilpotent cone of $\ld{\g}_{\Z'}$, then $\ld{\cN}^\reg \cong \ld{G}_{\Z'}/\spec \H_\ast(\Omega G; \Z')$. So if $\pi_\ast \MU[\Omega G]$ is isomorphic to $\H_\ast(\Omega G; \Z') \otimes_{\Z'} \pi_\ast(\MU)$ as Hopf algebras, one is led to the expectation that there is an isomorphism over $\pi_\ast(\MU)$
    $$\ld{G}^\reg_{\GG_\univ}/\ld{G} \times_{\ld{G}_{\GG_\univ}\mmod \ld{G}} \{0\} \cong \ld{\cN}^\reg/\ld{G}.$$
    When $G = \SU(n)$, we will describe $\pi_\ast \MU[\Omega G]$ (with $\MU$ replaced by the sphere spectrum) in the next section. It turns out that this can be identified with $\H_\ast(\Omega G; \Z') \otimes_{\Z'} \pi_\ast(\MU)$ as Hopf algebras, and so $\spec \pi_\ast \MU[\Omega \SU(n)]$ is isomorphic to the centralizer $Z_e(\PGL_n)$ of a regular nilpotent element in $\fr{pgl}_n$. Therefore, $\Loc^{\gr}(\Gr_{\SL_n}; \MU)$ should be identified with $\Perf^\gr(\spec(\pi_\ast \MU) \times \ld{\cN}^\reg/\PGL_n)$.
\end{remark} 
\begin{remark}
    It would be interesting to generalize \cref{expect: sphere satake}, or even \cref{thm: ku derived satake}, to the non-simply-laced case.
\end{remark}
\begin{remark}\label{rmk: not even good for elliptic}
    The right-hand side of \cref{expect: sphere satake} is likely not literally correct as written, at least if one wants to generalize away from complex-oriented rings. To explain this, observe that if the construction of the left-hand side of \cref{expect: sphere satake} is sufficiently well-behaved, it should imply that if $A$ is any $2$-periodic $\Eoo$-ring equipped with an oriented group scheme $\GG$ over $A$ (with underlying group scheme $\GG_0$ over $\pi_0(A)$) and a commutative diagram
    $$\xymatrix{
    \GG \ar[d] \ar[r] & \spec(\MU_T^T) \ar[d] \\
    \spec(A) \ar[r] & \spec(\MU)
    }$$
    of $\Eoo$-rings, then there is an equivalence of $\pi_0(A_G)$-linear $\infty$-categories
    $$\Shv_{G\pw{t}/Z(G)\pw{t}}^{c,\gr}(\Gr_G; A) \simeq \Perf^\gr(\ld{G}_{\GG_0}/\ld{G}).$$
    Note that the right-hand side is indeed well-defined by \cref{eq: HGmod H}, since $\GG_0$ is an honest $1$-dimensional linear algebraic group. Such an equivalence would in particular imply an equivalence between the localizing subcategory of the left-hand side spanned by \textit{locally constant sheaves} (i.e., local systems) and the localization of the right-hand side given by restriction to the regular locus. This more limited equivalence was proved in \cite{grg-reg} in the cases of $2$-periodic rational cohomology and periodic complex K-theory. 
    
    However, the case of elliptic cohomology is more subtle. Suppose $\GG_0$ is an elliptic curve $E$ over $\pi_0(A)$, so that $\GG_0$ is \textit{not} affine. One could nevertheless try to define $\ld{G}_{\GG_0}/\ld{G}$ as $\Map(\ul{\Hom}(E, B\GG_m), B\ld{G})$, so that it is isomorphic to the moduli stack $\Bun_{\ld{G}}(E^\vee)$ of $\ld{G}$-bundles on the dual elliptic curve $E^\vee$. This, however, is \textit{not} the object appearing in \cite{grg-reg}: instead, the relevant object is the substack of degree zero semistable $\ld{G}$-bundles on $E^\vee$. This subtlety suggests that the correct form of \cref{expect: univ probing} and \cref{expect: sphere satake} should have similar adornments on the dual side, but it is difficult to guess what these should be.
\end{remark}
\begin{remark}
    It should also be possible to extend \cref{conj: bzsv} to the case of coefficients in $\MU$ (once \cref{expect: shvsyn G equiv} is appropriately resolved), although at this stage of development, there is very little that can be said about the general theory of relative Langlands duality with such coefficients. One might hope that $\ld{G}_{\GG_\univ}/\ld{G}$ admits a $1$-shifted symplectic structure. The analogue of Hamiltonian $\ld{G}$-spaces in this new context should then be given by Lagrangian morphisms to $\ld{G}_{\GG_\univ}/\ld{G}$.

    This expected generalization of \cref{conj: bzsv} is in some ways not entirely unreasonable, since (as mentioned at the end of \cref{expect: sphere satake}) it can be viewed as concerned with the $G$-equivariant stable homotopy type of $\cL(G/H)$, or equivalently with the $H$-equivariant stable homotopy type of $\Omega(G/H)$. The proof of \cref{thm: rk 1 bzsv is true}, and in particular \cref{table: topology and dualities for rank 1 spherical varieties}, shows that these spectral decompositions already have manifestations at the level of spaces themselves, and not just at the level of ($\ku$-)chains.
\end{remark}

%% file: ku-hamiltonian/langlands-over-sphere.tex
\subsection{A calculation over the sphere}\label{subsec: all the expectations}

Let us end with a calculation in the \textit{non}equivariant setting. One could interpret this rather simple calculation as either a step towards \cref{expect: sphere satake}, or as evidence that equivariance is the most interesting part of \cref{expect: sphere satake}. In any case, recall the following simple consequence (and historical antecedent) of \cref{thm: derived satake}:
\begin{prop}[{Ginzburg, \cite[Proposition 1.7.2]{ginzburg-langlands}}]\label{prop: ginzburg coh}
    Let $G$ be a connected semisimple algebraic group over $\cc$; we will abusively also write $G$ to denote the maximal compact subgroup of $G(\cc)$.
    Let $\ld{\g}^e(2)$ denote the centralizer of a principal nilpotent element $e$ in $\ld{\g}$. Then there is an isomorphism $\H^\ast(\Omega G; \cc) \cong U(\ld{\g}^e(2))$ of Hopf algebras over $\cc$, such that the Chern class $c_1(\det)$ of the determinant line bundle\footnote{This is simply defined to be the line bundle classified by the $\E{2}$-map $\Omega^2(BG \to K(\Z,4))$, where the map $BG \to K(\Z,4)$ detects the Killing form in $\H^4(BG; \Z) \cong \Sym^2(\bX^\ast(T))^W$.} over $\Omega G$ is sent to the element $e \in U(\ld{\g}^e(2))$. 
    
    Moreover, if $\lambda$ is a dominant weight of $\ld{G}$ with associated highest weight representation $\ld{V}_\lambda$, and $\IC_\lambda$ is the $\IC$-sheaf associated to $\Gr_G^{\leq \lambda} \subseteq \Gr_G$, there is an isomorphism $\H^\ast(\Gr_G; \IC_\lambda) \cong \ld{V}_\lambda$ of $\H^\ast(\Omega G; \cc) \cong U(\ld{\g}^e(2))$-comodules.
\end{prop}
\begin{remark}
    By \cite[Theorem 1.10.3]{ginzburg-perverse}, if $G$ is simply-connected and $\lambda,\mu$ are dominant coweights of $G$ with associated irreducible representations $\ld{V}_\lambda, \ld{V}_\mu$ of $\ld{G}$, there is an isomorphism
    \begin{equation}\label{eq: grg coh and homs}
        \Ext^\bull_{\Shv^c(\Gr_G; \cc)}(\IC_\lambda, \IC_\mu) \cong \Hom^\bull_{\Rep(Z_e(\ld{G}))}(\ld{V}_\lambda, \ld{V}_\mu).
    \end{equation}
    Indeed, using \cref{thm: full faithful}, taking cohomology defines an isomorphism
    $$\Ext^\bull_{\Shv^c(\Gr_G; \cc)}(\IC_\lambda, \IC_\mu) \xar{\sim} \Hom^\bull_{\H^\ast(\Omega G; \cc)}(\H^\ast(\Omega G; \IC_\lambda), \H^\ast(\Omega G; \IC_\mu)).$$
    The desired isomorphism \cref{eq: grg coh and homs} then follows from the graded isomorphism $\H^\ast(\Omega G; \IC_\lambda) \cong \ld{V}_\lambda$ via \cref{prop: ginzburg coh}.
    
    In \cite[Proposition 1.10.4]{ginzburg-perverse}, the isomorphism \cref{eq: grg coh and homs} is rephrased without appeal to a particular choice of regular nilpotent element as follows.
    Let $\ld{\cN}$ denote the nilpotent cone of $\ld{\g}$ (equipped with its natural $\GG_m$-action), and for any $\ld{G}$-representation $\ld{V}$, let $\ld{\cV}$ denote the associated vector bundle $\co_{\ld{\cN}} \otimes_\cc \ld{V}$. If $G$ is simply-connected and $\lambda,\mu$ are dominant coweights of $G$, there is an isomorphism
    $$\Ext^\bull_{\Shv^c(\Gr_G; \cc)}(\IC_\lambda, \IC_\mu) \cong \Hom^\bull_{\Coh(\ld{\cN}/\ld{G})}(\ld{\cV}_\lambda, \ld{\cV}_\mu).$$
    This can be deduced from \cref{eq: grg coh and homs} using the fact that the regular nilpotent orbit $\ld{\cN}^\reg \subseteq \ld{\cN}$ is isomorphic to $\ld{G}/Z_e(\ld{G})$, and has complement of codimension $\geq 2$.
\end{remark}
\begin{remark}
    The isomorphism $\H^\ast(\Omega G; \cc) \cong U(\ld{\g}^e(2))$ is defined via a map $\H^\ast(\Omega G; \cc) \to U(\ld{\g}(2))$ of Hopf algebras, which is constructed using \cref{thm: abelian satake}. Namely, taking cohomology defines a functor $\Perv_{G(\co)}(\Gr_G; \cc) \simeq \Rep(\ld{G}) \to \Mod_{\H^\ast(\Omega G; \cc)}$, and hence a map of Hopf algebras as desired.
\end{remark}
\begin{remark}\label{rmk: yun zhu centralizer}
    Let $\ell_G$ denote the square of the ratio of the lengths of long roots and the short roots of $G$.
    If we replace $\cc$ by $\Z[1/\ell_G]$, there is an isomorphism of Hopf algebras between $\H^\ast(\Omega G; \Z[1/\ell_G])$ and the divided power Hopf algebra $U^\sharp(\ld{\g}^e(2))$ (i.e., distributions on $Z_{\ld{G}}(e)$). In fact, there is an isomorphism $\spec \H_\ast(\Omega G; \Z[1/\ell_G]) \cong Z_e(\ld{G})$ of group schemes over $\Z[1/\ell_G]$; see \cite[Theorem 6.1]{homology-langlands}. After rationalization (or even just inverting $|W|$), this follows from \cref{thm: classical homology loops G}.

    Let us mention how this isomorphism can be deduced after inverting $|W|$ using the analogue of \cref{thm: derived satake} for sheaves with coefficients in $\Z[1/|W|]$. Namely, this equivalence states that there is an equivalence $\Loc_G(\Omega G; \Z[1/|W|]) \simeq \Perf(\ld{\g}^{\ast,\reg}[2]/\ld{G})$, so that 
    $$\Loc(\Omega G; \Z[1/|W|]) \simeq \Perf(\ld{\g}^{\ast,\reg}[2]/\ld{G} \times_{\ld{\g}^\ast[2]\mmod \ld{G}} \{0\}).$$
    But the fiber product $\ld{\g}^{\ast,\reg}[2] \times_{\ld{\g}^\ast[2]\mmod \ld{G}} \{0\}$ is precisely the regular locus $\ld{\cN}^\reg$ in the (shifted) nilpotent cone $\ld{\cN}$, so that $\ld{\g}^{\ast,\reg}[2]/\ld{G} \times_{\ld{\g}^\ast[2]\mmod \ld{G}} \{0\} \cong \ld{\cN}^\reg/\ld{G}$. By \cite{kostant-lie-group-reps}, there is a unique regular $\ld{G}$-orbit in $\ld{\cN}$, and so $\ld{\cN}^\reg/\ld{G} \cong BZ_e(\ld{G})$. This implies that there is an equivalence
    $$\Loc(\Omega G; \Z[1/|W|]) \simeq \Perf(\sh^{1/2} BZ_e(\ld{G})).$$
    This equivalence sends the skyscraper sheaf at the basepoint of $\Omega G$ to the pushforward of the structure sheaf along the map $\spec(\Z[1/|W|]) \to BZ_e(\ld{G})$, and the constant sheaf on $\Omega G$ to the structure sheaf of $BZ_e(\ld{G})$. For instance, taking endomorphisms of the constant sheaf on the left-hand side, we get an isomorphism $\H_\ast(\Omega G; \Z[1/|W|]) \cong \co_{Z_e(\ld{G})}$, as expected from \cite[Theorem 6.1]{homology-langlands}.
    
    Moreover, as in \cref{rmk: Loops2 and BGa}, the above equivalence of $\infty$-categories implies (by computing endomorphisms of the skyscraper sheaf at the basepoint of $\Omega G$) that there is an isomorphism
    $$\H_\ast(\Omega^2 G; \Z[1/|W|]) \cong \H^\ast(\sh^{1/2} BZ_e(\ld{G}); \co).$$
    For instance, if $G = \SU(n)$ (in which case one does not need to invert $n!$), the centralizer $Z_e(\PGL_n)$ is isomorphic to the group scheme $\WW_{n-1}$ of length $n$ Witt vectors (see \cite[Example 4.1.8]{grg-reg}) where the $j$th ghost coordinate lives in weight $2j$: indeed, both group schemes can be identified with the group of matrices of the form
    $$Z_e(\PGL_n) \cong \left\{
    \begin{pmatrix}
        1 & x_1 & x_2 & x_3 & \cdots & x_{n-1} \\
        & 1 & x_1 & x_2 & \cdots & x_{n-2}\\
        & & \ddots & \vdots & \cdots & \vdots \\
        & & & 1 & x_1 & x_2 \\
        & & & & 1 & x_1 \\
        & & & & & 1
    \end{pmatrix}
    \right\} \cong \WW_{n-1}.$$
    There is then an isomorphism
    $$\H_\ast(\Omega^2 \SU(n); \Z) \cong \H^\ast(\sh^{1/2} B\WW_{n-1}; \co);$$
    moreover, this isomorphism is compatible as $n$ varies. After base-changing to $\FF_p$, this recovers \cite[Theorem A]{ravenel-homology-Loops2-SUn} and the main result of \cite{yamaguchi-double-loops-stiefel}.
    Note, also, that when $n = \infty$, we get the statement that $\H_\ast(\Omega^2 \SU; \Z)$ is isomorphic to $\H^\ast(\sh^{1/2} B\WW; \co)$.\footnote{This can be seen directly as follows (we will ignore gradings for simplicity). Bott periodicity gives an equivalence $\Omega^2 \SU \simeq \Omega \BU = \U$, so it suffices to compute $\H_\ast(\U; \Z)$. This in turn can be understood via the colimit over $n$ of $\H_\ast(\U(n); \Z)$; moreover, $\H_\ast(\U(n); \Z)$ is the cobar construction on $\H_\ast(\BU(n); \Z)$. However, it is easy to see that $\H^\ast(\BU(n); \Z) \cong \co_{\WW_n}$; this implies that $\H_\ast(\BU(n); \Z) \cong \co_{\WW_n^\vee}$, where $\WW_n^\vee$ is the Cartier dual of $\WW_n$. In particular, $\H_\ast(\U(n); \Z) \cong \H^\ast(B\WW_n^\vee; \co)$. In the limit, we find that $\H_\ast(\U; \Z) \cong \H^\ast(B\WW^\vee; \co)$; the identification with $\H^\ast(B\WW; \co)$ now follows from the Cartier self-duality of the (graded) Witt ring scheme (see, e.g., \cite[Remark C.2]{thh-xn} for a reference). One way to interpret this discussion is that the algebraic manifestation of the Bott periodicity equivalence $\Omega \SU \simeq \BU$ is that $\WW$ is Cartier self-dual, up to reversing the grading, as a graded commutative group scheme.}
\end{remark}
\begin{remark}
    The following discussion is a digression, and can be ignored by the uninterested reader.
    Recall that the regular nilpotent $e\in \ld{\g}$ is the Chern class of the determinant line bundle on $\Gr_G$.  It is therefore natural to wonder if there is an analogue of the determinant line bundle for arbitrary {spherical varieties} $X$, which reduces to the determinant line bundle on $\Gr_G$ in the group case for $G$. If $\Gr^X$ denotes the relative Grassmannian of \cite[Section 8.2]{bzsv}, there is a tautological map $\Gr^X \to \Gr_G$, so one can simply consider the pullback of the determinant line bundle to $\Gr^X$.
    
    However, this construction seems lacking in at least two ways. First, in the group case for $G$ (so $X = (G\times G)/G^\mathrm{diag}$), one recovers the tensor \textit{square} of the determinant line bundle on $\Gr_G$, and not the determinant line bundle itself. Second, in some sense, it misses the \textit{role} of the determinant line bundle: namely, it captures the first nontrivial integral cohomology class on $\Gr_G$, and hence the principal nilpotent $e\in \ld{\g}$, but the above line bundle on $\Gr^X$ does not capture anything ``new''.
    
    In order to understand the new phenomena which appear in the relative setting, suppose $X = G/H$ with $H\subseteq G$ being a connected reductive subgroup. Throughout this article, we have emphasized that the quotient $\Omega(G/H) \simeq \Gr_G/\Gr_H$ captures the homotopically interesting content of the local relative geometric Langlands conjectures, so it is natural to ask whether $\Omega(G/H)$ carries an analogue of the determinant line bundle. The answer to this particular question is ``no'' for a very na\"ive reason: if $G/H$ is $3$-connected, as is often the case (e.g., $\GL_n/\GL_{n-1}$ for $n>2$), then $\H^2(\Omega(G/H); \Z) = 0$, and so all complex line bundles over $\Omega(G/H)$ are trivial.
    
    However, if one relaxes the notion of a line bundle by allowing for \textit{categorification}, the answer seems to be ``yes''. Consider, for instance, the case of $\Spin_{2n+2}/\Spin_{2n+1} \simeq S^{2n+1}$, so that $\H^\ast(\Omega S^{2n+1}; \Z) \cong \Z\pdb{x}$ with $x$ in cohomological degree $2n$. As we have seen in \cref{thm: bzsv for Dn}, the class $x$ plays the role of the principal nilpotent $e \in \ld{\g}_X$ (just as in \cref{prop: ginzburg coh}).
    In any case, $x$ may be viewed as a map $\Omega S^{2n+1} \to K(\Z, 2n)$, and hence it classifies a nontrivial $(2n-2)$-$\GG_m$-gerbe over $\Omega S^{2n+1}$. When $n=1$, this is simply the determinant line bundle over $\Omega S^3 = \Gr_{\SL_2}$.
    
    For larger $n$, one can wonder whether this $(2n-2)$-$\GG_m$-gerbe arises from a vector bundle over $\Omega S^{2n+1}$, i.e., whether the map $\Omega S^{2n+1} \to K(\Z,2n)$ factors as a composite
    \begin{equation}\label{eq: Omega S2n+1 and KZ2n}
        \Omega S^{2n+1} \xar{f} \BSO(2n) \xar{c_n} K(\Z,2n).
    \end{equation}
    Regarding this map, we have the following result, which suggests that, in general, a geometric and representation-theoretic construction/meaning of such putative determinant $\GG_m$-gerbes on $\Omega(G/H)$ seems very interesting.
    \begin{lemma}
        Such a factorization \cref{eq: Omega S2n+1 and KZ2n} is impossible for $n\geq 2$.
    \end{lemma}
    \begin{proof}
        Let us begin with the cases $n\neq 1,2,4$. When restricted to the bottom cell $S^{2n} \subseteq \Omega S^{2n+1}$, the map $S^{2n} \to \BSO(2n)$ will classify a vector bundle $\xi$ over $S^{2n}$ whose Euler class is a generator of $\Z$. If such a $\xi$ exists, the attaching map $S^{4n-1} \to S^{2n}$ of the top cell of the Thom space of $\xi$ will be a solution to the Hopf invariant one problem; so $n$ must be $1,2,4$, giving a contradiction.

        We have already seen that the factorization \cref{eq: Omega S2n+1 and KZ2n} exists when $n=1$, so let us now look at the case $n=2$.
        Let $S^4 \to \BSU(2) = \HHP^\infty$ denote the inclusion of the bottom cell, so that the map $\BSU(2) \to \BSO(4)$ defines a map $S^4 \to \BSO(4)$ which induces the map $(0,1):\Z \to \Z^2$ on $\H^4(-;\Z)$. 
        One can identify $\pi_j(\BSO(4)) \cong \pi_{j-1}(\SO(3) \times S^3)$, and our map $S^4 \to \BSO(4)$ induces the map $(0,1): \Z \to \pi_3 \SO(3) \oplus \Z$ on $\pi_4$.
        To get the desired map $f$, we would need the attaching map of the $8$-skeleton of $\Omega S^5$ to be $2\nu \in \pi_7(S^4)$. However, this is false; see \cite{mathoverflow-OmegaS5}. In brief, the attaching map of the $8$-skeleton is given by the Whitehead bracket $[\iota_4, \iota_4] \in \pi_7(S^4)$, which was computed in \cite[Equation 5.8]{toda} to be $\pm (2\nu - \Sigma \nu')$. (Here, $\nu': S^3 \wedge S^3 \cong S^6 \to S^3$ is the ``Blakers-Massey'' map which generates $\pi_6(S^3)$. It is characterized by the property that the composite $S^3 \times S^3 \to S^6 \to S^3$ sends two unit quaternions $(x,y)$ to the commutator $xyx^{-1} y^{-1}$.) Since $\Sigma \nu'\neq 0\in \pi_7(S^4)$, the desired factorization cannot exist.
    
        The only remaining case is $n=4$. Here, again, an argument similar to that of \cite{mathoverflow-OmegaS5} works to show that the desired extension does not exist. Namely, let $S^8 \to \BSO(8)$ denote the vector bundle whose sphere bundle is the octonionic Hopf bundle $S^{15} \to S^8$. As in \cite[Page 133]{mimura-homotopy-low-rank}, one can identify $\pi_j \SO(8) \cong \pi_j(\SO(7) \times S^7)$; and our map $S^8 \to \BSO(8)$ induces the map $(0,1): \Z \to \pi_7(\SO(7)) \oplus \Z$ on $\pi_8$.
        In order to extend this to a map $\Omega S^9 \to \BSO(8)$, we would need the attaching map of the $16$-skeleton of $\Omega S^9$ to be $2\sigma \in \pi_{15} S^8$. However, the attaching map of the $16$-skeleton is given by the Whitehead bracket $[\iota_8, \iota_8] \in \pi_{15}(S^8)$, which was computed in \cite[Equation 5.17]{toda} to be $\pm (2\sigma - \Sigma \sigma')$. Here, $\sigma$ is the Hopf map and $\sigma' \in \pi_{14}(S^7)$ is an octonionic version of the Blakers-Massey map (defined, as in the preceding paragraph, using the commutator of two unit octonions). Since $\Sigma \sigma'\neq 0\in \pi_{15}(S^8)$, the desired factorization cannot exist.
    \end{proof}
\end{remark}
Our goal is to prove an analogue of \cref{prop: ginzburg coh} over the sphere spectrum in the simplest case of $G = \SL_n$, and (a piece) of the second part of \textit{loc. cit.} when $n=2$. Since our goal is to illustrate certain phenomena, as opposed to proving the most general statement, we will only stick to this simple case.

First, we need a construction describing a passage from homotopy theory to algebraic geometry. In previous sections, we implicitly used the construction taking in a commutative ring spectrum $R$ and producing the (graded) affine scheme $\spec \pi_\ast(R)$. This construction is well-behaved if $R$ is concentrated in even degrees, but generally not so otherwise; for example, if $R$ is the sphere spectrum, the scheme $\spec \pi_\ast(S^0)$ is an absolute nightmare, and so something like $\spec \pi_\ast(\Omega \SU(n)_+)$ would be even more complicated. 
One important insight suggested by chromatic homotopy theory is that the algebro-geometric object associated to a ring spectrum which is not even should instead be a (graded) \textit{stack}. This perspective is made precise (for $\Eoo$-rings) in \cite{even-filtr} and \cite{prismatization-even-filtr}. Let us briefly recall the relevant construction.
\begin{recall}
    An $\Eoo$-ring $A$ will be called \textit{even} if $\pi_\ast(A)$ is concentrated in even weights.
    If $R$ is an $\Eoo$-ring spectrum, let $\cM_R^\fil$ denote the filtered stack
    $$\cM_R^\fil := \colim_{R \to A} \spec \tau_{\geq \star}(A),$$
    where the colimit is taken over all $\Eoo$-ring maps $R \to A$ with $A$ being even. The structure sheaf of $\cM_R^\fil$ defines a filtered $\Eoo$-ring $\F^\star_\ev R := \lim_{R \to A} \tau_{\geq \star} A$. If $M$ is an $R$-module, let $\cf_M^\fil$ denote the quasicoherent sheaf over $\cM_R^\fil$ defined by $\lim_{R \to A} \tau_{\geq \star}(A \otimes_R M)$.
    Of course, if $R$ is already even, $\cM_R^\fil \cong \spec \tau_{\geq \star}(R)$. The unit map $S^0 \to R$ defines a map $\cM_{R}^\fil \to \cM_{S^0}^\fil$.
    
    A map $R \to A$ is called an \textit{eff cover} (for ``evenly faithfully flat'') if for every even $\Eoo$-$R$-algebra $B$, the base-change $A \otimes_R B$ is even and the map $\pi_\ast(A) \to \pi_\ast(A \otimes_R B)$ is faithfully flat. If $R$ admits an eff cover $R \to A$ by an even $\Eoo$-ring $A$, there is an equivalence
    $$\cM_R^\fil \cong \colim_{\Deltab^\op} \spec \tau_{\geq \star}(A^{\otimes_R \bull+1}).$$
    If $M$ is an $R$-module, $\cf_M^\fil$ can be identified with $\lim_{\Deltab} \tau_{\geq \star}(M \otimes_R A^{\otimes_R \bull+1})$.
    Each of these filtered stacks and quasicoherent sheaves defines graded stacks and quasicoherent sheaves, which we will simply denote by $\cM_R$ and $\cf_M$. For example,
    $$\cM_R = \colim_{R \to A} \spec \pi_\bull(A),$$
    where $\bull$ denotes the grading.
\end{recall}
\begin{example}\label{ex: even-filtr on sphere}
    The map $S^0 \to \MU$ is an eff cover: if $B$ is an even $\Eoo$-ring, it admits a complex orientation, and so $\MU \otimes B \cong B[\BU]$; but $\BU$ has even cells, so that $\pi_\ast(\MU \otimes B)$ is a free $\pi_\ast(B)$-module on classes in even weights, and hence is itself concentrated in even weights. Furthermore, results of Quillen, Araki, Landweber, and Novikov (see \cite{quillen-formal-gps, araki-typical-formal-gps, landweber-cooperations, novikov-cobordism}) identify $\cM_{S^0}$ with the moduli stack $\Msfg$ of graded ($1$-dimensional) spin formal groups in the sense of \cite{miller-comodules}. Explicitly, if $R$ is a commutative ring, an $R$-point of $\Msfg$ is the data of a line bundle $\cL$ over $\spec(R)$, a formal group $\Ghat$ over $\spec(R)$, and isomorphism $\omega_{\Ghat} \cong \cL^{\otimes 2}$. We will write $\Ghat_\univ$ to denote the universal spin formal group over $\Msfg$, and $\omega$ to denote the line bundle over $\Msfg$ given by $\omega_{\Ghat_\univ}$. Note that there is a canonical square root $\omega^{1/2}$; this plays the role of the ``weight $1$'' line bundle $\co(1)$ which appears throughout this article.
    
    Any spectrum $M$ therefore defines a quasicoherent sheaf $\cf_M$ on $\Msfg$. If $X$ is a space, let $\cf(X; S^0)$ denote the ind-coherent sheaf $\cf_{(S^0)^{X_+}}$ associated to the spherical cochains $(S^0)^{X_+}$, so that the pullback of $\cf(X; S^0)$ along the map $\spec \pi_\ast(\MU) \to \Msfg$ can be identified with $\pi_\ast \MU^{X_+}$. Note that the diagonal on $X$ equips $\cf(X; S^0)$ with the structure of an $\Eoo$-algebra in $\IndCoh(\Msfg)$. For instance, let $q: \Ghat_\univ \to \Msfg$ denote the universal spin formal group over $\Msfg$; then $\cf(\CP^\infty; S^0) \cong \co_{\Ghat_\univ}$. If $X$ is a finite space, let $\cf(X; S^0)^\vee$ denote the $\co_{\Msfg}$-linear dual of $\cf(X; S^0)$.
\end{example}

Let us now return to the calculation at hand. 
\begin{construction}
    Let $X$ be a scheme (or even a stack), let $\cL$ be a line bundle over $X$, and let $\cL^\times$ denote the associated $\GG_m$-torsor. Let $\WW$ denote the Witt ring scheme, so that it has a $\GG_m$-action where the $n$th ghost coordinate has weight $2n$. Let $\WW(\cL) = \WW \times^{\GG_m} \cL^\times$ denote the associated \textit{Witt group scheme} over $X$. Quotienting by the Verschiebung, one obtains the length $n$ Witt group scheme $\WW_n(\cL)$. Its Cartier dual will be denoted $\WW_n^\vee(\cL^{-1})$.
\end{construction}
Using that $\Omega \SU(n)$ has even cells, we may identify $\cM_{\Omega \SU(n)_+}$ with $\spec_{\Msfg} \cf(\Omega \SU(n); S^0)^\vee$, which is described in the following result. It can be viewed as a refinement of \cite[Example 4.1.8]{grg-reg}, and an analogue of \cref{prop: ginzburg coh} for $G = \SL_n$.
\begin{lemma}\label{lem: sphere and Loops SUn}
    Let $\Ghat_\univ$ denote the universal spin formal group over $\Msfg$, and let $\Ghat_\univ^\vee$ denote its Cartier dual.
    There is an isomorphism 
    $$\spec_{\Msfg} \cf(\Omega \SU(n); S^0)^\vee \cong \WW_{n-1}(\omega)$$
    of group stacks over $\Msfg$, and which induces an isomorphism 
    $$\spec_{\Msfg} \cf(\Omega \SU(n); S^0) \cong \WW_{n-1}^\vee(\omega^{-1}) \cong \Hom(\Ghat_\univ^\vee, \WW_{n-1}^\vee)$$
    of group stacks over $\Msfg$.
\end{lemma}
\begin{proof}
    The natural map $\CP^{n-1} \to \Omega \SU(n)$ is Bott's generating complex from \cite{bott-space-of-loops}, so that there is an isomorphism 
    $$\spec_{\Msfg} \cf(\Omega \SU(n); S^0)^\vee \cong \spec_{\Msfg} \Sym_{\Msfg}(\cf(\CP^{n-1}; S^0)^\vee),$$
    the latter being the total space $\bV(\cf(\CP^{n-1}; S^0))$ of the vector bundle $\cf(\CP^{n-1}; S^0)$.
    It therefore suffices to observe that there is an isomorphism of group schemes
    $$\bV(\cf(\CP^{n-1}; S^0)) \cong \WW_{n-1}(\omega).$$
    Since $\spec_{\Msfg} \cf(\Omega \SU(n); S^0)^\vee$ is Cartier dual to $\spec_{\Msfg} \cf(\Omega \SU(n); S^0)$, this implies that the latter is isomorphic to $\WW_{n-1}^\vee(\omega^{-1})$ (being Cartier dual to $\WW_{n-1}(\omega)$).
    \cref{lem: Lie and hom to Ga} implies that there is an isomorphism
    $$\Hom_{\Msfg}(\Ghat_\univ^\vee, \WW_{n-1}) \cong \WW_{n-1}(\omega^{-1}),$$
    since $\Lie(\Ghat_\univ) \cong \omega^{-1}$; this in turn implies that there is an isomorphism
    $$\Hom_{\Msfg}(\Ghat_\univ^\vee, \WW_{n-1}^\vee) \cong \WW_{n-1}^\vee(\omega^{-1}),$$
    as desired.
\end{proof}
\begin{remark}
    The group scheme $\WW_{n-1}(\omega)$ over $\Msfg$ can be identified with the fiber product $\WW_{n-1}(\co(2)) \times_{B\GG_m} \Msfg$, where the map $\Msfg \to B\GG_m$ classifies $\omega^{1/2}$.
\end{remark}
\begin{remark}
    The $\E{2}$-map $\Omega S^3 \to \CP^\infty$ classifying the determinant line bundle induces a map 
    $$\spec_{\Msfg} \cf(\CP^\infty; S^0)^\vee \to \spec_{\Msfg} \cf(\Omega S^3; S^0)^\vee,$$
    which, under \cref{lem: sphere and Loops SUn} for $n=2$ can be identified with the canonical map $\Ghat^\vee \to \bV(\omega)$ classifying the tautological section of $\Lie(\Ghat_\univ) \otimes \omega \cong \co_{\Msfg}$. The canonical map 
    $$\WW_1^\vee(\omega^{-1}) \cong \hat{\bV}(\omega^{-1})^\sharp \to \Ghat_\univ$$
    classifies the first Chern class of the determinant line bundle, and can be regarded as an analogue of the map $\spec \H^\ast(\Omega \SU(n); \Z) \to \AA^1$ detecting the principal nilpotent element of $\ld{\g} = \pgl_n$ under \cref{prop: ginzburg coh}.
\end{remark}
According to \cref{lem: sphere and Loops SUn}, $\WW_{n-1}(\omega)$ is to be understood as the analogue for the sphere of the centralizer $Z_{\PGL_n}(e)$.\footnote{Note that there is an isomorphism $Z_{\PGL_n}(e) \cong \WW_{n-1}(\co_{\spec(\Z)})$ of group schemes over $\Z$; see \cite[Example 4.1.8]{grg-reg}. Therefore, pulling back the isomorphism of \cref{lem: sphere and Loops SUn} along the map $\spec(\Z) \to \Msfg$ classifying the additive formal group precisely recovers the isomorphism $\spec \H_\ast(\Omega \SU(n); \Z) \cong Z_{\PGL_n}(e)$.}
In this setting, the grading on $Z_{\PGL_n}(e)$ translates into tensor powers of the line bundle $\omega$.
As a consequence of the above discussion, $\WW_{n-1}^\vee(\omega^{-1})$ is to be understood as the analogue for the sphere of the divided power enveloping algebra $U^\sharp(\fr{pgl}_n^e)$. 

Let us now illustrate an analogue of the second part of \cref{prop: ginzburg coh} for $G = \SL_2$.\footnote{To a representation theorist, this might seem like an especially trivial case, and likewise for a homotopy theorist; but perhaps for somewhat different reasons (either $\ld{G} = \PGL_2$ is too simple, or $\Omega G = \Omega S^3$ is too simple). For a geometric representation theorist, this example is trivial for two reasons!} In this case, we have the following (well-known) calculation.
\begin{lemma}\label{lem: schubert for GL2}
    Consider the dominant (co)weight $(i,j)$ of $\GL_2$, so that $i\geq j$. Then there is a cell decomposition
    $$\Gr_{\GL_2}^{\leq (i,j)} \cong \begin{cases}
        \coprod_{0\leq k \leq (i-j)/2} \Gr_{\GL_2}^{(i-k, j+k)} & i-j\text{ even}, \\
        \coprod_{0\leq k \leq (i-j-1)/2} \Gr_{\GL_2}^{(i-k, j+k)} & i-j\text{ odd}.
    \end{cases}$$
    If $i>j$, there is an isomorphism
    $$\Gr_{\GL_2}^{(i,j)} \cong \Map(\spec \cc[\epsilon]/\epsilon^{i-j}, \PP^1),$$
    and $\Gr_{\GL_2}^{(i,i)} = \spec(\cc)$.
\end{lemma}
\begin{proof}[Proof sketch]
    The only nontrivial claim is the calculation of $\Gr_{\GL_2}^{(i,j)}$ for $i>j$. Since there is an isomorphism $\Gr_{\GL_2}^{(i,j)} \cong \Gr_{\GL_2}^{(i-1,j-1)}$, we may assume that $j = 0$. In this case, the $\GL_2(\cc\pw{t})$-orbit of the loop $\begin{psmallmatrix}
        t^i & 0\\
        0 & t^j
    \end{psmallmatrix}$ consists of equivalence classes of matrices of the form $\begin{psmallmatrix}
        t^i a & b\\
        t^i c & d
    \end{psmallmatrix}$ with $a,b,c,d\in \cc\pw{t}$. Note that both $b$ and $d$ cannot both be nonconstant, so $\Gr_{\GL_2}^{(i,0)}$ is covered by the open loci where the constant term of $b$ (resp. of $d$) is invertible. Let $A = \begin{psmallmatrix}
        a & b \\
        c & d
    \end{psmallmatrix}$. Over the locus where the constant term of $b$ is invertible, we have the relation
    $$\tfrac{1}{\det(A)} \begin{psmallmatrix}
        t^i a & b\\
        t^i c & d
    \end{psmallmatrix} \begin{psmallmatrix}
        b & 0\\
        -at^i & adb^{-1} - c
    \end{psmallmatrix} = \begin{psmallmatrix}
        0 & 1\\
        -t^i & db^{-1}
    \end{psmallmatrix}.$$
    Multiplying on the right by some strictly upper-triangular matrix in $\GL_2(\cc\pw{t})$, we can further assume that the bottom-right corner lies in $\cc[t]/t^i$. This gives an isomorphism between $\{b(0) \neq 0\} \subseteq \Gr_{\GL_2}^{(i,0)}$ and $\cc[t]/t^i$. One can similarly obtain an isomorphism between $\{d(0) \neq 0\} \subseteq \Gr_{\GL_2}^{(i,0)}$ and $\cc[t]/t^i$ via matrices of the form $\begin{psmallmatrix}
        t^i & bd^{-1}\\
        0 & 1
    \end{psmallmatrix}$. The intersection of these two opens is precisely given by sending $db^{-1} \mapsto bd^{-1}$, i.e., is obtained by gluing $\cc[t]/t^i$ to itself along the inversion automorphism of $(\cc[t]/t^i)^\times$. But this is precisely the mapping scheme $\Map(\spec \cc[\epsilon]/\epsilon^{i-j}, \PP^1)$.
\end{proof}
\begin{corollary}\label{cor: partial james GrSL2}
    Let $i\geq 0$ be an integer.
    If $J_i(S^2) \hookrightarrow \Omega S^3$ denotes the partial James construction (with top cell in dimension $2i$), \cref{lem: schubert for GL2} gives a homotopy equivalence
    $$\Gr_{\SL_2}^{\leq i}(\cc) \cong \Gr_{\GL_2}^{\leq (i,-i)}(\cc) \simeq J_{2i}(S^2)$$
    which is compatible with the homotopy equivalence $\Gr_{\SL_2}(\cc) \simeq \Omega S^3$ from \cref{thm: quillen}.
\end{corollary}
Importantly, \cref{lem: schubert for GL2} implies that the Schubert varieties for $\SL_2$ are rationally smooth, and so $\IC_\lambda$ is simply the pushforward of the constant sheaf $\ul{\cc}$ on $\Gr_{\SL_2}^{\leq \lambda}$. Therefore, if $i\geq 0$ is an integer (viewed as a dominant coweight of $\SL_2$), \cref{cor: partial james GrSL2} implies that we may identify 
$$\H^\ast(\Omega \SU(2); \IC_i) \cong \H^\ast(J_{2i}(S^2); \cc),$$
and hence \cref{prop: ginzburg coh} says (in particular) that
\begin{equation}\label{eq: coh of james and Vi}
    \H^\ast(J_{2i}(S^2); \cc[-2i]) \cong \ld{V}_i
\end{equation}
with $\ld{V}_i$ being the $(2i+1)$-dimensional irreducible representation of $\PGL_2$.
To prove an analogue of the second half of \cref{prop: ginzburg coh}, we therefore need to compute $\cf(J_{2i}(S^2); S^0)$ as a $\cf(\Omega S^3; S^0)$-comodule (i.e., by \cref{lem: sphere and Loops SUn}, as a $\bV(\omega)$-representation). This is quite simple\footnote{Ha!}:
\begin{lemma}\label{lem: coh of J2i}
    Let $\Gamma^j_{\Msfg}(\omega) = (\omega^{\otimes j})^{\Sigma_j}$ denote the $j$th divided power of $\omega$, and let $\Gamma_{\Msfg}^{\leq 2i}(\omega) = \bigoplus_{0\leq j\leq 2i} \Gamma^j_{\Msfg}(\omega)$. Then there is an isomorphism 
    $$\cf(J_{2i}(S^2); S^0[-2i]) \cong \Gamma_{\Msfg}^{\leq 2i}(\omega) \otimes \omega^{\otimes -i},$$
    and the inclusion $J_{2i}(S^2) \to \Omega S^3$ induces the canonical map 
    $$\co_{\hat{\bV}(\omega^{-1})^\sharp} \cong \Gamma_{\Msfg}^\ast(\omega) \to \Gamma_{\Msfg}^{\leq 2i}(\omega).$$
\end{lemma}
Since the representation $\ld{V}_i$ of $\PGL_2$ (over $\Z$) can be identified with $\Gamma^{\leq 2i}(\Z\cdot e)$, \cref{lem: coh of J2i} can be viewed as an analogue of \cref{eq: coh of james and Vi} over the sphere spectrum. 
\begin{remark}
    Note that \cref{lem: coh of J2i} (or even \cref{eq: coh of james and Vi}) only describes the action of the centralizer $Z_{\PGL_2}(e)$, which is isomorphic to the unipotent radical of the upper-triangular Borel subgroup of $\PGL_2$, on the representation $\ld{V}_i$. However, it is well-known that the action of the entirety of $\PGL_2$ is completely encoded in natural structures present in this setting (this is a part of \cref{thm: abelian satake}). Namely, in the classical setting of \cref{prop: ginzburg coh}, the action of the Cartan subgroup $\GG_m \subseteq \PGL_2$ is encoded by the natural grading on $\H^\ast(\Omega \SU(2); \IC_i)$. Moreover, the action of the unipotent radical of the opposite Borel is encoded by the action of $Z_{\PGL_2}(e)$ under ``Poincar\'e duality'' on $J_{2i}(S^2)$. Although $J_{2i}(S^2)$ is generally not a smooth manifold (and is often not even a Poincar\'e duality complex unless $(2i)!$ is inverted in the ring of coefficients), its integral homology and cohomology \textit{groups} are dual to each other, with a shift of $4i$.

    Concretely, the abelian group underlying $\ld{V}_i$ is $\Z\{v_0, \cdots, v_{2i}\}$, where one should understand $v_j$ as $\frac{e^j}{j!}$ placed in weight $2j-2i$ (so $v_0$ is in weight $-2i$, and $v_{2i}$ is in weight $2i$). This specifies the action of the diagonal torus of $\PGL_2$. The action of $Z_{\PGL_2}(e) \cong \begin{psmallmatrix}
        1 & \ast\\
        0 & 1
    \end{psmallmatrix}$ can be described through the action of $U^\sharp(\fr{pgl}_2^e) = \Gamma_\Z(e)$: the element $\frac{e^n}{n!} \in U^\sharp(\fr{pgl}_2^e)$ sends $v_j \mapsto \binom{n+j}{n} v_{n+j}$.
    Finally, the action of the strictly lower-triangular matrices $\ld{N}^- = \begin{psmallmatrix}
        1 & 0\\
        \ast & 1
    \end{psmallmatrix} \subseteq \PGL_2$ is described by ``Poincar\'e duality'': this exchanges $v_j$ with $v_{2i-j}$, and therefore $\frac{f^n}{n!} \in U^\sharp(\fr{\ld{n}}^-) \cong \Gamma_\Z(f)$ sends $v_j \mapsto \binom{n+(2i-j)}{n} v_{j-n}$.
    
    The important point in the setting of \cref{lem: coh of J2i} is that the action of the Cartan subgroup is encoded by the tensor powers of the line bundle $\omega$ appearing in $\cf(J_{2i}(S^2); S^0[-2i])$. But as usual, the action of the unipotent radical of the opposite Borel is still encoded by ``Poincar\'e duality'' on $J_{2i}(S^2)$, which gives an isomorphism of $\co_{\Msfg}$-modules
    $$(\Gamma_{\Msfg}^{\leq 2i}(\omega) \otimes \omega^{\otimes -i})^\vee \cong \Gamma_{\Msfg}^{\leq 2i}(\omega) \otimes \omega^{\otimes i}.$$
\end{remark}
As a first step towards \cref{expect: sphere satake}, it would be interesting and important to prove an analogue of \cref{prop: ginzburg coh} with coefficients in the sphere spectrum for arbitrary (simply-connected and simply-laced) $G$, and not just for $\SL_2$.

%% file: appendices/shv-loops.tex
\section{Sheaves on loop spaces}

The proof of \cref{thm: ordinary homology criterion satake} relied on a comparison to positive characteristic. In the finite-type case, this comparison is provided by \cref{thm: changing coefficients and bases}, but the relevant comparison turns out to be more subtle in the infinite-type situation.

\begin{notation}
    Let $R$ be a commutative ring, and let $X$ be an $R$-scheme. Write $X\ls{t}$ to denote the prestack sending an $R$-algebra $S$ to $X(S\ls{t})$. Similarly, write $X\pw{t}$ to denote the prestack sending an $R$-algebra $S$ to $X(S\pw{t})$, and let $X\pw{t}/t^n$ denote the prestack sending an $R$-algebra $S$ to $X(S[t]/t^n)$.
\end{notation}
\begin{definition}\label{def: collected action}
    Let $X$ be an affine scheme defined over a commutative ring $R$ equipped with an action of a linear algebraic group $G$ over $R$. The $G\pw{t}$-action on $X\ls{t}$ is called \textit{placid} if:
    \begin{itemize}
        \item there is a presentation $X\ls{t} = \varinjlim_j X^j$, where each $X^j$ is an inverse limit $\varprojlim_n X^j_n$ with each $X^j_n$ being a $G\pw{t}$-equivariant scheme of finite type;
        \item the action of $G\pw{t}$ on $X^j_n$ factors through $G\pw{t}/t^{m_n}$ for some $m_n \gg 0$ (compatibly in $n$).
        \item The transition maps $X^j_n \to X^j_{n'}$ are $G\pw{t}/t^{m_n}$-equivariant affine smooth surjections.
    \end{itemize}
    Note that this is weaker than the condition that $G\pw{t}$-action on $X\ls{t}$ being placid in the sense of \cite[Section 7.3.1]{bzsv}: there, it is required that the transition maps $X^j_n \to X^j_{n-1}$ also be torsors for a unipotent group scheme.
\end{definition}
\begin{construction}\label{cstr: placid sheaves}
    In the above setup, let $\Shv^\et_{G\pw{t}}(X^j_n; \olQell)$ denote the $\infty$-category $\Shv^\et_{G\pw{t}/t^{m_n}}(X^j_n; \olQell)$. Note that since the kernel of the surjection $G\pw{t} \twoheadrightarrow G\pw{t}/t^{m_n}$ is unipotent, and the action of $G\pw{t}$ factors through this surjection, the $\infty$-category $\Shv^\et_{G\pw{t}/t^{m_n}}(X^j_n; \olQell)$ would be unchanged if we replace $m_n$ by any $m\geq m_n$.

    Let $\Shv^{\et}_{G\pw{t}}(X^j; \olQell)$ denote the direct limit
    $$\Shv^\et_{G\pw{t}}(X^j; \olQell) = \varinjlim_{f_{j,n}^!} \Shv^\et_{G\pw{t}}(X^j_n; \olQell)$$
    of the $\infty$-categories $\Shv^\et_{G\pw{t}}(X^j_n; \olQell)$ along $!$-pullbacks. Finally, define
    $$\Shv^\et_{G\pw{t}}(X\ls{t}; \olQell) = \varinjlim_{g^j_!} \Shv^\et_{G\pw{t}}(X^j; \olQell)$$
    to be the direct limit of the $\infty$-categories $\Shv^\et_{G\pw{t}}(X^j; \olQell)$ along the $!$-pushforward functors associated to the maps $g^j: X^j \to X^{j+1}$.
    
    Suppose that there are only countably many $G\pw{t}$-orbits on $X\ls{t}$. (If $X$ is affine and $G$ is reductive, \cref{thm: spherical gaitsgory nadler} says that this is the case if and only if $X$ is a spherical $G$-variety.) Then there are only finitely many $G\pw{t}/t^{m_n}$-orbits on $X^j_n$, and the maps $f_{j,n}: X^j_n \to X^j_{n-1}$ are $G\pw{t}$-equivariant and respect the stratifications on $X^j_n$ and $X^j_{n-1}$.
    Define $\Shv^{c,\et}_{G\pw{t}}(X\ls{t}; \olQell)$ to be the $\infty$-category obtained via the above procedure, except where $\Shv^\et_{G\pw{t}/t^{m_n}}(X^j_n; \olQell)$ is replaced by the $\infty$-category $\Shv^{c,\et}_{G\pw{t}/t^{m_n}}(X^j_n; \olQell)$ of $G\pw{t}/t^{m_n}$-equivariant \'etale sheaves on $X^j_n$ which are constructible with respect to the orbit stratification on $X^j_n$.
\end{construction}
\begin{theorem}\label{thm: changing coefficients and bases}
    Let $q \gg 0$ be a large prime power. Fix a prime $\ell \neq p$, and choose an isomorphism $\iota: \olQell \xar{\sim} \cc$.
    In \cref{setup: changing coefficients and bases}, suppose that the $H(R\pw{t})$-action on $X\ls{t}$ is placid. Then is there is a localization $R \subseteq R'$ such that for any $k$-point $R' \twoheadrightarrow k$ with $k$ being a finite field, there is a natural equivalence
    $$\Shv^{c,\et}_{H_{\ol{\FF}_q}\pw{t}}(X\ls{t}_{\ol{\FF}_q}; \ol{\QQ_\ell}) \xar{\sim} \Shv^{c}_{H(\cc\pw{t})}(X(\cc\ls{t}); \cc).$$
\end{theorem}
\begin{proof}
    By definition of the $\infty$-categories involved, it suffices to show that for each $j$ and $n$, there are compatible equivalences
    $$\Shv^{c,\et}_{H_{\ol{\FF}_q}\pw{t}/t^{m_n}}(X^j_{n,\ol{\FF}_q}; \ol{\QQ_\ell}) \xar{\sim} \Shv^{c}_{H_\cc\pw{t}/t^{m_n}}(X^j_{n,\cc}; \cc).$$
    This in turn is a consequence of \cref{thm: finite type changing coefficients and bases} applied to the group scheme $H\pw{t}/t^{m_n}$ over $R$.
\end{proof}
\begin{remark}
    \cref{thm: changing coefficients and bases} implicitly uses \cref{rmk: coeff change infty cat}: it is necessary to treat the equivalence in the proof of \cref{thm: changing coefficients and bases} as those of $\infty$-categories, so that the (co)limit constructions of \cref{cstr: placid sheaves} are legal.
\end{remark}
We will often use \cref{thm: changing coefficients and bases} in the case when $X$ is an affine homogeneous spherical $H$-variety. For this, we need:
\begin{conjecture}\label{conj: placidity}
    If $X$ is an affine homogeneous $H$-variety, the $H\pw{t}$-action on $X\ls{t}$ is placid (so that we can apply \cref{thm: changing coefficients and bases}).
\end{conjecture}
However, not being a specialist in the relevant technical details, I have not been able to verify this. It is quite likely that the assumption of placidity is not necessary to prove \cref{thm: ordinary homology criterion satake}, since it is only used to conclude formality of a certain $\Ext$-algebra. 
\begin{remark}
    When $X = \GL_n/\O_n$ and $X = \GL_{2n}/\Sp_{2n}$, \cref{conj: placidity} was verified as \cite[Proposition 22]{chen-yi-formality} (and the desired formality mentioned above was deduced as \cite[Theorem 23]{chen-yi-formality}). The stratification on $X\ls{t}$ defined in \cite{chen-yi-formality} induces one on the based loop space $\Omega X$, where it has appeared previously in \cite[Discussion after Proposition 1.4]{crabb-mitchell}.
\end{remark}

%% file: appendices/comparison-cohen-jones-yan.tex
\section{String topology for spherical homogeneous spaces}

\cref{prop: equiv coh as hochschild} allows us to use \cref{conj: bzsv} (rather, an integral refinement of it) to give a description of the (integral) homology of the free loop space of a compact homogeneous space equipped with the Chas-Sullivan product in Langlands dual terms. This section is mainly an exercise in taking derived quotients in simple examples.
\begin{assume}\label{assume: split form}
    There is a sufficiently large localization $\Z'$ of $\Z$ such that there is a Kostant slice $\kappa: \ld{\g}_{\Z'}^\ast(2)\mmod \ld{G}_{\Z'} \to \ld{\g}_{\Z'}^\ast(2-2\rho)$, and the equivalence of \cref{conj: bzsv} satisfies the following properties:
    \begin{enumerate}
        \item It admits a refinement to a $\Z'$-linear equivalence
        $$\Shv^{c,\Sat}_{G\pw{t}}(X\ls{t}; \Z') \simeq \Perf(\sh^{1/2} \ld{M}_{\Z'}/\ld{G}_{\Z'}(-2\rho)),$$
        where $\ld{M}_{\Z'}$ is a flat lift of $\ld{M}$ along $\Z' \to \QQ$, and $\ld{G}_{\Z'}$ is the Chevalley split form of $\ld{G}$. 
        \item This equivalence sends $\IC_0$ to $\sh^{1/2} \co_{\ld{M}_{\Z'}}$. 
        \item Suppose $X = G/H$ for a spherical subgroup $H\subseteq G$. As in \cref{rmk: suggestion kostant for duals}, there is an isomorphism $\ld{M}_{\Z'} \mmod \ld{G}_{\Z'}(-2\rho) \cong \ld{\fr{h}}_{\Z'}^\ast(2)\mmod \ld{H}_{\Z'}$ and a ``Kostant section'' $\kappa_{\ld{M}}: \ld{\fr{h}}_{\Z'}^\ast(2)\mmod \ld{H}_{\Z'} \to \ld{M}_{\Z'}$ which makes the following square commute:
        $$\xymatrix{
            \ld{\fr{h}}_{\Z'}^\ast(2)\mmod \ld{H}_{\Z'} \ar[r]^-{\kappa_{\ld{M}}} \ar[d] & \ld{M}_{\Z'}/\ld{G}_{\Z'}(-2\rho) \ar[d]^-\mu \\
            \ld{\g}_{\Z'}^\ast(2)\mmod \ld{G}_{\Z'} \ar[r]^-{\kappa} & \ld{\g}_{\Z'}^\ast(2-2\rho)/\ld{G}_{\Z'}(-2\rho).
        }$$
    \end{enumerate}
    See also \cite[Section 5.3]{bzsv} for related discussion.
\end{assume}
\begin{recall}
    Fix the same localization $\Z'$ of $\Z$ as in \cref{assume: split form}. 
    If $X$ is a closed oriented $n$-manifold, the homology $\cC_\ast(\cL X; \Z') := \Sigma^{n} C_{\ast}(\cL X; \Z')$ admits the structure of an algebra via the Chas-Sullivan product. This product is given by push-pull along the span $\cL X \times \cL X \leftarrow \Map(S^1 \vee S^1, X) \to \cL X$. If $X = G/H$, one can identify
    $$\cC_\ast(\cL X; \Z') = C_\ast(\Omega X; \Z')^{hH} \otimes_{\Z^{'hG}} \Z'.$$
    Up to completion, $\cC_\ast(\cL X; \Z')$ can be identified with the Hochschild cohomology $\HC(C^\ast(X; \Z')/\Z')$; in particular, $\cC_\ast(\cL X; \Z')$ is naturally an $\E{2}$-$\Z'$-algebra by the Deligne conjecture. See \cite[Theorem 1.0.2]{brav-rozenblyum}. If $X$ is simply-connected, one can also view $\cC_\ast(\cL X; \Z)$ as the symplectic cohomology $\SH(T^\ast X; \Z)$ by Viterbo's theorem (see \cite{abouzaid-viterbo}).
\end{recall}
Of course, the reader only interested in string topology with rational coefficients need not use \cref{assume: split form} in any of this discussion.
\begin{theorem}\label{prop: string top of homogeneous space}
    Fix a spherical subgroup $H \subseteq G$. Let $X \simeq (G/H)(\cc)$ denote the associated compact homogeneous manifold for the maximal compact subgroup of $G(\cc)$.
    Under \cref{assume: split form}, let
    $$\ld{J}_X' = \ld{\fr{h}}_{\Z'}^\ast(2)\mmod \ld{H}_{\Z'} \times_{\ld{M}_{\Z'}/\ld{G}_{\Z'}(-2\rho)} \ld{\fr{h}}_{\Z'}^\ast(2)\mmod \ld{H}_{\Z'},$$
    so that there is a canonical morphism $\ld{J}_X' \to \ld{\fr{h}}_{\Z'}^\ast(2)\mmod \ld{H}_{\Z'}$, and hence a canonical composite
    $$\ld{J}_X' \to \ld{\fr{h}}_{\Z'}^\ast(2)\mmod \ld{H}_{\Z'} \to \ld{\fr{g}}_{\Z'}^\ast(2)\mmod \ld{G}_{\Z'}.$$
    Then there is an equivalence
    $$\cC_\ast(\cL X; \Z') \simeq \sh^{1/2} \Gamma\left(\ld{J}_X' \times_{\ld{\fr{g}}_{\Z'}^\ast(2)\mmod \ld{G}_{\Z'}} \{0\}; \co\right).$$
\end{theorem}
\begin{proof}
    It follows from \cref{assume: split form}(b) and (an integral variant of) \cref{thm: ordinary homology criterion satake} that $\sh^{1/2} \ld{J}_X' \simeq \spec C^H_\ast(\Omega X; \Z')$, so that
    \begin{align*}
        \sh^{1/2} \ld{J}_X' \times_{\ld{\fr{g}}_{\Z'}^\ast(2)\mmod \ld{G}_{\Z'}} \{0\} & \simeq \spec \left( C^H_\ast(\Omega X; \Z') \otimes_{C_G^\ast(\ast; \Z')} \Z' \right) \\
        & \simeq \spec \left( C_\ast(\Omega X; \Z')^{hH} \otimes_{\Z^{'hG}} \Z' \right),
    \end{align*}
    whose ring of functions identifies with $\cC_\ast(\cL X; \Z')$.
\end{proof}
Note that the full strength of \cref{assume: split form} is not really necessary; all one needs is an integral and derived refinement of the isomorphism $\ld{J}_X' \simeq \spec \H^H_\ast(\Omega X; \Z')$ from \cref{thm: ordinary homology criterion satake}.
\begin{remark}
    Combined with \cref{thm: ordinary homology criterion satake}, \cref{prop: string top of homogeneous space} says that if $G/H$ is simply-connected, one could view 
    \begin{equation}\label{eq: reg centr and G equiv symplectic}
        \H^H_\ast(\Omega(G/H); \Z) \cong \SH_G^\ast(T^\ast(G/H); \Z),
    \end{equation}
    the latter being some appropriate version of $G$-equivariant symplectic cohomology (i.e., a $G$-equivariant analogue of the string topology algebra $\cC_\ast(\cL(G/H); \Z)$). In particular, there should be a map of algebras
    $$\H^G_\ast(\Omega G; \Z) \to \SH_G^\ast(T^\ast(G/H); \Z)$$
    which is ``coisotropic'' as explained in \cref{cor: desired map is coisotropic}. See \cite{coulomb-branch-symplectic-cohomology} for such a construction with coefficients in $\cc$.
    There should also be an equivalence (motivated by \cref{thm: ordinary homology criterion satake})
    $$\spec \SH_G^\ast(T^\ast(G/H); \Z') \cong \ld{J}_X'.$$
    Moreover, the isomorphism \cref{eq: reg centr and G equiv symplectic} should hold much more generally with ordinary homology replaced by equivariant (connective) K-theory or equivariant elliptic cohomology, and equivariant symplectic cohomology replaced by the appropriate variant.
\end{remark}
\begin{example}[Group case]
    Suppose $G = H \times H$, with $H$ embedded diagonally. Then $X$ is the maximal compact subgroup of $H(\cc)$, so we will simply write $X = H$ for notational simplicity. We will also omit grading shifts. In this case, $\ld{M} = T^\ast \ld{H}$.
    Assume a version of \cref{thm: derived satake} with $\Z'$ coefficients. Then \cref{prop: string top of homogeneous space} says that there is an equivalence
    \begin{align*}
        \cC_\ast(\cL H; \Z') & \simeq \sh^{1/2} \Gamma \left( \ld{J}_{\ld{H}_{\Z'}} \times_{\ld{\fr{h}}_{\Z'}^\ast(2)\mmod \ld{H}_{\Z'} \times \ld{\fr{h}}_{\Z'}^\ast(2)\mmod \ld{H}_{\Z'}} \{0\}; \co \right) \\
        & \simeq \sh^{1/2} \Gamma \left( \left(\ld{J}_{\ld{H}_{\Z'}} \times_{\ld{\fr{h}}_{\Z'}^\ast(2)\mmod \ld{H}_{\Z'}} \{0\}\right) \times \left(\{0\} \times_{\ld{\fr{h}}_{\Z'}^\ast(2)\mmod \ld{H}_{\Z'}} \{0\}\right); \co\right) \\
        & \simeq \sh^{1/2} \Gamma \left( Z_e(\ld{H}_{\Z'})(-2\rho) \times (\ld{\fr{h}}_{\Z'}^\ast(2)\mmod \ld{H}_{\Z'})[-1]; \co\right)
    \end{align*}
    The final term is simply the shearing of the tensor product of $\co_{Z_e(\ld{H}_{\Z'})(-2\rho)} \cong \H_\ast(\Omega H; \Z')$ (this isomorphism being provided by \cref{prop: ginzburg coh} and \cref{rmk: yun zhu centralizer}) with an exterior algebra (which identifies with $\H^\ast(H; \Z')$ by the theory of the Kostant slice). One therefore recovers \cite[Theorem 1]{hepworth-string-topology-lie-groups}, which identifies $\pi_\ast \cC_\ast(\cL H; \Z')$ with $\H_\ast(\Omega H; \Z') \otimes \H^\ast(H; \Z')$.
\end{example}
\cref{prop: string top of homogeneous space} also says that the calculations in \cref{sec: case by case} can be used to recover the calculations of \cite{cohen-jones-yan}, and more generally of the string topology of compact Riemannian symmetric spaces $X$ of rank one. For simplicity, we will assume $X$ is simply-connected, so $X$ is $S^n$, $\CP^n$, $\HHP^n$, or $\OP^2$. These cases fall under the purview of the following lemma.
\begin{lemma}\label{lem: derived quotient}
    Let $\Z[x,b]/bx^k$ denote the graded ring with $x$ in weight $-2i$ and $b$ in weight $2j$. Then the quotient $(\Z[x,b]/bx^k)/x^{k+\ell}$ has homotopy groups given by
    $$\pi_\ast (\Z[x,b]/bx^k)/x^{k+\ell} \cong \Z[x,b,\sigma(bx^{k+\ell})]/(x^{k+\ell}, bx^k, \sigma(bx^{k+\ell})^2, x^k \sigma(bx^{k+\ell})),$$
    where $\sigma(bx^{k+\ell})$ lives in weight $2(j-(k+\ell)i)$ and degree $1$. In particular, there is a graded isomorphism
    $$\pi_\ast \sh^{1/2} (\Z[x,b]/bx^k)/x^{k+\ell} \cong \Z[x,b,\sigma(bx^{k+\ell})]/(x^{k+\ell}, bx^k, \sigma(bx^{k+\ell})^2, x^k \sigma(bx^{k+\ell})),$$
    where $\sigma(bx^{k+\ell})$ lives in weight $2(j-(k+\ell)i)+1$.
\end{lemma}
\begin{proof}
    Let $n \geq k+\ell$; then, $bx^n = (bx^k)x^{n-k} = b(x^n)$ admits two nullhomotopies in this quotient, so we obtain a class $\sigma(bx^n)$ which lives in weight $2(j-ni)$ and degree $1$; this class lies in the Massey product/Toda bracket $\pdb{x^{n-k}, x^k, b}$. Note that $\sigma(bx^n) = x^{n-k-\ell} \sigma(bx^{k+\ell})$. The relation $x^{k+\ell} = 0$ implies that $x^{n-k} = 0$ for $n-k\geq k+\ell$; so $x^k \sigma(bx^{k+\ell}) = 0$.
\end{proof}
\begin{example}\label{ex: string top for spheres}
    Let $\Z' = \Z[\tfrac{1}{2}]$, let $G = \SO_{n+1}$, and let $H = \SO_n$. Using the calculations in the proof of \cref{thm: bzsv for Bn} and \cref{thm: bzsv for Dn} (or rather, its variant for $\SO_{n+1}/\SO_n$), we obtain isomorphisms
    \begin{align*}
        \ld{J}_{\SO_{n+1}/\SO_n} & \cong \spec \H_\ast^{\SO_n}(\Omega S^n; \Z') \\
        & \cong \begin{cases}
            \spec \Z'\left[p_1, \cdots, p_j, (a + a^{-1})^2, \tfrac{(a-a^{-1})^2}{p_j}, \tfrac{a^2-a^{-2}}{p_j^{1/2}}\right] & n = 2j+1, \\
            \spec \Z'[p_1, \cdots, p_{j-1}, c_j, b]/bc_j & n = 2j;
        \end{cases}
    \end{align*}
    here, $a$ lives in weight $0$, and $b$ lives in weight $4j-2$. See also \cite[Remark B.4]{grg-reg} for the first isomorphism. This scheme lives over $\ld{\g}_{\Z'}^\ast(2)\mmod \ld{G}_{\Z'}$, and we find:
    \begin{itemize}
        \item If $n = 2j+1$, the map $\ld{\fr{h}}_{\Z'}^\ast(2)\mmod \ld{H}_{\Z'} \to \ld{\g}_{\Z'}^\ast(2)\mmod \ld{G}_{\Z'}$ is given by the map
        $$\Z'[p_1,\cdots,p_j, c_{j+1}] \to \Z'[p_1,\cdots,p_j], \ c_{j+1} \mapsto 0.$$
        It follows from \cref{prop: string top of homogeneous space} that
        \begin{align*}
            \cC_\ast(\cL S^{2j+1}; \Z') & \simeq \sh^{1/2}\left(\Z'\left[p_1, \cdots, p_j, (a + a^{-1})^2, \tfrac{(a-a^{-1})^2}{p_j}, \tfrac{a^2-a^{-2}}{p_j^{1/2}}\right]\right)/(p_1, \cdots, p_j, c_{j+1}) \\
            & \simeq \sh^{1/2}\left(\Z'\left[(a + a^{-1})^2, \tfrac{(a-a^{-1})^2}{p_j}, \tfrac{a^2-a^{-2}}{p_j^{1/2}}\right]\right)/(p_j, c_{j+1}).
        \end{align*}
        Though this looks complicated, the relations simplify dramatically: we find that the homotopy of this algebra is simply
        $$\pi_\ast \cC_\ast(\cL S^{2j+1}; \Z') \cong \Z'\left[\tfrac{a^2-a^{-2}}{p_j^{1/2}}, \sigma(c_{j+1})\right]/\sigma(c_{j+1})^2,$$
        where $\tfrac{a^2-a^{-2}}{p_j^{1/2}}$ lives in weight $2j$ and $\sigma(c_{j+1})$ (coming from the two nullhomotopies of $c_{j+1}$) lives in weight $-2(j+1)+1 = -2j-1$. As expected, this is precisely \cite[Theorem 2(1)]{cohen-jones-yan}.
        \item If $n = 2j$, the map $\ld{\fr{h}}_{\Z'}^\ast(2)\mmod \ld{H}_{\Z'} \to \ld{\g}_{\Z'}^\ast(2)\mmod \ld{G}_{\Z'}$ is given by the map
        $$\Z'[p_1,\cdots,p_j] \to \Z'[p_1,\cdots,p_{j-1},c_j], \ p_j \mapsto c_j^2.$$
        It follows from \cref{prop: string top of homogeneous space} that 
        \begin{align*}
            \cC_\ast(\cL S^{2j}; \Z') & \simeq \sh^{1/2}(\Z'[p_1, \cdots, p_{j-1}, c_j, b]/bc_j)/(p_1, \cdots, p_{j-1}, c_j^2) \\
            & \simeq \sh^{1/2}(\Z'[c_j, b]/bc_j)/c_j^2.
        \end{align*}
        It follows from \cref{lem: derived quotient} that there is a graded isomorphism
        $$\pi_\ast \cC_\ast(\cL S^{2j}; \Z') \cong \Z'[c_j, b, \sigma(bc_j^2)]/(bc_j, c_j^2, \sigma(bc_j^2)^2, c_j \sigma(bc_j^2)),$$
        where $\sigma(bc_j^2)$ lives in weight $-1$.
        This is precisely \cite[Theorem 2(2)]{cohen-jones-yan}, as expected.
    \end{itemize}
\end{example}
\begin{example}\label{ex: string top for CPn}
    Let $\Z' = \Z[\frac{1}{n+1}]$. Using the calculations in the proof of (an integral version of) \cref{thm: bzsv for An}, we obtain isomorphisms
    \begin{align*}
        \ld{J}_{\PGL_{n+1}/\GL_n} & \cong \spec \H_\ast^{\GL_n}(\Omega \CP^n; \Z') \\
        & \cong \spec \Z'[c_1, \cdots, c_{n-1}, c_n, b]/bc_n.
    \end{align*}
    Here, $b$ lives in weight $2n$. This scheme lives over $\ld{\g}_{\Z'}^\ast(2)\mmod \ld{G}_{\Z'}$. Taking its fiber over the origin kills the ideal $(c_2, \cdots, c_{n-1}, c_1^n-c_n, c_1 c_n)$. Note that $c_1 c_n = c_1^{n+1}$. It follows from \cref{prop: string top of homogeneous space} and \cref{lem: derived quotient} that there is a graded isomorphism
    $$\pi_\ast \cC_\ast(\cL \CP^n; \Z') \cong \Z'[c_1, b, \sigma(bc_1^{n+1})]/(bc_1^n, c_1^{n+1}, \sigma(bc_1^{n+1})^2, c_1^n \sigma(bc_1^{n+1})),$$
    where $\sigma(bc_1^{n+1})$ lives in weight $-1$. As expected, this is \cite[Theorem 3]{cohen-jones-yan}.
\end{example}
\begin{remark}
    \cref{lem: derived quotient} shows that using \cref{thm: bzsv for Cn} and \cref{thm: bzsv for F4}, one obtains graded isomorphisms
    \begin{align*}
        \pi_\ast \cC_\ast(\cL \HHP^n; \Z') & \cong \Z'[p_1, b, \sigma(bp_1^{n+1})]/(bp_1^n, p_1^{n+1}, \sigma(bp_1^{n+1})^2, p_1^n \sigma(bp_1^{n+1})), \\
        \pi_\ast \cC_\ast(\cL \OP^2; \Z') & \cong \Z'[p_2, b, \sigma(bp_2^3)]/(bp_2^2, p_2^3, \sigma(bp_2^3)^2, p_2^2 \sigma(bp_2^3)).
    \end{align*}
    Here, $b$ lives in weight $4n+2$ in the first line, and in weight $22$ in the second line. In the first line, one can take $\Z' = \Z[\tfrac{1}{n+1}]$; and in the second line, one can take $\Z' = \Z[\tfrac{1}{3}]$; see \cite{string-top-quaternion}.\footnote{These denominators precisely encode the ``very good primes'' for $\ld{G}$; for instance, a prime $p$ is very good for $\PGL_{n+1}$ if $p\nmid (n+1)$, and is very good for $\F_4$ if $p\nmid 3$. If a prime $p$ is very good for $\ld{G}$, it admits a quasi-logarithm $\ld{G} \to \ld{\g}$ in the sense of \cite{kazhdan-varshavsky-quasilog} (see \cite[Corollary 6.4]{friedlander-negron}), which should imply that the usual statement of \cref{thm: derived satake} holds with coefficients in $\FF_p$. See \cref{rmk: remarks about derived satake}.}
\end{remark}

%% file: appendices/open-qns.tex
\section{Questions/further directions}

The work presented in this article is clearly far from being a complete story. There are numerous questions left open by our discussion, of varying levels of difficulty. I hope that the discussion in this article is compelling enough to motivate further study of these problems. 

Let us begin with some broad questions. First, some conceptual questions about the place of this story in mathematics and physics:
\begin{enumerate}
    \item Having been raised a topologist, it is inspiring to see the appearance of \textit{unstable} homotopy-theoretic structures in \cref{table: topology and dualities for rank 1 spherical varieties} (such as EHP sequences and Hopf fibrations). It suggests that there might be an unstable analogue of \cref{expect: sphere satake} (the latter is already not well-defined!) giving a spectral decomposition for ``$G$-equivariant sheaves of \textit{spaces} on $\Omega G$''. There are some indications that this might be possible, and many that it might be impossible --- but it is perhaps better to be optimistic!
    
    This question seems exceptionally difficult, even in the simplest case of $G$ being trivial. In this case, one can interpret the question as asking for an analogue of the theory of synthetic spectra \`a la \cite{piotr-synthetic} for spaces/homotopy types/anima. Considerations from the theory of power operations suggests that one should replace the fpqc stack $\Msfg$ of groupoids classifying spin formal groups and isomorphisms between them by some sort of fpqc stack of categories classifying spin formal groups and \textit{homomorphisms} (in particular, including isogenies) between them.
    \item As mentioned in \cite[Heuristic B.5]{grg-reg}, an unpublished conjecture of Gaiotto (which I learned about from Nakajima) says that the Coulomb branch of $4$d $\cN = 2$ pure gauge theory over $\RR^3 \times S^1$ with a generic choice of complex structure and gauge group $G$ can be modeled via the periodic $G$-equivariant complex K-theory of $\Omega G$. The calculations of this paper suggest that perhaps one can modify this proposal to use \textit{connective} equivariant complex K-theory instead, and that the Bott class is related to the radius of the circle $S^1$ (see also \cref{rmk: ku-ham and radius}). The Bott class being sent to zero then might correspond to the degeneration of $4$d $\cN = 2$ pure gauge theory into $3$d $\cN = 4$ pure gauge theory. For instance, in \cref{rmk: S1rot and torus}, we find that if $T = S^1$ (say), going from integral to $\ku$-homology has the effect of \textit{deforming} the relation $[p,a] = \hbar$ in the Weyl algebra of $\GG_a = \spec \Z[a]$ over $\Z$ to the relation $[p,a] = \hbar(1 + \beta a p)$ over $\Z[\beta]$. (The symbol $p$ corresponds to $\tfrac{1}{a} x$ in \cref{rmk: S1rot and torus}.) Does such a deformation have any precedence in physics?
    
    Here is a more optimistic question. Kapustin and Witten (among others) proposed viewing the derived geometric Satake equivalence of \cref{thm: derived satake} as an equivalence of categories stemming from S-duality for $4$d $\cN=4$ gauge theory. Does \cref{thm: ku derived satake} have any relationship to $5$d maximally supersymmetric Yang-Mills theory compactified on a circle of finite radius (related to the Bott class $\beta$)?
    \item One glaring omission in this article is the study of equivariance for loop rotation on the topological/A-side; as explained in \cite{nek-shat, ben-zvi-susy} via the $\Omega$-deformation and \cite[Remark B.2]{grg-reg}, this corresponds to deformation quantizing the spectral/B-side. From one point of view, it can be described as the study of the homology of the operad $\E{3} \rtimes S^1$ associated to rotations about a line inside $\RR^3$. Is there a good theory of quantizations of $\ku$-Hamiltonian spaces (and a classification thereof)?
    
    Motivated by this, and considerations involving Frobenius-constant quantizations \`a la Bezrukavnikov-Kaledin \cite{bezrukav-kaledin}, we are planning to discuss the $\U(n)$-equivariant cohomology of the $\E{2n+1}$-operad with \textit{arbitrary} complex-oriented coefficients in a later article \cite{Eodd-and-quantizations}. This discussion is heavily motivated by \cite{generalized-n-series} (in particular, \cite[Section 4.4]{generalized-n-series} can be regarded as the simplest case of this story, namely the case of tori). Part of the reason for not discussing loop rotation equivariance in this article is our reliance on \cref{prop: going from G of LG/H to H of Loops G/H}, which inherently breaks the natural $S^1$-action on the free loop space.
    \item Prove \cref{conj: generalized kostant slice}. In every example I have studied, the construction of the relevant Kostant slice $\kappa_{\ld{M}}: \ld{M}\mmod \ld{G} \to \ld{M}/\ld{G}$ relies on classical invariant theory, and proving a unified statement would be very interesting and important. One interesting feature in the known examples is that the construction of such slices involves constructing a ``stronger'' slice $\ld{M}\mmod \ld{G} \to \ld{M}$.  However, this depends on some noncanonical choices, i.e., the map $\ld{M}\mmod \ld{G} \to \ld{M}$ is not canonical. It would therefore also be interesting to show that $\ld{G}$ acts transitively on the set of all such slices $\ld{M} \mmod \ld{G} \to \ld{M}$ (in other words, that such a slice, if it exists, is \textit{unique} if we ask it to land in $\ld{M}/\ld{G}$ instead of $\ld{M}$). For the usual Kostant slice $\ld{\g}\mmod \ld{G} \to \ld{\g}$, this was shown in \cite{friedman-morgan-kostant-slice}.
    \item Give an explanation for the observation in \cref{rmk: minimal nilpotent orbits}.
    \item Can one prove \textit{ramified} versions of the local equivalences studied in this article? See \cref{conj: iwahori-ramification} for the analogue of \cite[Conjecture 7.5.1]{bzsv} with tame ramification.
    \item In many calculations, one gets the sense that the condition of sphericity should not be crucial in relative geometric Langlands. For instance, the criteria (a) and (b) of \cref{thm: ordinary homology criterion satake} do not use sphericity of $H_\cc \subseteq G_\cc$ at all (and as in \cref{rmk: a composite map in homology}, \cref{rmk: expected G mod T}, and \cref{conj: spherical levi}, there should be many interesting examples of relative Langlands duality for non-spherical subgroups). One of the immediate difficulties encountered when working with non-spherical subgroups is that $G_\cc(\cc\pw{t})$-orbits on $G_\cc(\cc\ls{t})/H_\cc(\cc\ls{t})$ are not parametrized by a discrete set; so it is hard to do ``combinatorics''. Nevertheless, it would be interesting to collect other non-spherical examples satisfying conditions (a) and (b) of \cref{thm: ordinary homology criterion satake}; a sufficient supply of examples might suggest a way to understand this generalization of relative Langlands duality.
    \item In the proof (from \cite{bf-derived-satake}) of \cref{thm: derived satake}, many calculations can be reduced to the case of semisimple rank $1$, essentially by localization on the affine space $\ld{\fr{t}}^\ast(2)\mmod W = \spec \H^\ast_G(\ast; \Z')$. Is it possible to tackle \cref{conj: bzsv} for affine homogeneous spherical varieties of higher rank by reduction to the rank $1$ case (in which case \cref{thm: rk 1 bzsv is true} can be applied)? The technique of \textit{localization} of spherical varieties, studied in \cite{luna-big-cells, knop-localization-spherical}, should be crucial here.
    \item Let $C_2$ be the cyclic group of order $2$. In the case of symmetric spaces, \cref{conj: bzsv} should be closely related to $C_2$-equivariant derived algebraic geometry as studied by Mike Hill and his collaborators (see, e.g., \cite{mike-hill-nsf-2014}). Let $\sigma$ denote the sign representation of $C_2$ on $\RR$, and let $\varrho = 1 + \sigma$ denote the regular representation of $C_2$ on $\cc$.\footnote{The notation is intended to distinguish this $C_2$ from the $\Z/2$ appearing elsewhere in this article, as well as $\varrho$ from the half-sum of positive roots (which is usually denoted $\rho$). On the topic of notation: it would have been great if $\varrho$ was used to denote the \textit{r(h)o}tation representation of $C_2$ (which is frequently denoted by $\lambda = 2\sigma$), instead of the regular representation of $C_2$.} If $Y$ is a $C_2$-space, let $\cL^\sigma Y$ denote the space of maps $\Map(S^\sigma, Y)$ equipped with its natural $C_2$-action.
    
    Let $G$ be a connected compact Lie group, let $\theta$ be a conjugate-linear involution on $G_\cc$ preserving $G$, and let $G^\theta$ denote the maximal compact subgroup of the fixed subgroup $G_\cc^\theta$.
    Then, there is a $C_2$-equivariant homotopy equivalence of orbispaces
    $$G\backslash \cL^\sigma G/G \simeq \Bun_{(G,\theta)}(\CP^1_\RR),$$
    which gives an equivalence
    $$(G\backslash \cL^\sigma G/G)^{C_2} \simeq G\backslash \cL(G/G^\theta).$$
    Note that $\CP^1_\RR = S^\varrho$ (which one can think of as the unit sphere in $\varrho \oplus \RR = \sigma \oplus \RR^2$) is the standard real structure on $\CP^1$.

    Let $R$ denote a normed algebra in genuine $C_2$-spectra, and let $\Phi^{C_2} R$ denote its geometric fixed points (for instance, if $R = \ul{\Z}$ is the constant Mackey functor, $\Phi^{C_2} R$ is the connective cover of the Tate construction $\FF_2^{tS^1} \simeq \Z^{tC_2}$). We expect that there is an $\infty$-category $\Shv^{c,\Sat}(\Bun_{(G,\theta)}(\CP^1_\RR); R)$ of ``$C_2$-equivariant sheaves of $R$-modules on $\Bun_{(G,\theta)}(\CP^1_\RR)$''. This $\infty$-category should satisfy an analogue of \cref{thm: derived satake}. Taking geometric fixed points (in some appropriate categorical sense) of this putative spectral decomposition of $\Shv^{c,\Sat}(\Bun_{(G,\theta)}(\CP^1_\RR); R)$ should produce an analogue/special case of \cref{conj: bzsv} which describes $\Shv^{c,\Sat}_{G\pw{t}}(G\ls{t}/G^\theta\ls{t}; \Phi^{C_2} R)$. For instance, one can try to access $\Shv^{c,\Sat}_{G\pw{t}}(G\ls{t}/G^\theta\ls{t}; \Phi^{C_2} R)$ via an analogue of Smith theory. The case of $G = \GL_n$ equipped with the orthogonal involution, and $R = \ul{\Z}$, is especially interesting, and we hope to study this question in future work. Developing this theory should lead to important interactions between geometric Langlands for symmetric spaces and $C_2$-equivariant homotopy theory.

    \item What is the underlying reason/broader context for the similarities evident in \cref{table: topology and dualities for rank 1 spherical varieties} between the dual varieties to the rank $1$ affine homogeneous varieties? In our story, these reflect the similitude between the based loop spaces of odd/even-dimensional spheres, but there are \textit{many} other examples of such similarities in higher rank. We will address this question in future work; see \cref{rmk: similarities between spheres} for a brief comment. The recent article \cite{leslie-endoscopic-symmetric-varieties} by Leslie is also concerned with this topic.
\end{enumerate}
Some of the most important issues in this article on the topological/A-side come from defining well-behaved sheaf theories.
\begin{enumerate}[resume]
    \item Verify \cref{hypothesis: rank 1 weakly placid}, and more generally \cref{conj: placidity}.
    \item Let $G$ be a compact Lie group, and let $X$ be an (ind-)finite stratified $G$-space. Is there a good $\sh^{1/2} \Z[\beta]$-linear $\infty$-category $\Shv_G^c(X; \ku)$ which agrees with the \textit{ad hoc} construction of \cref{cstr: Shv-Sat LG/H with ku}?
    \item Is there a good sheaf theory of synthetic $G$-equivariant constructible sheaves of spectra as in \cref{expect: shvsyn G equiv}? As a first step, it seems important to study $G$-equivariant analogues of synthetic spectra, as well as nonequivariant synthetic analogues of constructible sheaves of spectra.
\end{enumerate}
The discussion in this article also suggests several (less lofty, and presumably more approachable) questions.
\begin{enumerate}[resume]
    \item Can one define $\ld{G}_{\GG_\univ}$ as in \cref{expect: univ probing}?
    \item Is there an \textit{a priori} reason that if $\ld{M}_\beta$ is an $\ld{G}$-variety over $\Z[\beta]$ as in \cref{prop: ku rk 1 reg centr --> thm}, the map $\ld{M}_\beta/\ld{G} \to \ld{G}^\mathrm{sc}_\beta/\ld{G}$ admits a Lagrangian structure? See \cref{cor: desired map is coisotropic} for some progress in this direction.
    \item As in \cref{rmk: general map complex reflection}, let $W_1 \to W_2$ be a homomorphism of reflection groups acting on vector spaces $V_1 \to V_2$ over a field $k$ (possibly of nonzero characteristic), so that there is a map $V_1\mmod W_1 \to V_2\mmod W_2$. The Hochschild homology $\HH(V_1\mmod W_1 / V_2\mmod W_2)$ should be an interesting invariant associated to homomorphisms of reflection groups; what can one prove about it? In the case that the map of reflection groups comes from an inclusion of root data, this Hochschild homology is the content of \cref{prop: equiv coh as hochschild}, and therefore plays an important role in the relative Langlands program.
    \item Is there an analogue of \cref{thm: bzsv for Dn} for $\ku$? The difficulty lies in finding a $\beta$-deformation of $\ld{M}$.
    \item Can one extend \cref{thm: ku derived satake} to the non-simply-laced case? Similarly, can one extend \cref{prop: ginzburg coh} to coefficients with the sphere spectrum? Following \cref{lem: sphere and Loops SUn}, a first step will be an understanding of $\MU_\ast(\Omega G)$ as a Hopf algebra; the cobordism groups are known by \cite{petrie-MU-lie-group}. A description of the $\MU_\ast \MU$-comodule structure would then give the desired generalization of \cref{prop: ginzburg coh}. 
    \item Fix a prime $p$. Can the results of this article (even only the ones concerning coefficients in ordinary commutative rings) be extended to the setting of $p$-compact groups? This question was suggested by Haynes Miller. See \cite{grodal-icm} for a survey of the theory of $p$-compact groups.

    Let us make some brief comments about this question. First, if the input into the machinery of geometric Langlands duality is the \textit{$p$-completion} of a compact torus $T$, the ``Langlands dual'' appears to be the $p$-divisible subgroup $\ld{T}[p^\infty]$ of the Langlands dual torus $\ld{T}$. (In the rank one case, this reduces to the statement that the Cartier dual of $\Omega (S^1)^\wedge_p = \Z_p$, if interpreted suitably, is $\mu_{p^\infty}$.) Second, the simplest exotic case of the Sullivan sphere $(S^{2n-1})^\wedge_p = \Omega (B^2\Z_p)_{h\Cn}$ for $n|\#\FF_p^\times$, which behaves like a compact Lie group of rank $1$, already seems to exhibit an interesting Langlands duality: its ``dual group'' is built out of the $p$-divisible torus $\mu_{p^\infty}$ and the Weyl group $\Cn$ acting via $\Cn \subseteq \Z_p^\times \circlearrowright \mu_{p^\infty}$.
    Third, in \cite[Remark 3.13]{triple-product-cayley}, I suggest that the Dwyer-Wilkerson exotic $2$-compact group from \cite{dwyer-wilkerson-2-compact} might act on a framed $30$-manifold with Kervaire invariant $1$ (just as $\PGL_2$ acts on $\PGL_2^{\times 3}/\PGL_2 \cong \RP^3 \times \RP^3$ and $\G_2$ acts on $\SO_8/\mu_2 \cdot \G_2 \cong \RP^7 \times \RP^7$), and that the ``regular centralizer'' group scheme for this action should be related to the relative Langlands program. See \cite{triple-product-cayley} for (brief) further discussion.
    \item Finding other examples of the $\ku$-theoretic analogue of \cite[Conjecture 7.5.1]{bzsv} is an important goal. For instance, what is the $\ku$-theoretic dual of the spherical $\GL_n$-variety $\GL_n/(\GL_j \times \GL_{n-j})$ (which is homotopy equivalent to $\Gr_j(\cc^n)$)? In the arithmetic case, this was studied by Jacquet-Rallis in \cite{jacquet-rallis}, and is described in \cref{conj: spherical levi} in the geometric case.
    \item Is there an analogue of the theory of synthetic spectra for global homotopy theory?
    \item Is there a structure theory for $\ku$-Hamiltonian spaces? Namely, is there a well-behaved theory of ``hyperspherical'' $\ku$-Hamiltonian spaces which generalizes the notion of hyperspherical Hamiltonian spaces introduced in \cite[Section 3]{bzsv}, and an analogue of \cite[Theorem 3.6.1]{bzsv}?
    \item Related to the preceding point, can one develop a theory of $\beta$-deformed cotangent bundles for certain $G$-varieties? In the case of ``coisotropic subgroups'' as in \cref{rmk: subgroup coisotropic}, this would be a $\beta$-deformation of \cite{balibanu-mayrand} (or rather, of the specialization of their results to quasi-Hamiltonian spaces).
\end{enumerate}